\documentclass[10pt]{amsbook}
\usepackage{float}
\usepackage{amsmath}
\usepackage{amssymb,hyperref}
\usepackage{amsfonts}

\theoremstyle{plain}

\newtheorem{corollary}{Corollary}

\newtheorem{definition}{Definition}

\newtheorem{lemma}{Lemma}

\newtheorem{proposition}{Proposition}
\newtheorem{remark}{Remark}

\newtheorem{theorem}{Theorem}
\numberwithin{equation}{chapter} \numberwithin{section}{chapter}
\input{tcilatex}
\makeindex

\begin{document}
\frontmatter
\title[Advances on Inequalities of the Schwarz, Triangle and Heisenberg Type]%
{\textbf{\ Advances on Inequalities of the Schwarz, Triangle and Heisenberg
Type in Inner Product Spaces }}
\author[S.S. Dragomir]{Sever Silvestru Dragomir}
\address{School of Computer Science \& Mathematics\\
Victoria University, Melbourne\\
Victoria, Australia}
\email{sever.dragomir@vu.edu.au}
\urladdr{http://rgmia.vu.edu.au/SSDragomirWeb.html}

\begin{abstract}
The purpose of this book is to give a comprehensive introduction to several
inequalities in Inner Product Spaces that have important applications in
various topics of Contemporary Mathematics such as: Linear Operators Theory,
Partial Differential Equations, Nonlinear Analysis, Approximation Theory,
Optimization Theory, Numerical Analysis, Probability Theory, Statistics and
other fields.
\end{abstract}

\maketitle
\tableofcontents

\chapter*{Preface}

The purpose of this book, that can be seen as a continuation of the previous
one entitled "\textit{Advances on Inequalities of the Schwarz, Gr\"{u}ss and
Bessel Type in Inner Product Spaces}" (Nova Science Publishers, NY, 2005),
is to give a comprehensive introduction to other classes of inequalities in
Inner Product Spaces that have important applications in various topics of
Contemporary Mathematics such as: Linear Operators Theory, Partial
Differential Equations, Nonlinear Analysis, Approximation Theory,
Optimization Theory, Numerical Analysis, Probability Theory, Statistics and
other fields.

The monograph is intended for use by both researchers in various fields of
Mathematical Inequalities, domains which have grown exponentially in the
last decade, as well as by postgraduate students and scientists applying
inequalities in their specific areas.

The aim of Chapter \ref{ch1} is to present some fundamental analytic
properties concerning Hermitian forms defined on real or complex linear
spaces. The basic inequalities as well as various properties of
superadditivity and monotonicity for the diverse functionals that can be
naturally associated with the quantities involved in the Schwarz inequality
are given. Applications for orthonormal families, Gram determinants, linear
operators defined on Hilbert spaces and sequences of vectors are also
pointed out.

In Chapter \ref{ch2}, classical and recent refinements and reverse
inequalities for the Schwarz and the triangle inequalities are presented.
Further on, the inequalities obtained by Buzano, Richards, Precupanu and
Moore and their extensions and generalizations for orthonormal families of
vectors in both real and complex inner product spaces are outlined. Recent
results concerning the classical refinement of Schwarz inequality due to
Kurepa for the complexification of real inner product spaces are also
reviewed. Various applications for integral inequalities including a version
of Heisenberg inequality for vector valued functions in Hilbert spaces are
provided as well.

The aim of Chapter \ref{ch3} is to survey various recent reverses for the
generalised triangle inequality in both its simple form, that are closely
related to the Diaz-Metcalf results, or in the equivalent quadratic form
that maybe be of interest in the Geometry of Inner product Spaces.
Applications for vector valued integral inequalities and for complex numbers
are given as well.

Further on, in Chapter \ref{ch4}, some recent reverses of the continuous
triangle inequality for Bochner integrable functions with values in Hilbert
spaces and defined on a compact interval $\left[ a,b\right] \subset \mathbb{R%
}$ are surveyed. Applications for Lebesgue integrable complex-valued
functions that generalise and extend the classical result of Karamata are
provided as well.

In Chapter \ref{ch5} some reverses of the Cauchy-Buniakovsky-Schwarz
vector-valued integral inequalities under various assumptions of boundedness
for the functions involved are given. Natural applications for the
Heisenberg inequality for vector-valued functions in Hilbert spaces are also
provided.

The last chapter, Chapter \ref{ch6}, is a potpourri of other inequalities in
inner product spaces. The aim of the first section is to point out some
upper bounds for the distance $d\left( x,M\right) $ from a vector $x$ to a
finite dimensional subspace $M$ in terms of the linearly independent vectors 
$\left\{ x_{1},\dots ,x_{n}\right\} $ that span $M$. As a by-product of this
endeavour, some refinements of the generalisations for Bessel's inequality
due to several authors including: Boas, Bellman and Bombieri are obtained.
Refinements for the well known Hadamard's inequality for Gram determinants
are also derived.

In the second and third sections of this last chapter, several reverses for
the Cauchy-Bunyakovsky-Schwarz (CBS) inequality for sequences of vectors in
Hilbert spaces are obtained. Applications for bounding the distance to a
finite-dimensional subspace and in reversing the generalised triangle
inequality are also given.

For the sake of completeness, all the results presented are completely
proved and the original references where they have been firstly obtained are
mentioned. The chapters are relatively independent and can be read
separately.

\ 

\textit{The Author},

March, 2005.

\mainmatter

\pagenumbering{arabic}


\chapter{Inequalities for Hermitian Forms}\label{ch1}

\section{Introduction}

Let $\mathbb{K}$ be the field of real or complex numbers, i.e., $\mathbb{K=R}
$ or $\mathbb{C}$ and $X$ be a linear space over $\mathbb{K}$.

\begin{definition}
\label{d1.1.1}A functional $\left( \cdot ,\cdot \right) :X\times
X\rightarrow \mathbb{K}$ is said to be a Hermitian form on $X$ if

\begin{enumerate}
\item[(H1)] $\left( ax+by,z\right) =a\left( x,z\right) +b\left( y,z\right) $
for $a,b\in \mathbb{K}$ and $x,y,z\in X;$

\item[(H2)] $\left( x,y\right) =\overline{\left( y,x\right) }$ for all $%
x,y\in X.$
\end{enumerate}
\end{definition}

The functional $\left( \cdot ,\cdot \right) $ is said to be \textit{positive
semi-definite }on a subspace $Y$ of $X$ if

\begin{enumerate}
\item[(H3)] $\left( y,y\right) \geq 0$ for every $y\in Y,$
\end{enumerate}

and \textit{positive definite} on $Y$ if it is positive semi-definite on $Y$
and

\begin{enumerate}
\item[(H4)] $\left( y,y\right) =0,$ $y\in Y$ implies $y=0.$
\end{enumerate}

The functional $\left( \cdot ,\cdot \right) $ is said to be \textit{definite
}on $Y$ provided that either $\left( \cdot ,\cdot \right) $ or $-\left(
\cdot ,\cdot \right) $ is positive semi-definite on $Y.$

When a Hermitian functional $\left( \cdot ,\cdot \right) $ is
positive-definite on the whole space $X,$ then, as usual, we will call it an
\textit{inner product }on $X$ and will denote it by $\left\langle \cdot
,\cdot \right\rangle .$

The aim of this chapter is to present some fundamental analytic properties
concerning Hermitian forms defined on real or complex linear spaces. The
basic inequalities as well as various properties of superadditivity and
monotonicity for diverse functionals that can be naturally associated with
the quantities involved in the Schwarz inequality are given. Applications
for orthonormal families, Gram determinants, linear operators defined on
Hilbert spaces and sequences of vectors are also pointed out. The results
are completely proved and the original references where they have been
firstly obtained are mentioned.

\section{Hermitian Forms, Fundamental Properties\label{s1}}

\subsection{\label{ss1.1}Schwarz's Inequality}

We use the following notations related to a given Hermitian form $\left(
\cdot ,\cdot \right) $ on $X:$%
\begin{align*}
X_{0}& :=\left\{ x\in X|\left( x,x\right) =0\right\} , \\
K& :=\left\{ x\in X|\left( x,x\right) <0\right\}
\end{align*}%
and, for a given $z\in X,$%
\begin{equation*}
X^{\left( z\right) }:=\left\{ x\in X|\left( x,z\right) =0\right\} \quad
\text{and\quad }L\left( z\right) :=\left\{ az|a\in \mathbb{K}\right\} .
\end{equation*}

The following fundamental facts concerning Hermitian forms hold \cite{IHFKU}:

\begin{theorem}[Kurepa, 1968]
\label{t1.1.1}Let $X$ and $\left( \cdot ,\cdot \right) $ be as above.

\begin{enumerate}
\item If $e\in X$ is such that $\left( e,e\right) \neq 0,$ then we have the
decomposition%
\begin{equation}
X=L\left( e\right) \bigoplus X^{\left( e\right) },  \label{1.1.1}
\end{equation}%
where $\bigoplus $ denotes the direct sum of the linear subspaces $X^{\left(
e\right) }$ and $L\left( e\right) ;$

\item If the functional $\left( \cdot ,\cdot \right) $ is positive
semi-definite on $X^{\left( e\right) }$ for at least one $e\in K,$ then $%
\left( \cdot ,\cdot \right) $ is positive semi-definite on $X^{\left(
f\right) }$ for each $f\in K;$

\item The functional $\left( \cdot ,\cdot \right) $ is positive
semi-definite on $X^{\left( e\right) }$ with $e\in K$ if and only if the
inequality%
\begin{equation}
\left\vert \left( x,y\right) \right\vert ^{2}\geq \left( x,x\right) \left(
y,y\right)  \label{1.1.2}
\end{equation}%
holds for all $x\in K$ and all $y\in X;$

\item The functional $\left( \cdot ,\cdot \right) $ is semi-definite on $X$
if and only if the Schwarz's inequality%
\begin{equation}
\left\vert \left( x,y\right) \right\vert ^{2}\leq \left( x,x\right) \left(
y,y\right)  \label{1.1.3}
\end{equation}%
holds for all $x,y\in X;$

\item The case of equality holds in (\ref{1.1.3}) for $x,y\in X$ and in (\ref%
{1.1.2}), for $x\in K,$ $y\in X,$ respectively; if and only if there exists
a scalar $a\in \mathbb{K}$ such that%
\begin{equation*}
y-ax\in X_{0}^{\left( x\right) }:=X_{0}\cap X^{\left( x\right) }.
\end{equation*}
\end{enumerate}
\end{theorem}

\begin{proof}
We follow the argument in \cite{IHFKU}.

If $\left( e,e\right) \neq 0,$ then the element%
\begin{equation*}
x:=y-\frac{\left( y,e\right) }{\left( e,e\right) }e
\end{equation*}%
has the property that $\left( x,e\right) =0,$ i.e., $x\in X^{\left( e\right)
}.$ This proves that $X$ is a sum of the subspaces $L\left( e\right) $ and $%
X^{\left( s\right) }.$ The fact that the sum is direct is obvious.

Suppose that $\left( e,e\right) \neq 0$ and that $\left( \cdot ,\cdot
\right) $ is positive semi-definite on $X.$ Then for each $y\in X$ we have $%
y=ae+z$ with $a\in \mathbb{K}$ and $z\in X^{\left( e\right) },$ from where
we get%
\begin{equation}
\left\vert \left( e,y\right) \right\vert ^{2}-\left( e,e\right) \left(
y,y\right) =-\left( e,e\right) \left( z,z\right) .  \label{1.1.4}
\end{equation}%
From (\ref{1.1.4}) we get the inequality (\ref{1.1.3}), with $x=e,$ in the
case that $\left( e,e\right) >0$ and (\ref{1.1.2}) in the case that $\left(
e,e\right) <0.$ In addition to this, from (\ref{1.1.4}) we observe that the
case of equality holds in (\ref{1.1.2}) or in (\ref{1.1.3}) if and only if $%
\left( z,z\right) =0,$ i.e., if and only if $y-ae\in X_{0}^{\left( e\right)
}.$

Conversely, if (\ref{1.1.3}) holds for all $x,y\in X,$ then $\left(
x,x\right) $ has the same sign over the whole of $X,$ i.e., $\left( \cdot
,\cdot \right) $ is semi-definite on $X.$ In the same manner, from (\ref%
{1.1.2}), for $y\in X^{\left( e\right) },$ we get $\left( e,e\right) \cdot
\left( y,y\right) \leq 0,$ which implies $\left( y,y\right) \geq 0,$ i.e., $%
\left( \cdot ,\cdot \right) $ is positive semi-definite on $X^{\left(
e\right) }.$

Now, suppose that $\left( \cdot ,\cdot \right) $ is positive semi-definite
on $X^{\left( e\right) }$ for at least one $e\in K$. Let us prove that $%
\left( \cdot ,\cdot \right) $ is positive semi-definite on $X^{\left(
f\right) }$ for each $f\in K.$

For a given $f\in K,$ consider the vector%
\begin{equation}
e^{\prime }:=e-\frac{\left( e,f\right) }{\left( f,f\right) }f.  \label{1.1.5}
\end{equation}%
Now,%
\begin{equation*}
\left( e^{\prime },e^{\prime }\right) =\left( e^{\prime },e\right) =\frac{%
\left( e,e\right) \left( f,f\right) -\left\vert \left( e,f\right)
\right\vert ^{2}}{\left( f,f\right) },\quad \left( e^{\prime },f\right) =0
\end{equation*}%
and together with%
\begin{equation*}
\left\vert \left( e,y\right) \right\vert ^{2}\geq \left( e,e\right) \left(
y,y\right) \quad \text{for any}\quad y\in X
\end{equation*}%
imply $\left( e^{\prime },e^{\prime }\right) \geq 0.$

There are two cases to be considered: $\left( e^{\prime },e^{\prime }\right)
>0$ and $\left( e^{\prime },e^{\prime }\right) =0.$

If $\left( e^{\prime },e^{\prime }\right) >0,$ then for any $x\in X^{\left(
f\right) },$ the vector%
\begin{equation*}
x^{\prime }:=x-ae^{\prime }\quad \text{with}\quad a=\frac{\left( x,e^{\prime
}\right) }{\left( e^{\prime },e^{\prime }\right) }
\end{equation*}%
satisfies the conditions%
\begin{equation*}
\left( x^{\prime },e\right) =0\quad \text{and}\quad \left( x^{\prime
},f\right) =0
\end{equation*}%
which implies%
\begin{equation*}
x^{\prime }\in X^{\left( e\right) }\quad \text{and}\quad \left( x,x\right)
=\left\vert a\right\vert ^{2}\left( e^{\prime },e^{\prime }\right) +\left(
x^{\prime },x^{\prime }\right) \geq 0.
\end{equation*}%
Therefore $\left( \cdot ,\cdot \right) $ is a positive semi-definite
functional on $X^{\left( f\right) }.$

From the parallelogram identity:%
\begin{equation}
\left( x+y,x+y\right) +\left( x-y,x-y\right) =2\left[ \left( x,x\right)
+\left( y,y\right) \right] ,\quad x,y\in X  \label{1.1.6}
\end{equation}%
we conclude that the set $X_{0}^{\left( e\right) }=X_{0}\cap X^{\left(
e\right) }$ is a linear subspace of $X.$

Since%
\begin{equation}
\left( x,y\right) =\frac{1}{4}\left[ \left( x+y,x+y\right) +\left(
x-y,x-y\right) \right] ,\quad x,y\in X  \label{1.1.7}
\end{equation}%
in the case of real spaces, and%
\begin{multline}
\left( x,y\right) =\frac{1}{4}\left[ \left( x+y,x+y\right) +\left(
x-y,x-y\right) \right]  \label{1.1.8} \\
+\frac{i}{4}\left[ \left( x+iy,x+iy\right) -\left( x-iy,x-iy\right) \right]
,\quad x,y\in X
\end{multline}%
in the case of complex spaces, hence $\left( x,y\right) =0$ provided that $x$
and $y$ belong to $X_{0}^{\left( e\right) }.$

If $\left( e^{\prime },e^{\prime }\right) =0,$ then $\left( e^{\prime
},e\right) =\left( e^{\prime },e^{\prime }\right) =0$ and then we can
conclude that $e^{\prime }\in X_{0}^{\left( e\right) }.$ Also, since $\left(
e^{\prime },e^{\prime }\right) =0$ implies $\left( e,f\right) \neq 0,$ hence
we have%
\begin{equation*}
f=b\left( e-e^{\prime }\right) \quad \text{with}\quad b=\frac{\left(
f,f\right) }{\left( e,f\right) }.
\end{equation*}%
Now write%
\begin{equation*}
X^{\left( e\right) }=X_{0}^{\left( e\right) }\bigoplus X_{+}^{\left(
e\right) },
\end{equation*}%
where $X_{+}^{\left( e\right) }$ is any direct complement of $X_{0}^{\left(
e\right) }$ in the space $X^{\left( e\right) }.$ If $y\neq 0,$ then $y\in
X_{+}^{\left( e\right) }$ implies $\left( y,y\right) >0.$ For such a vector $%
y,$ the vector%
\begin{equation*}
y^{\prime }:=e^{\prime }-\frac{\left( e^{\prime },y\right) }{\left(
y,y\right) }\cdot y.
\end{equation*}%
is in $X^{\left( e\right) }$ and therefore $\left( y^{\prime },y^{\prime
}\right) \geq 0.$

On the other hand%
\begin{equation*}
\left( y^{\prime },y^{\prime }\right) =\left( e^{\prime },y^{\prime }\right)
=-\frac{\left\vert \left( e^{\prime },y\right) \right\vert ^{2}}{\left(
y,y\right) }.
\end{equation*}%
Hence $y\in X_{+}^{\left( e\right) }$ implies that $\left( e^{\prime
},y\right) =0,$ i.e.,%
\begin{equation*}
\left( e,y\right) =\frac{\left( e,f\right) }{\left( f,f\right) }\left(
f,y\right) ,
\end{equation*}%
which together with $y\in X^{\left( e\right) }$ leads to $\left( f,y\right)
=0.$ Thus $y\in X_{+}^{\left( e\right) }$ implies $y\in X^{\left( f\right)
}. $

On the other hand $x\in X_{0}^{\left( e\right) }$ and $f=b\left( e-e^{\prime
}\right) $ imply $\left( f,x\right) =-b\left( e^{\prime },x\right) =0$ due
to the fact that $e^{\prime },x\in X_{0}^{\left( e\right) }.$

Hence $x\in X_{0}^{\left( e\right) }$ implies $\left( x,f\right) =0,$ i.e., $%
x\in X^{\left( f\right) }.$

From $X_{0}^{\left( e\right) }\subseteq X^{\left( f\right) }$ and $%
X_{+}^{\left( e\right) }\subseteq X^{\left( f\right) }$ we get $X^{\left(
e\right) }\subseteq X^{\left( f\right) }.$ Since $e\notin X^{\left( f\right)
}$ and $X=L\left( e\right) \bigoplus X^{\left( e\right) },$ we deduce $%
X^{\left( e\right) }=X^{\left( f\right) }$ and then $\left( \cdot ,\cdot
\right) $ is positive semi-definite on $X^{\left( f\right) }.$

The theorem is completely proved.
\end{proof}

In the case of complex linear spaces we may state the following result as
well \cite{IHFKU}:

\begin{theorem}[Kurepa, 1968]
\label{t1.1.2}Let $X$ be a complex linear space and $\left( \cdot ,\cdot
\right) $ a hermitian functional on $X.$

\begin{enumerate}
\item The functional $\left( \cdot ,\cdot \right) $ is semi-definite on $X$
if and only if there exists at least one vector $e\in X$ with $\left(
e,e\right) \neq 0$ such that%
\begin{equation}
\left[ \func{Re}\left( e,y\right) \right] ^{2}\leq \left( e,e\right) \left(
y,y\right) ,  \label{1.1.9}
\end{equation}%
for all $y\in X;$

\item There is no nonzero Hermitian functional $\left( \cdot ,\cdot \right) $
such that the inequality%
\begin{equation}
\left[ \func{Re}\left( e,y\right) \right] ^{2}\geq \left( e,e\right) \left(
y,y\right) ,\quad \left( e,e\right) \neq 0,  \label{1.1.10}
\end{equation}%
holds for all $y\in X$ and for an $e\in X.$
\end{enumerate}
\end{theorem}

\begin{proof}
We follow the proof in \cite{IHFKU}.

Let $\sigma $ and $\tau $ be real numbers and $x\in X^{\left( e\right) }$ a
given vector. For $y:=\left( \sigma +i\tau \right) e+x$ we get%
\begin{equation}
\left[ \func{Re}\left( e,y\right) \right] ^{2}-\left( e,e\right) \left(
y,y\right) =-\tau ^{2}\left( e,e\right) ^{2}-\left( e,e\right) \left(
x,x\right) .  \label{1.1.11}
\end{equation}%
If $\left( \cdot ,\cdot \right) $ is semi-definite on $X,$ then (\ref{1.1.11}%
) implies (\ref{1.1.9}).

Conversely, if (\ref{1.1.9}) holds for all $y\in X$ and for at least one $%
e\in X,$ then $\left( \cdot ,\cdot \right) $ is semi-definite on $X^{\left(
e\right) }.$ But (\ref{1.1.9}) and (\ref{1.1.11}) for $\tau =0$ lead to $%
-\left( e,e\right) \left( x,x\right) \leq 0$ from which it follows that $%
\left( e,e\right) $ and $\left( x,x\right) $ are of the same sign so that $%
\left( \cdot ,\cdot \right) $ is semi-definite on $X.$

Suppose that $\left( \cdot ,\cdot \right) \neq 0$ and that (\ref{1.1.10})
holds. We can assume that $\left( e,e\right) <0.$ Then (\ref{1.1.10})
implies that $\left( \cdot ,\cdot \right) $ is positive semi-definite on $%
X^{\left( e\right) }.$ On the other hand, if $\tau $ is such that%
\begin{equation*}
\tau ^{2}>-\frac{\left( x,x\right) }{\left( e,e\right) },
\end{equation*}%
then (\ref{1.1.11}) leads to $\left[ \func{Re}\left( e,y\right) \right]
^{2}<\left( e,e\right) \left( y,y\right) $, contradicting (\ref{1.1.10}).

Hence, if a Hermitian functional $\left( \cdot ,\cdot \right) $ is not
semi-definite and if $-\left( e,e\right) \neq 0,$ then the function $%
y\longmapsto \left[ \func{Re}\left( e,y\right) \right] ^{2}-\left(
e,e\right) \left( y,y\right) $ takes both positive and negative values.

The theorem is completely proved.
\end{proof}

\subsection{\label{ss1.2}Schwarz's Inequality for the Complexification of a
Real Space}

Let $X$ be a real linear space. The \textit{complexification }$X_{\mathbb{C}%
} $ of $X$ is defined as a complex linear space $X\times X$ of all ordered
pairs $\left\{ x,y\right\} $ $\left( x,y\in X\right) $ endowed with the
operations:%
\begin{align*}
\left\{ x,y\right\} +\left\{ x^{\prime },y^{\prime }\right\} & :=\left\{
x+x^{\prime },y+y^{\prime }\right\} , \\
\left( \sigma +i\tau \right) \cdot \left\{ x,y\right\} & :=\left\{ \sigma
x-\tau y,\sigma x+\tau y\right\} ,
\end{align*}%
where $x,y,x^{\prime },y^{\prime }\in X$ and $\sigma ,\tau \in \mathbb{R}$
(see for instance \cite{IHFKU2}).

If $z=\left\{ x,y\right\} ,$ then we can define the conjugate vector $\bar{z}
$ of $z$ by $\bar{z}:=\left\{ x,-y\right\} .$ Similarly, with the scalar
case, we denote%
\begin{equation*}
\func{Re}z=\left\{ x,0\right\} \quad \text{and\quad }\func{Im}z:=\left\{
0,y\right\} .
\end{equation*}%
Formally, we can write $z=x+iy=\func{Re}z+i\func{Im}z$ and $\bar{z}=x-iy=%
\func{Re}z-i\func{Im}z.$

Now, let $\left( \cdot ,\cdot \right) $ be a Hermitian functional on $X.$ We
may define on the complexification $X_{\mathbb{C}}$ of $X,$ the \textit{%
complexification }of $\left( \cdot ,\cdot \right) ,$ denoted by $\left(
\cdot ,\cdot \right) _{\mathbb{C}}$ and defined by:%
\begin{equation*}
\left( x+iy,x^{\prime }+iy^{\prime }\right) _{\mathbb{C}}:=\left(
x,x^{\prime }\right) +\left( y,y^{\prime }\right) +i\left[ \left(
y,x^{\prime }\right) -\left( x,y^{\prime }\right) \right] ,
\end{equation*}%
for $x,y,x^{\prime },y^{\prime }\in X.$

The following result may be stated \cite{IHFKU}:

\begin{theorem}[Kurepa, 1968]
\label{t1.1.3}Let $X,$ $X_{\mathbb{C}},$ $\left( \cdot ,\cdot \right) $ and $%
\left( \cdot ,\cdot \right) _{\mathbb{C}}$ be as above. An inequality of
type (\ref{1.1.2}) and (\ref{1.1.3}) holds for the functional $\left( \cdot
,\cdot \right) _{\mathbb{C}}$ in the space $X_{\mathbb{C}}$ if and only if
the same type of inequality holds for the functional $\left( \cdot ,\cdot
\right) $ in the space $X.$
\end{theorem}

\begin{proof}
We follow the proof in \cite{IHFKU}.

Firstly, observe that $\left( \cdot ,\cdot \right) $ is semi-definite if and
only if $\left( \cdot ,\cdot \right) _{\mathbb{C}}$ is semi-definite.

Now, suppose that $e\in X$ is such that%
\begin{equation*}
\left\vert \left( e,y\right) \right\vert ^{2}\geq \left( e,e\right) \left(
y,y\right) ,\quad \left( e,e\right) <0
\end{equation*}%
for all $y\in X.$ Then for $x,y\in X$ we have%
\begin{align*}
\left\vert \left( e,x+iy\right) _{\mathbb{C}}\right\vert ^{2}& =\left[
\left( e,x\right) \right] ^{2}+\left[ \left( e,y\right) \right] ^{2} \\
& \geq \left( e,e\right) \left[ \left( x,x\right) +\left( y,y\right) \right]
\\
& =\left( e,e\right) \left( x+iy,x+iy\right) _{\mathbb{C}}.
\end{align*}%
Hence, if for the functional $\left( \cdot ,\cdot \right) $ on $X$ an
inequality of type (\ref{1.1.2}) holds, then the same type of inequality
holds in $X_{\mathbb{C}}$ for the corresponding functional $\left( \cdot
,\cdot \right) _{\mathbb{C}}.$

Conversely, suppose that $e,f\in X$ are such that%
\begin{equation}
\left\vert \left( e+if,x+iy\right) _{\mathbb{C}}\right\vert ^{2}\geq \left(
e+if,e+if\right) _{\mathbb{C}}\left( x+iy,x+iy\right) _{\mathbb{C}}
\label{1.1.12}
\end{equation}%
holds for all $x,y\in X$ and that%
\begin{equation}
\left( e+if,e+if\right) _{\mathbb{C}}=\left( e,e\right) +\left( f,f\right)
<0.  \label{1.1.13}
\end{equation}

If $e=af$ with a real number $a,$ then (\ref{1.1.13}) implies that $\left(
f,f\right) <0$ and (\ref{1.1.12}) for $y=0$ leads to%
\begin{equation*}
\left[ \left( f,x\right) \right] ^{2}\geq \left( f,f\right) \left(
x,x\right) ,
\end{equation*}%
for all $x\in X.$ Hence, in this case, we have an inequality of type (\ref%
{1.1.2}) for the functional $\left( \cdot ,\cdot \right) $ in $X.$

Suppose that $e$ and $g$ are linearly independent and by $Y=L\left(
e,f\right) $ let us denote the subspace of $X$ consisting of all linear
combinations of $e$ and $f.$ On $Y$ we define a hermitian functional $D$ by
setting $D\left( x,y\right) =\left( x,y\right) $ for $x,y\in Y.$ Let $D_{%
\mathbb{C}}$ be the complexification of $D.$ Then (\ref{1.1.12}) implies:%
\begin{multline}
 \left\vert D_{\mathbb{C}}\left( e+if,x+iy\right) \right\vert
^{2}
\label{1.1.14} \\
\geq D_{\mathbb{C}}\left( e+if,e+if\right) D_{\mathbb{C}}\left(
x+iy,x+iy\right) ,\quad x,y\in X
\end{multline}%
and (\ref{1.1.13}) implies%
\begin{equation}
D\left( e,e\right) +D\left( f,f\right) <0.  \label{1.1.15}
\end{equation}%
Further, consider in $Y$ a base consisting of the two vectors $\left\{
u_{1},u_{2}\right\} $ on which $D$ is diagonal, i.e., $D$ satisfies%
\begin{equation*}
D\left( x,y\right) =\lambda _{1}x_{1}y_{1}+\lambda _{2}x_{2}y_{2},
\end{equation*}%
where%
\begin{equation*}
x=x_{1}u_{1}+x_{2}u_{2},\quad y=y_{1}u_{1}+y_{2}u_{2},
\end{equation*}%
and%
\begin{equation*}
\lambda _{1}=D\left( u_{1},u_{1}\right) ,\quad \lambda _{2}=D\left(
u_{2},u_{2}\right) .
\end{equation*}%
Since for the functional $D$ we have the relations (\ref{1.1.15}) and (\ref%
{1.1.14}), we conclude that $D$ is not a semi-definite functional on $Y.$
Hence $\lambda _{1}\cdot \lambda _{2}<0,$ so we can take $\lambda _{1}<0$
and $\lambda _{2}>0.$

Set%
\begin{equation*}
X^{+}:=\left\{ x|\left( x,e\right) =\left( x,f\right) =0,\ x\in X\right\} .
\end{equation*}%
Obviously, $\left( x,e\right) =\left( x,f\right) =0$ if and only if $\left(
x_{1}u_{1}\right) =\left( x_{2}u_{2}\right) =0.$

Now, if $y\in X,$ then the vector%
\begin{equation}
x:=y-\frac{\left( y,u_{1}\right) }{\left( u_{1},u_{1}\right) }u_{1}-\frac{%
\left( y,u_{2}\right) }{\left( u_{2},u_{2}\right) }u_{2}  \label{1.1.16}
\end{equation}%
belongs to $X^{+}.$ From this it follows that%
\begin{equation*}
X=L\left( e,f\right) \bigoplus X^{+}.
\end{equation*}%
Now, replacing in (\ref{1.1.12}) the vector $x+iy$ with $z\in X^{+},$ we get
from (\ref{1.1.13}) that%
\begin{equation*}
\left[ \left( e,e\right) +\left( f,f\right) \right] \left( z,z\right) \leq 0,
\end{equation*}%
which, together with (\ref{1.1.13}) leads to $\left( z,z\right) \geq 0.$
Therefore the functional $\left( \cdot ,\cdot \right) $ is positive
semi-definite on $X^{+}.$

Now, since any $y\in X$ is of the form (\ref{1.1.16}), hence for $y\in
X^{\left( u_{1}\right) }$ we get%
\begin{equation*}
\left( y,y\right) =\left( x,x\right) +\frac{\left[ \left( y,u_{2}\right) %
\right] ^{2}}{\lambda _{2}},
\end{equation*}%
which is a nonnegative number. Thus, $\left( \cdot ,\cdot \right) $ is
positive semi-definite on the space $X^{\left( u_{1}\right) }.$ Since $%
\left( u_{1},u_{1}\right) <0$ we have $\left[ \left( u_{1},y\right) \right]
^{2}\geq \left( u_{1},u_{1}\right) \left( y,y\right) $ for any $y\in X$ and
the theorem is completely proved.
\end{proof}

\section{Superadditivity and Monotonicity\label{s2}}

\subsection{\label{ss2.1}The Convex Cone of Nonnegative Hermitian Forms}

Let $X$ be a linear space over the real or complex number field $\mathbb{K}$
and let us denote by $\mathcal{H}\left( X\right) $ the class of all positive
semi-definite Hermitian forms on $X,$ or, for simplicity, \textit{%
nonnegative }forms on $X,$ i.e., the mapping $\left( \cdot ,\cdot \right)
:X\times X\rightarrow \mathbb{K}$ belongs to $\mathcal{H}\left( X\right) $
if it satisfies the conditions

\begin{enumerate}
\item[(i)] $\left( x,x\right) \geq 0$ for all $x$ in $X;$

\item[(ii)] $\left( \alpha x+\beta y,z\right) =\alpha \left( x,z\right)
+\beta \left( y,z\right) $ for all $x,y\in X$ and $\alpha ,\beta \in \mathbb{%
K}$

\item[(iii)] $\left( y,x\right) =\overline{\left( x,y\right) }$ for all $%
x,y\in X.$
\end{enumerate}

If $\left( \cdot ,\cdot \right) \in \mathcal{H}\left( X\right) ,$ then the
functional $\left\Vert \cdot \right\Vert =\left( \cdot ,\cdot \right) ^{%
\frac{1}{2}}$ is a \textit{semi-norm} on $X$ and the following equivalent
versions of Schwarz's inequality hold:%
\begin{equation}
\left\Vert x\right\Vert ^{2}\left\Vert y\right\Vert ^{2}\geq \left\vert
\left( x,y\right) \right\vert ^{2}\quad \text{or\quad }\left\Vert
x\right\Vert \left\Vert y\right\Vert \geq \left\vert \left( x,y\right)
\right\vert  \label{1.2.1}
\end{equation}%
for any $x,y\in X.$

Now, let us observe that $\mathcal{H}\left( X\right) $ is a \textit{convex
cone }in the linear space of all mappings defined on $X^{2}$ with values in $%
\mathbb{K}$, i.e.,

\begin{enumerate}
\item[(e)] $\left( \cdot ,\cdot \right) _{1},\left( \cdot ,\cdot \right)
_{2}\in \mathcal{H}\left( X\right) $ implies that $\left( \cdot ,\cdot
\right) _{1}+\left( \cdot ,\cdot \right) _{2}\in \mathcal{H}\left( X\right)
; $

\item[(ee)] $\alpha \geq 0$ and $\left( \cdot ,\cdot \right) \in \mathcal{H}%
\left( X\right) $ implies that $\alpha \left( \cdot ,\cdot \right) \in
\mathcal{H}\left( X\right) .$
\end{enumerate}

We can introduce on $\mathcal{H}\left( X\right) $ the following binary
relation \cite{IHFDM}:%
\begin{equation}
\left( \cdot ,\cdot \right) _{2}\geq \left( \cdot ,\cdot \right) _{1}\quad
\text{if and only if\quad }\left\Vert x\right\Vert _{2}\geq \left\Vert
x\right\Vert _{1}\quad \text{for all \ }x\in X.  \label{1.2.2}
\end{equation}%
We observe that the following properties hold:

\begin{enumerate}
\item[(b)] $\left( \cdot ,\cdot \right) _{2}\geq \left( \cdot ,\cdot \right)
_{1}$ for all $\left( \cdot ,\cdot \right) \in \mathcal{H}\left( X\right) ;$

\item[(bb)] $\left( \cdot ,\cdot \right) _{3}\geq \left( \cdot ,\cdot
\right) _{2}$ and $\left( \cdot ,\cdot \right) _{2}\geq \left( \cdot ,\cdot
\right) _{1}$ implies that $\left( \cdot ,\cdot \right) _{3}\geq \left(
\cdot ,\cdot \right) _{1};$

\item[(bbb)] $\left( \cdot ,\cdot \right) _{2}\geq \left( \cdot ,\cdot
\right) _{1}$ and $\left( \cdot ,\cdot \right) _{1}\geq \left( \cdot ,\cdot
\right) _{2}$ implies that $\left( \cdot ,\cdot \right) _{2}=\left( \cdot
,\cdot \right) _{1};$
\end{enumerate}

i.e., the binary relation defined by (\ref{1.2.2}) is \textit{an order
relation }on $\mathcal{H}\left( X\right) .$

While (b) and (bb) are obvious from the definition, we should remark, for
(bbb), that if $\left( \cdot ,\cdot \right) _{2}\geq \left( \cdot ,\cdot
\right) _{1}$ and $\left( \cdot ,\cdot \right) _{1}\geq \left( \cdot ,\cdot
\right) _{2},$ then obviously $\left\Vert x\right\Vert _{2}=\left\Vert
x\right\Vert _{1}$ for all $x\in X,$ which implies, by the following well
known identity:%
\begin{multline}
\quad\left( x,y\right) _{k}\\:=\frac{1}{4}\left[ \left\Vert
x+y\right\Vert _{k}^{2}-\left\Vert x-y\right\Vert _{k}^{2}+i\left(
\left\Vert x+iy\right\Vert _{k}^{2}-\left\Vert x-iy\right\Vert
_{k}^{2}\right) \right] \label{1.2.3}\quad
\end{multline}%
with $x,y\in X$ and $k\in \left\{ 1,2\right\} $, that $\left(
x,y\right) _{2}=\left( x,y\right) _{1}$ for all $x,y\in X.$

\subsection{\label{ss2.2}The Superadditivity and Monotonicity of $\protect%
\sigma -$Mapping}

Let us consider the following mapping \cite{IHFDM}:%
\begin{equation*}
\sigma :\mathcal{H}\left( X\right) \times X^{2}\rightarrow \mathbb{R}%
_{+},\quad \sigma \left( \left( \cdot ,\cdot \right) ;x,y\right)
:=\left\Vert x\right\Vert \left\Vert y\right\Vert -\left\vert \left(
x,y\right) \right\vert ,
\end{equation*}%
which is closely related to Schwarz's inequality (\ref{1.2.1}).

The following simple properties of $\sigma $ are obvious:

\begin{enumerate}
\item[(s)] $\sigma \left( \alpha \left( \cdot ,\cdot \right) ;x,y\right)
=\alpha \sigma \left( \left( \cdot ,\cdot \right) ;x,y\right) ;$

\item[(ss)] $\sigma \left( \left( \cdot ,\cdot \right) ;y,x\right) =\sigma
\left( \left( \cdot ,\cdot \right) ;x,y\right) ;$

\item[(sss)] $\sigma \left( \left( \cdot ,\cdot \right) ;x,y\right) \geq 0$
(Schwarz's inequality);
\end{enumerate}

for any $\alpha \geq 0,$ $\left( \cdot ,\cdot \right) \in \mathcal{H}\left(
X\right) $ and $x,y\in X.$

The following result concerning the functional properties of $\sigma $ as a
function depending on the nonnegative hermitian form $\left( \cdot ,\cdot
\right) $ has been obtained in \cite{IHFDM}:

\begin{theorem}[Dragomir-Mond, 1994]
\label{t1.2.1}The mapping $\sigma $ satisfies the following statements:

\begin{enumerate}
\item[(i)] For every $\left( \cdot ,\cdot \right) _{i}\in \mathcal{H}\left(
X\right) $ $\left( i=1,2\right) $ one has the inequality%
\begin{multline}
\qquad \sigma \left( \left( \cdot ,\cdot \right) _{1}+\left( \cdot ,\cdot
\right) _{2};x,y\right)  \label{1.2.4} \\
\geq \sigma \left( \left( \cdot ,\cdot \right) _{1};x,y\right) +\sigma
\left( \left( \cdot ,\cdot \right) _{2};x,y\right) \qquad \left( \geq
0\right) \qquad
\end{multline}%
for all $x,y\in X$, i.e., the mapping $\sigma \left( \cdot ;x,y\right) $ is
superadditive on $\mathcal{H}\left( X\right) ;$

\item[(ii)] For every $\left( \cdot ,\cdot \right) _{i}\in \mathcal{H}\left(
X\right) $ $\left( i=1,2\right) $ with $\left( \cdot ,\cdot \right) _{2}\geq
\left( \cdot ,\cdot \right) _{1}$ one has%
\begin{equation}
\sigma \left( \left( \cdot ,\cdot \right) _{2};x,y\right) \geq \sigma \left(
\left( \cdot ,\cdot \right) _{1};x,y\right) \qquad \left( \geq 0\right)
\label{1.2.5}
\end{equation}%
for all $x,y\in X,$ i.e., the mapping $\sigma \left( \cdot ;x,y\right) $ is
nondecreasing on $\mathcal{H}\left( X\right) .$
\end{enumerate}
\end{theorem}

\begin{proof}
We follow the proof in \cite{IHFDM}.

\begin{enumerate}
\item[(i)] By the Cauchy-Bunyakovsky-Schwarz inequality for real numbers, we
have%
\begin{equation*}
\left( a^{2}+b^{2}\right) ^{\frac{1}{2}}\left( c^{2}+d^{2}\right) ^{\frac{1}{%
2}}\geq ac+bd;\quad a,b,c,d\geq 0.
\end{equation*}%
Therefore,%
\begin{align*}
& \sigma \left( \left( \cdot ,\cdot \right) _{1}+\left( \cdot ,\cdot \right)
_{2};x,y\right) \\
& =\left( \left\Vert x\right\Vert _{1}^{2}+\left\Vert x\right\Vert
_{2}^{2}\right) ^{\frac{1}{2}}\left( \left\Vert y\right\Vert
_{1}^{2}+\left\Vert y\right\Vert _{2}^{2}\right) ^{\frac{1}{2}}-\left\vert
\left( x,y\right) _{1}+\left( x,y\right) _{2}\right\vert \\
& \geq \left\Vert x\right\Vert _{1}\left\Vert y\right\Vert _{1}+\left\Vert
x\right\Vert _{2}\left\Vert y\right\Vert _{2}-\left\vert \left( x,y\right)
_{1}\right\vert -\left\vert \left( x,y\right) _{2}\right\vert \\
& =\sigma \left( \left( \cdot ,\cdot \right) _{1};x,y\right) +\sigma \left(
\left( \cdot ,\cdot \right) _{2};x,y\right) ,
\end{align*}%
for all $\left( \cdot ,\cdot \right) _{i}\in \mathcal{H}\left( X\right) $ $%
\left( i=1,2\right) $ and $x,y\in X,$ and the statement is proved.

\item[(ii)] Suppose that $\left( \cdot ,\cdot \right) _{2}\geq \left( \cdot
,\cdot \right) _{1}$ and define $\left( \cdot ,\cdot \right) _{2,1}:=\left(
\cdot ,\cdot \right) _{2}-\left( \cdot ,\cdot \right) _{1}.$ It is obvious
that $\left( \cdot ,\cdot \right) _{2,1}$ is a nonnegative hermitian form
and thus, by the above property one has,%
\begin{align*}
\sigma \left( \left( \cdot ,\cdot \right) _{2};x,y\right) & \geq \sigma
\left( \left( \cdot ,\cdot \right) _{2,1}+\left( \cdot ,\cdot \right)
_{1};x,y\right) \\
& \geq \sigma \left( \left( \cdot ,\cdot \right) _{2,1};x,y\right) +\sigma
\left( \left( \cdot ,\cdot \right) _{1};x,y\right)
\end{align*}%
from where we get:%
\begin{equation*}
\sigma \left( \left( \cdot ,\cdot \right) _{2};x,y\right) -\sigma \left(
\left( \cdot ,\cdot \right) _{1};x,y\right) \geq \sigma \left( \left( \cdot
,\cdot \right) _{2,1};x,y\right) \geq 0
\end{equation*}%
and the proof of the theorem is completed.
\end{enumerate}
\end{proof}

\begin{remark}
\label{r1.2.1}If we consider the related mapping \cite{IHFDM}%
\begin{equation*}
\sigma _{r}\left( \left( \cdot ,\cdot \right) ;x,y\right) :=\left\Vert
x\right\Vert \left\Vert y\right\Vert -\func{Re}\left( x,y\right) ,
\end{equation*}%
then we can show, as above, that $\sigma \left( \cdot ;x,y\right) $ is%
\textbf{\ superadditive and nondecreasing} on $\mathcal{H}\left( X\right) .$

Moreover, if we introduce another mapping, namely, \cite{IHFDM}%
\begin{equation*}
\tau :\mathcal{H}\left( X\right) \times X^{2}\rightarrow \mathbb{R}%
_{+},\quad \tau \left( \left( \cdot ,\cdot \right) ;x,y\right) :=\left(
\left\Vert x\right\Vert +\left\Vert y\right\Vert \right) ^{2}-\left\Vert
x+y\right\Vert ^{2},
\end{equation*}%
which is connected with the \textbf{triangle inequality}%
\begin{equation}
\left\Vert x+y\right\Vert \leq \left\Vert x\right\Vert +\left\Vert
y\right\Vert \quad \text{for any\quad }x,y\in X  \label{1.2.6}
\end{equation}%
then we observe that%
\begin{equation}
\tau \left( \left( \cdot ,\cdot \right) ;x,y\right) =2\sigma _{r}\left(
\left( \cdot ,\cdot \right) ;x,y\right)  \label{1.2.7}
\end{equation}%
for all $\left( \cdot ,\cdot \right) \in \mathcal{H}\left( X\right) $ and $%
x,y\in X,$ therefore $\sigma \left( \cdot ;x,y\right) $ is in its turn a
\textbf{superadditive} and \textbf{nondecreasing} functional on $\mathcal{H}%
\left( X\right) .$
\end{remark}

\subsection{\label{ss2.3}The Superadditivity and Monotonicity of $\protect%
\delta -$Mapping}

Now consider another mapping naturally associated to Schwarz's inequality,
namely \cite{IHFDM}%
\begin{equation*}
\delta :\mathcal{H}\left( X\right) \times X^{2}\rightarrow \mathbb{R}%
_{+},\quad \delta \left( \left( \cdot ,\cdot \right) ;x,y\right)
:=\left\Vert x\right\Vert ^{2}\left\Vert y\right\Vert ^{2}-\left\vert \left(
x,y\right) \right\vert ^{2}.
\end{equation*}%
It is obvious that the following properties are valid:

\begin{enumerate}
\item[(i)] $\delta \left( \left( \cdot ,\cdot \right) ;x,y\right) \geq 0$
(Schwarz's inequality);

\item[(ii)] $\delta \left( \left( \cdot ,\cdot \right) ;x,y\right) =\delta
\left( \left( \cdot ,\cdot \right) ;y,x\right) ;$

\item[(iii)] $\delta \left( \alpha \left( \cdot ,\cdot \right) ;x,y\right)
=\alpha ^{2}\delta \left( \left( \cdot ,\cdot \right) ;x,y\right) $
\end{enumerate}

for all $x,y\in X,$ $\alpha \geq 0$ and $\left( \cdot ,\cdot \right) \in
\mathcal{H}\left( X\right) .$

The following theorem incorporates some further properties of this
functional \cite{IHFDM}:

\begin{theorem}[Dragomir-Mond, 1994]
\label{t1.2.2}With the above assumptions, we have:

\begin{enumerate}
\item[(i)] If $\left( \cdot ,\cdot \right) _{i}\in \mathcal{H}\left(
X\right) $ $\left( i=1,2\right) ,$ then%
\begin{multline}
\delta \left( \left( \cdot ,\cdot \right) _{1}+\left( \cdot ,\cdot \right)
_{2};x,y\right) -\delta \left( \left( \cdot ,\cdot \right) _{1};x,y\right)
-\delta \left( \left( \cdot ,\cdot \right) _{2};x,y\right)  \label{1.2.8} \\
\geq \left( \det \left[
\begin{array}{ll}
\left\Vert x\right\Vert _{1} & \left\Vert y\right\Vert _{1} \\
\left\Vert x\right\Vert _{2} & \left\Vert y\right\Vert _{2}%
\end{array}%
\right] \right) ^{2}\qquad \left( \geq 0\right) ;
\end{multline}%
i.e., the mapping $\delta \left( \cdot ;x,y\right) $ is strong superadditive
on $\mathcal{H}\left( X\right) .$

\item[(ii)] If $\left( \cdot ,\cdot \right) _{i}\in \mathcal{H}\left(
X\right) $ $\left( i=1,2\right) ,$ with $\left( \cdot ,\cdot \right)
_{2}\geq \left( \cdot ,\cdot \right) _{1},$ then%
\begin{multline}
\delta \left( \left( \cdot ,\cdot \right) _{2};x,y\right) -\delta \left(
\left( \cdot ,\cdot \right) _{1};x,y\right)  \label{1.2.9} \\
\geq \left( \det \left[
\begin{array}{cc}
\left\Vert x\right\Vert _{1} & \left\Vert y\right\Vert _{1} \\
\left( \left\Vert x\right\Vert _{2}^{2}-\left\Vert x\right\Vert
_{1}^{2}\right) ^{\frac{1}{2}} & \left( \left\Vert y\right\Vert
_{2}^{2}-\left\Vert y\right\Vert _{1}^{2}\right) ^{\frac{1}{2}}%
\end{array}%
\right] \right) ^{2}\qquad \left( \geq 0\right) ;
\end{multline}%
i.e., the mapping $\delta \left( \cdot ;x,y\right) $ is strong nondecreasing
on $\mathcal{H}\left( X\right) .$
\end{enumerate}
\end{theorem}

\begin{proof}
(i) For all $\left( \cdot ,\cdot \right) _{i}\in \mathcal{H}\left( X\right) $
$\left( i=1,2\right) $ and $x,y\in X$ we have%
\begin{align}
& \delta \left( \left( \cdot ,\cdot \right) _{1}+\left( \cdot ,\cdot \right)
_{2};x,y\right)  \label{1.2.10} \\
& =\left( \left\Vert x\right\Vert _{2}^{2}-\left\Vert x\right\Vert
_{1}^{2}\right) \left( \left\Vert y\right\Vert _{2}^{2}-\left\Vert
y\right\Vert _{1}^{2}\right) -\left\vert \left( x,y\right) _{2}+\left(
x,y\right) _{1}\right\vert ^{2}  \notag \\
& \geq \left\Vert x\right\Vert _{2}^{2}\left\Vert y\right\Vert
_{2}^{2}+\left\Vert x\right\Vert _{1}^{2}\left\Vert y\right\Vert
_{1}^{2}+\left\Vert x\right\Vert _{1}^{2}\left\Vert y\right\Vert
_{2}^{2}+\left\Vert x\right\Vert _{2}^{2}\left\Vert y\right\Vert _{1}^{2}
\notag \\
& \qquad \qquad -\left( \left\vert \left( x,y\right) _{2}\right\vert
+\left\vert \left( x,y\right) _{1}\right\vert \right) ^{2}  \notag \\
& =\delta \left( \left( \cdot ,\cdot \right) _{2};x,y\right) +\delta \left(
\left( \cdot ,\cdot \right) _{1};x,y\right)  \notag \\
& \qquad \qquad +\left\Vert x\right\Vert _{1}^{2}\left\Vert y\right\Vert
_{2}^{2}+\left\Vert x\right\Vert _{2}^{2}\left\Vert y\right\Vert
_{1}^{2}-2\left\vert \left( x,y\right) _{2}\left( x,y\right) _{1}\right\vert
.  \notag
\end{align}%
By Schwarz's inequality we have%
\begin{equation}
\left\vert \left( x,y\right) _{2}\left( x,y\right) _{1}\right\vert \leq
\left\Vert x\right\Vert _{1}\left\Vert y\right\Vert _{1}\left\Vert
x\right\Vert _{2}\left\Vert y\right\Vert _{2},  \label{1.2.11}
\end{equation}%
therefore, by (\ref{1.2.10}) and (\ref{1.2.11}), we can state that%
\begin{align*}
& \delta \left( \left( \cdot ,\cdot \right) _{1}+\left( \cdot ,\cdot \right)
_{2};x,y\right) -\delta \left( \left( \cdot ,\cdot \right) _{1};x,y\right)
-\delta \left( \left( \cdot ,\cdot \right) _{2};x,y\right) \\
& \geq \left\Vert x\right\Vert _{1}^{2}\left\Vert y\right\Vert
_{2}^{2}+\left\Vert x\right\Vert _{2}^{2}\left\Vert y\right\Vert
_{1}^{2}-2\left\Vert x\right\Vert _{1}\left\Vert y\right\Vert _{1}\left\Vert
x\right\Vert _{2}\left\Vert y\right\Vert _{2} \\
& =\left( \left\Vert x\right\Vert _{1}\left\Vert y\right\Vert
_{2}-\left\Vert x\right\Vert _{2}\left\Vert y\right\Vert _{1}\right) ^{2}
\end{align*}%
and the inequality (\ref{1.2.8}) is proved.

(ii) Suppose that $\left( \cdot ,\cdot \right) _{2}\geq \left( \cdot ,\cdot
\right) _{1}$ and, as in Theorem \ref{t1.2.1}, define $\left( \cdot ,\cdot
\right) _{2,1}:=\left( \cdot ,\cdot \right) _{2}-\left( \cdot ,\cdot \right)
_{1}.$ Then $\left( \cdot ,\cdot \right) _{2,1}$ is a nonnegative hermitian
form and by (i) we have%
\begin{align*}
& \delta \left( \left( \cdot ,\cdot \right) _{2,1};x,y\right) -\delta \left(
\left( \cdot ,\cdot \right) _{1};x,y\right) \\
& =\delta \left( \left( \cdot ,\cdot \right) _{2,1}+\left( \cdot ,\cdot
\right) _{1};x,y\right) -\delta \left( \left( \cdot ,\cdot \right)
_{1};x,y\right) \\
& \geq \delta \left( \left( \cdot ,\cdot \right) _{2,1};x,y\right) +\left(
\det \left[
\begin{array}{ll}
\left\Vert x\right\Vert _{1} & \left\Vert y\right\Vert _{1} \\
\left\Vert x\right\Vert _{2,1} & \left\Vert y\right\Vert _{2,1}%
\end{array}%
\right] \right) ^{2} \\
& \geq \left( \det \left[
\begin{array}{ll}
\left\Vert x\right\Vert _{1} & \left\Vert y\right\Vert _{1} \\
\left\Vert x\right\Vert _{2,1} & \left\Vert y\right\Vert _{2,1}%
\end{array}%
\right] \right) ^{2}.
\end{align*}%
Since $\left\Vert z\right\Vert _{2,1}=\left( \left\Vert z\right\Vert
_{2}^{2}-\left\Vert z\right\Vert _{1}^{2}\right) ^{\frac{1}{2}}$ for $z\in
X, $ hence the inequality (\ref{1.2.9}) is proved.
\end{proof}

\begin{remark}
\label{r1.2.2}If we consider the functional $\delta _{r}\left( \left( \cdot
,\cdot \right) ;x,y\right) :=\left\Vert x\right\Vert ^{2}\left\Vert
y\right\Vert ^{2}-\left[ \func{Re}\left( x,y\right) \right] ^{2},$ then we
can state similar properties for it. We omit the details.
\end{remark}

\subsection{\label{ss2.4}Superadditivity and Monotonicity of $\protect\beta %
- $Mapping}

Consider the functional $\beta :\mathcal{H}\left( X\right) \times
X^{2}\rightarrow \mathbb{R}$ defined by \cite{IHFDM2}
\begin{equation}
\beta \left( \left( \cdot ,\cdot \right) ;x,y\right) =\left( \left\Vert
x\right\Vert ^{2}\left\Vert y\right\Vert ^{2}-\left\vert \left( x,y\right)
\right\vert ^{2}\right) ^{\frac{1}{2}}.  \label{1.2.12}
\end{equation}%
It is obvious that $\beta \left( \left( \cdot ,\cdot \right) ;x,y\right) =%
\left[ \delta \left( \left( \cdot ,\cdot \right) ;x,y\right) \right] ^{\frac{%
1}{2}}$ and thus it is monotonic nondecreasing on $\mathcal{H}\left(
X\right) .$ Before we prove that $\beta \left( \cdot ;x,y\right) $ is also
superadditive, which apparently does not follow from the properties of $%
\delta $ pointed out in the subsection above, we need the following simple
lemma:

\begin{lemma}
\label{l1.2.1}If $\left( \cdot ,\cdot \right) $ is a nonnegative Hermitian
form on $X,$ $x,y\in X$ and $\left\Vert y\right\Vert \neq 0,$ then%
\begin{equation}
\inf_{\lambda \in \mathbb{K}}\left\Vert x-\lambda y\right\Vert ^{2}=\frac{%
\left\Vert x\right\Vert ^{2}\left\Vert y\right\Vert ^{2}-\left\vert \left(
x,y\right) \right\vert ^{2}}{\left\Vert y\right\Vert ^{2}}.  \label{1.2.13}
\end{equation}
\end{lemma}

\begin{proof}
Observe that%
\begin{equation*}
\left\Vert x-\lambda y\right\Vert ^{2}=\left\Vert x\right\Vert ^{2}-2\func{Re%
}\left[ \lambda \left( x,y\right) \right] +\left\vert \lambda \right\vert
^{2}\left\Vert y\right\Vert ^{2}
\end{equation*}%
and, for $\left\Vert y\right\Vert \neq 0,$
\begin{multline*}
\frac{\left\Vert x\right\Vert ^{2}\left\Vert y\right\Vert ^{2}-\left\vert
\left( x,y\right) \right\vert ^{2}+\left\vert \mu \left\Vert y\right\Vert
^{2}-\left( x,y\right) \right\vert ^{2}}{\left\Vert y\right\Vert ^{2}} \\
=\left\Vert x\right\Vert ^{2}-2\func{Re}\left[ \mu \overline{\left(
x,y\right) }\right] +\left\vert \mu \right\vert ^{2}\left\Vert y\right\Vert
^{2},
\end{multline*}%
and since $\func{Re}\left[ \bar{\lambda}\left( x,y\right) \right] =\func{Re}%
\overline{\left[ \bar{\lambda}\left( x,y\right) \right] }=\func{Re}\left[
\lambda \overline{\left( x,y\right) }\right] ,$ we deduce the equality%
\begin{equation}
\left\Vert x-\lambda y\right\Vert ^{2}=\frac{\left\Vert x\right\Vert
^{2}\left\Vert y\right\Vert ^{2}-\left\vert \left( x,y\right) \right\vert
^{2}+\left\vert \mu \left\Vert y\right\Vert ^{2}-\left( x,y\right)
\right\vert ^{2}}{\left\Vert y\right\Vert ^{2}},  \label{1.2.14}
\end{equation}%
for any $x,y\in X$ with $\left\Vert y\right\Vert \neq 0.$

Taking the infimum over $\lambda \in \mathbb{K}$ in (\ref{1.2.14}), we
deduce the desired result (\ref{1.2.13}).
\end{proof}

For the subclass $\mathcal{JP}\left( X\right) ,$ of all inner products
defined on $X,$ of $\mathcal{H}\left( X\right) $ \ and $y\neq 0,$ we may
define%
\begin{align*}
\gamma \left( \left( \cdot ,\cdot \right) ;x,y\right) & =\frac{\left\Vert
x\right\Vert ^{2}\left\Vert y\right\Vert ^{2}-\left\vert \left( x,y\right)
\right\vert ^{2}}{\left\Vert y\right\Vert ^{2}} \\
& =\frac{\delta \left( \left( \cdot ,\cdot \right) ;x,y\right) }{\left\Vert
y\right\Vert ^{2}}.
\end{align*}

The following result may be stated (see also \cite{IHFDM2}):

\begin{theorem}[Dragomir-Mond, 1996]
\label{t1.2.3}The functional $\gamma \left( \cdot ;x,y\right) $ is
superadditive and monotonic nondecreasing on $\mathcal{JP}\left( X\right) $
for any $x,y\in X$ with $y\neq 0.$
\end{theorem}

\begin{proof}
Let $\left( \cdot ,\cdot \right) _{1},\left( \cdot ,\cdot \right) _{2}\in
\mathcal{JP}\left( X\right) .$ Then%
\begin{align}
& \gamma \left( \left( \cdot ,\cdot \right) _{1}+\left( \cdot ,\cdot \right)
_{2};x,y\right)  \label{1.2.15} \\
& =\frac{\left( \left\Vert x\right\Vert _{1}^{2}+\left\Vert x\right\Vert
_{2}^{2}\right) \left( \left\Vert y\right\Vert _{1}^{2}+\left\Vert
y\right\Vert _{2}^{2}\right) -\left\vert \left( x,y\right) _{1}+\left(
x,y\right) _{2}\right\vert ^{2}}{\left\Vert y\right\Vert _{1}^{2}\left\Vert
y\right\Vert _{2}^{2}}  \notag \\
& =\inf_{\lambda \in \mathbb{K}}\left[ \left\Vert x-\lambda y\right\Vert
_{1}^{2}+\left\Vert x-\lambda y\right\Vert _{2}^{2}\right] ,  \notag
\end{align}%
and for the last equality we have used Lemma \ref{l1.2.1}.

Also,%
\begin{align}
\gamma \left( \left( \cdot ,\cdot \right) _{i};x,y\right) & =\frac{%
\left\Vert x\right\Vert _{i}^{2}\left\Vert y\right\Vert _{i}^{2}-\left\vert
\left( x,y\right) _{i}\right\vert ^{2}}{\left\Vert y\right\Vert _{i}^{2}}
\label{1.2.16} \\
& =\inf_{\lambda \in \mathbb{K}}\left\Vert x-\lambda y\right\Vert
_{i}^{2},\qquad i=1,2.  \notag
\end{align}%
Utilising the infimum property that
\begin{equation*}
\inf_{\lambda \in \mathbb{K}}\left( f\left( \lambda \right) +g\left( \lambda
\right) \right) \geq \inf_{\lambda \in \mathbb{K}}f\left( \lambda \right)
+\inf_{\lambda \in \mathbb{K}}g\left( \lambda \right) ,
\end{equation*}
we can write that%
\begin{equation*}
\inf_{\lambda \in \mathbb{K}}\left[ \left\Vert x-\lambda y\right\Vert
_{1}^{2}+\left\Vert x-\lambda y\right\Vert _{2}^{2}\right] \geq
\inf_{\lambda \in \mathbb{K}}\left\Vert x-\lambda y\right\Vert
_{1}^{2}+\inf_{\lambda \in \mathbb{K}}\left\Vert x-\lambda y\right\Vert
_{2}^{2},
\end{equation*}%
which proves the superadditivity of $\gamma \left( \cdot ;x,y\right) .$

The monotonicity follows by the superadditivity property and the theorem is
completely proved.
\end{proof}

\begin{corollary}
\label{c1.2.1}If $\left( \cdot ,\cdot \right) _{i}\in \mathcal{JP}\left(
X\right) $ with $\left( \cdot ,\cdot \right) _{2}\geq \left( \cdot ,\cdot
\right) _{1}$ and $x,y\in X$ are such that $x,y\neq 0,$ then:%
\begin{align}
\delta \left( \left( \cdot ,\cdot \right) _{2};x,y\right) & \geq \max
\left\{ \frac{\left\Vert y\right\Vert _{2}^{2}}{\left\Vert y\right\Vert
_{1}^{2}},\frac{\left\Vert x\right\Vert _{2}^{2}}{\left\Vert x\right\Vert
_{1}^{2}}\right\} \delta \left( \left( \cdot ,\cdot \right) _{1};x,y\right)
\label{1.2.17} \\
& \left( \geq \delta \left( \left( \cdot ,\cdot \right) _{1};x,y\right)
\right)  \notag
\end{align}%
or equivalently, \cite{IHFDM2}%
\begin{multline}
\delta \left( \left( \cdot ,\cdot \right) _{2};x,y\right) -\delta \left(
\left( \cdot ,\cdot \right) _{1};x,y\right)  \label{1.2.18} \\
\geq \max \left\{ \frac{\left\Vert y\right\Vert _{2}^{2}-\left\Vert
y\right\Vert _{1}^{2}}{\left\Vert y\right\Vert _{1}^{2}},\frac{\left\Vert
x\right\Vert _{2}^{2}-\left\Vert x\right\Vert _{1}^{2}}{\left\Vert
x\right\Vert _{1}^{2}}\right\} \delta \left( \left( \cdot ,\cdot \right)
_{1};x,y\right) .
\end{multline}
\end{corollary}

The following strong superadditivity property of $\delta \left( \cdot
;x,y\right) $ that is different from the one in Subsection \ref{ss2.2} holds
\cite{IHFDM2}:

\begin{corollary}[Dragomir-Mond, 1996]
\label{c1.2.2}If $\left( \cdot ,\cdot \right) _{i}\in \mathcal{JP}\left(
X\right) $ and $x,y\in X$ with $x,y\neq 0,$ then%
\begin{multline}
\delta \left( \left( \cdot ,\cdot \right) _{1}+\left( \cdot ,\cdot \right)
_{2};x,y\right) -\delta \left( \left( \cdot ,\cdot \right) _{1};x,y\right)
-\delta \left( \left( \cdot ,\cdot \right) _{2};x,y\right)  \label{1.2.19} \\
\geq \max \left\{ \left( \frac{\left\Vert y\right\Vert _{2}}{\left\Vert
y\right\Vert _{1}}\right) ^{2}\delta \left( \left( \cdot ,\cdot \right)
_{1};x,y\right) +\left( \frac{\left\Vert y\right\Vert _{1}}{\left\Vert
y\right\Vert _{2}}\right) ^{2}\delta \left( \left( \cdot ,\cdot \right)
_{2};x,y\right) ;\right. \\
\left. \left( \frac{\left\Vert x\right\Vert _{2}}{\left\Vert x\right\Vert
_{1}}\right) ^{2}\delta \left( \left( \cdot ,\cdot \right) _{1};x,y\right)
+\left( \frac{\left\Vert x\right\Vert _{1}}{\left\Vert x\right\Vert _{2}}%
\right) ^{2}\delta \left( \left( \cdot ,\cdot \right) _{2};x,y\right)
\right\} \qquad \left( \geq 0\right) .
\end{multline}
\end{corollary}

\begin{proof}
Utilising the identities (\ref{1.2.15}) and (\ref{1.2.16}) and taking into
account that $\gamma \left( \cdot ;x,y\right) $ is superadditive, we can
state that%
\begin{align}
& \delta \left( \left( \cdot ,\cdot \right) _{1}+\left( \cdot ,\cdot \right)
_{2};x,y\right)  \label{1.2.20} \\
& \geq \frac{\left\Vert y\right\Vert _{1}^{2}+\left\Vert y\right\Vert
_{2}^{2}}{\left\Vert y\right\Vert _{1}^{2}}\delta \left( \left( \cdot ,\cdot
\right) _{1};x,y\right) +\frac{\left\Vert y\right\Vert _{1}^{2}+\left\Vert
y\right\Vert _{2}^{2}}{\left\Vert y\right\Vert _{2}^{2}}\delta \left( \left(
\cdot ,\cdot \right) _{2};x,y\right)  \notag \\
& =\delta \left( \left( \cdot ,\cdot \right) _{1};x,y\right) +\delta \left(
\left( \cdot ,\cdot \right) _{2};x,y\right)  \notag \\
& \qquad \qquad +\left( \frac{\left\Vert y\right\Vert _{2}}{\left\Vert
y\right\Vert _{1}}\right) ^{2}\delta \left( \left( \cdot ,\cdot \right)
_{1};x,y\right) +\left( \frac{\left\Vert y\right\Vert _{1}}{\left\Vert
y\right\Vert _{2}}\right) ^{2}\delta \left( \left( \cdot ,\cdot \right)
_{2};x,y\right)  \notag
\end{align}%
and a similar inequality with $x$ instead of $y.$ These show that the
desired inequality (\ref{1.2.19}) holds true.
\end{proof}

\begin{remark}
\label{r1.2.3}Obviously, all the inequalities above remain true if $\left(
\cdot ,\cdot \right) _{i},$ $i=1,2$ are nonnegative Hermitian forms for
which we have $\left\Vert x\right\Vert _{i},$ $\left\Vert y\right\Vert
_{i}\neq 0.$
\end{remark}

Finally, we may state and prove the superadditivity result for the mapping $%
\beta $ (see \cite{IHFDM2}):

\begin{theorem}[Dragomir-Mond, 1996]
\label{t1.2.4}The mapping $\beta $ defined by (\ref{1.2.12}) is
superadditive on $\mathcal{H}\left( X\right) .$
\end{theorem}

\begin{proof}
Without loss of generality, if $\left( \cdot ,\cdot \right) _{i}\in \mathcal{%
H}\left( X\right) $ and $x,y\in X,$ we may assume, for instance, that $%
\left\Vert y\right\Vert _{i}\neq 0,$ $i=1,2.$

If so, then%
\begin{multline*}
\qquad \left( \frac{\left\Vert y\right\Vert _{2}}{\left\Vert y\right\Vert
_{1}}\right) ^{2}\delta \left( \left( \cdot ,\cdot \right) _{1};x,y\right)
+\left( \frac{\left\Vert y\right\Vert _{1}}{\left\Vert y\right\Vert _{2}}%
\right) ^{2}\delta \left( \left( \cdot ,\cdot \right) _{2};x,y\right) \\
\geq 2\left[ \delta \left( \left( \cdot ,\cdot \right) _{1};x,y\right)
\delta \left( \left( \cdot ,\cdot \right) _{2};x,y\right) \right] ^{\frac{1}{%
2}},\qquad
\end{multline*}%
and by making use of (\ref{1.2.20}) we get:%
\begin{equation*}
\delta \left( \left( \cdot ,\cdot \right) _{1}+\left( \cdot ,\cdot \right)
_{2};x,y\right) \geq \left\{ \left[ \delta \left( \left( \cdot ,\cdot
\right) _{1};x,y\right) \right] ^{\frac{1}{2}}+\left[ \delta \left( \left(
\cdot ,\cdot \right) _{2};x,y\right) \right] ^{\frac{1}{2}}\right\} ^{2},
\end{equation*}%
which is exactly the superadditivity property for $\beta .$
\end{proof}

\section{Applications for General Inner Product Spaces\label{s3}}

\subsection{\label{ss3.1}Inequalities for Orthonormal Families}

Let $\left( H;\left\langle \cdot ,\cdot \right\rangle \right) $ be an inner
product space over the real or complex number field $\mathbb{K}.$ The family
of vectors $E:=\left\{ e_{i}\right\} _{i\in I}$ ($I$ is a finite or
infinite) is an \textit{orthonormal family} of vectors if $\left\langle
e_{i},e_{j}\right\rangle =\delta _{ij}$ for $i,j\in I,$ where $\delta _{ij}$
is Kronecker's delta.

The following inequality is well known in the literature as Bessel's
inequality:%
\begin{equation}
\sum_{i\in F}\left\vert \left\langle x,e_{i}\right\rangle \right\vert
^{2}\leq \left\Vert x\right\Vert ^{2}  \label{1.3.1}
\end{equation}%
for any $F$ a finite part of $I$ and $x$ a vector in $H.$

If by $\mathcal{F}\left( I\right) $ we denote the family of all finite parts
of $I$ (including the empty set $\varnothing $), then for any $F\in \mathcal{%
F}\left( I\right) \backslash \left\{ \varnothing \right\} $ the functional $%
\left( \cdot ,\cdot \right) _{F}:H\times H\rightarrow \mathbb{K}$ given by%
\begin{equation}
\left( x,y\right) _{F}:=\sum_{i\in F}\left\langle x,e_{i}\right\rangle
\left\langle e_{i},y\right\rangle  \label{1.3.2}
\end{equation}%
is a Hermitian form on $H.$

It is obvious that if $F_{1},F_{2}\in \mathcal{F}\left( I\right) \backslash
\left\{ \varnothing \right\} $ and $F_{1}\cap F_{2}=\varnothing ,$ then $%
\left( \cdot ,\cdot \right) _{F_{1}\cup F_{2}}=\left( \cdot ,\cdot \right)
_{F_{1}}+\left( \cdot ,\cdot \right) _{F_{2}}.$

We can define the functional $\sigma :\mathcal{F}\left( I\right) \times
H^{2}\rightarrow \mathbb{R}_{+}$ by%
\begin{equation}
\sigma \left( F;x,y\right) :=\left\Vert x\right\Vert _{F}\left\Vert
y\right\Vert _{F}-\left\vert \left( x,y\right) _{F}\right\vert ,
\label{1.3.3}
\end{equation}%
where%
\begin{equation*}
\left\Vert x\right\Vert _{F}:=\left( \sum_{i\in F}\left\vert \left\langle
x,e_{i}\right\rangle \right\vert ^{2}\right) ^{\frac{1}{2}}=\left[ \left(
x,x\right) _{F}\right] ^{\frac{1}{2}},\qquad x\in H.
\end{equation*}

The following proposition may be stated (see also \cite{IHFDM2}):

\begin{proposition}[Dragomir-Mond, 1995]
\label{p1.3.1}The mapping $\sigma $ satisfies the following

\begin{enumerate}
\item[(i)] If $F_{1},F_{2}\in \mathcal{F}\left( I\right) \backslash \left\{
\varnothing \right\} $ with $F_{1}\cap F_{2}=\varnothing ,$ then%
\begin{equation*}
\sigma \left( F_{1}\cup F_{2};x,y\right) \geq \sigma \left( F_{1};x,y\right)
+\sigma \left( F_{2};x,y\right) \qquad \left( \geq 0\right)
\end{equation*}%
for any $x,y\in H,$ i.e., the mapping $\sigma \left( \cdot ;x,y\right) $ is
an index set superadditive mapping on $\mathcal{F}\left( I\right) ;$

\item[(ii)] If $\varnothing \neq F_{1}\subseteq F_{2},$ $F_{1},F_{2}\in
\mathcal{F}\left( I\right) ,$ then%
\begin{equation*}
\sigma \left( F_{2};x,y\right) \geq \sigma \left( F_{1};x,y\right) \qquad
\left( \geq 0\right) ,
\end{equation*}%
i.e., the mapping $\sigma \left( \cdot ;x,y\right) $ is an index set
monotonic mapping on $\mathcal{F}\left( I\right) .$
\end{enumerate}
\end{proposition}

The proof is obvious by Theorem \ref{t1.2.1} and we omit the details.

We can also define the mapping $\sigma _{r}\left( \cdot ;\cdot ,\cdot
\right) :\mathcal{F}\left( I\right) \times H^{2}\rightarrow \mathbb{R}_{+}$
by%
\begin{equation*}
\sigma _{r}\left( F;x,y\right) :=\left\Vert x\right\Vert _{F}\left\Vert
y\right\Vert _{F}-\func{Re}\left( x,y\right) _{F},
\end{equation*}%
which also has the properties (i) and (ii) of Proposition \ref{p1.3.1}.

Since, by Bessel's inequality the hermitian form $\left( \cdot ,\cdot
\right) _{F}\leq \left\langle \cdot ,\cdot \right\rangle $ in the sense of
definition (\ref{1.2.2}) then by Theorem \ref{t1.2.1} we may state the
following \textit{refinements} of Schwarz's inequality \cite{IHFDM}:

\begin{proposition}[Dragomir-Mond, 1994]
\label{p1.3.2}For any $F\in \mathcal{F}\left( I\right) \backslash \left\{
0\right\} ,$ we have the inequalities%
\begin{multline}
\left\Vert x\right\Vert \left\Vert y\right\Vert -\left\vert \left\langle
x,y\right\rangle \right\vert  \label{1.3.4} \\
\geq \left( \sum_{i\in F}\left\vert \left\langle x,e_{i}\right\rangle
\right\vert ^{2}\right) ^{\frac{1}{2}}\left( \sum_{i\in F}\left\vert
\left\langle y,e_{i}\right\rangle \right\vert ^{2}\right) ^{\frac{1}{2}%
}-\left\vert \sum_{i\in F}\left\langle x,e_{i}\right\rangle \left\langle
e_{i},y\right\rangle \right\vert
\end{multline}%
and%
\begin{multline}
\left\Vert x\right\Vert \left\Vert y\right\Vert -\left\vert \left\langle
x,y\right\rangle \right\vert  \label{1.3.5} \\
\geq \left( \left\Vert x\right\Vert ^{2}-\sum_{i\in F}\left\vert
\left\langle x,e_{i}\right\rangle \right\vert ^{2}\right) ^{\frac{1}{2}%
}\left( \left\Vert y\right\Vert ^{2}-\sum_{i\in F}\left\vert \left\langle
y,e_{i}\right\rangle \right\vert ^{2}\right) ^{\frac{1}{2}} \\
-\left\vert \left\langle x,y\right\rangle -\sum_{i\in F}\left\langle
x,e_{i}\right\rangle \left\langle e_{i},y\right\rangle \right\vert
\end{multline}%
and the corresponding versions on replacing $\left\vert \cdot \right\vert $
by $\func{Re}\left( \cdot \right) ,$ where $x,y$ are vectors in $H.$
\end{proposition}

\begin{remark}
\label{r1.3.4}Note that the inequality (\ref{1.3.4}) and its version for $%
\func{Re}\left( \cdot \right) $ has been established for the first time and
utilising a different argument by Dragomir and S\'{a}ndor in 1994 (see \cite[%
Theorem 5 and Remark 2]{IHFDS}).
\end{remark}

If we now define the mapping $\delta :\mathcal{F}\left( I\right) \times
H^{2}\rightarrow \mathbb{R}_{+}$ by%
\begin{equation*}
\delta \left( F;x,y\right) :=\left\Vert x\right\Vert _{F}^{2}\left\Vert
y\right\Vert _{F}^{2}-\left\vert \left( x,y\right) _{F}\right\vert ^{2}
\end{equation*}%
and making use of Theorem \ref{t1.2.2}, we may state the following result
\cite{IHFDM2}.

\begin{proposition}[Dragomir-Mond, 1995]
\label{p1.3.3}The mapping $\delta $ satisfies the following properties:

\begin{enumerate}
\item[(i)] If $F_{1},F_{2}\in \mathcal{F}\left( I\right) \backslash \left\{
\varnothing \right\} $ with $F_{1}\cap F_{2}=\varnothing ,$ then%
\begin{multline}
\delta \left( F_{1}\cup F_{2};x,y\right) -\delta \left( F_{1};x,y\right)
-\delta \left( F_{2};x,y\right)  \label{1.3.6} \\
\geq \left( \det \left[
\begin{array}{ll}
\left\Vert x\right\Vert _{F_{1}} & \left\Vert y\right\Vert _{F_{1}} \\
\left\Vert x\right\Vert _{F_{2}} & \left\Vert y\right\Vert _{F_{2}}%
\end{array}%
\right] \right) ^{2}\qquad \left( \geq 0\right) ,
\end{multline}%
i.e., the mapping $\delta \left( \cdot ;x,y\right) $ is strong superadditive
as an index set mapping;

\item[(ii)] If $\varnothing \neq F_{1}\subseteq F_{2},$ $F_{1},F_{2}\in
\mathcal{F}\left( I\right) ,$ then%
\begin{multline}
\delta \left( F_{2};x,y\right) -\delta \left( F_{1};x,y\right)  \label{1.3.7}
\\
\geq \left( \det \left[
\begin{array}{cc}
\left\Vert x\right\Vert _{F_{1}} & \left\Vert y\right\Vert _{F_{1}} \\
\left( \left\Vert x\right\Vert _{F_{2}}^{2}-\left\Vert x\right\Vert
_{F_{1}}^{2}\right) ^{\frac{1}{2}} & \left( \left\Vert y\right\Vert
_{F_{2}}^{2}-\left\Vert y\right\Vert _{F_{1}}^{2}\right) ^{\frac{1}{2}}%
\end{array}%
\right] \right) ^{2}\qquad \left( \geq 0\right) ,
\end{multline}%
i.e., the mapping $\delta \left( \cdot ;x,y\right) $ is strong nondecreasing
as an index set mapping.
\end{enumerate}
\end{proposition}

On applying the same general result in Theorem \ref{t1.2.2}, (ii) for the
hermitian functionals $\left( \cdot ,\cdot \right) _{F}$ $\left( F\in
\mathcal{F}\left( I\right) \backslash \left\{ \varnothing \right\} \right) $
and $\left\langle \cdot ,\cdot \right\rangle $ for which, by Bessel's
inequality we know that $\left( \cdot ,\cdot \right) _{F}\leq \left\langle
\cdot ,\cdot \right\rangle ,$ we may state the following result as well,
which provides refinements for the Schwarz inequality.

\begin{proposition}[Dragomir-Mond, 1994]
\label{p1.3.4}For any $F\in \mathcal{F}\left( I\right) \backslash \left\{
\varnothing \right\} $, we have the inequalities:%
\begin{multline}
\left\Vert x\right\Vert ^{2}\left\Vert y\right\Vert ^{2}-\left\vert
\left\langle x,y\right\rangle \right\vert ^{2}  \label{1.3.8} \\
\geq \sum_{i\in F}\left\vert \left\langle x,e_{i}\right\rangle \right\vert
^{2}\sum_{i\in F}\left\vert \left\langle y,e_{i}\right\rangle \right\vert
^{2}-\left\vert \sum_{i\in F}\left\langle x,e_{i}\right\rangle \left\langle
e_{i},y\right\rangle \right\vert ^{2}\qquad \left( \geq 0\right)
\end{multline}%
and%
\begin{multline}
\left\Vert x\right\Vert ^{2}\left\Vert y\right\Vert ^{2}-\left\vert
\left\langle x,y\right\rangle \right\vert ^{2}  \label{1.3.9} \\
\geq \left( \left\Vert x\right\Vert ^{2}-\sum_{i\in F}\left\vert
\left\langle x,e_{i}\right\rangle \right\vert ^{2}\right) \left( \left\Vert
y\right\Vert ^{2}-\sum_{i\in F}\left\vert \left\langle y,e_{i}\right\rangle
\right\vert ^{2}\right) \\
-\left\vert \left\langle x,y\right\rangle -\sum_{i\in F}\left\langle
x,e_{i}\right\rangle \left\langle e_{i},y\right\rangle \right\vert
^{2}\qquad \left( \geq 0\right) ,
\end{multline}%
for any $x,y\in H.$
\end{proposition}

On utilising Corollary \ref{c1.2.2} we may state the following different
superadditivity property for the mapping $\delta \left( \cdot ;x,y\right) .$

\begin{proposition}
\label{p1.3.5}If $F_{1},F_{2}\in \mathcal{F}\left( I\right) \backslash
\left\{ \varnothing \right\} $ with $F_{1}\cap F_{2}=\varnothing ,$ then%
\begin{multline}
\delta \left( F_{1}\cup F_{2};x,y\right) -\delta \left( F_{1};x,y\right)
-\delta \left( F_{2};x,y\right)  \label{1.3.10} \\
\geq \max \left\{ \left( \frac{\left\Vert y\right\Vert _{F_{2}}}{\left\Vert
y\right\Vert _{F_{1}}}\right) ^{2}\delta \left( F_{1};x,y\right) +\left(
\frac{\left\Vert y\right\Vert _{F_{1}}}{\left\Vert y\right\Vert _{F_{2}}}%
\right) ^{2}\delta \left( F_{2};x,y\right) ;\right. \\
\left. \left( \frac{\left\Vert x\right\Vert _{F_{2}}}{\left\Vert
x\right\Vert _{F_{1}}}\right) ^{2}\delta \left( F_{1};x,y\right) +\left(
\frac{\left\Vert x\right\Vert _{F_{1}}}{\left\Vert x\right\Vert _{F_{2}}}%
\right) ^{2}\delta \left( F_{2};x,y\right) \right\} \qquad \left( \geq
0\right)
\end{multline}%
for any $x,y\in H\backslash \left\{ 0\right\} .$
\end{proposition}

Further, for $y\notin M^{\perp }$ where $M=Sp\left\{ e_{i}\right\} _{i\in I}$
is the linear space spanned by $E=\left\{ e_{i}\right\} _{i\in I},$ we can
also consider the functional $\gamma :\mathcal{F}\left( I\right) \times
H^{2}\rightarrow \mathbb{R}_{+}$ defined by%
\begin{equation*}
\gamma \left( F;x,y\right) :=\frac{\delta \left( F;x,y\right) }{\left\Vert
y\right\Vert _{F}^{2}}=\frac{\left\Vert x\right\Vert _{F}^{2}\left\Vert
y\right\Vert _{F}^{2}-\left\vert \left( x,y\right) _{F}\right\vert ^{2}}{%
\left\Vert y\right\Vert _{F}^{2}},
\end{equation*}%
where $x\in H$ and $F\neq \varnothing .$

Utilising Theorem \ref{t1.2.3}, we may state the following result concerning
the properties of the functional $\gamma \left( \cdot ;x,y\right) $ with $x$
and $y$ as above.

\begin{proposition}
\label{p1.3.7}For any $x\in H$ and $y\in H\backslash M^{\perp },$ the
functional $\gamma \left( \cdot ;x,y\right) $ is superadditive and monotonic
nondecreasing as an index set mapping on $\mathcal{F}\left( I\right) .$
\end{proposition}

Since $\left\langle \cdot ,\cdot \right\rangle \geq \left( \cdot ,\cdot
\right) _{F}$ for any $F\in \mathcal{F}\left( I\right) ,$ on making use of
Corollary \ref{c1.2.1}, we may state the following refinement of Schwarz's
inequality:

\begin{proposition}
\label{p1.3.8}Let $x\in H$ and $y\in H\backslash M_{F}^{\perp },$ where $%
M_{F}:=Sp\left\{ e_{i}\right\} _{i\in I}$ and $F\in \mathcal{F}\left(
I\right) \backslash \left\{ \varnothing \right\} $ is given. Then%
\begin{multline}
\left\Vert x\right\Vert ^{2}\left\Vert y\right\Vert ^{2}-\left\vert
\left\langle x,y\right\rangle \right\vert ^{2}\geq \max \left\{ \frac{%
\left\Vert y\right\Vert ^{2}}{\sum_{i\in F}\left\vert \left\langle
y,e_{i}\right\rangle \right\vert ^{2}},\frac{\left\Vert x\right\Vert ^{2}}{%
\sum_{i\in F}\left\vert \left\langle x,e_{i}\right\rangle \right\vert ^{2}}%
\right\}  \label{1.3.11} \\
\times \left( \sum_{i\in F}\left\vert \left\langle x,e_{i}\right\rangle
\right\vert ^{2}\sum_{i\in F}\left\vert \left\langle y,e_{i}\right\rangle
\right\vert ^{2}-\left\vert \sum_{i\in F}\left\langle x,e_{i}\right\rangle
\left\langle e_{i},y\right\rangle \right\vert ^{2}\right) \\
\left( \geq \sum_{i\in F}\left\vert \left\langle x,e_{i}\right\rangle
\right\vert ^{2}\sum_{i\in F}\left\vert \left\langle y,e_{i}\right\rangle
\right\vert ^{2}-\left\vert \sum_{i\in F}\left\langle x,e_{i}\right\rangle
\left\langle e_{i},y\right\rangle \right\vert ^{2}\right) ,
\end{multline}%
which is a refinement of (\ref{1.3.9}) in the case that $y\in H\backslash
M_{F}^{\perp }.$
\end{proposition}

Finally, consider the functional $\beta :\mathcal{F}\left( I\right) \times
H^{2}\rightarrow \mathbb{R}_{+}$ given by%
\begin{equation*}
\beta \left( F;x,y\right) :=\left[ \delta \left( F;x,y\right) \right] ^{%
\frac{1}{2}}=\left( \left\Vert x\right\Vert _{F}^{2}\left\Vert y\right\Vert
_{F}^{2}-\left\vert \left( x,y\right) _{F}\right\vert ^{2}\right) ^{\frac{1}{%
2}}.
\end{equation*}%
Utilising Theorem \ref{t1.2.4}, we may state the following:

\begin{proposition}
\label{p1.3.9}The functional $\beta \left( \cdot ;x,y\right) $ is
superadditive as an index set mapping on $\mathcal{F}\left( I\right) $ for
each $x,y\in H.$
\end{proposition}

As a dual approach, one may also consider the following form $\left( \cdot
,\cdot \right) _{C,F}:H\times H\rightarrow \mathbb{R}$ given by:%
\begin{equation}
\left( x,y\right) _{C,F}:=\left\langle x,y\right\rangle -\left( x,y\right)
_{F}=\left\langle x,y\right\rangle -\sum_{i\in F}\left\langle
x,e_{i}\right\rangle \left\langle e_{i},y\right\rangle .  \label{1.3.12}
\end{equation}%
By Bessel's inequality, we observe that $\left( \cdot ,\cdot \right) _{C,F}$
is a nonnegative hermitian form and, obviously%
\begin{equation*}
\left( \cdot ,\cdot \right) _{I}+\left( \cdot ,\cdot \right)
_{C,F}=\left\langle \cdot ,\cdot \right\rangle .
\end{equation*}%
Utilising the superadditivity properties from Section \ref{s2}, one may
state the following refinement of the Schwarz inequality:%
\begin{multline}
\left\Vert x\right\Vert \left\Vert y\right\Vert -\left\vert \left\langle
x,y\right\rangle \right\vert  \label{1.3.13} \\
\geq \left( \sum_{i\in F}\left\vert \left\langle x,e_{i}\right\rangle
\right\vert ^{2}\sum_{i\in F}\left\vert \left\langle y,e_{i}\right\rangle
\right\vert ^{2}\right) ^{\frac{1}{2}}-\left\vert \sum_{i\in F}\left\langle
x,e_{i}\right\rangle \left\langle e_{i},y\right\rangle \right\vert \\
+\left( \left\Vert x\right\Vert ^{2}-\sum_{i\in F}\left\vert \left\langle
x,e_{i}\right\rangle \right\vert ^{2}\right) ^{\frac{1}{2}}\left( \left\Vert
y\right\Vert ^{2}-\sum_{i\in F}\left\vert \left\langle y,e_{i}\right\rangle
\right\vert ^{2}\right) ^{\frac{1}{2}} \\
-\left\vert \left\langle x,y\right\rangle -\sum_{i\in F}\left\langle
x,e_{i}\right\rangle \left\langle e_{i},y\right\rangle \right\vert \qquad
\left( \geq 0\right) ,
\end{multline}%
\begin{multline}
\left\Vert x\right\Vert ^{2}\left\Vert y\right\Vert ^{2}-\left\vert
\left\langle x,y\right\rangle \right\vert ^{2}  \label{1.3.14} \\
\geq \sum_{i\in F}\left\vert \left\langle x,e_{i}\right\rangle \right\vert
^{2}\sum_{i\in F}\left\vert \left\langle y,e_{i}\right\rangle \right\vert
^{2}-\left\vert \sum_{i\in F}\left\langle x,e_{i}\right\rangle \left\langle
e_{i},y\right\rangle \right\vert ^{2} \\
+\left( \left\Vert x\right\Vert ^{2}-\sum_{i\in F}\left\vert \left\langle
x,e_{i}\right\rangle \right\vert ^{2}\right) \left( \left\Vert y\right\Vert
^{2}-\sum_{i\in F}\left\vert \left\langle y,e_{i}\right\rangle \right\vert
^{2}\right) \\
-\left\vert \left\langle x,y\right\rangle -\sum_{i\in F}\left\langle
x,e_{i}\right\rangle \left\langle e_{i},y\right\rangle \right\vert
^{2}\qquad \left( \geq 0\right)
\end{multline}%
and
\begin{multline}
\left( \left\Vert x\right\Vert ^{2}\left\Vert y\right\Vert ^{2}-\left\vert
\left\langle x,y\right\rangle \right\vert ^{2}\right) ^{\frac{1}{2}}
\label{1.3.15} \\
\geq \left[ \sum_{i\in F}\left\vert \left\langle x,e_{i}\right\rangle
\right\vert ^{2}\sum_{i\in F}\left\vert \left\langle y,e_{i}\right\rangle
\right\vert ^{2}-\left\vert \sum_{i\in F}\left\langle x,e_{i}\right\rangle
\left\langle e_{i},y\right\rangle \right\vert ^{2}\right] ^{\frac{1}{2}} \\
+\left[ \left( \left\Vert x\right\Vert ^{2}-\sum_{i\in F}\left\vert
\left\langle x,e_{i}\right\rangle \right\vert ^{2}\right) \left( \left\Vert
y\right\Vert ^{2}-\sum_{i\in F}\left\vert \left\langle y,e_{i}\right\rangle
\right\vert ^{2}\right) \right. \\
-\left. \left\vert \left\langle x,y\right\rangle -\sum_{i\in F}\left\langle
x,e_{i}\right\rangle \left\langle e_{i},y\right\rangle \right\vert ^{2}
\right] ^{\frac{1}{2}}\qquad \left( \geq 0\right) ,
\end{multline}%
for any $x,y\in H$ and $F\in \mathcal{F}\left( I\right) \backslash \left\{
\varnothing \right\} .$

\subsection{Inequalities for Gram Determinants}

Let $\left\{ x_{1},\dots ,x_{n}\right\} $ be vectors in the inner product
space $\left( H,\left\langle \cdot ,\cdot \right\rangle \right) $ over the
real or complex number field $\mathbb{K}$. Consider the \textit{gram matrix}
associated to the above vectors:%
\begin{equation*}
G\left( x_{1},\dots ,x_{n}\right) :=\left[
\begin{array}{cccc}
\left\langle x_{1},x_{1}\right\rangle & \left\langle x_{1},x_{2}\right\rangle
& \cdots & \left\langle x_{1},x_{n}\right\rangle \\
\left\langle x_{2},x_{1}\right\rangle &  & \cdots & \left\langle
x_{2},x_{n}\right\rangle \\
\cdots &  & \cdots & \cdots \\
\left\langle x_{n},x_{1}\right\rangle & \left\langle x_{n},x_{2}\right\rangle
& \cdots & \left\langle x_{n},x_{n}\right\rangle%
\end{array}%
\right] .
\end{equation*}%
The determinant%
\begin{equation*}
\Gamma \left( x_{1},\dots ,x_{n}\right) :=\det G\left( x_{1},\dots
,x_{n}\right)
\end{equation*}%
is called the Gram determinant associated to the system $\left\{ x_{1},\dots
,x_{n}\right\} .$

If $\left\{ x_{1},\dots ,x_{n}\right\} $ does not contain the null vector $%
0, $ then \cite{IHFFD}%
\begin{equation}
0\leq \Gamma \left( x_{1},\dots ,x_{n}\right) \leq \left\Vert
x_{1}\right\Vert ^{2}\left\Vert x_{2}\right\Vert ^{2}\cdots \left\Vert
x_{n}\right\Vert ^{2}.  \label{1.3.16}
\end{equation}%
The equality holds on the left (respectively right) side of (\ref{1.3.16})
if and only if $\left\{ x_{1},\dots ,x_{n}\right\} $ is linearly dependent
(respectively orthogonal). The first inequality in (\ref{1.3.16}) is known
in the literature as \textit{Gram's inequality} while the second one is
known as \textit{Hadamard's inequality.}

The following result obtained in \cite{IHFDS} may be regarded as a
refinement of Gram's inequality:

\begin{theorem}[Dragomir-S\'{a}ndor, 1994]
\label{t1.3.8}Let $\left\{ x_{1},\dots ,x_{n}\right\} $ be a system of
nonzero vectors in $H.$ Then for any $x,y\in H$ one has:%
\begin{equation}
\Gamma \left( x,x_{1},\dots ,x_{n}\right) \Gamma \left( y,x_{1},\dots
,x_{n}\right) \geq \left\vert \Gamma \left( x_{1},\dots ,x_{n}\right) \left(
x,y\right) \right\vert ^{2},  \label{1.3.17}
\end{equation}%
where $\Gamma \left( x_{1},\dots ,x_{n}\right) \left( x,y\right) $ is
defined by:%
\begin{multline*}
\Gamma \left( x_{1},\dots ,x_{n}\right) \left( x,y\right) \\
:=\det \left[
\begin{array}{cccc}
\left\langle x,y\right\rangle & \left\langle x,x_{1}\right\rangle & \cdots &
\left\langle x,x_{n}\right\rangle \\
\left\langle x_{1},y\right\rangle &  &  &  \\
\cdots &  & G\left( x_{1},\dots ,x_{n}\right) &  \\
\left\langle x_{n},y\right\rangle &  &  &
\end{array}%
\right] .
\end{multline*}
\end{theorem}

\begin{proof}
We follow the proof from \cite{IHFDS}.

Let us consider the mapping $p:H\times H\rightarrow \mathbb{K}$ given by%
\begin{equation*}
p\left( x,y\right) =\Gamma \left( x_{1},\dots ,x_{n}\right) \left(
x,y\right) .
\end{equation*}%
Utilising the properties of determinants, we notice that%
\begin{align*}
p\left( x,y\right) & =\Gamma \left( x,x_{1},\dots ,x_{n}\right) \geq 0, \\
p\left( x+y,z\right) & =\Gamma \left( x_{1},\dots ,x_{n}\right) \left(
x+y,z\right) \\
& =\Gamma \left( x_{1},\dots ,x_{n}\right) \left( x,z\right) +\Gamma \left(
x_{1},\dots ,x_{n}\right) \left( y,z\right) \\
& =p\left( x,z\right) +p\left( y,z\right) , \\
p\left( \alpha x,y\right) & =\alpha p\left( x,y\right) , \\
p\left( y,x\right) & =\overline{p\left( x,y\right) },
\end{align*}%
for any $x,y,z\in H$ and $\alpha \in \mathbb{K}$, showing that $p\left(
\cdot ,\cdot \right) $ is a nonnegative hermitian from on $X.$ Writing
Schwarz's inequality for $p\left( \cdot ,\cdot \right) $ we deduce the
desired result (\ref{1.3.17}).
\end{proof}

In a similar manner, if we define $q:H\times H\rightarrow \mathbb{K}$ by%
\begin{align*}
q\left( x,y\right) & :=\left( x,y\right) \prod_{i=1}^{n}\left\Vert
x_{i}\right\Vert ^{2}-p\left( x,y\right) \\
& =\left( x,y\right) \prod_{i=1}^{n}\left\Vert x_{i}\right\Vert ^{2}-\Gamma
\left( x_{1},\dots ,x_{n}\right) \left( x,y\right) ,
\end{align*}%
then, using Hadamard's inequality, we conclude that $q\left( \cdot ,\cdot
\right) $ is also a nonnegative hermitian form. Therefore, by Schwarz's
inequality applied for $q\left( \cdot ,\cdot \right) ,$ we can state the
following result as well \cite{IHFDS}:

\begin{theorem}[Dragomir-S\'{a}ndor, 1994]
\label{t1.3.9}With the assumptions of Theorem \ref{t1.3.8}, we have:%
\begin{multline}
\left[ \left\Vert x\right\Vert ^{2}\prod_{i=1}^{n}\left\Vert
x_{i}\right\Vert ^{2}-\Gamma \left( x,x_{1},\dots ,x_{n}\right) \right]
\label{1.3.18} \\
\times \left[ \left\Vert y\right\Vert ^{2}\prod_{i=1}^{n}\left\Vert
x_{i}\right\Vert ^{2}-\Gamma \left( y,x_{1},\dots ,x_{n}\right) \right] \\
\geq \left\vert \left\langle x,y\right\rangle \prod_{i=1}^{n}\left\Vert
x_{i}\right\Vert ^{2}-\Gamma \left( x_{1},\dots ,x_{n}\right) \left(
x,y\right) \right\vert ^{2},
\end{multline}%
for each $x,y\in H.$
\end{theorem}

Observing that, for a given set of nonzero vectors $\left\{ x_{1},\dots
,x_{n}\right\} ,$%
\begin{equation*}
p\left( x,y\right) +q\left( x,y\right) =\left( x,y\right)
\prod_{i=1}^{n}\left\Vert x_{i}\right\Vert ^{2},
\end{equation*}%
for any $x,y\in H,$ then, on making use of the superadditivity properties of
the various functionals defined in Section \ref{s2}, we can state the
following refinements of the Schwarz inequality in inner product spaces:%
\begin{multline}
\left[ \left\Vert x\right\Vert \left\Vert y\right\Vert -\left\vert
\left\langle x,y\right\rangle \right\vert \right] \prod_{i=1}^{n}\left\Vert
x_{i}\right\Vert ^{2}  \label{1.3.19} \\
\geq \left[ \Gamma \left( x,x_{1},\dots ,x_{n}\right) \Gamma \left(
y,x_{1},\dots ,x_{n}\right) \right] ^{\frac{1}{2}}-\left\vert \Gamma \left(
x_{1},\dots ,x_{n}\right) \left( x,y\right) \right\vert \\
+\left[ \left\Vert x\right\Vert ^{2}\prod_{i=1}^{n}\left\Vert
x_{i}\right\Vert ^{2}-\Gamma \left( x,x_{1},\dots ,x_{n}\right) \right] ^{%
\frac{1}{2}} \\
\times \left[ \left\Vert y\right\Vert ^{2}\prod_{i=1}^{n}\left\Vert
x_{i}\right\Vert ^{2}-\Gamma \left( y,x_{1},\dots ,x_{n}\right) \right] ^{%
\frac{1}{2}} \\
-\left\vert \left\langle x,y\right\rangle \prod_{i=1}^{n}\left\Vert
x_{i}\right\Vert ^{2}-\Gamma \left( x_{1},\dots ,x_{n}\right) \left(
x,y\right) \right\vert \qquad \left( \geq 0\right) ,
\end{multline}%
\begin{multline}
\left[ \left\Vert x\right\Vert ^{2}\left\Vert y\right\Vert ^{2}-\left\vert
\left\langle x,y\right\rangle \right\vert ^{2}\right] \prod_{i=1}^{n}\left%
\Vert x_{i}\right\Vert ^{4}  \label{1.3.20} \\
\Gamma \left( x,x_{1},\dots ,x_{n}\right) \Gamma \left( y,x_{1},\dots
,x_{n}\right) -\left\vert \Gamma \left( x_{1},\dots ,x_{n}\right) \left(
x,y\right) \right\vert ^{2} \\
+\left[ \left\Vert x\right\Vert ^{2}\prod_{i=1}^{n}\left\Vert
x_{i}\right\Vert ^{2}-\Gamma \left( x,x_{1},\dots ,x_{n}\right) \right]\displaybreak\\
\times \left[ \left\Vert y\right\Vert
^{2}\prod_{i=1}^{n}\left\Vert
x_{i}\right\Vert ^{2}-\Gamma \left( y,x_{1},\dots ,x_{n}\right) \right] \\
-\left\vert \left\langle x,y\right\rangle \prod_{i=1}^{n}\left\Vert
x_{i}\right\Vert ^{2}-\Gamma \left( x_{1},\dots ,x_{n}\right) \left(
x,y\right) \right\vert ^{2}\qquad \left( \geq 0\right) ,
\end{multline}%
and%
\begin{multline}
\left[ \left\Vert x\right\Vert \left\Vert y\right\Vert -\left\vert
\left\langle x,y\right\rangle \right\vert \right] ^{\frac{1}{2}%
}\prod_{i=1}^{n}\left\Vert x_{i}\right\Vert ^{2}  \label{1.3.21} \\
\geq \left[ \Gamma \left( x,x_{1},\dots ,x_{n}\right) \Gamma \left(
y,x_{1},\dots ,x_{n}\right) -\left\vert \Gamma \left( x_{1},\dots
,x_{n}\right) \left( x,y\right) \right\vert ^{2}\right] ^{\frac{1}{2}} \\
+\left\{ \left[ \left\Vert x\right\Vert ^{2}\prod_{i=1}^{n}\left\Vert
x_{i}\right\Vert ^{2}-\Gamma \left( x,x_{1},\dots ,x_{n}\right) \right]
\right. \\
\times \left[ \left\Vert y\right\Vert ^{2}\prod_{i=1}^{n}\left\Vert
x_{i}\right\Vert ^{2}-\Gamma \left( y,x_{1},\dots ,x_{n}\right) \right] \\
-\left. \left\vert \left\langle x,y\right\rangle \prod_{i=1}^{n}\left\Vert
x_{i}\right\Vert ^{2}-\Gamma \left( x_{1},\dots ,x_{n}\right) \left(
x,y\right) \right\vert ^{2}\right\} ^{\frac{1}{2}}\qquad \left( \geq
0\right) .
\end{multline}

\subsection{Inequalities for Linear Operators}

Let $A:H\rightarrow H$ be a linear bounded operator and%
\begin{equation*}
\left\Vert A\right\Vert :=\sup \left\{ \left\Vert Ax\right\Vert ,\left\Vert
x\right\Vert <1\right\}
\end{equation*}%
its norm.

If we consider the hermitian forms $\left( \cdot ,\cdot \right) _{2},$ $%
\left( \cdot ,\cdot \right) _{1}:H\rightarrow H$ defined by%
\begin{equation*}
\left( x,y\right) _{1}:=\left\langle Ax,Ay\right\rangle ,\qquad \left(
x,y\right) _{2}:=\left\Vert A\right\Vert ^{2}\left\langle x,y\right\rangle
\end{equation*}%
then obviously $\left( \cdot ,\cdot \right) _{2}\geq \left( \cdot ,\cdot
\right) _{1}$ in the sense of definition (\ref{1.2.2}) and utilising the
monotonicity properties of the functional considered in Section \ref{s2}, we
may state the following inequalities:%
\begin{equation}
\left\Vert A\right\Vert ^{2}\left[ \left\Vert x\right\Vert \left\Vert
y\right\Vert -\left\vert \left\langle x,y\right\rangle \right\vert \right]
\geq \left\Vert Ax\right\Vert \left\Vert Ay\right\Vert -\left\vert
\left\langle Ax,Ay\right\rangle \right\vert \qquad \left( \geq 0\right) ,
\label{1.3.22}
\end{equation}%
\begin{multline}
\left\Vert A\right\Vert ^{4}\left[ \left\Vert x\right\Vert ^{2}\left\Vert
y\right\Vert ^{2}-\left\vert \left\langle x,y\right\rangle \right\vert ^{2}%
\right]  \label{1.3.23} \\
\geq \left\Vert Ax\right\Vert ^{2}\left\Vert Ay\right\Vert ^{2}-\left\vert
\left\langle Ax,Ay\right\rangle \right\vert ^{2}\qquad \left( \geq 0\right)
\end{multline}%
for any $x,y\in H$, and the corresponding versions on replacing $\left\vert
\cdot \right\vert $ by $\func{Re}\left( \cdot \right) .$

The results (\ref{1.3.22}) and (\ref{1.3.23}) have been obtained by Dragomir
and Mond in \cite{IHFDM}.

On using Corollary \ref{c1.2.1}, we may deduce the following inequality as
well:%
\begin{multline}
\left\Vert A\right\Vert ^{2}\left[ \left\Vert x\right\Vert ^{2}\left\Vert
y\right\Vert ^{2}-\left\vert \left\langle x,y\right\rangle \right\vert ^{2}%
\right]  \label{1.3.24} \\
\geq \max \left\{ \frac{\left\Vert x\right\Vert ^{2}}{\left\Vert
Ax\right\Vert ^{2}},\frac{\left\Vert y\right\Vert ^{2}}{\left\Vert
Ay\right\Vert ^{2}}\right\} \left[ \left\Vert Ax\right\Vert ^{2}\left\Vert
Ay\right\Vert ^{2}-\left\vert \left\langle Ax,Ay\right\rangle \right\vert
^{2}\right] \qquad \left( \geq 0\right)
\end{multline}%
for any $x,y\in H$ with $Ax,Ay\neq 0;$ which improves (\ref{1.3.23}) for $%
x,y $ specified before.

Similarly, if $B:H\rightarrow H$ is a linear operator satisfying the
condition%
\begin{equation}
\left\Vert Bx\right\Vert \geq m\left\Vert x\right\Vert \quad \text{for any}%
\quad x\in H,  \label{1.3.25}
\end{equation}%
where $m>0$ is given, then the hermitian forms $\left[ x,y\right]
_{2}:=\left\langle Bx,By\right\rangle ,$ $\left[ x,y\right]
_{1}:=m^{2}\left\langle x,y\right\rangle ,$ have the property that $\left[
\cdot ,\cdot \right] _{2}\geq \left[ \cdot ,\cdot \right] _{1}.$ Therefore,
from the monotonicity results established in Section \ref{s2}, we can state
that%
\begin{equation}
\left\Vert Bx\right\Vert \left\Vert By\right\Vert -\left\vert \left\langle
Bx,By\right\rangle \right\vert \geq m^{2}\left[ \left\Vert x\right\Vert
\left\Vert y\right\Vert -\left\vert \left\langle x,y\right\rangle
\right\vert \right] \qquad \left( \geq 0\right) ,  \label{1.3.26}
\end{equation}%
\begin{multline}
\qquad \left\Vert Bx\right\Vert ^{2}\left\Vert By\right\Vert ^{2}-\left\vert
\left\langle Bx,By\right\rangle \right\vert ^{2}  \label{1.3.27} \\
\geq m^{4}\left[ \left\Vert x\right\Vert ^{2}\left\Vert y\right\Vert
^{2}-\left\vert \left\langle x,y\right\rangle \right\vert ^{2}\right] \quad
\left( \geq 0\right) \qquad
\end{multline}%
for any $x,y\in H,$ and the corresponding results on replacing $\left\vert
\cdot \right\vert $ by $\func{Re}\left( \cdot \right) .$

The same Corollary \ref{c1.2.1}, would give the inequality%
\begin{multline}
\left\Vert Bx\right\Vert ^{2}\left\Vert By\right\Vert ^{2}-\left\vert
\left\langle Bx,By\right\rangle \right\vert ^{2}  \label{1.3.28} \\
\geq m^{2}\max \left\{ \frac{\left\Vert Bx\right\Vert ^{2}}{\left\Vert
x\right\Vert ^{2}},\frac{\left\Vert By\right\Vert ^{2}}{\left\Vert
y\right\Vert ^{2}}\right\} \left[ \left\Vert x\right\Vert ^{2}\left\Vert
y\right\Vert ^{2}-\left\vert \left\langle x,y\right\rangle \right\vert ^{2}%
\right]
\end{multline}%
for $x,y\neq 0,$ which is an improvement of (\ref{1.3.27}).

We recall that a linear self-adjoint operator $P:H\rightarrow H$ is \textit{%
nonnegative} if $\left\langle Px,x\right\rangle \geq 0$ for any $x\in H.$ $P$
is called \textit{positive} if $\left\langle Px,x\right\rangle =0$ and
\textit{positive definite with the constant }$\gamma >0$ if $\left\langle
Px,x\right\rangle \geq \gamma \left\Vert x\right\Vert ^{2}$ for any $x\in H.$

If $A,B:H\rightarrow H$ are two linear self-adjoint operators such that $%
A\geq B$ (this means that $A-B$ is nonnegative), then the corresponding
hermitian forms $\left( x,y\right) _{A}:=\left\langle Ax,y\right\rangle $
and $\left( x,y\right) _{B}:=\left\langle Bx,y\right\rangle $ satisfies the
property that $\left( \cdot ,\cdot \right) _{A}\geq \left( \cdot ,\cdot
\right) _{B}.$

If by $\mathcal{P}\left( H\right) $ we denote the \textit{cone }\ of all
linear self-adjoint and nonnegative operators defined in the Hilbert space $%
H,$ then, on utilising the results of Section \ref{s2}, we may state that
the functionals $\sigma _{0},\delta _{0},\beta _{0}:\mathcal{P}\left(
H\right) \times H^{2}\rightarrow \left[ 0,\infty \right] $ given by%
\begin{align*}
\sigma _{0}\left( P;x,y\right) & :=\left\langle Ax,x\right\rangle ^{\frac{1}{%
2}}\left\langle Py,y\right\rangle ^{\frac{1}{2}}-\left\vert \left\langle
Px,y\right\rangle \right\vert , \\
\delta _{0}\left( P;x,y\right) & :=\left\langle Px,x\right\rangle
\left\langle Py,y\right\rangle -\left\vert \left\langle Px,y\right\rangle
\right\vert ^{2}, \\
\beta _{0}\left( P;x,y\right) & :=\left[ \left\langle Px,x\right\rangle
\left\langle Py,y\right\rangle -\left\vert \left\langle Px,y\right\rangle
\right\vert ^{2}\right] ^{\frac{1}{2}},
\end{align*}%
are \textit{superadditive }and \textit{monotonic decreasing }on $\mathcal{P}%
\left( H\right) ,$ i.e.,%
\begin{equation*}
\gamma _{0}\left( P+Q;x,y\right) \geq \gamma _{0}\left( P;x,y\right) +\gamma
_{0}\left( Q;x,y\right) \qquad \left( \geq 0\right)
\end{equation*}%
for any $P,Q\in \mathcal{P}\left( H\right) $ and $x,y\in H,$ and%
\begin{equation*}
\gamma _{0}\left( P;x,y\right) \geq \gamma _{0}\left( Q;x,y\right) \qquad
\left( \geq 0\right)
\end{equation*}%
for any $P,Q$ with $P\geq Q\geq 0$ and $x,y\in H,$ where $\gamma \in \left\{
\sigma ,\delta ,\beta \right\} .$

The superadditivity and monotonicity properties of $\sigma _{0}$ and $\delta
_{0}$ have been noted by Dragomir and Mond in \cite{IHFDM}.

If $u\in \mathcal{P}\left( H\right) $ is such that $I\geq U\geq 0,$ where $I$
is the identity operator, then on using the superadditivity property of the
functionals $\sigma _{0},\delta _{0}$ and $\beta _{0}$ one may state the
following refinements for the Schwarz inequality:%
\begin{multline}
\left\Vert x\right\Vert \left\Vert y\right\Vert -\left\vert \left\langle
x,y\right\rangle \right\vert \geq \left\langle Ux,x\right\rangle ^{\frac{1}{2%
}}\left\langle Uy,y\right\rangle ^{\frac{1}{2}}-\left\vert \left\langle
Ux,y\right\rangle \right\vert  \label{1.3.29} \\
+\left\langle \left( I-U\right) x,x\right\rangle ^{\frac{1}{2}}\left\langle
\left( I-U\right) y,y\right\rangle ^{\frac{1}{2}}-\left\vert \left\langle
\left( I-U\right) x,y\right\rangle \right\vert \qquad \left( \geq 0\right) ,
\end{multline}%
\begin{multline}
\left\Vert x\right\Vert ^{2}\left\Vert y\right\Vert ^{2}-\left\vert
\left\langle x,y\right\rangle \right\vert ^{2}\geq \left\langle
Ux,x\right\rangle \left\langle Uy,y\right\rangle -\left\vert \left\langle
Ux,y\right\rangle \right\vert ^{2}  \label{1.3.30} \\
+\left\langle \left( I-U\right) x,x\right\rangle \left\langle \left(
I-U\right) y,y\right\rangle -\left\vert \left\langle \left( I-U\right)
x,y\right\rangle \right\vert ^{2}\qquad \left( \geq 0\right) ,
\end{multline}%
and%
\begin{multline}
\left( \left\Vert x\right\Vert ^{2}\left\Vert y\right\Vert ^{2}-\left\vert
\left\langle x,y\right\rangle \right\vert ^{2}\right) ^{\frac{1}{2}}\geq
\left( \left\langle Ux,x\right\rangle \left\langle Uy,y\right\rangle
-\left\vert \left\langle Ux,y\right\rangle \right\vert ^{2}\right) ^{\frac{1%
}{2}}  \label{1.3.31} \\
+\left( \left\langle \left( I-U\right) x,x\right\rangle \left\langle \left(
I-U\right) y,y\right\rangle -\left\vert \left\langle \left( I-U\right)
x,y\right\rangle \right\vert ^{2}\right) ^{\frac{1}{2}}\qquad \left( \geq
0\right)
\end{multline}%
for any $x,y\in H.$

Note that (\ref{1.3.31}) is a better result than (\ref{1.3.30}).

Finally, if we assume that $D\in \mathcal{P}\left( H\right) $ with $D\geq
\gamma I$, where $\gamma >0,$ i.e., $D$ is positive definite on $H,$ then we
may state the following inequalities%
\begin{equation}
\left\langle Dx,x\right\rangle ^{\frac{1}{2}}\left\langle Dy,y\right\rangle
^{\frac{1}{2}}-\left\vert \left\langle Dx,y\right\rangle \right\vert \geq
\gamma \left[ \left\Vert x\right\Vert \left\Vert y\right\Vert -\left\vert
\left\langle x,y\right\rangle \right\vert \right] \quad \left( \geq 0\right)
,  \label{1.3.32}
\end{equation}%
\begin{multline}
\qquad \left\langle Dx,x\right\rangle \left\langle Dy,y\right\rangle
-\left\vert \left\langle Dx,y\right\rangle \right\vert ^{2}  \label{1.3.33}
\\
\geq \gamma ^{2}\left[ \left\Vert x\right\Vert ^{2}\left\Vert y\right\Vert
^{2}-\left\vert \left\langle x,y\right\rangle \right\vert ^{2}\right] \quad
\left( \geq 0\right) ,\qquad
\end{multline}%
for any $x,y\in H$ and
\begin{multline}
\left\langle Dx,x\right\rangle \left\langle Dy,y\right\rangle -\left\vert
\left\langle Dx,y\right\rangle \right\vert ^{2}  \label{1.3.34} \\
\geq \gamma \max \left\{ \frac{\left\langle Dx,x\right\rangle }{\left\Vert
x\right\Vert ^{2}},\frac{\left\langle Dy,y\right\rangle }{\left\Vert
y\right\Vert ^{2}}\right\} \left[ \left\Vert x\right\Vert ^{2}\left\Vert
y\right\Vert ^{2}-\left\vert \left\langle x,y\right\rangle \right\vert ^{2}%
\right] \quad \left( \geq 0\right)
\end{multline}%
for any $x,y\in H\backslash \left\{ 0\right\} .$

The results (\ref{1.3.32}) and (\ref{1.3.33}) have been obtained by Dragomir
and Mond in \cite{IHFDM}.

Note that (\ref{1.3.34}) is a better result than (\ref{1.3.33}).

The above results (\ref{1.3.29}) -- (\ref{1.3.34}) also hold for $\func{Re}%
\left( \cdot \right) $ instead of $\left\vert \cdot \right\vert .$

\section{Applications for Sequences of Vectors}

\subsection{The Case of Mapping $\protect\sigma $}

Let $\mathcal{P}_{f}\left( \mathbb{N}\right) $ be the family of finite parts
of the natural number set $\mathbb{N}$, $\mathcal{S}_{+}\left( \mathbb{R}%
\right) $ the cone of nonnegative real sequences and for a given inner
product space $\left( H;\left\langle \cdot ,\cdot \right\rangle \right) $
over the real or complex number field $\mathbb{K}$, $\mathcal{S}\left(
H\right) $ the linear space of all sequences of vectors from $H,$ i.e.,%
\begin{equation*}
\mathcal{S}\left( H\right) :=\left\{ \mathbf{x}|\mathbf{x}=\left(
x_{i}\right) _{i\in \mathbb{N}},\ x_{i}\in H,\ i\in \mathbb{N}\right\} .
\end{equation*}

Consider $\left\langle \cdot ,\cdot \right\rangle _{\mathbf{p},I}:\mathcal{S}%
\left( H\right) \times \mathcal{S}\left( H\right) \rightarrow \mathbb{R}$
defined by%
\begin{equation*}
\left\langle \mathbf{x},\mathbf{y}\right\rangle _{\mathbf{p},I}:=\sum_{i\in
I}p_{i}\left\langle x_{i},y_{i}\right\rangle .
\end{equation*}

We may define the mapping $\sigma $ by%
\begin{equation}
\sigma \left( \mathbf{p},I,\mathbf{x},\mathbf{y}\right) :=\left( \sum_{i\in
I}p_{i}\left\Vert x_{i}\right\Vert ^{2}\sum_{i\in I}p_{i}\left\Vert
y_{i}\right\Vert ^{2}\right) ^{\frac{1}{2}}-\left\vert \sum_{i\in
I}p_{i}\left\langle x_{i},y_{i}\right\rangle \right\vert ,  \label{1.4.1}
\end{equation}%
where $\mathbf{p}\in \mathcal{S}_{+}\left( \mathbb{R}\right) ,$ $I\in
\mathcal{P}_{f}\left( \mathbb{N}\right) $ and $\mathbf{x},\mathbf{y}\in
\mathcal{S}\left( H\right) .$

We observe that, for a $I\in \mathcal{P}_{f}\left( \mathbb{N}\right)
\backslash \left\{ \varnothing \right\} ,$ the functional $\left\langle
\cdot ,\cdot \right\rangle _{\mathbf{p},I}\geq \left\langle \cdot ,\cdot
\right\rangle _{\mathbf{q},I},$ provided $\mathbf{p}\geq \mathbf{q}\geq
\mathbf{0}.$

Using Theorem \ref{t1.2.1}, we may state the following result.

\begin{proposition}
\label{p1.4.10}Let $I\in \mathcal{P}_{f}\left( \mathbb{N}\right) \backslash
\left\{ \varnothing \right\} $, $\mathbf{x},\mathbf{y}\in \mathcal{S}\left(
H\right) .$ Then the functional $\sigma \left( \cdot ,I,\mathbf{x},\mathbf{y}%
\right) $ is superadditive and monotonic nondecreasing on $\mathcal{S}%
_{+}\left( \mathbb{R}\right) .$
\end{proposition}

If $I,J\in \mathcal{P}_{f}\left( \mathbb{N}\right) \backslash \left\{
\varnothing \right\} ,$ with $I\cap J=\varnothing $, for a given $\mathbf{p}%
\in \mathcal{S}_{+}\left( \mathbb{R}\right) ,$ we observe that%
\begin{equation}
\left\langle \cdot ,\cdot \right\rangle _{\mathbf{p},I\cup J}=\left\langle
\cdot ,\cdot \right\rangle _{\mathbf{p},I}+\left\langle \cdot ,\cdot
\right\rangle _{\mathbf{p},J}.  \label{1.4.2}
\end{equation}%
Taking into account this property and on making use of Theorem \ref{t1.2.1},
we may state the following result.

\begin{proposition}
\label{p1.4.11}Let $\mathbf{p}\in \mathcal{S}_{+}\left( \mathbb{R}\right) $
and $\mathbf{x},\mathbf{y}\in \mathcal{S}\left( H\right) .$

\begin{enumerate}
\item[(i)] For any $I,J\in \mathcal{P}_{f}\left( \mathbb{N}\right)
\backslash \left\{ \varnothing \right\} ,$ with $I\cap J=\varnothing ,$ we
have%
\begin{equation}
\sigma \left( \mathbf{p},I\cup J,\mathbf{x},\mathbf{y}\right) \geq \sigma
\left( \mathbf{p},I,\mathbf{x},\mathbf{y}\right) +\sigma \left( \mathbf{p},J,%
\mathbf{x},\mathbf{y}\right) \qquad \left( \geq 0\right) ,  \label{1.4.3}
\end{equation}%
i.e., $\sigma \left( \mathbf{p},\cdot ,\mathbf{x},\mathbf{y}\right) $ is
superadditive as an index set mapping on $\mathcal{P}_{f}\left( \mathbb{N}%
\right) .$

\item[(ii)] If $\varnothing \neq J\subseteq I,$ $I,J\in \mathcal{P}%
_{f}\left( \mathbb{N}\right) ,$ then%
\begin{equation}
\sigma \left( \mathbf{p},I,\mathbf{x},\mathbf{y}\right) \geq \sigma \left(
\mathbf{p},J,\mathbf{x},\mathbf{y}\right) \qquad \left( \geq 0\right) ,
\label{1.4.4}
\end{equation}%
i.e., $\sigma \left( \mathbf{p},\cdot ,\mathbf{x},\mathbf{y}\right) $ is
monotonic nondecreasing as an index set mapping on $\mathcal{S}_{+}\left(
\mathbb{R}\right) .$
\end{enumerate}
\end{proposition}

It is well known that the following Cauchy-Bunyakovsky-Schwarz (CBS) type
inequality for sequences of vectors in an inner product space holds true:%
\begin{equation}
\sum_{i\in I}p_{i}\left\Vert x_{i}\right\Vert ^{2}\sum_{i\in
I}p_{i}\left\Vert y_{i}\right\Vert ^{2}\geq \left\vert \sum_{i\in
I}p_{i}\left\langle x_{i},y_{i}\right\rangle \right\vert ^{2}  \label{1.4.5}
\end{equation}%
for $I\in \mathcal{P}_{f}\left( \mathbb{N}\right) \backslash \left\{
\varnothing \right\} ,$ $\mathbf{p}\in \mathcal{S}_{+}\left( \mathbb{R}%
\right) $ and $\mathbf{x},\mathbf{y}\in \mathcal{S}\left( H\right) .$

If $p_{i}>0$ for all $i\in I,$ then equality holds in (\ref{1.4.5}) if and
only if there exists a scalar $\lambda \in \mathbb{K}$ such that $%
x_{i}=\lambda y_{i},$ $i\in I.$

Utilising the above results for the functional $\sigma $, we may state the
following inequalities related to the (CBS)-inequality (\ref{1.4.5}).

\begin{enumerate}
\item[(1)] Let $\alpha _{i}\in \mathbb{R}$, $x_{i},y_{i}\in H,$ $i\in
\left\{ 1,\dots ,n\right\} .$ Then one has the inequality:%
\begin{multline}
\sum_{i=1}^{n}\left\Vert x_{i}\right\Vert ^{2}\sum_{i=1}^{n}\left\Vert
y_{i}\right\Vert ^{2}-\left\vert \sum_{i=1}^{n}\left\langle
x_{i},y_{i}\right\rangle \right\vert  \label{1.4.6} \\
\geq \left( \sum_{i=1}^{n}\left\Vert x_{i}\right\Vert ^{2}\sin ^{2}\alpha
_{i}\sum_{i=1}^{n}\left\Vert y_{i}\right\Vert ^{2}\sin ^{2}\alpha
_{i}\right) ^{\frac{1}{2}}-\left\vert \sum_{i=1}^{n}\left\langle
x_{i},y_{i}\right\rangle \sin ^{2}\alpha _{i}\right\vert \\
+\left( \sum_{i=1}^{n}\left\Vert x_{i}\right\Vert ^{2}\cos ^{2}\alpha
_{i}\sum_{i=1}^{n}\left\Vert y_{i}\right\Vert ^{2}\cos ^{2}\alpha
_{i}\right) ^{\frac{1}{2}} \\
-\left\vert \sum_{i=1}^{n}\left\langle x_{i},y_{i}\right\rangle \cos
^{2}\alpha _{i}\right\vert \geq 0.
\end{multline}

\item[(2)] Denote $S_{n}\left( \mathbf{1}\right) :=\left\{ \mathbf{p}\in
\mathcal{S}_{+}\left( \mathbb{R}\right) |p_{i}\leq 1\ \text{for all }i\in
\left\{ 1,\dots ,n\right\} \right\} .$ Then for all $x_{i},y_{i}\in H,$ $%
i\in \left\{ 1,\dots ,n\right\} ,$ we have the bound:%
\begin{multline}
\left( \sum_{i=1}^{n}\left\Vert x_{i}\right\Vert
^{2}\sum_{i=1}^{n}\left\Vert y_{i}\right\Vert ^{2}\right) ^{\frac{1}{2}%
}-\left\vert \sum_{i=1}^{n}\left\langle x_{i},y_{i}\right\rangle \right\vert
\label{1.4.7} \\
=\sup_{\mathbf{p}\in S_{n}\left( \mathbf{1}\right) }\left[ \left(
\sum_{i=1}^{n}p_{i}\left\Vert x_{i}\right\Vert
^{2}\sum_{i=1}^{n}p_{i}\left\Vert y_{i}\right\Vert ^{2}\right) ^{\frac{1}{2}%
}-\left\vert \sum_{i=1}^{n}p_{i}\left\langle x_{i},y_{i}\right\rangle
\right\vert \right] \geq 0.
\end{multline}

\item[(3)] Let $p_{i}\geq 0,$ $x_{i},y_{i}\in H,$ $i\in \left\{ 1,\dots
,n\right\} .$ Then we have the inequality:%
\begin{multline}
\left( \sum_{i=1}^{2n}p_{i}\left\Vert x_{i}\right\Vert
^{2}\sum_{i=1}^{2n}p_{i}\left\Vert y_{i}\right\Vert ^{2}\right) ^{\frac{1}{2}%
}-\left\vert \sum_{i=1}^{2n}p_{i}\left\langle x_{i},y_{i}\right\rangle
\right\vert  \label{1.4.8} \\
\geq \left( \sum_{k=1}^{n}p_{2k}\left\Vert x_{2k}\right\Vert
^{2}\sum_{k=1}^{n}p_{2k}\left\Vert y_{2k}\right\Vert ^{2}\right) ^{\frac{1}{2%
}}-\left\vert \sum_{k=1}^{n}p_{2k}\left\langle x_{2k},y_{2k}\right\rangle
\right\vert \\
+\left( \sum_{k=1}^{n}p_{2k-1}\left\Vert x_{2k-1}\right\Vert
^{2}\sum_{k=1}^{n}p_{2k-1}\left\Vert y_{2k-1}\right\Vert ^{2}\right) ^{\frac{%
1}{2}} \\
-\left\vert \sum_{k=1}^{n}p_{2k-1}\left\langle
x_{2k-1},y_{2k-1}\right\rangle \right\vert \qquad \left( \geq 0\right) .
\end{multline}

\item[(4)] We have the bound:%
\begin{multline}
\left[ \sum_{i=1}^{n}p_{i}\left\Vert x_{i}\right\Vert
^{2}\sum_{i=1}^{n}p_{i}\left\Vert y_{i}\right\Vert ^{2}\right] ^{\frac{1}{2}%
}-\left\vert \sum_{i=1}^{n}p_{i}\left\langle x_{i},y_{i}\right\rangle
\right\vert  \label{1.4.9} \\
=\sup_{\varnothing \neq I\subseteq \left\{ 1,\dots ,n\right\} }\left( \left[
\sum_{i\in I}p_{i}\left\Vert x_{i}\right\Vert ^{2}\sum_{i\in
I}p_{i}\left\Vert y_{i}\right\Vert ^{2}\right] ^{\frac{1}{2}}-\left\vert
\sum_{i\in I}p_{i}\left\langle x_{i},y_{i}\right\rangle \right\vert \right)
\geq 0.
\end{multline}

\item[(5)] The sequence $S_{n}$ given by%
\begin{equation*}
S_{n}:=\left( \sum_{i=1}^{n}p_{i}\left\Vert x_{i}\right\Vert
^{2}\sum_{i=1}^{n}p_{i}\left\Vert y_{i}\right\Vert ^{2}\right) ^{\frac{1}{2}%
}-\left\vert \sum_{i=1}^{n}p_{i}\left\langle x_{i},y_{i}\right\rangle
\right\vert
\end{equation*}%
is nondecreasing, i.e.,%
\begin{equation}
S_{k+1}\geq S_{k},\quad k\geq 2  \label{1.4.10}
\end{equation}%
and we have the bound%
\begin{multline}
S_{n}\geq \max_{1\leq i<j\leq n}\left\{ \left( p_{i}\left\Vert
x_{i}\right\Vert ^{2}+p_{j}\left\Vert x_{j}\right\Vert ^{2}\right) ^{\frac{1%
}{2}}\left( p_{i}\left\Vert y_{i}\right\Vert ^{2}+p_{j}\left\Vert
y_{j}\right\Vert ^{2}\right) ^{\frac{1}{2}}\right.  \label{1.4.11} \\
-\left\vert p_{i}\left\langle x_{i},y_{i}\right\rangle +p_{j}\left\langle
x_{j},y_{j}\right\rangle \right\vert \big\}\geq 0,
\end{multline}%
for $n\geq 2$ and $x_{i},y_{i}\in H,$ $i\in \left\{ 1,\dots ,n\right\} .$
\end{enumerate}

\begin{remark}
\label{r1.4.5}The results in this subsection have been obtained by Dragomir
and Mond in \cite{IHFDM} for the particular case of scalar sequences $%
\mathbf{x} $ and $\mathbf{y}.$
\end{remark}

\subsection{The Case of Mapping $\protect\delta $}

Under the assumptions of the above subsection, we can define the following
functional%
\begin{equation*}
\delta \left( \mathbf{p},I,\mathbf{x},\mathbf{y}\right) :=\sum_{i\in
I}p_{i}\left\Vert x_{i}\right\Vert ^{2}\sum_{i\in I}p_{i}\left\Vert
y_{i}\right\Vert ^{2}-\left\vert \sum_{i\in I}p_{i}\left\langle
x_{i},y_{i}\right\rangle \right\vert ^{2},
\end{equation*}%
where $\mathbf{p}\in \mathcal{S}_{+}\left( \mathbb{R}\right) ,$ $I\in
\mathcal{P}_{f}\left( \mathbb{N}\right) \backslash \left\{ \varnothing
\right\} $ and $\mathbf{x},\mathbf{y}\in \mathcal{S}\left( H\right) .$

Utilising Theorem \ref{t1.2.2}, we may state the following results.

\begin{proposition}
\label{p1.4.12}We have

\begin{enumerate}
\item[(i)] For any $\mathbf{p},\mathbf{q}\in \mathcal{S}_{+}\left( \mathbb{R}%
\right) ,$ $I\in \mathcal{P}_{f}\left( \mathbb{N}\right) \backslash \left\{
\varnothing \right\} $ and $\mathbf{x},\mathbf{y}\in \mathcal{S}\left(
H\right) $ we have%
\begin{multline}
\delta \left( \mathbf{p}+\mathbf{q},I,\mathbf{x},\mathbf{y}\right) -\delta
\left( \mathbf{p},I,\mathbf{x},\mathbf{y}\right) -\delta \left( \mathbf{q},I,%
\mathbf{x},\mathbf{y}\right)  \label{1.4.12} \\
\geq \left( \det \left[
\begin{array}{ll}
\left( \sum\limits_{i\in I}p_{i}\left\Vert x_{i}\right\Vert ^{2}\right) ^{%
\frac{1}{2}} & \left( \sum\limits_{i\in I}p_{i}\left\Vert y_{i}\right\Vert
^{2}\right) ^{\frac{1}{2}} \\
\left( \sum\limits_{i\in I}q_{i}\left\Vert x_{i}\right\Vert ^{2}\right) ^{%
\frac{1}{2}} & \left( \sum\limits_{i\in I}q_{i}\left\Vert y_{i}\right\Vert
^{2}\right) ^{\frac{1}{2}}%
\end{array}%
\right] \right) ^{2}\geq 0.
\end{multline}

\item[(ii)] If $\mathbf{p}\geq \mathbf{q}\geq \mathbf{0},$ then%
\begin{multline}
\delta \left( \mathbf{p},I,\mathbf{x},\mathbf{y}\right) -\delta \left(
\mathbf{q},I,\mathbf{x},\mathbf{y}\right)  \label{1.4.13} \\
\geq \left( \det \left[
\begin{array}{cc}
\left( \sum\limits_{i\in I}p_{i}\left\Vert x_{i}\right\Vert ^{2}\right) ^{%
\frac{1}{2}} & \left( \sum\limits_{i\in I}p_{i}\left\Vert y_{i}\right\Vert
^{2}\right) ^{\frac{1}{2}} \\
\left( \sum\limits_{i\in I}\left( p_{i}-q_{i}\right) \left\Vert
x_{i}\right\Vert ^{2}\right) ^{\frac{1}{2}} & \left( \sum\limits_{i\in
I}\left( p_{i}-q_{i}\right) \left\Vert y_{i}\right\Vert ^{2}\right) ^{\frac{1%
}{2}}%
\end{array}%
\right] \right) ^{2}\geq 0.
\end{multline}
\end{enumerate}
\end{proposition}

\begin{proposition}
\label{p1.4.13}We have

\begin{enumerate}
\item[(i)] For any $I,J\in \mathcal{P}_{f}\left( \mathbb{N}\right) ,$ with $%
I\cap J=\varnothing $ and $\mathbf{p}\in \mathcal{S}_{+}\left( \mathbb{R}%
\right) ,\ \mathbf{x},\mathbf{y}\in \mathcal{S}\left( H\right) ,$ we have%
\begin{multline}
\delta \left( \mathbf{p},I\cup J,\mathbf{x},\mathbf{y}\right) -\delta \left(
\mathbf{p},I,\mathbf{x},\mathbf{y}\right) -\delta \left( \mathbf{p},J,%
\mathbf{x},\mathbf{y}\right)  \label{1.4.14} \\
\geq \left( \det \left[
\begin{array}{ll}
\left( \sum\limits_{i\in I}p_{i}\left\Vert x_{i}\right\Vert ^{2}\right) ^{%
\frac{1}{2}} & \left( \sum\limits_{i\in I}p_{i}\left\Vert y_{i}\right\Vert
^{2}\right) ^{\frac{1}{2}} \\
\left( \sum\limits_{i\in J}p_{i}\left\Vert x_{i}\right\Vert ^{2}\right) ^{%
\frac{1}{2}} & \left( \sum\limits_{i\in J}p_{i}\left\Vert y_{i}\right\Vert
^{2}\right) ^{\frac{1}{2}}%
\end{array}%
\right] \right) ^{2}\geq 0.
\end{multline}

\item[(ii)] If $\varnothing \neq J\subseteq I,$ $I\neq J,$ $I,J\in \mathcal{P%
}_{f}\left( \mathbb{N}\right) ,$ then we have%
\begin{multline}
\delta \left( \mathbf{p},I,\mathbf{x},\mathbf{y}\right) -\delta \left(
\mathbf{p},J,\mathbf{x},\mathbf{y}\right)  \label{1.4.15} \\
\geq \left( \det \left[
\begin{array}{ll}
\left( \sum\limits_{i\in I}p_{i}\left\Vert x_{i}\right\Vert ^{2}\right) ^{%
\frac{1}{2}} & \left( \sum\limits_{i\in I}p_{i}\left\Vert y_{i}\right\Vert
^{2}\right) ^{\frac{1}{2}} \\
\left( \sum\limits_{i\in I\backslash J}p_{i}\left\Vert x_{i}\right\Vert
^{2}\right) ^{\frac{1}{2}} & \left( \sum\limits_{i\in I\backslash
J}p_{i}\left\Vert y_{i}\right\Vert ^{2}\right) ^{\frac{1}{2}}%
\end{array}%
\right] \right) ^{2}\geq 0.
\end{multline}
\end{enumerate}
\end{proposition}

The following particular instances that provide refinements for the
(CBS)-inequality may be stated as well:%
\begin{align}
& \sum_{i\in I}\left\Vert x_{i}\right\Vert ^{2}\sum_{i\in I}\left\Vert
y_{i}\right\Vert ^{2}-\left\vert \sum_{i\in I}\left\langle
x_{i},y_{i}\right\rangle \right\vert ^{2}  \label{1.4.16} \\
& \geq \sum_{i\in I}\left\Vert x_{i}\right\Vert ^{2}\sin ^{2}\alpha
_{i}\sum_{i\in I}\left\Vert y_{i}\right\Vert ^{2}\sin ^{2}\alpha
_{i}-\left\vert \sum_{i\in I}\left\langle x_{i},y_{i}\right\rangle \sin
^{2}\alpha _{i}\right\vert ^{2}  \notag \\
& \qquad +\sum_{i\in I}\left\Vert x_{i}\right\Vert ^{2}\cos ^{2}\alpha
_{i}\sum_{i\in I}\left\Vert y_{i}\right\Vert ^{2}\cos ^{2}\alpha _{i}  \notag
\\
& \qquad -\left\vert \sum_{i\in I}\left\langle x_{i},y_{i}\right\rangle \cos
^{2}\alpha _{i}\right\vert ^{2}  \notag \\
& \geq \left( \det \left[
\begin{array}{ll}
\left( \sum\limits_{i\in I}\left\Vert x_{i}\right\Vert ^{2}\sin ^{2}\alpha
_{i}\right) ^{\frac{1}{2}} & \left( \sum\limits_{i\in I}\left\Vert
y_{i}\right\Vert ^{2}\sin ^{2}\alpha _{i}\right) ^{\frac{1}{2}} \\
\left( \sum\limits_{i\in I}\left\Vert x_{i}\right\Vert ^{2}\cos ^{2}\alpha
_{i}\right) ^{\frac{1}{2}} & \left( \sum\limits_{i\in I}\left\Vert
y_{i}\right\Vert ^{2}\cos ^{2}\alpha _{i}\right) ^{\frac{1}{2}}%
\end{array}%
\right] \right) ^{2}  \notag \\
& \geq 0,  \notag
\end{align}%
where $x_{i},y_{i}\in H,$ $\alpha _{i}\in \mathbb{R}$, $i\in I$ and $I\in
\mathcal{P}_{f}\left( \mathbb{N}\right) \backslash \left\{ \varnothing
\right\} .$

Suppose that $p_{i}\geq 0,$ $x_{i},y_{i}\in H,$ $i\in \left\{ 1,\dots
,2n\right\} .$ Then%
\begin{align}
& \sum_{i=1}^{2n}p_{i}\left\Vert x_{i}\right\Vert
^{2}\sum_{i=1}^{2n}p_{i}\left\Vert y_{i}\right\Vert ^{2}-\left\vert
\sum_{i=1}^{2n}p_{i}\left\langle x_{i},y_{i}\right\rangle \right\vert ^{2}
\label{1.4.17} \\
& \geq \sum_{k=1}^{n}p_{2k}\left\Vert x_{2k}\right\Vert
^{2}\sum_{k=1}^{n}p_{2k}\left\Vert y_{2k}\right\Vert ^{2}-\left\vert
\sum_{k=1}^{n}p_{2k}\left\langle x_{2k},y_{2k}\right\rangle \right\vert ^{2}
\notag \\
& \qquad +\sum_{k=1}^{n}p_{2k-1}\left\Vert x_{2k-1}\right\Vert
^{2}\sum_{k=1}^{n}p_{2k-1}\left\Vert y_{2k-1}\right\Vert ^{2}  \notag \\
& \qquad -\left\vert \sum_{k=1}^{n}p_{2k-1}\left\langle
x_{2k-1},y_{2k-1}\right\rangle \right\vert ^{2}  \notag \\
& \geq \hspace{-0.02in}\left( \hspace{-0.02in}\det \left[ \hspace{-0.02in}%
\begin{array}{ll}
\left( \sum\limits_{k=1}^{n}p_{2k}\left\Vert x_{2k}\right\Vert ^{2}\right) ^{%
\frac{1}{2}} & \hspace{-0.02in}\left( \sum\limits_{k=1}^{n}p_{2k}\left\Vert
y_{2k}\right\Vert ^{2}\right) ^{\frac{1}{2}} \\
\left( \sum\limits_{k=1}^{n}p_{2k-1}\left\Vert x_{2k-1}\right\Vert
^{2}\right) ^{\frac{1}{2}}\hspace{-0.02in} & \hspace{-0.02in}\left(
\sum\limits_{k=1}^{n}p_{2k-1}\left\Vert y_{2k-1}\right\Vert ^{2}\right) ^{%
\frac{1}{2}}%
\end{array}%
\right] \hspace{-0.02in}\right) ^{2}  \notag \\
& \geq 0.  \notag
\end{align}

\begin{remark}
\label{r1.4.6}The above results (\ref{1.4.12}) -- (\ref{1.4.17}) have been
obtained for the case where $\mathbf{x}$ and $\mathbf{y}$ are real or
complex numbers by Dragomir and Mond \cite{IHFDM}.
\end{remark}

Further, if we use Corollaries \ref{c1.2.2} and \ref{c1.2.1}, then we can
state the following propositions as well.

\begin{proposition}
\label{p1.4.14}We have

\begin{enumerate}
\item[(i)] For any $\mathbf{p},\mathbf{q}\in \mathcal{S}_{+}\left( \mathbb{R}%
\right) ,$ $I\in \mathcal{P}_{f}\left( \mathbb{N}\right) \backslash \left\{
\varnothing \right\} $ and $\mathbf{x},\mathbf{y}\in \mathcal{S}\left(
H\right) \backslash \left\{ 0\right\} $ we have%
\begin{multline}
\delta \left( \mathbf{p}+\mathbf{q},I,\mathbf{x},\mathbf{y}\right) -\delta
\left( \mathbf{p},I,\mathbf{x},\mathbf{y}\right) -\delta \left( \mathbf{q},I,%
\mathbf{x},\mathbf{y}\right)   \label{1.4.18} \\
\geq \max \left\{ \frac{\sum_{i\in I}p_{i}\left\Vert x_{i}\right\Vert ^{2}}{%
\sum_{i\in I}q_{i}\left\Vert x_{i}\right\Vert ^{2}}\delta \left( \mathbf{q}%
,I,\mathbf{x},\mathbf{y}\right) +\frac{\sum_{i\in I}q_{i}\left\Vert
x_{i}\right\Vert ^{2}}{\sum_{i\in I}p_{i}\left\Vert x_{i}\right\Vert ^{2}}%
\delta \left( \mathbf{p},I,\mathbf{x},\mathbf{y}\right) ,\right.  \\
\left. \frac{\sum_{i\in I}p_{i}\left\Vert y_{i}\right\Vert ^{2}}{\sum_{i\in
I}q_{i}\left\Vert y_{i}\right\Vert ^{2}}\delta \left( \mathbf{q},I,\mathbf{x}%
,\mathbf{y}\right) +\frac{\sum_{i\in I}q_{i}\left\Vert y_{i}\right\Vert ^{2}%
}{\sum_{i\in I}p_{i}\left\Vert y_{i}\right\Vert ^{2}}\delta \left( \mathbf{p}%
,I,\mathbf{x},\mathbf{y}\right) \right\} \geq 0.
\end{multline}

\item[(ii)] If $\mathbf{p}\geq \mathbf{q}\geq \mathbf{0}$ and $I\in \mathcal{%
P}_{f}\left( \mathbb{N}\right) \backslash \left\{ \varnothing \right\} $, $%
\mathbf{x},\mathbf{y}\in \mathcal{S}\left( H\right) \backslash \left\{
0\right\} ,$ then:%
\begin{multline}
\delta \left( \mathbf{p},I,\mathbf{x},\mathbf{y}\right) -\delta \left(
\mathbf{q},I,\mathbf{x},\mathbf{y}\right)   \label{1.4.19} \\
\geq \max \left\{ \frac{\sum_{i\in I}\left( p_{i}-q_{i}\right) \left\Vert
x_{i}\right\Vert ^{2}}{\sum_{i\in I}p_{i}\left\Vert x_{i}\right\Vert ^{2}},%
\frac{\sum_{i\in I}\left( p_{i}-q_{i}\right) \left\Vert y_{i}\right\Vert ^{2}%
}{\sum_{i\in I}p_{i}\left\Vert y_{i}\right\Vert ^{2}}\right\} \delta \left(
\mathbf{p},I,\mathbf{x},\mathbf{y}\right) \geq 0.
\end{multline}
\end{enumerate}
\end{proposition}

\begin{proposition}
\label{p1.4.15}We have

\begin{enumerate}
\item[(i)] For any $I,J\in \mathcal{P}_{f}\left( \mathbb{N}\right)
\backslash \left\{ \varnothing \right\} ,$ with $I\cap J=\varnothing $ and $%
\mathbf{p}\in \mathcal{S}_{+}\left( \mathbb{R}\right) ,\ \mathbf{x},\mathbf{y%
}\in \mathcal{S}\left( H\right) \backslash \left\{ 0\right\} ,$ we have%
\begin{multline}
\delta \left( \mathbf{p},I\cup J,\mathbf{x},\mathbf{y}\right) -\delta \left(
\mathbf{p},I,\mathbf{x},\mathbf{y}\right) -\delta \left( \mathbf{p},J,%
\mathbf{x},\mathbf{y}\right)   \label{1.4.20} \\
\geq \max \left\{ \frac{\sum_{i\in I}p_{i}\left\Vert x_{i}\right\Vert ^{2}}{%
\sum_{j\in J}p_{j}\left\Vert x_{j}\right\Vert ^{2}}\delta \left( \mathbf{p}%
,J,\mathbf{x},\mathbf{y}\right) +\frac{\sum_{j\in J}p_{j}\left\Vert
x_{j}\right\Vert ^{2}}{\sum_{i\in I}p_{i}\left\Vert x_{i}\right\Vert ^{2}}%
\delta \left( \mathbf{p},I,\mathbf{x},\mathbf{y}\right) ,\right.  \\
\left. \frac{\sum_{i\in I}p_{i}\left\Vert y_{i}\right\Vert ^{2}}{\sum_{j\in
J}p_{j}\left\Vert y_{j}\right\Vert ^{2}}\delta \left( \mathbf{p},J,\mathbf{x}%
,\mathbf{y}\right) +\frac{\sum_{j\in J}p_{j}\left\Vert y_{j}\right\Vert ^{2}%
}{\sum_{i\in I}p_{i}\left\Vert y_{i}\right\Vert ^{2}}\delta \left( \mathbf{p}%
,I,\mathbf{x},\mathbf{y}\right) \right\} \geq 0.
\end{multline}

\item[(ii)] If $\varnothing \neq J\subseteq I,$ $I\neq J,$ $I,J\in \mathcal{P%
}_{f}\left( \mathbb{N}\right) \backslash \left\{ \varnothing \right\} $ and $%
\mathbf{p}\in \mathcal{S}_{+}\left( \mathbb{R}\right) \backslash \left\{
0\right\} ,$ $\mathbf{x},\mathbf{y}\in \mathcal{S}\left( H\right) \backslash
\left\{ 0\right\} ,$ then%
\begin{multline}
\delta \left( \mathbf{p},I,\mathbf{x},\mathbf{y}\right) -\delta \left(
\mathbf{p},J,\mathbf{x},\mathbf{y}\right)   \label{1.4.21} \\
\geq \max \left\{ \frac{\sum_{k\in I\backslash J}p_{k}\left\Vert
x_{k}\right\Vert ^{2}}{\sum_{i\in I}p_{i}\left\Vert x_{i}\right\Vert ^{2}},%
\frac{\sum_{k\in I\backslash J}p_{k}\left\Vert y_{k}\right\Vert ^{2}}{%
\sum_{i\in I}p_{i}\left\Vert y_{i}\right\Vert ^{2}}\right\} \delta \left(
\mathbf{p},J,\mathbf{x},\mathbf{y}\right) \geq 0.
\end{multline}
\end{enumerate}
\end{proposition}

\begin{remark}
\label{r1.4.7}The results in Proposition \ref{p1.4.14} have been obtained by
Dragomir and Mond in \cite{IHFDM2} for the case of scalar sequences $\mathbf{%
x}$ and $\mathbf{y}.$
\end{remark}

\subsection{The Case of Mapping $\protect\beta $}

With the assumptions in the first subsections, we can define the following
functional%
\begin{align*}
\beta \left( \mathbf{p},I,\mathbf{x},\mathbf{y}\right) & :=\left[ \delta
\left( \mathbf{p},I,\mathbf{x},\mathbf{y}\right) \right] ^{\frac{1}{2}} \\
& =\left[ \sum_{i\in I}p_{i}\left\Vert x_{i}\right\Vert ^{2}\sum_{i\in
I}p_{i}\left\Vert y_{i}\right\Vert ^{2}-\left\vert \sum_{i\in
I}p_{i}\left\langle x_{i},y_{i}\right\rangle \right\vert ^{2}\right] ^{\frac{%
1}{2}},
\end{align*}%
where $\mathbf{p}\in \mathcal{S}_{+}\left( \mathbb{R}\right) ,$ $I\in
\mathcal{P}_{f}\left( \mathbb{N}\right) \backslash \left\{ \varnothing
\right\} $ and $\mathbf{x},\mathbf{y}\in \mathcal{S}\left( H\right) .$

Utilising Theorem \ref{t1.2.4}, we can state the following results:

\begin{proposition}
\label{p1.4.16}We have

\begin{enumerate}
\item[(i)] The functional $\beta \left( \cdot ,I,\mathbf{x},\mathbf{y}%
\right) $ is superadditive on $\mathcal{S}_{+}\left( \mathbb{R}\right) $ for
any $I\in \mathcal{P}_{f}\left( \mathbb{N}\right) \backslash \left\{
\varnothing \right\} $ and $\mathbf{x},\mathbf{y}\in \mathcal{S}\left(
H\right) .$

\item[(ii)] The functional $\beta \left( \mathbf{p},\cdot ,\mathbf{x},%
\mathbf{y}\right) $ is superadditive as an index set mapping on $\mathcal{P}%
_{f}\left( \mathbb{N}\right) $ and $\mathbf{x},\mathbf{y}\in \mathcal{S}%
\left( H\right) .$
\end{enumerate}
\end{proposition}

As simple consequences of the above proposition, we may state the following
refinements of the (CBS)-inequality.

\begin{enumerate}
\item[(a)] If $\mathbf{x},\mathbf{y}\in \mathcal{S}\left( H\right) $ and $%
\alpha _{i}\in \mathbb{R}$, $i\in I$ with $I\in \mathcal{P}_{f}\left(
\mathbb{N}\right) \backslash \left\{ \varnothing \right\} ,$ then%
\begin{multline}
\left( \sum_{i\in I}\left\Vert x_{i}\right\Vert ^{2}\sum_{i\in I}\left\Vert
y_{i}\right\Vert ^{2}-\left\vert \sum_{i\in I}\left\langle
x_{i},y_{i}\right\rangle \right\vert ^{2}\right) ^{\frac{1}{2}}
\label{1.4.22} \\
\geq \left( \sum_{i\in I}\left\Vert x_{i}\right\Vert ^{2}\sin ^{2}\alpha
_{i}\sum_{i\in I}\left\Vert y_{i}\right\Vert ^{2}\sin ^{2}\alpha
_{i}-\left\vert \sum_{i\in I}\left\langle x_{i},y_{i}\right\rangle \sin
^{2}\alpha _{i}\right\vert ^{2}\right) ^{\frac{1}{2}} \\
+\left( \sum_{i\in I}\left\Vert x_{i}\right\Vert ^{2}\cos ^{2}\alpha
_{i}\sum_{i\in I}\left\Vert y_{i}\right\Vert ^{2}\cos ^{2}\alpha
_{i}-\left\vert \sum_{i\in I}\left\langle x_{i},y_{i}\right\rangle \cos
^{2}\alpha _{i}\right\vert ^{2}\right) ^{\frac{1}{2}}\geq 0.
\end{multline}

\item[(b)] If $x_{i},y_{i}\in H,$ $p_{i}>0,$ $i\in \left\{ 1,\dots
,2n\right\} ,$ then%
\begin{multline}
\left( \sum_{i=1}^{2n}p_{i}\left\Vert x_{i}\right\Vert
^{2}\sum_{i=1}^{2n}p_{i}\left\Vert y_{i}\right\Vert ^{2}-\left\vert
\sum_{i=1}^{2n}p_{i}\left\langle x_{i},y_{i}\right\rangle \right\vert
^{2}\right) ^{\frac{1}{2}}  \label{1.4.23} \\
\geq \left( \sum_{k=1}^{n}p_{2k}\left\Vert x_{2k}\right\Vert
^{2}\sum_{k=1}^{n}p_{2k}\left\Vert y_{2k}\right\Vert ^{2}-\left\vert
\sum_{k=1}^{n}p_{2k}\left\langle x_{2k},y_{2k}\right\rangle \right\vert
^{2}\right) ^{\frac{1}{2}} \\
+\left( \sum_{k=1}^{n}p_{2k-1}\left\Vert x_{2k-1}\right\Vert
^{2}\sum_{k=1}^{n}p_{2k-1}\left\Vert y_{2k-1}\right\Vert ^{2}\right.  \\
-\left. \left\vert \sum_{k=1}^{n}p_{2k-1}\left\langle
x_{2k-1},y_{2k-1}\right\rangle \right\vert ^{2}\right) ^{\frac{1}{2}}\qquad
\left( \geq 0\right) .
\end{multline}
\end{enumerate}

\begin{remark}
\label{r1.4.8}Part (i) of Proposition \ref{p1.4.16} and the inequality (\ref%
{1.4.21}) have been obtained by Dragomir and Mond in \cite{IHFDM2} for the
case of scalar sequences $\mathbf{x}$ and $\mathbf{y}.$
\end{remark}

%


\chapter{Schwarz Related Inequalities}

\label{ch2}

\section{Introduction}

Let $H$ be a linear space over the real or complex number field $\mathbb{K}$%
. The functional $\left\langle \cdot ,\cdot \right\rangle :H\times
H\rightarrow \mathbb{K}$ is called an \textit{inner product }on $H$ if it
satisfies the conditions

\begin{enumerate}
\item[(i)] $\left\langle x,x\right\rangle \geq 0$ for any $x\in H$ and $%
\left\langle x,x\right\rangle =0$ iff $x=0;$

\item[(ii)] $\left\langle \alpha x+\beta y,z\right\rangle =\alpha
\left\langle x,z\right\rangle +\beta \left\langle y,z\right\rangle $ for any 
$\alpha ,\beta \in \mathbb{K}$ and $x,y,z\in H;$

\item[(iii)] $\left\langle y,x\right\rangle =\overline{\left\langle
x,y\right\rangle }$ for any $x,y\in H.$
\end{enumerate}

A first fundamental consequence of the properties (i)-(iii) above, is the 
\textit{Schwarz inequality:}%
\begin{equation}
\left\vert \left\langle x,y\right\rangle \right\vert ^{2}\leq \left\langle
x,x\right\rangle \left\langle y,y\right\rangle ,  \label{CBS}
\end{equation}%
for any $x,y\in H.$ The equality holds in (\ref{CBS}) if and only if the
vectors $x$ and $y$ are \textit{linearly dependent,} i.e., there exists a
nonzero constant $\alpha \in \mathbb{K}$ so that $x=\alpha y.$

If we denote $\left\Vert x\right\Vert :=\sqrt{\left\langle x,x\right\rangle }%
,x\in H,$ then one may state the following properties

\begin{enumerate}
\item[(n)] $\left\Vert x\right\Vert \geq 0$ for any $x\in H$ and $\left\Vert
x\right\Vert =0$ iff $x=0;$

\item[(nn)] $\left\Vert \alpha x\right\Vert =\left\vert \alpha \right\vert
\left\Vert x\right\Vert $ for any $\alpha \in \mathbb{K}$ and $x\in H;$

\item[(nnn)] $\left\Vert x+y\right\Vert \leq \left\Vert x\right\Vert
+\left\Vert y\right\Vert $ for any $x,y\in H$ (the triangle inequality);
\end{enumerate}

i.e., $\left\Vert \cdot \right\Vert $ is a \textit{norm} on $H.$

In this chapter we present some classical and recent refinements and reverse
inequalities for the Schwarz and the triangle inequalities. More precisely,
we point out upper bounds or positive lower bounds for the nonnegative
quantities%
\begin{equation*}
\left\Vert x\right\Vert \left\Vert y\right\Vert -\left\vert \left\langle
x,y\right\rangle \right\vert ,\text{ \ }\left\Vert x\right\Vert
^{2}\left\Vert y\right\Vert ^{2}-\left\vert \left\langle x,y\right\rangle
\right\vert ^{2}
\end{equation*}%
and 
\begin{equation*}
\left\Vert x\right\Vert +\left\Vert y\right\Vert -\left\Vert x+y\right\Vert
\end{equation*}%
under various assumptions for the vectors $x,y\in H.$

If the vectors $x,y\in H$ are not \textit{orthogonal}, i.e., $\left\langle
x,y\right\rangle \neq 0,$ then some upper and lower bounds for the
supra-unitary quantities%
\begin{equation*}
\frac{\left\Vert x\right\Vert \left\Vert y\right\Vert }{\left\vert
\left\langle x,y\right\rangle \right\vert },\ \ \frac{\left\Vert
x\right\Vert ^{2}\left\Vert y\right\Vert ^{2}}{\left\vert \left\langle
x,y\right\rangle \right\vert ^{2}}
\end{equation*}%
under appropriate restrictions for the vectors $x$ and $y$ are provided as
well.

The inequalities obtained by Buzano, Richards, Precupanu and Moore and their
extensions and generalizations for orthonormal families of vectors in both
real and complex inner product spaces are presented. Recent results
concerning the classical refinement of Schwarz inequality due to Kurepa for
the complexification of real inner product spaces are also reviewed. Various
applications for integral inequalities including a version of Heisenberg
inequality for vector valued functions in Hilbert spaces are provided as
well.

\section{Inequalities Related to Schwarz's One}

\subsection{Some Refinements}

The following result holds \cite[Theorem 1]{SSD3} (see also \cite[Theorem 2]%
{DS}).

\begin{theorem}[Dragomir, 1985]
\label{t2.1.1}Let $\left( H,\left\langle \cdot ,\cdot \right\rangle \right) $
be a real or complex inner product space. Then%
\begin{multline}
\quad \left( \left\Vert x\right\Vert ^{2}\left\Vert y\right\Vert
^{2}-\left\vert \left\langle x,y\right\rangle \right\vert ^{2}\right) \left(
\left\Vert y\right\Vert ^{2}\left\Vert z\right\Vert ^{2}-\left\vert
\left\langle y,z\right\rangle \right\vert ^{2}\right)  \label{2.1.1} \\
\geq \left\vert \left\langle x,z\right\rangle \left\Vert y\right\Vert
^{2}-\left\langle x,y\right\rangle \left\langle y,z\right\rangle \right\vert
^{2}\quad
\end{multline}%
for any $x,y,z\in H.$
\end{theorem}

\begin{proof}
We follow the proof in \cite{SSD3}.

Let us consider the mapping%
\begin{equation*}
p_{y}:H\times H\rightarrow \mathbb{K},\qquad p_{y}\left( x,z\right)
=\left\langle x,z\right\rangle \left\Vert y\right\Vert ^{2}-\left\langle
x,y\right\rangle \left\langle y,z\right\rangle
\end{equation*}%
for each $y\in H\backslash \left\{ 0\right\} .$

It is easy to see that $p_{y}\left( \cdot ,\cdot \right) $ is a nonnegative
Hermitian form and then on writing Schwarz's inequality%
\begin{equation*}
\left\vert p_{y}\left( x,z\right) \right\vert ^{2}\leq p_{y}\left(
x,x\right) p_{y}\left( z,z\right) ,\qquad x,z\in H
\end{equation*}%
we obtain the desired inequality (\ref{2.1.1}).
\end{proof}

\begin{remark}
\label{r2.1.1}From (\ref{2.1.1}) it follows that \cite[Corollary 1]{SSD3}
(see also \cite[Corollary 2.1]{DS})%
\begin{multline}
\left( \left\Vert x+z\right\Vert ^{2}\left\Vert y\right\Vert ^{2}-\left\vert
\left\langle x+z,y\right\rangle \right\vert ^{2}\right) ^{\frac{1}{2}}
\label{2.1.2} \\
\leq \left( \left\Vert x\right\Vert ^{2}\left\Vert y\right\Vert
^{2}-\left\vert \left\langle x,y\right\rangle \right\vert ^{2}\right) ^{%
\frac{1}{2}}+\left( \left\Vert y\right\Vert ^{2}\left\Vert z\right\Vert
^{2}-\left\vert \left\langle y,z\right\rangle \right\vert ^{2}\right) ^{%
\frac{1}{2}}
\end{multline}%
for every $x,y,z\in H.$

Putting $z=\lambda y$ in (\ref{2.1.2}), we get:%
\begin{align}
0& \leq \left\Vert x+\lambda y\right\Vert ^{2}\left\Vert y\right\Vert
^{2}-\left\vert \left\langle x+\lambda y,y\right\rangle \right\vert ^{2}
\label{2.1.3} \\
& \leq \left\Vert x\right\Vert ^{2}\left\Vert y\right\Vert ^{2}-\left\vert
\left\langle x,y\right\rangle \right\vert ^{2}  \notag
\end{align}%
and, in particular,%
\begin{equation}
0\leq \left\Vert x\pm y\right\Vert ^{2}\left\Vert y\right\Vert
^{2}-\left\vert \left\langle x\pm y,y\right\rangle \right\vert ^{2}\leq
\left\Vert x\right\Vert ^{2}\left\Vert y\right\Vert ^{2}-\left\vert
\left\langle x,y\right\rangle \right\vert ^{2}  \label{2.1.4}
\end{equation}%
for every $x,y\in H.$
\end{remark}

Both inequalities (\ref{2.1.3}) and (\ref{2.1.4}) have been obtained in \cite%
{SSD3}.

We note here that the inequality (\ref{2.1.3}) is in fact equivalent to the
following statement%
\begin{equation}
\sup_{\lambda \in \mathbb{K}}\left[ \left\Vert x+\lambda y\right\Vert
^{2}\left\Vert y\right\Vert ^{2}-\left\vert \left\langle x+\lambda
y,y\right\rangle \right\vert ^{2}\right] =\left\Vert x\right\Vert
^{2}\left\Vert y\right\Vert ^{2}-\left\vert \left\langle x,y\right\rangle
\right\vert ^{2}  \label{2.1.5}
\end{equation}%
for each $x,y\in H.$

The following corollary may be stated \cite[Corollary 2]{SSD3} (see also 
\cite[Corollary 2.2]{DS}):

\begin{corollary}[Dragomir, 1985]
\label{c2.1.1}For any $x,y,z\in H\backslash \left\{ 0\right\} $ we have the
inequality%
\begin{equation}
\left\vert \frac{\left\langle x,y\right\rangle }{\left\Vert x\right\Vert
\left\Vert y\right\Vert }\right\vert ^{2}+\left\vert \frac{\left\langle
y,z\right\rangle }{\left\Vert y\right\Vert \left\Vert z\right\Vert }%
\right\vert ^{2}+\left\vert \frac{\left\langle z,x\right\rangle }{\left\Vert
z\right\Vert \left\Vert x\right\Vert }\right\vert ^{2}\leq 1+2\left\vert 
\frac{\left\langle x,y\right\rangle \left\langle y,z\right\rangle
\left\langle z,x\right\rangle }{\left\Vert x\right\Vert ^{2}\left\Vert
y\right\Vert ^{2}\left\Vert z\right\Vert ^{2}}\right\vert .  \label{2.1.6}
\end{equation}
\end{corollary}

\begin{proof}
By the modulus properties we obviously have%
\begin{equation*}
\left\vert \left\langle x,z\right\rangle \left\Vert y\right\Vert
^{2}-\left\langle x,y\right\rangle \left\langle y,z\right\rangle \right\vert
\geq \left\vert \left\vert \left\langle x,z\right\rangle \right\vert
\left\Vert y\right\Vert ^{2}-\left\vert \left\langle x,y\right\rangle
\right\vert \left\vert \left\langle y,z\right\rangle \right\vert \right\vert
.
\end{equation*}%
Therefore, by (\ref{2.1.1}) we may state that%
\begin{multline*}
\left( \left\Vert x\right\Vert ^{2}\left\Vert y\right\Vert ^{2}-\left\vert
\left\langle x,y\right\rangle \right\vert ^{2}\right) \left( \left\Vert
y\right\Vert ^{2}\left\Vert z\right\Vert ^{2}-\left\vert \left\langle
y,z\right\rangle \right\vert ^{2}\right) \\
\geq \left\vert \left\langle x,z\right\rangle \right\vert ^{2}\left\Vert
y\right\Vert ^{4}-2\left\vert \left\langle x,y\right\rangle \left\langle
y,z\right\rangle \left\langle z,x\right\rangle \right\vert \left\Vert
y\right\Vert ^{2}+\left\vert \left\langle x,y\right\rangle \right\vert
^{2}\left\vert \left\langle y,z\right\rangle \right\vert ^{2},
\end{multline*}%
which, upon elementary calculation, is equivalent to (\ref{2.1.6}).
\end{proof}

\begin{remark}
\label{r2.1.2}If we utilise the elementary inequality $a^{2}+b^{2}+c^{2}\geq
3abc$ when $a,b,c\geq 0,$ then one can state the following inequality%
\begin{equation}
3\left\vert \frac{\left\langle x,y\right\rangle \left\langle
y,z\right\rangle \left\langle z,x\right\rangle }{\left\Vert x\right\Vert
^{2}\left\Vert y\right\Vert ^{2}\left\Vert z\right\Vert ^{2}}\right\vert
\leq \left\vert \frac{\left\langle x,y\right\rangle }{\left\Vert
x\right\Vert \left\Vert y\right\Vert }\right\vert ^{2}+\left\vert \frac{%
\left\langle y,z\right\rangle }{\left\Vert y\right\Vert \left\Vert
z\right\Vert }\right\vert ^{2}+\left\vert \frac{\left\langle
z,x\right\rangle }{\left\Vert z\right\Vert \left\Vert x\right\Vert }%
\right\vert ^{2}  \label{2.1.7}
\end{equation}%
for any $x,y,z\in H\backslash \left\{ 0\right\} .$ Therefore, the inequality
(\ref{2.1.6}) may be regarded as a reverse inequality of (\ref{2.1.7}).
\end{remark}

The following refinement of the Schwarz inequality holds \cite[Theorem 2]%
{SSD3} (see also \cite[Corollary 1.1]{DS}):

\begin{theorem}[Dragomir, 1985]
\label{t2.1.2}For any $x,y\in H$ and $e\in H$ with $\left\Vert e\right\Vert
=1,$ the following refinement of the Schwarz inequality holds:%
\begin{equation}
\left\Vert x\right\Vert \left\Vert y\right\Vert \geq \left\vert \left\langle
x,y\right\rangle -\left\langle x,e\right\rangle \left\langle
e,y\right\rangle \right\vert +\left\vert \left\langle x,e\right\rangle
\left\langle e,y\right\rangle \right\vert \geq \left\vert \left\langle
x,y\right\rangle \right\vert .  \label{2.1.8}
\end{equation}
\end{theorem}

\begin{proof}
We follow the proof in \cite{SSD3}.

Applying the inequality (\ref{2.1.1}), we can state that%
\begin{equation}
\left( \left\Vert x\right\Vert ^{2}-\left\vert \left\langle x,e\right\rangle
\right\vert ^{2}\right) \left( \left\Vert y\right\Vert ^{2}-\left\vert
\left\langle y,e\right\rangle \right\vert ^{2}\right) \geq \left\vert
\left\langle x,y\right\rangle -\left\langle x,e\right\rangle \left\langle
e,y\right\rangle \right\vert ^{2}.  \label{2.1.9}
\end{equation}%
Utilising the elementary inequality for real numbers%
\begin{equation}
\left( m^{2}-n^{2}\right) \left( p^{2}-q^{2}\right) \leq \left( mp-nq\right)
^{2},  \label{2.1.10}
\end{equation}%
we can easily see that%
\begin{multline}
\quad \left( \left\Vert x\right\Vert \left\Vert y\right\Vert -\left\vert
\left\langle x,e\right\rangle \left\langle e,y\right\rangle \right\vert
\right) ^{2}  \label{2.1.11} \\
\geq \left( \left\Vert x\right\Vert ^{2}-\left\vert \left\langle
x,e\right\rangle \right\vert ^{2}\right) \left( \left\Vert y\right\Vert
^{2}-\left\vert \left\langle y,e\right\rangle \right\vert ^{2}\right) \quad
\end{multline}%
for any $x,y,e\in H$ with $\left\Vert e\right\Vert =1.$

Since, by Schwarz's inequality%
\begin{equation}
\left\Vert x\right\Vert \left\Vert y\right\Vert \geq \left\vert \left\langle
x,e\right\rangle \left\langle e,y\right\rangle \right\vert  \label{2.1.12}
\end{equation}%
hence, by (\ref{2.1.9}) and (\ref{2.1.11}) we deduce the first part of (\ref%
{2.1.9}).

The second part of (\ref{2.1.9}) is obvious.
\end{proof}

\begin{corollary}[Dragomir, 1985]
\label{c2.1.2}If $x,y,e\in H$ are such that $\left\Vert e\right\Vert =1$ and 
$x\perp y,$ then%
\begin{equation}
\left\Vert x\right\Vert \left\Vert y\right\Vert \geq 2\left\vert
\left\langle x,e\right\rangle \left\langle e,y\right\rangle \right\vert .
\label{2.1.13}
\end{equation}
\end{corollary}

\begin{remark}
\label{r2.1.3}Assume that $A:H\rightarrow H$ is a bounded linear operator on 
$H.$ For $x,e\in H$ with $\left\Vert x\right\Vert =\left\Vert e\right\Vert
=1,$ we have by (\ref{2.1.8}) that%
\begin{equation}
\left\Vert Ay\right\Vert \geq \left\vert \left\langle x,Ay\right\rangle
-\left\langle x,e\right\rangle \left\langle e,Ay\right\rangle \right\vert
+\left\vert \left\langle x,e\right\rangle \left\langle e,Ay\right\rangle
\right\vert \geq \left\vert \left\langle x,Ay\right\rangle \right\vert
\label{2.1.14}
\end{equation}%
for any $y\in H.$

Taking the supremum over $x\in H,$ $\left\Vert x\right\Vert =1$ in (\ref%
{2.1.14}) and noting that $\left\Vert Ay\right\Vert =\sup\limits_{\left\Vert
x\right\Vert =1}\left\vert \left\langle x,Ay\right\rangle \right\vert ,$ we
deduce the representation%
\begin{equation}
\left\Vert Ay\right\Vert =\sup_{\left\Vert x\right\Vert =1}\left\{
\left\vert \left\langle x,Ay\right\rangle -\left\langle x,e\right\rangle
\left\langle e,Ay\right\rangle \right\vert +\left\vert \left\langle
x,e\right\rangle \left\langle e,Ay\right\rangle \right\vert \right\}
\label{2.1.15}
\end{equation}%
for any $y\in H.$ Finally, on taking the supremum over $y\in H,$ $\left\Vert
y\right\Vert =1$ in (\ref{2.1.15}) we get%
\begin{equation}
\left\Vert A\right\Vert =\sup_{\left\Vert y\right\Vert =1,\left\Vert
x\right\Vert =1}\left\{ \left\vert \left\langle x,Ay\right\rangle
-\left\langle x,e\right\rangle \left\langle e,Ay\right\rangle \right\vert
+\left\vert \left\langle x,e\right\rangle \left\langle e,Ay\right\rangle
\right\vert \right\}  \label{2.1.16}
\end{equation}%
for any $e\in H,$ $\left\Vert e\right\Vert =1,$ a representation that has
been obtained in \cite[Eq. 9]{SSD3}.
\end{remark}

\begin{remark}
\label{r2.1.4}Let $\left( H;\left\langle \cdot ,\cdot \right\rangle \right) $
be a Hilbert space. Then for any continuous linear functional $%
f:H\rightarrow \mathbb{K}$, $f\neq 0,$ there exists, by the Riesz
representation theorem a unique vector $e\in H\backslash \left\{ 0\right\} $
such that $f\left( x\right) =\left\langle x,e\right\rangle $ for $x\in H$
and $\left\Vert f\right\Vert =\left\Vert e\right\Vert .$

If $E$ is a nonzero linear subspace of $H$ and if we denote by $E^{\perp }$
its \ orthogonal complement, i.e., we recall that $E^{\perp }:=\left\{ y\in
H|y\perp x\right\} $ then for any $x\in E$ and $y\in E^{\perp },$ by (\ref%
{2.1.13}) we may state that%
\begin{equation*}
\left\Vert x\right\Vert \left\Vert y\right\Vert \geq 2\left\vert
\left\langle x,\frac{e}{\left\Vert x\right\Vert }\right\rangle \left\langle
y,\frac{e}{\left\Vert y\right\Vert }\right\rangle \right\vert ,
\end{equation*}%
giving, for $x,y\neq 0,$ that%
\begin{equation}
\left\Vert f\right\Vert ^{2}\geq 2\left\vert \left\langle x,e\right\rangle
\left\langle y,e\right\rangle \right\vert =2\left\vert f\left( x\right)
\right\vert \left\vert f\left( y\right) \right\vert  \label{2.1.17}
\end{equation}%
for any $x\in E$ and $y\in E^{\perp }.$

If by $\left\Vert f\right\Vert _{E}$ we denote the norm of the functional $f$
restricted to $E,$ i.e., $\left\Vert f\right\Vert _{E}=\sup_{x\in
E\backslash \left\{ 0\right\} }\frac{\left\vert f\left( x\right) \right\vert 
}{\left\Vert x\right\Vert },$ then, on taking the supremum over $x\in E$ and 
$y\in E^{\perp }$ in (\ref{2.1.17}) we deduce%
\begin{equation}
\left\Vert f\right\Vert ^{2}\geq 2\left\Vert f\right\Vert _{E}\cdot
\left\Vert f\right\Vert _{E^{\perp }}  \label{2.1.18}
\end{equation}%
for any $E$ a nonzero linear subspace of the Hilbert space $H$ and a given
functional $f\in H^{\ast }\backslash \left\{ 0\right\} .$

We note that the inequality (\ref{2.1.18}) has been obtained in \cite[Eq. 10]%
{SSD3}.
\end{remark}

\subsection{A Conditional Inequality}

The following result providing a lower bound for the norm product under
suitable conditions holds \cite{DS1} (see also \cite[Theorem 1]{DS}):

\begin{theorem}[Dragomir-S\'{a}ndor, 1986]
\label{t2.2.3}Let $x,y,a,b\in H,$ where $\left( H;\left\langle \cdot ,\cdot
\right\rangle \right) $ is an inner product space, be such that%
\begin{equation}
\left\Vert a\right\Vert ^{2}\leq 2\func{Re}\left\langle x,a\right\rangle
\qquad \text{and\qquad }\left\Vert y\right\Vert ^{2}\leq 2\func{Re}%
\left\langle y,b\right\rangle  \label{2.2.1}
\end{equation}%
holds true. Then%
\begin{multline}
\left\Vert x\right\Vert \left\Vert y\right\Vert \geq \left( 2\func{Re}%
\left\langle x,a\right\rangle -\left\Vert a\right\Vert ^{2}\right) ^{\frac{1%
}{2}}\left( 2\func{Re}\left\langle y,b\right\rangle -\left\Vert b\right\Vert
^{2}\right) ^{\frac{1}{2}}  \label{2.2.2} \\
+\left\vert \left\langle x,y\right\rangle -\left\langle x,b\right\rangle
-\left\langle a,y\right\rangle +\left\langle a,b\right\rangle \right\vert .
\end{multline}
\end{theorem}

\begin{proof}
We follow the proof in \cite{DS1}.

Observe that%
\begin{align}
& \left\vert \left\langle x,y\right\rangle -\left\langle x,b\right\rangle
-\left\langle a,y\right\rangle +\left\langle a,b\right\rangle \right\vert
\label{2.2.3} \\
& =\left\vert \left\langle x-a,y-b\right\rangle \right\vert ^{2}  \notag \\
& \leq \left\Vert x-a\right\Vert ^{2}\left\Vert y-b\right\Vert ^{2}  \notag
\\
& =\left[ \left\Vert x\right\Vert ^{2}-\left( 2\func{Re}\left\langle
x,a\right\rangle -\left\Vert a\right\Vert ^{2}\right) \right] \left[
\left\Vert y\right\Vert ^{2}-\left( 2\func{Re}\left\langle y,b\right\rangle
-\left\Vert b\right\Vert ^{2}\right) \right] .  \notag
\end{align}%
Applying the elementary inequality (\ref{2.1.10}) we have%
\begin{multline}
\left\{ \left\Vert x\right\Vert ^{2}-\left[ \left( 2\func{Re}\left\langle
x,a\right\rangle -\left\Vert a\right\Vert ^{2}\right) ^{\frac{1}{2}}\right]
^{2}\right\}  \label{2.2.4} \\
\times \left\{ \left\Vert y\right\Vert ^{2}-\left[ \left( 2\func{Re}%
\left\langle y,b\right\rangle -\left\Vert b\right\Vert ^{2}\right) ^{\frac{1%
}{2}}\right] ^{2}\right\} \\
\leq \left[ \left\Vert x\right\Vert \left\Vert y\right\Vert -\left( 2\func{Re%
}\left\langle x,a\right\rangle -\left\Vert a\right\Vert ^{2}\right) ^{\frac{1%
}{2}}\left( 2\func{Re}\left\langle y,b\right\rangle -\left\Vert b\right\Vert
^{2}\right) ^{\frac{1}{2}}\right] .
\end{multline}%
Since%
\begin{align*}
0& \leq \left( 2\func{Re}\left\langle x,a\right\rangle -\left\Vert
a\right\Vert ^{2}\right) ^{\frac{1}{2}}\leq \left\Vert x\right\Vert \quad 
\text{and} \\
0& \leq \left( 2\func{Re}\left\langle y,b\right\rangle -\left\Vert
b\right\Vert ^{2}\right) ^{\frac{1}{2}}\leq \left\Vert y\right\Vert
\end{align*}%
hence%
\begin{equation*}
\left\Vert x\right\Vert \left\Vert y\right\Vert \geq \left( 2\func{Re}%
\left\langle x,a\right\rangle -\left\Vert a\right\Vert ^{2}\right) ^{\frac{1%
}{2}}\left( 2\func{Re}\left\langle y,b\right\rangle -\left\Vert b\right\Vert
^{2}\right) ^{\frac{1}{2}}
\end{equation*}%
and by (\ref{2.2.3}) and (\ref{2.2.4}) we deduce the desired result (\ref%
{2.2.2}).
\end{proof}

\begin{remark}
\label{r2.2.4}As pointed out in \cite{DS1}, if we consider $a=\left\langle
x,e\right\rangle e,$ $b=\left\langle y,e\right\rangle e$ with $e\in H,$ $%
\left\Vert e\right\Vert =1$, then the condition (\ref{2.2.1}) is obviously
satisfied and the inequality (\ref{2.2.2}) becomes%
\begin{align}
\left\Vert x\right\Vert \left\Vert y\right\Vert & \geq \left\vert
\left\langle x,e\right\rangle \left\langle e,y\right\rangle \right\vert
+\left\vert \left\langle x,y\right\rangle -\left\langle x,e\right\rangle
\left\langle e,y\right\rangle \right\vert  \label{2.2.5} \\
(& \geq \left\vert \left\langle x,y\right\rangle \right\vert ),  \notag
\end{align}%
which is the refinement of the Schwarz inequality incorporated in (\ref%
{2.1.8}).
\end{remark}

For vectors located in a closed ball centered at $0$ and of radius $\sqrt{2}%
, $ one can state the following corollary as well \cite[Corollary 1.2]{DS}.

\begin{corollary}
\label{c2.2.3}Let $x,y\in H$ such that $\left\Vert x\right\Vert ,\left\Vert
y\right\Vert \leq \sqrt{2}.$ Then%
\begin{multline}
\left\Vert x\right\Vert \left\Vert y\right\Vert \geq \left\vert \left\langle
x,y\right\rangle \right\vert ^{2}\left( 2-\left\Vert x\right\Vert
^{2}\right) ^{\frac{1}{2}}\left( 2-\left\Vert y\right\Vert ^{2}\right) ^{%
\frac{1}{2}}  \label{2.2.6} \\
+\left\vert \left\langle x,y\right\rangle \right\vert \left\vert
1-\left\Vert x\right\Vert ^{2}-\left\Vert y\right\Vert ^{2}+\left\vert
\left\langle x,y\right\rangle \right\vert ^{2}\right\vert .
\end{multline}
\end{corollary}

\begin{proof}
Follows by Theorem \ref{t2.2.3} on choosing $a=\left\langle x,y\right\rangle
y,$ $b=\left\langle y,x\right\rangle x.$ We omit the details.
\end{proof}

\subsection{A Refinement for Orthonormal Families}

The following result provides a generalisation for a refinement of the
Schwarz inequality incorporated in (\ref{2.1.8}) \cite[Theorem 3]{SSD3} (see
also \cite[Theorem]{SSD4} or \cite[Theorem 3]{DS}):

\begin{theorem}[Dragomir, 1985]
\label{t2.3.4}Let $\left( H;\left\langle \cdot ,\cdot \right\rangle \right) $
be an inner product space over the real or complex number field $\mathbb{K}$
and $\left\{ e_{i}\right\} _{i\in I}$ an orthonormal family in $I.$ For any $%
F$ a nonempty finite part of $I$ we have the following refinement of
Schwarz's inequality:%
\begin{align}
\left\Vert x\right\Vert \left\Vert y\right\Vert & \geq \left\vert
\left\langle x,y\right\rangle -\sum_{i\in F}\left\langle
x,e_{i}\right\rangle \left\langle e_{i},y\right\rangle \right\vert
+\sum_{i\in F}\left\vert \left\langle x,e_{i}\right\rangle \left\langle
e_{i},y\right\rangle \right\vert  \label{2.3.1} \\
& \geq \left\vert \left\langle x,y\right\rangle -\sum_{i\in F}\left\langle
x,e_{i}\right\rangle \left\langle e_{i},y\right\rangle \right\vert
+\left\vert \sum_{i\in F}\left\langle x,e_{i}\right\rangle \left\langle
e_{i},y\right\rangle \right\vert  \notag \\
& \geq \left\vert \left\langle x,y\right\rangle \right\vert ,  \notag
\end{align}%
where $x,y\in H.$
\end{theorem}

\begin{proof}
We follow the proof in \cite{SSD3}.

We apply the Schwarz inequality to obtain%
\begin{multline}
\left\vert \left\langle x-\sum_{i\in F}\left\langle x,e_{i}\right\rangle
e_{i},\ \ y-\sum_{i\in F}\left\langle y,e_{i}\right\rangle
e_{i}\right\rangle \right\vert ^{2}  \label{2.3.2} \\
\leq \left\Vert x-\sum_{i\in F}\left\langle x,e_{i}\right\rangle
e_{i}\right\Vert ^{2}\left\Vert y-\sum_{i\in F}\left\langle
y,e_{i}\right\rangle e_{i}\right\Vert ^{2}.
\end{multline}%
Since a simple calculation with orthonormal vectors shows that%
\begin{align*}
\left\Vert x-\sum_{i\in F}\left\langle x,e_{i}\right\rangle e_{i}\right\Vert
^{2}& =\left\Vert x\right\Vert ^{2}-\sum_{i\in F}\left\vert \left\langle
x,e_{i}\right\rangle \right\vert ^{2}, \\
\left\Vert y-\sum_{i\in F}\left\langle y,e_{i}\right\rangle e_{i}\right\Vert
^{2}& =\left\Vert y\right\Vert ^{2}-\sum_{i\in F}\left\vert \left\langle
y,e_{i}\right\rangle \right\vert ^{2},
\end{align*}%
and%
\begin{equation*}
\left\langle x-\sum_{i\in F}\left\langle x,e_{i}\right\rangle e_{i},\ \
y-\sum_{i\in F}\left\langle y,e_{i}\right\rangle e_{i}\right\rangle
=\left\langle x,y\right\rangle -\sum_{i\in F}\left\langle
x,e_{i}\right\rangle \left\langle e_{i},y\right\rangle ,
\end{equation*}%
hence (\ref{2.3.2}) is equivalent to%
\begin{multline}
\left\vert \left\langle x,y\right\rangle -\sum_{i\in F}\left\langle
x,e_{i}\right\rangle \left\langle e_{i},y\right\rangle \right\vert ^{2}
\label{2.3.3} \\
\leq \left( \left\Vert x\right\Vert ^{2}-\sum_{i\in F}\left\vert
\left\langle x,e_{i}\right\rangle \right\vert ^{2}\right) \left( \left\Vert
y\right\Vert ^{2}-\sum_{i\in F}\left\vert \left\langle y,e_{i}\right\rangle
\right\vert ^{2}\right)
\end{multline}%
for any $x,y\in H.$

Further, we need the following Acz\'{e}l type inequality%
\begin{equation}
\left( \alpha ^{2}-\sum_{i\in F}\alpha _{i}^{2}\right) \left( \beta
^{2}-\sum_{i\in F}\beta _{i}^{2}\right) \leq \left( \alpha \beta -\sum_{i\in
F}\alpha _{i}\beta _{i}\right) ^{2},  \label{2.3.4}
\end{equation}%
provided that $\alpha ^{2}\geq \sum_{i\in F}\alpha _{i}^{2}$ and $\beta
^{2}\geq \sum_{i\in F}\beta _{i}^{2}$, where $\alpha ,\beta ,\alpha
_{i},\beta _{i}\in \mathbb{R}$, $i\in F.$

For an Acz\'{e}l inequality that holds under slightly weaker conditions and
a different proof based on polynomials, see \cite[p. 57]{MIT}.

For the sake of completeness, we give here a direct proof of (\ref{2.3.4}).

Utilising the elementary inequality (\ref{2.1.10}), we can write%
\begin{multline}
\left( \alpha ^{2}-\left[ \left( \sum_{i\in F}\alpha _{i}^{2}\right) ^{\frac{%
1}{2}}\right] ^{2}\right) \left( \beta ^{2}-\left[ \left( \sum_{i\in F}\beta
_{i}^{2}\right) ^{\frac{1}{2}}\right] ^{2}\right)  \label{2.3.5} \\
\leq \left[ \left\vert \alpha \beta \right\vert -\left( \sum_{i\in F}\alpha
_{i}^{2}\right) ^{\frac{1}{2}}\left( \sum_{i\in F}\beta _{i}^{2}\right) ^{%
\frac{1}{2}}\right] ^{2}.
\end{multline}%
Since $\left\vert \alpha \right\vert \geq \left( \sum_{i\in F}\alpha
_{i}^{2}\right) ^{\frac{1}{2}}$ and $\left\vert \beta \right\vert \geq
\left( \sum_{i\in F}\beta _{i}^{2}\right) ^{\frac{1}{2}},$ then%
\begin{equation*}
\left\vert \alpha \beta \right\vert \geq \left( \sum_{i\in F}\alpha
_{i}^{2}\right) ^{\frac{1}{2}}\left( \sum_{i\in F}\beta _{i}^{2}\right) ^{%
\frac{1}{2}}.
\end{equation*}%
Therefore, by the Cauchy-Bunyakovsky-Schwarz inequality, we have that%
\begin{align*}
\left\vert \left\vert \alpha \beta \right\vert -\left( \sum_{i\in F}\alpha
_{i}^{2}\right) ^{\frac{1}{2}}\left( \sum_{i\in F}\beta _{i}^{2}\right) ^{%
\frac{1}{2}}\right\vert & =\left\vert \alpha \beta \right\vert -\left(
\sum_{i\in F}\alpha _{i}^{2}\right) ^{\frac{1}{2}}\left( \sum_{i\in F}\beta
_{i}^{2}\right) ^{\frac{1}{2}} \\
& \leq \left\vert \alpha \beta \right\vert -\left\vert \sum_{i\in F}\alpha
_{i}\beta _{i}\right\vert \\
& =\left\vert \left\vert \alpha \beta \right\vert -\left\vert \sum_{i\in
F}\alpha _{i}\beta _{i}\right\vert \right\vert \\
& \leq \left\vert \alpha \beta -\sum_{i\in F}\alpha _{i}\beta
_{i}\right\vert ,
\end{align*}%
showing that%
\begin{equation}
\left[ \left\vert \alpha \beta \right\vert -\left( \sum_{i\in F}\alpha
_{i}^{2}\right) ^{\frac{1}{2}}\left( \sum_{i\in F}\beta _{i}^{2}\right) ^{%
\frac{1}{2}}\right] ^{2}\leq \left( \alpha \beta -\sum_{i\in F}\alpha
_{i}\beta _{i}\right) ^{2}  \label{2.3.6}
\end{equation}%
and then, by (\ref{2.3.5}) and (\ref{2.3.6}) we deduce the desired result (%
\ref{2.3.4}).

By Bessel's inequality we obviously have that%
\begin{equation*}
\left\Vert x\right\Vert ^{2}\geq \sum_{i\in F}\left\vert \left\langle
x,e_{i}\right\rangle \right\vert ^{2}\qquad \text{and\qquad }\left\Vert
y\right\Vert ^{2}\geq \sum_{i\in F}\left\vert \left\langle
y,e_{i}\right\rangle \right\vert ^{2},
\end{equation*}%
therefore, on applying the inequality (\ref{2.3.4}) we deduce that%
\begin{multline}
\left( \left\Vert x\right\Vert ^{2}-\sum_{i\in F}\left\vert \left\langle
x,e_{i}\right\rangle \right\vert ^{2}\right) \left( \left\Vert y\right\Vert
^{2}-\sum_{i\in F}\left\vert \left\langle y,e_{i}\right\rangle \right\vert
^{2}\right)  \label{2.3.7} \\
\leq \left( \left\Vert x\right\Vert \left\Vert y\right\Vert -\sum_{i\in
F}\left\vert \left\langle x,e_{i}\right\rangle \left\langle
e_{i},y\right\rangle \right\vert \right) ^{2}.
\end{multline}%
Since $\left\Vert x\right\Vert \left\Vert y\right\Vert -\sum_{i\in
F}\left\vert \left\langle x,e_{i}\right\rangle \left\langle
e_{i},y\right\rangle \right\vert \geq 0,$ hence by (\ref{2.3.3}) and (\ref%
{2.3.7}) we deduce the first part of (\ref{2.3.1}).

The second and third parts are obvious.
\end{proof}

When the vectors are orthogonal, the following result may be stated \cite%
{SSD4} (see also \cite[Corollary 3.1]{DS}).

\begin{corollary}
\label{c2.3.3}If $\left\{ e_{i}\right\} _{i\in I}$ is an orthonormal family
in $\left( H,\left\langle \cdot ,\cdot \right\rangle \right) $ and $x,y\in H$
with $x\perp y,$ then we have the inequality:%
\begin{align}
\left\Vert x\right\Vert \left\Vert y\right\Vert & \geq \left\vert \sum_{i\in
F}\left\langle x,e_{i}\right\rangle \left\langle e_{i},y\right\rangle
\right\vert +\sum_{i\in F}\left\vert \left\langle x,e_{i}\right\rangle
\left\langle e_{i},y\right\rangle \right\vert  \label{2.3.8} \\
& \geq 2\left\vert \sum_{i\in F}\left\langle x,e_{i}\right\rangle
\left\langle e_{i},y\right\rangle \right\vert ,  \notag
\end{align}%
for any nonempty finite part of $I.$
\end{corollary}

\section{Kurepa Type Refinements for the Schwarz Inequality}

\subsection{Kurepa's Inequality}

In 1960, N.G. de Bruijn proved the following refinement of the celebrated
Cauchy-Bunyakovsky-Schwarz (CBS) inequality for a sequence of real numbers
and the second of complex numbers, see \cite{BR} or \cite[p. 48]{SSD1}:

\begin{theorem}[de Bruijn, 1960]
\label{t3.1.1}Let $\left( a_{1},\dots ,a_{n}\right) $ be an $n-$tuple of
real numbers and $\left( z_{1},\dots ,z_{n}\right) $ an $n-$tuple of complex
numbers. Then%
\begin{align}
\left\vert \sum_{k=1}^{n}a_{k}z_{k}\right\vert ^{2}& \leq \frac{1}{2}%
\sum_{k=1}^{n}a_{k}^{2}\left[ \sum_{k=1}^{n}\left\vert z_{k}\right\vert
^{2}+\left\vert \sum_{k=1}^{n}z_{k}^{2}\right\vert \right]  \label{3.1.1} \\
& \left( \leq \sum_{k=1}^{n}a_{k}^{2}\cdot \sum_{k=1}^{n}\left\vert
z_{k}\right\vert ^{2}\right) .  \notag
\end{align}%
Equality holds in (\ref{3.1.1}) if and only if, for $k\in \left\{ 1,\dots
,n\right\} ,$ $a_{k}=\func{Re}\left( \lambda z_{k}\right) ,$ where $\lambda $
is a complex number such that $\lambda ^{2}\sum_{k=1}^{n}z_{n}^{2}$ is a
nonnegative real number.
\end{theorem}

In 1966, in an effort to extend this result to inner products, Kurepa \cite%
{KU} obtained the following refinement for the complexification of a real
inner product space $\left( H;\left\langle \cdot ,\cdot \right\rangle
\right) :$

\begin{theorem}[Kurepa, 1966]
\label{t3.1.2}Let $\left( H;\left\langle \cdot ,\cdot \right\rangle \right) $
be a real inner product space and $\left( H_{\mathbb{C}},\left\langle \cdot
,\cdot \right\rangle _{\mathbb{C}}\right) $ its complexification. For any $%
a\in H$ and $z\in H_{\mathbb{C}}$ we have the inequality:%
\begin{align}
\left\vert \left\langle z,a\right\rangle _{\mathbb{C}}\right\vert ^{2}& \leq 
\frac{1}{2}\left\Vert a\right\Vert ^{2}\left[ \left\Vert z\right\Vert _{%
\mathbb{C}}^{2}+\left\vert \left\langle z,\bar{z}\right\rangle _{\mathbb{C}%
}\right\vert \right]  \label{3.1.2} \\
& \left( \leq \left\Vert a\right\Vert ^{2}\left\Vert z\right\Vert _{\mathbb{C%
}}^{2}\right) .  \notag
\end{align}
\end{theorem}

To be comprehensive, we define in the following the concept of
complexification for a real inner product space.

Let $H$ be a real inner product space with the scalar product $\left\langle
\cdot ,\cdot \right\rangle $ and the norm $\left\Vert \cdot \right\Vert .$
The \textit{complexification }$H_{\mathbb{C}}$ of $H$ is defined as a
complex linear space $H\times H$ of all ordered pairs $\left( x,y\right) $ $%
\left( x,y\in H\right) $ endowed with the operations%
\begin{align*}
\left( x,y\right) +\left( x^{\prime },y^{\prime }\right) & :=\left(
x+x^{\prime },y+y^{\prime }\right) ,\qquad x,x^{\prime },y,y^{\prime }\in H;
\\
\left( \sigma +i\tau \right) \cdot \left( x,y\right) & :=\left( \sigma
x-\tau y,\tau x+\sigma y\right) ,\qquad x,y\in H\text{ \ and \ }\sigma ,\tau
\in \mathbb{R}.
\end{align*}%
On $H_{\mathbb{C}}$ one can canonically consider the \textit{scalar product} 
$\left\langle \cdot ,\cdot \right\rangle _{\mathbb{C}}$ defined by:%
\begin{equation*}
\left\langle z,z^{\prime }\right\rangle _{\mathbb{C}}:=\left\langle
x,x^{\prime }\right\rangle +\left\langle y,y^{\prime }\right\rangle +i\left[
\left\langle y,x^{\prime }\right\rangle -\left\langle x,y^{\prime
}\right\rangle \right]
\end{equation*}%
where $z=\left( x,y\right) ,$ $z^{\prime }=\left( x^{\prime },y^{\prime
}\right) \in H_{\mathbb{C}}.$ Obviously,%
\begin{equation*}
\left\Vert z\right\Vert _{\mathbb{C}}^{2}=\left\Vert x\right\Vert
^{2}+\left\Vert y\right\Vert ^{2},
\end{equation*}%
where $z=\left( x,y\right) .$

The conjugate of a vector $z=\left( x,y\right) \in H_{\mathbb{C}}$ is
defined by $\bar{z}:=\left( x,-y\right) .$

It is easy to see that the elements of $H_{\mathbb{C}}$ under defined
operations behave as formal \textquotedblleft complex\textquotedblright\
combinations $x+iy$ with $x,y\in H.$ Because of this, we may write $z=x+iy$
instead of $z=\left( x,y\right) .$ Thus, $\bar{z}=x-iy.$

\subsection{A Generalisation of Kurepa's Inequality}

The following lemma is of interest \cite{SSD2}.

\begin{lemma}
\label{l3.1.1}Let $f:\left[ 0,2\pi \right] \rightarrow \mathbb{R}$ given by%
\begin{equation}
f\left( \alpha \right) =\lambda \sin ^{2}\alpha +2\beta \sin \alpha \cos
\alpha +\alpha \cos ^{2}\alpha ,  \label{3.1.3}
\end{equation}%
where $\lambda ,\beta ,\gamma \in \mathbb{R}$. Then%
\begin{equation}
\sup_{\alpha \in \left[ 0,2\pi \right] }f\left( \alpha \right) =\frac{1}{2}%
\left( \lambda +\gamma \right) +\frac{1}{2}\left[ \left( \gamma -\lambda
\right) ^{2}+4\beta ^{2}\right] ^{\frac{1}{2}}.  \label{3.1.4}
\end{equation}
\end{lemma}

\begin{proof}
Since%
\begin{equation*}
\sin ^{2}\alpha =\frac{1-\cos 2\alpha }{2},\quad \cos ^{2}\alpha =\frac{%
1+\cos 2\alpha }{2},\quad 2\sin \alpha \cos \alpha =\sin 2\alpha ,
\end{equation*}%
hence $f$ may be written as%
\begin{equation}
f\left( \alpha \right) =\frac{1}{2}\left( \lambda +\gamma \right) +\frac{1}{2%
}\left( \gamma -\lambda \right) \cos 2\alpha +\beta \sin 2\alpha .
\label{3.1.5}
\end{equation}%
If $\beta =0,$ then (\ref{3.1.5}) becomes%
\begin{equation*}
f\left( \alpha \right) =\frac{1}{2}\left( \lambda +\gamma \right) +\frac{1}{2%
}\left( \gamma -\lambda \right) \cos 2\alpha .
\end{equation*}%
Obviously, in this case%
\begin{equation*}
\sup_{\alpha \in \left[ 0,2\pi \right] }f\left( \alpha \right) =\frac{1}{2}%
\left( \lambda +\gamma \right) +\frac{1}{2}\left\vert \gamma -\lambda
\right\vert =\max \left\{ \gamma ,\lambda \right\} .
\end{equation*}%
If $\beta \neq 0,$ then (\ref{3.1.5}) becomes%
\begin{equation*}
f\left( \alpha \right) =\frac{1}{2}\left( \lambda +\gamma \right) +\beta %
\left[ \sin 2\alpha +\frac{\left( \gamma -\lambda \right) }{\beta }\cos
2\alpha \right] .
\end{equation*}%
Let $\varphi \in \left( -\frac{\pi }{2},\frac{\pi }{2}\right) $ for which $%
\tan \varphi =\frac{\gamma -\lambda }{2\beta }.$ Then $f$ can be written as%
\begin{equation*}
f\left( \alpha \right) =\frac{1}{2}\left( \lambda +\gamma \right) +\frac{%
\beta }{\cos \varphi }\sin \left( 2\alpha +\varphi \right) .
\end{equation*}%
For this function, obviously%
\begin{equation}
\sup_{\alpha \in \left[ 0,2\pi \right] }f\left( \alpha \right) =\frac{1}{2}%
\left( \lambda +\gamma \right) +\frac{\left\vert \beta \right\vert }{%
\left\vert \cos \varphi \right\vert }.  \label{3.1.6}
\end{equation}%
Since%
\begin{equation*}
\frac{\sin ^{2}\varphi }{\cos ^{2}\varphi }=\frac{\left( \gamma -\lambda
\right) ^{2}}{4\beta ^{2}},
\end{equation*}%
hence,%
\begin{equation*}
\frac{1}{\left\vert \cos \varphi \right\vert }=\frac{\left[ \left( \gamma
-\lambda \right) ^{2}+4\beta ^{2}\right] ^{\frac{1}{2}}}{2\left\vert \beta
\right\vert },
\end{equation*}%
and from (\ref{3.1.6}) we deduce the desired result (\ref{3.1.4}).
\end{proof}

The following result holds \cite{SSD2}.

\begin{theorem}[Dragomir, 2004]
\label{t3.1.3}Let $\left( H;\left\langle \cdot ,\cdot \right\rangle \right) $
be a complex inner product space. If $x,y,z\in H$ are such that%
\begin{equation}
\func{Im}\left\langle x,z\right\rangle =\func{Im}\left\langle
y,z\right\rangle =0,  \label{3.1.7}
\end{equation}%
then we have the inequality:%
\begin{align}
& \func{Re}^{2}\left\langle x,z\right\rangle +\func{Re}^{2}\left\langle
y,z\right\rangle  \label{3.1.8} \\
& =\left\vert \left\langle x+iy,z\right\rangle \right\vert ^{2}  \notag \\
& \leq \frac{1}{2}\left\{ \left\Vert x\right\Vert ^{2}+\left\Vert
y\right\Vert ^{2}+\left[ \left( \left\Vert x\right\Vert ^{2}-\left\Vert
y\right\Vert ^{2}\right) ^{2}-4\func{Re}^{2}\left\langle x,y\right\rangle %
\right] ^{\frac{1}{2}}\right\} \left\Vert z\right\Vert ^{2}  \notag \\
& \leq \left( \left\Vert x\right\Vert ^{2}+\left\Vert y\right\Vert
^{2}\right) \left\Vert z\right\Vert ^{2}.  \notag
\end{align}
\end{theorem}

\begin{proof}
Obviously, by (\ref{3.1.7}), we have%
\begin{equation*}
\left\langle x+iy,z\right\rangle =\func{Re}\left\langle x,z\right\rangle +i%
\func{Re}\left\langle y,z\right\rangle
\end{equation*}%
and the first part of (\ref{3.1.8}) holds true.

Now, let $\varphi \in \left[ 0,2\pi \right] $ be such that $\left\langle
x+iy,z\right\rangle =e^{i\varphi }\left\vert \left\langle
x+iy,z\right\rangle \right\vert .$ Then%
\begin{equation*}
\left\vert \left\langle x+iy,z\right\rangle \right\vert =e^{-i\varphi
}\left\langle x+iy,z\right\rangle =\left\langle e^{-i\varphi }\left(
x+iy\right) ,z\right\rangle .
\end{equation*}%
Utilising the above identity, we can write:%
\begin{align*}
\left\vert \left\langle x+iy,z\right\rangle \right\vert & =\func{Re}%
\left\langle e^{-i\varphi }\left( x+iy\right) ,z\right\rangle \\
& =\func{Re}\left\langle \left( \cos \varphi -i\sin \varphi \right) \left(
x+iy\right) ,z\right\rangle \\
& =\func{Re}\left\langle \cos \varphi \cdot x+\sin \varphi \cdot y-i\sin
\varphi \cdot x+i\cos \varphi \cdot y,z\right\rangle \\
& =\func{Re}\left\langle \cos \varphi \cdot x+\sin \varphi \cdot
y,z\right\rangle +\func{Im}\left\langle \sin \varphi \cdot x-\cos \varphi
\cdot y,z\right\rangle \\
& =\func{Re}\left\langle \cos \varphi \cdot x+\sin \varphi \cdot
y,z\right\rangle +\sin \varphi \func{Im}\left\langle x,z\right\rangle -\cos
\varphi \func{Im}\left\langle y,z\right\rangle \\
& =\func{Re}\left\langle \cos \varphi \cdot x+\sin \varphi \cdot
y,z\right\rangle ,
\end{align*}%
and for the last equality we have used the assumption (\ref{3.1.7}).

Taking the square and using the Schwarz inequality for the inner product $%
\left\langle \cdot ,\cdot \right\rangle ,$ we have%
\begin{align}
\left\vert \left\langle x+iy,z\right\rangle \right\vert ^{2}& =\left[ \func{%
Re}\left\langle \cos \varphi \cdot x+\sin \varphi \cdot y,z\right\rangle %
\right] ^{2}  \label{3.1.9} \\
& \leq \left\Vert \cos \varphi \cdot x+\sin \varphi \cdot y\right\Vert
^{2}\left\Vert z\right\Vert ^{2}.  \notag
\end{align}%
On making use of Lemma \ref{l3.1.1}, we have%
\begin{align*}
& \sup_{\alpha \in \left[ 0,2\pi \right] }\left\Vert \cos \varphi \cdot
x+\sin \varphi \cdot y\right\Vert ^{2} \\
& =\sup_{\alpha \in \left[ 0,2\pi \right] }\left[ \left\Vert x\right\Vert
^{2}\cos ^{2}\varphi +2\func{Re}\left\langle x,y\right\rangle \sin \varphi
\cos \varphi +\left\Vert y\right\Vert ^{2}\sin {}^{2}\varphi \right] \\
& =\frac{1}{2}\left\{ \left\Vert x\right\Vert ^{2}+\left\Vert y\right\Vert
^{2}+\left[ \left( \left\Vert x\right\Vert ^{2}-\left\Vert y\right\Vert
^{2}\right) ^{2}+4\func{Re}^{2}\left\langle x,y\right\rangle \right] ^{\frac{%
1}{2}}\right\}
\end{align*}%
and the first inequality in (\ref{3.1.8}) is proved.

Observe that%
\begin{align*}
& \left( \left\Vert x\right\Vert ^{2}-\left\Vert y\right\Vert ^{2}\right)
^{2}+4\func{Re}^{2}\left\langle x,y\right\rangle \\
& =\left( \left\Vert x\right\Vert ^{2}+\left\Vert y\right\Vert ^{2}\right)
^{2}-4\left[ \left\Vert x\right\Vert ^{2}\left\Vert y\right\Vert ^{2}-\func{%
Re}^{2}\left\langle x,y\right\rangle \right] \\
& \leq \left( \left\Vert x\right\Vert ^{2}+\left\Vert y\right\Vert
^{2}\right) ^{2}
\end{align*}%
and the last part of (\ref{3.1.8}) is proved.
\end{proof}

\begin{remark}
\label{r3.1.4}Observe that if $\left( H,\left\langle \cdot ,\cdot
\right\rangle \right) $ is a real inner product space, then for any $%
x,y,z\in H$ one has:%
\begin{align}
& \left\langle x,z\right\rangle ^{2}+\left\langle y,z\right\rangle ^{2}
\label{3.1.10} \\
& \leq \frac{1}{2}\left\{ \left\Vert x\right\Vert ^{2}+\left\Vert
y\right\Vert ^{2}+\left[ \left( \left\Vert x\right\Vert ^{2}-\left\Vert
y\right\Vert ^{2}\right) ^{2}+4\left\langle x,y\right\rangle ^{2}\right]
\right\} ^{\frac{1}{2}}\left\Vert z\right\Vert ^{2}  \notag \\
& \leq \left( \left\Vert x\right\Vert ^{2}+\left\Vert y\right\Vert
^{2}\right) \left\Vert z\right\Vert ^{2}.  \notag
\end{align}
\end{remark}

\begin{remark}
\label{r3.1.5}If $H$ is a real space, $\left\langle \cdot ,\cdot
\right\rangle $ the real inner product, $H_{\mathbb{C}}$ its
complexification and $\left\langle \cdot ,\cdot \right\rangle _{\mathbb{C}}$
the corresponding complexification for $\left\langle \cdot ,\cdot
\right\rangle $, then for $x,y\in H$ and $w:=x+iy\in H_{\mathbb{C}}$ and for 
$e\in H$ we have%
\begin{equation*}
\func{Im}\left\langle x,e\right\rangle _{\mathbb{C}}=\func{Im}\left\langle
y,e\right\rangle _{\mathbb{C}}=0,
\end{equation*}%
\begin{equation*}
\left\Vert w\right\Vert _{\mathbb{C}}^{2}=\left\Vert x\right\Vert
^{2}+\left\Vert y\right\Vert ^{2},\qquad \left\vert \left\langle w,\bar{w}%
\right\rangle _{\mathbb{C}}\right\vert =\left( \left\Vert x\right\Vert
^{2}-\left\Vert y\right\Vert ^{2}\right) ^{2}+4\left\langle x,y\right\rangle
^{2},
\end{equation*}%
where $\bar{w}=x-iy\in H_{\mathbb{C}}.$

Applying Theorem \ref{t3.1.3} for the complex space $H_{\mathbb{C}}$ and
complex inner product $\left\langle \cdot ,\cdot \right\rangle _{\mathbb{C}%
}, $ we deduce%
\begin{equation}
\left\vert \left\langle w,e\right\rangle _{\mathbb{C}}\right\vert ^{2}\leq 
\frac{1}{2}\left\Vert e\right\Vert ^{2}\left[ \left\Vert w\right\Vert _{%
\mathbb{C}}^{2}+\left\vert \left\langle w,\bar{w}\right\rangle _{\mathbb{C}%
}\right\vert \right] \leq \left\Vert e\right\Vert ^{2}\left\Vert
w\right\Vert _{\mathbb{C}}^{2},  \label{3.1.11}
\end{equation}%
which is Kurepa's inequality (\ref{3.1.2}).
\end{remark}

\begin{corollary}
\label{c3.1.1}Let $x,y,z$ be as in Theorem \ref{t3.1.3}. In addition, if $%
\func{Re}\left\langle x,y\right\rangle =0,$ then%
\begin{equation}
\left[ \func{Re}^{2}\left\langle x,z\right\rangle +\func{Re}^{2}\left\langle
y,z\right\rangle \right] ^{\frac{1}{2}}\leq \left\Vert z\right\Vert \cdot
\max \left\{ \left\Vert x\right\Vert ,\left\Vert y\right\Vert \right\} .
\label{3.1.12}
\end{equation}
\end{corollary}

\begin{remark}
\label{r3.1.6}If $H$ is a real space and $\left\langle \cdot ,\cdot
\right\rangle $ a real inner product on $H,$ then for any $x,y,z\in H$ with $%
\left\langle x,y\right\rangle =0$ we have%
\begin{equation}
\left[ \left\langle x,z\right\rangle ^{2}+\left\langle y,z\right\rangle ^{2}%
\right] ^{\frac{1}{2}}\leq \left\Vert z\right\Vert \cdot \max \left\{
\left\Vert x\right\Vert ,\left\Vert y\right\Vert \right\} .  \label{3.1.13}
\end{equation}
\end{remark}

\subsection{A Related Result}

Utilising Lemma \ref{l3.1.1}, we may state and prove the following result as
well.

\begin{theorem}[Dragomir, 2004]
\label{t3.1.4}Let $\left( H,\left\langle \cdot ,\cdot \right\rangle \right) $
\ be a real or complex inner product space. Then we have the inequalities:%
\begin{align}
& \frac{1}{2}\bigg\{\left\vert \left\langle v,t\right\rangle \right\vert
^{2}+\left\vert \left\langle w,t\right\rangle \right\vert ^{2}+\left[ \left(
\left\vert \left\langle v,t\right\rangle \right\vert ^{2}-\left\vert
\left\langle w,t\right\rangle \right\vert ^{2}\right) ^{2}\right.
\label{3.1.14} \\
& \quad \left. +\left. 4\left( \func{Re}\left\langle v,t\right\rangle \func{%
Re}\left\langle w,t\right\rangle +\func{Im}\left\langle v,t\right\rangle 
\func{Im}\left\langle w,t\right\rangle \right) ^{2}\right] ^{\frac{1}{2}%
}\right\}  \notag \\
& \leq \frac{1}{2}\left\Vert t\right\Vert ^{2}\left\{ \left\Vert
v\right\Vert ^{2}+\left\Vert w\right\Vert ^{2}+\left[ \left( \left\Vert
v\right\Vert ^{2}-\left\Vert w\right\Vert ^{2}\right) ^{2}+4\func{Re}%
^{2}\left( v,w\right) \right] ^{\frac{1}{2}}\right\}  \notag \\
& \leq \left( \left\Vert v\right\Vert ^{2}+\left\Vert w\right\Vert
^{2}\right) \left\Vert t\right\Vert ^{2},  \notag
\end{align}%
for all $v,w,t\in H.$
\end{theorem}

\begin{proof}
Observe that, by Schwarz's inequality%
\begin{equation}
\left\vert \left( \cos \varphi \cdot v+\sin \varphi \cdot w,z\right)
\right\vert ^{2}\leq \left\Vert \cos \varphi \cdot v+\sin \varphi \cdot
w\right\Vert ^{2}\left\Vert z\right\Vert ^{2}  \label{3.1.15}
\end{equation}%
for any $\varphi \in \left[ 0,2\pi \right] .$

Since%
\begin{align*}
I\left( \varphi \right) & :=\left\Vert \cos \varphi \cdot v+\sin \varphi
\cdot w\right\Vert ^{2} \\
& =\cos ^{2}\varphi \left\Vert v\right\Vert ^{2}+2\func{Re}\left( v,w\right)
\sin \varphi \cos \varphi +\left\Vert w\right\Vert ^{2}\sin ^{2}\varphi ,
\end{align*}%
hence, as in Theorem \ref{t3.1.3},%
\begin{equation*}
\sup_{\varphi \in \left[ 0,2\pi \right] }I\left( \varphi \right) =\frac{1}{2}%
\left\{ \left\Vert v\right\Vert ^{2}+\left\Vert w\right\Vert ^{2}+\left[
\left( \left\Vert v\right\Vert ^{2}-\left\Vert w\right\Vert ^{2}\right)
^{2}+4\func{Re}^{2}\left( v,w\right) \right] ^{\frac{1}{2}}\right\} .
\end{equation*}%
Also, denoting%
\begin{align*}
J\left( \varphi \right) & :=\left\vert \cos \varphi \left\langle
v,z\right\rangle +\sin \varphi \left\langle w,z\right\rangle \right\vert \\
& =\cos ^{2}\varphi \left\vert \left\langle v,z\right\rangle \right\vert
^{2}+2\sin \varphi \cos \varphi \func{Re}\left[ \left\langle
v,z\right\rangle \overline{\left\langle w,z\right\rangle }\right] +\sin
^{2}\varphi \left\vert \left\langle w,z\right\rangle \right\vert ^{2},
\end{align*}%
then, on applying Lemma \ref{l3.1.1}, we deduce that%
\begin{multline*}
\sup_{\varphi \in \left[ 0,2\pi \right] }J\left( \varphi \right) =\frac{1}{2}%
\bigg\{\left\vert \left\langle v,t\right\rangle \right\vert ^{2}+\left\vert
\left\langle w,t\right\rangle \right\vert ^{2} \\
+\left. \left[ \left( \left\vert \left\langle v,t\right\rangle \right\vert
^{2}-\left\vert \left\langle w,t\right\rangle \right\vert ^{2}\right) ^{2}+4%
\func{Re}^{2}\left[ \left\langle v,z\right\rangle \overline{\left\langle
w,z\right\rangle }\right] \right] ^{\frac{1}{2}}\right\}
\end{multline*}%
and, since%
\begin{equation*}
\func{Re}\left[ \left\langle v,z\right\rangle \overline{\left\langle
w,z\right\rangle }\right] =\func{Re}\left\langle v,t\right\rangle \func{Re}%
\left\langle w,t\right\rangle +\func{Im}\left\langle v,t\right\rangle \func{%
Im}\left\langle w,t\right\rangle ,
\end{equation*}%
hence, on taking the supremum in the inequality (\ref{3.1.15}), we deduce
the desired inequality (\ref{3.1.14}).
\end{proof}

\begin{remark}
\label{r3.1.7}In the real case, (\ref{3.1.14}) provides the same inequality
we obtained in (\ref{3.1.10}).

In the complex case, if we assume that \ $v,w,t\in H$ are such that%
\begin{equation*}
\func{Re}\left\langle v,t\right\rangle \func{Re}\left\langle
w,t\right\rangle =-\func{Im}\left\langle v,t\right\rangle \func{Im}%
\left\langle w,t\right\rangle ,
\end{equation*}%
then (\ref{3.1.14}) becomes:%
\begin{multline}
\max \left\{ \left\vert \left\langle v,t\right\rangle \right\vert
^{2},\left\vert \left\langle w,t\right\rangle \right\vert ^{2}\right\}
\label{3.1.1.6} \\
\leq \frac{1}{2}\left\Vert t\right\Vert ^{2}\left\{ \left\Vert v\right\Vert
^{2}+\left\Vert w\right\Vert ^{2}+\left[ \left( \left\Vert v\right\Vert
^{2}-\left\Vert w\right\Vert ^{2}\right) ^{2}+4\func{Re}^{2}\left(
v,w\right) \right] ^{\frac{1}{2}}\right\} .
\end{multline}
\end{remark}

\section{Refinements of Buzano's and Kurepa's Inequalities}

\subsection{Introduction}

In \cite{B}, M.L. Buzano obtained the following extension of the celebrated
Schwarz's inequality in a real or complex inner product space $\left(
H;\left\langle \cdot ,\cdot \right\rangle \right) :$%
\begin{equation}
\left\vert \left\langle a,x\right\rangle \left\langle x,b\right\rangle
\right\vert \leq \frac{1}{2}\left[ \left\Vert a\right\Vert \cdot \left\Vert
b\right\Vert +\left\vert \left\langle a,b\right\rangle \right\vert \right]
\left\Vert x\right\Vert ^{2},  \label{ch3.1.1}
\end{equation}%
for any $a,b,x\in H.$

It is clear that for $a=b,$ the above inequality becomes the standard
Schwarz inequality%
\begin{equation}
\left\vert \left\langle a,x\right\rangle \right\vert ^{2}\leq \left\Vert
a\right\Vert ^{2}\left\Vert x\right\Vert ^{2},\qquad a,x\in H;
\label{ch3.1.2}
\end{equation}%
with equality if and only if there exists a scalar $\lambda \in \mathbb{K}$ $%
\left( \mathbb{K}=\mathbb{R}\text{ \ or }\mathbb{C}\right) $ such that $%
x=\lambda a.$

As noted by M. Fujii and F. Kubo in \cite{FK}, where they provided a simple
proof of (\ref{ch3.1.1}) by utilising orthogonal projection arguments, the
case of equality holds in (\ref{ch3.1.1}) if%
\begin{equation*}
x=\left\{ 
\begin{array}{ll}
\alpha \left( \frac{a}{\left\Vert a\right\Vert }+\frac{\left\langle
a,b\right\rangle }{\left\vert \left\langle a,b\right\rangle \right\vert }%
\cdot \frac{b}{\left\Vert b\right\Vert }\right) , & \text{when \ }%
\left\langle a,b\right\rangle \neq 0 \\ 
&  \\ 
\alpha \left( \frac{a}{\left\Vert a\right\Vert }+\beta \cdot \frac{b}{%
\left\Vert b\right\Vert }\right) , & \text{when\ }\left\langle
a,b\right\rangle =0,%
\end{array}%
\right.
\end{equation*}%
where $\alpha ,\beta \in \mathbb{K}$.

It might be useful to observe that, out of (\ref{ch3.1.1}), one may get the
following discrete inequality:%
\begin{multline}
\left\vert \sum_{i=1}^{n}p_{i}a_{i}\overline{x_{i}}\sum_{i=1}^{n}p_{i}x_{i}%
\overline{b_{i}}\right\vert  \label{ch3.1.3} \\
\leq \frac{1}{2}\left[ \left( \sum_{i=1}^{n}p_{i}\left\vert a_{i}\right\vert
^{2}\sum_{i=1}^{n}p_{i}\left\vert b_{i}\right\vert ^{2}\right) ^{\frac{1}{2}%
}+\left\vert \sum_{i=1}^{n}p_{i}a_{i}\overline{b_{i}}\right\vert \right]
\sum_{i=1}^{n}p_{i}\left\vert x_{i}\right\vert ^{2},
\end{multline}%
where $p_{i}\geq 0,$ $a_{i},x_{i},b_{i}\in \mathbb{C}$, $i\in \left\{
1,\dots ,n\right\} .$

If one takes in (\ref{ch3.1.3}) $b_{i}=\overline{a_{i}}$ for $i\in \left\{
1,\dots ,n\right\} ,$ then one obtains%
\begin{equation}
\left\vert \sum_{i=1}^{n}p_{i}a_{i}\overline{x_{i}}%
\sum_{i=1}^{n}p_{i}a_{i}x_{i}\right\vert \leq \frac{1}{2}\left[
\sum_{i=1}^{n}p_{i}\left\vert a_{i}\right\vert ^{2}+\left\vert
\sum_{i=1}^{n}p_{i}a_{i}^{2}\right\vert \right] \sum_{i=1}^{n}p_{i}\left%
\vert x_{i}\right\vert ^{2},  \label{ch3.1.4}
\end{equation}%
for any $p_{i}\geq 0,$ $a_{i},x_{i},b_{i}\in \mathbb{C}$, $i\in \left\{
1,\dots ,n\right\} .$

Note that, if $x_{i},$ $i\in \left\{ 1,\dots ,n\right\} $ are real numbers,
then out of (\ref{ch3.1.4}), we may deduce the de Bruijn refinement of the
celebrated Cauchy-Bunyakovsky-Schwarz inequality \cite{BR}%
\begin{equation}
\left\vert \sum_{i=1}^{n}p_{i}x_{i}z_{i}\right\vert ^{2}\leq \frac{1}{2}%
\sum_{i=1}^{n}p_{i}x_{i}^{2}\left[ \sum_{i=1}^{n}p_{i}\left\vert
z_{i}\right\vert ^{2}+\left\vert \sum_{i=1}^{n}p_{i}z_{i}^{2}\right\vert %
\right] ,  \label{ch3.1.5}
\end{equation}%
where $z_{i}\in \mathbb{C}$, $i\in \left\{ 1,\dots ,n\right\} .$ In this
way, Buzano's result may be regarded as a generalisation of de Bruijn's
inequality.

Similar comments obviously apply for integrals, but, for the sake of brevity
we do not mention them here.

The aim of the present section is to establish some related results as well
as a refinement of Buzano's inequality for real or complex inner product
spaces. An improvement of Kurepa's inequality for the complexification of a
real inner product and the corresponding applications for discrete and
integral inequalities are also provided.

\subsection{Some Buzano Type Inequalities}

The following result may be stated \cite{DRAG1}.

\begin{theorem}[Dragomir, 2004]
\label{ch3.t2.1}Let $\left( H;\left\langle \cdot ,\cdot \right\rangle
\right) $ be an inner product space over the real or complex number field $%
\mathbb{K}$. For all $\alpha \in \mathbb{K}\backslash \left\{ 0\right\} $
and $x,a,b\in H,$ $\alpha \neq 0,$ one has the inequality%
\begin{multline}
\left\vert \frac{\left\langle a,x\right\rangle \left\langle x,b\right\rangle 
}{\left\Vert x\right\Vert ^{2}}-\frac{\left\langle a,b\right\rangle }{\alpha 
}\right\vert  \label{ch3.2.1} \\
\leq \frac{\left\Vert b\right\Vert }{\left\vert \alpha \right\vert
\left\Vert x\right\Vert }\left[ \left\vert \alpha -1\right\vert
^{2}\left\vert \left\langle a,x\right\rangle \right\vert ^{2}+\left\Vert
x\right\Vert ^{2}\left\Vert a\right\Vert ^{2}-\left\vert \left\langle
a,x\right\rangle \right\vert ^{2}\right] .
\end{multline}%
The case of equality holds in (\ref{ch3.2.1}) if and only if there exists a
scalar $\lambda \in \mathbb{K}$ so that%
\begin{equation}
\alpha \cdot \frac{\left\langle a,x\right\rangle }{\left\Vert x\right\Vert
^{2}}x=a+\lambda b.  \label{ch3.2.1.a}
\end{equation}
\end{theorem}

\begin{proof}
We follow the proof in \cite{DRAG1}.

Using Schwarz's inequality, we have that%
\begin{equation}
\left\vert \left\langle \alpha \cdot \frac{\left\langle a,x\right\rangle }{%
\left\Vert x\right\Vert ^{2}}x-a,b\right\rangle \right\vert ^{2}\leq
\left\Vert \alpha \cdot \frac{\left\langle a,x\right\rangle }{\left\Vert
x\right\Vert ^{2}}x-a\right\Vert ^{2}\left\Vert b\right\Vert ^{2}
\label{ch3.2.2}
\end{equation}%
and since%
\begin{align*}
\left\Vert \alpha \cdot \frac{\left\langle a,x\right\rangle }{\left\Vert
x\right\Vert ^{2}}x-a\right\Vert ^{2}& =\left\vert \alpha \right\vert ^{2}%
\frac{\left\vert \left\langle a,x\right\rangle \right\vert ^{2}}{\left\Vert
x\right\Vert ^{2}}-2\frac{\left\vert \left\langle a,x\right\rangle
\right\vert ^{2}}{\left\Vert x\right\Vert ^{2}}\func{Re}\alpha +\left\Vert
a\right\Vert ^{2} \\
& =\frac{\left\vert \alpha -1\right\vert ^{2}\left\vert \left\langle
a,x\right\rangle \right\vert ^{2}+\left\Vert x\right\Vert ^{2}\left\Vert
a\right\Vert ^{2}-\left\vert \left\langle a,x\right\rangle \right\vert ^{2}}{%
\left\Vert x\right\Vert ^{2}}
\end{align*}%
and%
\begin{equation*}
\left\langle \alpha \cdot \frac{\left\langle a,x\right\rangle }{\left\Vert
x\right\Vert ^{2}}x-a,b\right\rangle =\alpha \left[ \frac{\left\langle
a,x\right\rangle \left\langle x,b\right\rangle }{\left\Vert x\right\Vert ^{2}%
}-\frac{\left\langle a,b\right\rangle }{\alpha }\right] ,
\end{equation*}%
hence by (\ref{ch3.2.1}) we deduce the desired inequality (\ref{ch3.2.1}).

The case of equality is obvious from the above considerations related to the
Schwarz's inequality (\ref{ch3.1.2}).
\end{proof}

\begin{remark}
\label{ch3.r2.2}Using the continuity property of the modulus, i.e., $%
\left\vert \left\vert z\right\vert -\left\vert u\right\vert \right\vert \leq
\left\vert z-u\right\vert ,$ $z,u\in \mathbb{K}$, we have:%
\begin{equation}
\left\vert \frac{\left\vert \left\langle a,x\right\rangle \left\langle
x,b\right\rangle \right\vert }{\left\Vert x\right\Vert ^{2}}-\frac{%
\left\vert \left\langle a,b\right\rangle \right\vert }{\left\vert \alpha
\right\vert }\right\vert \leq \left\vert \frac{\left\langle a,x\right\rangle
\left\langle x,b\right\rangle }{\left\Vert x\right\Vert ^{2}}-\frac{%
\left\langle a,b\right\rangle }{\alpha }\right\vert .  \label{ch3.2.3}
\end{equation}%
Therefore, by (\ref{ch3.2.1}) and (\ref{ch3.2.3}), one may deduce the
following double inequality:%
\begin{align}
& \frac{1}{\left\vert \alpha \right\vert }\left[ \left\vert \left\langle
a,b\right\rangle \right\vert -\frac{\left\Vert b\right\Vert }{\left\Vert
x\right\Vert }\right.  \label{ch3.2.4} \\
& \qquad \times \left. \left[ \left( \left\vert \alpha -1\right\vert
^{2}\left\vert \left\langle x,a\right\rangle \right\vert ^{2}+\left\Vert
x\right\Vert ^{2}\left\Vert a\right\Vert ^{2}-\left\vert \left\langle
a,x\right\rangle \right\vert ^{2}\right) ^{\frac{1}{2}}\right] \right] 
\notag \\
& \leq \frac{\left\vert \left\langle a,x\right\rangle \left\langle
x,b\right\rangle \right\vert }{\left\Vert x\right\Vert ^{2}}  \notag \\
& \leq \frac{1}{\left\vert \alpha \right\vert }\left[ \left\vert
\left\langle a,b\right\rangle \right\vert +\frac{\left\Vert b\right\Vert }{%
\left\Vert x\right\Vert }\right]  \notag \\
& \qquad \times \left[ \left( \left\vert \alpha -1\right\vert ^{2}\left\vert
\left\langle x,a\right\rangle \right\vert ^{2}+\left\Vert x\right\Vert
^{2}\left\Vert a\right\Vert ^{2}-\left\vert \left\langle x,a\right\rangle
\right\vert ^{2}\right) ^{\frac{1}{2}}\right] ,  \notag
\end{align}%
for each $\alpha \in \mathbb{K}\backslash \left\{ 0\right\} ,$ $a,b,x\in H$
and $x\neq 0.$
\end{remark}

It is obvious that, out of (\ref{ch3.2.1}), we can obtain various particular
inequalities. We mention in the following a class of these which is related
to Buzano's result (\ref{ch3.1.1}) \cite{DRAG1}.

\begin{corollary}[Dragomir, 2004]
\label{ch3.c2.3}Let $a,b,x\in H,$ $x\neq 0$ and $\eta \in \mathbb{K}$ with $%
\left\vert \eta \right\vert =1,$ $\func{Re}\eta \neq -1.$ Then we have the
inequality:%
\begin{equation}
\left\vert \frac{\left\langle a,x\right\rangle \left\langle x,b\right\rangle 
}{\left\Vert x\right\Vert ^{2}}-\frac{\left\langle a,b\right\rangle }{1+\eta 
}\right\vert \leq \frac{\left\Vert a\right\Vert \left\Vert b\right\Vert }{%
\sqrt{2}\sqrt{1+\func{Re}\eta }},  \label{ch3.2.5}
\end{equation}%
and, in particular, for $\eta =1,$ the inequality:%
\begin{equation}
\left\vert \frac{\left\langle a,x\right\rangle \left\langle x,b\right\rangle 
}{\left\Vert x\right\Vert ^{2}}-\frac{\left\langle a,b\right\rangle }{2}%
\right\vert \leq \frac{\left\Vert a\right\Vert \left\Vert b\right\Vert }{2}.
\label{ch3.2.6}
\end{equation}
\end{corollary}

\begin{proof}
It follows by Theorem \ref{ch3.t2.1} on choosing $\alpha =1+\eta $ and we
omit the details.
\end{proof}

\begin{remark}
Using the continuity property of modulus, we get from (\ref{ch3.2.5}) that:%
\begin{equation*}
\frac{\left\vert \left\langle a,x\right\rangle \left\langle x,b\right\rangle
\right\vert }{\left\Vert x\right\Vert ^{2}}\leq \frac{\left\vert
\left\langle a,b\right\rangle \right\vert +\left\Vert a\right\Vert
\left\Vert b\right\Vert }{\sqrt{2}\sqrt{1+\func{Re}\eta }},\qquad \left\vert
\eta \right\vert =1,\ \ \func{Re}\eta \neq -1,
\end{equation*}%
which provides, as the best possible inequality, the above result due to
Buzano (\ref{ch3.1.1}).
\end{remark}

\begin{remark}
\label{ch3.r2.3}If the space is real, then the inequality (\ref{ch3.2.1}) is
obviously equivalent to:%
\begin{align}
& \frac{\left\langle a,b\right\rangle }{\alpha }-\frac{\left\Vert
b\right\Vert }{\left\vert \alpha \right\vert \left\Vert x\right\Vert }\left[
\left( \alpha -1\right) ^{2}\left\langle a,x\right\rangle ^{2}+\left\Vert
x\right\Vert ^{2}\left\Vert a\right\Vert ^{2}-\left\langle a,x\right\rangle
^{2}\right] ^{\frac{1}{2}}  \label{ch3.2.7} \\
& \leq \frac{\left\langle a,x\right\rangle \left\langle x,b\right\rangle }{%
\left\Vert x\right\Vert ^{2}}  \notag \\
& \leq \frac{\left\langle a,b\right\rangle }{\alpha }+\frac{\left\Vert
b\right\Vert }{\left\vert \alpha \right\vert \left\Vert x\right\Vert }\left[
\left( \alpha -1\right) ^{2}\left\langle a,x\right\rangle ^{2}+\left\Vert
x\right\Vert ^{2}\left\Vert a\right\Vert ^{2}-\left\langle a,x\right\rangle
^{2}\right] ^{\frac{1}{2}}  \notag
\end{align}%
for any $\alpha \in \mathbb{R}\backslash \left\{ 0\right\} $ and $a,b,x\in
H, $ $x\neq 0.$

If in (\ref{ch3.2.7}) we take $\alpha =2,$ then we get%
\begin{align}
\frac{1}{2}\left[ \left\langle a,b\right\rangle -\left\Vert a\right\Vert
\left\Vert b\right\Vert \right] \left\Vert x\right\Vert ^{2}& \leq
\left\langle a,x\right\rangle \left\langle x,b\right\rangle  \label{ch3.2.8}
\\
& \leq \frac{1}{2}\left[ \left\langle a,b\right\rangle +\left\Vert
a\right\Vert \left\Vert y\right\Vert \right] \left\Vert x\right\Vert ^{2}, 
\notag
\end{align}%
which apparently, as mentioned by T. Precupanu in \cite{P}, has been
obtained independently of Buzano, by U. Richard in \cite{R}.

In \cite{PE}, Pe\v{c}ari\'{c} gave a simple direct proof of (\ref{ch3.2.8})
without mentioning the work of either Buzano or Richard, but tracked down
the result, in a particular form, to an earlier paper due to C. Blatter \cite%
{BL}.
\end{remark}

Obviously, the following refinement of Buzano's result may be stated \cite%
{DRAG1}.

\begin{corollary}[Dragomir, 2004]
\label{ch3.c2.4}Let $\left( H;\left\langle \cdot ,\cdot \right\rangle
\right) $ be a real or complex inner product space and $a,b,x\in H.$ Then%
\begin{align}
\left\vert \left\langle a,x\right\rangle \left\langle x,b\right\rangle
\right\vert & \leq \left\vert \left\langle a,x\right\rangle \left\langle
x,b\right\rangle -\frac{1}{2}\left\langle a,b\right\rangle \left\Vert
x\right\Vert ^{2}\right\vert +\frac{1}{2}\left\vert \left\langle
a,b\right\rangle \right\vert \left\Vert x\right\Vert ^{2}  \label{ch3.2.9} \\
& \leq \frac{1}{2}\left[ \left\Vert a\right\Vert \left\Vert b\right\Vert
+\left\vert \left\langle a,b\right\rangle \right\vert \right] \left\Vert
x\right\Vert ^{2}.  \notag
\end{align}
\end{corollary}

\begin{proof}
The first inequality in (\ref{ch3.2.9}) follows by the triangle inequality
for the modulus $\left\vert \cdot \right\vert .$ The second inequality is
merely (\ref{ch3.2.6}) in which we added the same quantity to both sides.
\end{proof}

\begin{remark}
\label{ch3.r2.5}For $\alpha =1,$ we deduce from (\ref{ch3.2.1}) the
following inequality:%
\begin{equation}
\left\vert \frac{\left\langle a,x\right\rangle \left\langle x,b\right\rangle 
}{\left\Vert x\right\Vert ^{2}}-\left\langle a,b\right\rangle \right\vert
\leq \frac{\left\Vert b\right\Vert }{\left\Vert x\right\Vert }\left[
\left\Vert x\right\Vert ^{2}\left\Vert a\right\Vert ^{2}-\left\vert
\left\langle a,x\right\rangle \right\vert ^{2}\right] ^{\frac{1}{2}}
\label{ch3.2.10}
\end{equation}%
for any $a,b,x\in H$ with $x\neq 0.$

If the space is real, then (\ref{ch3.2.10}) is equivalent to%
\begin{align}
& \left\langle a,b\right\rangle -\frac{\left\Vert b\right\Vert }{\left\Vert
x\right\Vert }\left[ \left\Vert x\right\Vert ^{2}\left\Vert a\right\Vert
^{2}-\left\vert \left\langle a,x\right\rangle \right\vert ^{2}\right] ^{%
\frac{1}{2}}  \label{ch3.2.11} \\
& \leq \frac{\left\langle a,x\right\rangle \left\langle x,b\right\rangle }{%
\left\Vert x\right\Vert ^{2}}  \notag \\
& \leq \frac{\left\Vert b\right\Vert }{\left\Vert x\right\Vert }\left[
\left\Vert x\right\Vert ^{2}\left\Vert a\right\Vert ^{2}-\left\vert
\left\langle a,x\right\rangle \right\vert ^{2}\right] ^{\frac{1}{2}%
}+\left\langle a,b\right\rangle ,  \notag
\end{align}%
which is similar to Richard's inequality (\ref{ch3.2.8}).
\end{remark}

\subsection{Applications to Kurepa's Inequality}

In 1960, N.G. de Bruijn \cite{BR} obtained the following refinement of the
Cauchy-Bunyakovsky-Schwarz inequality:%
\begin{equation}
\left\vert \sum_{i=1}^{n}a_{i}z_{i}\right\vert ^{2}\leq \frac{1}{2}%
\sum_{i=1}^{n}a_{i}^{2}\left[ \sum_{i=1}^{n}\left\vert z_{i}\right\vert
^{2}+\left\vert \sum_{i=1}^{n}z_{i}^{2}\right\vert \right] ,  \label{ch3.3.1}
\end{equation}%
provided that $a_{i}$ are real numbers while $z_{i}$ are complex for each $%
i\in \left\{ 1,...,n\right\} .$

In \cite{KU}, S. Kurepa proved the following generalisation of the de Bruijn
result:

\begin{theorem}[Kurepa, 1966]
\label{ch3.t3.1}Let $\left( H;\left\langle \cdot ,\cdot \right\rangle
\right) $ be a real inner product space and $\left( H_{\mathbb{C}%
},\left\langle \cdot ,\cdot \right\rangle _{\mathbb{C}}\right) $ its
complexification. Then for any $a\in H$ and $z\in H_{\mathbb{C}},$ one has
the following refinement of Schwarz's inequality%
\begin{equation}
\left\vert \left\langle a,z\right\rangle _{\mathbb{C}}\right\vert ^{2}\leq 
\frac{1}{2}\left\Vert a\right\Vert ^{2}\left[ \left\Vert z\right\Vert _{%
\mathbb{C}}^{2}+\left\vert \left\langle z,\bar{z}\right\rangle _{\mathbb{C}%
}\right\vert \right] \leq \left\Vert a\right\Vert ^{2}\left\Vert
z\right\Vert _{\mathbb{C}}^{2},  \label{ch3.3.2}
\end{equation}%
where $\bar{z}$ denotes the conjugate of $z\in H_{\mathbb{C}}.$
\end{theorem}

As consequences of this general result, Kurepa noted the following integral,
respectively, discrete inequality:

\begin{corollary}[Kurepa, 1966]
\label{ch3.c3.2}Let $\left( S,\Sigma ,\mu \right) $ be a positive measure
space and $a,z\in L_{2}\left( S,\Sigma ,\mu \right) ,$ the Hilbert space of
complex-valued $2-\mu -$integrable functions defined on $S.$ If $a$ is a
real function and $z$ is a complex function, then 
\begin{multline}
\left\vert \int_{S}a\left( t\right) z\left( t\right) d\mu \left( t\right)
\right\vert ^{2}  \label{ch3.3.3} \\
\leq \frac{1}{2}\cdot \int_{S}a^{2}\left( t\right) d\mu \left( t\right) %
\left[ \int_{S}\left\vert z\left( t\right) \right\vert ^{2}d\mu \left(
t\right) +\left\vert \int_{S}z^{2}\left( t\right) d\mu \left( t\right)
\right\vert \right] .
\end{multline}
\end{corollary}

\begin{corollary}[Kurepa, 1966]
\label{ch3.c3.3}If $a_{1},\dots ,a_{n}$ are real numbers, $z_{1},\dots
,z_{n} $ are complex numbers and $\left( A_{ij}\right) $ is a positive
definite real matrix of dimension $n\times n$, then%
\begin{equation}
\left\vert \sum_{i,j=1}^{n}A_{ij}a_{i}z_{j}\right\vert ^{2}\leq \frac{1}{2}%
\sum_{i,j=1}^{n}A_{ij}a_{i}a_{j}\left[ \sum_{i,j=1}^{n}A_{ij}z_{i}\overline{%
z_{j}}+\left\vert \sum_{i,j=1}^{n}A_{ij}z_{i}\overline{z_{j}}\right\vert %
\right] .  \label{ch3.3.4}
\end{equation}
\end{corollary}

The following refinement of Kurepa's result may be stated \cite{DRAG1}.

\begin{theorem}[Dragomir, 2004]
\label{ch3.t3.3}Let $\left( H;\left\langle \cdot ,\cdot \right\rangle
\right) $ be a real inner product space and $\left( H_{\mathbb{C}%
},\left\langle \cdot ,\cdot \right\rangle _{\mathbb{C}}\right) $ its
complexification. Then for any $e\in H$ and $w\in H_{\mathbb{C}},$ one has
the inequality:%
\begin{align}
\left\vert \left\langle w,e\right\rangle _{\mathbb{C}}\right\vert ^{2}& \leq
\left\vert \left\langle w,e\right\rangle _{\mathbb{C}}^{2}-\frac{1}{2}%
\left\langle w,\bar{w}\right\rangle _{\mathbb{C}}\left\Vert e\right\Vert
^{2}\right\vert +\frac{1}{2}\left\vert \left\langle w,\bar{w}\right\rangle _{%
\mathbb{C}}\right\vert \left\Vert e\right\Vert ^{2}  \label{ch3.3.5} \\
& \leq \frac{1}{2}\left\Vert e\right\Vert ^{2}\left[ \left\Vert w\right\Vert
_{\mathbb{C}}^{2}+\left\vert \left\langle w,\bar{w}\right\rangle _{\mathbb{C}%
}\right\vert \right] .  \notag
\end{align}
\end{theorem}

\begin{proof}
We follow the proof in \cite{DRAG1}.

If we apply Corollary \ref{ch3.c3.3} for $\left( H_{\mathbb{C}},\left\langle
\cdot ,\cdot \right\rangle _{\mathbb{C}}\right) $ and $x=e\in H,$ $a=w$ and $%
b=\bar{w},$ then we have%
\begin{align}
& \left\vert \left\langle w,e\right\rangle _{\mathbb{C}}\left\langle e,\bar{w%
}\right\rangle _{\mathbb{C}}\right\vert  \label{ch3.3.6} \\
& \leq \left\vert \left\langle w,e\right\rangle _{\mathbb{C}}\left\langle e,%
\bar{w}\right\rangle _{\mathbb{C}}-\frac{1}{2}\left\langle w,\bar{w}%
\right\rangle _{\mathbb{C}}\left\Vert e\right\Vert ^{2}\right\vert +\frac{1}{%
2}\left\vert \left\langle w,\bar{w}\right\rangle _{\mathbb{C}}\right\vert
\left\Vert e\right\Vert ^{2}  \notag \\
& \leq \frac{1}{2}\left\Vert e\right\Vert ^{2}\left[ \left\Vert w\right\Vert
_{\mathbb{C}}\left\Vert \bar{w}\right\Vert _{\mathbb{C}}+\left\vert
\left\langle w,\bar{w}\right\rangle _{\mathbb{C}}\right\vert \right] . 
\notag
\end{align}%
Now, if we assume that $w=\left( x,y\right) \in H_{\mathbb{C}},$ then, by
the definition of $\left\langle \cdot ,\cdot \right\rangle _{\mathbb{C}},$
we have%
\begin{align*}
\left\langle w,e\right\rangle _{\mathbb{C}}& =\left\langle \left( x,y\right)
,\left( e,0\right) \right\rangle _{\mathbb{C}} \\
& =\left\langle x,e\right\rangle +\left\langle y,0\right\rangle +i\left[
\left\langle y,e\right\rangle -\left\langle x,0\right\rangle \right] \\
& =\left\langle e,x\right\rangle +i\left\langle e,y\right\rangle ,
\end{align*}%
\begin{align*}
\left\langle e,\bar{w}\right\rangle _{\mathbb{C}}& =\left\langle \left(
e,0\right) ,\left( x,-y\right) \right\rangle _{\mathbb{C}} \\
& =\left\langle e,x\right\rangle +\left\langle 0,-y\right\rangle +i\left[
\left\langle 0,x\right\rangle -\left\langle e,-y\right\rangle \right] \\
& =\left\langle e,x\right\rangle +i\left\langle e,y\right\rangle
=\left\langle w,e\right\rangle _{\mathbb{C}}
\end{align*}%
and%
\begin{equation*}
\left\Vert \bar{w}\right\Vert _{\mathbb{C}}^{2}=\left\Vert x\right\Vert
^{2}+\left\Vert y\right\Vert ^{2}=\left\Vert w\right\Vert _{\mathbb{C}}^{2}.
\end{equation*}%
Therefore, by (\ref{ch3.3.6}), we deduce the desired result (\ref{ch3.3.5}).
\end{proof}

Denote by $\ell _{\rho }^{2}\left( \mathbb{C}\right) $ the Hilbert space of
all complex sequences $z=\left( z_{i}\right) _{i\in \mathbb{N}}$ with the
property that for $\rho _{i}\geq 0$ with $\sum_{i=1}^{\infty }\rho _{i}=1$
we have $\sum_{i=1}^{\infty }\rho _{i}\left\vert z_{i}\right\vert
^{2}<\infty .$ If $a=\left( a_{i}\right) _{i\in \mathbb{N}}$ is a sequence
of real numbers such that $a\in \ell _{\rho }^{2}\left( \mathbb{C}\right) ,$
then for any $z\in \ell _{\rho }^{2}\left( \mathbb{C}\right) $ we have the
inequality:%
\begin{align}
& \left\vert \sum_{i=1}^{\infty }\rho _{i}a_{i}z_{i}\right\vert ^{2}
\label{ch3.3.7} \\
& \leq \left\vert \left( \sum_{i=1}^{\infty }\rho _{i}a_{i}z_{i}\right) ^{2}-%
\frac{1}{2}\sum_{i=1}^{\infty }\rho _{i}a_{i}^{2}\sum_{i=1}^{\infty }\rho
_{i}z_{i}^{2}\right\vert +\frac{1}{2}\sum_{i=1}^{\infty }\rho
_{i}a_{i}^{2}\left\vert \sum_{i=1}^{\infty }\rho _{i}z_{i}^{2}\right\vert 
\notag \\
& \leq \frac{1}{2}\sum_{i=1}^{\infty }\rho _{i}a_{i}^{2}\left[
\sum_{i=1}^{\infty }\rho _{i}\left\vert z_{i}\right\vert ^{2}+\left\vert
\sum_{i=1}^{\infty }\rho _{i}z_{i}^{2}\right\vert \right] .  \notag
\end{align}

Similarly, if by $L_{\rho }^{2}\left( S,\Sigma ,\mu \right) $ we understand
the Hilbert space of all complex-valued functions $f:S\rightarrow \mathbb{C}$
with the property that for the $\mu -$measurable function $\rho \geq 0$ with 
$\int_{S}\rho \left( t\right) d\mu \left( t\right) =1$ we have%
\begin{equation*}
\int_{S}\rho \left( t\right) \left\vert f\left( t\right) \right\vert
^{2}d\mu \left( t\right) <\infty ,
\end{equation*}%
then for a real function $a\in L_{\rho }^{2}\left( S,\Sigma ,\mu \right) $
and any $f\in L_{\rho }^{2}\left( S,\Sigma ,\mu \right) ,$ we have the
inequalities%
\begin{align}
& \left\vert \int_{S}\rho \left( t\right) a\left( t\right) f\left( t\right)
d\mu \left( t\right) \right\vert ^{2}  \label{ch3.3.8} \\
& \leq \left\vert \left( \int_{S}\rho \left( t\right) a\left( t\right)
f\left( t\right) d\mu \left( t\right) \right) ^{2}\right.  \notag \\
& \qquad \qquad -\frac{1}{2}\left. \int_{S}\rho \left( t\right) f^{2}\left(
t\right) d\mu \left( t\right) \int_{S}\rho \left( t\right) a^{2}\left(
t\right) d\mu \left( t\right) \right\vert  \notag \\
& \qquad \qquad +\frac{1}{2}\left\vert \int_{S}\rho \left( t\right)
f^{2}\left( t\right) d\mu \left( t\right) \right\vert \int_{S}\rho \left(
t\right) a^{2}\left( t\right) d\mu \left( t\right)  \notag \\
& \leq \frac{1}{2}\int_{S}\rho \left( t\right) a^{2}\left( t\right) d\mu
\left( t\right)  \notag \\
& \qquad \qquad \times \left[ \int_{S}\rho \left( t\right) \left\vert
f\left( t\right) \right\vert ^{2}d\mu \left( t\right) +\left\vert
\int_{S}\rho \left( t\right) f^{2}\left( t\right) d\mu \left( t\right)
\right\vert \right] .  \notag
\end{align}

\section{Inequalities for Orthornormal Families}

\subsection{\label{ch4.s1}Introduction}

In \cite{B}, M.L. Buzano obtained the following extension of the celebrated
Schwarz's inequality in a real or complex inner product space $\left(
H;\left\langle \cdot ,\cdot \right\rangle \right) :$%
\begin{equation}
\left\vert \left\langle a,x\right\rangle \left\langle x,b\right\rangle
\right\vert \leq \frac{1}{2}\left[ \left\Vert a\right\Vert \left\Vert
b\right\Vert +\left\vert \left\langle a,b\right\rangle \right\vert \right]
\left\Vert x\right\Vert ^{2},  \label{ch4.1.1}
\end{equation}%
for any $a,b,x\in H.$

It is clear that the above inequality becomes, for $a=b,$ the Schwarz's
inequality%
\begin{equation}
\left\vert \left\langle a,x\right\rangle \right\vert ^{2}\leq \left\Vert
a\right\Vert ^{2}\left\Vert x\right\Vert ^{2},\quad a,x\in H;
\label{ch4.1.2}
\end{equation}%
in which the equality holds if and only if there exists a scalar $\lambda
\in \mathbb{K}$ $\left( \mathbb{R},\mathbb{C}\right) $ so that $x=\lambda a.$

As noted by T. Precupanu in \cite{P}, independently of Buzano, U. Richard 
\cite{R} obtained the following similar inequality holding in real inner
product spaces:%
\begin{align}
\frac{1}{2}\left\Vert x\right\Vert ^{2}\left[ \left\langle a,b\right\rangle
-\left\Vert a\right\Vert \left\Vert b\right\Vert \right] & \leq \left\langle
a,x\right\rangle \left\langle x,b\right\rangle  \label{ch4.1.4} \\
& \leq \frac{1}{2}\left\Vert x\right\Vert ^{2}\left[ \left\langle
a,b\right\rangle +\left\Vert a\right\Vert \left\Vert b\right\Vert \right] . 
\notag
\end{align}

The main aim of the present section is to obtain similar results for
families of orthonormal vectors in $\left( H;\left\langle \cdot ,\cdot
\right\rangle \right) ,$ real or complex space, that are naturally connected
with the celebrated Bessel inequality and improve the results of Busano,
Richard and Kurepa.

\subsection{A Generalisation for Orthonormal Families\label{ch4.s2}}

We say that the finite family $\left\{ e_{i}\right\} _{i\in I}$ \ ($I$ is
finite) of vectors is \textit{orthonormal} if $\left\langle
e_{i},e_{j}\right\rangle =0$ if $i,j\in I$ with $i\neq j$ and $\left\Vert
e_{i}\right\Vert =1$ for each $i\in I.$ The following result may be stated 
\cite{DRAG2}:

\begin{theorem}[Dragomir, 2004]
\label{ch4.t2.1}Let $\left( H;\left\langle \cdot ,\cdot \right\rangle
\right) $ be an inner product space over the real or complex number field $%
\mathbb{K}$ and $\left\{ e_{i}\right\} _{i\in I}$ a finite orthonormal
family in $H.$ Then for any $a,b\in H,$ one has the inequality:%
\begin{equation}
\left\vert \sum_{i\in I}\left\langle a,e_{i}\right\rangle \left\langle
e_{i},b\right\rangle -\frac{1}{2}\left\langle a,b\right\rangle \right\vert
\leq \frac{1}{2}\left\Vert a\right\Vert \left\Vert b\right\Vert .
\label{ch4.2.1}
\end{equation}%
The case of equality holds in (\ref{ch4.2.1}) if and only if%
\begin{equation}
\sum_{i\in I}\left\langle a,e_{i}\right\rangle e_{i}=\frac{1}{2}a+\left(
\sum_{i\in I}\left\langle a,e_{i}\right\rangle \left\langle
e_{i},b\right\rangle -\frac{1}{2}\left\langle a,b\right\rangle \right) \cdot 
\frac{b}{\left\Vert b\right\Vert ^{2}}.  \label{ch4.2.2}
\end{equation}
\end{theorem}

\begin{proof}
We follow the proof in \cite{DRAG2}.

It is well known that, for $e\neq 0$ and $f\in H,$ the following identity
holds:%
\begin{equation}
\frac{\left\Vert f\right\Vert ^{2}\left\Vert e\right\Vert ^{2}-\left\vert
\left\langle f,e\right\rangle \right\vert ^{2}}{\left\Vert e\right\Vert ^{2}}%
=\left\Vert f-\frac{\left\langle f,e\right\rangle e}{\left\Vert e\right\Vert
^{2}}\right\Vert ^{2}.  \label{ch4.2.3}
\end{equation}%
Therefore, in Schwarz's inequality%
\begin{equation}
\left\vert \left\langle f,e\right\rangle \right\vert ^{2}\leq \left\Vert
f\right\Vert ^{2}\left\Vert e\right\Vert ^{2},\quad f,e\in H;
\label{ch4.2.4}
\end{equation}%
the case of equality, for $e\neq 0,$ holds if and only if%
\begin{equation*}
f=\frac{\left\langle f,e\right\rangle e}{\left\Vert e\right\Vert ^{2}}.
\end{equation*}%
Let $f:=2\sum_{i\in I}\left\langle a,e_{i}\right\rangle e_{i}-a$ and $e:=b.$
Then, by Schwarz's inequality (\ref{ch4.2.4}), we may state that%
\begin{equation}
\left\vert \left\langle 2\sum_{i\in I}\left\langle a,e_{i}\right\rangle
e_{i}-a,b\right\rangle \right\vert ^{2}\leq \left\Vert 2\sum_{i\in
I}\left\langle a,e_{i}\right\rangle e_{i}-a\right\Vert ^{2}\left\Vert
b\right\Vert ^{2}  \label{ch4.2.5}
\end{equation}%
with equality, for $b\neq 0,$ if and only if%
\begin{equation}
2\sum_{i\in I}\left\langle a,e_{i}\right\rangle e_{i}-a=\left\langle
2\sum_{i\in I}\left\langle a,e_{i}\right\rangle e_{i}-a,b\right\rangle \frac{%
b}{\left\Vert b\right\Vert ^{2}}.  \label{ch4.2.6}
\end{equation}%
Since%
\begin{equation*}
\left\langle 2\sum_{i\in I}\left\langle a,e_{i}\right\rangle
e_{i}-a,b\right\rangle =2\sum_{i\in I}\left\langle a,e_{i}\right\rangle
\left\langle e_{i},b\right\rangle -\left\langle a,b\right\rangle
\end{equation*}%
and%
\begin{align*}
& \left\Vert 2\sum_{i\in I}\left\langle a,e_{i}\right\rangle
e_{i}-a\right\Vert ^{2} \\
& =4\left\Vert \sum_{i\in I}\left\langle a,e_{i}\right\rangle
e_{i}\right\Vert ^{2}-4\func{Re}\left\langle \sum_{i\in I}\left\langle
a,e_{i}\right\rangle e_{i},a\right\rangle +\left\Vert a\right\Vert ^{2} \\
& =4\sum_{i\in I}\left\vert \left\langle a,e_{i}\right\rangle \right\vert
^{2}-4\sum_{i\in I}\left\vert \left\langle a,e_{i}\right\rangle \right\vert
^{2}+\left\Vert a\right\Vert ^{2} \\
& =\left\Vert a\right\Vert ^{2},
\end{align*}%
hence by (\ref{ch4.2.5}) we deduce the desired inequality (\ref{ch4.2.1}).

Finally, as (\ref{ch4.2.2}) is equivalent to%
\begin{equation*}
\sum_{i\in I}\left\langle a,e_{i}\right\rangle e_{i}-\frac{a}{2}=\left(
\sum_{i\in I}\left\langle a,e_{i}\right\rangle \left\langle
e_{i},b\right\rangle -\frac{1}{2}\left\langle a,b\right\rangle \right) \frac{%
b}{\left\Vert b\right\Vert ^{2}},
\end{equation*}%
hence the equality holds in (\ref{ch4.2.1}) if and only if (\ref{ch4.2.2})
is valid.
\end{proof}

The following result is well known in the literature as Bessel's inequality%
\begin{equation}
\sum_{i\in I}\left\vert \left\langle x,e_{i}\right\rangle \right\vert
^{2}\leq \left\Vert x\right\Vert ^{2},\quad x\in H,  \label{ch4.2.7}
\end{equation}%
where, as above, $\left\{ e_{i}\right\} _{i\in I}$ is a finite orthonormal
family in the inner product space $\left( H;\left\langle \cdot ,\cdot
\right\rangle \right) .$

If one chooses $a=b=x$ in (\ref{ch4.2.1}), then one gets the inequality%
\begin{equation*}
\left\vert \sum_{i\in I}\left\vert \left\langle x,e_{i}\right\rangle
\right\vert ^{2}-\frac{1}{2}\left\Vert x\right\Vert ^{2}\right\vert \leq 
\frac{1}{2}\left\Vert x\right\Vert ^{2},
\end{equation*}%
which is obviously equivalent to Bessel's inequality (\ref{ch4.2.7}).
Therefore, the inequality (\ref{ch4.2.1}) may be regarded as a
generalisation of Bessel's inequality as well.

Utilising the Bessel and Cauchy-Bunyakovsky-Schwarz inequalities, one may
state that%
\begin{equation}
\left\vert \sum_{i\in I}\left\langle a,e_{i}\right\rangle \left\langle
e_{i},b\right\rangle \right\vert \leq \left[ \sum_{i\in I}\left\vert
\left\langle a,e_{i}\right\rangle \right\vert ^{2}\sum_{i\in I}\left\vert
\left\langle b,e_{i}\right\rangle \right\vert ^{2}\right] ^{\frac{1}{2}}\leq
\left\Vert a\right\Vert \left\Vert b\right\Vert  \label{ch4.2.8}
\end{equation}

A different refinement of the inequality between the first and the last term
in (\ref{ch4.2.8}) is incorporated in the following \cite{DRAG2}:

\begin{corollary}[Dragomir, 2004]
\label{ch4.c2.2}With the assumption of Theorem \ref{ch4.t2.1}, we have%
\begin{align}
\left\vert \sum_{i\in I}\left\langle a,e_{i}\right\rangle \left\langle
e_{i},b\right\rangle \right\vert & \leq \left\vert \sum_{i\in I}\left\langle
a,e_{i}\right\rangle \left\langle e_{i},b\right\rangle -\frac{1}{2}%
\left\langle a,b\right\rangle \right\vert +\frac{1}{2}\left\vert
\left\langle a,b\right\rangle \right\vert  \label{ch4.2.9} \\
& \leq \frac{1}{2}\left[ \left\Vert a\right\Vert \left\Vert b\right\Vert
+\left\vert \left\langle a,b\right\rangle \right\vert \right]  \notag \\
& \leq \left\Vert a\right\Vert \left\Vert b\right\Vert .  \notag
\end{align}
\end{corollary}

\begin{remark}
\label{ch4.r2.3}If the space $\left( H;\left\langle \cdot ,\cdot
\right\rangle \right) $ is real, then, obviously, (\ref{ch4.2.1}) is
equivalent to:%
\begin{equation}
\frac{1}{2}\left( \left\langle a,b\right\rangle -\left\Vert a\right\Vert
\left\Vert b\right\Vert \right) \leq \sum_{i\in I}\left\langle
a,e_{i}\right\rangle \left\langle e_{i},b\right\rangle \leq \frac{1}{2}\left[
\left\Vert a\right\Vert \left\Vert b\right\Vert +\left\langle
a,b\right\rangle \right] .  \label{ch4.2.10}
\end{equation}
\end{remark}

\begin{remark}
\label{ch4.r2.4}It is obvious that if the family comprises of only a single
element $e=\frac{x}{\left\Vert x\right\Vert },$ $x\in H,$ $x\neq 0,$ then
from (\ref{ch4.2.9}) we recapture the refinement of Buzano's inequality
incorporated in (\ref{ch4.1.1}) while from (\ref{ch4.2.10}) we deduce
Richard's result from (\ref{ch4.1.4}).
\end{remark}

The following corollary of Theorem \ref{ch4.t2.1} is of interest as well 
\cite{DRAG2}:

\begin{corollary}[Dragomir, 2004]
\label{ch4.c2.3}Let $\left\{ e_{i}\right\} _{i\in I}$ be a finite
orthonormal family in $\left( H;\left\langle \cdot ,\cdot \right\rangle
\right) .$ If $x,y\in H\backslash \left\{ 0\right\} $ are such that there
exists the constants $m_{i},n_{i},$ $M_{i},$ $N_{i}\in \mathbb{R}$, $i\in I$
such that:%
\begin{equation}
-1\leq m_{i}\leq \frac{\func{Re}\left\langle x,e_{i}\right\rangle }{%
\left\Vert x\right\Vert }\cdot \frac{\func{Re}\left\langle
y,e_{i}\right\rangle }{\left\Vert y\right\Vert }\leq M_{i}\leq 1,\quad i\in I
\label{ch4.2.11}
\end{equation}%
and%
\begin{equation}
-1\leq n_{i}\leq \frac{\func{Im}\left\langle x,e_{i}\right\rangle }{%
\left\Vert x\right\Vert }\cdot \frac{\func{Im}\left\langle
y,e_{i}\right\rangle }{\left\Vert y\right\Vert }\leq N_{i}\leq 1,\quad i\in I
\label{ch4.2.12}
\end{equation}%
then%
\begin{equation}
2\sum_{i\in I}\left( m_{i}+n_{i}\right) -1\leq \frac{\func{Re}\left\langle
x,y\right\rangle }{\left\Vert x\right\Vert \left\Vert y\right\Vert }\leq
1+2\sum_{i\in I}\left( M_{i}+N_{i}\right) .  \label{ch4.2.13}
\end{equation}
\end{corollary}

\begin{proof}
We follow the proof in \cite{DRAG2}.

Using Theorem \ref{ch4.t2.1} and the fact that for any complex number $z,$ $%
\left\vert z\right\vert \geq \left\vert \func{Re}z\right\vert ,$ we have%
\begin{align}
& \left\vert \sum_{i\in I}\func{Re}\left[ \left\langle x,e_{i}\right\rangle
\left\langle e_{i},y\right\rangle \right] -\frac{1}{2}\func{Re}\left\langle
x,y\right\rangle \right\vert  \label{ch4.2.14} \\
& \leq \left\vert \sum_{i\in I}\left\langle x,e_{i}\right\rangle
\left\langle e_{i},y\right\rangle -\frac{1}{2}\left\langle x,y\right\rangle
\right\vert  \notag \\
& \leq \frac{1}{2}\left\Vert x\right\Vert \left\Vert y\right\Vert .  \notag
\end{align}%
Since%
\begin{equation*}
\func{Re}\left[ \left\langle x,e_{i}\right\rangle \left\langle
e_{i},y\right\rangle \right] =\func{Re}\left\langle x,e_{i}\right\rangle 
\func{Re}\left\langle y,e_{i}\right\rangle +\func{Im}\left\langle
x,e_{i}\right\rangle \func{Im}\left\langle y,e_{i}\right\rangle ,
\end{equation*}%
hence by (\ref{ch4.2.14}) we have:%
\begin{align}
& -\frac{1}{2}\left\Vert x\right\Vert \left\Vert y\right\Vert +\frac{1}{2}%
\func{Re}\left\langle x,y\right\rangle  \label{ch4.2.15} \\
& \leq \sum_{i\in I}\func{Re}\left\langle x,e_{i}\right\rangle \func{Re}%
\left\langle y,e_{i}\right\rangle +\sum_{i\in I}\func{Im}\left\langle
x,e_{i}\right\rangle \func{Im}\left\langle y,e_{i}\right\rangle  \notag \\
& \leq \frac{1}{2}\left\Vert x\right\Vert \left\Vert y\right\Vert +\frac{1}{2%
}\func{Re}\left\langle x,y\right\rangle .  \notag
\end{align}%
Utilising the assumptions (\ref{ch4.2.11}) and (\ref{ch4.2.12}), we have%
\begin{equation}
\sum_{i\in I}m_{i}\leq \sum_{i\in I}\frac{\func{Re}\left\langle
x,e_{i}\right\rangle \func{Re}\left\langle y,e_{i}\right\rangle }{\left\Vert
x\right\Vert \left\Vert y\right\Vert }\leq \sum_{i\in I}M_{i}
\label{ch4.2.16}
\end{equation}%
and 
\begin{equation}
\sum_{i\in I}n_{i}\leq \sum_{i\in I}\frac{\func{Im}\left\langle
x,e_{i}\right\rangle \func{Im}\left\langle y,e_{i}\right\rangle }{\left\Vert
x\right\Vert \left\Vert y\right\Vert }\leq \sum_{i\in I}N_{i}.
\label{ch4.2.17}
\end{equation}%
Finally, on making use of (\ref{ch4.2.15}) -- (\ref{ch4.2.17}), we deduce
the desired result (\ref{ch4.2.13}).
\end{proof}

\begin{remark}
By Schwarz's inequality, is it obvious that, in general,%
\begin{equation*}
-1\leq \frac{\func{Re}\left\langle x,y\right\rangle }{\left\Vert
x\right\Vert \left\Vert y\right\Vert }\leq 1.
\end{equation*}%
Consequently, the left inequality in (\ref{ch4.2.13}) is of interest when $%
\sum_{i\in I}\left( m_{i}+n_{i}\right) >0,$ while the right inequality in (%
\ref{ch4.2.13}) is of interest when $\sum_{i\in I}\left( M_{i}+N_{i}\right)
<0.$
\end{remark}

\subsection{Refinements of Kurepa's Inequality\label{ch4.s3}}

The following result holds \cite{DRAG2}.

\begin{theorem}[Dragomir, 2004]
\label{ch4.t3.1}Let $\left\{ e_{j}\right\} _{j\in I}$ be a finite
orthonormal family in the real inner product space $\left( H;\left\langle
\cdot ,\cdot \right\rangle \right) .$ Then for any $w\in H_{\mathbb{C}},$
where $\left( H_{\mathbb{C}};\left\langle \cdot ,\cdot \right\rangle _{%
\mathbb{C}}\right) $ is the complexification of $\left( H;\left\langle \cdot
,\cdot \right\rangle \right) ,$ one has the following Bessel's type
inequality:%
\begin{align}
\left\vert \sum_{j\in I}\left\langle w,e_{j}\right\rangle _{\mathbb{C}%
}^{2}\right\vert & \leq \left\vert \sum_{j\in I}\left\langle
w,e_{j}\right\rangle _{\mathbb{C}}^{2}-\frac{1}{2}\left\langle w,\bar{w}%
\right\rangle _{\mathbb{C}}\right\vert +\frac{1}{2}\left\vert \left\langle w,%
\bar{w}\right\rangle _{\mathbb{C}}\right\vert  \label{ch4.3.3} \\
& \leq \frac{1}{2}\left[ \left\Vert w\right\Vert _{\mathbb{C}%
}^{2}+\left\vert \left\langle w,\bar{w}\right\rangle _{\mathbb{C}%
}\right\vert \right] \leq \left\Vert w\right\Vert _{\mathbb{C}}^{2}.  \notag
\end{align}
\end{theorem}

\begin{proof}
We follow the proof in \cite{DRAG2}.

Define $f_{j}\in H_{\mathbb{C}},$ $f_{j}:=\left( e_{j},0\right) ,$ $j\in I.$
For any $k,j\in I$ we have%
\begin{equation*}
\left\langle f_{i},f_{j}\right\rangle _{\mathbb{C}}=\left\langle \left(
e_{k},0\right) ,\left( e_{j},0\right) \right\rangle _{\mathbb{C}%
}=\left\langle e_{k},e_{j}\right\rangle =\delta _{kj},
\end{equation*}%
therefore $\left\{ f_{j}\right\} _{j\in I}$ is an orthonormal family in $%
\left( H_{\mathbb{C}};\left\langle \cdot ,\cdot \right\rangle _{\mathbb{C}%
}\right) .$

If we apply Theorem \ref{ch4.t2.1} for $\left( H_{\mathbb{C}};\left\langle
\cdot ,\cdot \right\rangle _{\mathbb{C}}\right) ,$ $a=w,$ $b=\bar{w},$ we
may write:%
\begin{equation}
\left\vert \sum_{j\in I}\left\langle w,e_{j}\right\rangle _{\mathbb{C}%
}\left\langle e_{j},\bar{w}\right\rangle _{\mathbb{C}}-\frac{1}{2}%
\left\langle w,\bar{w}\right\rangle _{\mathbb{C}}\right\vert \leq \frac{1}{2}%
\left\Vert w\right\Vert _{\mathbb{C}}\left\Vert \bar{w}\right\Vert _{\mathbb{%
C}}.  \label{ch4.3.4}
\end{equation}%
However, for $w:=\left( x,y\right) \in H_{\mathbb{C}}$, we have $\bar{w}%
=\left( x,-y\right) $ and%
\begin{equation*}
\left\langle e_{j},\bar{w}\right\rangle _{\mathbb{C}}=\left\langle \left(
e_{j},0\right) ,\left( x,-y\right) \right\rangle _{\mathbb{C}}=\left\langle
e_{j},x\right\rangle -i\left\langle e_{j},-y\right\rangle =\left\langle
e_{j},x\right\rangle +i\left\langle e_{j},y\right\rangle
\end{equation*}%
and%
\begin{equation*}
\left\langle w,e_{j}\right\rangle _{\mathbb{C}}=\left\langle \left(
x,y\right) ,\left( e_{j},0\right) \right\rangle _{\mathbb{C}}=\left\langle
e_{j},x\right\rangle -i\left\langle e_{j},-y\right\rangle =\left\langle
x,e_{j}\right\rangle +i\left\langle e_{j},y\right\rangle
\end{equation*}%
for any $j\in I.$ Thus $\left\langle e_{j},\bar{w}\right\rangle
=\left\langle w,e_{j}\right\rangle $ for each $j\in I$ and since%
\begin{equation*}
\left\Vert w\right\Vert _{\mathbb{C}}=\left\Vert \bar{w}\right\Vert _{%
\mathbb{C}}=\left( \left\Vert x\right\Vert ^{2}+\left\Vert y\right\Vert
^{2}\right) ^{\frac{1}{2}},
\end{equation*}%
we get from (\ref{ch4.3.4}) that%
\begin{equation}
\left\vert \sum_{j\in I}\left\langle w,e_{j}\right\rangle _{\mathbb{C}}^{2}-%
\frac{1}{2}\left\langle w,\bar{w}\right\rangle _{\mathbb{C}}\right\vert \leq 
\frac{1}{2}\left\Vert w\right\Vert _{\mathbb{C}}^{2}.  \label{ch4.3.5}
\end{equation}%
Now, observe that the first inequality in (\ref{ch4.3.3}) follows by the
triangle inequality, the second is an obvious consequence of (\ref{ch4.3.5})
and the last one is derived from Schwarz's result.
\end{proof}

\begin{remark}
\label{ch4.r3.2}If the family $\left\{ e_{j}\right\} _{j\in I}$ contains
only a single element $e=\frac{x}{\left\Vert x\right\Vert },$ $x\in H,$ $%
x\neq 0,$ then from (\ref{ch4.3.3}) we deduce (\ref{ch3.3.6}), which, in its
turn, provides a refinement of Kurepa's inequality (\ref{ch3.3.2}).
\end{remark}

\subsection{\label{ch4.s4}An Application for $L_{2}\left[ -\protect\pi ,%
\protect\pi \right] $}

It is well known that in the Hilbert space $L_{2}\left[ -\pi ,\pi \right] $
of all functions $f:\left[ -\pi ,\pi \right] \rightarrow \mathbb{C}$ with
the property that $f$ is Lebesgue measurable on $\left[ -\pi ,\pi \right] $
and $\int_{-\pi }^{\pi }\left\vert f\left( t\right) \right\vert
^{2}dt<\infty ,$ the set of functions%
\begin{equation*}
\left\{ \frac{1}{\sqrt{2\pi }},\frac{1}{\sqrt{\pi }}\cos t,\frac{1}{\sqrt{%
\pi }}\sin t,\dots ,\frac{1}{\sqrt{\pi }}\cos nt,\frac{1}{\sqrt{\pi }}\sin
nt,\dots \right\}
\end{equation*}%
is orthonormal.

If by $\limfunc{trig}t,$ we denote either $\sin t$ or $\cos t$, $t\in \left[
-\pi ,\pi \right] ,$ then on using the results from Sections \ref{ch4.s2}
and \ref{ch4.s3}, we may state the following inequality:%
\begin{multline}
\left\vert \frac{1}{\pi }\sum_{k=1}^{n}\int_{-\pi }^{\pi }f\left( t\right) 
\limfunc{trig}\left( kt\right) dt\cdot \int_{-\pi }^{\pi }\overline{g\left(
t\right) }\limfunc{trig}\left( kt\right) dt\right.  \label{ch4.4.1} \\
\left. -\frac{1}{2}\int_{-\pi }^{\pi }f\left( t\right) \overline{g\left(
t\right) }dt\right\vert ^{2} \\
\leq \frac{1}{4}\int_{-\pi }^{\pi }\left\vert f\left( t\right) \right\vert
^{2}dt\int_{-\pi }^{\pi }\left\vert g\left( t\right) \right\vert ^{2}dt,
\end{multline}%
where all $\limfunc{trig}\left( kt\right) $ is either $\sin kt$ or $\cos kt,$
$k\in \left\{ 1,\dots ,n\right\} $ and $f\in L_{2}\left[ -\pi ,\pi \right] .$

This follows by Theorem \ref{ch4.t2.1}.

If one uses Corollary \ref{ch4.c2.2}, then one can state the following chain
of inequalities%
\begin{align}
& \left\vert \frac{1}{\pi }\sum_{k=1}^{n}\int_{-\pi }^{\pi }f\left( t\right) 
\limfunc{trig}\left( kt\right) dt\cdot \int_{-\pi }^{\pi }\overline{g\left(
t\right) }\limfunc{trig}\left( kt\right) dt\right\vert  \label{ch4.4.2} \\
& \leq \left\vert \frac{1}{\pi }\sum_{k=1}^{n}\int_{-\pi }^{\pi }f\left(
t\right) \limfunc{trig}\left( kt\right) dt\cdot \int_{-\pi }^{\pi }\overline{%
g\left( t\right) }\limfunc{trig}\left( kt\right) dt\right.  \notag \\
& \qquad \qquad \left. -\frac{1}{2}\int_{-\pi }^{\pi }f\left( t\right) 
\overline{g\left( t\right) }dt\right\vert +\frac{1}{2}\left\vert \int_{-\pi
}^{\pi }f\left( t\right) \overline{g\left( t\right) }dt\right\vert  \notag \\
& \leq \frac{1}{2}\left[ \left( \int_{-\pi }^{\pi }\left\vert f\left(
t\right) \right\vert ^{2}dt\int_{-\pi }^{\pi }\left\vert g\left( t\right)
\right\vert ^{2}dt\right) ^{\frac{1}{2}}+\left\vert \int_{-\pi }^{\pi
}f\left( t\right) \overline{g\left( t\right) }dt\right\vert \right]  \notag
\\
& \leq \left( \int_{-\pi }^{\pi }\left\vert f\left( t\right) \right\vert
^{2}dt\int_{-\pi }^{\pi }\left\vert g\left( t\right) \right\vert
^{2}dt\right) ^{\frac{1}{2}},  \notag
\end{align}%
where $f\in L_{2}\left[ -\pi ,\pi \right] .$

Finally, by employing Theorem \ref{ch4.t3.1}, we may state:%
\begin{align*}
& \frac{1}{\pi }\left\vert \sum_{k=1}^{n}\left[ \int_{-\pi }^{\pi }f\left(
t\right) \limfunc{trig}\left( kt\right) dt\right] ^{2}\right\vert \\
& \leq \left\vert \frac{1}{\pi }\sum_{k=1}^{n}\left[ \int_{-\pi }^{\pi
}f\left( t\right) \limfunc{trig}\left( kt\right) dt\right] ^{2}-\frac{1}{2}%
\int_{-\pi }^{\pi }f^{2}\left( t\right) dt\right\vert +\frac{1}{2}\left\vert
\int_{-\pi }^{\pi }f^{2}\left( t\right) dt\right\vert \\
& \leq \frac{1}{2}\left[ \int_{-\pi }^{\pi }\left\vert f\left( t\right)
\right\vert ^{2}dt+\left\vert \int_{-\pi }^{\pi }f^{2}\left( t\right)
dt\right\vert \right] \leq \int_{-\pi }^{\pi }\left\vert f\left( t\right)
\right\vert ^{2}dt,
\end{align*}%
where $f\in L_{2}\left[ -\pi ,\pi \right] .$

\section{Generalizations of Precupanu's Inequality}

\subsection{\label{ch5.s1}Introduction}

In 1976, T. Precupanu \cite{P} obtained the following result related to the
Schwarz inequality in a real inner product space $\left( H;\left\langle
\cdot ,\cdot \right\rangle \right) :$

\begin{theorem}[Precupanu, 1976]
\label{ch5.t1.1}For any $a\in H,$ $x,y\in H\backslash \left\{ 0\right\} ,$
we have the inequality:%
\begin{align}
& \frac{-\left\Vert a\right\Vert \left\Vert b\right\Vert +\left\langle
a,b\right\rangle }{2}  \label{ch5.1.1} \\
& \leq \frac{\left\langle x,a\right\rangle \left\langle x,b\right\rangle }{%
\left\Vert x\right\Vert ^{2}}+\frac{\left\langle y,a\right\rangle
\left\langle y,b\right\rangle }{\left\Vert y\right\Vert ^{2}}-2\cdot \frac{%
\left\langle x,a\right\rangle \left\langle y,b\right\rangle \left\langle
x,y\right\rangle }{\left\Vert x\right\Vert ^{2}\left\Vert y\right\Vert ^{2}}
\notag \\
& \leq \frac{\left\Vert a\right\Vert \left\Vert b\right\Vert +\left\langle
a,b\right\rangle }{2}.  \notag
\end{align}%
In the right-hand side or in the left-hand side of (\ref{ch5.1.1}) we have
equality if and only if there are $\lambda ,\mu \in \mathbb{R}$ such that%
\begin{equation}
\lambda \frac{\left\langle x,a\right\rangle }{\left\Vert x\right\Vert ^{2}}%
\cdot x+\mu \frac{\left\langle y,b\right\rangle }{\left\Vert y\right\Vert
^{2}}\cdot y=\frac{1}{2}\left( \lambda a+\mu b\right) .  \label{ch5.1.2}
\end{equation}
\end{theorem}

Note for instance that \cite{P}, if $y\perp b,$ i.e., $\left\langle
y,b\right\rangle =0,$ then by (\ref{ch5.1.1}) one may deduce:%
\begin{equation}
\frac{-\left\Vert a\right\Vert \left\Vert b\right\Vert +\left\langle
a,b\right\rangle }{2}\left\Vert x\right\Vert ^{2}\leq \left\langle
x,a\right\rangle \left\langle x,b\right\rangle \leq \frac{\left\Vert
a\right\Vert \left\Vert b\right\Vert +\left\langle a,b\right\rangle }{2}%
\left\Vert x\right\Vert ^{2}  \label{ch5.1.3}
\end{equation}%
for any $a,b,x\in H,$ an inequality that has been obtained previously by U.
Richard \cite{R}. The case of equality in the right-hand side or in the
left-hand side of (\ref{ch5.1.3}) holds if and only if there are $\lambda
,\mu \in \mathbb{R}$ with%
\begin{equation}
2\lambda \left\langle x,a\right\rangle x=\left( \lambda a+\mu b\right)
\left\Vert x\right\Vert ^{2}.  \label{ch5.1.4}
\end{equation}%
For $a=b,$ we may obtain from (\ref{ch5.1.1}) the following inequality \cite%
{P}%
\begin{equation}
0\leq \frac{\left\langle x,a\right\rangle ^{2}}{\left\Vert x\right\Vert ^{2}}%
+\frac{\left\langle y,a\right\rangle ^{2}}{\left\Vert y\right\Vert ^{2}}%
-2\cdot \frac{\left\langle x,a\right\rangle \left\langle y,a\right\rangle
\left\langle x,y\right\rangle }{\left\Vert x\right\Vert ^{2}\left\Vert
y\right\Vert ^{2}}\leq \left\Vert a\right\Vert ^{2}.  \label{ch5.1.5}
\end{equation}%
This inequality implies \cite{P}:%
\begin{equation}
\frac{\left\langle x,y\right\rangle }{\left\Vert x\right\Vert \left\Vert
y\right\Vert }\geq \frac{1}{2}\left[ \frac{\left\langle x,a\right\rangle }{%
\left\Vert x\right\Vert \left\Vert a\right\Vert }+\frac{\left\langle
y,a\right\rangle }{\left\Vert y\right\Vert \left\Vert a\right\Vert }\right]
^{2}-\frac{3}{2}.  \label{ch5.1.6}
\end{equation}

In \cite{M}, M.H. Moore pointed out the following reverse of the Schwarz
inequality%
\begin{equation}
\left\vert \left\langle y,z\right\rangle \right\vert \leq \left\Vert
y\right\Vert \left\Vert z\right\Vert ,\qquad y,z\in H,  \label{ch5.1.7}
\end{equation}%
where some information about a third vector $x$ is known:

\begin{theorem}[Moore, 1973]
\label{ch5.t1.2}Let $\left( H;\left\langle \cdot ,\cdot \right\rangle
\right) $ be an inner product space over the real field $\mathbb{R}$ and $%
x,y,z\in H$ such that:%
\begin{equation}
\left\vert \left\langle x,y\right\rangle \right\vert \geq \left(
1-\varepsilon \right) \left\Vert x\right\Vert \left\Vert y\right\Vert
,\qquad \left\vert \left\langle x,z\right\rangle \right\vert \geq \left(
1-\varepsilon \right) \left\Vert x\right\Vert \left\Vert z\right\Vert ,
\label{ch5.1.8}
\end{equation}%
where $\varepsilon $ is a positive real number, reasonably small. Then%
\begin{equation}
\left\vert \left\langle y,z\right\rangle \right\vert \geq \max \left\{
1-\varepsilon -\sqrt{2\varepsilon },1-4\varepsilon ,0\right\} \left\Vert
y\right\Vert \left\Vert z\right\Vert .  \label{ch5.1.9}
\end{equation}
\end{theorem}

Utilising Richard's inequality (\ref{ch5.1.3}) written in the following
equivalent form:%
\begin{equation}
2\cdot \frac{\left\langle x,a\right\rangle \left\langle x,b\right\rangle }{%
\left\Vert x\right\Vert ^{2}}-\left\Vert a\right\Vert \left\Vert
b\right\Vert \leq \left\langle a,b\right\rangle \leq 2\cdot \frac{%
\left\langle x,a\right\rangle \left\langle x,b\right\rangle }{\left\Vert
x\right\Vert ^{2}}+\left\Vert a\right\Vert \left\Vert b\right\Vert
\label{ch5.1.10}
\end{equation}%
for any $a,b\in H$ and $a\in H\backslash \left\{ 0\right\} ,$ Precupanu has
obtained the following Moore's type result:

\begin{theorem}[Precupanu, 1976]
\label{ch5.t1.3}Let $\left( H;\left\langle \cdot ,\cdot \right\rangle
\right) $ be a real inner product space. If $a,b,x\in H$ and $0<\varepsilon
_{1}<\varepsilon _{2}$ are such that:%
\begin{align}
\varepsilon _{1}\left\Vert x\right\Vert \left\Vert a\right\Vert & \leq
\left\langle x,a\right\rangle \leq \varepsilon _{2}\left\Vert x\right\Vert
\left\Vert a\right\Vert ,  \label{ch5.1.11} \\
\varepsilon _{1}\left\Vert x\right\Vert \left\Vert b\right\Vert & \leq
\left\langle x,b\right\rangle \leq \varepsilon _{2}\left\Vert x\right\Vert
\left\Vert b\right\Vert ,  \notag
\end{align}%
then%
\begin{equation}
\left( 2\varepsilon _{1}^{2}-1\right) \left\Vert a\right\Vert \left\Vert
b\right\Vert \leq \left\langle a,b\right\rangle \leq \left( 2\varepsilon
_{1}^{2}+1\right) \left\Vert a\right\Vert \left\Vert b\right\Vert .
\label{ch5.1.12}
\end{equation}
\end{theorem}

Remark that the right inequality is always satisfied, since by Schwarz's
inequality, we have $\left\langle a,b\right\rangle \leq \left\Vert
a\right\Vert \left\Vert b\right\Vert $. The left inequality may be useful
when one assumes that $\varepsilon _{1}\in (0,1].$ In that case, from (\ref%
{ch5.1.12}), we obtain%
\begin{equation}
-\left\Vert a\right\Vert \left\Vert b\right\Vert \leq \left( 2\varepsilon
_{1}^{2}-1\right) \left\Vert a\right\Vert \left\Vert b\right\Vert \leq
\left\langle a,b\right\rangle  \label{ch5.1.13}
\end{equation}%
provided $\varepsilon _{1}\left\Vert x\right\Vert \left\Vert a\right\Vert
\leq \left\langle x,a\right\rangle $ and $\varepsilon _{1}\left\Vert
x\right\Vert \left\Vert b\right\Vert \leq \left\langle x,b\right\rangle ,$
which is a refinement of Schwarz's inequality%
\begin{equation*}
-\left\Vert a\right\Vert \left\Vert b\right\Vert \leq \left\langle
a,b\right\rangle .
\end{equation*}

In the complex case, apparently independent of Richard, M.L. Buzano obtained
in \cite{B} the following inequality%
\begin{equation}
\left\vert \left\langle x,a\right\rangle \left\langle x,b\right\rangle
\right\vert \leq \frac{\left\Vert a\right\Vert \left\Vert b\right\Vert
+\left\vert \left\langle a,b\right\rangle \right\vert }{2}\cdot \left\Vert
x\right\Vert ^{2},  \label{ch5.1.14}
\end{equation}%
provided $x,a,b$ are vectors in the complex inner product space $\left(
H;\left\langle \cdot ,\cdot \right\rangle \right) .$

In the same paper \cite{P}, Precupanu, without mentioning Buzano's name in
relation to the inequality (\ref{ch5.1.14}), observed that, on utilising (%
\ref{ch5.1.14}), one may obtain the following result of Moore type:

\begin{theorem}[Precupanu, 1976]
\label{ch5.t1.4}Let $\left( H;\left\langle \cdot ,\cdot \right\rangle
\right) $ be a (real or) complex inner product space. If $x,a,b\in H$ are
such that%
\begin{equation}
\left\vert \left\langle x,a\right\rangle \right\vert \geq \left(
1-\varepsilon \right) \left\Vert x\right\Vert \left\Vert a\right\Vert
,\qquad \left\vert \left\langle x,b\right\rangle \right\vert \geq \left(
1-\varepsilon \right) \left\Vert x\right\Vert \left\Vert b\right\Vert ,
\label{ch5.1.15}
\end{equation}%
then%
\begin{equation}
\left\vert \left\langle a,b\right\rangle \right\vert \geq \left(
1-4\varepsilon +2\varepsilon ^{2}\right) \left\Vert a\right\Vert \left\Vert
b\right\Vert .  \label{ch5.1.16}
\end{equation}
\end{theorem}

Note that the above theorem is useful when, for $\varepsilon \in (0,1],$ the
quantity $1-4\varepsilon +2\varepsilon ^{2}>0,$ i.e., $\varepsilon \in
\left( 0,1-\frac{\sqrt{2}}{2}\right] .$

\begin{remark}
When the space is real, the inequality (\ref{ch5.1.16}) provides a better
lower bound for $\left\vert \left\langle a,b\right\rangle \right\vert $ than
the second bound in Moore's result (\ref{ch5.1.9}). However, it is not known
if the first bound in (\ref{ch5.1.9}) remains valid for the case of complex
spaces. From Moore's original proof, apparently, the fact that the space $%
\left( H;\left\langle \cdot ,\cdot \right\rangle \right) $ is real plays an
essential role.
\end{remark}

Before we point out some new results for orthonormal families of vectors in
real or complex inner product spaces, we state the following result that
complements the Moore type results outlined above for real spaces \cite%
{DRAG3}:

\begin{theorem}[Dragomir, 2004]
\label{ch5.t1.5}Let $\left( H;\left\langle \cdot ,\cdot \right\rangle
\right) $ be a real inner product space and $a,b,x,y\in H\backslash \left\{
0\right\} .$

\begin{enumerate}
\item[(i)] If there exist $\delta _{1},\delta _{2}\in (0,1]$ such that%
\begin{equation*}
\frac{\left\langle x,a\right\rangle }{\left\Vert x\right\Vert \left\Vert
a\right\Vert }\geq \delta _{1},\qquad \frac{\left\langle y,a\right\rangle }{%
\left\Vert y\right\Vert \left\Vert a\right\Vert }\geq \delta _{2}
\end{equation*}%
and $\delta _{1}+\delta _{2}\geq 1,$ then%
\begin{equation}
\frac{\left\langle x,y\right\rangle }{\left\Vert x\right\Vert \left\Vert
y\right\Vert }\geq \frac{1}{2}\left( \delta _{1}+\delta _{2}\right) ^{2}-%
\frac{3}{2}\qquad \left( \geq -1\right) .  \label{ch5.1.17}
\end{equation}

\item[(ii)] If there exist $\mu _{1}\left( \mu _{2}\right) \in \mathbb{R}$
such that%
\begin{equation*}
\mu _{1}\left\Vert a\right\Vert \left\Vert b\right\Vert \leq \frac{%
\left\langle x,a\right\rangle \left\langle x,b\right\rangle }{\left\Vert
x\right\Vert ^{2}}\left( \leq \mu _{2}\left\Vert a\right\Vert \left\Vert
b\right\Vert \right)
\end{equation*}%
and $1\geq \mu _{1}\geq 0$ $\left( -1\leq \mu _{2}\leq 0\right) ,$ then%
\begin{equation}
\left[ -1\leq \right] 2\mu _{1}-1\leq \frac{\left\langle a,b\right\rangle }{%
\left\Vert a\right\Vert \left\Vert b\right\Vert }\left( \leq 2\mu _{2}+1%
\left[ \leq 1\right] \right) .  \label{ch5.1.18}
\end{equation}
\end{enumerate}
\end{theorem}

The proof is obvious by the inequalities (\ref{ch5.1.6}) and (\ref{ch5.1.10}%
). We omit the details.

\subsection{Inequalities for Orthonormal Families}

The following result may be stated \cite{DRAG3}.

\begin{theorem}[Dragomir, 2004]
\label{ch5.t2.1}Let $\left\{ e_{i}\right\} _{i\in I}$ and $\left\{
f_{j}\right\} _{j\in J}$ be two finite families of orthonormal vectors in $%
\left( H;\left\langle \cdot ,\cdot \right\rangle \right) .$ For any $x,y\in
H\backslash \left\{ 0\right\} $ one has the inequality%
\begin{multline}
\left\vert \sum_{i\in I}\left\langle x,e_{i}\right\rangle \left\langle
e_{i},y\right\rangle +\sum_{j\in J}\left\langle x,f_{j}\right\rangle
\left\langle f_{j},y\right\rangle \right.  \label{ch5.2.1} \\
-\left. 2\sum_{i\in I,j\in J}\left\langle x,e_{i}\right\rangle \left\langle
f_{j},y\right\rangle \left\langle e_{i},f_{j}\right\rangle -\frac{1}{2}%
\left\langle x,y\right\rangle \right\vert \leq \frac{1}{2}\left\Vert
x\right\Vert \left\Vert y\right\Vert .
\end{multline}%
The case of equality holds in (\ref{ch5.2.1}) if and only if there exists a $%
\lambda \in \mathbb{K}$ such that%
\begin{equation}
x-\lambda y=2\left( \sum_{i\in I}\left\langle x,e_{i}\right\rangle
e_{i}-\lambda \sum_{j\in J}\left\langle y,f_{j}\right\rangle f_{j}\right) .
\label{ch5.2.2}
\end{equation}
\end{theorem}

\begin{proof}
We follow the proof in \cite{DRAG3}.

We know that, if $u,v\in H,$ $v\neq 0,$ then%
\begin{equation}
\left\Vert u-\frac{\left\langle u,v\right\rangle }{\left\Vert v\right\Vert
^{2}}\cdot v\right\Vert ^{2}=\frac{\left\Vert u\right\Vert ^{2}\left\Vert
v\right\Vert ^{2}-\left\vert \left\langle u,v\right\rangle \right\vert ^{2}}{%
\left\Vert v\right\Vert ^{2}}  \label{ch5.2.3}
\end{equation}%
showing that, in Schwarz's inequality%
\begin{equation}
\left\vert \left\langle u,v\right\rangle \right\vert ^{2}\leq \left\Vert
u\right\Vert ^{2}\left\Vert v\right\Vert ^{2},  \label{ch5.2.4}
\end{equation}%
the case of equality, for $v\neq 0,$ holds if and only if%
\begin{equation}
u=\frac{\left\langle u,v\right\rangle }{\left\Vert v\right\Vert ^{2}}\cdot
v,\   \label{ch5.2.5}
\end{equation}%
i.e. there exists a $\lambda \in \mathbb{R}$ \ such that $u=\lambda v.$

Now, let $u:=2\sum_{i\in I}\left\langle x,e_{i}\right\rangle e_{i}-x$ and $%
v:=2\sum_{j\in J}\left\langle y,f_{j}\right\rangle f_{j}-y.$

Observe that%
\begin{align*}
\left\Vert u\right\Vert ^{2}& =\left\Vert 2\sum_{i\in I}\left\langle
x,e_{i}\right\rangle e_{i}\right\Vert ^{2}-4\func{Re}\left\langle \sum_{i\in
I}\left\langle x,e_{i}\right\rangle e_{i},x\right\rangle +\left\Vert
x\right\Vert ^{2} \\
& =4\sum_{i\in I}\left\vert \left\langle x,e_{i}\right\rangle \right\vert
^{2}-4\sum_{i\in I}\left\vert \left\langle x,e_{i}\right\rangle \right\vert
^{2}+\left\Vert x\right\Vert ^{2}=\left\Vert x\right\Vert ^{2},
\end{align*}%
and, similarly%
\begin{equation*}
\left\Vert v\right\Vert ^{2}=\left\Vert y\right\Vert ^{2}.
\end{equation*}%
Also,%
\begin{multline*}
\left\langle u,v\right\rangle =4\sum_{i\in I,j\in J}\left\langle
x,e_{i}\right\rangle \left\langle f_{j},y\right\rangle \left\langle
e_{i},f_{j}\right\rangle +\left\langle x,y\right\rangle \\
-2\sum_{i\in I}\left\langle x,e_{i}\right\rangle \left\langle
e_{i},y\right\rangle -2\sum_{j\in J}\left\langle x,f_{j}\right\rangle
\left\langle f_{j},y\right\rangle .
\end{multline*}%
Therefore, by Schwarz's inequality (\ref{ch5.2.4}) we deduce the desired
inequality (\ref{ch5.2.1}). By (\ref{ch5.2.5}), the case of equality holds
in (\ref{ch5.2.1}) if and only if there exists a $\lambda \in \mathbb{K}$
such that%
\begin{equation*}
2\sum_{i\in I}\left\langle x,e_{i}\right\rangle e_{i}-x=\lambda \left(
2\sum_{j\in J}\left\langle y,f_{j}\right\rangle f_{j}-y\right) ,
\end{equation*}%
which is equivalent to (\ref{ch5.2.2}).
\end{proof}

\begin{remark}
If in (\ref{ch5.2.2}) we choose $x=y,$ then we get the inequality:%
\begin{multline}
\left\vert \sum_{i\in I}\left\vert \left\langle x,e_{i}\right\rangle
\right\vert ^{2}+\sum_{j\in J}\left\vert \left\langle x,f_{j}\right\rangle
\right\vert ^{2}\right.  \label{ch5.2.6} \\
-2\left. \sum_{i\in I,j\in J}\left\langle x,e_{i}\right\rangle \left\langle
f_{j},x\right\rangle \left\langle e_{i},f_{j}\right\rangle -\frac{1}{2}%
\left\Vert x\right\Vert ^{2}\right\vert \leq \frac{1}{2}\left\Vert
x\right\Vert ^{2}
\end{multline}%
for any $x\in H.$

If in the above theorem we assume that $I=J$ and $f_{i}=e_{i},$ $i\in I,$
then we get from (\ref{ch5.2.1}) the Schwarz inequality $\left\vert
\left\langle x,y\right\rangle \right\vert \leq \left\Vert x\right\Vert
\left\Vert y\right\Vert .$

If $I\cap J=\varnothing ,$ $I\cup J=K,$ $g_{k}=e_{k},$ $k\in I,$ $%
g_{k}=f_{k},$ $k\in J$ and $\left\{ g_{k}\right\} _{k\in K}$ is orthonormal,
then from (\ref{ch5.2.1}) we get:%
\begin{equation}
\left\vert \sum_{k\in K}\left\langle x,g_{k}\right\rangle \left\langle
g_{k},y\right\rangle -\frac{1}{2}\left\langle x,y\right\rangle \right\vert
\leq \frac{1}{2}\left\Vert x\right\Vert \left\Vert y\right\Vert ,\qquad
x,y\in H  \label{ch5.2.7}
\end{equation}%
which has been obtained earlier by the author in \cite{DRAG1}.
\end{remark}

If $I$ and $J$ reduce to one element, namely $e_{1}=\frac{e}{\left\Vert
e\right\Vert },$ $f_{1}=\frac{f}{\left\Vert f\right\Vert }$ with $e,f\neq 0,$
then from (\ref{ch5.2.1}) we get%
\begin{multline}
\left\vert \frac{\left\langle x,e\right\rangle \left\langle e,y\right\rangle 
}{\left\Vert e\right\Vert ^{2}}+\frac{\left\langle x,f\right\rangle
\left\langle f,y\right\rangle }{\left\Vert f\right\Vert ^{2}}-2\cdot \frac{%
\left\langle x,e\right\rangle \left\langle f,y\right\rangle \left\langle
e,f\right\rangle }{\left\Vert e\right\Vert ^{2}\left\Vert f\right\Vert ^{2}}-%
\frac{1}{2}\left\langle x,y\right\rangle \right\vert  \label{ch5.2.8} \\
\leq \frac{1}{2}\left\Vert x\right\Vert \left\Vert y\right\Vert ,\qquad
x,y\in H
\end{multline}%
which is the corresponding complex version of Precupanu's inequality (\ref%
{ch5.1.1}).

If in (\ref{ch5.2.8}) we assume that $x=y,$ then we get%
\begin{multline}
\left\vert \frac{\left\vert \left\langle x,e\right\rangle \right\vert ^{2}}{%
\left\Vert e\right\Vert ^{2}}+\frac{\left\vert \left\langle x,f\right\rangle
\right\vert ^{2}}{\left\Vert f\right\Vert ^{2}}-2\cdot \frac{\left\langle
x,e\right\rangle \left\langle f,e\right\rangle \left\langle e,f\right\rangle 
}{\left\Vert e\right\Vert ^{2}\left\Vert f\right\Vert ^{2}}-\frac{1}{2}%
\left\Vert x\right\Vert ^{2}\right\vert  \label{ch5.2.9} \\
\leq \frac{1}{2}\left\Vert x\right\Vert ^{2}.
\end{multline}

The following corollary may be stated \cite{DRAG3}:

\begin{corollary}[Dragomir, 2004]
\label{ch5.c2.3}With the assumptions of Theorem \ref{ch5.t2.1}, we have:%
\begin{align}
& \left\vert \sum_{i\in I}\left\langle x,e_{i}\right\rangle \left\langle
e_{i},y\right\rangle +\sum_{j\in J}\left\langle x,f_{j}\right\rangle
\left\langle f_{j},y\right\rangle \right.  \label{ch5.2.10} \\
& \qquad \qquad -\left. 2\sum_{i\in I,j\in J}\left\langle
x,e_{i}\right\rangle \left\langle f_{j},y\right\rangle \left\langle
e_{i},f_{j}\right\rangle \right\vert  \notag \\
& \leq \frac{1}{2}\left\vert \left\langle x,y\right\rangle \right\vert
+\left\vert \sum_{i\in I}\left\langle x,e_{i}\right\rangle \left\langle
e_{i},y\right\rangle +\sum_{j\in J}\left\langle x,f_{j}\right\rangle
\left\langle f_{j},y\right\rangle \right.  \notag \\
& \qquad \qquad \qquad -\left. 2\sum_{i\in I,j\in J}\left\langle
x,e_{i}\right\rangle \left\langle f_{j},y\right\rangle \left\langle
e_{i},f_{j}\right\rangle -\frac{1}{2}\left\vert \left\langle
x,y\right\rangle \right\vert \right\vert  \notag \\
& \leq \frac{1}{2}\left[ \left\vert \left\langle x,y\right\rangle
\right\vert +\left\Vert x\right\Vert \left\Vert y\right\Vert \right] . 
\notag
\end{align}
\end{corollary}

\begin{proof}
The first inequality follows by the triangle inequality for the modulus. The
second inequality follows by (\ref{ch5.2.1}) on adding the quantity $\frac{1%
}{2}\left\vert \left\langle x,y\right\rangle \right\vert $ on both sides.
\end{proof}

\begin{remark}

\begin{enumerate}
\item \label{ch5.r2.4}If we choose in (\ref{ch5.2.10}), $x=y,$ then we get:%
\begin{align}
& \left\vert \sum_{i\in I}\left\vert \left\langle x,e_{i}\right\rangle
\right\vert ^{2}+\sum_{j\in J}\left\vert \left\langle x,f_{j}\right\rangle
\right\vert ^{2}\right.  \label{ch5.2.11} \\
& \qquad -\left. 2\sum_{i\in I,j\in J}\left\langle x,e_{i}\right\rangle
\left\langle f_{j},x\right\rangle \left\langle e_{i},f_{j}\right\rangle
\right\vert  \notag \\
& \leq \left\vert \sum_{i\in I}\left\vert \left\langle x,e_{i}\right\rangle
\right\vert ^{2}+\sum_{j\in J}\left\vert \left\langle x,f_{j}\right\rangle
\right\vert ^{2}\right.  \notag \\
& \qquad -\left. 2\sum_{i\in I,j\in J}\left\langle x,e_{i}\right\rangle
\left\langle f_{j},x\right\rangle \left\langle e_{i},f_{j}\right\rangle -%
\frac{1}{2}\left\Vert x\right\Vert ^{2}\right\vert +\frac{1}{2}\left\Vert
x\right\Vert ^{2}  \notag \\
& \leq \left\Vert x\right\Vert ^{2}.  \notag
\end{align}%
We observe that (\ref{ch5.2.11}) will generate Bessel's inequality if $%
\left\{ e_{i}\right\} _{i\in I},$ $\left\{ f_{j}\right\} _{j\in J}$ are
disjoint parts of a larger orthonormal family.

\item From (\ref{ch5.2.8}) one can obtain:%
\begin{multline}
\left\vert \frac{\left\langle x,e\right\rangle \left\langle e,y\right\rangle 
}{\left\Vert e\right\Vert ^{2}}+\frac{\left\langle x,f\right\rangle
\left\langle f,y\right\rangle }{\left\Vert f\right\Vert ^{2}}-2\cdot \frac{%
\left\langle x,e\right\rangle \left\langle f,y\right\rangle \left\langle
e,f\right\rangle }{\left\Vert e\right\Vert ^{2}\left\Vert f\right\Vert ^{2}}%
\right\vert  \label{ch5.2.12} \\
\leq \frac{1}{2}\left[ \left\Vert x\right\Vert \left\Vert y\right\Vert
+\left\vert \left\langle x,y\right\rangle \right\vert \right]
\end{multline}%
and in particular 
\begin{equation}
\left\vert \frac{\left\vert \left\langle x,e\right\rangle \right\vert ^{2}}{%
\left\Vert e\right\Vert ^{2}}+\frac{\left\vert \left\langle x,f\right\rangle
\right\vert ^{2}}{\left\Vert f\right\Vert ^{2}}-2\cdot \frac{\left\langle
x,e\right\rangle \left\langle f,e\right\rangle \left\langle e,f\right\rangle 
}{\left\Vert e\right\Vert ^{2}\left\Vert f\right\Vert ^{2}}\right\vert \leq
\left\Vert x\right\Vert ^{2},  \label{ch5.2.13}
\end{equation}%
for any $x,y\in H.$
\end{enumerate}
\end{remark}

The case of real inner products will provide a natural genearlization for
Precupanu's inequality (\ref{ch5.1.1}) \cite{DRAG3}:

\begin{corollary}[Dragomir, 2004]
\label{ch5.c2.5}Let $\left( H;\left\langle \cdot ,\cdot \right\rangle
\right) $ be a real inner product space and $\left\{ e_{i}\right\} _{i\in
I}, $ $\left\{ f_{j}\right\} _{j\in J}$ two finite families of orthonormal
vectors in $\left( H;\left\langle \cdot ,\cdot \right\rangle \right) .$ For
any $x,y\in H\backslash \left\{ 0\right\} $ one has the double inequality:%
\begin{align}
\frac{1}{2}\left[ \left\vert \left\langle x,y\right\rangle \right\vert
-\left\Vert x\right\Vert \left\Vert y\right\Vert \right] & \leq \sum_{i\in
I}\left\langle x,e_{i}\right\rangle \left\langle y,e_{i}\right\rangle
+\sum_{j\in J}\left\langle x,f_{j}\right\rangle \left\langle
y,f_{j}\right\rangle  \label{ch5.2.14} \\
& \qquad \qquad -2\sum_{i\in I,j\in J}\left\langle x,e_{i}\right\rangle
\left\langle y,f_{j}\right\rangle \left\langle e_{i},f_{j}\right\rangle 
\notag \\
& \leq \frac{1}{2}\left[ \left\Vert x\right\Vert \left\Vert y\right\Vert
+\left\vert \left\langle x,y\right\rangle \right\vert \right] .  \notag
\end{align}%
In particular, we have%
\begin{align}
0& \leq \sum_{i\in I}\left\langle x,e_{i}\right\rangle ^{2}+\sum_{j\in
J}\left\langle x,f_{j}\right\rangle ^{2}-2\sum_{i\in I,j\in J}\left\langle
x,e_{i}\right\rangle \left\langle x,f_{j}\right\rangle \left\langle
e_{i},f_{j}\right\rangle  \label{ch5.2.15} \\
& \leq \left\Vert x\right\Vert ^{2},  \notag
\end{align}%
for any $x\in H.$
\end{corollary}

\begin{remark}
Similar particular inequalities to those incorporated in (\ref{ch5.2.7}) -- (%
\ref{ch5.2.13}) may be stated, but we omit them.
\end{remark}

\subsection{Refinements of Kurepa's Inequality\label{ch5.s3}}

The following result may be stated \cite{DRAG3}.

\begin{theorem}[Dragomir, 2004]
\label{ch5.t3.1}Let $\left( H;\left\langle \cdot ,\cdot \right\rangle
\right) $ be a real inner product space and $\left\{ e_{i}\right\} _{i\in
I},\left\{ f_{j}\right\} _{j\in J}$ two finite families in $H.$ If $\left(
H_{\mathbb{C}};\left\langle \cdot ,\cdot \right\rangle _{\mathbb{C}}\right) $
is the complexification of $\left( H;\left\langle \cdot ,\cdot \right\rangle
\right) ,$ then for any $w\in H_{\mathbb{C}},$ we have the inequalities%
\begin{align}
& \left\vert \sum_{i\in I}\left\langle w,e_{i}\right\rangle _{\mathbb{C}%
}^{2}+\sum_{j\in J}\left\langle w,f_{j}\right\rangle _{\mathbb{C}%
}^{2}-2\sum_{i\in I,j\in J}\left\langle w,e_{i}\right\rangle _{\mathbb{C}%
}\left\langle w,f_{j}\right\rangle _{\mathbb{C}}\left\langle
e_{i},f_{j}\right\rangle \right\vert  \label{ch5.3.3} \\
& \leq \frac{1}{2}\left\vert \left\langle w,\bar{w}\right\rangle _{\mathbb{C}%
}\right\vert +\left\vert \sum_{i\in I}\left\langle w,e_{i}\right\rangle _{%
\mathbb{C}}^{2}+\sum_{j\in J}\left\langle w,f_{j}\right\rangle _{\mathbb{C}%
}^{2}\right.  \notag \\
& \qquad \qquad \qquad -\left. 2\sum_{i\in I,j\in J}\left\langle
w,e_{i}\right\rangle _{\mathbb{C}}\left\langle w,f_{j}\right\rangle _{%
\mathbb{C}}\left\langle e_{i},f_{j}\right\rangle -\frac{1}{2}\left\langle w,%
\bar{w}\right\rangle _{\mathbb{C}}\right\vert  \notag \\
& \leq \frac{1}{2}\left[ \left\Vert w\right\Vert _{\mathbb{C}%
}^{2}+\left\vert \left\langle w,\bar{w}\right\rangle _{\mathbb{C}%
}\right\vert \right] \leq \left\Vert w\right\Vert _{\mathbb{C}}^{2}.  \notag
\end{align}
\end{theorem}

\begin{proof}
Define $g_{j}\in H_{\mathbb{C}},$ $g_{j}:=\left( e_{j},0\right) ,$ $j\in I.$
For any $k,j\in I$ we have%
\begin{equation*}
\left\langle g_{k},g_{j}\right\rangle _{\mathbb{C}}=\left\langle \left(
e_{k},0\right) ,\left( e_{j},0\right) \right\rangle _{\mathbb{C}%
}=\left\langle e_{k},e_{j}\right\rangle =\delta _{kj},
\end{equation*}%
therefore $\left\{ g_{j}\right\} _{j\in I}$ is an orthonormal family in $%
\left( H_{\mathbb{C}};\left\langle \cdot ,\cdot \right\rangle _{\mathbb{C}%
}\right) .$

If we apply Corollary \ref{ch5.c2.3} for $\left( H_{\mathbb{C}};\left\langle
\cdot ,\cdot \right\rangle _{\mathbb{C}}\right) ,$ $x=w,$ $y=\bar{w},$ we
may write:%
\begin{align}
& \left\vert \sum_{i\in I}\left\langle w,e_{i}\right\rangle _{\mathbb{C}%
}\left\langle e_{i},\bar{w}\right\rangle _{\mathbb{C}}+\sum_{j\in
J}\left\langle w,f_{j}\right\rangle \left\langle f_{j},\bar{w}\right\rangle
\right.  \label{ch5.3.4} \\
& \qquad \qquad \qquad -\left. 2\sum_{i\in I,j\in J}\left\langle
w,e_{i}\right\rangle _{\mathbb{C}}\left\langle f_{j},\overline{w}%
\right\rangle _{\mathbb{C}}\left\langle e_{i},f_{j}\right\rangle \right\vert
\notag \\
& \leq \frac{1}{2}\left\Vert w\right\Vert _{\mathbb{C}}\left\Vert \bar{w}%
\right\Vert _{\mathbb{C}}+\left\vert \sum_{i\in I}\left\langle
w,e_{i}\right\rangle _{\mathbb{C}}\left\langle e_{i},\bar{w}\right\rangle _{%
\mathbb{C}}+\sum_{j\in J}\left\langle w,f_{j}\right\rangle \left\langle
f_{j},\bar{w}\right\rangle \right.  \notag \\
& \qquad \qquad \qquad -\left. 2\sum_{i\in I,j\in J}\left\langle
w,e_{i}\right\rangle _{\mathbb{C}}\left\langle f_{j},\overline{w}%
\right\rangle _{\mathbb{C}}\left\langle e_{i},f_{j}\right\rangle -\frac{1}{2}%
\left\langle w,\bar{w}\right\rangle _{\mathbb{C}}\right\vert  \notag \\
& \leq \frac{1}{2}\left[ \left\vert \left\langle w,\bar{w}\right\rangle _{%
\mathbb{C}}\right\vert +\left\Vert w\right\Vert _{\mathbb{C}}\left\Vert \bar{%
w}\right\Vert _{\mathbb{C}}\right] .  \notag
\end{align}%
However, for $w:=\left( x,y\right) \in H_{\mathbb{C}}$, we have $\bar{w}%
=\left( x,-y\right) $ and%
\begin{equation*}
\left\langle e_{j},\bar{w}\right\rangle _{\mathbb{C}}=\left\langle \left(
e_{j},0\right) ,\left( x,-y\right) \right\rangle _{\mathbb{C}}=\left\langle
e_{j},x\right\rangle +i\left\langle e_{j},y\right\rangle
\end{equation*}%
and%
\begin{equation*}
\left\langle w,e_{j}\right\rangle _{\mathbb{C}}=\left\langle \left(
x,y\right) ,\left( e_{j},0\right) \right\rangle _{\mathbb{C}}=\left\langle
x,e_{j}\right\rangle +i\left\langle e_{j},y\right\rangle
\end{equation*}%
showing that $\left\langle e_{j},\bar{w}\right\rangle =\left\langle
w,e_{j}\right\rangle $ for any $j\in I$. A similar relation is true for $%
f_{j}$ and since%
\begin{equation*}
\left\Vert w\right\Vert _{\mathbb{C}}=\left\Vert \bar{w}\right\Vert _{%
\mathbb{C}}=\left( \left\Vert x\right\Vert ^{2}+\left\Vert y\right\Vert
^{2}\right) ^{\frac{1}{2}},
\end{equation*}%
hence from (\ref{ch5.3.4}) we deduce the desired inequality (\ref{ch5.3.3}).
\end{proof}

\begin{remark}
It is obvious that, if one family, say $\left\{ f_{j}\right\} _{j\in J}$ is
empty, then, on observing that all sums $\sum_{j\in J}$ should be zero, from
(\ref{ch5.3.3}) one would get \cite{DRAG1}%
\begin{align}
& \left\vert \sum_{i\in I}\left\langle w,e_{i}\right\rangle _{\mathbb{C}%
}^{2}\right\vert  \label{ch5.3.5} \\
& \leq \frac{1}{2}\left\vert \left\langle w,\bar{w}\right\rangle _{\mathbb{C}%
}\right\vert +\left\vert \sum_{i\in I}\left\langle w,e_{i}\right\rangle _{%
\mathbb{C}}^{2}-\frac{1}{2}\left\langle w,\bar{w}\right\rangle _{\mathbb{C}%
}\right\vert  \notag \\
& \leq \frac{1}{2}\left[ \left\Vert w\right\Vert _{\mathbb{C}%
}^{2}+\left\vert \left\langle w,\bar{w}\right\rangle _{\mathbb{C}%
}\right\vert \right] \leq \left\Vert w\right\Vert _{\mathbb{C}}^{2}.  \notag
\end{align}%
If in (\ref{ch5.3.5}) one assumes that the family $\left\{ e_{i}\right\}
_{i\in I}$ contains only one element $e=\frac{a}{\left\Vert a\right\Vert }%
,a\neq 0,$ then by selecting $w=z,$ one would deduce (\ref{ch3.3.5}), which
is a refinement for Kurepa's inequality.
\end{remark}

\section{Some New Refinements of the Schwarz Inequality}

\subsection{Refinements}

The following result holds \cite{DRAG01}.

\begin{theorem}[Dragomir, 2004]
\label{ch2.t2.1}Let $\left( H;\left\langle \cdot ,\cdot \right\rangle
\right) $ be an inner product space over the real or complex number field $%
\mathbb{K}$ and $r_{1},r_{2}>0.$ If $x,y\in H$ are with the property that 
\begin{equation}
\left\Vert x-y\right\Vert \geq r_{2}\geq r_{1}\geq \left\vert \left\Vert
x\right\Vert -\left\Vert y\right\Vert \right\vert ,  \label{ch2.2.1}
\end{equation}%
then we have the following refinement of Schwarz's inequality 
\begin{equation}
\left\Vert x\right\Vert \left\Vert y\right\Vert -\func{Re}\left\langle
x,y\right\rangle \geq \frac{1}{2}\left( r_{2}^{2}-r_{1}^{2}\right) \left(
\geq 0\right) .  \label{ch2.2.2}
\end{equation}%
The constant $\frac{1}{2}$ \ is best possible in the sense that it cannot be
replaced by a larger quantity.
\end{theorem}

\begin{proof}
From the first inequality in (\ref{ch2.2.1}) we have 
\begin{equation}
\left\| x\right\| ^{2}+\left\| y\right\| ^{2}\geq r_{2}^{2}+2\func{Re}%
\left\langle x,y\right\rangle .  \label{ch2.2.3}
\end{equation}
Subtracting in (\ref{ch2.2.3}) the quantity $2\left\| x\right\| \left\|
y\right\| ,$ we get 
\begin{equation}
\left( \left\| x\right\| -\left\| y\right\| \right) ^{2}\geq
r_{2}^{2}-2\left( \left\| x\right\| \left\| y\right\| -\func{Re}\left\langle
x,y\right\rangle \right) .  \label{ch2.2.4}
\end{equation}
Since, by the second inequality in (\ref{ch2.2.1}) we have 
\begin{equation}
r_{1}^{2}\geq \left( \left\| x\right\| -\left\| y\right\| \right) ^{2},
\label{ch2.2.5}
\end{equation}
hence from (\ref{ch2.2.4}) and (\ref{ch2.2.5}) we deduce the desired
inequality (\ref{ch2.2.2}).

To prove the sharpness of the constant $\frac{1}{2}$ in (\ref{ch2.2.2}), let
us assume that there is a constant $C>0$ such that 
\begin{equation}
\left\| x\right\| \left\| y\right\| -\func{Re}\left\langle x,y\right\rangle
\geq C\left( r_{2}^{2}-r_{1}^{2}\right) ,  \label{ch2.2.6}
\end{equation}
provided that $x$ and $y$ satisfy (\ref{ch2.2.1}).

Let $e\in H$ with $\left\Vert e\right\Vert =1$ and for $r_{2}>r_{1}>0,$
define 
\begin{equation}
x=\frac{r_{2}+r_{1}}{2}\cdot e\text{ \ and \ }y=\frac{r_{1}-r_{2}}{2}\cdot e.
\label{ch2.2.7}
\end{equation}%
Then 
\begin{equation*}
\left\Vert x-y\right\Vert =r_{2}\text{ and }\left\vert \left\Vert
x\right\Vert -\left\Vert y\right\Vert \right\vert =r_{1},
\end{equation*}%
showing that the condition (\ref{ch2.2.1}) is fulfilled with equality.

If we replace $x$ and $y$ as defined in (\ref{ch2.2.7}) into the inequality (%
\ref{ch2.2.6}), then we get 
\begin{equation*}
\frac{r_{2}^{2}-r_{1}^{2}}{2}\geq C\left( r_{2}^{2}-r_{1}^{2}\right) ,
\end{equation*}%
which implies that $C\leq \frac{1}{2},$ and the theorem is completely proved.
\end{proof}

The following corollary holds.

\begin{corollary}
\label{ch2.c2.2}With the assumptions of Theorem \ref{ch2.t2.1}, we have the
inequality: 
\begin{equation}
\left\Vert x\right\Vert +\left\Vert y\right\Vert -\frac{\sqrt{2}}{2}%
\left\Vert x+y\right\Vert \geq \frac{\sqrt{2}}{2}\sqrt{r_{2}^{2}-r_{1}^{2}}.
\label{ch2.2.8}
\end{equation}
\end{corollary}

\begin{proof}
We have, by (\ref{ch2.2.2}), that 
\begin{equation*}
\left( \left\| x\right\| +\left\| y\right\| \right) ^{2}-\left\| x+y\right\|
^{2}=2\left( \left\| x\right\| \left\| y\right\| -\func{Re}\left\langle
x,y\right\rangle \right) \geq r_{2}^{2}-r_{1}^{2}\geq 0
\end{equation*}
which gives 
\begin{equation}
\left( \left\| x\right\| +\left\| y\right\| \right) ^{2}\geq \left\|
x+y\right\| ^{2}+\left( \sqrt{r_{2}^{2}-r_{1}^{2}}\right) ^{2}.
\label{ch2.2.9}
\end{equation}
By making use of the elementary inequality 
\begin{equation*}
2\left( \alpha ^{2}+\beta ^{2}\right) \geq \left( \alpha +\beta \right)
^{2},\qquad \alpha ,\beta \geq 0;
\end{equation*}
we get 
\begin{equation}
\left\| x+y\right\| ^{2}+\left( \sqrt{r_{2}^{2}-r_{1}^{2}}\right) ^{2}\geq 
\frac{1}{2}\left( \left\| x+y\right\| +\sqrt{r_{2}^{2}-r_{1}^{2}}\right)
^{2}.  \label{ch2.2.10}
\end{equation}
Utilising (\ref{ch2.2.9}) and (\ref{ch2.2.10}), we deduce the desired
inequality (\ref{ch2.2.8}).
\end{proof}

If $\left( H;\left\langle \cdot ,\cdot \right\rangle \right) $ is a Hilbert
space and $\left\{ e_{i}\right\} _{i\in I}$ is an orthornormal family in $H,$
i.e., we recall that $\left\langle e_{i},e_{j}\right\rangle =\delta _{ij}$
for any $i,j\in I,$ where $\delta _{ij}$ is Kronecker's delta, then we have
the following inequality which is well known in the literature as \textit{%
Bessel's inequality} 
\begin{equation}
\sum_{i\in I}\left\vert \left\langle x,e_{i}\right\rangle \right\vert
^{2}\leq \left\Vert x\right\Vert ^{2}\text{ \ \ for each \ }x\in H.
\label{ch2.2.11}
\end{equation}%
Here, the meaning of the sum is 
\begin{equation*}
\sum_{i\in I}\left\vert \left\langle x,e_{i}\right\rangle \right\vert
^{2}=\sup_{F\subset I}\left\{ \sum_{i\in F}\left\vert \left\langle
x,e_{i}\right\rangle \right\vert ^{2},\ F\text{ is a finite part of }%
I\right\} .
\end{equation*}

The following result providing a refinement of the Bessel inequality (\ref%
{ch2.2.11}) holds \cite{DRAG01}.

\begin{theorem}[Dragomir, 2004]
\label{ch2.t2.3}Let $\left( H;\left\langle \cdot ,\cdot \right\rangle
\right) $ be a Hilbert space and $\left\{ e_{i}\right\} _{i\in I}$ an
orthonormal family in $H.$ If $x\in H,$ $x\neq 0,$ and $r_{2},r_{1}>0$ are
such that: 
\begin{multline}
\left\Vert x-\sum_{i\in I}\left\langle x,e_{i}\right\rangle e_{i}\right\Vert
\\
\geq r_{2}\geq r_{1}\geq \left\Vert x\right\Vert -\left( \sum_{i\in
I}\left\vert \left\langle x,e_{i}\right\rangle \right\vert ^{2}\right) ^{%
\frac{1}{2}}\left( \geq 0\right) ,  \label{ch2.2.12}
\end{multline}%
then we have the inequality 
\begin{equation}
\left\Vert x\right\Vert -\left( \sum_{i\in I}\left\vert \left\langle
x,e_{i}\right\rangle \right\vert ^{2}\right) ^{\frac{1}{2}}\geq \frac{1}{2}%
\cdot \frac{r_{2}^{2}-r_{1}^{2}}{\left( \sum_{i\in I}\left\vert \left\langle
x,e_{i}\right\rangle \right\vert ^{2}\right) ^{\frac{1}{2}}}\left( \geq
0\right) .  \label{ch2.2.13}
\end{equation}%
The constant $\frac{1}{2}$ is best possible.
\end{theorem}

\begin{proof}
Consider $y:=\sum_{i\in I}\left\langle x,e_{i}\right\rangle e_{i}.$
Obviously, since $H$ is a Hilbert space, $y\in H.$ We also note that 
\begin{equation*}
\left\Vert y\right\Vert =\left\Vert \sum_{i\in I}\left\langle
x,e_{i}\right\rangle e_{i}\right\Vert =\sqrt{\left\Vert \sum_{i\in
I}\left\langle x,e_{i}\right\rangle e_{i}\right\Vert ^{2}}=\sqrt{\sum_{i\in
I}\left\vert \left\langle x,e_{i}\right\rangle \right\vert ^{2}},
\end{equation*}
and thus (\ref{ch2.2.12}) is in fact (\ref{ch2.2.1}) of Theorem \ref%
{ch2.t2.1}.

Since 
\begin{align*}
\left\Vert x\right\Vert \left\Vert y\right\Vert -\func{Re}\left\langle
x,y\right\rangle & =\left\Vert x\right\Vert \left( \sum_{i\in I}\left\vert
\left\langle x,e_{i}\right\rangle \right\vert ^{2}\right) ^{\frac{1}{2}}-%
\func{Re}\left\langle x,\sum_{i\in I}\left\langle x,e_{i}\right\rangle
e_{i}\right\rangle \\
& =\left( \sum_{i\in I}\left\vert \left\langle x,e_{i}\right\rangle
\right\vert ^{2}\right) ^{\frac{1}{2}}\left[ \left\Vert x\right\Vert -\left(
\sum_{i\in I}\left\vert \left\langle x,e_{i}\right\rangle \right\vert
^{2}\right) ^{\frac{1}{2}}\right] ,
\end{align*}%
hence, by (\ref{ch2.2.2}), we deduce the desired result (\ref{ch2.2.13}).

We will prove the sharpness of the constant for the case of one element,
i.e., $I=\left\{ 1\right\} ,$ $e_{1}=e\in H,$ $\left\Vert e\right\Vert =1.$
For this, assume that there exists a constant $D>0$ such that 
\begin{equation}
\left\Vert x\right\Vert -\left\vert \left\langle x,e\right\rangle
\right\vert \geq D\cdot \frac{r_{2}^{2}-r_{1}^{2}}{\left\vert \left\langle
x,e\right\rangle \right\vert }  \label{ch2.2.14}
\end{equation}%
provided $x\in H\backslash \left\{ 0\right\} $ satisfies the condition 
\begin{equation}
\left\Vert x-\left\langle x,e\right\rangle e\right\Vert \geq r_{2}\geq
r_{1}\geq \left\Vert x\right\Vert -\left\vert \left\langle x,e\right\rangle
\right\vert .  \label{ch2.2.15}
\end{equation}%
Assume that $x=\lambda e+\mu f$ with $e,f\in H,$ $\left\Vert e\right\Vert
=\left\Vert f\right\Vert =1$ and $e\perp f.$ We wish to see if there exists
positive numbers $\lambda ,\mu $ such that 
\begin{equation}
\left\Vert x-\left\langle x,e\right\rangle e\right\Vert
=r_{2}>r_{1}=\left\Vert x\right\Vert -\left\vert \left\langle
x,e\right\rangle \right\vert .  \label{ch2.2.16}
\end{equation}%
Since (for $\lambda ,\mu >0$) 
\begin{equation*}
\left\Vert x-\left\langle x,e\right\rangle e\right\Vert =\mu
\end{equation*}%
and 
\begin{equation*}
\left\Vert x\right\Vert -\left\vert \left\langle x,e\right\rangle
\right\vert =\sqrt{\lambda ^{2}+\mu ^{2}}-\lambda
\end{equation*}%
hence, by (\ref{ch2.2.16}), we get $\mu =r_{2}$ and 
\begin{equation*}
\sqrt{\lambda ^{2}+r_{2}^{2}}-\lambda =r_{1}
\end{equation*}%
giving 
\begin{equation*}
\lambda ^{2}+r_{2}^{2}=\lambda ^{2}+2\lambda r_{1}+r_{1}^{2}
\end{equation*}%
from where we get 
\begin{equation*}
\lambda =\frac{r_{2}^{2}-r_{1}^{2}}{2r_{1}}>0.
\end{equation*}%
With these values for $\lambda $ and $\mu ,$ we have 
\begin{equation*}
\left\Vert x\right\Vert -\left\vert \left\langle x,e\right\rangle
\right\vert =r_{1},\qquad \left\vert \left\langle x,e\right\rangle
\right\vert =\frac{r_{2}^{2}-r_{1}^{2}}{2r_{1}}
\end{equation*}%
and thus, from (\ref{ch2.2.14}), we deduce 
\begin{equation*}
r_{1}\geq D\cdot \frac{r_{2}^{2}-r_{1}^{2}}{\frac{r_{2}^{2}-r_{1}^{2}}{2r_{1}%
}},
\end{equation*}%
giving $D\leq \frac{1}{2}.$ This proves the theorem.
\end{proof}

The following corollary is obvious.

\begin{corollary}
\label{ch2.c2.4}Let $x,y\in H$ with $\left\langle x,y\right\rangle \neq 0$
and $r_{2}\geq r_{1}>0$ such that 
\begin{align}
\left\Vert \left\Vert y\right\Vert x-\frac{\left\langle x,y\right\rangle }{%
\left\Vert y\right\Vert }\cdot y\right\Vert & \geq r_{2}\left\Vert
y\right\Vert \geq r_{1}\left\Vert y\right\Vert  \label{ch2.2.17} \\
& \geq \left\Vert x\right\Vert \left\Vert y\right\Vert -\left\vert
\left\langle x,y\right\rangle \right\vert \left( \geq 0\right) .  \notag
\end{align}%
Then we have the following refinement of the Schwarz's inequality: 
\begin{equation}
\left\Vert x\right\Vert \left\Vert y\right\Vert -\left\vert \left\langle
x,y\right\rangle \right\vert \geq \frac{1}{2}\left(
r_{2}^{2}-r_{1}^{2}\right) \frac{\left\Vert y\right\Vert ^{2}}{\left\vert
\left\langle x,y\right\rangle \right\vert }\left( \geq 0\right) .
\label{ch2.2.18}
\end{equation}%
The constant $\frac{1}{2}$ is best possible.
\end{corollary}

The following lemma holds \cite{DRAG01}.

\begin{lemma}[Dragomir, 2004]
\label{ch2.l2.5}Let $\left( H;\left\langle \cdot ,\cdot \right\rangle
\right) $ be an inner product space and $R\geq 1.$ For $x,y\in H,$ the
subsequent statements are equivalent:

\begin{enumerate}
\item[(i)] The following refinement of the triangle inequality holds: 
\begin{equation}
\left\Vert x\right\Vert +\left\Vert y\right\Vert \geq R\left\Vert
x+y\right\Vert ;  \label{ch2.2.19}
\end{equation}

\item[(ii)] The following refinement of the Schwarz inequality holds: 
\begin{equation}
\left\Vert x\right\Vert \left\Vert y\right\Vert -\func{Re}\left\langle
x,y\right\rangle \geq \frac{1}{2}\left( R^{2}-1\right) \left\Vert
x+y\right\Vert ^{2}.  \label{ch2.2.20}
\end{equation}
\end{enumerate}
\end{lemma}

\begin{proof}
Taking the square in (\ref{ch2.2.19}), we have 
\begin{equation}
2\left\Vert x\right\Vert \left\Vert y\right\Vert \geq \left( R^{2}-1\right)
\left\Vert x\right\Vert ^{2}+2R^{2}\func{Re}\left\langle x,y\right\rangle
+\left( R^{2}-1\right) \left\Vert y\right\Vert ^{2}.  \label{ch2.2.21}
\end{equation}
Subtracting from both sides of (\ref{ch2.2.21}) the quantity $2\func{Re}%
\left\langle x,y\right\rangle ,$ we obtain 
\begin{align*}
2\left( \left\Vert x\right\Vert \left\Vert y\right\Vert -\func{Re}%
\left\langle x,y\right\rangle \right) & \geq \left( R^{2}-1\right) \left[
\left\Vert x\right\Vert ^{2}+2\func{Re}\left\langle x,y\right\rangle
+\left\Vert y\right\Vert ^{2}\right] \\
& =\left( R^{2}-1\right) \left\Vert x+y\right\Vert ^{2},
\end{align*}
which is clearly equivalent to (\ref{ch2.2.20}).
\end{proof}

By the use of the above lemma, we may now state the following theorem
concerning another refinement of the Schwarz inequality \cite{DRAG01}.

\begin{theorem}[Dragomir, 2004]
\label{ch2.t2.6}Let $\left( H;\left\langle \cdot ,\cdot \right\rangle
\right) $ be an inner product space over the real or complex number field
and $R\geq 1,$ $r\geq 0.$ If $x,y\in H$ are such that 
\begin{equation}
\frac{1}{R}\left( \left\Vert x\right\Vert +\left\Vert y\right\Vert \right)
\geq \left\Vert x+y\right\Vert \geq r,  \label{ch2.2.22}
\end{equation}%
then we have the following refinement of the Schwarz inequality 
\begin{equation}
\left\Vert x\right\Vert \left\Vert y\right\Vert -\func{Re}\left\langle
x,y\right\rangle \geq \frac{1}{2}\left( R^{2}-1\right) r^{2}.
\label{ch2.2.23}
\end{equation}%
The constant $\frac{1}{2}$ is best possible in the sense that it cannot be
replaced by a larger quantity.
\end{theorem}

\begin{proof}
The inequality (\ref{ch2.2.23}) follows easily from Lemma \ref{ch2.l2.5}. We
need only prove that $\frac{1}{2}$ is the best possible constant in (\ref%
{ch2.2.23}).

Assume that there exists a $C>0$ such that 
\begin{equation}
\left\Vert x\right\Vert \left\Vert y\right\Vert -\func{Re}\left\langle
x,y\right\rangle \geq C\left( R^{2}-1\right) r^{2}  \label{ch2.2.24}
\end{equation}%
provided $x,y,$ $R$ and $r$ satisfy (\ref{ch2.2.22}).

Consider $r=1,$ $R>1$ and choose $x=\frac{1-R}{2}e,$ $y=\frac{1+R}{2}e$ with 
$e\in H,$ $\left\Vert e\right\Vert =1.$ Then 
\begin{equation*}
x+y=e,\ \ \ \frac{\left\Vert x\right\Vert +\left\Vert y\right\Vert }{R}=1
\end{equation*}%
and thus (\ref{ch2.2.22}) holds with equality on both sides.

From (\ref{ch2.2.24}), for the above choices, we have $\frac{1}{2}\left(
R^{2}-1\right) \geq C\left( R^{2}-1\right) ,$ which shows that $C\leq \frac{1%
}{2}.$
\end{proof}

Finally, the following result also holds \cite{DRAG01}.

\begin{theorem}[Dragomir, 2004]
\label{ch2.t2.7}Let $\left( H;\left\langle \cdot ,\cdot \right\rangle
\right) $ be an inner product space over the real or complex number field $%
\mathbb{K}$ and $r\in (0,1].$ For $x,y\in H,$ the following statements are
equivalent:

\begin{enumerate}
\item[(i)] We have the inequality 
\begin{equation}
\left\vert \left\Vert x\right\Vert -\left\Vert y\right\Vert \right\vert \leq
r\left\Vert x-y\right\Vert ;  \label{ch2.2.25}
\end{equation}

\item[(ii)] We have the following refinement of the Schwarz inequality 
\begin{equation}
\left\Vert x\right\Vert \left\Vert y\right\Vert -\func{Re}\left\langle
x,y\right\rangle \geq \frac{1}{2}\left( 1-r^{2}\right) \left\Vert
x-y\right\Vert ^{2}.  \label{ch2.2.26}
\end{equation}
The constant $\frac{1}{2}$ in (\ref{ch2.2.26}) is best possible.
\end{enumerate}
\end{theorem}

\begin{proof}
Taking the square in (\ref{ch2.2.25}), we have 
\begin{equation*}
\left\Vert x\right\Vert ^{2}-2\left\Vert x\right\Vert \left\Vert
y\right\Vert +\left\Vert y\right\Vert ^{2}\leq r^{2}\left( \left\Vert
x\right\Vert ^{2}-2\func{Re}\left\langle x,y\right\rangle +\left\Vert
y\right\Vert ^{2}\right)
\end{equation*}
which is clearly equivalent to 
\begin{equation*}
\left( 1-r^{2}\right) \left[ \left\Vert x\right\Vert ^{2}-2\func{Re}%
\left\langle x,y\right\rangle +\left\Vert y\right\Vert ^{2}\right] \leq
2\left( \left\Vert x\right\Vert \left\Vert y\right\Vert -\func{Re}%
\left\langle x,y\right\rangle \right)
\end{equation*}
or with (\ref{ch2.2.26}).

Now, assume that (\ref{ch2.2.26}) holds with a constant $E>0,$ i.e., 
\begin{equation}
\left\| x\right\| \left\| y\right\| -\func{Re}\left\langle x,y\right\rangle
\geq E\left( 1-r^{2}\right) \left\| x-y\right\| ^{2},  \label{ch2.2.27}
\end{equation}
provided (\ref{ch2.2.25}) holds.

Define $x=\frac{r+1}{2}e,$ $y=\frac{r-1}{2}e$ with $e\in H,$ $\left\Vert
e\right\Vert =1.$ Then 
\begin{equation*}
\left\vert \left\Vert x\right\Vert -\left\Vert y\right\Vert \right\vert
=r,\quad \left\Vert x-y\right\Vert =1
\end{equation*}%
showing that (\ref{ch2.2.25}) holds with equality.

If we replace $x$ and $y$ in (\ref{ch2.2.27}), then we get $E\left(
1-r^{2}\right) \leq \frac{1}{2}\left( 1-r^{2}\right) ,$ implying that $E\leq 
\frac{1}{2}.$
\end{proof}

\subsection{Discrete Inequalities}

Assume that $\left( K;\left( \cdot ,\cdot \right) \right) $ is a Hilbert
space over the real or complex number field. Assume also that $p_{i}\geq 0$, 
$i\in H$ with $\sum_{i=1}^{\infty }p_{i}=1$ and define 
\begin{equation*}
\ell _{p}^{2}\left( K\right) :=\left\{ \mathbf{x}:=\left. \left(
x_{i}\right) _{i\in \mathbb{N}}\right| x_{i}\in \mathbb{K},\ i\in \mathbb{N}%
\text{ \ and \ }\sum_{i=1}^{\infty }p_{i}\left\| x_{i}\right\| ^{2}<\infty
\right\} .
\end{equation*}
It is well known that $\ell _{p}^{2}\left( K\right) $ endowed with the inner
product $\left\langle \cdot ,\cdot \right\rangle _{p}$ defined by 
\begin{equation*}
\left\langle \mathbf{x},\mathbf{y}\right\rangle _{p}:=\sum_{i=1}^{\infty
}p_{i}\left( x_{i},y_{i}\right)
\end{equation*}
and generating the norm 
\begin{equation*}
\left\| \mathbf{x}\right\| _{p}:=\left( \sum_{i=1}^{\infty }p_{i}\left\|
x_{i}\right\| ^{2}\right) ^{\frac{1}{2}}
\end{equation*}
is a Hilbert space over $\mathbb{K}$.

We may state the following discrete inequality improving the
Cauchy-Bunyakovsky-Schwarz classical result \cite{DRAG01}.

\begin{proposition}
\label{ch2.p3.1}Let $\left( K;\left( \cdot ,\cdot \right) \right) $ be a
Hilbert space and $p_{i}\geq 0$ $\left( i\in \mathbb{N}\right) $ with $%
\sum_{i=1}^{\infty }p_{i}=1.$ Assume that $\mathbf{x},\mathbf{y}\in \ell
_{p}^{2}\left( K\right) $ and $r_{1},r_{2}>0$ satisfy the condition 
\begin{equation}
\left\Vert x_{i}-y_{i}\right\Vert \geq r_{2}\geq r_{1}\geq \left\vert
\left\Vert x_{i}\right\Vert -\left\Vert y_{i}\right\Vert \right\vert
\label{ch2.3.1}
\end{equation}%
for each $i\in \mathbb{N}$. Then we have the following refinement of the
Cauchy-Bunyakovsky-Schwarz inequality 
\begin{multline}
\left( \sum_{i=1}^{\infty }p_{i}\left\Vert x_{i}\right\Vert
^{2}\sum_{i=1}^{\infty }p_{i}\left\Vert y_{i}\right\Vert ^{2}\right) ^{\frac{%
1}{2}}-\sum_{i=1}^{\infty }p_{i}\func{Re}\left( x_{i},y_{i}\right)
\label{ch2.3.2} \\
\geq \frac{1}{2}\left( r_{2}^{2}-r_{1}^{2}\right) \geq 0.
\end{multline}%
The constant $\frac{1}{2}$ is best possible.
\end{proposition}

\begin{proof}
From the condition (\ref{ch2.3.1}) we simply deduce 
\begin{align}
\sum_{i=1}^{\infty }p_{i}\left\Vert x_{i}-y_{i}\right\Vert ^{2}& \geq
r_{2}^{2}\geq r_{1}^{2}\geq \sum_{i=1}^{\infty }p_{i}\left( \left\Vert
x_{i}\right\Vert -\left\Vert y_{i}\right\Vert \right) ^{2}  \label{ch2.3.3}
\\
& \geq \left[ \left( \sum_{i=1}^{\infty }p_{i}\left\Vert x_{i}\right\Vert
^{2}\right) ^{\frac{1}{2}}-\left( \sum_{i=1}^{\infty }p_{i}\left\Vert
y_{i}\right\Vert ^{2}\right) ^{\frac{1}{2}}\right] ^{2}.  \notag
\end{align}
In terms of the norm $\left\Vert \cdot \right\Vert _{p},$ the inequality (%
\ref{ch2.3.3}) may be written as 
\begin{equation}
\left\Vert \mathbf{x}-\mathbf{y}\right\Vert _{p}\geq r_{2}\geq r_{1}\geq
\left\vert \left\Vert \mathbf{x}\right\Vert _{p}-\left\Vert \mathbf{y}%
\right\Vert _{p}\right\vert .  \label{ch2.3.4}
\end{equation}%
Utilising Theorem \ref{ch2.t2.1} for the Hilbert space $\left( \ell
_{p}^{2}\left( K\right) ,\left\langle \cdot ,\cdot \right\rangle _{p}\right)
,$ we deduce the desired inequality (\ref{ch2.3.2}).

For $n=1$ $\left( p_{1}=1\right) ,$ the inequality (\ref{ch2.3.2}) reduces
to (\ref{ch2.2.2}) for which we have shown that $\frac{1}{2}$ is the best
possible constant.
\end{proof}

By the use of Corollary \ref{ch2.c2.2}, we may state the following result as
well.

\begin{corollary}
\label{ch2.c3.2}With the assumptions of Proposition \ref{ch2.p3.1}, we have
the inequality 
\begin{multline}
\left( \sum_{i=1}^{\infty }p_{i}\left\Vert x_{i}\right\Vert ^{2}\right) ^{%
\frac{1}{2}}+\left( \sum_{i=1}^{\infty }p_{i}\left\Vert y_{i}\right\Vert
^{2}\right) ^{\frac{1}{2}}  \label{ch2.3.5} \\
-\frac{\sqrt{2}}{2}\left( \sum_{i=1}^{\infty }p_{i}\left\Vert
x_{i}+y_{i}\right\Vert ^{2}\right) ^{\frac{1}{2}}\geq \frac{\sqrt{2}}{2}%
\sqrt{r_{2}^{2}-r_{1}^{2}}.
\end{multline}
\end{corollary}

The following proposition also holds \cite{DRAG01}.

\begin{proposition}
\label{ch2.p3.3}Let $\left( K;\left( \cdot ,\cdot \right) \right) $ be a
Hilbert space and $p_{i}\geq 0$ $\left( i\in \mathbb{N}\right) $ with $%
\sum_{i=1}^{\infty }p_{i}=1.$ Assume that $\mathbf{x},\mathbf{y}\in \ell
_{p}^{2}\left( K\right) $ and $R\geq 1,$ $r\geq 0$ satisfy the condition 
\begin{equation}
\frac{1}{R}\left( \left\Vert x_{i}\right\Vert +\left\Vert y_{i}\right\Vert
\right) \geq \left\Vert x_{i}+y_{i}\right\Vert \geq r  \label{ch2.3.6}
\end{equation}%
for each $i\in \mathbb{N}$. Then we have the following refinement of the
Schwarz inequality 
\begin{multline}
\left( \sum_{i=1}^{\infty }p_{i}\left\Vert x_{i}\right\Vert
^{2}\sum_{i=1}^{\infty }p_{i}\left\Vert y_{i}\right\Vert ^{2}\right) ^{\frac{%
1}{2}}-\sum_{i=1}^{\infty }p_{i}\func{Re}\left( x_{i},y_{i}\right)
\label{ch2.3.7} \\
\geq \frac{1}{2}\left( R^{2}-1\right) r^{2}.
\end{multline}%
The constant $\frac{1}{2}$ is best possible in the sense that it cannot be
replaced by a larger quantity.
\end{proposition}

\begin{proof}
By (\ref{ch2.3.6}) we deduce 
\begin{equation}
\frac{1}{R}\left[ \sum_{i=1}^{\infty }p_{i}\left( \left\Vert
x_{i}\right\Vert +\left\Vert y_{i}\right\Vert \right) ^{2}\right] ^{\frac{1}{%
2}}\geq \left( \sum_{i=1}^{\infty }p_{i}\left\Vert x_{i}+y_{i}\right\Vert
^{2}\right) ^{\frac{1}{2}}\geq r.  \label{ch2.3.8}
\end{equation}%
By the classical Minkowsky inequality for nonnegative numbers, we have 
\begin{multline}
\quad \left( \sum_{i=1}^{\infty }p_{i}\left\Vert x_{i}\right\Vert
^{2}\right) ^{\frac{1}{2}}+\left( \sum_{i=1}^{\infty }p_{i}\left\Vert
y_{i}\right\Vert ^{2}\right) ^{\frac{1}{2}}  \label{ch2.3.9} \\
\geq \left[ \sum_{i=1}^{\infty }p_{i}\left( \left\Vert x_{i}\right\Vert
+\left\Vert y_{i}\right\Vert \right) ^{2}\right] ^{\frac{1}{2}},\quad
\end{multline}%
and thus, by utilising (\ref{ch2.3.8}) and (\ref{ch2.3.9}), we may state in
terms of $\left\Vert \cdot \right\Vert _{p}$ the following inequality 
\begin{equation}
\frac{1}{R}\left( \left\Vert \mathbf{x}\right\Vert _{p}+\left\Vert \mathbf{y}%
\right\Vert _{p}\right) \geq \left\Vert \mathbf{x}+\mathbf{y}\right\Vert
_{p}\geq r.  \label{ch2.3.10}
\end{equation}%
Employing Theorem \ref{ch2.t2.6} for the Hilbert space $\ell _{p}^{2}\left(
K\right) $ and the inequality (\ref{ch2.3.10}), we deduce the desired result
(\ref{ch2.3.7}).

Since, for $p=1,$ $n=1,$ (\ref{ch2.3.7}) reduced to (\ref{ch2.2.23}) for
which we have shown that $\frac{1}{2}$ is the best constant, we conclude
that $\frac{1}{2}$ is the best constant in (\ref{ch2.3.7}) as well.
\end{proof}

Finally, we may state and prove the following result \cite{DRAG01}
incorporated in

\begin{proposition}
\label{ch2.p3.4}Let $\left( K;\left( \cdot ,\cdot \right) \right) $ be a
Hilbert space and $p_{i}\geq 0$ $\left( i\in \mathbb{N}\right) $ with $%
\sum_{i=1}^{\infty }p_{i}=1.$ Assume that $\mathbf{x},\mathbf{y}\in \ell
_{p}^{2}\left( K\right) $ and $r\in (0,1]$ such that 
\begin{equation}
\left\vert \left\Vert x_{i}\right\Vert -\left\Vert y_{i}\right\Vert
\right\vert \leq r\left\Vert x_{i}-y_{i}\right\Vert \text{ \ for each }i\in 
\mathbb{N},  \label{ch2.3.11}
\end{equation}
holds true. Then we have the following refinement of the Schwarz inequality 
\begin{multline}
\left( \sum_{i=1}^{\infty }p_{i}\left\Vert x_{i}\right\Vert
^{2}\sum_{i=1}^{\infty }p_{i}\left\Vert y_{i}\right\Vert ^{2}\right) ^{\frac{%
1}{2}}-\sum_{i=1}^{\infty }p_{i}\func{Re}\left( x_{i},y_{i}\right)
\label{ch2.3.12} \\
\geq \frac{1}{2}\left( 1-r^{2}\right) \sum_{i=1}^{\infty }p_{i}\left\Vert
x_{i}-y_{i}\right\Vert ^{2}.
\end{multline}
The constant $\frac{1}{2}$ is best possible in (\ref{ch2.3.12}).
\end{proposition}

\begin{proof}
From (\ref{ch2.3.11}) we have 
\begin{equation*}
\left[ \sum_{i=1}^{\infty }p_{i}\left( \left\Vert x_{i}\right\Vert
-\left\Vert y_{i}\right\Vert \right) ^{2}\right] ^{\frac{1}{2}}\leq r\left[
\sum_{i=1}^{\infty }p_{i}\left\Vert x_{i}-y_{i}\right\Vert ^{2}\right] ^{%
\frac{1}{2}}.
\end{equation*}%
Utilising the following elementary result 
\begin{equation*}
\left\vert \left( \sum_{i=1}^{\infty }p_{i}\left\Vert x_{i}\right\Vert
^{2}\right) ^{\frac{1}{2}}-\left( \sum_{i=1}^{\infty }p_{i}\left\Vert
y_{i}\right\Vert ^{2}\right) ^{\frac{1}{2}}\right\vert \leq \left(
\sum_{i=1}^{\infty }p_{i}\left( \left\Vert x_{i}\right\Vert -\left\Vert
y_{i}\right\Vert \right) ^{2}\right) ^{\frac{1}{2}},
\end{equation*}%
we may state that 
\begin{equation*}
\left\vert \left\Vert \mathbf{x}\right\Vert _{p}-\left\Vert \mathbf{y}%
\right\Vert _{p}\right\vert \leq r\left\Vert \mathbf{x}-\mathbf{y}%
\right\Vert _{p}.
\end{equation*}%
Now, by making use of Theorem \ref{ch2.t2.7}, we deduce the desired
inequality (\ref{ch2.3.12}) and the fact that $\frac{1}{2}$ is the best
possible constant. We omit the details.
\end{proof}

\subsection{Integral Inequalities}

Assume that $\left( K;\left( \cdot ,\cdot \right) \right) $ is a Hilbert
space over the real or complex number field $\mathbb{K}$. If $\rho :\left[
a,b\right] \subset \mathbb{R}\rightarrow \lbrack 0,\infty )$ is a Lebesgue
integrable function with $\int_{a}^{b}\rho \left( t\right) dt=1,$ then we
may consider the space $L_{\rho }^{2}\left( \left[ a,b\right] ;K\right) $ of
all functions $f:\left[ a,b\right] \rightarrow K,$ that are Bochner
measurable and $\int_{a}^{b}\rho \left( t\right) \left\Vert f\left( t\right)
\right\Vert ^{2}dt<\infty $. It is known that $L_{\rho }^{2}\left( \left[ a,b%
\right] ;K\right) $ endowed with the inner product $\left\langle \cdot
,\cdot \right\rangle _{\rho }$ defined by 
\begin{equation*}
\left\langle f,g\right\rangle _{\rho }:=\int_{a}^{b}\rho \left( t\right)
\left( f\left( t\right) ,g\left( t\right) \right) dt
\end{equation*}%
and generating the norm 
\begin{equation*}
\left\Vert f\right\Vert _{\rho }:=\left( \int_{a}^{b}\rho \left( t\right)
\left\Vert f\left( t\right) \right\Vert ^{2}dt\right) ^{\frac{1}{2}}
\end{equation*}%
is a Hilbert space over $\mathbb{K}$.

Now we may state and prove the first refinement of the
Cauchy-Bunyakovsky-Schwarz integral inequality \cite{DRAG01}.

\begin{proposition}
\label{ch2.p4.1}Assume that $f,g\in L_{\rho }^{2}\left( \left[ a,b\right]
;K\right) $ and $r_{2},r_{1}>0$ satisfy the condition 
\begin{equation}
\left\Vert f\left( t\right) -g\left( t\right) \right\Vert \geq r_{2}\geq
r_{1}\geq \left\vert \left\Vert f\left( t\right) \right\Vert -\left\Vert
g\left( t\right) \right\Vert \right\vert  \label{ch2.4.1}
\end{equation}%
for a.e. $t\in \left[ a,b\right] .$ Then we have the inequality 
\begin{multline}
\left( \int_{a}^{b}\rho \left( t\right) \left\Vert f\left( t\right)
\right\Vert ^{2}dt\int_{a}^{b}\rho \left( t\right) \left\Vert g\left(
t\right) \right\Vert ^{2}dt\right) ^{\frac{1}{2}}  \label{ch2.4.2} \\
-\int_{a}^{b}\rho \left( t\right) \func{Re}\left( f\left( t\right) ,g\left(
t\right) \right) dt\geq \frac{1}{2}\left( r_{2}^{2}-r_{1}^{2}\right) \left(
\geq 0\right) .
\end{multline}%
The constant $\frac{1}{2}$ is best possible in (\ref{ch2.4.2}).
\end{proposition}

\begin{proof}
Integrating (\ref{ch2.4.1}), we get 
\begin{multline}
\left( \int_{a}^{b}\rho \left( t\right) \left( \left\Vert f\left( t\right)
-g\left( t\right) \right\Vert \right) ^{2}dt\right) ^{\frac{1}{2}}
\label{ch2.4.3} \\
\geq r_{2}\geq r_{1}\geq \left( \int_{a}^{b}\rho \left( t\right) \left(
\left\Vert f\left( t\right) \right\Vert -\left\Vert g\left( t\right)
\right\Vert \right) ^{2}dt\right) ^{\frac{1}{2}}.
\end{multline}%
Utilising the obvious fact 
\begin{multline}
\left[ \int_{a}^{b}\rho \left( t\right) \left( \left\Vert f\left( t\right)
\right\Vert -\left\Vert g\left( t\right) \right\Vert \right) ^{2}dt\right] ^{%
\frac{1}{2}}  \label{ch2.4.4} \\
\geq \left\vert \left( \int_{a}^{b}\rho \left( t\right) \left\Vert f\left(
t\right) \right\Vert ^{2}dt\right) ^{\frac{1}{2}}-\left( \int_{a}^{b}\rho
\left( t\right) \left\Vert g\left( t\right) \right\Vert ^{2}dt\right) ^{%
\frac{1}{2}}\right\vert ,
\end{multline}
we can state the following inequality in terms of the $\left\Vert \cdot
\right\Vert _{\rho }$ norm: 
\begin{equation}
\left\Vert f-g\right\Vert _{\rho }\geq r_{2}\geq r_{1}\geq \left\vert
\left\Vert f\right\Vert _{\rho }-\left\Vert g\right\Vert _{\rho }\right\vert
.  \label{ch2.4.5}
\end{equation}%
Employing Theorem \ref{ch2.t2.1} for the Hilbert space $L_{\rho }^{2}\left( %
\left[ a,b\right] ;K\right) ,$ we deduce the desired inequality (\ref%
{ch2.4.2}).

To prove the sharpness of $\frac{1}{2}$ in (\ref{ch2.4.2}), we choose $a=0,$ 
$b=1,$ $f\left( t\right) =1,$ $t\in \left[ 0,1\right] $ and $f\left(
t\right) =x,$ $g\left( t\right) =y,$ $t\in \left[ a,b\right] ,$ $x,y\in K.$
Then (\ref{ch2.4.2}) becomes 
\begin{equation*}
\left\| x\right\| \left\| y\right\| -\func{Re}\left\langle x,y\right\rangle
\geq \frac{1}{2}\left( r_{2}^{2}-r_{1}^{2}\right)
\end{equation*}
provided 
\begin{equation*}
\left\| x-y\right\| \geq r_{2}\geq r_{1}\geq \left| \left\| x\right\|
-\left\| y\right\| \right| ,
\end{equation*}
which, by Theorem \ref{ch2.t2.1} has the quantity $\frac{1}{2}$ as the best
possible constant.
\end{proof}

The following corollary holds.

\begin{corollary}
\label{ch2.c4.2}With the assumptions of Proposition \ref{ch2.p4.1}, we have
the inequality 
\begin{multline}
\left( \int_{a}^{b}\rho \left( t\right) \left\Vert f\left( t\right)
\right\Vert ^{2}dt\right) ^{\frac{1}{2}}+\left( \int_{a}^{b}\rho \left(
t\right) \left\Vert g\left( t\right) \right\Vert ^{2}dt\right) ^{\frac{1}{2}}
\label{ch2.4.6} \\
-\frac{\sqrt{2}}{2}\left( \int_{a}^{b}\rho \left( t\right) \left\Vert
f\left( t\right) +g\left( t\right) \right\Vert ^{2}dt\right) ^{\frac{1}{2}%
}\geq \frac{\sqrt{2}}{2}\sqrt{r_{2}^{2}-r_{1}^{2}}.
\end{multline}
\end{corollary}

The following two refinements of the Cauchy-Bunyakovsky-Schwarz (CBS)
integral inequality also hold.

\begin{proposition}
\label{ch2.p4.3}If $f,g\in L_{\rho }^{2}\left( \left[ a,b\right] ;K\right) $
and $R\geq 1,r\geq 0$ satisfy the condition 
\begin{equation}
\frac{1}{R}\left( \left\Vert f\left( t\right) \right\Vert +\left\Vert
g\left( t\right) \right\Vert \right) \geq \left\Vert f\left( t\right)
+g\left( t\right) \right\Vert \geq r  \label{ch2.4.7}
\end{equation}%
for a.e. $t\in \left[ a,b\right] ,$ then we have the inequality 
\begin{multline}
\left( \int_{a}^{b}\rho \left( t\right) \left\Vert f\left( t\right)
\right\Vert ^{2}dt\int_{a}^{b}\rho \left( t\right) \left\Vert g\left(
t\right) \right\Vert ^{2}dt\right) ^{\frac{1}{2}}  \label{ch2.4.8} \\
-\int_{a}^{b}\rho \left( t\right) \func{Re}\left( f\left( t\right) ,g\left(
t\right) \right) dt\geq \frac{1}{2}\left( R^{2}-1\right) r^{2}.
\end{multline}%
The constant $\frac{1}{2}$ is best possible in (\ref{ch2.4.8}).
\end{proposition}

The proof follows by Theorem \ref{ch2.t2.6} and we omit the details.

\begin{proposition}
\label{ch2.p4.4}If $f,g\in L_{\rho }^{2}\left( \left[ a,b\right] ;K\right) $
and $\zeta \in (0,1]$ satisfy the condition 
\begin{equation}
\left\vert \left\Vert f\left( t\right) \right\Vert -\left\Vert g\left(
t\right) \right\Vert \right\vert \leq \zeta \left\Vert f\left( t\right)
-g\left( t\right) \right\Vert  \label{ch2.4.9}
\end{equation}%
for a.e. $t\in \left[ a,b\right] ,$ then we have the inequality 
\begin{multline}
\left( \int_{a}^{b}\rho \left( t\right) \left\Vert f\left( t\right)
\right\Vert ^{2}dt\int_{a}^{b}\rho \left( t\right) \left\Vert g\left(
t\right) \right\Vert ^{2}dt\right) ^{\frac{1}{2}}  \label{ch2.4.10} \\
-\int_{a}^{b}\rho \left( t\right) \func{Re}\left( f\left( t\right) ,g\left(
t\right) \right) dt \\
\geq \frac{1}{2}\left( 1-\zeta ^{2}\right) \int_{a}^{b}\rho \left( t\right)
\left\Vert f\left( t\right) -g\left( t\right) \right\Vert ^{2}dt.
\end{multline}%
The constant $\frac{1}{2}$ is best possible in (\ref{ch2.4.10}).
\end{proposition}

The proof follows by Theorem \ref{ch2.t2.7} and we omit the details.

\subsection{Refinements of the Heisenberg Inequality}

It is well known that if $\left( H;\left\langle \cdot ,\cdot \right\rangle
\right) $ is a real or complex Hilbert space and $f:\left[ a,b\right]
\subset \mathbb{R\rightarrow }H$ is an \textit{absolutely continuous
vector-valued }function, then $f$ is differentiable almost everywhere on $%
\left[ a,b\right] ,$ the derivative $f^{\prime }:\left[ a,b\right]
\rightarrow H$ is Bochner integrable on $\left[ a,b\right] $ and 
\begin{equation}
f\left( t\right) =\int_{a}^{t}f^{\prime }\left( s\right) ds\qquad \text{for
any \ }t\in \left[ a,b\right] .  \label{ch2.5.1}
\end{equation}

The following theorem provides a version of the Heisenberg inequalities in
the general setting of Hilbert spaces \cite{DRAG01}.

\begin{theorem}[Dragomir, 2004]
\label{ch2.t5.1}Let $\varphi :\left[ a,b\right] \rightarrow H$ be an
absolutely continuous function with the property that $b\left\Vert \varphi
\left( b\right) \right\Vert ^{2}=a\left\Vert \varphi \left( a\right)
\right\Vert ^{2}.$ Then we have the inequality: 
\begin{equation}
\left( \int_{a}^{b}\left\Vert \varphi \left( t\right) \right\Vert
^{2}dt\right) ^{2}\leq 4\int_{a}^{b}t^{2}\left\Vert \varphi \left( t\right)
\right\Vert ^{2}dt\cdot \int_{a}^{b}\left\Vert \varphi ^{\prime }\left(
t\right) \right\Vert ^{2}dt.  \label{ch2.5.2}
\end{equation}%
The constant $4$ is best possible in the sense that it cannot be replaced by
a smaller constant.
\end{theorem}

\begin{proof}
Integrating by parts, we have successively 
\begin{align}
& \int_{a}^{b}\left\Vert \varphi \left( t\right) \right\Vert ^{2}dt
\label{ch2.5.3} \\
& =t\left\Vert \varphi \left( t\right) \right\Vert ^{2}\bigg|%
_{a}^{b}-\int_{a}^{b}t\frac{d}{dt}\left( \left\Vert \varphi \left( t\right)
\right\Vert ^{2}\right) dt  \notag \\
& =b\left\Vert \varphi \left( b\right) \right\Vert ^{2}-a\left\Vert \varphi
\left( a\right) \right\Vert ^{2}-\int_{a}^{b}t\frac{d}{dt}\left\langle
\varphi \left( t\right) ,\varphi \left( t\right) \right\rangle dt  \notag \\
& =-\int_{a}^{b}t\left[ \left\langle \varphi ^{\prime }\left( t\right)
,\varphi \left( t\right) \right\rangle +\left\langle \varphi \left( t\right)
,\varphi ^{\prime }\left( t\right) \right\rangle \right] dt  \notag \\
& =-2\int_{a}^{b}t\func{Re}\left\langle \varphi ^{\prime }\left( t\right)
,\varphi \left( t\right) \right\rangle dt  \notag \\
& =2\int_{a}^{b}\func{Re}\left\langle \varphi ^{\prime }\left( t\right)
,\left( -t\right) \varphi \left( t\right) \right\rangle dt.  \notag
\end{align}%
If we apply the Cauchy-Bunyakovsky-Schwarz integral inequality 
\begin{equation*}
\int_{a}^{b}\func{Re}\left\langle g\left( t\right) ,h\left( t\right)
\right\rangle dt\leq \left( \int_{a}^{b}\left\Vert g\left( t\right)
\right\Vert ^{2}dt\int_{a}^{b}\left\Vert h\left( t\right) \right\Vert
^{2}dt\right) ^{\frac{1}{2}}
\end{equation*}%
for $g\left( t\right) =\varphi ^{\prime }\left( t\right) ,$ $h\left(
t\right) =-t\varphi \left( t\right) ,$ $t\in \left[ a,b\right] ,$ then we
deduce the desired inequality (\ref{ch2.4.2}).

The fact that $4$ is the best constant in (\ref{ch2.4.2}) follows from the
fact that in the (CBS) inequality, the case of equality holds iff $g\left(
t\right) =\lambda h\left( t\right) $ for a.e. $t\in \left[ a,b\right] $ and $%
\lambda $ a given scalar in $\mathbb{K}$. We omit the details.
\end{proof}

For details on the classical Heisenberg inequality, see, for instance, \cite%
{HLP}.

Utilising Proposition \ref{ch2.p4.1}, we can state the following refinement 
\cite{DRAG01} of the Heisenberg inequality obtained above in (\ref{ch2.5.2}):

\begin{proposition}
\label{ch2.p5.1}Assume that $\varphi :\left[ a,b\right] \rightarrow H$ is as
in the hypothesis of Theorem \ref{ch2.t5.1}. In addition, if there exists $%
r_{2},r_{1}>0$ so that 
\begin{equation*}
\left\Vert \varphi ^{\prime }\left( t\right) +t\varphi \left( t\right)
\right\Vert \geq r_{2}\geq r_{1}\geq \left\vert \left\Vert \varphi ^{\prime
}\left( t\right) \right\Vert -\left\vert t\right\vert \left\Vert \varphi
\left( t\right) \right\Vert \right\vert
\end{equation*}%
for a.e. $t\in \left[ a,b\right] ,$ then we have the inequality 
\begin{multline*}
\left( \int_{a}^{b}t^{2}\left\Vert \varphi \left( t\right) \right\Vert
^{2}dt\cdot \int_{a}^{b}\left\Vert \varphi ^{\prime }\left( t\right)
\right\Vert ^{2}dt\right) ^{\frac{1}{2}}-\frac{1}{2}\int_{a}^{b}\left\Vert
\varphi \left( t\right) \right\Vert ^{2}dt \\
\geq \frac{1}{2}\left( b-a\right) \left( r_{2}^{2}-r_{1}^{2}\right) \left(
\geq 0\right) .
\end{multline*}
\end{proposition}

The proof follows by Proposition \ref{ch2.p4.1} on choosing $f\left(
t\right) =\varphi ^{\prime }\left( t\right) ,g\left( t\right) =-t\varphi
\left( t\right) $ and $\rho \left( t\right) =\frac{1}{b-a},t\in \left[ a,b%
\right] .$

On utilising the Proposition \ref{ch2.p4.3} for the same choices of $f,g$
and $\rho ,$ we may state the following results as well \cite{DRAG01}:

\begin{proposition}
\label{ch2.p5.2}Assume that $\varphi :\left[ a,b\right] \rightarrow H$ is as
in the hypothesis of Theorem \ref{ch2.t5.1}. In addition, if there exists $%
R\geq 1$ and $r>0$ so that 
\begin{equation*}
\frac{1}{R}\left( \left\Vert \varphi ^{\prime }\left( t\right) \right\Vert
+\left\vert t\right\vert \left\Vert \varphi \left( t\right) \right\Vert
\right) \geq \left\Vert \varphi ^{\prime }\left( t\right) -t\varphi \left(
t\right) \right\Vert \geq r
\end{equation*}%
for a.e. $t\in \left[ a,b\right] ,$ then we have the inequality 
\begin{multline*}
\left( \int_{a}^{b}t^{2}\left\Vert \varphi \left( t\right) \right\Vert
^{2}dt\cdot \int_{a}^{b}\left\Vert \varphi ^{\prime }\left( t\right)
\right\Vert ^{2}dt\right) ^{\frac{1}{2}}-\frac{1}{2}\int_{a}^{b}\left\Vert
\varphi \left( t\right) \right\Vert ^{2}dt \\
\geq \frac{1}{2}\left( b-a\right) \left( R^{2}-1\right) r^{2}\left( \geq
0\right) .
\end{multline*}
\end{proposition}

Finally, we can state

\begin{proposition}
\label{ch2.p5.3}Let $\varphi :\left[ a,b\right] \rightarrow H$ be as in the
hypothesis of Theorem \ref{ch2.t5.1}. In addition, if there exists $\zeta
\in (0,1]$ so that 
\begin{equation*}
\left\vert \left\Vert \varphi ^{\prime }\left( t\right) \right\Vert
-\left\vert t\right\vert \left\Vert \varphi \left( t\right) \right\Vert
\right\vert \leq \zeta \left\Vert \varphi ^{\prime }\left( t\right)
+t\varphi \left( t\right) \right\Vert
\end{equation*}%
for a.e. $t\in \left[ a,b\right] ,$ then we have the inequality 
\begin{multline*}
\left( \int_{a}^{b}t^{2}\left\Vert \varphi \left( t\right) \right\Vert
^{2}dt\cdot \int_{a}^{b}\left\Vert \varphi ^{\prime }\left( t\right)
\right\Vert ^{2}dt\right) ^{\frac{1}{2}}-\frac{1}{2}\int_{a}^{b}\left\Vert
\varphi \left( t\right) \right\Vert ^{2}dt \\
\geq \frac{1}{2}\left( 1-\zeta ^{2}\right) \int_{a}^{b}\left\Vert \varphi
^{\prime }\left( t\right) +t\varphi \left( t\right) \right\Vert ^{2}dt\left(
\geq 0\right) .
\end{multline*}
\end{proposition}

This follows by Proposition \ref{ch2.p4.4} and we omit the details.

\section{More Schwarz Related Inequalities}

\subsection{Introduction}

In practice, one may need reverses of the Schwarz inequality, namely, upper
bounds for the quantities%
\begin{equation*}
\left\Vert x\right\Vert \left\Vert y\right\Vert -\func{Re}\left\langle
x,y\right\rangle ,\qquad \left\Vert x\right\Vert ^{2}\left\Vert y\right\Vert
^{2}-\left( \func{Re}\left\langle x,y\right\rangle \right) ^{2}
\end{equation*}%
and%
\begin{equation*}
\frac{\left\Vert x\right\Vert \left\Vert y\right\Vert }{\func{Re}%
\left\langle x,y\right\rangle }
\end{equation*}%
or the corresponding expressions where $\func{Re}\left\langle
x,y\right\rangle $ is replaced by either $\left\vert \func{Re}\left\langle
x,y\right\rangle \right\vert $ or $\left\vert \left\langle x,y\right\rangle
\right\vert ,$ under suitable assumptions for the vectors $x,y$ in an inner
product space $\left( H;\left\langle \cdot ,\cdot \right\rangle \right) $
over the real or complex number field $\mathbb{K}$.

In this class of results, we mention the following recent reverses of the
Schwarz inequality due to the present author, that can be found, for
instance, in the survey work \cite{SSD3a}, where more specific references
are provided:

\begin{theorem}[Dragomir, 2004]
\label{ch2x.t1.4}Let $\left( H;\left\langle \cdot ,\cdot \right\rangle
\right) $ be an inner product space over $\mathbb{K}$ $\left( \mathbb{K}=%
\mathbb{C},\mathbb{R}\right) .$ If $a,A\in \mathbb{K}$ and $x,y\in H$ are
such that either%
\begin{equation}
\func{Re}\left\langle Ay-x,x-ay\right\rangle \geq 0,  \label{ch2x.1.8}
\end{equation}%
or, equivalently,%
\begin{equation}
\left\Vert x-\frac{A+a}{2}y\right\Vert \leq \frac{1}{2}\left\vert
A-a\right\vert \left\Vert y\right\Vert ,  \label{ch2x.1.9}
\end{equation}%
then the following reverse for the quadratic form of the Schwarz inequality%
\begin{align}
(0& \leq )\left\Vert x\right\Vert ^{2}\left\Vert y\right\Vert
^{2}-\left\vert \left\langle x,y\right\rangle \right\vert ^{2}
\label{ch2x.1.10} \\
& \leq \left\{ 
\begin{array}{l}
\frac{1}{4}\left\vert A-a\right\vert ^{2}\left\Vert y\right\Vert
^{4}-\left\vert \frac{A+a}{2}\left\Vert y\right\Vert ^{2}-\left\langle
x,y\right\rangle \right\vert ^{2} \\ 
\\ 
\frac{1}{4}\left\vert A-a\right\vert ^{2}\left\Vert y\right\Vert
^{4}-\left\Vert y\right\Vert ^{2}\func{Re}\left\langle Ay-x,x-ay\right\rangle%
\end{array}%
\right.  \notag \\
& \leq \frac{1}{4}\left\vert A-a\right\vert ^{2}\left\Vert y\right\Vert ^{4}
\notag
\end{align}%
holds.

If in addition, we have $\func{Re}\left( A\bar{a}\right) >0,$ then%
\begin{equation}
\left\Vert x\right\Vert \left\Vert y\right\Vert \leq \frac{1}{2}\cdot \frac{%
\func{Re}\left[ \left( \bar{A}+\bar{a}\right) \left\langle x,y\right\rangle %
\right] }{\sqrt{\func{Re}\left( A\bar{a}\right) }}\leq \frac{1}{2}\cdot 
\frac{\left\vert A+a\right\vert }{\sqrt{\func{Re}\left( A\bar{a}\right) }}%
\left\vert \left\langle x,y\right\rangle \right\vert ,  \label{ch2x.1.11}
\end{equation}%
and%
\begin{equation}
(0\leq )\left\Vert x\right\Vert ^{2}\left\Vert y\right\Vert ^{2}-\left\vert
\left\langle x,y\right\rangle \right\vert ^{2}\leq \frac{1}{4}\cdot \frac{%
\left\vert A-a\right\vert ^{2}}{\func{Re}\left( A\bar{a}\right) }\left\vert
\left\langle x,y\right\rangle \right\vert ^{2}.  \label{ch2x.1.12}
\end{equation}%
Also, if (\ref{ch2x.1.8}) or (\ref{ch2x.1.9}) are valid and $A\neq -a,$ then
we have the reverse for the simple form of Schwarz inequality%
\begin{align}
(0& \leq )\left\Vert x\right\Vert \left\Vert y\right\Vert -\left\vert
\left\langle x,y\right\rangle \right\vert \leq \left\Vert x\right\Vert
\left\Vert y\right\Vert -\left\vert \func{Re}\left[ \frac{\bar{A}+\bar{a}}{%
\left\vert A+a\right\vert }\left\langle x,y\right\rangle \right] \right\vert
\label{ch2x.1.13} \\
& \leq \left\Vert x\right\Vert \left\Vert y\right\Vert -\func{Re}\left[ 
\frac{\bar{A}+\bar{a}}{\left\vert A+a\right\vert }\left\langle
x,y\right\rangle \right] \leq \frac{1}{4}\cdot \frac{\left\vert
A-a\right\vert ^{2}}{\left\vert A+a\right\vert }\left\Vert y\right\Vert ^{2}.
\notag
\end{align}%
The multiplicative constants $\frac{1}{4}$ and $\frac{1}{2}$ above are best
possible in the sense that they cannot be replaced by a smaller quantity.
\end{theorem}

For some classical results related to Schwarz inequality, see \cite{B}, \cite%
{FK}, \cite{PE}, \cite{P}, \cite{R} and the references therein.

The main aim of the present section is to point out other results in
connection with both the quadratic and simple forms of the Schwarz
inequality. As applications, some reverse results for the generalised
triangle inequality, i.e., upper bounds for the quantity%
\begin{equation*}
(0\leq )\sum_{i=1}^{n}\left\Vert x_{i}\right\Vert -\left\Vert
\sum_{i=1}^{n}x_{i}\right\Vert
\end{equation*}%
under various assumptions for the vectors $x_{i}\in H,$ $i\in \left\{
1,\dots ,n\right\} ,$ are established.

\subsection{Refinements and Reverses\label{ch2x.s2}}

The following result holds \cite{DRAG02}.

\begin{proposition}
\label{ch2x.p2.1}Let $\left( H;\left\langle \cdot ,\cdot \right\rangle
\right) $ be an inner product space over the real or complex number field $%
\mathbb{K}$. The subsequent statements are equivalent.

\begin{enumerate}
\item[(i)] The following inequality holds%
\begin{equation}
\left\Vert \frac{x}{\left\Vert x\right\Vert }-\frac{y}{\left\Vert
y\right\Vert }\right\Vert \leq \left( \geq \right) r;  \label{ch2x.2.1}
\end{equation}

\item[(ii)] The following reverse (improvement) of Schwarz's inequality holds%
\begin{equation}
\left\Vert x\right\Vert \left\Vert y\right\Vert -\func{Re}\left\langle
x,y\right\rangle \leq \left( \geq \right) \frac{1}{2}r^{2}\left\Vert
x\right\Vert \left\Vert y\right\Vert .  \label{ch2x.2.2}
\end{equation}%
The constant $\frac{1}{2}$ is best possible in (\ref{ch2x.2.2}) in the sense
that it cannot be replaced by a larger (smaller) quantity.
\end{enumerate}
\end{proposition}

\begin{remark}
\label{ch2x.r2.2}Since%
\begin{align*}
\left\Vert \left\Vert y\right\Vert x-\left\Vert x\right\Vert y\right\Vert &
=\left\Vert \left\Vert y\right\Vert \left( x-y\right) +\left( \left\Vert
y\right\Vert -\left\Vert x\right\Vert \right) y\right\Vert \\
& \leq \left\Vert y\right\Vert \left\Vert x-y\right\Vert +\left\vert
\left\Vert y\right\Vert -\left\Vert x\right\Vert \right\vert \left\Vert
y\right\Vert \\
& \leq 2\left\Vert y\right\Vert \left\Vert x-y\right\Vert
\end{align*}%
hence a sufficient condition for (\ref{ch2x.2.1}) to hold is%
\begin{equation}
\left\Vert x-y\right\Vert \leq \frac{r}{2}\left\Vert x\right\Vert .
\label{ch2x.2.3}
\end{equation}
\end{remark}

\begin{remark}
\label{ch2x.r2.3}Utilising the Dunkl-Williams inequality \cite{DW}%
\begin{equation}
\left\Vert a-b\right\Vert \geq \frac{1}{2}\left( \left\Vert a\right\Vert
+\left\Vert b\right\Vert \right) \left\Vert \frac{a}{\left\Vert a\right\Vert 
}-\frac{b}{\left\Vert b\right\Vert }\right\Vert ,\ \ a,b\in H\backslash
\left\{ 0\right\}  \label{ch2x.2.4}
\end{equation}%
with equality if and only if either $\left\Vert a\right\Vert =\left\Vert
b\right\Vert $ or $\left\Vert a\right\Vert +\left\Vert b\right\Vert
=\left\Vert a-b\right\Vert ,$ we can state the following inequality%
\begin{equation}
\frac{\left\Vert x\right\Vert \left\Vert y\right\Vert -\func{Re}\left\langle
x,y\right\rangle }{\left\Vert x\right\Vert \left\Vert y\right\Vert }\leq
2\left( \frac{\left\Vert x-y\right\Vert }{\left\Vert x\right\Vert
+\left\Vert y\right\Vert }\right) ^{2},\ \ x,y\in H\backslash \left\{
0\right\} .  \label{ch2x.2.5}
\end{equation}%
Obviously, if $x,y\in H\backslash \left\{ 0\right\} $ are such that%
\begin{equation}
\left\Vert x-y\right\Vert \leq \eta \left( \left\Vert x\right\Vert
+\left\Vert y\right\Vert \right) ,  \label{ch2x.2.6}
\end{equation}%
with $\eta \in (0,1],$ then one has the following reverse of the Schwarz
inequality%
\begin{equation}
\left\Vert x\right\Vert \left\Vert y\right\Vert -\func{Re}\left\langle
x,y\right\rangle \leq 2\eta ^{2}\left\Vert x\right\Vert \left\Vert
y\right\Vert  \label{ch2x.2.7}
\end{equation}%
that is similar to (\ref{ch2x.2.2}).
\end{remark}

The following result may be stated as well \cite{DRAG02}.

\begin{proposition}
\label{ch2x.t2.4}If $x,y\in H\backslash \left\{ 0\right\} $ and $\rho >0$
are such that%
\begin{equation}
\left\Vert \frac{x}{\left\Vert y\right\Vert }-\frac{y}{\left\Vert
x\right\Vert }\right\Vert \leq \rho ,  \label{ch2x.2.8}
\end{equation}%
then we have the following reverse of Schwarz's inequality%
\begin{align}
\left( 0\leq \right) \left\Vert x\right\Vert \left\Vert y\right\Vert
-\left\vert \left\langle x,y\right\rangle \right\vert & \leq \left\Vert
x\right\Vert \left\Vert y\right\Vert -\func{Re}\left\langle x,y\right\rangle
\label{ch2x.2.9} \\
& \leq \frac{1}{2}\rho ^{2}\left\Vert x\right\Vert \left\Vert y\right\Vert .
\notag
\end{align}%
The case of equality holds in the last inequality in (\ref{ch2x.2.9}) if and
only if%
\begin{equation}
\left\Vert x\right\Vert =\left\Vert y\right\Vert \qquad \text{and}\qquad
\left\Vert x-y\right\Vert =\rho .  \label{ch2x.2.10}
\end{equation}%
The constant $\frac{1}{2}$ in (\ref{ch2x.2.9}) cannot be replaced by a
smaller quantity.
\end{proposition}

\begin{proof}
Taking the square in (\ref{ch2x.2.8}), we get%
\begin{equation}
\frac{\left\Vert x\right\Vert ^{2}}{\left\Vert y\right\Vert ^{2}}-\frac{2%
\func{Re}\left\langle x,y\right\rangle }{\left\Vert x\right\Vert \left\Vert
y\right\Vert }+\frac{\left\Vert y\right\Vert ^{2}}{\left\Vert x\right\Vert
^{2}}\leq \rho ^{2}.  \label{ch2x.2.11}
\end{equation}%
Since, obviously%
\begin{equation}
2\leq \frac{\left\Vert x\right\Vert ^{2}}{\left\Vert y\right\Vert ^{2}}+%
\frac{\left\Vert y\right\Vert ^{2}}{\left\Vert x\right\Vert ^{2}}
\label{ch2x.2.12}
\end{equation}%
with equality iff $\left\Vert x\right\Vert =\left\Vert y\right\Vert ,$ hence
by (\ref{ch2x.2.11}) we deduce the second inequality in (\ref{ch2x.2.9}).

The case of equality and the best constant are obvious and we omit the
details.
\end{proof}

\begin{remark}
\label{ch2x.r2.5}In \cite{H}, Hile obtained the following inequality%
\begin{equation}
\left\Vert \left\Vert x\right\Vert ^{v}x-\left\Vert y\right\Vert
^{v}y\right\Vert \leq \frac{\left\Vert x\right\Vert ^{v+1}-\left\Vert
y\right\Vert ^{v+1}}{\left\Vert x\right\Vert -\left\Vert y\right\Vert }%
\left\Vert x-y\right\Vert  \label{ch2x.2.12a}
\end{equation}%
provided $v>0$ and $\left\Vert x\right\Vert \neq \left\Vert y\right\Vert .$

If in (\ref{ch2x.2.12a}) we choose $v=1$ and take the square, then we get%
\begin{equation}
\left\Vert x\right\Vert ^{4}-2\left\Vert x\right\Vert \left\Vert
y\right\Vert \func{Re}\left\langle x,y\right\rangle +\left\Vert y\right\Vert
^{4}\leq \left( \left\Vert x\right\Vert +\left\Vert y\right\Vert \right)
^{2}\left\Vert x-y\right\Vert ^{2}.  \label{ch2x.2.13}
\end{equation}%
Since,%
\begin{equation*}
\left\Vert x\right\Vert ^{4}+\left\Vert y\right\Vert ^{4}\geq 2\left\Vert
x\right\Vert ^{2}\left\Vert y\right\Vert ^{2},
\end{equation*}%
hence, by (\ref{ch2x.2.13}) we deduce%
\begin{equation}
\left( 0\leq \right) \left\Vert x\right\Vert \left\Vert y\right\Vert -\func{%
Re}\left\langle x,y\right\rangle \leq \frac{1}{2}\cdot \frac{\left(
\left\Vert x\right\Vert +\left\Vert y\right\Vert \right) ^{2}\left\Vert
x-y\right\Vert ^{2}}{\left\Vert x\right\Vert \left\Vert y\right\Vert },
\label{ch2x.2.13a}
\end{equation}%
provided $x,y\in H\backslash \left\{ 0\right\} .$
\end{remark}

The following inequality is due to Goldstein, Ryff and Clarke \cite[p. 309]%
{GRC}:%
\begin{multline}
\left\Vert x\right\Vert ^{2r}+\left\Vert y\right\Vert ^{2r}-2\left\Vert
x\right\Vert ^{r}\left\Vert y\right\Vert ^{r}\cdot \frac{\func{Re}%
\left\langle x,y\right\rangle }{\left\Vert x\right\Vert \left\Vert
y\right\Vert }  \label{ch2x.2.13.a} \\
\leq \left\{ 
\begin{array}{ll}
r^{2}\left\Vert x\right\Vert ^{2r-2}\left\Vert x-y\right\Vert ^{2} & \text{%
if \ }r\geq 1 \\ 
&  \\ 
\left\Vert y\right\Vert ^{2r-2}\left\Vert x-y\right\Vert ^{2} & \text{if \ }%
r<1%
\end{array}%
\right.
\end{multline}%
provided $r\in \mathbb{R}$ and $x,y\in H$ with $\left\Vert x\right\Vert \geq
\left\Vert y\right\Vert .$

Utilising (\ref{ch2x.2.13.a}) we may state the following proposition
containing a different reverse of the Schwarz inequality in inner product
spaces \cite{DRAG02}.

\begin{proposition}
\label{ch2x.p2.5.a}Let $\left( H;\left\langle \cdot ,\cdot \right\rangle
\right) $ be an inner product space over the real or complex number field $%
\mathbb{K}$. If $x,y\in H\backslash \left\{ 0\right\} $ and $\left\Vert
x\right\Vert \geq \left\Vert y\right\Vert ,$ then we have%
\begin{align}
0& \leq \left\Vert x\right\Vert \left\Vert y\right\Vert -\left\vert
\left\langle x,y\right\rangle \right\vert \leq \left\Vert x\right\Vert
\left\Vert y\right\Vert -\func{Re}\left\langle x,y\right\rangle
\label{ch2x.2.13.b} \\
& \leq \left\{ 
\begin{array}{ll}
\frac{1}{2}r^{2}\left( \frac{\left\Vert x\right\Vert }{\left\Vert
y\right\Vert }\right) ^{r-1}\left\Vert x-y\right\Vert ^{2} & \text{if \ }%
r\geq 1, \\ 
&  \\ 
\frac{1}{2}\left( \frac{\left\Vert x\right\Vert }{\left\Vert y\right\Vert }%
\right) ^{1-r}\left\Vert x-y\right\Vert ^{2} & \text{if \ }r<1.%
\end{array}%
\right.  \notag
\end{align}
\end{proposition}

\begin{proof}
It follows from (\ref{ch2x.2.13.a}), on dividing by $\left\Vert x\right\Vert
^{r}\left\Vert y\right\Vert ^{r},$ that%
\begin{multline}
\left( \frac{\left\Vert x\right\Vert }{\left\Vert y\right\Vert }\right)
^{r}+\left( \frac{\left\Vert y\right\Vert }{\left\Vert x\right\Vert }\right)
^{r}-2\cdot \frac{\func{Re}\left\langle x,y\right\rangle }{\left\Vert
x\right\Vert \left\Vert y\right\Vert }  \label{ch2x.2.13.c} \\
\leq \left\{ 
\begin{array}{ll}
r^{2}\cdot \frac{\left\Vert x\right\Vert ^{r-2}}{\left\Vert y\right\Vert ^{r}%
}\left\Vert x-y\right\Vert ^{2} & \text{if \ }r\geq 1, \\ 
&  \\ 
\frac{\left\Vert y\right\Vert ^{r-2}}{\left\Vert x\right\Vert ^{r}}%
\left\Vert x-y\right\Vert ^{2} & \text{if \ }r<1.%
\end{array}%
\right.
\end{multline}%
Since%
\begin{equation*}
\left( \frac{\left\Vert x\right\Vert }{\left\Vert y\right\Vert }\right)
^{r}+\left( \frac{\left\Vert y\right\Vert }{\left\Vert x\right\Vert }\right)
^{r}\geq 2,
\end{equation*}%
hence, by (\ref{ch2x.2.13.c}) one has%
\begin{equation*}
2-2\cdot \frac{\func{Re}\left\langle x,y\right\rangle }{\left\Vert
x\right\Vert \left\Vert y\right\Vert }\leq \left\{ 
\begin{array}{ll}
r^{2}\frac{\left\Vert x\right\Vert ^{r-2}}{\left\Vert y\right\Vert ^{r}}%
\left\Vert x-y\right\Vert ^{2} & \text{if \ }r\geq 1, \\ 
&  \\ 
\frac{\left\Vert y\right\Vert ^{r-2}}{\left\Vert x\right\Vert ^{r}}%
\left\Vert x-y\right\Vert ^{2} & \text{if \ }r<1.%
\end{array}%
\right.
\end{equation*}%
Dividing this inequality by 2 and multiplying with $\left\Vert x\right\Vert
\left\Vert y\right\Vert ,$ we deduce the desired result in (\ref{ch2x.2.13.b}%
).
\end{proof}

Another result providing a different additive reverse (refinement) of the
Schwarz inequality may be stated \cite{DRAG02}.

\begin{proposition}
\label{ch2x.p2.6}Let $x,y\in H$ with $y\neq 0$ and $r>0.$ The subsequent
statements are equivalent:

\begin{enumerate}
\item[(i)] The following inequality holds:%
\begin{equation}
\left\Vert x-\frac{\left\langle x,y\right\rangle }{\left\Vert y\right\Vert
^{2}}\cdot y\right\Vert \leq \left( \geq \right) r;  \label{ch2x.2.14}
\end{equation}

\item[(ii)] The following reverse (refinement) of the quadratic Schwarz
inequality holds:%
\begin{equation}
\left\Vert x\right\Vert ^{2}\left\Vert y\right\Vert ^{2}-\left\vert
\left\langle x,y\right\rangle \right\vert ^{2}\leq \left( \geq \right)
r^{2}\left\Vert y\right\Vert ^{2}.  \label{ch2x.2.15}
\end{equation}
\end{enumerate}
\end{proposition}

The proof is obvious on taking the square in (\ref{ch2x.2.14}) and
performing the calculation.

\begin{remark}
\label{ch2x.r2.7}Since%
\begin{align*}
\left\Vert \left\Vert y\right\Vert ^{2}x-\left\langle x,y\right\rangle
y\right\Vert & =\left\Vert \left\Vert y\right\Vert ^{2}\left( x-y\right)
-\left\langle x-y,y\right\rangle y\right\Vert \\
& \leq \left\Vert y\right\Vert ^{2}\left\Vert x-y\right\Vert +\left\vert
\left\langle x-y,y\right\rangle \right\vert \left\Vert y\right\Vert \\
& \leq 2\left\Vert x-y\right\Vert \left\Vert y\right\Vert ^{2},
\end{align*}%
hence a sufficient condition for the inequality (\ref{ch2x.2.14}) to hold is
that%
\begin{equation}
\left\Vert x-y\right\Vert \leq \frac{r}{2}.  \label{ch2x.2.16}
\end{equation}
\end{remark}

The following proposition may give a complementary approach \cite{DRAG02}:

\begin{proposition}
\label{ch2x.p2.6.a}Let $x,y\in H$ with $\left\langle x,y\right\rangle \neq 0$
and $\rho >0.$ If 
\begin{equation}
\left\Vert x-\frac{\left\langle x,y\right\rangle }{\left\vert \left\langle
x,y\right\rangle \right\vert }\cdot y\right\Vert \leq \rho ,
\label{ch2x.2.16.a}
\end{equation}%
then%
\begin{equation}
(0\leq )\left\Vert x\right\Vert \left\Vert y\right\Vert -\left\vert
\left\langle x,y\right\rangle \right\vert \leq \frac{1}{2}\rho ^{2}.
\label{ch2x.2.16.b}
\end{equation}%
The multiplicative constant $\frac{1}{2}$ is best possible in (\ref%
{ch2x.2.16.b}).
\end{proposition}

The proof is similar to the ones outlined above and we omit it.

For the case of complex inner product spaces, we may state the following
result \cite{DRAG02}.

\begin{proposition}
\label{ch2x.t2.8}Let $\left( H;\left\langle \cdot ,\cdot \right\rangle
\right) $ be a complex inner product space and $\alpha \in \mathbb{C}$ a
given complex number with $\func{Re}\alpha ,$ $\func{Im}\alpha >0.$ If $%
x,y\in H$ are such that%
\begin{equation}
\left\Vert x-\frac{\func{Im}\alpha }{\func{Re}\alpha }\cdot y\right\Vert
\leq r,  \label{ch2x.2.17}
\end{equation}%
then we have the inequality%
\begin{align}
(0& \leq )\left\Vert x\right\Vert \left\Vert y\right\Vert -\left\vert
\left\langle x,y\right\rangle \right\vert \leq \left\Vert x\right\Vert
\left\Vert y\right\Vert -\func{Re}\left\langle x,y\right\rangle
\label{ch2x.2.18} \\
& \leq \frac{1}{2}\cdot \frac{\func{Re}\alpha }{\func{Im}\alpha }\cdot r^{2}.
\notag
\end{align}%
The equality holds in the second inequality in (\ref{ch2x.2.18}) if and only
if the case of equality holds in (\ref{ch2x.2.17}) and $\func{Re}\alpha
\cdot \left\Vert x\right\Vert =\func{Im}\alpha \cdot \left\Vert y\right\Vert
.$
\end{proposition}

\begin{proof}
Observe that the condition (\ref{ch2x.2.17}) is equivalent to%
\begin{equation}
\left[ \func{Re}\alpha \right] ^{2}\left\Vert x\right\Vert ^{2}+\left[ \func{%
Im}\alpha \right] ^{2}\left\Vert y\right\Vert ^{2}\leq 2\func{Re}\alpha 
\func{Im}\alpha \func{Re}\left\langle x,y\right\rangle +\left[ \func{Re}%
\alpha \right] ^{2}r^{2}.  \label{ch2x.2.19}
\end{equation}%
On the other hand, on utilising the elementary inequality%
\begin{equation}
2\func{Re}\alpha \func{Im}\alpha \left\Vert x\right\Vert \left\Vert
y\right\Vert \leq \left[ \func{Re}\alpha \right] ^{2}\left\Vert x\right\Vert
^{2}+\left[ \func{Im}\alpha \right] ^{2}\left\Vert y\right\Vert ^{2},
\label{ch2x.2.20}
\end{equation}%
with equality if and only if $\func{Re}\alpha \cdot \left\Vert x\right\Vert =%
\func{Im}\alpha \cdot \left\Vert y\right\Vert ,$ we deduce from (\ref%
{ch2x.2.19}) that%
\begin{equation}
2\func{Re}\alpha \func{Im}\alpha \left\Vert x\right\Vert \left\Vert
y\right\Vert \leq 2\func{Re}\alpha \func{Im}\alpha \func{Re}\left\langle
x,y\right\rangle +r^{2}\left[ \func{Re}\alpha \right] ^{2}  \label{ch2x.2.21}
\end{equation}%
giving the desired inequality (\ref{ch2x.2.18}).

The case of equality follows from the above and we omit the details.
\end{proof}

The following different reverse for the Schwarz inequality that holds for
both real and complex inner product spaces may be stated as well \cite%
{DRAG02}.

\begin{theorem}[Dragomir, 2004]
\label{ch2x.t2.9}Let $\left( H;\left\langle \cdot ,\cdot \right\rangle
\right) $ be an inner product space over $\mathbb{K}$, $\mathbb{K}=\mathbb{C}%
,\mathbb{R}.$ If $\alpha \in \mathbb{K}\backslash \left\{ 0\right\} ,$ then%
\begin{align}
0& \leq \left\Vert x\right\Vert \left\Vert y\right\Vert -\left\vert
\left\langle x,y\right\rangle \right\vert \leq \left\Vert x\right\Vert
\left\Vert y\right\Vert -\func{Re}\left[ \frac{\alpha ^{2}}{\left\vert
\alpha \right\vert ^{2}}\left\langle x,y\right\rangle \right]
\label{ch2x.2.22} \\
& \leq \frac{1}{2}\cdot \frac{\left[ \left\vert \func{Re}\alpha \right\vert
\left\Vert x-y\right\Vert +\left\vert \func{Im}\alpha \right\vert \left\Vert
x+y\right\Vert \right] ^{2}}{\left\vert \alpha \right\vert ^{2}}\leq \frac{1%
}{2}\cdot I^{2},  \notag
\end{align}%
where%
\begin{equation}
I:=\left\{ 
\begin{array}{l}
\max \left\{ \left\vert \func{Re}\alpha \right\vert ,\left\vert \func{Im}%
\alpha \right\vert \right\} \left( \left\Vert x-y\right\Vert +\left\Vert
x+y\right\Vert \right) ; \\ 
\\ 
\left( \left\vert \func{Re}\alpha \right\vert ^{p}+\left\vert \func{Im}%
\alpha \right\vert ^{p}\right) ^{\frac{1}{p}}\left( \left\Vert
x-y\right\Vert ^{q}+\left\Vert x+y\right\Vert ^{q}\right) ^{\frac{1}{q}}, \\ 
\hfill p>1,\ \frac{1}{p}+\frac{1}{q}=1; \\ 
\max \left\{ \left\Vert x-y\right\Vert ,\left\Vert x+y\right\Vert \right\}
\left( \left\vert \func{Re}\alpha \right\vert +\left\vert \func{Im}\alpha
\right\vert \right) .%
\end{array}%
\right.  \label{ch2x.2.23}
\end{equation}
\end{theorem}

\begin{proof}
Observe, for $\alpha \in \mathbb{K}\backslash \left\{ 0\right\} ,$ that%
\begin{align*}
\left\Vert \alpha x-\bar{\alpha}y\right\Vert ^{2}& =\left\vert \alpha
\right\vert ^{2}\left\Vert x\right\Vert ^{2}-2\func{Re}\left\langle \alpha x,%
\bar{\alpha}y\right\rangle +\left\vert \alpha \right\vert ^{2}\left\Vert
y\right\Vert ^{2} \\
& =\left\vert \alpha \right\vert ^{2}\left( \left\Vert x\right\Vert
^{2}+\left\Vert y\right\Vert ^{2}\right) -2\func{Re}\left[ \alpha
^{2}\left\langle x,y\right\rangle \right] .
\end{align*}%
Since $\left\Vert x\right\Vert ^{2}+\left\Vert y\right\Vert ^{2}\geq
2\left\Vert x\right\Vert \left\Vert y\right\Vert ,$ hence%
\begin{equation}
\left\Vert \alpha x-\bar{\alpha}y\right\Vert ^{2}\geq 2\left\vert \alpha
\right\vert ^{2}\left\{ \left\Vert x\right\Vert \left\Vert y\right\Vert -%
\func{Re}\left[ \frac{\alpha ^{2}}{\left\vert \alpha \right\vert ^{2}}%
\left\langle x,y\right\rangle \right] \right\} .  \label{ch2x.2.24}
\end{equation}%
On the other hand, we have%
\begin{align}
\left\Vert \alpha x-\bar{\alpha}y\right\Vert & =\left\Vert \left( \func{Re}%
\alpha +i\func{Im}\alpha \right) x-\left( \func{Re}\alpha -i\func{Im}\alpha
\right) y\right\Vert  \label{ch2x.2.25} \\
& =\left\Vert \func{Re}\alpha \left( x-y\right) +i\func{Im}\alpha \left(
x+y\right) \right\Vert  \notag \\
& \leq \left\vert \func{Re}\alpha \right\vert \left\Vert x-y\right\Vert
+\left\vert \func{Im}\alpha \right\vert \left\Vert x+y\right\Vert .  \notag
\end{align}%
Utilising (\ref{ch2x.2.24}) and (\ref{ch2x.2.25}) we deduce the third
inequality in (\ref{ch2x.2.22}).

For the last inequality we use the following elementary inequality%
\begin{equation}
\alpha a+\beta b\leq \left\{ 
\begin{array}{ll}
\max \left\{ \alpha ,\beta \right\} \left( a+b\right) &  \\ 
&  \\ 
\left( \alpha ^{p}+\beta ^{p}\right) ^{\frac{1}{p}}\left( a^{q}+b^{q}\right)
^{\frac{1}{q}}, & p>1,\ \frac{1}{p}+\frac{1}{q}=1,%
\end{array}%
\right.  \label{ch2x.2.26}
\end{equation}%
provided $\alpha ,\beta ,a,b\geq 0.$
\end{proof}

The following result may be stated \cite{DRAG02}.

\begin{proposition}
\label{ch2x.p2.11}Let $\left( H;\left\langle \cdot ,\cdot \right\rangle
\right) $ be an inner product over$\mathbb{\ K}$ and $e\in H,$ $\left\Vert
e\right\Vert =1.$ If $\lambda \in \left( 0,1\right) ,$ then%
\begin{multline}
\func{Re}\left[ \left\langle x,y\right\rangle -\left\langle x,e\right\rangle
\left\langle e,y\right\rangle \right]  \label{ch2x.2.28} \\
\leq \frac{1}{4}\cdot \frac{1}{\lambda \left( 1-\lambda \right) }\left[
\left\Vert \lambda x+\left( 1-\lambda \right) y\right\Vert ^{2}-\left\vert
\left\langle \lambda x+\left( 1-\lambda \right) y,e\right\rangle \right\vert
^{2}\right] .
\end{multline}%
The constant $\frac{1}{4}$ is best possible.
\end{proposition}

\begin{proof}
Firstly, note that the following equality holds true%
\begin{equation*}
\left\langle x-\left\langle x,e\right\rangle e,y-\left\langle
y,e\right\rangle e\right\rangle =\left\langle x,y\right\rangle -\left\langle
x,e\right\rangle \left\langle e,y\right\rangle .
\end{equation*}%
Utilising the elementary inequality%
\begin{equation*}
\func{Re}\left\langle z,w\right\rangle \leq \frac{1}{4}\left\Vert
z+w\right\Vert ^{2},\qquad z,w\in H
\end{equation*}%
we have%
\begin{align*}
& \func{Re}\left\langle x-\left\langle x,e\right\rangle e,y-\left\langle
y,e\right\rangle e\right\rangle \\
& =\frac{1}{\lambda \left( 1-\lambda \right) }\func{Re}\left\langle \lambda
x-\left\langle \lambda x,e\right\rangle e,\left( 1-\lambda \right)
y-\left\langle \left( 1-\lambda \right) y,e\right\rangle e\right\rangle \\
& \leq \frac{1}{4}\cdot \frac{1}{\lambda \left( 1-\lambda \right) }\left[
\left\Vert \lambda x+\left( 1-\lambda \right) y\right\Vert ^{2}-\left\vert
\left\langle \lambda x+\left( 1-\lambda \right) y,e\right\rangle \right\vert
^{2}\right] ,
\end{align*}%
proving the desired inequality (\ref{ch2x.2.28}).
\end{proof}

\begin{remark}
\label{ch2x.r2.12}For $\lambda =\frac{1}{2},$ we get the simpler inequality:%
\begin{equation}
\func{Re}\left[ \left\langle x,y\right\rangle -\left\langle x,e\right\rangle
\left\langle e,y\right\rangle \right] \leq \left\Vert \frac{x+y}{2}%
\right\Vert ^{2}-\left\vert \left\langle \frac{x+y}{2},e\right\rangle
\right\vert ^{2},  \label{ch2x.2.29}
\end{equation}%
that has been obtained in \cite[p. 46]{SSD3a}, for which the sharpness of
the inequality was established.
\end{remark}

The following result may be stated as well \cite{DRAG02}.

\begin{theorem}[Dragomir, 2004]
\label{ch2x.t2.13}Let $\left( H;\left\langle \cdot ,\cdot \right\rangle
\right) $ be an inner product space over $\mathbb{K}$ and $p\geq 1.$ Then
for any $x,y\in H$ we have%
\begin{align}
0& \leq \left\Vert x\right\Vert \left\Vert y\right\Vert -\left\vert
\left\langle x,y\right\rangle \right\vert \leq \left\Vert x\right\Vert
\left\Vert y\right\Vert -\func{Re}\left\langle x,y\right\rangle
\label{ch2x.2.30} \\
& \leq \frac{1}{2}\times \left\{ 
\begin{array}{l}
\left[ \left( \left\Vert x\right\Vert +\left\Vert y\right\Vert \right)
^{2p}-\left\Vert x+y\right\Vert ^{2p}\right] ^{\frac{1}{p}}, \\ 
\\ 
\left[ \left\Vert x-y\right\Vert ^{2p}-\left\vert \left\Vert x\right\Vert
-\left\Vert y\right\Vert \right\vert ^{2p}\right] ^{\frac{1}{p}}.%
\end{array}%
\right.  \notag
\end{align}
\end{theorem}

\begin{proof}
Firstly, observe that%
\begin{equation*}
2\left( \left\Vert x\right\Vert \left\Vert y\right\Vert -\func{Re}%
\left\langle x,y\right\rangle \right) =\left( \left\Vert x\right\Vert
+\left\Vert y\right\Vert \right) ^{2}-\left\Vert x+y\right\Vert ^{2}.
\end{equation*}%
Denoting $D:=\left\Vert x\right\Vert \left\Vert y\right\Vert -\func{Re}%
\left\langle x,y\right\rangle ,$ then we have%
\begin{equation}
2D+\left\Vert x+y\right\Vert ^{2}=\left( \left\Vert x\right\Vert +\left\Vert
y\right\Vert \right) ^{2}.  \label{ch2x.2.31}
\end{equation}%
Taking in (\ref{ch2x.2.31}) the power $p\geq 1$ and using the elementary
inequality 
\begin{equation}
\left( a+b\right) ^{p}\geq a^{p}+b^{p};a,b\geq 0,  \label{ch2x.2.31.a}
\end{equation}%
we have%
\begin{equation*}
\left( \left\Vert x\right\Vert +\left\Vert y\right\Vert \right) ^{2p}=\left(
2D+\left\Vert x+y\right\Vert ^{2}\right) ^{p}\geq 2^{p}D^{p}+\left\Vert
x+y\right\Vert ^{2p}
\end{equation*}%
giving%
\begin{equation*}
D^{p}\leq \frac{1}{2^{p}}\left[ \left( \left\Vert x\right\Vert +\left\Vert
y\right\Vert \right) ^{2p}-\left\Vert x+y\right\Vert ^{2p}\right] ,
\end{equation*}%
which is clearly equivalent to the first branch of the third inequality in (%
\ref{ch2x.2.30}).

With the above notation, we also have%
\begin{equation}
2D+\left( \left\Vert x\right\Vert -\left\Vert y\right\Vert \right)
^{2}=\left\Vert x-y\right\Vert ^{2}.  \label{ch2x.2.32}
\end{equation}%
Taking the power $p\geq 1$ in (\ref{ch2x.2.32}) and using the inequality (%
\ref{ch2x.2.31.a}) we deduce%
\begin{equation*}
\left\Vert x-y\right\Vert ^{2p}\geq 2^{p}D^{p}+\left\vert \left\Vert
x\right\Vert -\left\Vert y\right\Vert \right\vert ^{2p},
\end{equation*}%
from where we get the last part of (\ref{ch2x.2.30}).
\end{proof}

\subsection{More Schwarz Related Inequalities}

Before we point out other inequalities related to the Schwarz inequality, we
need the following identity that is interesting in itself \cite{DRAG02}.

\begin{lemma}[Dragomir, 2004]
\label{ch2x.l2.14}Let $\left( H;\left\langle \cdot ,\cdot \right\rangle
\right) $ be an inner product space over the real or complex number field $%
\mathbb{K}$, $e\in H,$ $\left\Vert e\right\Vert =1,$ $\alpha \in H$ and $%
\gamma ,\Gamma \in \mathbb{K}$. Then we have the identity:%
\begin{multline}
\left\Vert x\right\Vert ^{2}-\left\vert \left\langle x,e\right\rangle
\right\vert ^{2}  \label{ch2x.2.33} \\
=\left( \func{Re}\Gamma -\func{Re}\left\langle x,e\right\rangle \right)
\left( \func{Re}\left\langle x,e\right\rangle -\func{Re}\gamma \right)  \\
+\left( \func{Im}\Gamma -\func{Im}\left\langle x,e\right\rangle \right)
\left( \func{Im}\left\langle x,e\right\rangle -\func{Im}\gamma \right)  \\
+\left\Vert x-\frac{\gamma +\Gamma }{2}e\right\Vert ^{2}-\frac{1}{4}%
\left\vert \Gamma -\gamma \right\vert ^{2}.
\end{multline}
\end{lemma}

\begin{proof}
We start with the following known equality (see for instance \cite[eq. (2.6)]%
{SSD1a})%
\begin{multline}
\left\Vert x\right\Vert ^{2}-\left\vert \left\langle x,e\right\rangle
\right\vert ^{2}  \label{ch2x.2.34} \\
=\func{Re}\left[ \left( \Gamma -\left\langle x,e\right\rangle \right) \left( 
\overline{\left\langle x,e\right\rangle }-\bar{\gamma}\right) \right] -\func{%
Re}\left\langle \Gamma e-x,x-\gamma e\right\rangle
\end{multline}%
holding for $x\in H,$ $e\in H,$ $\left\Vert e\right\Vert =1$ and $\gamma
,\Gamma \in \mathbb{K}$.

We also know that (see for instance \cite{SSD1b})%
\begin{equation}
-\func{Re}\left\langle \Gamma e-x,x-\gamma e\right\rangle =\left\Vert x-%
\frac{\gamma +\Gamma }{2}e\right\Vert ^{2}-\frac{1}{4}\left\vert \Gamma
-\gamma \right\vert ^{2}.  \label{ch2x.2.35}
\end{equation}%
Since%
\begin{multline}
\func{Re}\left[ \left( \Gamma -\left\langle x,e\right\rangle \right) \left( 
\overline{\left\langle x,e\right\rangle }-\bar{\gamma}\right) \right] 
\label{ch2x.2.36} \\
=\left( \func{Re}\Gamma -\func{Re}\left\langle x,e\right\rangle \right)
\left( \func{Re}\left\langle x,e\right\rangle -\func{Re}\gamma \right)  \\
+\left( \func{Im}\Gamma -\func{Im}\left\langle x,e\right\rangle \right)
\left( \func{Im}\left\langle x,e\right\rangle -\func{Im}\gamma \right) ,
\end{multline}%
hence, by (\ref{ch2x.2.34}) -- (\ref{ch2x.2.36}), we deduce the desired
identity (\ref{ch2x.2.33}).
\end{proof}

The following general result providing a reverse of the Schwarz inequality
may be stated \cite{DRAG02}.

\begin{proposition}
\label{ch2x.t2.15}Let $\left( H;\left\langle \cdot ,\cdot \right\rangle
\right) $ be an inner product space over $\mathbb{K},$ $e\in H,$ $\left\Vert
e\right\Vert =1,$ $x\in H$ and $\gamma ,\Gamma \in \mathbb{K}$. Then we have
the inequality:%
\begin{equation}
\left( 0\leq \right) \left\Vert x\right\Vert ^{2}-\left\vert \left\langle
x,e\right\rangle \right\vert ^{2}\leq \left\Vert x-\frac{\gamma +\Gamma }{2}%
\cdot e\right\Vert ^{2}.  \label{ch2x.2.37}
\end{equation}%
The constant $\frac{1}{2}$ is best possible in (\ref{ch2x.2.37}). The case
of equality holds in (\ref{ch2x.2.37}) if and only if%
\begin{equation}
\func{Re}\left\langle x,e\right\rangle =\func{Re}\left( \frac{\gamma +\Gamma 
}{2}\right) ,\qquad \func{Im}\left\langle x,e\right\rangle =\func{Im}\left( 
\frac{\gamma +\Gamma }{2}\right) .  \label{ch2x.2.38}
\end{equation}
\end{proposition}

\begin{proof}
Utilising the elementary inequality for real numbers%
\begin{equation*}
\alpha \beta \leq \frac{1}{4}\left( \alpha +\beta \right) ^{2},\qquad \alpha
,\beta \in \mathbb{R};
\end{equation*}%
with equality iff $\alpha =\beta ,$ we have%
\begin{equation}
\left( \func{Re}\Gamma -\func{Re}\left\langle x,e\right\rangle \right)
\left( \func{Re}\left\langle x,e\right\rangle -\func{Re}\gamma \right) \leq 
\frac{1}{4}\left( \func{Re}\Gamma -\func{Re}\gamma \right) ^{2}
\label{ch2x.2.39}
\end{equation}%
and%
\begin{equation}
\left( \func{Im}\Gamma -\func{Im}\left\langle x,e\right\rangle \right)
\left( \func{Im}\left\langle x,e\right\rangle -\func{Im}\gamma \right) \leq 
\frac{1}{4}\left( \func{Im}\Gamma -\func{Im}\gamma \right) ^{2}
\label{ch2x.2.40}
\end{equation}%
with equality if and only if%
\begin{equation*}
\func{Re}\left\langle x,e\right\rangle =\frac{\func{Re}\Gamma +\func{Re}%
\gamma }{2}\qquad \text{and\qquad }\func{Im}\left\langle x,e\right\rangle =%
\frac{\func{Im}\Gamma +\func{Im}\gamma }{2}.
\end{equation*}%
Finally, on making use of (\ref{ch2x.2.39}), (\ref{ch2x.2.40}) and the
identity (\ref{ch2x.2.33}), we deduce the desired result (\ref{ch2x.2.37}).
\end{proof}

The following result may be stated as well \cite{DRAG02}.

\begin{proposition}
\label{ch2x.t2.16}Let $\left( H;\left\langle \cdot ,\cdot \right\rangle
\right) $ be an inner product space over $\mathbb{K},$ $e\in H,$ $\left\Vert
e\right\Vert =1,$ $x\in H$ and $\gamma ,\Gamma \in \mathbb{K}$. If $x\in H$
is such that%
\begin{equation}
\func{Re}\gamma \leq \func{Re}\left\langle x,e\right\rangle \leq \func{Re}%
\Gamma \qquad \text{and\qquad }\func{Im}\gamma \leq \func{Im}\left\langle
x,e\right\rangle \leq \func{Im}\Gamma ,  \label{ch2x.2.41}
\end{equation}%
then we have the inequality%
\begin{equation}
\left\Vert x\right\Vert ^{2}-\left\vert \left\langle x,e\right\rangle
\right\vert ^{2}\geq \left\Vert x-\frac{\gamma +\Gamma }{2}e\right\Vert ^{2}-%
\frac{1}{4}\left\vert \Gamma -\gamma \right\vert ^{2}.  \label{ch2x.2.42}
\end{equation}%
The constant $\frac{1}{4}$ is best possible in (\ref{ch2x.2.42}). The case
of equality holds in (\ref{ch2x.2.42}) if and only if%
\begin{equation*}
\func{Re}\left\langle x,e\right\rangle =\func{Re}\Gamma \text{ or }\func{Re}%
\left\langle x,e\right\rangle =\func{Re}\gamma 
\end{equation*}%
and%
\begin{equation*}
\func{Im}\left\langle x,e\right\rangle =\func{Im}\Gamma \text{ or }\func{Im}%
\left\langle x,e\right\rangle =\func{Im}\gamma .
\end{equation*}
\end{proposition}

\begin{proof}
From the hypothesis we obviously have%
\begin{equation*}
\left( \func{Re}\Gamma -\func{Re}\left\langle x,e\right\rangle \right)
\left( \func{Re}\left\langle x,e\right\rangle -\func{Re}\gamma \right) \geq 0
\end{equation*}%
and%
\begin{equation*}
\left( \func{Im}\Gamma -\func{Im}\left\langle x,e\right\rangle \right)
\left( \func{Im}\left\langle x,e\right\rangle -\func{Im}\gamma \right) \geq
0.
\end{equation*}%
Utilising the identity (\ref{ch2x.2.33}) we deduce the desired result (\ref%
{ch2x.2.42}). The case of equality is obvious.
\end{proof}

Further on, we can state the following reverse of the quadratic Schwarz
inequality \cite{DRAG02}:

\begin{proposition}
\label{ch2x.t2.17}Let $\left( H;\left\langle \cdot ,\cdot \right\rangle
\right) $ be an inner product space over $\mathbb{K},$ $e\in H,$ $\left\Vert
e\right\Vert =1.$ If $\gamma ,\Gamma \in \mathbb{K}$ and $x\in H$ are such
that either%
\begin{equation}
\func{Re}\left\langle \Gamma e-x,x-\gamma e\right\rangle \geq 0
\label{ch2x.2.43}
\end{equation}%
or, equivalently,%
\begin{equation}
\left\Vert x-\frac{\gamma +\Gamma }{2}e\right\Vert \leq \frac{1}{2}%
\left\vert \Gamma -\gamma \right\vert ,  \label{ch2x.2.44}
\end{equation}%
then%
\begin{align}
(0& \leq )\left\Vert x\right\Vert ^{2}-\left\vert \left\langle
x,e\right\rangle \right\vert ^{2}  \label{ch2x.2.45} \\
& \leq \left( \func{Re}\Gamma -\func{Re}\left\langle x,e\right\rangle
\right) \left( \func{Re}\left\langle x,e\right\rangle -\func{Re}\gamma
\right)  \notag \\
& \qquad \qquad \qquad +\left( \func{Im}\Gamma -\func{Im}\left\langle
x,e\right\rangle \right) \left( \func{Im}\left\langle x,e\right\rangle -%
\func{Im}\gamma \right)  \notag \\
& \leq \frac{1}{4}\left\vert \Gamma -\gamma \right\vert ^{2}.  \notag
\end{align}%
The case of equality holds in (\ref{ch2x.2.45}) if it holds either in (\ref%
{ch2x.2.43}) or (\ref{ch2x.2.44}).
\end{proposition}

The proof is obvious by Lemma \ref{ch2x.l2.14} and we omit the details.

\begin{remark}
\label{ch2x.r2.18}We remark that the inequality (\ref{ch2x.2.45}) may also
be used to get, for instance, the following result%
\begin{multline}
\left\Vert x\right\Vert ^{2}-\left\vert \left\langle x,e\right\rangle
\right\vert ^{2}  \label{ch2x.2.46} \\
\leq \left[ \left( \func{Re}\Gamma -\func{Re}\left\langle x,e\right\rangle
\right) ^{2}+\left( \func{Im}\Gamma -\func{Im}\left\langle x,e\right\rangle
\right) ^{2}\right] ^{\frac{1}{2}} \\
\times \left[ \left( \func{Re}\left\langle x,e\right\rangle -\func{Re}\gamma
\right) ^{2}+\left( \func{Im}\left\langle x,e\right\rangle -\func{Im}\gamma
\right) ^{2}\right] ^{\frac{1}{2}},
\end{multline}%
that provides a different bound than $\frac{1}{4}\left\vert \Gamma -\gamma
\right\vert ^{2}$ for the quantity $\left\Vert x\right\Vert ^{2}-\left\vert
\left\langle x,e\right\rangle \right\vert ^{2}.$
\end{remark}

The following result may be stated as well \cite{DRAG02}.

\begin{theorem}[Dragomir, 2004]
\label{ch2x.t2.19}Let $\left( H;\left\langle \cdot ,\cdot \right\rangle
\right) $ be an inner product space over $\mathbb{K}$ and $\alpha ,\gamma
>0, $ $\beta \in \mathbb{K}$ with $\left\vert \beta \right\vert ^{2}\geq
\alpha \gamma .$ If $x,a\in H$ are such that $a\neq 0$ and%
\begin{equation}
\left\Vert x-\frac{\beta }{\alpha }a\right\Vert \leq \frac{\left( \left\vert
\beta \right\vert ^{2}-\alpha \gamma \right) ^{\frac{1}{2}}}{\alpha }%
\left\Vert a\right\Vert ,  \label{ch2x.2.47}
\end{equation}%
then we have the following reverses of Schwarz's inequality%
\begin{align}
\left\Vert x\right\Vert \left\Vert a\right\Vert & \leq \frac{\func{Re}\beta
\cdot \func{Re}\left\langle x,a\right\rangle +\func{Im}\beta \cdot \func{Im}%
\left\langle x,a\right\rangle }{\sqrt{\alpha \gamma }}  \label{ch2x.2.48} \\
& \leq \frac{\left\vert \beta \right\vert \left\vert \left\langle
x,a\right\rangle \right\vert }{\sqrt{\alpha \gamma }}  \notag
\end{align}%
and%
\begin{equation}
\left( 0\leq \right) \left\Vert x\right\Vert ^{2}\left\Vert a\right\Vert
^{2}-\left\vert \left\langle x,a\right\rangle \right\vert ^{2}\leq \frac{%
\left\vert \beta \right\vert ^{2}-\alpha \gamma }{\alpha \gamma }\left\vert
\left\langle x,a\right\rangle \right\vert ^{2}.  \label{ch2x.2.49}
\end{equation}
\end{theorem}

\begin{proof}
Taking the square in (\ref{ch2x.2.47}), it becomes equivalent to%
\begin{equation*}
\left\Vert x\right\Vert ^{2}-\frac{2}{\alpha }\func{Re}\left[ \bar{\beta}%
\left\langle x,a\right\rangle \right] +\frac{\left\vert \beta \right\vert
^{2}}{\alpha ^{2}}\left\Vert a\right\Vert ^{2}\leq \frac{\left\vert \beta
\right\vert ^{2}-\alpha \gamma }{\alpha ^{2}}\left\Vert a\right\Vert ^{2},
\end{equation*}%
which is clearly equivalent to%
\begin{align}
\alpha \left\Vert x\right\Vert ^{2}+\gamma \left\Vert a\right\Vert ^{2}&
\leq 2\func{Re}\left[ \bar{\beta}\left\langle x,a\right\rangle \right]
\label{ch2x.2.50} \\
& =2\left[ \func{Re}\beta \cdot \func{Re}\left\langle x,a\right\rangle +%
\func{Im}\beta \cdot \func{Im}\left\langle x,a\right\rangle \right] .  \notag
\end{align}%
On the other hand, since%
\begin{equation}
2\sqrt{\alpha \gamma }\left\Vert x\right\Vert \left\Vert a\right\Vert \leq
\alpha \left\Vert x\right\Vert ^{2}+\gamma \left\Vert a\right\Vert ^{2},
\label{ch2x.2.51}
\end{equation}%
hence by (\ref{ch2x.2.50}) and (\ref{ch2x.2.51}) we deduce the first
inequality in (\ref{ch2x.2.48}).

The other inequalities are obvious.
\end{proof}

\begin{remark}
\label{ch2x.r2.20}The above inequality (\ref{ch2x.2.48}) contains in
particular the reverse (\ref{ch2x.1.11}) of the Schwarz inequality. Indeed,
if we assume that $\alpha =1,$ $\beta =\frac{\delta +\Delta }{2},$ $\delta
,\Delta \in \mathbb{K}$, with $\gamma =\func{Re}\left( \Delta \bar{\gamma}%
\right) >0,$ then the condition $\left\vert \beta \right\vert ^{2}\geq
\alpha \gamma $ is equivalent to $\left\vert \delta +\Delta \right\vert
^{2}\geq 4\func{Re}\left( \Delta \bar{\gamma}\right) $ which is actually $%
\left\vert \Delta -\delta \right\vert ^{2}\geq 0.$ With this assumption, (%
\ref{ch2x.2.47}) becomes%
\begin{equation*}
\left\Vert x-\frac{\delta +\Delta }{2}\cdot a\right\Vert \leq \frac{1}{2}%
\left\vert \Delta -\delta \right\vert \left\Vert a\right\Vert ,
\end{equation*}%
which implies the reverse of the Schwarz inequality%
\begin{align*}
\left\Vert x\right\Vert \left\Vert a\right\Vert & \leq \frac{\func{Re}\left[
\left( \bar{\Delta}+\bar{\delta}\right) \left\langle x,a\right\rangle \right]
}{2\sqrt{\func{Re}\left( \Delta \bar{\delta}\right) }} \\
& \leq \frac{\left\vert \Delta +\delta \right\vert }{2\sqrt{\func{Re}\left(
\Delta \bar{\delta}\right) }}\left\vert \left\langle x,a\right\rangle
\right\vert ,
\end{align*}%
which is (\ref{ch2x.1.11}).
\end{remark}

The following particular case of Theorem \ref{ch2x.t2.19} may be stated \cite%
{DRAG02}:

\begin{corollary}
\label{ch2x.c2.21}Let $\left( H;\left\langle \cdot ,\cdot \right\rangle
\right) $ be an inner product space over $\mathbb{K},$ $\varphi \in \lbrack
0,2\pi ),$ $\theta \in \left( 0,\frac{\pi }{2}\right) .$ If $x,a\in H$ are
such that $a\neq 0$ and%
\begin{equation}
\left\Vert x-\left( \cos \varphi +i\sin \varphi \right) a\right\Vert \leq
\cos \theta \left\Vert a\right\Vert ,  \label{ch2x.2.52}
\end{equation}%
then we have the reverses of the Schwarz inequality%
\begin{equation}
\left\Vert x\right\Vert \left\Vert a\right\Vert \leq \frac{\cos \varphi 
\func{Re}\left\langle x,a\right\rangle +\sin \varphi \func{Im}\left\langle
x,a\right\rangle }{\sin \theta }.  \label{ch2x.2.52a}
\end{equation}%
In particular, if%
\begin{equation*}
\left\Vert x-a\right\Vert \leq \cos \theta \left\Vert a\right\Vert ,
\end{equation*}%
then%
\begin{equation*}
\left\Vert x\right\Vert \left\Vert a\right\Vert \leq \frac{1}{\cos \theta }%
\func{Re}\left\langle x,a\right\rangle ;
\end{equation*}%
and if%
\begin{equation*}
\left\Vert x-ia\right\Vert \leq \cos \theta \left\Vert a\right\Vert ,
\end{equation*}%
then%
\begin{equation*}
\left\Vert x\right\Vert \left\Vert a\right\Vert \leq \frac{1}{\cos \theta }%
\func{Im}\left\langle x,a\right\rangle .
\end{equation*}
\end{corollary}

\subsection{Reverses of the Generalised Triangle Inequality}

In \cite{SSD4a}, the author obtained the following reverse result for the
generalised triangle inequality%
\begin{equation}
\sum_{i=1}^{n}\left\Vert x_{i}\right\Vert \geq \left\Vert
\sum_{i=1}^{n}x_{i}\right\Vert ,  \label{ch2x.3.1}
\end{equation}%
provided $x_{i}\in H,$ $i\in \left\{ 1,\dots ,n\right\} $ are vectors in a
real or complex inner product $\left( H;\left\langle \cdot ,\cdot
\right\rangle \right) :$

\begin{theorem}[Dragomir, 2004]
\label{ch2x.t3.1}Let $e,x_{i}\in H,$ $i\in \left\{ 1,\dots ,n\right\} $ with 
$\left\Vert e\right\Vert =1.$ If $k_{i}\geq 0,$ $i\in \left\{ 1,\dots
,n\right\} $ are such that%
\begin{equation}
\left( 0\leq \right) \left\Vert x_{i}\right\Vert -\func{Re}\left\langle
e,x_{i}\right\rangle \leq k_{i}\qquad \text{for each \qquad }i\in \left\{
1,\dots ,n\right\} ,  \label{ch2x.3.2}
\end{equation}%
then we have the inequality%
\begin{equation}
\left( 0\leq \right) \sum_{i=1}^{n}\left\Vert x_{i}\right\Vert -\left\Vert
\sum_{i=1}^{n}x_{i}\right\Vert \leq \sum_{i=1}^{n}k_{i}.  \label{ch2x.3.3}
\end{equation}%
The equality holds in (\ref{ch2x.3.3}) if and only if%
\begin{equation}
\sum_{i=1}^{n}\left\Vert x_{i}\right\Vert \geq \sum_{i=1}^{n}k_{i}
\label{ch2x.3.4}
\end{equation}%
and%
\begin{equation}
\sum_{i=1}^{n}x_{i}=\left( \sum_{i=1}^{n}\left\Vert x_{i}\right\Vert
-\sum_{i=1}^{n}k_{i}\right) e.  \label{ch2x.3.5}
\end{equation}
\end{theorem}

By utilising some of the results obtained in Section \ref{ch2x.s2}, we point
out several reverses of the generalised triangle inequality (\ref{ch2x.3.1})
that are corollaries of the above Theorem \ref{ch2x.t3.1} \cite{DRAG02}.

\begin{corollary}
\label{ch2x.c3.2}Let $e,$ $x_{i}\in H\backslash \left\{ 0\right\} ,$ $i\in
\left\{ 1,\dots ,n\right\} $ with $\left\Vert e\right\Vert =1.$ If%
\begin{equation}
\left\Vert \frac{x_{i}}{\left\Vert x_{i}\right\Vert }-e\right\Vert \leq
r_{i}\qquad \text{for each \qquad }i\in \left\{ 1,\dots ,n\right\} ,
\label{ch2x.3.6}
\end{equation}%
then%
\begin{align}
(0& \leq )\sum_{i=1}^{n}\left\Vert x_{i}\right\Vert -\left\Vert
\sum_{i=1}^{n}x_{i}\right\Vert  \label{ch2x.3.7} \\
& \leq \frac{1}{2}\sum_{i=1}^{n}r_{i}^{2}\left\Vert x_{i}\right\Vert  \notag
\\
& \leq \frac{1}{2}\times \left\{ 
\begin{array}{ll}
\left( \max\limits_{1\leq i\leq n}r_{i}\right)
^{2}\sum\limits_{i=1}^{n}\left\Vert x_{i}\right\Vert ; &  \\ 
&  \\ 
\left( \sum\limits_{i=1}^{n}r_{i}^{2p}\right) ^{\frac{1}{p}}\left(
\sum\limits_{i=1}^{n}\left\Vert x_{i}\right\Vert ^{q}\right) ^{\frac{1}{q}},
& p>1,\ \frac{1}{p}+\frac{1}{q}=1; \\ 
&  \\ 
\max\limits_{1\leq i\leq n}\left\Vert x_{i}\right\Vert
\sum\limits_{i=1}^{n}r_{i}^{2}. & 
\end{array}%
\right.  \notag
\end{align}
\end{corollary}

\begin{proof}
The first part follows from Proposition \ref{ch2x.p2.1} on choosing $%
x=x_{i}, $ $y=e$ and applying Theorem \ref{ch2x.t3.1}. The last part is
obvious by H\"{o}lder's inequality.
\end{proof}

\begin{remark}
\label{ch2x.r3.3}One would obtain the same reverse inequality (\ref{ch2x.3.7}%
) if one were to use Theorem \ref{ch2x.t2.4}. In this case, the assumption (%
\ref{ch2x.3.6}) should be replaced by%
\begin{equation}
\left\Vert \left\Vert x_{i}\right\Vert x_{i}-e\right\Vert \leq
r_{i}\left\Vert x_{i}\right\Vert \qquad \text{for each \qquad }i\in \left\{
1,\dots ,n\right\} .  \label{ch2x.3.8}
\end{equation}
\end{remark}

On utilising the inequalities (\ref{ch2x.2.5}) and (\ref{ch2x.2.13.a}) one
may state the following corollary of Theorem \ref{ch2x.t3.1} \cite{DRAG02}.

\begin{corollary}
\label{ch2x.c3.4}Let $e,$ $x_{i}\in H\backslash \left\{ 0\right\} ,$ $i\in
\left\{ 1,\dots ,n\right\} $ with $\left\Vert e\right\Vert =1.$ Then we have
the inequality%
\begin{equation}
(0\leq )\sum_{i=1}^{n}\left\Vert x_{i}\right\Vert -\left\Vert
\sum_{i=1}^{n}x_{i}\right\Vert \leq \min \left\{ A,B\right\} ,
\label{ch2x.3.9}
\end{equation}%
where%
\begin{equation*}
A:=2\sum_{i=1}^{n}\left\Vert x_{i}\right\Vert \left( \frac{\left\Vert
x_{i}-e\right\Vert }{\left\Vert x_{i}\right\Vert +1}\right) ^{2},
\end{equation*}%
and%
\begin{equation*}
B:=\frac{1}{2}\sum_{i=1}^{n}\frac{\left( \left\Vert x_{i}\right\Vert
+1\right) ^{2}\left\Vert x_{i}-e\right\Vert ^{2}}{\left\Vert
x_{i}\right\Vert }.
\end{equation*}
\end{corollary}

For vectors located outside the closed unit ball $\bar{B}\left( 0,1\right)
:=\left\{ z\in H|\left\Vert z\right\Vert \leq 1\right\} ,$ we may state the
following result \cite{DRAG02}.

\begin{corollary}
\label{ch2x.c3.5}Assume that $x_{i}\notin \bar{B}\left( 0,1\right) ,$ $i\in
\left\{ 1,\dots ,n\right\} $ and $e\in H,$ $\left\Vert e\right\Vert =1.$
Then we have the inequality:%
\begin{align}
(0& \leq )\sum_{i=1}^{n}\left\Vert x_{i}\right\Vert -\left\Vert
\sum_{i=1}^{n}x_{i}\right\Vert  \label{ch2x.3.10} \\
& \leq \left\{ 
\begin{array}{ll}
\dfrac{1}{2}p^{2}\sum\limits_{i=1}^{n}\left\Vert x_{i}\right\Vert
^{p-1}\left\Vert x_{i}-e\right\Vert ^{2}, & \text{if \ }p\geq 1 \\ 
&  \\ 
\dfrac{1}{2}\sum\limits_{i=1}^{n}\left\Vert x_{i}\right\Vert
^{1-p}\left\Vert x_{i}-e\right\Vert ^{2}, & \text{if \ }p<1.%
\end{array}%
\right.  \notag
\end{align}
\end{corollary}

The proof follows by Proposition \ref{ch2x.p2.5.a} and Theorem \ref%
{ch2x.t3.1}.

For complex spaces one may state the following result as well \cite{DRAG02}.

\begin{corollary}
\label{ch2x.c3.6}Let $\left( H;\left\langle \cdot ,\cdot \right\rangle
\right) $ be a complex inner product space and $\alpha _{i}\in \mathbb{C}$
with $\func{Re}\alpha _{i},$ $\func{Im}\alpha _{i}>0,$ $i\in \left\{ 1,\dots
,n\right\} . $ If $x_{i},e\in H,$ $i\in \left\{ 1,\dots ,n\right\} $ with $%
\left\Vert e\right\Vert =1$ and%
\begin{equation}
\left\Vert x_{i}-\frac{\func{Im}\alpha _{i}}{\func{Re}\alpha _{i}}\cdot
e\right\Vert \leq d_{i},\qquad i\in \left\{ 1,\dots ,n\right\} ,
\label{ch2x.3.11}
\end{equation}%
then%
\begin{equation}
(0\leq )\sum_{i=1}^{n}\left\Vert x_{i}\right\Vert -\left\Vert
\sum_{i=1}^{n}x_{i}\right\Vert \leq \frac{1}{2}\sum_{i=1}^{n}\frac{\func{Re}%
\alpha _{i}}{\func{Im}\alpha _{i}}\cdot d_{i}^{2}.  \label{ch2x.3.12}
\end{equation}
\end{corollary}

The proof follows by Theorems \ref{ch2x.t2.8} and \ref{ch2x.t3.1} and the
details are omitted.

Finally, by the use of Theorem \ref{ch2x.t2.13}, we can state \cite{DRAG02}:

\begin{corollary}
\label{ch2x.c3.7}If $x_{i},e\in H,$ $i\in \left\{ 1,\dots ,n\right\} $ with $%
\left\Vert e\right\Vert =1$ and $p\geq 1,$ then we have the inequalities:%
\begin{align}
(0& \leq )\sum_{i=1}^{n}\left\Vert x_{i}\right\Vert -\left\Vert
\sum_{i=1}^{n}x_{i}\right\Vert  \label{ch2x.3.13} \\
& \leq \frac{1}{2}\times \left\{ 
\begin{array}{l}
\sum\limits_{i=1}^{n}\left[ \left( \left\Vert x_{i}\right\Vert +1\right)
^{2p}-\left\Vert x_{i}+e\right\Vert ^{2p}\right] ^{\frac{1}{p}}, \\ 
\\ 
\sum\limits_{i=1}^{n}\left[ \left\Vert x_{i}-e\right\Vert ^{2p}-\left\vert
\left\Vert x_{i}\right\Vert -1\right\vert ^{2p}\right] ^{\frac{1}{p}}.%
\end{array}%
\right.  \notag
\end{align}
\end{corollary}

%


\chapter{Reverses for the Triangle Inequality}\label{ch3}

\section{Introduction}

The following reverse of the generalised triangle inequality%
\begin{equation*}
\cos \theta \sum_{k=1}^{n}\left\vert z_{k}\right\vert \leq \left\vert
\sum_{k=1}^{n}z_{k}\right\vert ,
\end{equation*}%
provided the complex numbers $z_{k},$ $k\in \left\{ 1,\dots ,n\right\} $
satisfy the assumption%
\begin{equation*}
a-\theta \leq \arg \left( z_{k}\right) \leq a+\theta ,\ \ \text{for any \ }%
k\in \left\{ 1,\dots ,n\right\} ,
\end{equation*}%
where $a\in \mathbb{R}$ and $\theta \in \left( 0,\frac{\pi }{2}\right) $ was
first discovered by M. Petrovich in 1917, \cite{RTxP} (see \cite[p. 492]%
{RTxMPF}) and subsequently was rediscovered by other authors, including J.
Karamata \cite[p. 300 -- 301]{RTxKA}, H.S. Wilf \cite{RTxW}, and in an
equivalent form by M. Marden \cite{RTxMA}.

In 1966, J.B. Diaz and F.T. Metcalf \cite{RTxDM} proved the following
reverse of the triangle inequality:

\begin{theorem}[Diaz-Metcalf, 1966]
\label{RT0.ta} Let $a$ be a unit vector in the inner product space $\left(
H;\left\langle \cdot ,\cdot \right\rangle \right) $ over the real or complex
number field $\mathbb{K}$. Suppose that the vectors $x_{i}\in H\backslash
\left\{ 0\right\} ,$ $i\in \left\{ 1,\dots ,n\right\} $ satisfy%
\begin{equation}
0\leq r\leq \frac{\func{Re}\left\langle x_{i},a\right\rangle }{\left\Vert
x_{i}\right\Vert },\ \ \ \ \ i\in \left\{ 1,\dots ,n\right\} .
\label{RT0.1.1}
\end{equation}%
Then%
\begin{equation}
r\sum_{i=1}^{n}\left\Vert x_{i}\right\Vert \leq \left\Vert
\sum_{i=1}^{n}x_{i}\right\Vert ,  \label{RT0.1.2}
\end{equation}%
where equality holds if and only if%
\begin{equation}
\sum_{i=1}^{n}x_{i}=r\left( \sum_{i=1}^{n}\left\Vert x_{i}\right\Vert
\right) a.  \label{RT0.1.3}
\end{equation}
\end{theorem}

A generalisation of this result for orthonormal families is incorporated in
the following result \cite{RTxDM}.

\begin{theorem}[Diaz-Metcalf, 1966]
\label{RT0.tb}Let $a_{1},\dots ,a_{n}$ be orthonormal vectors in $H.$
Suppose the vectors $x_{1},\dots ,x_{n}\in H\backslash \left\{ 0\right\} $
satisfy%
\begin{equation}
0\leq r_{k}\leq \frac{\func{Re}\left\langle x_{i},a_{k}\right\rangle }{%
\left\Vert x_{i}\right\Vert },\ \ \ \ \ i\in \left\{ 1,\dots ,n\right\} ,\
k\in \left\{ 1,\dots ,m\right\} .  \label{RT0.1.4}
\end{equation}%
Then%
\begin{equation}
\left( \sum_{k=1}^{m}r_{k}^{2}\right) ^{\frac{1}{2}}\sum_{i=1}^{n}\left\Vert
x_{i}\right\Vert \leq \left\Vert \sum_{i=1}^{n}x_{i}\right\Vert ,
\label{RT0.1.5}
\end{equation}%
where equality holds if and only if%
\begin{equation}
\sum_{i=1}^{n}x_{i}=\left( \sum_{i=1}^{n}\left\Vert x_{i}\right\Vert \right)
\sum_{k=1}^{m}r_{k}a_{k}.  \label{RT0.1.6}
\end{equation}
\end{theorem}

Similar results valid for semi-inner products may be found in \cite{RTxK}
and \cite{RTxM}.

For other classical inequalities related to the triangle inequality, see
Chapter XVII of the book \cite{RTxMPF} and the references therein.

The aim of the present chapter is to provide various recent reverses for the
generalised triangle inequality in both its simple form that are closely
related to the Diaz-Metcalf results mentioned above, or in the equivalent
quadratic form, i.e., upper bounds for

\begin{equation*}
\left( \sum_{i=1}^{n}\left\Vert x_{i}\right\Vert \right) ^{2}-\left\Vert
\sum_{i=1}^{n}x_{i}\right\Vert ^{2}
\end{equation*}%
and

\begin{equation*}
\frac{\left\Vert \sum_{i=1}^{n}x_{i}\right\Vert ^{2}}{\left(
\sum_{i=1}^{n}\left\Vert x_{i}\right\Vert \right) ^{2}}.
\end{equation*}

Applications for vector valued integral inequalities and for complex numbers
are given as well.

\section{Some Inequalities of Diaz-Metcalf Type}

\subsection{The Case of One Vector}

The following result with a natural geometrical meaning holds \cite{RTxSSD1}:

\begin{theorem}[Dragomir, 2004]
\label{RTt2.1.1}Let $a$ be a unit vector in the inner product space $\left(
H;\left\langle \cdot ,\cdot \right\rangle \right) $ and $\rho \in \left(
0,1\right) .$ If $x_{i}\in H,$ $i\in \left\{ 1,\dots ,n\right\} $ are such
that%
\begin{equation}
\left\Vert x_{i}-a\right\Vert \leq \rho \text{ \ for each \ }i\in \left\{
1,\dots ,n\right\} ,  \label{RT2.1.1}
\end{equation}%
then we have the inequality%
\begin{equation}
\sqrt{1-\rho ^{2}}\sum_{i=1}^{n}\left\Vert x_{i}\right\Vert \leq \left\Vert
\sum_{i=1}^{n}x_{i}\right\Vert ,  \label{RT2.2.1}
\end{equation}%
with equality if and only if%
\begin{equation}
\sum_{i=1}^{n}x_{i}=\sqrt{1-\rho ^{2}}\left( \sum_{i=1}^{n}\left\Vert
x_{i}\right\Vert \right) a.  \label{RT2.3.1}
\end{equation}
\end{theorem}

\begin{proof}
From (\ref{RT2.1.1}) we have%
\begin{equation*}
\left\Vert x_{i}\right\Vert ^{2}-2\func{Re}\left\langle x_{i},a\right\rangle
+1\leq \rho ^{2},
\end{equation*}%
giving%
\begin{equation}
\left\Vert x_{i}\right\Vert ^{2}+1-\rho ^{2}\leq 2\func{Re}\left\langle
x_{i},a\right\rangle ,  \label{RT2.4.1}
\end{equation}%
for each $i\in \left\{ 1,\dots ,n\right\} .$

Dividing by $\sqrt{1-\rho ^{2}}>0,$ we deduce%
\begin{equation}
\frac{\left\Vert x_{i}\right\Vert ^{2}}{\sqrt{1-\rho ^{2}}}+\sqrt{1-\rho ^{2}%
}\leq \frac{2\func{Re}\left\langle x_{i},a\right\rangle }{\sqrt{1-\rho ^{2}}}%
,  \label{RT2.5.1}
\end{equation}%
for each $i\in \left\{ 1,\dots ,n\right\} .$

On the other hand, by the elementary inequality%
\begin{equation}
\frac{p}{\alpha }+q\alpha \geq 2\sqrt{pq},\ \ \ p,q\geq 0,\ \alpha >0
\label{RT2.5.1.a}
\end{equation}%
we have%
\begin{equation}
2\left\Vert x_{i}\right\Vert \leq \frac{\left\Vert x_{i}\right\Vert ^{2}}{%
\sqrt{1-\rho ^{2}}}+\sqrt{1-\rho ^{2}}  \label{RT2.6.1}
\end{equation}%
and thus, by (\ref{RT2.5.1}) and (\ref{RT2.6.1}), we deduce%
\begin{equation*}
\frac{\func{Re}\left\langle x_{i},a\right\rangle }{\left\Vert
x_{i}\right\Vert }\geq \sqrt{1-\rho ^{2}},
\end{equation*}%
for each $i\in \left\{ 1,\dots ,n\right\} .$ Applying Theorem \ref{RT0.ta}
for $r=\sqrt{1-\rho ^{2}},$ we deduce the desired inequality (\ref{RT2.2.1}).
\end{proof}

The following results may be stated as well \cite{RTxSSD1}.

\begin{theorem}[Dragomir, 2004]
\label{RTt2.3.1}Let $a$ be a unit vector in the inner product space $\left(
H;\left\langle \cdot ,\cdot \right\rangle \right) $ and $M\geq m>0.$ If $%
x_{i}\in H,$ $i\in \left\{ 1,\dots ,n\right\} $ are such that either%
\begin{equation}
\func{Re}\left\langle Ma-x_{i},x_{i}-ma\right\rangle \geq 0  \label{RT2.14.1}
\end{equation}%
or, equivalently,%
\begin{equation}
\left\Vert x_{i}-\frac{M+m}{2}\cdot a\right\Vert \leq \frac{1}{2}\left(
M-m\right)  \label{RT2.15.1}
\end{equation}%
holds for each $i\in \left\{ 1,\dots ,n\right\} ,$ then we have the
inequality%
\begin{equation}
\frac{2\sqrt{mM}}{m+M}\sum_{i=1}^{n}\left\Vert x_{i}\right\Vert \leq
\left\Vert \sum_{i=1}^{n}x_{i}\right\Vert ,  \label{RT2.16.1}
\end{equation}%
or, equivalently,
\begin{equation}
\left( 0\leq \right) \sum_{i=1}^{n}\left\Vert x_{i}\right\Vert -\left\Vert
\sum_{i=1}^{n}x_{i}\right\Vert \leq \frac{\left( \sqrt{M}-\sqrt{m}\right)
^{2}}{2\sqrt{mM}}\left\Vert \sum_{i=1}^{n}x_{i}\right\Vert .
\label{RT2.17.1}
\end{equation}%
The equality holds in (\ref{RT2.16.1}) (or in (\ref{RT2.17.1})) if and only
if%
\begin{equation}
\sum_{i=1}^{n}x_{i}=\frac{2\sqrt{mM}}{m+M}\left( \sum_{i=1}^{n}\left\Vert
x_{i}\right\Vert \right) a.  \label{RT2.19.1}
\end{equation}
\end{theorem}

\begin{proof}
Firstly, we remark that if $x,z,Z\in H,$ then the following statements are
equivalent:

\begin{enumerate}
\item[(i)] $\func{Re}\left\langle Z-x,x-z\right\rangle \geq 0;$

\item[(ii)] $\left\Vert x-\frac{Z+z}{2}\right\Vert \leq \frac{1}{2}%
\left\Vert Z-z\right\Vert .$
\end{enumerate}

Using this fact, one may simply realize that (\ref{RT2.14.1}) and (\ref%
{RT2.15.1}) are equivalent.

Now, from (\ref{RT2.14.1}), we get%
\begin{equation*}
\left\Vert x_{i}\right\Vert ^{2}+mM\leq \left( M+m\right) \func{Re}%
\left\langle x_{i},a\right\rangle ,
\end{equation*}%
for any $i\in \left\{ 1,\dots ,n\right\} .$ Dividing this inequality by $%
\sqrt{mM}>0,$ we deduce the following inequality that will be used in the
sequel%
\begin{equation}
\frac{\left\Vert x_{i}\right\Vert ^{2}}{\sqrt{mM}}+\sqrt{mM}\leq \frac{M+m}{%
\sqrt{mM}}\func{Re}\left\langle x_{i},a\right\rangle ,  \label{RT2.20.1}
\end{equation}%
for each $i\in \left\{ 1,\dots ,n\right\} .$

Using the inequality (\ref{RT2.5.1.a}) from Theorem \ref{RTt2.1.1}, we also
have%
\begin{equation}
2\left\Vert x_{i}\right\Vert \leq \frac{\left\Vert x_{i}\right\Vert ^{2}}{%
\sqrt{mM}}+\sqrt{mM},  \label{RT2.21.1}
\end{equation}%
for each $i\in \left\{ 1,\dots ,n\right\} .$

Utilizing (\ref{RT2.20.1}) and (\ref{RT2.21.1}), we may conclude with the
following inequality%
\begin{equation*}
\left\Vert x_{i}\right\Vert \leq \frac{M+m}{\sqrt{mM}}\func{Re}\left\langle
x_{i},a\right\rangle ,
\end{equation*}%
which is equivalent to%
\begin{equation}
\frac{2\sqrt{mM}}{m+M}\leq \frac{\func{Re}\left\langle x_{i},a\right\rangle
}{\left\Vert x_{i}\right\Vert }  \label{RT2.22.1}
\end{equation}%
for any $i\in \left\{ 1,\dots ,n\right\} .$

Finally, on applying the Diaz-Metcalf result in Theorem \ref{RT0.ta} for $r=%
\frac{2\sqrt{mM}}{m+M}$, we deduce the desired conclusion.

The equivalence between (\ref{RT2.16.1}) and (\ref{RT2.17.1}) follows by
simple calculation and we omit the details.
\end{proof}

\subsection{The Case of $m$ Vectors}

In a similar manner to the one used in the proof of Theorem \ref{RTt2.1.1}
and by the use of the Diaz-Metcalf inequality incorporated in Theorem \ref%
{RT0.tb}, we can also prove the following result \cite{RTxSSD1} :

\begin{proposition}
\label{RTt2.2.1}Let $a_{1},\dots ,a_{n}$ be orthonormal vectors in $H.$
Suppose the vectors $x_{1},\dots ,x_{n}\in H\backslash \left\{ 0\right\} $
satisfy%
\begin{equation}
\left\Vert x_{i}-a_{k}\right\Vert \leq \rho _{k}\text{ \ for each}\ \ i\in
\left\{ 1,\dots ,n\right\} ,\ k\in \left\{ 1,\dots ,m\right\} ,
\label{RT2.11.1}
\end{equation}%
where $\rho _{k}\in \left( 0,1\right) ,k\in \left\{ 1,\dots ,m\right\} .$
Then we have the following reverse of the triangle inequality%
\begin{equation}
\left( m-\sum_{k=1}^{m}\rho _{k}^{2}\right) ^{\frac{1}{2}}\sum_{i=1}^{n}%
\left\Vert x_{i}\right\Vert \leq \left\Vert \sum_{i=1}^{n}x_{i}\right\Vert .
\label{RT2.12.1}
\end{equation}%
The equality holds in (\ref{RT2.12.1}) if and only if%
\begin{equation}
\sum_{i=1}^{n}x_{i}=\left( \sum_{i=1}^{n}\left\Vert x_{i}\right\Vert \right)
\sum_{k=1}^{m}\left( 1-\rho _{k}^{2}\right) ^{\frac{1}{2}}a_{k}.
\label{RT2.13.1}
\end{equation}
\end{proposition}

Finally, by the use of Theorem \ref{RT0.tb} and a similar technique to that
employed in the proof of Theorem \ref{RTt2.3.1}, we may state the following
result \cite{RTxSSD1}:

\begin{proposition}
\label{RTt2.4.1}Let $a_{1},\dots ,a_{n}$ be orthonormal vectors in $H.$
Suppose the vectors $x_{1},\dots ,x_{n}\in H\backslash \left\{ 0\right\} $
satisfy%
\begin{equation}
\func{Re}\left\langle M_{k}a_{k}-x_{i},x_{i}-\mu _{k}a_{k}\right\rangle \geq
0,  \label{RT2.23.1}
\end{equation}%
or, equivalently,%
\begin{equation}
\left\Vert x_{i}-\frac{M_{k}+\mu _{k}}{2}a_{k}\right\Vert \leq \frac{1}{2}%
\left( M_{k}-\mu _{k}\right) ,  \label{RT2.24.1}
\end{equation}%
for any $i\in \left\{ 1,\dots ,n\right\} $ and $k\in \left\{ 1,\dots
,m\right\} ,$ where $M_{k}\geq \mu _{k}>0$ for each $k\in \left\{ 1,\dots
,m\right\} .$

Then we have the inequality%
\begin{equation}
2\left( \sum_{k=1}^{m}\frac{\mu _{k}M_{k}}{\left( \mu _{k}+M_{k}\right) ^{2}}%
\right) ^{\frac{1}{2}}\sum_{i=1}^{n}\left\Vert x_{i}\right\Vert \leq
\left\Vert \sum_{i=1}^{n}x_{i}\right\Vert .  \label{RT2.25.1}
\end{equation}%
The equality holds in (\ref{RT2.25.1}) iff%
\begin{equation}
\sum_{i=1}^{n}x_{i}=2\left( \sum_{i=1}^{n}\left\Vert x_{i}\right\Vert
\right) \sum_{k=1}^{m}\frac{\sqrt{\mu _{k}M_{k}}}{\mu _{k}+M_{k}}a_{k}.
\label{RT2.26.1}
\end{equation}
\end{proposition}

\section{Additive Reverses for the Triangle Inequality}

\subsection{The Case of One Vector}

In this section we establish some additive reverses of the generalised
triangle inequality in real or complex inner product spaces.

The following result holds \cite{RTxSSD1}:

\begin{theorem}[Dragomir, 2004]
\label{RTt3.1.1}Let $\left( H;\left\langle \cdot ,\cdot \right\rangle
\right) $ be an inner product space over the real or complex number field $%
\mathbb{K}$ and $e,$ $x_{i}\in H,$ $i\in \left\{ 1,\dots ,n\right\} $ with $%
\left\Vert e\right\Vert =1.$ If $k_{i}\geq 0$, $i\in \left\{ 1,\dots
,n\right\} ,$ are such that%
\begin{equation}
\left\Vert x_{i}\right\Vert -\func{Re}\left\langle e,x_{i}\right\rangle \leq
k_{i}\text{\ \ for each}\ \ i\in \left\{ 1,\dots ,n\right\} ,
\label{RT3.1.1}
\end{equation}%
then we have the inequality%
\begin{equation}
\left( 0\leq \right) \sum_{i=1}^{n}\left\Vert x_{i}\right\Vert -\left\Vert
\sum_{i=1}^{n}x_{i}\right\Vert \leq \sum_{i=1}^{n}k_{i}.  \label{RT3.2.1}
\end{equation}%
The equality holds in (\ref{RT3.2.1}) if and only if%
\begin{equation}
\sum_{i=1}^{n}\left\Vert x_{i}\right\Vert \geq \sum_{i=1}^{n}k_{i}
\label{RT3.3.1}
\end{equation}%
and%
\begin{equation}
\sum_{i=1}^{n}x_{i}=\left( \sum_{i=1}^{n}\left\Vert x_{i}\right\Vert
-\sum_{i=1}^{n}k_{i}\right) e.  \label{RT3.4.1}
\end{equation}
\end{theorem}

\begin{proof}
If we sum in (\ref{RT3.1.1}) over $i$ from 1 to $n,$ then we get%
\begin{equation}
\sum_{i=1}^{n}\left\Vert x_{i}\right\Vert \leq \func{Re}\left\langle
e,\sum_{i=1}^{n}x_{i}\right\rangle +\sum_{i=1}^{n}k_{i}.  \label{RT3.5.1}
\end{equation}%
By Schwarz's inequality for $e$ and $\sum_{i=1}^{n}x_{i},$ we have%
\begin{align}
\func{Re}\left\langle e,\sum_{i=1}^{n}x_{i}\right\rangle & \leq \left\vert
\func{Re}\left\langle e,\sum_{i=1}^{n}x_{i}\right\rangle \right\vert
\label{RT3.6.1} \\
& \leq \left\vert \left\langle e,\sum_{i=1}^{n}x_{i}\right\rangle
\right\vert \leq \left\Vert e\right\Vert \left\Vert
\sum_{i=1}^{n}x_{i}\right\Vert =\left\Vert \sum_{i=1}^{n}x_{i}\right\Vert .
\notag
\end{align}%
Making use of (\ref{RT3.5.1}) and (\ref{RT3.6.1}), we deduce the desired
inequality (\ref{RT3.1.1}).

If (\ref{RT3.3.1}) and (\ref{RT3.4.1}) hold, then%
\begin{equation*}
\left\Vert \sum_{i=1}^{n}x_{i}\right\Vert =\left\vert
\sum_{i=1}^{n}\left\Vert x_{i}\right\Vert -\sum_{i=1}^{n}k_{i}\right\vert
\left\Vert e\right\Vert =\sum_{i=1}^{n}\left\Vert x_{i}\right\Vert
-\sum_{i=1}^{n}k_{i},
\end{equation*}%
and the equality in the second part of (\ref{RT3.2.1}) holds true.

Conversely, if the equality holds in (\ref{RT3.2.1}), then, obviously (\ref%
{RT3.3.1}) is valid and we need only to prove (\ref{RT3.4.1}).

Now, if the equality holds in (\ref{RT3.2.1}) then it must hold in (\ref%
{RT3.1.1}) for each $i\in \left\{ 1,\dots ,n\right\} $ and also must hold in
any of the inequalities in (\ref{RT3.6.1}).

It is well known that in Schwarz's inequality $\left\vert \left\langle
u,v\right\rangle \right\vert \leq \left\Vert u\right\Vert \left\Vert
v\right\Vert $ $\left( u,v\in H\right) $ the case of equality holds iff
there exists a $\lambda \in \mathbb{K}$ such that $u=\lambda v.$ We note
that in the weaker inequality $\func{Re}\left\langle u,v\right\rangle \leq
\left\Vert u\right\Vert \left\Vert v\right\Vert $ the case of equality holds
iff $\lambda \geq 0$ and $u=\lambda v.$

Consequently, the equality holds in all inequalities (\ref{RT3.6.1})
simultaneously iff there exists a $\mu \geq 0$ with%
\begin{equation}
\mu e=\sum_{i=1}^{n}x_{i}.  \label{RT3.7.1}
\end{equation}

If we sum the equalities in (\ref{RT3.1.1}) over $i$ from 1 to $n,$ then we
deduce%
\begin{equation}
\sum_{i=1}^{n}\left\Vert x_{i}\right\Vert -\func{Re}\left\langle
e,\sum_{i=1}^{n}x_{i}\right\rangle =\sum_{i=1}^{n}k_{i}.  \label{RT3.8.1}
\end{equation}%
Replacing $\sum_{i=1}^{n}\left\Vert x_{i}\right\Vert $ from (\ref{RT3.7.1})
into (\ref{RT3.8.1}), we deduce%
\begin{equation*}
\sum_{i=1}^{n}\left\Vert x_{i}\right\Vert -\mu \left\Vert e\right\Vert
^{2}=\sum_{i=1}^{n}k_{i},
\end{equation*}%
from where we get $\mu =\sum_{i=1}^{n}\left\Vert x_{i}\right\Vert
-\sum_{i=1}^{n}k_{i}.$ Using (\ref{RT3.7.1}), we deduce (\ref{RT3.4.1}) and
the theorem is proved.
\end{proof}

\subsection{The Case of $m$ Vectors}

If we turn our attention to the case of orthogonal families, then we may
state the following result as well \cite{RTxSSD1}.

\begin{theorem}[Dragomir, 2004]
\label{RTt3.2.1}Let $\left( H;\left\langle \cdot ,\cdot \right\rangle
\right) \ $be an inner product space over the real or complex number field $%
\mathbb{K}$, $\left\{ e_{k}\right\} _{k\in \left\{ 1,\dots ,m\right\} }$ a
family of orthonormal vectors in $H,$ $x_{i}\in H,$ $M_{i,k}\geq 0$ for $%
i\in \left\{ 1,\dots ,n\right\} $ and $k\in \left\{ 1,\dots ,m\right\} $
such that%
\begin{equation}
\left\Vert x_{i}\right\Vert -\func{Re}\left\langle e_{k},x_{i}\right\rangle
\leq M_{ik}\text{ }  \label{RT3.9.1}
\end{equation}%
\ for each$\ \ i\in \left\{ 1,\dots ,n\right\} ,\ k\in \left\{ 1,\dots
,m\right\} .$ Then we have the inequality%
\begin{equation}
\sum_{i=1}^{n}\left\Vert x_{i}\right\Vert \leq \frac{1}{\sqrt{m}}\left\Vert
\sum_{i=1}^{n}x_{i}\right\Vert +\frac{1}{m}\sum_{i=1}^{n}%
\sum_{k=1}^{m}M_{ik}.  \label{RT3.10.1}
\end{equation}%
The equality holds true in (\ref{RT3.10.1}) if and only if
\begin{equation}
\sum_{i=1}^{n}\left\Vert x_{i}\right\Vert \geq \frac{1}{m}%
\sum_{i=1}^{n}\sum_{k=1}^{m}M_{ik}  \label{RT3.11.1}
\end{equation}%
and%
\begin{equation}
\sum_{i=1}^{n}x_{i}=\left( \sum_{i=1}^{n}\left\Vert x_{i}\right\Vert -\frac{1%
}{m}\sum_{i=1}^{n}\sum_{k=1}^{m}M_{ik}\right) \sum_{k=1}^{m}e_{k}.
\label{RT3.12.1}
\end{equation}
\end{theorem}

\begin{proof}
If we sum over $i$ from 1 to $n$ in (\ref{RT3.9.1}), then we obtain%
\begin{equation*}
\sum_{i=1}^{n}\left\Vert x_{i}\right\Vert \leq \func{Re}\left\langle
e,\sum_{i=1}^{n}x_{i}\right\rangle +\sum_{i=1}^{n}M_{ik},
\end{equation*}%
for each $k\in \left\{ 1,\dots ,m\right\} .$ Summing these inequalities over
$k$ from 1 to $m,$ we deduce%
\begin{equation}
\sum_{i=1}^{n}\left\Vert x_{i}\right\Vert \leq \frac{1}{m}\func{Re}%
\left\langle \sum_{k=1}^{m}e_{k},\sum_{i=1}^{n}x_{i}\right\rangle +\frac{1}{m%
}\sum_{i=1}^{n}\sum_{k=1}^{m}M_{ik}.  \label{RT3.13.1}
\end{equation}%
By Schwarz's inequality for $\sum_{k=1}^{m}e_{k}$ and $\sum_{i=1}^{n}x_{i}$
we have%
\begin{align}
\func{Re}\left\langle \sum_{k=1}^{m}e_{k},\sum_{i=1}^{n}x_{i}\right\rangle &
\leq \left\vert \func{Re}\left\langle
\sum_{k=1}^{m}e_{k},\sum_{i=1}^{n}x_{i}\right\rangle \right\vert
\label{RT3.14.1} \\
& \leq \left\vert \left\langle
\sum_{k=1}^{m}e_{k},\sum_{i=1}^{n}x_{i}\right\rangle \right\vert  \notag \\
& \leq \left\Vert \sum_{k=1}^{m}e_{k}\right\Vert \left\Vert
\sum_{i=1}^{n}x_{i}\right\Vert  \notag \\
& =\sqrt{m}\left\Vert \sum_{i=1}^{n}x_{i}\right\Vert ,  \notag
\end{align}%
since, obviously,%
\begin{equation*}
\left\Vert \sum_{k=1}^{m}e_{k}\right\Vert =\sqrt{\left\Vert
\sum_{k=1}^{m}e_{k}\right\Vert ^{2}}=\sqrt{\sum_{k=1}^{m}\left\Vert
e_{k}\right\Vert ^{2}}=\sqrt{m}.
\end{equation*}%
Making use of (\ref{RT3.13.1}) and (\ref{RT3.14.1}), we deduce the desired
inequality (\ref{RT3.10.1}).

If (\ref{RT3.11.1}) and (\ref{RT3.12.1}) hold, then%
\begin{align*}
\frac{1}{\sqrt{m}}\left\Vert \sum_{i=1}^{n}x_{i}\right\Vert & =\left\vert
\sum_{i=1}^{n}\left\Vert x_{i}\right\Vert -\frac{1}{m}\sum_{i=1}^{n}%
\sum_{k=1}^{m}M_{ik}\right\vert \left\Vert \sum_{k=1}^{m}e_{k}\right\Vert \\
& =\frac{\sqrt{m}}{\sqrt{m}}\left( \sum_{i=1}^{n}\left\Vert x_{i}\right\Vert
-\frac{1}{m}\sum_{i=1}^{n}\sum_{k=1}^{m}M_{ik}\right) \\
& =\sum_{i=1}^{n}\left\Vert x_{i}\right\Vert -\frac{1}{m}\sum_{i=1}^{n}%
\sum_{k=1}^{m}M_{ik},
\end{align*}%
and the equality in (\ref{RT3.10.1}) holds true.

Conversely, if the equality holds in (\ref{RT3.10.1}), then, obviously (\ref%
{RT3.11.1}) is valid.

Now if the equality holds in (\ref{RT3.10.1}), then it must hold in (\ref%
{RT3.9.1}) for each $i\in \left\{ 1,\dots ,n\right\} $ and $k\in \left\{
1,\dots ,m\right\} $ and also must hold in any of the inequalities in (\ref%
{RT3.14.1}).

It is well known that in Schwarz's inequality $\func{Re}\left\langle
u,v\right\rangle \leq \left\Vert u\right\Vert \left\Vert v\right\Vert ,$ the
equality occurs iff $u=\lambda v$ with $\lambda \geq 0,$ consequently, the
equality holds in all inequalities (\ref{RT3.14.1}) simultaneously iff there
exists a $\mu \geq 0$ with%
\begin{equation}
\mu \sum_{k=1}^{m}e_{k}=\sum_{i=1}^{n}x_{i}.  \label{RT3.15.1}
\end{equation}%
If we sum the equality in (\ref{RT3.9.1}) over $i$ from 1 to $n$ and $k$
from 1 to $m,$ then we deduce%
\begin{equation}
m\sum_{i=1}^{n}\left\Vert x_{i}\right\Vert -\func{Re}\left\langle
\sum_{k=1}^{m}e_{k},\sum_{i=1}^{n}x_{i}\right\rangle
=\sum_{i=1}^{n}\sum_{k=1}^{m}M_{ik}.  \label{RT3.16.1}
\end{equation}%
Replacing $\sum_{i=1}^{n}x_{i}$ from (\ref{RT3.15.1}) into (\ref{RT3.16.1}),
we deduce%
\begin{equation*}
m\sum_{i=1}^{n}\left\Vert x_{i}\right\Vert -\mu \sum_{k=1}^{m}\left\Vert
e_{k}\right\Vert ^{2}=\sum_{i=1}^{n}\sum_{k=1}^{m}M_{ik}
\end{equation*}%
giving%
\begin{equation*}
\mu =\sum_{i=1}^{n}\left\Vert x_{i}\right\Vert -\frac{1}{m}%
\sum_{i=1}^{n}\sum_{k=1}^{m}M_{ik}.
\end{equation*}%
Using (\ref{RT3.15.1}), we deduce (\ref{RT3.12.1}) and the theorem is proved.
\end{proof}

\section{Further Additive Reverses}

\subsection{The Case of Small Balls}

In this section we point out different additive reverses of the generalised
triangle inequality under simpler conditions for the vectors involved.

The following result holds \cite{RTxSSD1}:

\begin{theorem}[Dragomir, 2004]
\label{RTt4.1.1}Let $\left( H;\left\langle \cdot ,\cdot \right\rangle
\right) $ be an inner product space over the real or complex number field $%
\mathbb{K}$ and $e,x_{i}\in H,$ $i\in \left\{ 1,\dots ,n\right\} $ with $%
\left\Vert e\right\Vert =1.$ If $\rho \in \left( 0,1\right) $ and $x_{i},$ $%
i\in \left\{ 1,\dots ,n\right\} $ are such that%
\begin{equation}
\left\Vert x_{i}-e\right\Vert \leq \rho \text{ \ for each}\ \ i\in \left\{
1,\dots ,n\right\} ,  \label{RT4.1.1}
\end{equation}%
then we have the inequality%
\begin{align}
& \left( 0\leq \right) \sum_{i=1}^{n}\left\Vert x_{i}\right\Vert -\left\Vert
\sum_{i=1}^{n}x_{i}\right\Vert  \label{RT4.2.1} \\
& \leq \frac{\rho ^{2}}{\sqrt{1-\rho ^{2}}\left( 1+\sqrt{1-\rho ^{2}}\right)
}\func{Re}\left\langle \sum_{i=1}^{n}x_{i},e\right\rangle  \notag \\
& \left( \leq \frac{\rho ^{2}}{\sqrt{1-\rho ^{2}}\left( 1+\sqrt{1-\rho ^{2}}%
\right) }\left\Vert \sum_{i=1}^{n}x_{i}\right\Vert \right) .  \notag
\end{align}%
The equality holds in (\ref{RT4.2.1}) if and only if%
\begin{equation}
\sum_{i=1}^{n}\left\Vert x_{i}\right\Vert \geq \frac{\rho ^{2}}{\sqrt{1-\rho
^{2}}\left( 1+\sqrt{1-\rho ^{2}}\right) }\func{Re}\left\langle
\sum_{i=1}^{n}x_{i},e\right\rangle  \label{RT4.3.1}
\end{equation}%
and
\begin{multline}
\sum_{i=1}^{n}x_{i}  \label{RT4.4.1} \\
=\left( \sum_{i=1}^{n}\left\Vert x_{i}\right\Vert -\frac{\rho ^{2}}{\sqrt{%
1-\rho ^{2}}\left( 1+\sqrt{1-\rho ^{2}}\right) }\func{Re}\left\langle
\sum_{i=1}^{n}x_{i},e\right\rangle \right) e.
\end{multline}
\end{theorem}

\begin{proof}
We know, from the proof of Theorem \ref{RTt3.1.1}, that, if (\ref{RT4.1.1})
is fulfilled, then we have the inequality%
\begin{equation*}
\left\Vert x_{i}\right\Vert \leq \frac{1}{\sqrt{1-\rho ^{2}}}\func{Re}%
\left\langle x_{i},e\right\rangle
\end{equation*}%
for each $i\in \left\{ 1,\dots ,n\right\} ,$ implying%
\begin{align}
\left\Vert x_{i}\right\Vert -\func{Re}\left\langle x_{i},e\right\rangle &
\leq \left( \frac{1}{\sqrt{1-\rho ^{2}}}-1\right) \func{Re}\left\langle
x_{i},e\right\rangle  \label{RT4.5.1} \\
& =\frac{\rho ^{2}}{\sqrt{1-\rho ^{2}}\left( 1+\sqrt{1-\rho ^{2}}\right) }%
\func{Re}\left\langle x_{i},e\right\rangle  \notag
\end{align}%
for each $i\in \left\{ 1,\dots ,n\right\} .$

Now, making use of Theorem \ref{RTt2.1.1}, for
\begin{equation*}
k_{i}:=\frac{\rho ^{2}}{\sqrt{1-\rho ^{2}}\left( 1+\sqrt{1-\rho ^{2}}\right)
}\func{Re}\left\langle x_{i},e\right\rangle ,\ \ \ i\in \left\{ 1,\dots
,n\right\} ,
\end{equation*}%
we easily deduce the conclusion of the theorem.

We omit the details.
\end{proof}

We may state the following result as well \cite{RTxSSD1}:

\begin{theorem}[Dragomir, 2004]
\label{RTt4.2.1}Let $\left( H;\left\langle \cdot ,\cdot \right\rangle
\right) $ be an inner product space and $e\in H,$ $M\geq m>0.$ If $x_{i}\in
H,$ $i\in \left\{ 1,\dots ,n\right\} $ are such that either%
\begin{equation}
\func{Re}\left\langle Me-x_{i},x_{i}-me\right\rangle \geq 0,  \label{RT4.6.1}
\end{equation}%
or, equivalently,%
\begin{equation}
\left\Vert x_{i}-\frac{M+m}{2}e\right\Vert \leq \frac{1}{2}\left( M-m\right)
\label{RT4.7.1}
\end{equation}%
holds for each $i\in \left\{ 1,\dots ,n\right\} ,$ then we have the
inequality%
\begin{align}
\left( 0\leq \right) \sum_{i=1}^{n}\left\Vert x_{i}\right\Vert -\left\Vert
\sum_{i=1}^{n}x_{i}\right\Vert & \leq \frac{\left( \sqrt{M}-\sqrt{m}\right)
^{2}}{2\sqrt{mM}}\func{Re}\left\langle \sum_{i=1}^{n}x_{i},e\right\rangle
\label{RT4.8.1} \\
& \left( \leq \frac{\left( \sqrt{M}-\sqrt{m}\right) ^{2}}{2\sqrt{mM}}%
\left\Vert \sum_{i=1}^{n}x_{i}\right\Vert \right) .  \notag
\end{align}%
The equality holds in (\ref{RT4.8.1}) if and only if%
\begin{equation}
\sum_{i=1}^{n}\left\Vert x_{i}\right\Vert \geq \frac{\left( \sqrt{M}-\sqrt{m}%
\right) ^{2}}{2\sqrt{mM}}\func{Re}\left\langle
\sum_{i=1}^{n}x_{i},e\right\rangle  \label{RT4.9.1}
\end{equation}%
and%
\begin{equation}
\sum_{i=1}^{n}x_{i}=\left( \sum_{i=1}^{n}\left\Vert x_{i}\right\Vert -\frac{%
\left( \sqrt{M}-\sqrt{m}\right) ^{2}}{2\sqrt{mM}}\func{Re}\left\langle
\sum_{i=1}^{n}x_{i},e\right\rangle \right) e.  \label{RT4.10.1}
\end{equation}
\end{theorem}

\begin{proof}
We know, from the proof of Theorem \ref{RTt2.3.1}, that if (\ref{RT4.6.1})
is fulfilled, then we have the inequality%
\begin{equation*}
\left\Vert x_{i}\right\Vert \leq \frac{M+m}{2\sqrt{mM}}\func{Re}\left\langle
x_{i},e\right\rangle
\end{equation*}%
for each $i\in \left\{ 1,\dots ,n\right\} .$ This is equivalent to%
\begin{equation*}
\left\Vert x_{i}\right\Vert -\func{Re}\left\langle x_{i},e\right\rangle \leq
\frac{\left( \sqrt{M}-\sqrt{m}\right) ^{2}}{2\sqrt{mM}}\func{Re}\left\langle
x_{i},e\right\rangle
\end{equation*}%
for each $i\in \left\{ 1,\dots ,n\right\} .$

Now, making use of Theorem \ref{RTt3.1.1}, we deduce the conclusion of the
theorem. We omit the details.
\end{proof}

\begin{remark}
If one uses Theorem \ref{RTt3.2.1} instead of Theorem \ref{RTt3.1.1} above,
then one can state the corresponding generalisation for families of
orthonormal vectors of the inequalities (\ref{RT4.2.1}) and (\ref{RT4.8.1})
respectively. We do not provide them here.
\end{remark}

\subsection{The Case of Arbitrary Balls}

Now, on utilising a slightly different approach, we may point out the
following result \cite{RTxSSD1}:

\begin{theorem}[Dragomir, 2004]
\label{RTt4.3.1}Let $\left( H;\left\langle \cdot ,\cdot \right\rangle
\right) $ be an inner product space over $\mathbb{K}$ and $e,$ $x_{i}\in H,$
$i\in \left\{ 1,\dots ,n\right\} $ with $\left\Vert e\right\Vert =1.$ If $%
r_{i}>0,$ $i\in \left\{ 1,\dots ,n\right\} $ are such that%
\begin{equation}
\left\Vert x_{i}-e\right\Vert \leq r_{i}\text{ \ for each}\ \ i\in \left\{
1,\dots ,n\right\} ,  \label{RT4.11.1}
\end{equation}%
then we have the inequality%
\begin{equation}
0\leq \sum_{i=1}^{n}\left\Vert x_{i}\right\Vert -\left\Vert
\sum_{i=1}^{n}x_{i}\right\Vert \leq \frac{1}{2}\sum_{i=1}^{n}r_{i}^{2}.
\label{RT4.12.1}
\end{equation}%
The equality holds in (\ref{RT4.12.1}) if and only if%
\begin{equation}
\sum_{i=1}^{n}\left\Vert x_{i}\right\Vert \geq \frac{1}{2}%
\sum_{i=1}^{n}r_{i}^{2}  \label{RT4.13.1}
\end{equation}%
and%
\begin{equation}
\sum_{i=1}^{n}x_{i}=\left( \sum_{i=1}^{n}\left\Vert x_{i}\right\Vert -\frac{1%
}{2}\sum_{i=1}^{n}r_{i}^{2}\right) e.  \label{RT4.14.1}
\end{equation}
\end{theorem}

\begin{proof}
The condition (\ref{RT4.11.1}) is clearly equivalent to%
\begin{equation}
\left\Vert x_{i}\right\Vert ^{2}+1\leq \func{Re}\left\langle
x_{i},e\right\rangle +r_{i}^{2}  \label{RT4.15.1}
\end{equation}%
for each $i\in \left\{ 1,\dots ,n\right\} .$

Using the elementary inequality%
\begin{equation}
2\left\Vert x_{i}\right\Vert \leq \left\Vert x_{i}\right\Vert ^{2}+1,
\label{RT4.16.1}
\end{equation}%
for each $i\in \left\{ 1,\dots ,n\right\} ,$ then, by (\ref{RT4.15.1}) and (%
\ref{RT4.16.1}), we deduce
\begin{equation*}
2\left\Vert x_{i}\right\Vert \leq 2\func{Re}\left\langle
x_{i},e\right\rangle +r_{i}^{2},
\end{equation*}%
giving%
\begin{equation}
\left\Vert x_{i}\right\Vert -\func{Re}\left\langle x_{i},e\right\rangle \leq
\frac{1}{2}r_{i}^{2}  \label{RT4.17.1}
\end{equation}%
for each $i\in \left\{ 1,\dots ,n\right\} .$

Now, utilising Theorem \ref{RTt3.1.1} for $k_{i}=\frac{1}{2}r_{i}^{2},$ $%
i\in \left\{ 1,\dots ,n\right\} ,$ we deduce the desired result. We omit the
details.
\end{proof}

Finally, we may state and prove the following result as well \cite{RTxSSD1}.

\begin{theorem}[Dragomir, 2004]
\label{RTt4.4.1}Let $\left( H;\left\langle \cdot ,\cdot \right\rangle
\right) $ be an inner product space over $\mathbb{K}$ and $e,$ $x_{i}\in H,$
$i\in \left\{ 1,\dots ,n\right\} $ with $\left\Vert e\right\Vert =1.$ If $%
M_{i}\geq m_{i}>0,$ $i\in \left\{ 1,\dots ,n\right\} ,$ are such that%
\begin{equation}
\left\Vert x_{i}-\frac{M_{i}+m_{i}}{2}e\right\Vert \leq \frac{1}{2}\left(
M_{i}-m_{i}\right) ,  \label{RT4.18.1}
\end{equation}%
or, equivalently,%
\begin{equation}
\func{Re}\left\langle M_{i}e-x,x-m_{i}e\right\rangle \geq 0  \label{RT4.19.1}
\end{equation}%
for each $i\in \left\{ 1,\dots ,n\right\} ,$ then we have the inequality%
\begin{equation}
\left( 0\leq \right) \sum_{i=1}^{n}\left\Vert x_{i}\right\Vert -\left\Vert
\sum_{i=1}^{n}x_{i}\right\Vert \leq \frac{1}{4}\sum_{i=1}^{n}\frac{\left(
M_{i}-m_{i}\right) ^{2}}{M_{i}+m_{i}}.  \label{RT4.20.1}
\end{equation}%
The equality holds in (\ref{RT4.20.1}) if and only if%
\begin{equation}
\sum_{i=1}^{n}\left\Vert x_{i}\right\Vert \geq \frac{1}{4}\sum_{i=1}^{n}%
\frac{\left( M_{i}-m_{i}\right) ^{2}}{M_{i}+m_{i}}  \label{RT4.21.1}
\end{equation}%
and
\begin{equation}
\sum_{i=1}^{n}x_{i}=\left( \sum_{i=1}^{n}\left\Vert x_{i}\right\Vert -\frac{1%
}{4}\sum_{i=1}^{n}\frac{\left( M_{i}-m_{i}\right) ^{2}}{M_{i}+m_{i}}\right)
e.  \label{RT4.22.1}
\end{equation}
\end{theorem}

\begin{proof}
The condition (\ref{RT4.18.1}) is equivalent to:%
\begin{equation*}
\left\Vert x_{i}\right\Vert ^{2}+\left( \frac{M_{i}+m_{i}}{2}\right)
^{2}\leq 2\func{Re}\left\langle x_{i},\frac{M_{i}+m_{i}}{2}e\right\rangle +%
\frac{1}{4}\left( M_{i}-m_{i}\right) ^{2}
\end{equation*}%
and since%
\begin{equation*}
2\left( \frac{M_{i}+m_{i}}{2}\right) \left\Vert x_{i}\right\Vert \leq
\left\Vert x_{i}\right\Vert ^{2}+\left( \frac{M_{i}+m_{i}}{2}\right) ^{2},
\end{equation*}%
then we get%
\begin{equation*}
2\left( \frac{M_{i}+m_{i}}{2}\right) \left\Vert x_{i}\right\Vert \leq 2\cdot
\frac{M_{i}+m_{i}}{2}\func{Re}\left\langle x_{i},e\right\rangle +\frac{1}{4}%
\left( M_{i}-m_{i}\right) ^{2},
\end{equation*}%
or, equivalently,%
\begin{equation*}
\left\Vert x_{i}\right\Vert -\func{Re}\left\langle x_{i},e\right\rangle \leq
\frac{1}{4}\cdot \frac{\left( M_{i}-m_{i}\right) ^{2}}{M_{i}+m_{i}}
\end{equation*}%
for each $i\in \left\{ 1,\dots ,n\right\} .$

Now, making use of Theorem \ref{RTt3.1.1} for $k_{i}:=\frac{1}{4}\cdot \frac{%
\left( M_{i}-m_{i}\right) ^{2}}{M_{i}+m_{i}},$ $i\in \left\{ 1,\dots
,n\right\} ,$ we deduce the desired result.
\end{proof}

\begin{remark}
If one uses Theorem \ref{RTt3.2.1} instead of Theorem \ref{RTt3.1.1} above,
then one can state the corresponding generalisation for families of
orthonormal vectors of the inequalities in (\ref{RT4.12.1}) and (\ref%
{RT4.20.1}) respectively. We omit the details.
\end{remark}

\section{Reverses of Schwarz Inequality}

In this section we outline a procedure showing how some of the above results
for triangle inequality may be employed to obtain reverses for the
celebrated Schwarz inequality.

For $a\in H,$ $\left\Vert a\right\Vert =1$ and $r\in \left( 0,1\right) $
define the closed ball
\begin{equation*}
\overline{D}\left( a,r\right) :=\left\{ x\in H,\left\Vert x-a\right\Vert
\leq r\right\} .
\end{equation*}%
The following reverse of the Schwarz inequality holds \cite{RTxSSD1}:

\begin{proposition}
\label{RTp5.1.1} If $x,y\in \overline{D}\left( a,r\right) $ with $a\in H,$ $%
\left\Vert a\right\Vert =1$ and $r\in \left( 0,1\right) ,$ then we have the
inequality%
\begin{equation}
\left( 0\leq \right) \frac{\left\Vert x\right\Vert \left\Vert y\right\Vert -%
\func{Re}\left\langle x,y\right\rangle }{\left( \left\Vert x\right\Vert
+\left\Vert y\right\Vert \right) ^{2}}\leq \frac{1}{2}r^{2}.  \label{RT5.1.1}
\end{equation}%
The constant $\frac{1}{2}$ in $\left( \ref{RT5.1.1}\right) $ is best
possible in the sense that it cannot be replaced by a smaller quantity.
\end{proposition}

\begin{proof}
Using Theorem \ref{RTt2.1.1} for $x_{1}=x,x_{2}=y,\rho =r,$ we have
\begin{equation}
\sqrt{1-r^{2}}\left( \left\Vert x\right\Vert +\left\Vert y\right\Vert
\right) \leq \left\Vert x+y\right\Vert .  \label{RT5.2.1}
\end{equation}%
Taking the square in (\ref{RT5.2.1}) we deduce
\begin{equation*}
\left( 1-r^{2}\right) \left( \left\Vert x\right\Vert ^{2}+2\left\Vert
x\right\Vert \left\Vert y\right\Vert +\left\Vert y\right\Vert ^{2}\right)
\leq \left\Vert x\right\Vert ^{2}+2\func{Re}\left\langle x,y\right\rangle
+\left\Vert y\right\Vert ^{2}
\end{equation*}%
which is clearly equivalent to (\ref{RT5.1.1}).

Now, assume that (\ref{RT5.1.1}) holds with a constant $C>0$ instead of $%
\frac{1}{2},i.e.,$%
\begin{equation}
\frac{\left\Vert x\right\Vert \left\Vert y\right\Vert -\func{Re}\left\langle
x,y\right\rangle }{\left( \left\Vert x\right\Vert +\left\Vert y\right\Vert
\right) ^{2}}\leq Cr^{2}  \label{RT5.3.1}
\end{equation}%
provided $x,y\in \overline{D}\left( a,r\right) $ with $a\in H,$ $\left\Vert
a\right\Vert =1$ and $r\in \left( 0,1\right) .$

Let $e\in H$ with $\left\Vert e\right\Vert =1$ and $e\perp a.$ Define $%
x=a+re,y=a-re.$ Then
\begin{equation*}
\left\Vert x\right\Vert =\sqrt{1+r^{2}}=\left\Vert y\right\Vert ,\text{ }%
\func{Re}\left\langle x,y\right\rangle =1-r^{2}
\end{equation*}%
and thus, from (\ref{RT5.3.1}), we have%
\begin{equation*}
\frac{1+r^{2}-\left( 1-r^{2}\right) }{\left( 2\sqrt{1+r^{2}}\right) ^{2}}%
\leq Cr^{2}
\end{equation*}%
giving
\begin{equation*}
\frac{1}{2}\leq \left( 1+r^{2}\right) C
\end{equation*}%
for any $r\in \left( 0,1\right) .$ If in this inequality we let $%
r\rightarrow 0+,$ then we get $C\geq \frac{1}{2}$ and the proposition is
proved.
\end{proof}

In a similar way, by the use of Theorem \ref{RTt2.3.1}, we may prove the
following reverse of the Schwarz inequality as well \cite{RTxSSD1}:

\begin{proposition}
\label{RTp5.2.1} If $a\in H,$ $\left\Vert a\right\Vert =1,$ $M\geq m>0$ and $%
x,y\in H$ are so that either%
\begin{equation*}
\func{Re}\left\langle Ma-x,x-ma\right\rangle ,\func{Re}\left\langle
Ma-y,y-ma\right\rangle \geq 0
\end{equation*}%
or, equivalently,%
\begin{equation*}
\left\Vert x-\frac{m+M}{2}a\right\Vert ,\left\Vert y-\frac{m+M}{2}%
a\right\Vert \leq \frac{1}{2}\left( M-m\right)
\end{equation*}%
hold, then
\begin{equation*}
\left( 0\leq \right) \frac{\left\Vert x\right\Vert \left\Vert y\right\Vert -%
\func{Re}\left\langle x,y\right\rangle }{\left( \left\Vert x\right\Vert
+\left\Vert y\right\Vert \right) ^{2}}\leq \frac{1}{2}\left( \frac{M-m}{M+m}%
\right) ^{2}.
\end{equation*}%
The constant $\frac{1}{2}$ cannot be replaced by a smaller quantity.
\end{proposition}

\begin{remark}
On utilising Theorem \ref{RTt2.2.1} and Theorem \ref{RTt2.4.1}, we may
deduce some similar reverses of Schwarz inequality provided $x,y\in \cap
_{k=1}^{m}\overline{D}\left( a_{k},\rho _{k}\right) ,$ assumed not to be
empty, where $a_{1},...,a_{n}$ are orthonormal vectors in $H$ and $\rho
_{k}\in \left( 0,1\right) $ for $k\in \left\{ 1,...,m\right\} .$ We omit the
details.
\end{remark}

\begin{remark}
For various different reverses of Schwarz inequality in inner product
spaces, see the recent survey \cite{RTxSSD}.
\end{remark}

\section{Quadratic Reverses of the Triangle Inequality}

\subsection{The General Case}

The following lemma holds \cite{RTxSSD2}:

\begin{lemma}[Dragomir, 2004]
\label{RTl2.1.2}Let $\left( H;\left\langle \cdot ,\cdot \right\rangle
\right) $ be an inner product space over the real or complex number field $%
\mathbb{K}$, $x_{i}\in H,$ $i\in \left\{ 1,\dots ,n\right\} $ and $k_{ij}>0$
for $1\leq i<j\leq n$ such that%
\begin{equation}
0\leq \left\Vert x_{i}\right\Vert \left\Vert x_{j}\right\Vert -\func{Re}%
\left\langle x_{i},x_{j}\right\rangle \leq k_{ij}  \label{RT2.1.2}
\end{equation}%
for $1\leq i<j\leq n.$ Then we have the following quadratic reverse of the
triangle inequality%
\begin{equation}
\left( \sum_{i=1}^{n}\left\Vert x_{i}\right\Vert \right) ^{2}\leq \left\Vert
\sum_{i=1}^{n}x_{i}\right\Vert ^{2}+2\sum_{1\leq i<j\leq n}k_{ij}.
\label{RT2.2.2}
\end{equation}%
The case of equality holds in (\ref{RT2.2.2}) if and only if it holds in (%
\ref{RT2.1.2}) for each $i,j$ with $1\leq i<j\leq n.$
\end{lemma}

\begin{proof}
We observe that the following identity holds:%
\begin{align}
& \left( \sum_{i=1}^{n}\left\Vert x_{i}\right\Vert \right) ^{2}-\left\Vert
\sum_{i=1}^{n}x_{i}\right\Vert ^{2}  \label{RT2.3.2} \\
& =\sum_{i,j=1}^{n}\left\Vert x_{i}\right\Vert \left\Vert x_{j}\right\Vert
-\left\langle \sum_{i=1}^{n}x_{i},\sum_{j=1}^{n}x_{j}\right\rangle  \notag \\
& =\sum_{i,j=1}^{n}\left\Vert x_{i}\right\Vert \left\Vert x_{j}\right\Vert
-\sum_{i,j=1}^{n}\func{Re}\left\langle x_{i},x_{j}\right\rangle  \notag \\
& =\sum_{i,j=1}^{n}\left[ \left\Vert x_{i}\right\Vert \left\Vert
x_{j}\right\Vert -\func{Re}\left\langle x_{i},x_{j}\right\rangle \right]
\notag \\
& =\sum_{1\leq i<j\leq n}\left[ \left\Vert x_{i}\right\Vert \left\Vert
x_{j}\right\Vert -\func{Re}\left\langle x_{i},x_{j}\right\rangle \right]
\notag \\
& +\sum_{1\leq j<i\leq n}\left[ \left\Vert x_{i}\right\Vert \left\Vert
x_{j}\right\Vert -\func{Re}\left\langle x_{i},x_{j}\right\rangle \right]
\notag \\
& =2\sum_{1\leq i<j\leq n}\left[ \left\Vert x_{i}\right\Vert \left\Vert
x_{j}\right\Vert -\func{Re}\left\langle x_{i},x_{j}\right\rangle \right] .
\notag
\end{align}%
Using the condition (\ref{RT2.1.2}), we deduce that%
\begin{equation*}
\sum_{1\leq i<j\leq n}\left[ \left\Vert x_{i}\right\Vert \left\Vert
x_{j}\right\Vert -\func{Re}\left\langle x_{i},x_{j}\right\rangle \right]
\leq \sum_{1\leq i<j\leq n}k_{ij},
\end{equation*}%
and by (\ref{RT2.3.2}), we get the desired inequality (\ref{RT2.2.2}).

The case of equality is obvious by the identity (\ref{RT2.3.2}) and we omit
the details.
\end{proof}

\begin{remark}
\label{RTr2.2.2}From (\ref{RT2.2.2}) one may deduce the coarser inequality
that might be useful in some applications:%
\begin{align*}
0& \leq \sum_{i=1}^{n}\left\Vert x_{i}\right\Vert -\left\Vert
\sum_{i=1}^{n}x_{i}\right\Vert  \\
& \leq \sqrt{2}\left( \sum_{1\leq i<j\leq n}k_{ij}\right) ^{\frac{1}{2}%
}\qquad \left( \leq \sqrt{2}\sum_{1\leq i<j\leq n}\sqrt{k_{ij}}\right) .
\end{align*}
\end{remark}

\begin{remark}
\label{RTr2.3.2}If the condition (\ref{RT2.1.2}) is replaced with the
following refinement of Schwarz's inequality:%
\begin{equation}
\left( 0\leq \right) \delta _{ij}\leq \left\Vert x_{i}\right\Vert \left\Vert
x_{j}\right\Vert -\func{Re}\left\langle x_{i},x_{j}\right\rangle \text{ for }%
1\leq i<j\leq n,  \label{RT2.4.2}
\end{equation}%
then the following refinement of the quadratic generalised triangle
inequality is valid:%
\begin{equation}
\left( \sum_{i=1}^{n}\left\Vert x_{i}\right\Vert \right) ^{2}\geq \left\Vert
\sum_{i=1}^{n}x_{i}\right\Vert ^{2}+2\sum_{1\leq i<j\leq n}\delta _{ij}\quad
\left( \geq \left\Vert \sum_{i=1}^{n}x_{i}\right\Vert ^{2}\right) .
\label{RT2.5.2}
\end{equation}%
The equality holds in the first part of (\ref{RT2.5.2}) iff the case of
equality holds in (\ref{RT2.4.2}) for each $1\leq i<j\leq n.$
\end{remark}

The following result holds \cite{RTxSSD2}.

\begin{proposition}
\label{RTp2.4.2}Let $\left( H;\left\langle \cdot ,\cdot \right\rangle
\right) $ be as above, $x_{i}\in H,$ $i\in \left\{ 1,\dots ,n\right\} $ and $%
r>0$ such that%
\begin{equation}
\left\Vert x_{i}-x_{j}\right\Vert \leq r  \label{RT2.6.2}
\end{equation}%
for $1\leq i<j\leq n.$ Then%
\begin{equation}
\left( \sum_{i=1}^{n}\left\Vert x_{i}\right\Vert \right) ^{2}\leq \left\Vert
\sum_{i=1}^{n}x_{i}\right\Vert ^{2}+\frac{n\left( n-1\right) }{2}r^{2}.
\label{RT2.7.2}
\end{equation}%
The case of equality holds in (\ref{RT2.7.2}) if and only if
\begin{equation}
\left\Vert x_{i}\right\Vert \left\Vert x_{j}\right\Vert -\func{Re}%
\left\langle x_{i},x_{j}\right\rangle =\frac{1}{2}r^{2}  \label{RT2.8.2}
\end{equation}%
for each $i,j$ with $1\leq i<j\leq n.$
\end{proposition}

\begin{proof}
The inequality (\ref{RT2.6.2}) is obviously equivalent to%
\begin{equation*}
\left\Vert x_{i}\right\Vert ^{2}+\left\Vert x_{j}\right\Vert ^{2}\leq 2\func{%
Re}\left\langle x_{i},x_{j}\right\rangle +r^{2}
\end{equation*}%
for $1\leq i<j\leq n.$ Since%
\begin{equation*}
2\left\Vert x_{i}\right\Vert \left\Vert x_{j}\right\Vert \leq \left\Vert
x_{i}\right\Vert ^{2}+\left\Vert x_{j}\right\Vert ^{2},\ \ 1\leq i<j\leq n;
\end{equation*}%
hence%
\begin{equation}
\left\Vert x_{i}\right\Vert \left\Vert x_{j}\right\Vert -\func{Re}%
\left\langle x_{i},x_{j}\right\rangle \leq \frac{1}{2}r^{2}  \label{RT2.9.2}
\end{equation}%
for any $i,j$ with $1\leq i<j\leq n.$

Applying Lemma \ref{RTl2.1.2} for $k_{ij}:=\frac{1}{2}r^{2}$ and taking into
account that
\begin{equation*}
\sum_{1\leq i<j\leq n}k_{ij}=\frac{n\left( n-1\right) }{4}r^{2},
\end{equation*}%
we deduce the desired inequality (\ref{RT2.7.2}). The case of equality is
also obvious by the above lemma and we omit the details.
\end{proof}

\subsection{Inequalities in Terms of the Forward Difference}

In the same spirit, and if some information about the forward difference $%
\Delta x_{k}:=x_{k+1}-x_{k}$ $\left( 1\leq k\leq n-1\right) $ are available,
then the following simple quadratic reverse of the generalised triangle
inequality may be stated \cite{RTxSSD2}.

\begin{corollary}
\label{RTc2.5.2}Let $\left( H;\left\langle \cdot ,\cdot \right\rangle
\right) $ be an inner product space and $x_{i}\in H,$ $i\in \left\{ 1,\dots
,n\right\} .$ Then we have the inequality%
\begin{equation}
\left( \sum_{i=1}^{n}\left\Vert x_{i}\right\Vert \right) ^{2}\leq \left\Vert
\sum_{i=1}^{n}x_{i}\right\Vert ^{2}+\frac{n\left( n-1\right) }{2}%
\sum_{k=1}^{n-1}\left\Vert \Delta x_{k}\right\Vert .  \label{RT2.10.2}
\end{equation}%
The constant $\frac{1}{2}$ is best possible in the sense that it cannot be
replaced in general by a smaller quantity.
\end{corollary}

\begin{proof}
Let $1\leq i<j\leq n.$ Then, obviously,
\begin{equation*}
\left\Vert x_{j}-x_{i}\right\Vert =\left\Vert \sum_{k=i}^{j-1}\Delta
x_{k}\right\Vert \leq \sum_{k=i}^{j-1}\left\Vert \Delta x_{k}\right\Vert
\leq \sum_{k=1}^{n-1}\left\Vert \Delta x_{k}\right\Vert .
\end{equation*}%
Applying Proposition \ref{RTp2.4.2} for $r:=\sum_{k=1}^{n-1}\left\Vert
\Delta x_{k}\right\Vert ,$ we deduce the desired result (\ref{RT2.10.2}).

To prove the sharpness of the constant $\frac{1}{2},$ assume that the
inequality (\ref{RT2.10.2}) holds with a constant $c>0,$ i.e.,
\begin{equation}
\left( \sum_{i=1}^{n}\left\Vert x_{i}\right\Vert \right) ^{2}\leq \left\Vert
\sum_{i=1}^{n}x_{i}\right\Vert ^{2}+cn\left( n-1\right)
\sum_{k=1}^{n-1}\left\Vert \Delta x_{k}\right\Vert  \label{RT2.11.2}
\end{equation}%
for $n\geq 2,$ $x_{i}\in H,$ $i\in \left\{ 1,\dots ,n\right\} .$

If we choose in (\ref{RT2.11.2}), $n=2,$ $x_{1}=-\frac{1}{2}e,$ $x_{2}=\frac{%
1}{2}e,$ $e\in H,$ $\left\Vert e\right\Vert =1,$ then we get $1\leq 2c,$
giving $c\geq \frac{1}{2}.$
\end{proof}

The following result providing a reverse of the quadratic generalised
triangle inequality in terms of the sup-norm of the forward differences also
holds \cite{RTxSSD2}.

\begin{proposition}
\label{RTp2.6.2}Let $\left( H;\left\langle \cdot ,\cdot \right\rangle
\right) $ be an inner product space and $x_{i}\in H,$ $i\in \left\{ 1,\dots
,n\right\} .$ Then we have the inequality%
\begin{equation}
\left( \sum_{i=1}^{n}\left\Vert x_{i}\right\Vert \right) ^{2}\leq \left\Vert
\sum_{i=1}^{n}x_{i}\right\Vert ^{2}+\frac{n^{2}\left( n^{2}-1\right) }{12}%
\max_{1\leq k\leq n-1}\left\Vert \Delta x_{k}\right\Vert ^{2}.
\label{RT2.12.2}
\end{equation}%
The constant $\frac{1}{12}$ is best possible in (\ref{RT2.12.2}).
\end{proposition}

\begin{proof}
As above, we have that%
\begin{equation*}
\left\Vert x_{j}-x_{i}\right\Vert \leq \sum_{k=i}^{j-1}\left\Vert \Delta
x_{k}\right\Vert \leq \left( j-i\right) \max_{1\leq k\leq n-1}\left\Vert
\Delta x_{k}\right\Vert ,
\end{equation*}%
for $1\leq i<j\leq n.$

Squaring the above inequality, we get%
\begin{equation*}
\left\Vert x_{j}\right\Vert ^{2}+\left\Vert x_{i}\right\Vert ^{2}\leq 2\func{%
Re}\left\langle x_{i},x_{j}\right\rangle +\left( j-i\right) ^{2}\max_{1\leq
k\leq n-1}\left\Vert \Delta x_{k}\right\Vert ^{2}
\end{equation*}%
for any $i,j$ with $1\leq i<j\leq n,$ and since%
\begin{equation*}
2\left\Vert x_{i}\right\Vert \left\Vert x_{j}\right\Vert \leq \left\Vert
x_{j}\right\Vert ^{2}+\left\Vert x_{i}\right\Vert ^{2},
\end{equation*}%
hence%
\begin{equation}
0\leq \left\Vert x_{i}\right\Vert \left\Vert x_{j}\right\Vert -\func{Re}%
\left\langle x_{i},x_{j}\right\rangle \leq \frac{1}{2}\left( j-i\right)
^{2}\max_{1\leq k\leq n-1}\left\Vert \Delta x_{k}\right\Vert ^{2}
\label{RT2.13.2}
\end{equation}%
for any $i,j$ with $1\leq i<j\leq n.$

Applying Lemma \ref{RTl2.1.2} for $k_{ij}:=\frac{1}{2}\left( j-i\right)
^{2}\max\limits_{1\leq k\leq n-1}\left\Vert \Delta x_{k}\right\Vert ^{2},$
we can state that%
\begin{equation*}
\left( \sum_{i=1}^{n}\left\Vert x_{i}\right\Vert \right) ^{2}\leq \left\Vert
\sum_{i=1}^{n}x_{i}\right\Vert ^{2}+\sum_{1\leq i<j\leq n}\left( j-i\right)
^{2}\max_{1\leq k\leq n-1}\left\Vert \Delta x_{k}\right\Vert ^{2}.
\end{equation*}%
However,%
\begin{align*}
\sum_{1\leq i<j\leq n}\left( j-i\right) ^{2}& =\frac{1}{2}%
\sum_{i,j=1}^{n}\left( j-i\right) ^{2}=n\sum_{k=1}^{n}k^{2}-\left(
\sum_{k=1}^{n}k\right) ^{2} \\
& =\frac{n^{2}\left( n^{2}-1\right) }{12}
\end{align*}%
giving the desired inequality.

To prove the sharpness of the constant, assume that (\ref{RT2.12.2}) holds
with a constant $D>0,$ i.e.,%
\begin{equation}
\left( \sum_{i=1}^{n}\left\Vert x_{i}\right\Vert \right) ^{2}\leq \left\Vert
\sum_{i=1}^{n}x_{i}\right\Vert ^{2}+Dn^{2}\left( n^{2}-1\right) \max_{1\leq
k\leq n-1}\left\Vert \Delta x_{k}\right\Vert ^{2}  \label{RT2.14.2}
\end{equation}%
for $n\geq 2,$ $x_{i}\in H,$ $i\in \left\{ 1,\dots ,n\right\} .$

If in (\ref{RT2.14.2}) we choose $n=2,$ $x_{1}=-\frac{1}{2}e,$ $x_{2}=\frac{1%
}{2}e,$ $e\in H,$ $\left\Vert e\right\Vert =1,$ then we get $1\leq 12D$
giving $D\geq \frac{1}{12}.$
\end{proof}

The following result may be stated as well \cite{RTxSSD2}.

\begin{proposition}
\label{RTp2.7.2}Let $\left( H;\left\langle \cdot ,\cdot \right\rangle
\right) $ be an inner product space and $x_{i}\in H,$ $i\in \left\{ 1,\dots
,n\right\} .$ Then we have the inequality:%
\begin{equation}
\left( \sum_{i=1}^{n}\left\Vert x_{i}\right\Vert \right) ^{2}\leq \left\Vert
\sum_{i=1}^{n}x_{i}\right\Vert ^{2}+\sum_{1\leq i<j\leq n}\left( j-i\right)
^{\frac{2}{q}}\left( \sum_{k=1}^{n-1}\left\Vert \Delta x_{k}\right\Vert
^{p}\right) ^{\frac{2}{p}},  \label{RT2.15.2}
\end{equation}%
where $p>1,$ $\frac{1}{p}+\frac{1}{q}=1.$

The constant $E=1$ in front of the double sum cannot generally be replaced
by a smaller constant.
\end{proposition}

\begin{proof}
Using H\"{o}lder's inequality, we have%
\begin{align*}
\left\Vert x_{j}-x_{i}\right\Vert & \leq \sum_{k=i}^{j-1}\left\Vert \Delta
x_{k}\right\Vert \leq \left( j-i\right) ^{\frac{1}{q}}\left(
\sum_{k=i}^{j-1}\left\Vert \Delta x_{k}\right\Vert ^{p}\right) ^{\frac{1}{p}}
\\
& \leq \left( j-i\right) ^{\frac{1}{q}}\left( \sum_{k=1}^{n-1}\left\Vert
\Delta x_{k}\right\Vert ^{p}\right) ^{\frac{1}{p}},
\end{align*}%
for $1\leq i<j\leq n.$

Squaring the previous inequality, we get%
\begin{equation*}
\left\Vert x_{j}\right\Vert ^{2}+\left\Vert x_{i}\right\Vert ^{2}\leq 2\func{%
Re}\left\langle x_{i},x_{j}\right\rangle +\left( j-i\right) ^{\frac{2}{q}%
}\left( \sum_{k=1}^{n-1}\left\Vert \Delta x_{k}\right\Vert ^{p}\right) ^{%
\frac{2}{p}},
\end{equation*}%
for $1\leq i<j\leq n.$

Utilising the same argument from the proof of Proposition \ref{RTp2.6.2}, we
deduce the desired inequality (\ref{RT2.15.2}).

Now assume that (\ref{RT2.15.2}) holds with a constant $E>0,$ i.e.,%
\begin{equation*}
\left( \sum_{i=1}^{n}\left\Vert x_{i}\right\Vert \right) ^{2}\leq \left\Vert
\sum_{i=1}^{n}x_{i}\right\Vert ^{2}+E\sum_{1\leq i<j\leq n}\left( j-i\right)
^{\frac{2}{q}}\left( \sum_{k=1}^{n-1}\left\Vert \Delta x_{k}\right\Vert
^{p}\right) ^{\frac{2}{p}},
\end{equation*}%
for $n\geq 2$ and $x_{i}\in H,$ $i\in \left\{ 1,\dots ,n\right\} ,$ $p>1,$ $%
\frac{1}{p}+\frac{1}{q}=1.$

For $n=2,$ $x_{1}=-\frac{1}{2}e,$ $x_{2}=\frac{1}{2}e,$ $\left\Vert
e\right\Vert =1,$ we get $1\leq E,$ showing the fact that the inequality (%
\ref{RT2.15.2}) is sharp.
\end{proof}

The particular case $p=q=2$ is of interest \cite{RTxSSD2}.

\begin{corollary}
\label{RTc2.8.2}Let $\left( H;\left\langle \cdot ,\cdot \right\rangle
\right) $ be an inner product space and $x_{i}\in H,$ $i\in \left\{ 1,\dots
,n\right\} .$ Then we have the inequality:%
\begin{equation}
\left( \sum_{i=1}^{n}\left\Vert x_{i}\right\Vert \right) ^{2}\leq \left\Vert
\sum_{i=1}^{n}x_{i}\right\Vert ^{2}+\frac{\left( n^{2}-1\right) n}{6}%
\sum_{k=1}^{n-1}\left\Vert \Delta x_{k}\right\Vert ^{2}.  \label{RT2.16.2}
\end{equation}%
The constant $\frac{1}{6}$ is best possible in (\ref{RT2.16.2}).
\end{corollary}

\begin{proof}
For $p=q=2,$ Proposition \ref{RTp2.7.2} provides the inequality%
\begin{equation*}
\left( \sum_{i=1}^{n}\left\Vert x_{i}\right\Vert \right) ^{2}\leq \left\Vert
\sum_{i=1}^{n}x_{i}\right\Vert ^{2}+\sum_{1\leq i<j\leq n}\left( j-i\right)
\sum_{k=1}^{n-1}\left\Vert \Delta x_{k}\right\Vert ^{2},
\end{equation*}%
and since%
\begin{align*}
& \sum\limits_{1\leq i<j\leq n}\left( j-i\right)  \\
& =1+\left( 1+2\right) +\left( 1+2+3\right) +\cdots +\left( 1+2+\cdots
+n-1\right)  \\
& =\sum\limits_{k=1}^{n-1}\left( 1+2+\cdots +k\right)
=\sum\limits_{k=1}^{n-1}\frac{k\left( k+1\right) }{2}=\frac{n\left(
n^{2}-1\right) }{6},
\end{align*}%
hence the inequality (\ref{RT2.15.2}) is proved. The best constant may be
shown in the same way as above but we omit the details.
\end{proof}

\subsection{A Different Quadratic Inequality}

Finally, we may state and prove the following different result \cite{RTxSSD2}%
.

\begin{theorem}[Dragomir, 2004]
\label{RTt2.9.2}Let $\left( H;\left\langle \cdot ,\cdot \right\rangle
\right) $ be an inner product space, $y_{i}\in H,$ $i\in \left\{ 1,\dots
,n\right\} $ and $M\geq m>0$ are such that either%
\begin{equation}
\func{Re}\left\langle My_{j}-y_{i},y_{i}-my_{j}\right\rangle \geq 0\text{ \
for }1\leq i<j\leq n,  \label{RT2.17.2}
\end{equation}%
or, equivalently,%
\begin{equation}
\left\Vert y_{i}-\frac{M+m}{2}y_{j}\right\Vert \leq \frac{1}{2}\left(
M-m\right) \left\Vert y_{j}\right\Vert \text{ \ for }1\leq i<j\leq n.
\label{RT2.18.2}
\end{equation}%
Then we have the inequality%
\begin{equation}
\left( \sum_{i=1}^{n}\left\Vert y_{i}\right\Vert \right) ^{2}\leq \left\Vert
\sum_{i=1}^{n}y_{i}\right\Vert ^{2}+\frac{1}{2}\cdot \frac{\left( M-m\right)
^{2}}{M+m}\sum\limits_{k=1}^{n-1}k\left\Vert y_{k+1}\right\Vert ^{2}.
\label{RT2.19.2}
\end{equation}%
The case of equality holds in (\ref{RT2.19.2}) if and only if%
\begin{equation}
\left\Vert y_{i}\right\Vert \left\Vert y_{j}\right\Vert -\func{Re}%
\left\langle y_{i},y_{j}\right\rangle =\frac{1}{4}\cdot \frac{\left(
M-m\right) ^{2}}{M+m}\left\Vert y_{j}\right\Vert ^{2}  \label{RT2.20.2}
\end{equation}%
for each $i,j$ with $1\leq i<j\leq n.$
\end{theorem}

\begin{proof}
Taking the square in (\ref{RT2.18.2}), we get%
\begin{multline*}
\quad \left\Vert y_{i}\right\Vert ^{2}+\frac{\left( M-m\right) ^{2}}{M+m}%
\left\Vert y_{j}\right\Vert ^{2} \\
\leq 2\func{Re}\left\langle y_{i},\frac{M+m}{2}y_{j}\right\rangle +\frac{1}{n%
}\left( M-m\right) ^{2}\left\Vert y_{j}\right\Vert ^{2}\quad
\end{multline*}%
for $1\leq i<j\leq n,$ and since, obviously,%
\begin{equation*}
2\left( \frac{M+m}{2}\right) \left\Vert y_{i}\right\Vert \left\Vert
y_{j}\right\Vert \leq \left\Vert y_{i}\right\Vert ^{2}+\frac{\left(
M-m\right) ^{2}}{M+m}\left\Vert y_{j}\right\Vert ^{2},
\end{equation*}%
hence%
\begin{multline*}
\quad  2\left( \frac{M+m}{2}\right) \left\Vert y_{i}\right\Vert
\left\Vert
y_{j}\right\Vert \\
\leq 2\func{Re}\left\langle y_{i},\frac{M+m}{2}y_{j}\right\rangle +\frac{1}{n%
}\left( M-m\right) ^{2}\left\Vert y_{j}\right\Vert ^{2},\qquad
\end{multline*}%
giving the much simpler inequality%
\begin{equation}
\left\Vert y_{i}\right\Vert \left\Vert y_{j}\right\Vert -\func{Re}%
\left\langle y_{i},y_{j}\right\rangle \leq \frac{1}{4}\cdot \frac{\left(
M-m\right) ^{2}}{M+m}\left\Vert y_{j}\right\Vert ^{2},  \label{RT2.21.2}
\end{equation}%
for $1\leq i<j\leq n.$

Applying Lemma \ref{RTl2.1.2} for $k_{ij}:=\frac{1}{4}\cdot \frac{\left(
M-m\right) ^{2}}{M+m}\left\Vert y_{j}\right\Vert ^{2},$ we deduce%
\begin{equation}
\left( \sum_{i=1}^{n}\left\Vert y_{i}\right\Vert \right) ^{2}\leq \left\Vert
\sum_{i=1}^{n}y_{i}\right\Vert ^{2}+\frac{1}{2}\frac{\left( M-m\right) ^{2}}{%
M+m}\sum_{1\leq i<j\leq n}\left\Vert y_{j}\right\Vert ^{2}  \label{RT2.22.2}
\end{equation}%
with equality if and only if (\ref{RT2.21.2}) holds for each $i,j$ with $%
1\leq i<j\leq n.$

Since%
\begin{align*}
\sum_{1\leq i<j\leq n}\left\Vert y_{j}\right\Vert ^{2}& =\sum_{1<j\leq
n}\left\Vert y_{j}\right\Vert ^{2}+\sum_{2<j\leq n}\left\Vert
y_{j}\right\Vert ^{2}+\cdots +\sum_{n-1<j\leq n}\left\Vert y_{j}\right\Vert
^{2} \\
& =\sum_{j=2}^{n}\left\Vert y_{j}\right\Vert ^{2}+\sum_{j=3}^{n}\left\Vert
y_{j}\right\Vert ^{2}+\cdots +\sum_{j=n-1}^{n}\left\Vert y_{j}\right\Vert
^{2}+\left\Vert y_{n}\right\Vert ^{2} \\
& =\sum_{j=2}^{n}\left( j-1\right) \left\Vert y_{j}\right\Vert
^{2}=\sum\limits_{k=1}^{n-1}k\left\Vert y_{k+1}\right\Vert ^{2},
\end{align*}%
hence the inequality (\ref{RT2.19.2}) is obtained.
\end{proof}

\section{Further Quadratic Refinements}

\subsection{The General Case}

The following lemma is of interest in itself as well \cite{RTxSSD2}.

\begin{lemma}[Dragomir, 2004]
\label{RTl3.1.2}Let $\left( H;\left\langle \cdot ,\cdot \right\rangle
\right) $ be an inner product space over the real or complex number field $%
\mathbb{K}$, $x_{i}\in H,$ $i\in \left\{ 1,\dots ,n\right\} $ and $k\geq 1$
with the property that:
\begin{equation}
\left\Vert x_{i}\right\Vert \left\Vert x_{j}\right\Vert \leq k\func{Re}%
\left\langle x_{i},x_{j}\right\rangle ,  \label{RT3.1.2}
\end{equation}%
for each $i,j$ with $1\leq i<j\leq n.$ Then%
\begin{equation}
\left( \sum_{i=1}^{n}\left\Vert x_{i}\right\Vert \right) ^{2}+\left(
k-1\right) \sum_{i=1}^{n}\left\Vert x_{i}\right\Vert ^{2}\leq k\left\Vert
\sum_{i=1}^{n}x_{i}\right\Vert ^{2}.  \label{RT3.2.2}
\end{equation}%
The equality holds in (\ref{RT3.2.2}) if and only if it holds in (\ref%
{RT3.1.2}) for each $i,j$ with $1\leq i<j\leq n.$
\end{lemma}

\begin{proof}
Firstly, let us observe that the following identity holds true:%
\begin{align}
& \left( \sum_{i=1}^{n}\left\Vert x_{i}\right\Vert \right) ^{2}-k\left\Vert
\sum_{i=1}^{n}x_{i}\right\Vert ^{2}  \label{RT3.3.2} \\
& =\sum_{i,j=1}^{n}\left\Vert x_{i}\right\Vert \left\Vert x_{j}\right\Vert
-k\left\langle \sum_{i=1}^{n}x_{i},\sum_{j=1}^{n}x_{j}\right\rangle  \notag
\\
& =\sum_{i,j=1}^{n}\left[ \left\Vert x_{i}\right\Vert \left\Vert
x_{j}\right\Vert -k\func{Re}\left\langle x_{i},x_{j}\right\rangle \right]
\notag \\
& =2\sum_{1\leq i<j\leq n}\left[ \left\Vert x_{i}\right\Vert \left\Vert
x_{j}\right\Vert -k\func{Re}\left\langle x_{i},x_{j}\right\rangle \right]
+\left( 1-k\right) \sum_{i=1}^{n}\left\Vert x_{i}\right\Vert ^{2},  \notag
\end{align}%
since, obviously, $\func{Re}\left\langle x_{i},x_{j}\right\rangle =\func{Re}%
\left\langle x_{j},x_{i}\right\rangle $ for any $i,j\in \left\{ 1,\dots
,n\right\} .$

Using the assumption (\ref{RT3.1.2}), we obtain%
\begin{equation*}
\sum_{1\leq i<j\leq n}\left[ \left\Vert x_{i}\right\Vert \left\Vert
x_{j}\right\Vert -k\func{Re}\left\langle x_{i},x_{j}\right\rangle \right]
\leq 0
\end{equation*}%
and thus, from (\ref{RT3.3.2}), we deduce the desired inequality (\ref%
{RT3.2.2}).

The case of equality is obvious by the identity (\ref{RT3.3.2}) and we omit
the details.
\end{proof}

\begin{remark}
\label{RTr3.4.2}The inequality (\ref{RT3.2.2}) provides the following
reverse of the quadratic generalised triangle inequality:%
\begin{equation}
0\leq \left( \sum_{i=1}^{n}\left\Vert x_{i}\right\Vert \right)
^{2}-\sum_{i=1}^{n}\left\Vert x_{i}\right\Vert ^{2}\leq k\left[ \left\Vert
\sum_{i=1}^{n}x_{i}\right\Vert ^{2}-\sum_{i=1}^{n}\left\Vert
x_{i}\right\Vert ^{2}\right] .  \label{RT3.4.2}
\end{equation}
\end{remark}

\begin{remark}
\label{RTr3.5.2}Since $k=1$ and $\sum_{i=1}^{n}\left\Vert x_{i}\right\Vert
^{2}\geq 0,$ hence by (\ref{RT3.2.2}) one may deduce the following reverse
of the triangle inequality%
\begin{equation}
\sum_{i=1}^{n}\left\Vert x_{i}\right\Vert \leq \sqrt{k}\left\Vert
\sum_{i=1}^{n}x_{i}\right\Vert ,  \label{RT3.5.2}
\end{equation}%
provided (\ref{RT3.1.2}) holds true for $1\leq i<j\leq n.$
\end{remark}

The following corollary providing a better bound for $\sum_{i=1}^{n}\left%
\Vert x_{i}\right\Vert ,$ holds \cite{RTxSSD2}.

\begin{corollary}
\label{RTc3.6.2}With the assumptions in Lemma \ref{RTl3.1.2}, one has the
inequality:%
\begin{equation}
\sum_{i=1}^{n}\left\Vert x_{i}\right\Vert \leq \sqrt{\frac{nk}{n+k-1}}%
\left\Vert \sum_{i=1}^{n}x_{i}\right\Vert .  \label{RT3.6.2}
\end{equation}
\end{corollary}

\begin{proof}
Using the Cauchy-Bunyakovsky-Schwarz inequality%
\begin{equation*}
n\sum_{i=1}^{n}\left\Vert x_{i}\right\Vert ^{2}\geq \left(
\sum_{i=1}^{n}\left\Vert x_{i}\right\Vert \right) ^{2}
\end{equation*}%
we get%
\begin{equation}
\left( k-1\right) \sum_{i=1}^{n}\left\Vert x_{i}\right\Vert ^{2}+\left(
\sum_{i=1}^{n}\left\Vert x_{i}\right\Vert \right) ^{2}\geq \left( \frac{k-1}{%
n}+1\right) \left( \sum_{i=1}^{n}\left\Vert x_{i}\right\Vert \right) ^{2}.
\label{RT3.7.2}
\end{equation}%
Consequently, by (\ref{RT3.7.2}) and (\ref{RT3.2.2}) we deduce%
\begin{equation*}
k\left\Vert \sum_{i=1}^{n}x_{i}\right\Vert ^{2}\geq \frac{n+k-1}{n}\left(
\sum_{i=1}^{n}\left\Vert x_{i}\right\Vert \right) ^{2}
\end{equation*}%
giving the desired inequality (\ref{RT3.6.2}).
\end{proof}

\subsection{Asymmetric Assumptions}

The following result may be stated as well \cite{RTxSSD2}.

\begin{theorem}[Dragomir, 2004]
\label{RTt3.7.2}Let $\left( H;\left\langle \cdot ,\cdot \right\rangle
\right) $ be an inner product space and $x_{i}\in H\backslash \left\{
0\right\} ,$ $i\in \left\{ 1,\dots ,n\right\} ,$ $\rho \in \left( 0,1\right)
,$ such that%
\begin{equation}
\left\Vert x_{i}-\frac{x_{j}}{\left\Vert x_{j}\right\Vert }\right\Vert \leq
\rho \ \ \text{ \ for }1\leq i<j\leq n.  \label{RT3.8.2}
\end{equation}%
Then we have the inequality%
\begin{multline}
\sqrt{1-\rho ^{2}}\left( \sum_{i=1}^{n}\left\Vert x_{i}\right\Vert \right)
^{2}+\left( 1-\sqrt{1-\rho ^{2}}\right) \sum_{i=1}^{n}\left\Vert
x_{i}\right\Vert ^{2}  \label{RT3.9.2} \\
\leq \left\Vert \sum_{i=1}^{n}x_{i}\right\Vert ^{2}.
\end{multline}%
The case of equality holds in (\ref{RT3.9.2}) iff%
\begin{equation}
\left\Vert x_{i}\right\Vert \left\Vert x_{j}\right\Vert =\frac{1}{\sqrt{%
1-\rho ^{2}}}\func{Re}\left\langle x_{i},x_{j}\right\rangle  \label{RT3.10.2}
\end{equation}%
for any $1\leq i<j\leq n.$
\end{theorem}

\begin{proof}
The condition (\ref{RT3.1.2}) is obviously equivalent to%
\begin{equation*}
\left\Vert x_{i}\right\Vert ^{2}+1-\rho ^{2}\leq 2\func{Re}\left\langle
x_{i},\frac{x_{j}}{\left\Vert x_{j}\right\Vert }\right\rangle
\end{equation*}%
for each $1\leq i<j\leq n.$

Dividing by $\sqrt{1-\rho ^{2}}>0,$ we deduce%
\begin{equation}
\frac{\left\Vert x_{i}\right\Vert ^{2}}{\sqrt{1-\rho ^{2}}}+\sqrt{1-\rho ^{2}%
}\leq \frac{2}{\sqrt{1-\rho ^{2}}}\func{Re}\left\langle x_{i},\frac{x_{j}}{%
\left\Vert x_{j}\right\Vert }\right\rangle ,  \label{RT3.11.2}
\end{equation}%
for $1\leq i<j\leq n.$

On the other hand, by the elementary inequality%
\begin{equation}
\frac{p}{\alpha }+q\alpha \geq 2\sqrt{pq},\ \ \ p,q\geq 0,\ \alpha >0
\label{RT3.12.2}
\end{equation}%
we have%
\begin{equation}
2\left\Vert x_{i}\right\Vert \leq \frac{\left\Vert x_{i}\right\Vert ^{2}}{%
\sqrt{1-\rho ^{2}}}+\sqrt{1-\rho ^{2}}.  \label{RT3.13.2}
\end{equation}%
Making use of (\ref{RT3.11.2}) and (\ref{RT3.13.2}), we deduce that
\begin{equation*}
\left\Vert x_{i}\right\Vert \left\Vert x_{j}\right\Vert \leq \frac{1}{\sqrt{%
1-\rho ^{2}}}\func{Re}\left\langle x_{i},x_{j}\right\rangle
\end{equation*}%
for $1\leq i<j\leq n.$

Now, applying Lemma \ref{RTl2.1.2} for $k=\frac{1}{\sqrt{1-\rho ^{2}}},$ we
deduce the desired result.
\end{proof}

\begin{remark}
\label{RTr3.8.2}If we assume that $\left\Vert x_{i}\right\Vert =1,$ $i\in
\left\{ 1,\dots ,n\right\} ,$ satisfying the simpler condition%
\begin{equation}
\left\Vert x_{j}-x_{i}\right\Vert \leq \rho \ \ \text{ \ for }1\leq i<j\leq
n,  \label{RT3.14.2}
\end{equation}%
then, from (\ref{RT3.9.2}), we deduce the following lower bound for $%
\left\Vert \sum_{i=1}^{n}x_{i}\right\Vert ,$ namely%
\begin{equation}
\left[ n+n\left( n-1\right) \sqrt{1-\rho ^{2}}\right] ^{\frac{1}{2}}\leq
\left\Vert \sum_{i=1}^{n}x_{i}\right\Vert .  \label{RT3.15.2}
\end{equation}%
The equality holds in (\ref{RT3.15.2}) iff $\sqrt{1-\rho ^{2}}=\func{Re}%
\left\langle x_{i},x_{j}\right\rangle $ for $1\leq i<j\leq n.$
\end{remark}

\begin{remark}
\label{RTr3.9.2}Under the hypothesis of Proposition \ref{RTp2.7.2}, we have
the coarser but simpler reverse of the triangle inequality%
\begin{equation}
\sqrt[4]{1-\rho ^{2}}\sum_{i=1}^{n}\left\Vert x_{i}\right\Vert \leq
\left\Vert \sum_{i=1}^{n}x_{i}\right\Vert .  \label{RT3.16.2}
\end{equation}%
Also, applying Corollary \ref{RTc3.6.2} for $k=\frac{1}{\sqrt{1-\rho ^{2}}},$
we can state that%
\begin{equation}
\sum_{i=1}^{n}\left\Vert x_{i}\right\Vert \leq \sqrt{\frac{n}{n\sqrt{1-\rho
^{2}}+1-\sqrt{1-\rho ^{2}}}}\left\Vert \sum_{i=1}^{n}x_{i}\right\Vert ,
\label{RT3.17.2}
\end{equation}%
provided $x_{i}\in H$ satisfy (\ref{RT3.8.2}) for $1\leq i<j\leq n.$
\end{remark}

In the same manner, we can state and prove the following reverse of the
quadratic generalised triangle inequality \cite{RTxSSD2}.

\begin{theorem}[Dragomir, 2004]
\label{RTt3.10.2}Let $\left( H;\left\langle \cdot ,\cdot \right\rangle
\right) $ be an inner product space over the real or complex number field $%
\mathbb{K}$, $x_{i}\in H,$ $i\in \left\{ 1,\dots ,n\right\} $ and $M\geq m>0$
such that either%
\begin{equation}
\func{Re}\left\langle Mx_{j}-x_{i},x_{i}-mx_{j}\right\rangle \geq 0\text{ \
for }1\leq i<j\leq n,  \label{RT3.18.2}
\end{equation}%
or, equivalently,
\begin{equation}
\left\Vert x_{i}-\frac{M+m}{2}x_{j}\right\Vert \leq \frac{1}{2}\left(
M-m\right) \left\Vert x_{j}\right\Vert \text{ \ for }1\leq i<j\leq n
\label{RT3.19.2}
\end{equation}%
hold. Then%
\begin{multline}
\frac{2\sqrt{mM}}{M+m}\left( \sum_{i=1}^{n}\left\Vert x_{i}\right\Vert
\right) ^{2}+\frac{\left( \sqrt{M}-\sqrt{m}\right) ^{2}}{M+m}%
\sum_{i=1}^{n}\left\Vert x_{i}\right\Vert ^{2}  \label{RT3.20.2} \\
\leq \left\Vert \sum_{i=1}^{n}x_{i}\right\Vert ^{2}.
\end{multline}%
The case of equality holds in (\ref{RT3.20.2}) if and only if%
\begin{equation}
\left\Vert x_{i}\right\Vert \left\Vert x_{j}\right\Vert =\frac{M+m}{2\sqrt{mM%
}}\func{Re}\left\langle x_{i},x_{j}\right\rangle \text{ \ for }1\leq i<j\leq
n.  \label{RT3.21.2}
\end{equation}
\end{theorem}

\begin{proof}
From (\ref{RT3.18.2}), observe that%
\begin{equation}
\left\Vert x_{i}\right\Vert ^{2}+Mm\left\Vert x_{j}\right\Vert ^{2}\leq
\left( M+m\right) \func{Re}\left\langle x_{i},x_{j}\right\rangle ,
\label{RT3.22.2}
\end{equation}%
for $1\leq i<j\leq n.$ Dividing (\ref{RT3.22.2}) by $\sqrt{mM}>0,$ we deduce%
\begin{equation*}
\frac{\left\Vert x_{i}\right\Vert ^{2}}{\sqrt{mM}}+\sqrt{mM}\left\Vert
x_{j}\right\Vert ^{2}\leq \frac{M+m}{\sqrt{mM}}\func{Re}\left\langle
x_{i},x_{j}\right\rangle ,
\end{equation*}%
and since, obviously%
\begin{equation*}
2\left\Vert x_{i}\right\Vert \left\Vert x_{j}\right\Vert \leq \frac{%
\left\Vert x_{i}\right\Vert ^{2}}{\sqrt{mM}}+\sqrt{mM}\left\Vert
x_{j}\right\Vert ^{2}
\end{equation*}%
hence%
\begin{equation*}
\left\Vert x_{i}\right\Vert \left\Vert x_{j}\right\Vert \leq \frac{M+m}{2%
\sqrt{mM}}\func{Re}\left\langle x_{i},x_{j}\right\rangle ,\text{ \ for }%
1\leq i<j\leq n.
\end{equation*}%
Applying Lemma \ref{RTl3.1.2} for $k=\frac{M+m}{2\sqrt{mM}}\geq 1,$ we
deduce the desired result.
\end{proof}

\begin{remark}
\label{RTr3.11.2}We also must note that a simpler but coarser inequality
that can be obtained from (\ref{RT3.20.2}) is
\begin{equation*}
\left( \frac{2\sqrt{mM}}{M+m}\right) ^{\frac{1}{2}}\sum_{i=1}^{n}\left\Vert
x_{i}\right\Vert \leq \left\Vert \sum_{i=1}^{n}x_{i}\right\Vert ,
\end{equation*}%
provided (\ref{RT3.18.2}) holds true.
\end{remark}

Finally, a different result related to the generalised triangle inequality
is incorporated in the following theorem \cite{RTxSSD2}.

\begin{theorem}[Dragomir, 2004]
\label{RTt3.12.2}Let $\left( H;\left\langle \cdot ,\cdot \right\rangle
\right) $ be an inner product space over $\mathbb{K}$, $\eta >0$ and $%
x_{i}\in H,$ $i\in \left\{ 1,\dots ,n\right\} $ with the property that%
\begin{equation}
\left\Vert x_{j}-x_{i}\right\Vert \leq \eta <\left\Vert x_{j}\right\Vert
\text{ \ for each \ }i,j\in \left\{ 1,\dots ,n\right\} .  \label{RT3.23.2}
\end{equation}%
Then we have the following reverse of the triangle inequality%
\begin{equation}
\frac{\sum_{i=1}^{n}\sqrt{\left\Vert x_{i}\right\Vert ^{2}-\eta ^{2}}}{%
\left\Vert \sum_{i=1}^{n}x_{i}\right\Vert }\leq \frac{\left\Vert
\sum_{i=1}^{n}x_{i}\right\Vert }{\sum_{i=1}^{n}\left\Vert x_{i}\right\Vert }.
\label{RT3.24.2}
\end{equation}%
The equality holds in (\ref{RT3.24.2}) iff%
\begin{equation}
\left\Vert x_{i}\right\Vert \sqrt{\left\Vert x_{j}\right\Vert ^{2}-\eta ^{2}}%
=\func{Re}\left\langle x_{i},x_{j}\right\rangle \text{ \ for each \ }i,j\in
\left\{ 1,\dots ,n\right\} .  \label{RT3.25.2}
\end{equation}
\end{theorem}

\begin{proof}
From (\ref{RT3.23.2}), we have%
\begin{equation*}
\left\Vert x_{i}\right\Vert ^{2}+\left\Vert x_{j}\right\Vert ^{2}-\eta
^{2}\leq 2\func{Re}\left\langle x_{i},x_{j}\right\rangle ,\ \ \ i,j\in
\left\{ 1,\dots ,n\right\} .
\end{equation*}%
On the other hand,
\begin{equation*}
2\left\Vert x_{i}\right\Vert \sqrt{\left\Vert x_{j}\right\Vert ^{2}-\eta ^{2}%
}\leq \left\Vert x_{i}\right\Vert ^{2}+\left\Vert x_{j}\right\Vert ^{2}-\eta
^{2},\ \ \ i,j\in \left\{ 1,\dots ,n\right\}
\end{equation*}%
and thus%
\begin{equation*}
\left\Vert x_{i}\right\Vert \sqrt{\left\Vert x_{j}\right\Vert ^{2}-\eta ^{2}}%
\leq \func{Re}\left\langle x_{i},x_{j}\right\rangle ,\ \ \ i,j\in \left\{
1,\dots ,n\right\} .
\end{equation*}%
Summing over $i,j\in \left\{ 1,\dots ,n\right\} ,$ we deduce the desired
inequality (\ref{RT3.24.2}).

The case of equality is also obvious from the above, and we omit the details.
\end{proof}

\section{Reverses for Complex Spaces}

\subsection{The Case of One Vector}

The following result holds \cite{RTxSSD3}.

\begin{theorem}[Dragomir, 2004]
\label{RTt2.1.3}Let $\left( H;\left\langle \cdot ,\cdot \right\rangle
\right) $ be a complex inner product space. Suppose that the vectors $%
x_{k}\in H,$ $k\in \left\{ 1,\dots ,n\right\} $ satisfy the condition%
\begin{equation}
0\leq r_{1}\left\Vert x_{k}\right\Vert \leq \func{Re}\left\langle
x_{k},e\right\rangle ,\quad 0\leq r_{2}\left\Vert x_{k}\right\Vert \leq
\func{Im}\left\langle x_{k},e\right\rangle  \label{RT2.1.3}
\end{equation}%
for each $k\in \left\{ 1,\dots ,n\right\} ,$ where $e\in H$ is such that $%
\left\Vert e\right\Vert =1$ and $r_{1},r_{2}\geq 0.$ Then we have the
inequality%
\begin{equation}
\sqrt{r_{1}^{2}+r_{2}^{2}}\sum_{k=1}^{n}\left\Vert x_{k}\right\Vert \leq
\left\Vert \sum_{k=1}^{n}x_{k}\right\Vert ,  \label{RT2.2.3}
\end{equation}%
where equality holds if and only if%
\begin{equation}
\sum_{k=1}^{n}x_{k}=\left( r_{1}+ir_{2}\right) \left(
\sum_{k=1}^{n}\left\Vert x_{k}\right\Vert \right) e.  \label{RT2.3.3}
\end{equation}
\end{theorem}

\begin{proof}
In view of the Schwarz inequality in the complex inner product space $\left(
H;\left\langle \cdot ,\cdot \right\rangle \right) ,$ we have%
\begin{align}
\left\Vert \sum_{k=1}^{n}x_{k}\right\Vert ^{2}& =\left\Vert
\sum_{k=1}^{n}x_{k}\right\Vert ^{2}\left\Vert e\right\Vert ^{2}\geq
\left\vert \left\langle \sum_{k=1}^{n}x_{k},e\right\rangle \right\vert ^{2}
\label{RT2.4.3} \\
& =\left\vert \left\langle \sum_{k=1}^{n}x_{k},e\right\rangle \right\vert
^{2}  \notag \\
& =\left\vert \sum_{k=1}^{n}\func{Re}\left\langle x_{k},e\right\rangle
+i\left( \sum_{k=1}^{n}\func{Im}\left\langle x_{k},e\right\rangle \right)
\right\vert ^{2}  \notag \\
& =\left( \sum_{k=1}^{n}\func{Re}\left\langle x_{k},e\right\rangle \right)
^{2}+\left( \sum_{k=1}^{n}\func{Im}\left\langle x_{k},e\right\rangle \right)
^{2}.  \notag
\end{align}%
Now, by hypothesis (\ref{RT2.1.3})%
\begin{equation}
\left( \sum_{k=1}^{n}\func{Re}\left\langle x_{k},e\right\rangle \right)
^{2}\geq r_{1}^{2}\left( \sum_{k=1}^{n}\left\Vert x_{k}\right\Vert \right)
^{2}  \label{RT2.5.3}
\end{equation}%
and%
\begin{equation}
\left( \sum_{k=1}^{n}\func{Im}\left\langle x_{k},e\right\rangle \right)
^{2}\geq r_{2}^{2}\left( \sum_{k=1}^{n}\left\Vert x_{k}\right\Vert \right)
^{2}.  \label{RT2.6.3}
\end{equation}%
If we add (\ref{RT2.5.3}) and (\ref{RT2.6.3}) and use (\ref{RT2.4.3}), then
we deduce the desired inequality (\ref{RT2.2.3}).

Now, if (\ref{RT2.3.3}) holds, then%
\begin{equation*}
\left\Vert \sum_{k=1}^{n}x_{k}\right\Vert =\left\vert
r_{1}+ir_{2}\right\vert \left( \sum_{k=1}^{n}\left\Vert x_{k}\right\Vert
\right) \left\Vert e\right\Vert =\sqrt{r_{1}^{2}+r_{2}^{2}}%
\sum_{k=1}^{n}\left\Vert x_{k}\right\Vert
\end{equation*}%
and the case of equality is valid in (\ref{RT2.2.3}).

Before we prove the reverse implication, let us observe that for $x\in H$
and $e\in H,$ $\left\Vert e\right\Vert =1,$ the following identity is true%
\begin{equation*}
\left\Vert x-\left\langle x,e\right\rangle e\right\Vert ^{2}=\left\Vert
x\right\Vert ^{2}-\left\vert \left\langle x,e\right\rangle \right\vert ^{2},
\end{equation*}%
therefore $\left\Vert x\right\Vert =\left\vert \left\langle x,e\right\rangle
\right\vert $ if and only if $x=\left\langle x,e\right\rangle e.$

If we assume that equality holds in (\ref{RT2.2.3}), then the case of
equality must hold in all the inequalities required in the argument used to
prove the inequality (\ref{RT2.2.3}), and we may state that%
\begin{equation}
\left\Vert \sum_{k=1}^{n}x_{k}\right\Vert =\left\vert \left\langle
\sum_{k=1}^{n}x_{k},e\right\rangle \right\vert ,  \label{RT2.7.3}
\end{equation}%
and
\begin{equation}
r_{1}\left\Vert x_{k}\right\Vert =\func{Re}\left\langle x_{k},e\right\rangle
,\quad r_{2}\left\Vert x_{k}\right\Vert =\func{Im}\left\langle
x_{k},e\right\rangle  \label{RT2.8.3}
\end{equation}%
for each $k\in \left\{ 1,\dots ,n\right\} .$

From (\ref{RT2.7.3}) we deduce%
\begin{equation}
\sum_{k=1}^{n}x_{k}=\left\langle \sum_{k=1}^{n}x_{k},e\right\rangle e
\label{RT2.9.3}
\end{equation}%
and from (\ref{RT2.8.3}), by multiplying the second equation with $i$ and
summing both equations over $k$ from $1$ to $n,$ we deduce%
\begin{equation}
\left( r_{1}+ir_{2}\right) \sum_{k=1}^{n}\left\Vert x_{k}\right\Vert
=\left\langle \sum_{k=1}^{n}x_{k},e\right\rangle .  \label{RT2.10.3}
\end{equation}%
Finally, by (\ref{RT2.10.3}) and (\ref{RT2.9.3}), we get the desired
equality (\ref{RT2.3.3}).
\end{proof}

The following corollary is of interest \cite{RTxSSD3}.

\begin{corollary}
\label{RTc2.2.3}Let $e$ a unit vector in the complex inner product space $%
\left( H;\left\langle \cdot ,\cdot \right\rangle \right) $ and $\rho
_{1},\rho _{2}\in \left( 0,1\right) .$ If $x_{k}\in H,$ $k\in \left\{
1,\dots ,n\right\} $ are such that%
\begin{equation}
\left\Vert x_{k}-e\right\Vert \leq \rho _{1},\ \ \ \left\Vert
x_{k}-ie\right\Vert \leq \rho _{2}\ \ \ \text{for each \ }k\in \left\{
1,\dots ,n\right\} ,  \label{RT2.11.3}
\end{equation}%
then we have the inequality%
\begin{equation}
\sqrt{2-\rho _{1}^{2}-\rho _{2}^{2}}\sum_{k=1}^{n}\left\Vert
x_{k}\right\Vert \leq \left\Vert \sum_{k=1}^{n}x_{k}\right\Vert ,
\label{RT2.12.3}
\end{equation}%
with equality if and only if%
\begin{equation}
\sum_{k=1}^{n}x_{k}=\left( \sqrt{1-\rho _{1}^{2}}+i\sqrt{1-\rho _{2}^{2}}%
\right) \left( \sum_{k=1}^{n}\left\Vert x_{k}\right\Vert \right) e.
\label{RT2.13.3}
\end{equation}
\end{corollary}

\begin{proof}
From the first inequality in (\ref{RT2.11.3}) we deduce%
\begin{equation}
0\leq \sqrt{1-\rho _{1}^{2}}\left\Vert x_{k}\right\Vert \leq \func{Re}%
\left\langle x_{k},e\right\rangle  \label{RT2.16.3}
\end{equation}%
for each $k\in \left\{ 1,\dots ,n\right\} .$

From the second inequality in (\ref{RT2.11.3}) we deduce%
\begin{equation*}
0\leq \sqrt{1-\rho _{2}^{2}}\left\Vert x_{k}\right\Vert \leq \func{Re}%
\left\langle x_{k},ie\right\rangle
\end{equation*}%
for each $k\in \left\{ 1,\dots ,n\right\} .$ Since%
\begin{equation*}
\func{Re}\left\langle x_{k},ie\right\rangle =\func{Im}\left\langle
x_{k},e\right\rangle ,
\end{equation*}%
hence%
\begin{equation}
0\leq \sqrt{1-\rho _{2}^{2}}\left\Vert x_{k}\right\Vert \leq \func{Im}%
\left\langle x_{k},e\right\rangle  \label{RT2.17.3}
\end{equation}%
for each $k\in \left\{ 1,\dots ,n\right\} .$

Now, observe from (\ref{RT2.16.3}) and (\ref{RT2.17.3}), that the condition (%
\ref{RT2.1.3}) of Theorem \ref{RTt2.1.3} is satisfied for $r_{1}=\sqrt{%
1-\rho _{1}^{2}},$ $r_{2}=\sqrt{1-\rho _{2}^{2}}\in \left( 0,1\right) ,$ and
thus the corollary is proved.
\end{proof}

The following corollary may be stated as well \cite{RTxSSD3}.

\begin{corollary}
\label{RTc2.3.3}Let $e$ be a unit vector in the complex inner product space $%
\left( H;\left\langle \cdot ,\cdot \right\rangle \right) $ and $M_{1}\geq
m_{1}>0,$ $M_{2}\geq m_{2}>0.$ If $x_{k}\in H,$ $k\in \left\{ 1,\dots
,n\right\} $ are such that either%
\begin{align}
\func{Re}\left\langle M_{1}e-x_{k},x_{k}-m_{1}e\right\rangle & \geq 0,\
\label{RT2.18.3} \\
\func{Re}\left\langle M_{2}ie-x_{k},x_{k}-m_{2}ie\right\rangle & \geq 0
\notag
\end{align}%
or, equivalently,%
\begin{align}
\left\Vert x_{k}-\frac{M_{1}+m_{1}}{2}e\right\Vert & \leq \frac{1}{2}\left(
M_{1}-m_{1}\right) ,  \label{RT2.19.3} \\
\left\Vert x_{k}-\frac{M_{2}+m_{2}}{2}ie\right\Vert & \leq \frac{1}{2}\left(
M_{2}-m_{2}\right) ,  \notag
\end{align}%
for each $k\in \left\{ 1,\dots ,n\right\} ,$ then we have the inequality%
\begin{equation}
2\left[ \frac{m_{1}M_{1}}{\left( M_{1}+m_{1}\right) ^{2}}+\frac{m_{2}M_{2}}{%
\left( M_{2}+m_{2}\right) ^{2}}\right] ^{\frac{1}{2}}\sum_{k=1}^{n}\left%
\Vert x_{k}\right\Vert \leq \left\Vert \sum_{k=1}^{n}x_{k}\right\Vert .
\label{RT2.20.3}
\end{equation}%
The equality holds in (\ref{RT2.20.3}) if and only if%
\begin{equation}
\sum_{k=1}^{n}x_{k}=2\left( \frac{\sqrt{m_{1}M_{1}}}{M_{1}+m_{1}}+i\frac{%
\sqrt{m_{2}M_{2}}}{M_{2}+m_{2}}\right) \left( \sum_{k=1}^{n}\left\Vert
x_{k}\right\Vert \right) e.  \label{RT2.21.3}
\end{equation}
\end{corollary}

\begin{proof}
From the first inequality in (\ref{RT2.18.3})%
\begin{equation}
0\leq \frac{2\sqrt{m_{1}M_{1}}}{M_{1}+m_{1}}\left\Vert x_{k}\right\Vert \leq
\func{Re}\left\langle x_{k},e\right\rangle  \label{RT2.24.3}
\end{equation}%
for each $k\in \left\{ 1,\dots ,n\right\} .$

Now, the proof follows the same path as the one of Corollary \ref{RTc2.2.3}
and we omit the details.
\end{proof}

\subsection{The Case of $m$ Orthonormal Vectors}

In \cite{RTxDM}, the authors have proved the following reverse of the
generalised triangle inequality in terms of orthonormal vectors \cite%
{RTxSSD3}.

\begin{theorem}[Diaz-Metcalf, 1966]
\label{RTt3.1.3}Let $e_{1},\dots ,e_{m}$ be orthonormal vectors in $\left(
H;\left\langle \cdot ,\cdot \right\rangle \right) $, i.e., we recall that $%
\left\langle e_{i},e_{j}\right\rangle =0$ if $i\neq j$ and $\left\Vert
e_{i}\right\Vert =1,$ $i,j\in \left\{ 1,\dots ,m\right\} .$ Suppose that the
vectors $x_{1},\dots ,x_{n}\in H$ satisfy%
\begin{equation}
0\leq r_{k}\left\Vert x_{j}\right\Vert \leq \func{Re}\left\langle
x_{j},e_{k}\right\rangle ,  \notag
\end{equation}%
$j\in \left\{ 1,\dots ,n\right\} ,\ k\in \left\{ 1,\dots ,m\right\} .$ Then%
\begin{equation}
\left( \sum_{k=1}^{m}r_{k}^{2}\right) ^{\frac{1}{2}}\sum_{j=1}^{n}\left\Vert
x_{j}\right\Vert \leq \left\Vert \sum_{j=1}^{n}x_{j}\right\Vert ,
\label{RT3.2.3}
\end{equation}%
where equality holds if and only if%
\begin{equation}
\sum_{j=1}^{n}x_{j}=\left( \sum_{j=1}^{n}\left\Vert x_{j}\right\Vert \right)
\sum_{k=1}^{m}r_{k}e_{k}.  \label{RT3.3.3}
\end{equation}
\end{theorem}

If the space $\left( H;\left\langle \cdot ,\cdot \right\rangle \right) $ is
complex and more information is available for the imaginary part, then the
following result may be stated as well \cite{RTxSSD3}.

\begin{theorem}[Dragomir, 2004]
\label{RTt3.2.3}Let $e_{1},\dots ,e_{m}\in H$ be an orthonormal family of
vectors in the complex inner product space $H.$ If the vectors $x_{1},\dots
,x_{n}\in H$ satisfy the conditions%
\begin{equation}
0\leq r_{k}\left\Vert x_{j}\right\Vert \leq \func{Re}\left\langle
x_{j},e_{k}\right\rangle ,\qquad 0\leq \rho _{k}\left\Vert x_{j}\right\Vert
\leq \func{Im}\left\langle x_{j},e_{k}\right\rangle  \label{RT3.4.3}
\end{equation}%
for each $j\in \left\{ 1,\dots ,n\right\} $ and $k\in \left\{ 1,\dots
,m\right\} ,$ then we have the following reverse of the generalised triangle
inequality;%
\begin{equation}
\left[ \sum_{k=1}^{m}\left( r_{k}^{2}+\rho _{k}^{2}\right) \right] ^{\frac{1%
}{2}}\sum_{j=1}^{n}\left\Vert x_{j}\right\Vert \leq \left\Vert
\sum_{j=1}^{n}x_{j}\right\Vert .  \label{RT3.5.3}
\end{equation}%
The equality holds in (\ref{RT3.5.3}) if and only if%
\begin{equation}
\sum_{j=1}^{n}x_{j}=\left( \sum_{j=1}^{n}\left\Vert x_{j}\right\Vert \right)
\sum_{k=1}^{m}\left( r_{k}+i\rho _{k}\right) e_{k}.  \label{RT3.6.3}
\end{equation}
\end{theorem}

\begin{proof}
Before we prove the theorem, let us recall that, if $x\in H$ and $%
e_{1},\dots ,e_{m}$ are orthogonal vectors, then the following identity
holds true:%
\begin{equation}
\left\Vert x-\sum_{k=1}^{m}\left\langle x,e_{k}\right\rangle
e_{k}\right\Vert ^{2}=\left\Vert x\right\Vert ^{2}-\sum_{k=1}^{n}\left\vert
\left\langle x,e_{k}\right\rangle \right\vert ^{2}.  \label{RT3.7.3}
\end{equation}%
As a consequence of this identity, we note the \textit{Bessel inequality}%
\begin{equation}
\sum_{k=1}^{m}\left\vert \left\langle x,e_{k}\right\rangle \right\vert
^{2}\leq \left\Vert x\right\Vert ^{2},x\in H.  \label{RT3.8.3}
\end{equation}%
The case of equality holds in (\ref{RT3.8.3}) if and only if (see (\ref%
{RT3.7.3}))%
\begin{equation}
x=\sum_{k=1}^{m}\left\langle x,e_{k}\right\rangle e_{k}.  \label{RT3.9.3}
\end{equation}%
Applying Bessel's inequality for $x=\sum_{j=1}^{n}x_{j},$ we have%
\begin{align}
\left\Vert \sum_{j=1}^{n}x_{j}\right\Vert ^{2}& \geq
\sum_{k=1}^{m}\left\vert \left\langle \sum_{j=1}^{n}x_{j},e_{k}\right\rangle
\right\vert ^{2}=\sum_{k=1}^{m}\left\vert \sum_{j=1}^{n}\left\langle
x_{j},e_{k}\right\rangle \right\vert ^{2}  \label{RT3.10.3} \\
& =\sum_{k=1}^{m}\left\vert \left( \sum_{j=1}^{n}\func{Re}\left\langle
x_{j},e_{k}\right\rangle \right) +i\left( \sum_{j=1}^{n}\func{Im}%
\left\langle x_{j},e_{k}\right\rangle \right) \right\vert ^{2}  \notag \\
& =\sum_{k=1}^{m}\left[ \left( \sum_{j=1}^{n}\func{Re}\left\langle
x_{j},e_{k}\right\rangle \right) ^{2}+\left( \sum_{j=1}^{n}\func{Im}%
\left\langle x_{j},e_{k}\right\rangle \right) ^{2}\right] .  \notag
\end{align}%
Now, by the hypothesis (\ref{RT3.4.3}) we have%
\begin{equation}
\left( \sum_{j=1}^{n}\func{Re}\left\langle x_{j},e_{k}\right\rangle \right)
^{2}\geq r_{k}^{2}\left( \sum_{j=1}^{n}\left\Vert x_{j}\right\Vert \right)
^{2}  \label{RT3.11}
\end{equation}%
and%
\begin{equation}
\left( \sum_{j=1}^{n}\func{Im}\left\langle x_{j},e_{k}\right\rangle \right)
^{2}\geq \rho _{k}^{2}\left( \sum_{j=1}^{n}\left\Vert x_{j}\right\Vert
\right) ^{2}.  \label{RT3.12.3}
\end{equation}%
Further, on making use of (\ref{RT3.10.3}) -- (\ref{RT3.12.3}), we deduce%
\begin{align*}
\left\Vert \sum_{j=1}^{n}x_{j}\right\Vert ^{2}& \geq \sum_{k=1}^{m}\left[
r_{k}^{2}\left( \sum_{j=1}^{n}\left\Vert x_{j}\right\Vert \right) ^{2}+\rho
_{k}^{2}\left( \sum_{j=1}^{n}\left\Vert x_{j}\right\Vert \right) ^{2}\right]
\\
& =\left( \sum_{j=1}^{n}\left\Vert x_{j}\right\Vert \right)
^{2}\sum_{k=1}^{m}\left( r_{k}^{2}+\rho _{k}^{2}\right) ,
\end{align*}%
which is clearly equivalent to (\ref{RT3.5.3}).

Now, if (\ref{RT3.6.3}) holds, then%
\begin{align*}
\left\Vert \sum_{j=1}^{n}x_{j}\right\Vert ^{2}& =\left(
\sum_{j=1}^{n}\left\Vert x_{j}\right\Vert \right) ^{2}\left\Vert
\sum_{k=1}^{m}\left( r_{k}+i\rho _{k}\right) e_{k}\right\Vert ^{2} \\
& =\left( \sum_{j=1}^{n}\left\Vert x_{j}\right\Vert \right)
^{2}\sum_{k=1}^{m}\left\vert r_{k}+i\rho _{k}\right\vert ^{2} \\
& =\left( \sum_{j=1}^{n}\left\Vert x_{j}\right\Vert \right)
^{2}\sum_{k=1}^{m}\left( r_{k}^{2}+\rho _{k}^{2}\right) ,
\end{align*}%
and the case of equality holds in (\ref{RT3.5.3}).

Conversely, if the equality holds in (\ref{RT3.5.3}), then it must hold in
all the inequalities used to prove (\ref{RT3.5.3}) and therefore we must have%
\begin{equation}
\left\Vert \sum_{j=1}^{n}x_{j}\right\Vert ^{2}=\sum_{k=1}^{m}\left\vert
\sum_{j=1}^{n}\left\langle x_{j},e_{k}\right\rangle \right\vert ^{2}
\label{RT3.13.3}
\end{equation}%
and%
\begin{equation}
r_{k}\left\Vert x_{j}\right\Vert =\func{Re}\left\langle
x_{j},e_{k}\right\rangle ,\qquad \rho _{k}\left\Vert x_{j}\right\Vert =\func{%
Im}\left\langle x_{j},e_{k}\right\rangle  \label{RT3.14.3}
\end{equation}%
for each $j\in \left\{ 1,\dots ,n\right\} $ and $k\in \left\{ 1,\dots
,m\right\} .$

Using the identity (\ref{RT3.7.3}), we deduce from (\ref{RT3.13.3}) that%
\begin{equation}
\sum_{j=1}^{n}x_{j}=\sum_{k=1}^{m}\left\langle
\sum_{j=1}^{n}x_{j},e_{k}\right\rangle e_{k}.  \label{RT3.15.3}
\end{equation}%
Multiplying the second equality in (\ref{RT3.14.3}) with the imaginary unit $%
i$ and summing the equality over $j$ from $1$ to $n,$ we deduce%
\begin{equation}
\left( r_{k}+i\rho _{k}\right) \sum_{j=1}^{n}\left\Vert x_{j}\right\Vert
=\left\langle \sum_{j=1}^{n}x_{j},e_{k}\right\rangle  \label{RT3.16.3}
\end{equation}%
for each $k\in \left\{ 1,\dots ,n\right\} .$

Finally, utilising (\ref{RT3.15.3}) and (\ref{RT3.16.3}), we deduce (\ref%
{RT3.6.3}) and the theorem is proved.
\end{proof}

The following corollaries are of interest \cite{RTxSSD3}.

\begin{corollary}
\label{RTc3.2.3}Let $e_{1},\dots ,e_{m}$ be orthonormal vectors in the
complex inner product space $\left( H;\left\langle \cdot ,\cdot
\right\rangle \right) $ and $\rho _{k},\eta _{k}\in \left( 0,1\right) ,$ $%
k\in \left\{ 1,\dots ,n\right\} .$ If $x_{1},\dots ,x_{n}\in H$ are such that%
\begin{equation*}
\left\Vert x_{j}-e_{k}\right\Vert \leq \rho _{k},\qquad \left\Vert
x_{j}-ie_{k}\right\Vert \leq \eta _{k}
\end{equation*}%
for each $j\in \left\{ 1,\dots ,n\right\} $ and $k\in \left\{ 1,\dots
,m\right\} ,$ then we have the inequality%
\begin{equation}
\left[ \sum_{k=1}^{m}\left( 2-\rho _{k}^{2}-\eta _{k}^{2}\right) \right] ^{%
\frac{1}{2}}\sum_{j=1}^{n}\left\Vert x_{j}\right\Vert \leq \left\Vert
\sum_{j=1}^{n}x_{j}\right\Vert .  \label{RT3.17.3}
\end{equation}%
The case of equality holds in (\ref{RT3.17.3}) if and only if%
\begin{equation}
\sum_{j=1}^{n}x_{j}=\left( \sum_{j=1}^{n}\left\Vert x_{j}\right\Vert \right)
\sum_{k=1}^{m}\left( \sqrt{1-\rho _{k}^{2}}+i\sqrt{1-\eta _{k}^{2}}\right)
e_{k}.  \label{RT3.18.3}
\end{equation}
\end{corollary}

The proof employs Theorem \ref{RTt3.2.3} and is similar to the one from
Corollary \ref{RTc2.2.3}. We omit the details.

\begin{corollary}
\label{RTc3.3.3}Let $e_{1},\dots ,e_{m}$ be as in Corollary \ref{RTc3.2.3}
and $M_{k}\geq m_{k}>0,$ $N_{k}\geq n_{k}>0,$ $k\in \left\{ 1,\dots
,m\right\} .$ If $x_{1},\dots ,x_{n}\in H$ are such that either%
\begin{equation*}
\func{Re}\left\langle M_{k}e_{k}-x_{j},x_{j}-m_{k}e_{k}\right\rangle \geq
0,\ \ \func{Re}\left\langle N_{k}ie_{k}-x_{j},x_{j}-n_{k}ie_{k}\right\rangle
\geq 0
\end{equation*}%
or, equivalently,%
\begin{align*}
\left\Vert x_{j}-\frac{M_{k}+m_{k}}{2}e_{k}\right\Vert & \leq \frac{1}{2}%
\left( M_{k}-m_{k}\right) ,\  \\
\left\Vert x_{j}-\frac{N_{k}+n_{k}}{2}ie_{k}\right\Vert & \leq \frac{1}{2}%
\left( N_{k}-n_{k}\right)
\end{align*}%
for each $j\in \left\{ 1,\dots ,n\right\} $ and $k\in \left\{ 1,\dots
,m\right\} ,$ then we have the inequality%
\begin{equation}
2\left\{ \sum_{k=1}^{m}\left[ \frac{m_{k}M_{k}}{\left( M_{k}+m_{k}\right)
^{2}}+\frac{n_{k}N_{k}}{\left( N_{k}+n_{k}\right) ^{2}}\right] \right\} ^{%
\frac{1}{2}}\sum_{j=1}^{n}\left\Vert x_{j}\right\Vert \leq \left\Vert
\sum_{j=1}^{n}x_{j}\right\Vert .  \label{RT3.19.3}
\end{equation}%
The case of equality holds in (\ref{RT3.19.3}) if and only if%
\begin{equation}
\sum_{j=1}^{n}x_{j}=2\left( \sum_{j=1}^{n}\left\Vert x_{j}\right\Vert
\right) \sum_{k=1}^{m}\left( \frac{\sqrt{m_{k}M_{k}}}{M_{k}+m_{k}}+i\frac{%
\sqrt{n_{k}N_{k}}}{N_{k}+n_{k}}\right) e_{k}.  \label{RT3.20.3}
\end{equation}
\end{corollary}

The proof employs Theorem \ref{RTt3.2.3} and is similar to the one in
Corollary \ref{RTc2.3.3}. We omit the details.

\section{Applications for Vector-Valued Integral Inequalities}

Let $\left( H;\left\langle \cdot ,\cdot \right\rangle \right) $ be a Hilbert
space over the real or complex number field, $\left[ a,b\right] $ a compact
interval in $\mathbb{R}$ and $\eta :\left[ a,b\right] \rightarrow \lbrack
0,\infty )$ a Lebesgue integrable function on $\left[ a,b\right] $ with the
property that $\int_{a}^{b}\eta \left( t\right) dt=1.$ If, by $L_{\eta
}\left( \left[ a,b\right] ;H\right) $ we denote the Hilbert space of all
Bochner measurable functions $f:\left[ a,b\right] \rightarrow H$ with the
property that $\int_{a}^{b}\eta \left( t\right) \left\Vert f\left( t\right)
\right\Vert ^{2}dt<\infty ,$ then the norm $\left\Vert \cdot \right\Vert
_{\eta }$ of this space is generated by the inner product $\left\langle
\cdot ,\cdot \right\rangle _{\eta }:H\times H\rightarrow \mathbb{K}$ defined
by
\begin{equation*}
\left\langle f,g\right\rangle _{\eta }:=\int_{a}^{b}\eta \left( t\right)
\left\langle f\left( t\right) ,g\left( t\right) \right\rangle dt.
\end{equation*}%
The following proposition providing a reverse of the integral generalised
triangle inequality may be stated \cite{RTxSSD1}.

\begin{proposition}
\label{RTp6.1.1} Let $\left( H;\left\langle \cdot ,\cdot \right\rangle
\right) $ be a Hilbert space and $\eta :\left[ a,b\right] \rightarrow
\lbrack 0,\infty )$ as above. If $g\in L_{\eta }\left( \left[ a,b\right]
;H\right) $ is so that $\int_{a}^{b}\eta \left( t\right) \left\Vert g\left(
t\right) \right\Vert ^{2}dt=1$ and $f_{i}\in L_{\eta }\left( \left[ a,b%
\right] ;H\right) ,i\in \left\{ 1,\dots ,n\right\} ,$ $\rho \in \left(
0,1\right) $ are so that%
\begin{equation}
\left\Vert f_{i}\left( t\right) -g\left( t\right) \right\Vert \leq \rho
\label{RT6.1.1}
\end{equation}%
for a.e. $t\in \left[ a,b\right] $ and each $i\in \left\{ 1,\dots ,n\right\}
,$ then we have the inequality%
\begin{multline}
\quad \sqrt{1-\rho ^{2}}\sum_{i=1}^{n}\left( \int_{a}^{b}\eta
\left( t\right) \left\Vert f_{i}\left( t\right) \right\Vert
^{2}dt\right) ^{\frac{1}{2}}
\label{RT6.2.1} \\
\leq \left( \int_{a}^{b}\eta \left( t\right) \left\Vert
\sum_{i=1}^{n}f_{i}\left( t\right) \right\Vert ^{2}dt\right)
^{\frac{1}{2}}.\quad
\end{multline}%
The case of equality holds in (\ref{RT6.2.1}) if and only if%
\begin{equation*}
\sum_{i=1}^{n}f_{i}\left( t\right) =\sqrt{1-\rho ^{2}}\sum_{i=1}^{n}\left(
\int_{a}^{b}\eta \left( t\right) \left\Vert f_{i}\left( t\right) \right\Vert
^{2}dt\right) ^{\frac{1}{2}}\cdot g\left( t\right)
\end{equation*}%
for a.e. $t\in \left[ a,b\right] .$
\end{proposition}

\begin{proof}
Observe, by (\ref{RT6.2.1}), that%
\begin{align*}
\left\Vert f_{i}-g\right\Vert _{\eta }& =\left( \int_{a}^{b}\eta \left(
t\right) \left\Vert f_{i}\left( t\right) -g\left( t\right) \right\Vert
^{2}dt\right) ^{\frac{1}{2}} \\
& \leq \left( \int_{a}^{b}\eta \left( t\right) \rho ^{2}dt\right) ^{\frac{1}{%
2}}=\rho
\end{align*}%
for each $i\in \left\{ 1,\dots ,n\right\} .$ Applying Theorem \ref{RTt2.1.1}
for the Hilbert space $L_{\eta }\left( \left[ a,b\right] ;H\right) ,$ we
deduce the desired result.
\end{proof}

The following result may be stated as well \cite{RTxSSD1}.

\begin{proposition}
\label{RTp6.2.1} Let $H,\eta ,g$ be as in Proposition \ref{RTp6.1.1}. If $%
f_{i}\in L_{\eta }\left( \left[ a,b\right] ;H\right) ,i\in \left\{ 1,\dots
,n\right\} $ and $M\geq m>0$ are so that either%
\begin{equation*}
\func{Re}\left\langle Mg\left( t\right) -f_{i}\left( t\right) ,f_{i}\left(
t\right) -mg\left( t\right) \right\rangle \geq 0
\end{equation*}%
or, equivalently,
\begin{equation*}
\left\Vert f_{i}\left( t\right) -\frac{m+M}{2}g\left( t\right) \right\Vert
\leq \frac{1}{2}\left( M-m\right)
\end{equation*}%
for a.e. $t\in \left[ a,b\right] $ and each $i\in \left\{ 1,\dots ,n\right\}
,$ then we have the inequality%
\begin{multline}
\quad  \frac{2\sqrt{mM}}{m+M}\sum_{i=1}^{n}\left( \int_{a}^{b}\eta
\left( t\right) \left\Vert f_{i}\left( t\right) \right\Vert
^{2}dt\right) ^{\frac{1}{2}}
\label{RT6.3.1} \\
\leq \left( \int_{a}^{b}\eta \left( t\right) \left\Vert
\sum_{i=1}^{n}f_{i}\left( t\right) \right\Vert ^{2}dt\right)
^{\frac{1}{2}}.\quad
\end{multline}%
The equality holds in (\ref{RT6.3.1}) if and only if%
\begin{equation*}
\sum_{i=1}^{n}f_{i}\left( t\right) =\frac{2\sqrt{mM}}{m+M}%
\sum_{i=1}^{n}\left( \int_{a}^{b}\eta \left( t\right) \left\Vert f_{i}\left(
t\right) \right\Vert ^{2}dt\right) ^{\frac{1}{2}}\cdot g\left( t\right) ,
\end{equation*}%
for a.e. $t\in \left[ a,b\right] .$
\end{proposition}

The following proposition providing a reverse of the integral generalised
triangle inequality may be stated \cite{RTxSSD2}.

\begin{proposition}
\label{RTp4.1.2} Let $\left( H;\left\langle \cdot ,\cdot \right\rangle
\right) $ be a Hilbert space and $\eta :\left[ a,b\right] \rightarrow
\lbrack 0,\infty )$ as above. If $g\in L_{\eta }\left( \left[ a,b\right]
;H\right) $ is so that $\int_{a}^{b}\eta \left( t\right) \left\Vert g\left(
t\right) \right\Vert ^{2}dt=1$ and $f_{i}\in L_{\eta }\left( \left[ a,b%
\right] ;H\right) ,i\in \left\{ 1,\dots ,n\right\} ,$ and $M\geq m>0$ are so
that either%
\begin{equation}
\func{Re}\left\langle Mf_{j}\left( t\right) -f_{i}\left( t\right)
,f_{i}\left( t\right) -mf_{j}\left( t\right) \right\rangle \geq 0
\label{RT4.1.2}
\end{equation}%
or, equivalently,
\begin{equation*}
\left\Vert f_{i}\left( t\right) -\frac{m+M}{2}f_{j}\left( t\right)
\right\Vert \leq \frac{1}{2}\left( M-m\right) \left\Vert f_{j}\left(
t\right) \right\Vert
\end{equation*}%
for a.e. $t\in \left[ a,b\right] $ and $1\leq i<j\leq n,$ then we have the
inequality%
\begin{multline}
\left[ \sum_{i=1}^{n}\left( \int_{a}^{b}\eta \left( t\right) \left\Vert
f_{i}\left( t\right) \right\Vert ^{2}dt\right) ^{\frac{1}{2}}\right] ^{2}
\label{RT4.2.2} \\
\leq \int_{a}^{b}\eta \left( t\right) \left\Vert \sum_{i=1}^{n}f_{i}\left(
t\right) \right\Vert ^{2}dt \\
+\frac{1}{2}\cdot \frac{\left( M-m\right) ^{2}}{m+M}\int_{a}^{b}\eta \left(
t\right) \left( \sum_{k=1}^{n-1}k\left\Vert f_{k+1}\left( t\right)
\right\Vert ^{2}\right) dt.
\end{multline}%
The case of equality holds in (\ref{RT4.2.2}) if and only if%
\begin{multline*}
\left( \int_{a}^{b}\eta \left( t\right) \left\Vert f_{i}\left( t\right)
\right\Vert ^{2}dt\right) ^{\frac{1}{2}}\left( \int_{a}^{b}\eta \left(
t\right) \left\Vert f_{j}\left( t\right) \right\Vert ^{2}dt\right) ^{\frac{1%
}{2}} \\
-\int_{a}^{b}\eta \left( t\right) \func{Re}\left\langle f_{i}\left( t\right)
,f_{j}\left( t\right) \right\rangle dt \\
=\frac{1}{4}\cdot \frac{\left( M-m\right) ^{2}}{m+M}\int_{a}^{b}\eta \left(
t\right) \left\Vert f_{j}\left( t\right) \right\Vert ^{2}dt
\end{multline*}%
for each $i,j$ with $1\leq i<j\leq n.$
\end{proposition}

\begin{proof}
We observe that
\begin{multline*}
\func{Re}\left\langle Mf_{j}-f_{i},f_{i}-mf_{j}\right\rangle _{\eta } \\
=\int_{a}^{b}\eta \left( t\right) \func{Re}\left\langle Mf_{j}\left(
t\right) -f_{i}\left( t\right) ,f_{i}\left( t\right) -mf_{j}\left( t\right)
\right\rangle dt\geq 0
\end{multline*}%
for any $i,j$ with $1\leq i<j\leq n.$

Applying Theorem \ref{RTt2.9.2} for the Hilbert space $L_{\eta }\left( \left[
a,b\right] ;H\right) $ and for $y_{i}=f_{i},i\in \left\{ 1,\dots ,n\right\}
, $ we deduce the desired result.
\end{proof}

Another integral inequality incorporated in the following proposition holds
\cite{RTxSSD2}:

\begin{proposition}
\label{RTp4.2.2} With the assumptions of Proposition \ref{RTp4.1.2}, we have
\begin{multline}
\frac{2\sqrt{mM}}{m+M}\left[ \sum_{i=1}^{n}\left( \int_{a}^{b}\eta \left(
t\right) \left\Vert f_{i}\left( t\right) \right\Vert ^{2}dt\right) ^{\frac{1%
}{2}}\right] ^{2}  \label{RT4.3.2} \\
+\frac{\left( \sqrt{M}-\sqrt{m}\right) ^{2}}{m+M}\sum_{i=1}^{n}\int_{a}^{b}%
\eta \left( t\right) \left\Vert f_{i}\left( t\right) \right\Vert ^{2}dt \\
\leq \int_{a}^{b}\eta \left( t\right) \left\Vert \sum_{i=1}^{n}f_{i}\left(
t\right) \right\Vert ^{2}dt.
\end{multline}%
The case of equality holds in (\ref{RT4.3.2}) if and only if
\begin{multline*}
\left( \int_{a}^{b}\eta \left( t\right) \left\Vert f_{i}\left( t\right)
\right\Vert ^{2}dt\right) ^{\frac{1}{2}}\left( \int_{a}^{b}\eta \left(
t\right) \left\Vert f_{j}\left( t\right) \right\Vert ^{2}dt\right) ^{\frac{1%
}{2}} \\
=\frac{M+m}{2\sqrt{mM}}\int_{a}^{b}\eta \left( t\right) \func{Re}%
\left\langle f_{i}\left( t\right) ,f_{j}\left( t\right) \right\rangle dt
\end{multline*}%
for any $i,j$ with $1\leq i<j\leq n.$
\end{proposition}

The proof is obvious by Theorem \ref{RTt3.10.2} and we omit the details.

\section{Applications for Complex Numbers}

The following reverse of the generalised triangle inequality with a clear
geometric meaning may be stated \cite{RTxSSD3}.

\begin{proposition}
\label{RTp4.1.3}Let $z_{1},\dots ,z_{n}$ be complex numbers with the
property that%
\begin{equation}
0\leq \varphi _{1}\leq \arg \left( z_{k}\right) \leq \varphi _{2}<\frac{\pi
}{2}  \label{RT4.1.3}
\end{equation}%
for each $k\in \left\{ 1,\dots ,n\right\} .$ Then we have the inequality%
\begin{equation}
\sqrt{\sin ^{2}\varphi _{1}+\cos ^{2}\varphi _{2}}\sum_{k=1}^{n}\left\vert
z_{k}\right\vert \leq \left\vert \sum_{k=1}^{n}z_{k}\right\vert .
\label{RT4.3.3}
\end{equation}%
The equality holds in (\ref{RT4.3.3}) if and only if%
\begin{equation}
\sum_{k=1}^{n}z_{k}=\left( \cos \varphi _{2}+i\sin \varphi _{1}\right)
\sum_{k=1}^{n}\left\vert z_{k}\right\vert .  \label{RT4.4.3}
\end{equation}
\end{proposition}

\begin{proof}
Let $z_{k}=a_{k}+ib_{k}.$ We may assume that $b_{k}\geq 0,$ $a_{k}>0,$ $k\in
\left\{ 1,\dots ,n\right\} ,$ since, by (\ref{RT4.1.3}), $\frac{b_{k}}{a_{k}}%
=\tan \left[ \arg \left( z_{k}\right) \right] \in \left[ 0,\frac{\pi }{2}%
\right) ,$ $k\in \left\{ 1,\dots ,n\right\} .$ By (\ref{RT4.1.3}), we
obviously have%
\begin{equation*}
0\leq \tan ^{2}\varphi _{1}\leq \frac{b_{k}^{2}}{a_{k}^{2}}\leq \tan
^{2}\varphi _{2},\qquad k\in \left\{ 1,\dots ,n\right\}
\end{equation*}%
from where we get%
\begin{equation*}
\frac{b_{k}^{2}+a_{k}^{2}}{a_{k}^{2}}\leq \frac{1}{\cos ^{2}\varphi _{2}}%
,\qquad k\in \left\{ 1,\dots ,n\right\} ,\ \varphi _{2}\in \left( 0,\frac{%
\pi }{2}\right)
\end{equation*}%
and%
\begin{equation*}
\frac{a_{k}^{2}+b_{k}^{2}}{a_{k}^{2}}\leq \frac{1+\tan ^{2}\varphi _{1}}{%
\tan ^{2}\varphi _{1}}=\frac{1}{\sin ^{2}\varphi _{1}},\qquad k\in \left\{
1,\dots ,n\right\} ,\ \varphi _{1}\in \left( 0,\frac{\pi }{2}\right)
\end{equation*}%
giving the inequalities%
\begin{equation*}
\left\vert z_{k}\right\vert \cos \varphi _{2}\leq \func{Re}\left(
z_{k}\right) ,\ \ \left\vert z_{k}\right\vert \sin \varphi _{1}\leq \func{Im}%
\left( z_{k}\right)
\end{equation*}%
for each $k\in \left\{ 1,\dots ,n\right\} .$

Now, applying Theorem \ref{RTt2.1.3} for the complex inner product $\mathbb{C%
}$ endowed with the inner product $\left\langle z,w\right\rangle =z\cdot
\bar{w} $ for $x_{k}=z_{k},$ $r_{1}=\cos \varphi _{2},$ $r_{2}=\sin \varphi
_{1}$ and $e=1,$ we deduce the desired inequality (\ref{RT4.3.3}). The case
of equality is also obvious by Theorem \ref{RTt2.1.3} and the proposition is
proven.
\end{proof}

Another result that has an obvious geometrical interpretation is the
following one.

\begin{proposition}
\label{RTp4.2.3}Let $c\in \mathbb{C}$ with $\left\vert z\right\vert =1$ and $%
\rho _{1},\rho _{2}\in \left( 0,1\right) .$ If $z_{k}\in \mathbb{C}$, $k\in
\left\{ 1,\dots ,n\right\} $ are such that%
\begin{equation}
\left\vert z_{k}-c\right\vert \leq \rho _{1},\ \ \left\vert
z_{k}-ic\right\vert \leq \rho _{2}\text{ \ for each \ }k\in \left\{ 1,\dots
,n\right\} ,  \label{RT4.5.3}
\end{equation}%
then we have the inequality%
\begin{equation}
\sqrt{2-\rho _{1}^{2}-\rho _{2}^{2}}\sum_{k=1}^{n}\left\vert
z_{k}\right\vert \leq \left\vert \sum_{k=1}^{n}z_{k}\right\vert ,
\label{RT4.6.3}
\end{equation}%
with equality if and only if%
\begin{equation}
\sum_{k=1}^{n}z_{k}=\left( \sqrt{1-\rho _{1}^{2}}+i\sqrt{1-\rho _{2}^{2}}%
\right) \left( \sum_{k=1}^{n}\left\vert z_{k}\right\vert \right) c.
\label{RT4.7.3}
\end{equation}
\end{proposition}

The proof is obvious by Corollary \ref{RTc2.2.3} applied for $H=\mathbb{C}$.

\begin{remark}
If we choose $e=1,$ and for $\rho _{1},\rho _{2}\in \left( 0,1\right) $ we
define $\bar{D}\left( 1,\rho _{1}\right) :=\left\{ z\in \mathbb{C}%
|\left\vert z-1\right\vert \leq \rho _{1}\right\} ,$ $\bar{D}\left( i,\rho
_{2}\right) :=\left\{ z\in \mathbb{C}|\left\vert z-i\right\vert \leq \rho
_{2}\right\} ,$ then obviously the intersection%
\begin{equation*}
S_{\rho _{1},\rho _{2}}:=\bar{D}\left( 1,\rho _{1}\right) \cap \bar{D}\left(
i,\rho _{2}\right)
\end{equation*}%
is nonempty if and only if $\rho _{1}+\rho _{2}\geq \sqrt{2}.$

If $z_{k}\in S_{\rho _{1},\rho _{2}}$ for $k\in \left\{ 1,\dots ,n\right\} ,$
then (\ref{RT4.6.3}) holds true. The equality holds in (\ref{RT4.6.3}) if
and only if%
\begin{equation*}
\sum_{k=1}^{n}z_{k}=\left( \sqrt{1-\rho _{1}^{2}}+i\sqrt{1-\rho _{2}^{2}}%
\right) \sum_{k=1}^{n}\left\vert z_{k}\right\vert .
\end{equation*}
\end{remark}

%


\chapter{ Reverses for the Continuous Triangle Inequality}\label{ch4}

\section{Introduction}

Let $f:\left[ a,b\right] \rightarrow \mathbb{K}$, $\mathbb{K}=\mathbb{C}$ or
$\mathbb{R}$ be a Lebesgue integrable function. The following inequality,
which is the continuous version of the \textit{triangle inequality}%
\begin{equation}
\left\vert \int_{a}^{b}f\left( x\right) dx\right\vert \leq
\int_{a}^{b}\left\vert f\left( x\right) \right\vert dx,  \label{RTxx0.1.1}
\end{equation}%
plays a fundamental role in Mathematical Analysis and its applications.

It appears, see \cite[p. 492]{RTCMPF}, that the first reverse inequality for
(\ref{RTxx0.1.1}) was obtained by J. Karamata in his book from 1949, \cite%
{RTCK}. It can be stated as%
\begin{equation}
\cos \theta \int_{a}^{b}\left\vert f\left( x\right) \right\vert dx\leq
\left\vert \int_{a}^{b}f\left( x\right) dx\right\vert  \label{RTxx0.1.2}
\end{equation}%
provided%
\begin{equation*}
-\theta \leq \arg f\left( x\right) \leq \theta ,\ \ x\in \left[ a,b\right]
\end{equation*}%
for given $\theta \in \left( 0,\frac{\pi }{2}\right) .$

This integral inequality is the continuous version of a reverse inequality
for the generalised triangle inequality%
\begin{equation}
\cos \theta \sum_{i=1}^{n}\left\vert z_{i}\right\vert \leq \left\vert
\sum_{i=1}^{n}z_{i}\right\vert ,  \label{RTxx0.1.3}
\end{equation}%
provided%
\begin{equation*}
a-\theta \leq \arg \left( z_{i}\right) \leq a+\theta ,\ \ \text{for \ }i\in
\left\{ 1,\dots ,n\right\} ,
\end{equation*}%
where $a\in \mathbb{R}$ and $\theta \in \left( 0,\frac{\pi }{2}\right) ,$
which, as pointed out in \cite[p. 492]{RTCMPF}, was first discovered by M.
Petrovich in 1917, \cite{RTCP}, and, subsequently rediscovered by other
authors, including J. Karamata \cite[p. 300 -- 301]{RTCK}, H.S. Wilf \cite%
{RTCW}, and in an equivalent form, by M. Marden \cite{RTCM}.

The first to consider the problem in the more general case of Hilbert and
Banach spaces, were J.B. Diaz and F.T. Metcalf \cite{RTCDM} who showed that,
in an inner product space $H$ over the real or complex number field, the
following reverse of the triangle inequality holds%
\begin{equation}
r\sum_{i=1}^{n}\left\Vert x_{i}\right\Vert \leq \left\Vert
\sum_{i=1}^{n}x_{i}\right\Vert ,  \label{RTxx0.1.4}
\end{equation}%
provided%
\begin{equation*}
0\leq r\leq \frac{\func{Re}\left\langle x_{i},a\right\rangle }{\left\Vert
x_{i}\right\Vert },\ \ \ \ \ i\in \left\{ 1,\dots ,n\right\} ,
\end{equation*}%
and $a\in H$ is a unit vector, i.e., $\left\Vert a\right\Vert =1.$ The case
of equality holds in (\ref{RTxx0.1.4}) if and only if%
\begin{equation}
\sum_{i=1}^{n}x_{i}=r\left( \sum_{i=1}^{n}\left\Vert x_{i}\right\Vert
\right) a.  \label{RTxx0.1.5}
\end{equation}

A generalisation of this result for orthonormal families is also well known
\cite{RTCDM}:

Let $a_{1},\dots ,a_{m}$ be $m$ orthonormal vectors in $H.$ Suppose the
vectors $x_{1},\dots ,x_{n}\in H\backslash \left\{ 0\right\} $ satisfy%
\begin{equation*}
0\leq r_{k}\leq \frac{\func{Re}\left\langle x_{i},a_{k}\right\rangle }{%
\left\Vert x_{i}\right\Vert },\ \ \ \ \ i\in \left\{ 1,\dots ,n\right\} ,\
k\in \left\{ 1,\dots ,m\right\} .
\end{equation*}%
Then%
\begin{equation*}
\left( \sum_{k=1}^{m}r_{k}^{2}\right) ^{\frac{1}{2}}\sum_{i=1}^{n}\left\Vert
x_{i}\right\Vert \leq \left\Vert \sum_{i=1}^{n}x_{i}\right\Vert ,
\end{equation*}%
where equality holds if and only if%
\begin{equation*}
\sum_{i=1}^{n}x_{i}=\left( \sum_{i=1}^{n}\left\Vert x_{i}\right\Vert \right)
\sum_{k=1}^{m}r_{k}a_{k}.
\end{equation*}%
The main aim of this chapter is to survey some recent reverses of the
triangle inequality for Bochner integrable functions $f$ with values in
Hilbert spaces and defined on a compact interval $\left[ a,b\right] \subset
\mathbb{R}$. Applications for Lebesgue integrable complex-valued functions
are provided as well.

\section{Multiplicative Reverses}

\subsection{Reverses for a Unit Vector}

We recall that $f\in L\left( \left[ a,b\right] ;H\right) ,$ the space of
Bochner integrable functions with values in a Hilbert space $H,$ if and only
if $f:\left[ a,b\right] \rightarrow H$ is Bochner measurable on $\left[ a,b%
\right] $ and the Lebesgue integral $\int_{a}^{b}\left\Vert f\left( t\right)
\right\Vert dt$ is finite.

The following result holds \cite{RTCSSD1}:

\begin{theorem}[Dragomir, 2004]
\label{RTxxt1.2.1}If $f\in L\left( \left[ a,b\right] ;H\right) $ is such
that there exists a constant $K\geq 1$ and a vector $e\in H,$ $\left\Vert
e\right\Vert =1$ with%
\begin{equation}
\left\Vert f\left( t\right) \right\Vert \leq K\func{Re}\left\langle f\left(
t\right) ,e\right\rangle \ \ \ \text{for a.e. }t\in \left[ a,b\right] ,
\label{RTxx1.2.1}
\end{equation}%
then we have the inequality:%
\begin{equation}
\int_{a}^{b}\left\Vert f\left( t\right) \right\Vert dt\leq K\left\Vert
\int_{a}^{b}f\left( t\right) dt\right\Vert .  \label{RTxx1.2.2}
\end{equation}%
The case of equality holds in (\ref{RTxx1.2.2}) if and only if%
\begin{equation}
\int_{a}^{b}f\left( t\right) dt=\frac{1}{K}\left( \int_{a}^{b}\left\Vert
f\left( t\right) \right\Vert dt\right) e.  \label{RTxx1.2.3}
\end{equation}
\end{theorem}

\begin{proof}
By the Schwarz inequality in inner product spaces, we have%
\begin{align}
\left\Vert \int_{a}^{b}f\left( t\right) dt\right\Vert & =\left\Vert
\int_{a}^{b}f\left( t\right) dt\right\Vert \left\Vert e\right\Vert
\label{RTxx1.2.4} \\
& \geq \left\vert \left\langle \int_{a}^{b}f\left( t\right)
dt,e\right\rangle \right\vert \geq \left\vert \func{Re}\left\langle
\int_{a}^{b}f\left( t\right) dt,e\right\rangle \right\vert  \notag \\
& \geq \func{Re}\left\langle \int_{a}^{b}f\left( t\right) dt,e\right\rangle
=\int_{a}^{b}\func{Re}\left\langle f\left( t\right) ,e\right\rangle dt.
\notag
\end{align}%
From the condition (\ref{RTxx1.2.1}), on integrating over $\left[ a,b\right]
,$ we deduce%
\begin{equation}
\int_{a}^{b}\func{Re}\left\langle f\left( t\right) ,e\right\rangle dt\geq
\frac{1}{K}\int_{a}^{b}\left\Vert f\left( t\right) \right\Vert dt,
\label{RTxx1.2.5}
\end{equation}%
and thus, on making use of (\ref{RTxx1.2.4}) and (\ref{RTxx1.2.5}), we
obtain the desired inequality (\ref{RTxx1.2.2}).

If (\ref{RTxx1.2.3}) holds true, then, obviously%
\begin{equation*}
K\left\Vert \int_{a}^{b}f\left( t\right) dt\right\Vert =\left\Vert
e\right\Vert \int_{a}^{b}\left\Vert f\left( t\right) \right\Vert
dt=\int_{a}^{b}\left\Vert f\left( t\right) \right\Vert dt,
\end{equation*}%
showing that (\ref{RTxx1.2.2}) holds with equality.

If we assume that the equality holds in (\ref{RTxx1.2.2}), then by the
argument provided at the beginning of our proof, we must have equality in
each of the inequalities from (\ref{RTxx1.2.4}) and (\ref{RTxx1.2.5}).

Observe that in Schwarz's inequality $\left\Vert x\right\Vert \left\Vert
y\right\Vert \geq \func{Re}\left\langle x,y\right\rangle ,$ $x,y\in H,$ the
case of equality holds if and only if there exists a positive scalar $\mu $
such that $x=\mu e.$ Therefore, equality holds in the first inequality in (%
\ref{RTxx1.2.4}) iff $\int_{a}^{b}f\left( t\right) dt=\lambda e$, with $%
\lambda \geq 0$ $.$

If we assume that a strict inequality holds in (\ref{RTxx1.2.1}) on a subset
of nonzero Lebesgue measure in $\left[ a,b\right] ,$ then $$%
\int_{a}^{b}\left\Vert f\left( t\right) \right\Vert dt<K\int_{a}^{b}\func{Re}%
\left\langle f\left( t\right) ,e\right\rangle dt,$$ and by
(\ref{RTxx1.2.4}) we deduce a strict inequality in
(\ref{RTxx1.2.2}), which contradicts the
assumption. Thus, we must have $\left\Vert f\left( t\right) \right\Vert =K%
\func{Re}\left\langle f\left( t\right) ,e\right\rangle $ for a.e. $t\in %
\left[ a,b\right] .$

If we integrate this equality, we deduce%
\begin{align*}
\int_{a}^{b}\left\Vert f\left( t\right) \right\Vert dt& =K\int_{a}^{b}\func{%
Re}\left\langle f\left( t\right) ,e\right\rangle dt=K\func{Re}\left\langle
\int_{a}^{b}f\left( t\right) dt,e\right\rangle \\
& =K\func{Re}\left\langle \lambda e,e\right\rangle =\lambda K
\end{align*}%
giving%
\begin{equation*}
\lambda =\frac{1}{K}\int_{a}^{b}\left\Vert f\left( t\right) \right\Vert dt,
\end{equation*}%
and thus the equality (\ref{RTxx1.2.3}) is necessary.

This completes the proof.
\end{proof}

A more appropriate result from an applications point of view is perhaps the
following result \cite{RTCSSD1}.

\begin{corollary}
\label{RTxxc1.2.2}Let $e$ be a unit vector in the Hilbert space $\left(
H;\left\langle \cdot ,\cdot \right\rangle \right) ,$ $\rho \in \left(
0,1\right) $ and $f\in L\left( \left[ a,b\right] ;H\right) $ so that%
\begin{equation}
\left\Vert f\left( t\right) -e\right\Vert \leq \rho \ \ \ \text{for a.e. }%
t\in \left[ a,b\right] .  \label{RTxx1.2.6}
\end{equation}%
Then we have the inequality%
\begin{equation}
\sqrt{1-\rho ^{2}}\int_{a}^{b}\left\Vert f\left( t\right) \right\Vert dt\leq
\left\Vert \int_{a}^{b}f\left( t\right) dt\right\Vert ,  \label{RTxx1.2.7}
\end{equation}%
with equality if and only if%
\begin{equation}
\int_{a}^{b}f\left( t\right) dt=\sqrt{1-\rho ^{2}}\left(
\int_{a}^{b}\left\Vert f\left( t\right) \right\Vert dt\right) e.
\label{RTxx1.2.8}
\end{equation}
\end{corollary}

\begin{proof}
From (\ref{RTxx1.2.6}), we have%
\begin{equation*}
\left\Vert f\left( t\right) \right\Vert ^{2}-2\func{Re}\left\langle f\left(
t\right) ,e\right\rangle +1\leq \rho ^{2},
\end{equation*}%
giving%
\begin{equation*}
\left\Vert f\left( t\right) \right\Vert ^{2}+1-\rho ^{2}\leq 2\func{Re}%
\left\langle f\left( t\right) ,e\right\rangle
\end{equation*}%
for a.e. $t\in \left[ a,b\right] .$

Dividing by $\sqrt{1-\rho ^{2}}>0,$ we deduce%
\begin{equation}
\frac{\left\Vert f\left( t\right) \right\Vert ^{2}}{\sqrt{1-\rho ^{2}}}+%
\sqrt{1-\rho ^{2}}\leq \frac{2\func{Re}\left\langle f\left( t\right)
,e\right\rangle }{\sqrt{1-\rho ^{2}}}  \label{RTxx1.2.9}
\end{equation}%
for a.e. $t\in \left[ a,b\right] .$

On the other hand, by the elementary inequality%
\begin{equation*}
\frac{p}{\alpha }+q\alpha \geq 2\sqrt{pq},\ \ \ p,q\geq 0,\ \alpha >0
\end{equation*}%
we have%
\begin{equation}
2\left\Vert f\left( t\right) \right\Vert \leq \frac{\left\Vert f\left(
t\right) \right\Vert ^{2}}{\sqrt{1-\rho ^{2}}}+\sqrt{1-\rho ^{2}}
\label{RTxx1.2.10}
\end{equation}%
for each $t\in \left[ a,b\right] .$

Making use of (\ref{RTxx1.2.9}) and (\ref{RTxx1.2.10}), we deduce%
\begin{equation*}
\left\Vert f\left( t\right) \right\Vert \leq \frac{1}{\sqrt{1-\rho ^{2}}}%
\func{Re}\left\langle f\left( t\right) ,e\right\rangle
\end{equation*}%
for a.e. $t\in \left[ a,b\right] .$

Applying Theorem \ref{RTxxt1.2.1} for $K=\frac{1}{\sqrt{1-\rho ^{2}}},$ we
deduce the desired inequality (\ref{RTxx1.2.7}).
\end{proof}

In the same spirit, we also have the following corollary \cite{RTCSSD1}.

\begin{corollary}
\label{RTxxc1.2.3}Let $e$ be a unit vector in $H$ and $M\geq m>0.$ If $f\in
L\left( \left[ a,b\right] ;H\right) $ is such that%
\begin{equation}
\func{Re}\left\langle Me-f\left( t\right) ,f\left( t\right) -me\right\rangle
\geq 0\ \ \   \label{RTxx1.2.9a}
\end{equation}%
or, equivalently,%
\begin{equation}
\left\Vert f\left( t\right) -\frac{M+m}{2}e\right\Vert \leq \frac{1}{2}%
\left( M-m\right) \ \ \   \label{RTxx1.2.10a}
\end{equation}%
for a.e. $t\in \left[ a,b\right] ,$ then we have the inequality%
\begin{equation}
\frac{2\sqrt{mM}}{M+m}\int_{a}^{b}\left\Vert f\left( t\right) \right\Vert
dt\leq \left\Vert \int_{a}^{b}f\left( t\right) dt\right\Vert ,
\label{RTxx1.2.11}
\end{equation}%
or, equivalently,%
\begin{align}
(0& \leq )\int_{a}^{b}\left\Vert f\left( t\right) \right\Vert dt-\left\Vert
\int_{a}^{b}f\left( t\right) dt\right\Vert  \label{RTxx1.2.12} \\
& \leq \frac{\left( \sqrt{M}-\sqrt{m}\right) ^{2}}{M+m}\left\Vert
\int_{a}^{b}f\left( t\right) dt\right\Vert .  \notag
\end{align}%
The equality holds in (\ref{RTxx1.2.11}) (or in the second part of (\ref%
{RTxx1.2.12})) if and only if%
\begin{equation}
\int_{a}^{b}f\left( t\right) dt=\frac{2\sqrt{mM}}{M+m}\left(
\int_{a}^{b}\left\Vert f\left( t\right) \right\Vert dt\right) e.
\label{RTxx1.2.13}
\end{equation}
\end{corollary}

\begin{proof}
Firstly, we remark that if $x,z,Z\in H,$ then the following statements are
equivalent

\begin{enumerate}
\item[(i)] $\func{Re}\left\langle Z-x,x-z\right\rangle \geq 0$

and

\item[(ii)] $\left\Vert x-\frac{Z+z}{2}\right\Vert \leq \frac{1}{2}%
\left\Vert Z-z\right\Vert .$
\end{enumerate}

Using this fact, we may simply realise that (\ref{RTxx1.2.9}) and (\ref%
{RTxx1.2.10}) are equivalent.

Now, from (\ref{RTxx1.2.9}), we obtain%
\begin{equation*}
\left\Vert f\left( t\right) \right\Vert ^{2}+mM\leq \left( M+m\right) \func{%
Re}\left\langle f\left( t\right) ,e\right\rangle
\end{equation*}%
for a.e. $t\in \left[ a,b\right] .$ Dividing this inequality with $\sqrt{mM}%
>0,$ we deduce the following inequality that will be used in the sequel%
\begin{equation}
\frac{\left\Vert f\left( t\right) \right\Vert ^{2}}{\sqrt{mM}}+\sqrt{mM}\leq
\frac{M+m}{\sqrt{mM}}\func{Re}\left\langle f\left( t\right) ,e\right\rangle
\label{RTxx1.2.14}
\end{equation}%
for a.e. $t\in \left[ a,b\right] .$

On the other hand%
\begin{equation}
2\left\Vert f\left( t\right) \right\Vert \leq \frac{\left\Vert f\left(
t\right) \right\Vert ^{2}}{\sqrt{mM}}+\sqrt{mM},  \label{RTxx1.2.15}
\end{equation}%
for any $t\in \left[ a,b\right] .$

Utilising (\ref{RTxx1.2.14}) and (\ref{RTxx1.2.15}), we may conclude with
the following inequality
\begin{equation*}
\left\Vert f\left( t\right) \right\Vert \leq \frac{M+m}{2\sqrt{mM}}\func{Re}%
\left\langle f\left( t\right) ,e\right\rangle ,
\end{equation*}%
for a.e. $t\in \left[ a,b\right] .$

Applying Theorem \ref{RTxxt1.2.1} for the constant $K:=\frac{m+M}{2\sqrt{mM}}%
\geq 1,$ we deduce the desired result.
\end{proof}

\subsection{Reverses for Orthonormal Families of Vectors}

The following result for orthonormal vectors in $H$ holds \cite{RTCSSD1}.

\begin{theorem}[Dragomir, 2004]
\label{RTxxt1.3.1}Let $\left\{ e_{1},\dots ,e_{n}\right\} $ be a family of
orthonormal vectors in $H,$ $k_{i}\geq 0,$ $i\in \left\{ 1,\dots ,n\right\} $
and $f\in L\left( \left[ a,b\right] ;H\right) $ such that%
\begin{equation}
k_{i}\left\Vert f\left( t\right) \right\Vert \leq \func{Re}\left\langle
f\left( t\right) ,e_{i}\right\rangle   \label{RTxx1.3.1}
\end{equation}%
for each $i\in \left\{ 1,\dots ,n\right\} $ and for a.e. $t\in \left[ a,b%
\right] .$

Then%
\begin{equation}
\left( \sum_{i=1}^{n}k_{i}^{2}\right) ^{\frac{1}{2}}\int_{a}^{b}\left\Vert
f\left( t\right) \right\Vert dt\leq \left\Vert \int_{a}^{b}f\left( t\right)
dt\right\Vert ,  \label{RTxx1.3.2}
\end{equation}%
where equality holds if and only if%
\begin{equation}
\int_{a}^{b}f\left( t\right) dt=\left( \int_{a}^{b}\left\Vert f\left(
t\right) \right\Vert dt\right) \sum_{i=1}^{n}k_{i}e_{i}.  \label{RTxx1.3.3}
\end{equation}
\end{theorem}

\begin{proof}
By Bessel's inequality applied for $\int_{a}^{b}f\left( t\right) dt$ and the
orthonormal vectors $\left\{ e_{1},\dots ,e_{n}\right\} ,$ we have%
\begin{align}
\left\Vert \int_{a}^{b}f\left( t\right) dt\right\Vert ^{2}& \geq
\sum_{i=1}^{n}\left\vert \left\langle \int_{a}^{b}f\left( t\right)
dt,e_{i}\right\rangle \right\vert ^{2}  \label{RTxx1.3.4} \\
& \geq \sum_{i=1}^{n}\left[ \func{Re}\left\langle \int_{a}^{b}f\left(
t\right) dt,e_{i}\right\rangle \right] ^{2}  \notag \\
& =\sum_{i=1}^{n}\left[ \int_{a}^{b}\func{Re}\left\langle f\left( t\right)
,e_{i}\right\rangle dt\right] ^{2}.  \notag
\end{align}%
Integrating (\ref{RTxx1.3.1}), we get for each $i\in \left\{ 1,\dots
,n\right\} $%
\begin{equation*}
0\leq k_{i}\int_{a}^{b}\left\Vert f\left( t\right) \right\Vert dt\leq
\int_{a}^{b}\func{Re}\left\langle f\left( t\right) ,e_{i}\right\rangle dt,
\end{equation*}%
implying%
\begin{equation}
\sum_{i=1}^{n}\left[ \int_{a}^{b}\func{Re}\left\langle f\left( t\right)
,e_{i}\right\rangle dt\right] ^{2}\geq \sum_{i=1}^{n}k_{i}^{2}\left(
\int_{a}^{b}\left\Vert f\left( t\right) \right\Vert dt\right) ^{2}.
\label{RTxx1.3.5}
\end{equation}%
On making use of (\ref{RTxx1.3.4}) and (\ref{RTxx1.3.5}), we deduce%
\begin{equation*}
\left\Vert \int_{a}^{b}f\left( t\right) dt\right\Vert ^{2}\geq
\sum_{i=1}^{n}k_{i}^{2}\left( \int_{a}^{b}\left\Vert f\left( t\right)
\right\Vert dt\right) ^{2},
\end{equation*}%
which is clearly equivalent to (\ref{RTxx1.3.2}).

If (\ref{RTxx1.3.3}) holds true, then%
\begin{align*}
\left\Vert \int_{a}^{b}f\left( t\right) dt\right\Vert & =\left(
\int_{a}^{b}\left\Vert f\left( t\right) \right\Vert dt\right) \left\Vert
\sum_{i=1}^{n}k_{i}e_{i}\right\Vert \\
& =\left( \int_{a}^{b}\left\Vert f\left( t\right) \right\Vert dt\right)
\left[ \sum_{i=1}^{n}k_{i}^{2}\left\Vert e_{i}\right\Vert ^{2}\right] ^{%
\frac{1}{2}} \\
& =\left( \sum_{i=1}^{n}k_{i}^{2}\right) ^{\frac{1}{2}}\int_{a}^{b}\left%
\Vert f\left( t\right) \right\Vert dt,
\end{align*}%
showing that (\ref{RTxx1.3.2}) holds with equality.

Now, suppose that there is an $i_{0}\in \left\{ 1,\dots ,n\right\} $ for
which%
\begin{equation*}
k_{i_{0}}\left\Vert f\left( t\right) \right\Vert <\func{Re}\left\langle
f\left( t\right) ,e_{i_{0}}\right\rangle
\end{equation*}%
on a subset of nonzero Lebesgue measure in $\left[ a,b\right] .$ Then
obviously%
\begin{equation*}
k_{i_{0}}\int_{a}^{b}\left\Vert f\left( t\right) \right\Vert dt<\int_{a}^{b}%
\func{Re}\left\langle f\left( t\right) ,e_{i_{0}}\right\rangle dt,
\end{equation*}%
and using the argument given above, we deduce%
\begin{equation*}
\left( \sum_{i=1}^{n}k_{i}^{2}\right) ^{\frac{1}{2}}\int_{a}^{b}\left\Vert
f\left( t\right) \right\Vert dt<\left\Vert \int_{a}^{b}f\left( t\right)
dt\right\Vert .
\end{equation*}%
Therefore, if the equality holds in (\ref{RTxx1.3.2}), we must have%
\begin{equation}
k_{i}\left\Vert f\left( t\right) \right\Vert =\func{Re}\left\langle f\left(
t\right) ,e_{i}\right\rangle   \label{RTxx1.3.6}
\end{equation}%
for each $i\in \left\{ 1,\dots ,n\right\} $ and a.e. $t\in \left[ a,b\right]
.$

Also, if the equality holds in (\ref{RTxx1.3.2}), then we must have equality
in all inequalities (\ref{RTxx1.3.4}), this means that%
\begin{equation}
\int_{a}^{b}f\left( t\right) dt=\sum_{i=1}^{n}\left\langle
\int_{a}^{b}f\left( t\right) dt,e_{i}\right\rangle e_{i}  \label{RTxx1.3.7}
\end{equation}%
and
\begin{equation}
\func{Im}\left\langle \int_{a}^{b}f\left( t\right) dt,e_{i}\right\rangle =0%
\text{ \ for each \ }i\in \left\{ 1,\dots ,n\right\} .  \label{RTxx1.3.8}
\end{equation}%
Using (\ref{RTxx1.3.6}) and (\ref{RTxx1.3.8}) in (\ref{RTxx1.3.7}), we deduce%
\begin{align*}
\int_{a}^{b}f\left( t\right) dt& =\sum_{i=1}^{n}\func{Re}\left\langle
\int_{a}^{b}f\left( t\right) dt,e_{i}\right\rangle e_{i} \\
& =\sum_{i=1}^{n}\int_{a}^{b}\func{Re}\left\langle f\left( t\right)
,e_{i}\right\rangle e_{i}dt \\
& =\sum_{i=1}^{n}\left( \int_{a}^{b}\left\Vert f\left( t\right) \right\Vert
dt\right) k_{i}e_{i} \\
& =\int_{a}^{b}\left\Vert f\left( t\right) \right\Vert
dt\sum_{i=1}^{n}k_{i}e_{i},
\end{align*}%
and the condition (\ref{RTxx1.3.3}) is necessary.

This completes the proof.
\end{proof}

The following two corollaries are of interest \cite{RTCSSD1}.

\begin{corollary}
\label{RTxxc1.3.2}Let $\left\{ e_{1},\dots ,e_{n}\right\} $ be a family of
orthonormal vectors in $H,$ $\rho _{i}\in \left( 0,1\right) ,$ $i\in \left\{
1,\dots ,n\right\} $ and $f\in L\left( \left[ a,b\right] ;H\right) $ such
that:%
\begin{equation}
\left\Vert f\left( t\right) -e_{i}\right\Vert \leq \rho _{i}\text{ \ for \ }%
i\in \left\{ 1,\dots ,n\right\} \text{ \ and a.e. }t\in \left[ a,b\right] .
\label{RTxx1.3.9}
\end{equation}%
Then we have the inequality%
\begin{equation*}
\left( n-\sum_{i=1}^{n}\rho _{i}^{2}\right) ^{\frac{1}{2}}\int_{a}^{b}\left%
\Vert f\left( t\right) \right\Vert dt\leq \left\Vert \int_{a}^{b}f\left(
t\right) dt\right\Vert ,
\end{equation*}%
with equality if and only if%
\begin{equation*}
\int_{a}^{b}f\left( t\right) dt=\int_{a}^{b}\left\Vert f\left( t\right)
\right\Vert dt\sum_{i=1}^{n}\left( 1-\rho _{i}^{2}\right) ^{1/2}e_{i}.
\end{equation*}
\end{corollary}

\begin{proof}
From the proof of Theorem \ref{RTxxt1.2.1}, we know that (\ref{RTxx1.3.3})
implies the inequality%
\begin{equation*}
\sqrt{1-\rho _{i}^{2}}\left\Vert f\left( t\right) \right\Vert \leq \func{Re}%
\left\langle f\left( t\right) ,e_{i}\right\rangle ,\ \ \ i\in \left\{
1,\dots ,n\right\} ,\text{ \ for a.e. }t\in \left[ a,b\right] .
\end{equation*}%
Now, applying Theorem \ref{RTxxt1.3.1} for $k_{i}:=\sqrt{1-\rho _{i}^{2}},$ $%
i\in \left\{ 1,\dots ,n\right\} $, we deduce the desired result.
\end{proof}

A different results is incorporated in (see \cite{RTCSSD1}):

\begin{corollary}
\label{RTxxc1.3.3}Let $\left\{ e_{1},\dots ,e_{n}\right\} $ be a family of
orthonormal vectors in $H,$ $M_{i}\geq m_{i}>0,$ $i\in \left\{ 1,\dots
,n\right\} $ and $f\in L\left( \left[ a,b\right] ;H\right) $ such that%
\begin{equation}
\func{Re}\left\langle M_{i}e_{i}-f\left( t\right) ,f\left( t\right)
-m_{i}e_{i}\right\rangle \geq 0  \label{RTxx1.3.10}
\end{equation}%
or, equivalently,%
\begin{equation*}
\left\Vert f\left( t\right) -\frac{M_{i}+m_{i}}{2}e_{i}\right\Vert \leq
\frac{1}{2}\left( M_{i}-m_{i}\right)
\end{equation*}%
for \ $i\in \left\{ 1,\dots ,n\right\} $ \ and a.e. $t\in \left[ a,b\right]
. $ Then we have the reverse of the continuous triangle inequality%
\begin{equation*}
\left[ \sum_{i=1}^{n}\frac{4m_{i}M_{i}}{\left( m_{i}+M_{i}\right) ^{2}}%
\right] ^{\frac{1}{2}}\int_{a}^{b}\left\Vert f\left( t\right) \right\Vert
dt\leq \left\Vert \int_{a}^{b}f\left( t\right) dt\right\Vert ,
\end{equation*}%
with equality if and only if%
\begin{equation*}
\int_{a}^{b}f\left( t\right) dt=\int_{a}^{b}\left\Vert f\left( t\right)
\right\Vert dt\left( \sum_{i=1}^{n}\frac{2\sqrt{m_{i}M_{i}}}{m_{i}+M_{i}}%
e_{i}\right) .
\end{equation*}
\end{corollary}

\begin{proof}
From the proof of Corollary \ref{RTxxc1.3.2}, we know (\ref{RTxx1.3.10})
implies that%
\begin{equation*}
\frac{2\sqrt{m_{i}M_{i}}}{m_{i}+M_{i}}\left\Vert f\left( t\right)
\right\Vert \leq \func{Re}\left\langle f\left( t\right) ,e_{i}\right\rangle
,\ \ \ i\in \left\{ 1,\dots ,n\right\} \text{ \ and a.e. }t\in \left[ a,b%
\right] .
\end{equation*}%
Now, applying Theorem \ref{RTxxt1.3.1} for $k_{i}:=\frac{2\sqrt{m_{i}M_{i}}}{%
m_{i}+M_{i}},$ $i\in \left\{ 1,\dots ,n\right\} ,$ we deduce the desired
result.
\end{proof}

\section{Some Additive Reverses}

\subsection{The Case of a Unit Vector}

The following result holds \cite{RTCSSD2}.

\begin{theorem}[Dragomir, 2004]
\label{RTxxt2.2.1}If $f\in L\left( \left[ a,b\right] ;H\right) $ is such
that there exists a vector $e\in H,$ $\left\Vert e\right\Vert =1$ and $k:%
\left[ a,b\right] \rightarrow \lbrack 0,\infty ),$ a Lebesgue integrable
function with%
\begin{equation}
\left\Vert f\left( t\right) \right\Vert -\func{Re}\left\langle f\left(
t\right) ,e\right\rangle \leq k\left( t\right) \ \ \ \text{for a.e. }t\in %
\left[ a,b\right] ,  \label{RTxx2.2.1}
\end{equation}%
then we have the inequality:%
\begin{equation}
\left( 0\leq \right) \int_{a}^{b}\left\Vert f\left( t\right) \right\Vert
dt-\left\Vert \int_{a}^{b}f\left( t\right) dt\right\Vert \leq
\int_{a}^{b}k\left( t\right) dt.  \label{RTxx2.2.2}
\end{equation}%
The equality holds in (\ref{RTxx2.2.2}) if and only if%
\begin{equation}
\int_{a}^{b}\left\Vert f\left( t\right) \right\Vert dt\geq
\int_{a}^{b}k\left( t\right) dt  \label{RTxx2.2.3}
\end{equation}%
and%
\begin{equation}
\int_{a}^{b}f\left( t\right) dt=\left( \int_{a}^{b}\left\Vert f\left(
t\right) \right\Vert dt-\int_{a}^{b}k\left( t\right) dt\right) e.
\label{RTxx2.2.4}
\end{equation}
\end{theorem}

\begin{proof}
If we integrate the inequality (\ref{RTxx2.2.1}), we get%
\begin{equation}
\int_{a}^{b}\left\Vert f\left( t\right) \right\Vert dt\leq \func{Re}%
\left\langle \int_{a}^{b}f\left( t\right) dt,e\right\rangle
+\int_{a}^{b}k\left( t\right) dt.  \label{RTxx2.2.5}
\end{equation}%
By Schwarz's inequality for $e$ and $\int_{a}^{b}f\left( t\right) dt,$ we
have%
\begin{align}
& \func{Re}\left\langle \int_{a}^{b}f\left( t\right) dt,e\right\rangle
\label{RTxx2.2.6} \\
& \leq \left\vert \func{Re}\left\langle \int_{a}^{b}f\left( t\right)
dt,e\right\rangle \right\vert \leq \left\vert \left\langle
\int_{a}^{b}f\left( t\right) dt,e\right\rangle \right\vert  \notag \\
& \leq \left\Vert \int_{a}^{b}f\left( t\right) dt\right\Vert \left\Vert
e\right\Vert =\left\Vert \int_{a}^{b}f\left( t\right) dt\right\Vert .  \notag
\end{align}%
Making use of (\ref{RTxx2.2.5}) and (\ref{RTxx2.2.6}), we deduce the desired
inequality (\ref{RTxx2.2.2}).

If (\ref{RTxx2.2.3}) and (\ref{RTxx2.2.4}) hold true, then%
\begin{align*}
\left\Vert \int_{a}^{b}f\left( t\right) dt\right\Vert & =\left\vert
\int_{a}^{b}\left\Vert f\left( t\right) \right\Vert dt-\int_{a}^{b}k\left(
t\right) dt\right\vert \left\Vert e\right\Vert \\
& =\int_{a}^{b}\left\Vert f\left( t\right) \right\Vert
dt-\int_{a}^{b}k\left( t\right) dt
\end{align*}%
and the equality holds true in (\ref{RTxx2.2.2}).

Conversely, if the equality holds in (\ref{RTxx2.2.2}), then, obviously (\ref%
{RTxx2.2.3}) is valid and we need only to prove (\ref{RTxx2.2.4}).

If $\left\Vert f\left( t\right) \right\Vert -\func{Re}\left\langle f\left(
t\right) ,e\right\rangle <k\left( t\right) $ on a subset of nonzero Lebesgue
measure in $\left[ a,b\right] ,$ then (\ref{RTxx2.2.5}) holds as a strict
inequality, implying that (\ref{RTxx2.2.2}) also holds as a strict
inequality. Therefore, if we assume that equality holds in (\ref{RTxx2.2.2}%
), then we must have%
\begin{equation}
\left\Vert f\left( t\right) \right\Vert =\func{Re}\left\langle f\left(
t\right) ,e\right\rangle +k\left( t\right) \ \ \text{for a.e.\ \ }t\in \left[
a,b\right] .  \label{RTxx2.2.7}
\end{equation}

It is well known that in Schwarz's inequality $\left\Vert x\right\Vert
\left\Vert y\right\Vert \geq \func{Re}\left\langle x,y\right\rangle $ the
equality holds iff there exists a $\lambda \geq 0$ such that $x=\lambda y.$
Therefore, if we assume that the equality holds in all of (\ref{RTxx2.2.6}),
then there exists a $\lambda \geq 0$ such that%
\begin{equation}
\int_{a}^{b}f\left( t\right) dt=\lambda e.  \label{RTxx2.2.8}
\end{equation}%
Integrating (\ref{RTxx2.2.7}) on $\left[ a,b\right] ,$ we deduce%
\begin{equation*}
\int_{a}^{b}\left\Vert f\left( t\right) \right\Vert dt=\func{Re}\left\langle
\int_{a}^{b}f\left( t\right) dt,e\right\rangle +\int_{a}^{b}k\left( t\right)
dt,
\end{equation*}%
and thus, by (\ref{RTxx2.2.8}), we get%
\begin{equation*}
\int_{a}^{b}\left\Vert f\left( t\right) \right\Vert dt=\lambda \left\Vert
e\right\Vert ^{2}+\int_{a}^{b}k\left( t\right) dt,
\end{equation*}%
giving $\lambda =\int_{a}^{b}\left\Vert f\left( t\right) \right\Vert
dt-\int_{a}^{b}k\left( t\right) dt.$

Using (\ref{RTxx2.2.8}), we deduce (\ref{RTxx2.2.4}) and the theorem is
completely proved.
\end{proof}

The following corollary may be useful for applications \cite{RTCSSD2}.

\begin{corollary}
\label{RTxxc2.2.2}If $f\in L\left( \left[ a,b\right] ;H\right) $ is such
that there exists a vector $e\in H,$ $\left\Vert e\right\Vert =1$ and $\rho
\in \left( 0,1\right) $ such that
\begin{equation}
\left\Vert f\left( t\right) -e\right\Vert \leq \rho \ \ \ \text{for a.e. }%
t\in \left[ a,b\right] ,  \label{RTxx2.2.9}
\end{equation}%
then we have the inequality%
\begin{align}
(0& \leq )\int_{a}^{b}\left\Vert f\left( t\right) \right\Vert dt-\left\Vert
\int_{a}^{b}f\left( t\right) dt\right\Vert  \label{RTxx2.2.10} \\
& \leq \frac{\rho ^{2}}{\sqrt{1-\rho ^{2}}\left( 1+\sqrt{1-\rho ^{2}}\right)
}\func{Re}\left\langle \int_{a}^{b}f\left( t\right) dt,e\right\rangle  \notag
\\
& \left( \leq \frac{\rho ^{2}}{\sqrt{1-\rho ^{2}}\left( 1+\sqrt{1-\rho ^{2}}%
\right) }\left\Vert \int_{a}^{b}f\left( t\right) dt\right\Vert \right) .
\notag
\end{align}%
The equality holds in (\ref{RTxx2.2.10}) if and only if%
\begin{equation}
\int_{a}^{b}\left\Vert f\left( t\right) \right\Vert dt\geq \frac{\rho ^{2}}{%
\sqrt{1-\rho ^{2}}\left( 1+\sqrt{1-\rho ^{2}}\right) }\func{Re}\left\langle
\int_{a}^{b}f\left( t\right) dt,e\right\rangle  \label{RTxx2.2.11}
\end{equation}%
and%
\begin{multline}
\int_{a}^{b}f\left( t\right) dt  \label{RTxx2.2.12} \\
=\left( \int_{a}^{b}\left\Vert f\left( t\right) \right\Vert dt-\frac{\rho
^{2}}{\sqrt{1-\rho ^{2}}\left( 1+\sqrt{1-\rho ^{2}}\right) }\func{Re}%
\left\langle \int_{a}^{b}f\left( t\right) dt,e\right\rangle \right) e.
\end{multline}
\end{corollary}

\begin{proof}
Firstly, note that (\ref{RTxx2.2.3}) is equivalent to%
\begin{equation*}
\left\Vert f\left( t\right) \right\Vert ^{2}+1-\rho ^{2}\leq 2\func{Re}%
\left\langle f\left( t\right) ,e\right\rangle ,
\end{equation*}%
giving%
\begin{equation*}
\frac{\left\Vert f\left( t\right) \right\Vert ^{2}}{\sqrt{1-\rho ^{2}}}+%
\sqrt{1-\rho ^{2}}\leq \frac{2\func{Re}\left\langle f\left( t\right)
,e\right\rangle }{\sqrt{1-\rho ^{2}}}
\end{equation*}%
for a.e. $t\in \left[ a,b\right] .$

Since, obviously%
\begin{equation*}
2\left\Vert f\left( t\right) \right\Vert \leq \frac{\left\Vert f\left(
t\right) \right\Vert ^{2}}{\sqrt{1-\rho ^{2}}}+\sqrt{1-\rho ^{2}}
\end{equation*}%
for any $t\in \left[ a,b\right] ,$ then we deduce the inequality%
\begin{equation*}
\left\Vert f\left( t\right) \right\Vert \leq \frac{\func{Re}\left\langle
f\left( t\right) ,e\right\rangle }{\sqrt{1-\rho ^{2}}}\ \ \ \text{for a.e. }%
t\in \left[ a,b\right] ,
\end{equation*}%
which is clearly equivalent to%
\begin{equation*}
\left\Vert f\left( t\right) \right\Vert -\func{Re}\left\langle f\left(
t\right) ,e\right\rangle \leq \frac{\rho ^{2}}{\sqrt{1-\rho ^{2}}\left( 1+%
\sqrt{1-\rho ^{2}}\right) }\func{Re}\left\langle f\left( t\right)
,e\right\rangle
\end{equation*}%
for a.e. $t\in \left[ a,b\right] .$

Applying Theorem \ref{RTxxt2.2.1} for $k\left( t\right) :=\frac{\rho ^{2}}{%
\sqrt{1-\rho ^{2}}\left( 1+\sqrt{1-\rho ^{2}}\right) }\func{Re}\left\langle
f\left( t\right) ,e\right\rangle ,$ we deduce the desired result.
\end{proof}

In the same spirit, we also have the following corollary \cite{RTCSSD2}.

\begin{corollary}
\label{RTxxc2.2.3}If $f\in L\left( \left[ a,b\right] ;H\right) $ is such
that there exists a vector $e\in H,$ $\left\Vert e\right\Vert =1$ and $M\geq
m>0$ such that either%
\begin{equation}
\func{Re}\left\langle Me-f\left( t\right) ,f\left( t\right) -me\right\rangle
\geq 0\ \ \   \label{RTxx2.2.13}
\end{equation}%
or, equivalently,%
\begin{equation}
\left\Vert f\left( t\right) -\frac{M+m}{2}e\right\Vert \leq \frac{1}{2}%
\left( M-m\right) \ \   \label{RTxx2.2.14}
\end{equation}%
$\ $for a.e. $t\in \left[ a,b\right] ,$ then we have the inequality%
\begin{align}
(0& \leq )\int_{a}^{b}\left\Vert f\left( t\right) \right\Vert dt-\left\Vert
\int_{a}^{b}f\left( t\right) dt\right\Vert  \label{RTxx2.2.15} \\
& \leq \frac{\left( \sqrt{M}-\sqrt{m}\right) ^{2}}{2\sqrt{mM}}\func{Re}%
\left\langle \int_{a}^{b}f\left( t\right) dt,e\right\rangle  \notag \\
& \left( \leq \frac{\left( \sqrt{M}-\sqrt{m}\right) ^{2}}{2\sqrt{mM}}%
\left\Vert \int_{a}^{b}f\left( t\right) dt\right\Vert \right) .  \notag
\end{align}%
The equality holds in (\ref{RTxx2.2.15}) if and only if%
\begin{equation*}
\int_{a}^{b}\left\Vert f\left( t\right) \right\Vert dt\geq \frac{\left(
\sqrt{M}-\sqrt{m}\right) ^{2}}{2\sqrt{mM}}\func{Re}\left\langle
\int_{a}^{b}f\left( t\right) dt,e\right\rangle
\end{equation*}%
and%
\begin{equation*}
\int_{a}^{b}f\left( t\right) dt=\left( \int_{a}^{b}\left\Vert f\left(
t\right) \right\Vert dt-\frac{\left( \sqrt{M}-\sqrt{m}\right) ^{2}}{2\sqrt{mM%
}}\func{Re}\left\langle \int_{a}^{b}f\left( t\right) dt,e\right\rangle
\right) e.
\end{equation*}
\end{corollary}

\begin{proof}
Observe that (\ref{RTxx2.2.13}) is clearly equivalent to%
\begin{equation*}
\left\Vert f\left( t\right) \right\Vert ^{2}+mM\leq \left( M+m\right) \func{%
Re}\left\langle f\left( t\right) ,e\right\rangle
\end{equation*}%
for a.e. $t\in \left[ a,b\right] ,$ giving the inequality%
\begin{equation*}
\frac{\left\Vert f\left( t\right) \right\Vert ^{2}}{\sqrt{mM}}+\sqrt{mM}\leq
\frac{M+m}{\sqrt{mM}}\func{Re}\left\langle f\left( t\right) ,e\right\rangle
\end{equation*}%
for a.e. $t\in \left[ a,b\right] .$

Since, obviously,%
\begin{equation*}
2\left\Vert f\left( t\right) \right\Vert \leq \frac{\left\Vert f\left(
t\right) \right\Vert ^{2}}{\sqrt{mM}}+\sqrt{mM}
\end{equation*}%
for any $t\in \left[ a,b\right] ,$ hence we deduce the inequality%
\begin{equation*}
\left\Vert f\left( t\right) \right\Vert \leq \frac{M+m}{\sqrt{mM}}\func{Re}%
\left\langle f\left( t\right) ,e\right\rangle \ \ \text{for a.e. }t\in \left[
a,b\right] ,
\end{equation*}%
which is clearly equivalent to%
\begin{equation*}
\left\Vert f\left( t\right) \right\Vert -\func{Re}\left\langle f\left(
t\right) ,e\right\rangle \leq \frac{\left( \sqrt{M}-\sqrt{m}\right) ^{2}}{2%
\sqrt{mM}}\func{Re}\left\langle f\left( t\right) ,e\right\rangle
\end{equation*}%
for a.e. $t\in \left[ a,b\right] .$

Finally, applying Theorem \ref{RTxxt2.2.1}, we obtain the desired result.
\end{proof}

We can state now (see also \cite{RTCSSD2}):

\begin{corollary}
\label{RTxxc2.2.4}If $f\in L\left( \left[ a,b\right] ;H\right) $ and $r\in
L_{2}\left( \left[ a,b\right] ;H\right) ,$ $e\in H,$ $\left\Vert
e\right\Vert =1$ are such that
\begin{equation}
\left\Vert f\left( t\right) -e\right\Vert \leq r\left( t\right) \ \ \ \text{%
for a.e. }t\in \left[ a,b\right] ,  \label{RTxx2.2.16}
\end{equation}%
then we have the inequality%
\begin{equation}
\left( 0\leq \right) \int_{a}^{b}\left\Vert f\left( t\right) \right\Vert
dt-\left\Vert \int_{a}^{b}f\left( t\right) dt\right\Vert \leq \frac{1}{2}%
\int_{a}^{b}r^{2}\left( t\right) dt.  \label{RTxx2.2.17}
\end{equation}%
The equality holds in (\ref{RTxx2.2.17}) if and only if%
\begin{equation*}
\int_{a}^{b}\left\Vert f\left( t\right) \right\Vert dt\geq \frac{1}{2}%
\int_{a}^{b}r^{2}\left( t\right) dt
\end{equation*}%
and%
\begin{equation*}
\int_{a}^{b}f\left( t\right) dt=\left( \int_{a}^{b}\left\Vert f\left(
t\right) \right\Vert dt-\frac{1}{2}\int_{a}^{b}r^{2}\left( t\right)
dt\right) e.
\end{equation*}
\end{corollary}

\begin{proof}
The condition (\ref{RTxx2.2.16}) is obviously equivalent to%
\begin{equation*}
\left\Vert f\left( t\right) \right\Vert ^{2}+1\leq 2\func{Re}\left\langle
f\left( t\right) ,e\right\rangle +r^{2}\left( t\right)
\end{equation*}%
for a.e. $t\in \left[ a,b\right] .$

Using the elementary inequality%
\begin{equation*}
2\left\Vert f\left( t\right) \right\Vert \leq \left\Vert f\left( t\right)
\right\Vert ^{2}+1,\ \ t\in \left[ a,b\right] ,
\end{equation*}%
we deduce%
\begin{equation*}
\left\Vert f\left( t\right) \right\Vert -\func{Re}\left\langle f\left(
t\right) ,e\right\rangle \leq \frac{1}{2}r^{2}\left( t\right)
\end{equation*}%
for a.e. $t\in \left[ a,b\right] .$

Applying Theorem \ref{RTxxt2.2.1} for $k\left( t\right) :=\frac{1}{2}%
r^{2}\left( t\right) ,$ $t\in \left[ a,b\right] $, we deduce the desired
result.
\end{proof}

Finally, we may state and prove the following result as well \cite{RTCSSD2}.

\begin{corollary}
\label{RTxxc2.2.5}If $f\in L\left( \left[ a,b\right] ;H\right) $, $e\in H,$ $%
\left\Vert e\right\Vert =1$ and $M,m:\left[ a,b\right] \rightarrow \lbrack
0,\infty )$ with $M\geq m$ a.e. on $\left[ a,b\right] ,$ are such that $%
\frac{\left( M-m\right) ^{2}}{M+m}\in L\left[ a,b\right] $ and either%
\begin{equation}
\left\Vert f\left( t\right) -\frac{M\left( t\right) +m\left( t\right) }{2}%
e\right\Vert \leq \frac{1}{2}\left[ M\left( t\right) -m\left( t\right) %
\right]  \label{RTxx2.2.18}
\end{equation}%
or, equivalently,%
\begin{equation}
\func{Re}\left\langle M\left( t\right) e-f\left( t\right) ,f\left( t\right)
-m\left( t\right) e\right\rangle \geq 0\ \ \   \label{RTxx2.2.19}
\end{equation}%
for a.e. $t\in \left[ a,b\right] ,$ then we have the inequality%
\begin{equation}
\left( 0\leq \right) \int_{a}^{b}\left\Vert f\left( t\right) \right\Vert
dt-\left\Vert \int_{a}^{b}f\left( t\right) dt\right\Vert \leq \frac{1}{4}%
\int_{a}^{b}\frac{\left[ M\left( t\right) -m\left( t\right) \right] ^{2}}{%
M\left( t\right) +m\left( t\right) }dt.  \label{RTxx2.2.20}
\end{equation}%
The equality holds in (\ref{RTxx2.2.20}) if and only if%
\begin{equation*}
\int_{a}^{b}\left\Vert f\left( t\right) \right\Vert dt\geq \frac{1}{4}%
\int_{a}^{b}\frac{\left[ M\left( t\right) -m\left( t\right) \right] ^{2}}{%
M\left( t\right) +m\left( t\right) }dt
\end{equation*}%
and%
\begin{equation*}
\int_{a}^{b}f\left( t\right) dt=\left( \int_{a}^{b}\left\Vert f\left(
t\right) \right\Vert dt-\frac{1}{4}\int_{a}^{b}\frac{\left[ M\left( t\right)
-m\left( t\right) \right] ^{2}}{M\left( t\right) +m\left( t\right) }%
dt\right) e.
\end{equation*}
\end{corollary}

\begin{proof}
The condition (\ref{RTxx2.2.18}) is equivalent to%
\begin{multline*}
\left\Vert f\left( t\right) \right\Vert ^{2}+\left( \frac{M\left( t\right)
+m\left( t\right) }{2}\right) ^{2} \\
\leq 2\left( \frac{M\left( t\right) +m\left( t\right) }{2}\right) \func{Re}%
\left\langle f\left( t\right) ,e\right\rangle +\frac{1}{4}\left[ M\left(
t\right) -m\left( t\right) \right] ^{2}
\end{multline*}%
for a.e. $t\in \left[ a,b\right] ,$ and since%
\begin{equation*}
2\left( \frac{M\left( t\right) +m\left( t\right) }{2}\right) \left\Vert
f\left( t\right) \right\Vert \leq \left\Vert f\left( t\right) \right\Vert
^{2}+\left( \frac{M\left( t\right) +m\left( t\right) }{2}\right) ^{2},\ \ \
t\in \left[ a,b\right]
\end{equation*}%
hence%
\begin{equation*}
\left\Vert f\left( t\right) \right\Vert -\func{Re}\left\langle f\left(
t\right) ,e\right\rangle \leq \frac{1}{4}\frac{\left[ M\left( t\right)
-m\left( t\right) \right] ^{2}}{M\left( t\right) +m\left( t\right) }
\end{equation*}%
for a.e. $t\in \left[ a,b\right] .$

Now, applying Theorem \ref{RTxxt2.2.1} for $k\left( t\right) :=\frac{1}{4}%
\frac{\left[ M\left( t\right) -m\left( t\right) \right] ^{2}}{M\left(
t\right) +m\left( t\right) },$ $t\in \left[ a,b\right] $, we deduce the
desired inequality.
\end{proof}

\subsection{Additive Reverses for Orthonormal Families}

The following reverse of the continuous triangle inequality for vector
valued integrals holds \cite{RTCSSD2}.

\begin{theorem}[Dragomir, 2004]
\label{RTxxt2.3.1}Let $f\in L\left( \left[ a,b\right] ;H\right) ,$ where $H$
is a Hilbert space over the real or complex number field $\mathbb{K}$, $%
\left\{ e_{i}\right\} _{i\in \left\{ 1,\dots ,n\right\} }$ an orthonormal
family in $H$ and $M_{i}\in L\left[ a,b\right] ,$ $i\in \left\{ 1,\dots
,n\right\} .$ If we assume that%
\begin{equation}
\left\Vert f\left( t\right) \right\Vert -\func{Re}\left\langle f\left(
t\right) ,e_{i}\right\rangle \leq M_{i}\left( t\right) \ \ \text{for a.e. }%
t\in \left[ a,b\right] ,  \label{RTxx2.3.1}
\end{equation}%
then we have the inequality%
\begin{equation}
\int_{a}^{b}\left\Vert f\left( t\right) \right\Vert dt\leq \frac{1}{\sqrt{n}}%
\left\Vert \int_{a}^{b}f\left( t\right) dt\right\Vert +\frac{1}{n}%
\sum_{i=1}^{n}\int_{a}^{b}M_{i}\left( t\right) dt.  \label{RTxx2.3.2}
\end{equation}%
The equality holds in (\ref{RTxx2.3.2}) if and only if%
\begin{equation}
\int_{a}^{b}\left\Vert f\left( t\right) \right\Vert dt\geq \frac{1}{n}%
\sum_{i=1}^{n}\int_{a}^{b}M_{i}\left( t\right) dt  \label{RTxx2.3.3}
\end{equation}%
and
\begin{equation}
\int_{a}^{b}f\left( t\right) dt=\left( \int_{a}^{b}\left\Vert f\left(
t\right) \right\Vert dt-\frac{1}{n}\sum_{i=1}^{n}\int_{a}^{b}M_{i}\left(
t\right) dt\right) \sum_{i=1}^{n}e_{i}.  \label{RTxx2.3.4}
\end{equation}
\end{theorem}

\begin{proof}
If we integrate the inequality (\ref{RTxx2.3.1}) on $\left[ a,b\right] ,$ we
get%
\begin{equation*}
\int_{a}^{b}\left\Vert f\left( t\right) \right\Vert dt\leq \func{Re}%
\left\langle \int_{a}^{b}f\left( t\right) dt,e_{i}\right\rangle
+\int_{a}^{b}M_{i}\left( t\right) dt
\end{equation*}%
for each $i\in \left\{ 1,\dots ,n\right\} .$ Summing these inequalities over
$i$ from $1$ to $n,$ we deduce%
\begin{equation}
\int_{a}^{b}\left\Vert f\left( t\right) \right\Vert dt\leq \frac{1}{n}\func{%
Re}\left\langle \int_{a}^{b}f\left( t\right)
dt,\sum_{i=1}^{n}e_{i}\right\rangle +\frac{1}{n}\sum_{i=1}^{n}%
\int_{a}^{b}M_{i}\left( t\right) dt.  \label{RTxx2.3.5}
\end{equation}%
By Schwarz's inequality for $\int_{a}^{b}f\left( t\right) dt$ and $%
\sum_{i=1}^{n}e_{i},$ we have%
\begin{align}
& \func{Re}\left\langle \int_{a}^{b}f\left( t\right)
dt,\sum_{i=1}^{n}e_{i}\right\rangle  \label{RTxx2.3.6} \\
& \leq \left\vert \func{Re}\left\langle \int_{a}^{b}f\left( t\right)
dt,\sum_{i=1}^{n}e_{i}\right\rangle \right\vert \leq \left\vert \left\langle
\int_{a}^{b}f\left( t\right) dt,\sum_{i=1}^{n}e_{i}\right\rangle \right\vert
\notag \\
& \leq \left\Vert \int_{a}^{b}f\left( t\right) dt\right\Vert \left\Vert
\sum_{i=1}^{n}e_{i}\right\Vert =\sqrt{n}\left\Vert \int_{a}^{b}f\left(
t\right) dt\right\Vert ,  \notag
\end{align}%
since%
\begin{equation*}
\left\Vert \sum_{i=1}^{n}e_{i}\right\Vert =\sqrt{\left\Vert
\sum_{i=1}^{n}e_{i}\right\Vert ^{2}}=\sqrt{\sum_{i=1}^{n}\left\Vert
e_{i}\right\Vert ^{2}}=\sqrt{n}.
\end{equation*}%
Making use of (\ref{RTxx2.3.5}) and (\ref{RTxx2.3.6}), we deduce the desired
inequality (\ref{RTxx2.3.2}).

If (\ref{RTxx2.3.3}) and (\ref{RTxx2.3.4}) hold, then
\begin{align*}
\frac{1}{\sqrt{n}}\left\Vert \int_{a}^{b}f\left( t\right) dt\right\Vert & =%
\frac{1}{\sqrt{n}}\left\vert \int_{a}^{b}\left\Vert f\left( t\right)
\right\Vert dt-\frac{1}{n}\sum_{i=1}^{n}\int_{a}^{b}M_{i}\left( t\right)
dt\right\vert \left\Vert \sum_{i=1}^{n}e_{i}\right\Vert \\
& =\left( \int_{a}^{b}\left\Vert f\left( t\right) \right\Vert dt-\frac{1}{n}%
\sum_{i=1}^{n}\int_{a}^{b}M_{i}\left( t\right) dt\right)
\end{align*}%
and the equality in (\ref{RTxx2.3.2}) holds true.

Conversely, if the equality holds in (\ref{RTxx2.3.2}), then, obviously, (%
\ref{RTxx2.3.3}) is valid.

Taking into account the argument presented above for the previous result (%
\ref{RTxx2.3.2}), it is obvious that, if the equality holds in (\ref%
{RTxx2.3.2}), then it must hold in (\ref{RTxx2.3.1}) for a.e. $t\in \left[
a,b\right] $ and for each $i\in \left\{ 1,\dots ,n\right\} $ and also the
equality must hold in any of the inequalities in (\ref{RTxx2.3.6}).

It is well known that in Schwarz's inequality $\func{Re}\left\langle
u,v\right\rangle \leq \left\Vert u\right\Vert \left\Vert v\right\Vert ,$ the
equality occurs if and only if $u=\lambda v$ with $\lambda \geq 0,$
consequently, the equality holds in all inequalities from (\ref{RTxx2.3.6})
simultaneously iff there exists a $\mu \geq 0$ with%
\begin{equation}
\mu \sum_{i=1}^{n}e_{i}=\int_{a}^{b}f\left( t\right) dt.  \label{RTxx2.3.7}
\end{equation}%
If we integrate the equality in (\ref{RTxx2.3.1}) and sum over $i,$ we deduce%
\begin{equation}
n\int_{a}^{b}f\left( t\right) dt=\func{Re}\left\langle \int_{a}^{b}f\left(
t\right) dt,\sum_{i=1}^{n}e_{i}\right\rangle
+\sum_{i=1}^{n}\int_{a}^{b}M_{i}\left( t\right) dt.  \label{RTxx2.3.8}
\end{equation}%
Replacing $\int_{a}^{b}f\left( t\right) dt$ from (\ref{RTxx2.3.7}) into (\ref%
{RTxx2.3.8}), we deduce%
\begin{align}
n\int_{a}^{b}f\left( t\right) dt& =\mu \left\Vert
\sum_{i=1}^{n}e_{i}\right\Vert ^{2}+\sum_{i=1}^{n}\int_{a}^{b}M_{i}\left(
t\right) dt  \label{RTxx2.3.9} \\
& =\mu n+\sum_{i=1}^{n}\int_{a}^{b}M_{i}\left( t\right) dt.  \notag
\end{align}

Finally, we note that (\ref{RTxx2.3.7}) and (\ref{RTxx2.3.9}) will produce
the required identity (\ref{RTxx2.3.4}), and the proof is complete.
\end{proof}

The following corollaries may be of interest for applications \cite{RTCSSD2}.

\begin{corollary}
\label{RTxxc2.3.2}Let $f\in L\left( \left[ a,b\right] ;H\right) ,$ $\left\{
e_{i}\right\} _{i\in \left\{ 1,\dots ,n\right\} }$ an orthonormal family in $%
H$ and $\rho _{i}\in \left( 0,1\right) ,$ $i\in \left\{ 1,\dots ,n\right\} $
such that%
\begin{equation}
\left\Vert f\left( t\right) -e_{i}\right\Vert \leq \rho _{i}\ \ \text{for
a.e. }t\in \left[ a,b\right] .  \label{RTxx2.3.10}
\end{equation}%
Then we have the inequalities:%
\begin{align}
& \int_{a}^{b}\left\Vert f\left( t\right) \right\Vert dt\leq \frac{1}{\sqrt{n%
}}\left\Vert \int_{a}^{b}f\left( t\right) dt\right\Vert  \label{RTxx2.3.11a}
\\
& +\func{Re}\left\langle \int_{a}^{b}f\left( t\right) dt,\frac{1}{n}%
\sum_{i=1}^{n}\frac{\rho _{i}^{2}}{\sqrt{1-\rho _{i}^{2}}\left( 1+\sqrt{%
1-\rho _{i}^{2}}\right) }e_{i}\right\rangle  \notag \\
& \leq \frac{1}{\sqrt{n}}\left\Vert \int_{a}^{b}f\left( t\right)
dt\right\Vert  \notag \\
& \times \left[ 1+\left( \frac{1}{n}\sum_{i=1}^{n}\frac{\rho _{i}^{2}}{\sqrt{%
1-\rho _{i}^{2}}\left( 1+\sqrt{1-\rho _{i}^{2}}\right) }\right) ^{\frac{1}{2}%
}\right] .  \notag
\end{align}%
The equality holds in the first inequality in (\ref{RTxx2.3.11a}) if and
only if%
\begin{equation*}
\int_{a}^{b}\left\Vert f\left( t\right) \right\Vert dt\geq \func{Re}%
\left\langle \int_{a}^{b}f\left( t\right) dt,\frac{1}{n}\sum_{i=1}^{n}\frac{%
\rho _{i}^{2}}{\sqrt{1-\rho _{i}^{2}}\left( 1+\sqrt{1-\rho _{i}^{2}}\right) }%
e_{i}\right\rangle
\end{equation*}%
and%
\begin{multline*}
\int_{a}^{b}f\left( t\right) dt \\
=\left( \int_{a}^{b}\left\Vert f\left( t\right) \right\Vert dt-\func{Re}%
\left\langle \int_{a}^{b}f\left( t\right) dt,\frac{1}{n}\sum_{i=1}^{n}\frac{%
\rho _{i}^{2}}{\sqrt{1-\rho _{i}^{2}}\left( 1+\sqrt{1-\rho _{i}^{2}}\right) }%
e_{i}\right\rangle \right) \\
\times \sum_{i=1}^{n}e_{i}.
\end{multline*}
\end{corollary}

\begin{proof}
As in the proof of Corollary \ref{RTxxc2.2.2}, the assumption (\ref%
{RTxx2.3.10}) implies%
\begin{equation*}
\left\Vert f\left( t\right) \right\Vert -\func{Re}\left\langle f\left(
t\right) ,e_{i}\right\rangle \leq \frac{\rho _{i}^{2}}{\sqrt{1-\rho _{i}^{2}}%
\left( \sqrt{1-\rho _{i}^{2}}+1\right) }\func{Re}\left\langle f\left(
t\right) ,e_{i}\right\rangle
\end{equation*}%
for a.e. $t\in \left[ a,b\right] $ and for each $i\in \left\{ 1,\dots
,n\right\} .$

Now, if we apply Theorem \ref{RTxxt2.3.1} for
\begin{equation*}
M_{i}\left( t\right) :=\frac{\rho _{i}^{2}\func{Re}\left\langle f\left(
t\right) ,e_{i}\right\rangle }{\sqrt{1-\rho _{i}^{2}}\left( \sqrt{1-\rho
_{i}^{2}}+1\right) },\ \ i\in \left\{ 1,\dots ,n\right\} ,\ \ t\in \left[ a,b%
\right] ,
\end{equation*}%
we deduce the first inequality in (\ref{RTxx2.3.11a}).

By Schwarz's inequality in $H,$ we have%
\begin{align*}
& \func{Re}\left\langle \int_{a}^{b}f\left( t\right) dt,\frac{1}{n}%
\sum_{i=1}^{n}\frac{\rho _{i}^{2}}{\sqrt{1-\rho _{i}^{2}}\left( 1+\sqrt{%
1-\rho _{i}^{2}}\right) }e_{i}\right\rangle \\
& \leq \left\Vert \int_{a}^{b}f\left( t\right) dt\right\Vert \left\Vert
\frac{1}{n}\sum_{i=1}^{n}\frac{\rho _{i}^{2}}{\sqrt{1-\rho _{i}^{2}}\left( 1+%
\sqrt{1-\rho _{i}^{2}}\right) }e_{i}\right\Vert \\
& =\frac{1}{n}\left\Vert \int_{a}^{b}f\left( t\right) dt\right\Vert \left(
\sum_{i=1}^{n}\left[ \frac{\rho _{i}^{2}}{\sqrt{1-\rho _{i}^{2}}\left( 1+%
\sqrt{1-\rho _{i}^{2}}\right) }\right] ^{2}\right) ^{\frac{1}{2}},
\end{align*}%
which implies the second inequality in (\ref{RTxx2.3.11a}).
\end{proof}

The second result is incorporated in \cite{RTCSSD2}:

\begin{corollary}
\label{RTxxc2.3.3}Let $f\in L\left( \left[ a,b\right] ;H\right) ,$ $\left\{
e_{i}\right\} _{i\in \left\{ 1,\dots ,n\right\} }$ an orthonormal family in $%
H$ and $M_{i}\geq m_{i}>0$ such that either%
\begin{equation}
\func{Re}\left\langle M_{i}e_{i}-f\left( t\right) ,f\left( t\right)
-m_{i}e_{i}\right\rangle \geq 0\   \label{RTxx2.3.12a}
\end{equation}%
or, equivalently,%
\begin{equation*}
\left\Vert f\left( t\right) -\frac{M_{i}+m_{i}}{2}\cdot e_{i}\right\Vert
\leq \frac{1}{2}\left( M_{i}-m_{i}\right) \
\end{equation*}%
for a.e. $t\in \left[ a,b\right] $ \ and each \ $i\in \left\{ 1,\dots
,n\right\} .$Then we have
\begin{align}
& \int_{a}^{b}\left\Vert f\left( t\right) \right\Vert dt\leq \frac{1}{\sqrt{n%
}}\left\Vert \int_{a}^{b}f\left( t\right) dt\right\Vert  \label{RTxx2.3.13}
\\
& +\func{Re}\left\langle \int_{a}^{b}f\left( t\right) dt,\frac{1}{n}%
\sum_{i=1}^{n}\frac{\left( \sqrt{M_{i}}-\sqrt{m_{i}}\right) ^{2}}{2\sqrt{%
m_{i}M_{i}}}e_{i}\right\rangle  \notag \\
& \leq \frac{1}{\sqrt{n}}\left\Vert \int_{a}^{b}f\left( t\right)
dt\right\Vert \left[ 1+\left( \frac{1}{n}\sum_{i=1}^{n}\frac{\left( \sqrt{%
M_{i}}-\sqrt{m_{i}}\right) ^{4}}{4m_{i}M_{i}}\right) ^{\frac{1}{2}}\right] .
\notag
\end{align}%
The equality holds in the first inequality in (\ref{RTxx2.3.13}) if and only
if%
\begin{equation*}
\int_{a}^{b}\left\Vert f\left( t\right) \right\Vert dt\geq \func{Re}%
\left\langle \int_{a}^{b}f\left( t\right) dt,\frac{1}{n}\sum_{i=1}^{n}\frac{%
\left( \sqrt{M_{i}}-\sqrt{m_{i}}\right) ^{2}}{2\sqrt{m_{i}M_{i}}}%
e_{i}\right\rangle
\end{equation*}%
and%
\begin{multline*}
\int_{a}^{b}f\left( t\right) dt \\
=\left( \int_{a}^{b}\left\Vert f\left( t\right) \right\Vert dt-\func{Re}%
\left\langle \int_{a}^{b}f\left( t\right) dt,\frac{1}{n}\sum_{i=1}^{n}\frac{%
\left( \sqrt{M_{i}}-\sqrt{m_{i}}\right) ^{2}}{2\sqrt{m_{i}M_{i}}}%
e_{i}\right\rangle \right) \\
\times \sum_{i=1}^{n}e_{i}.
\end{multline*}
\end{corollary}

\begin{proof}
As in the proof of Corollary \ref{RTxxc2.2.3}, from (\ref{RTxx2.3.12a}), we
have%
\begin{equation*}
\left\Vert f\left( t\right) \right\Vert -\func{Re}\left\langle f\left(
t\right) ,e_{i}\right\rangle \leq \frac{\left( \sqrt{M_{i}}-\sqrt{m_{i}}%
\right) ^{2}}{2\sqrt{m_{i}M_{i}}}\func{Re}\left\langle f\left( t\right)
,e_{i}\right\rangle
\end{equation*}%
for a.e. $t\in \left[ a,b\right] $ and $i\in \left\{ 1,\dots ,n\right\} .$

Applying Theorem \ref{RTxxt2.3.1} for%
\begin{equation*}
M_{i}\left( t\right) :=\frac{\left( \sqrt{M_{i}}-\sqrt{m_{i}}\right) ^{2}}{2%
\sqrt{m_{i}M_{i}}}\func{Re}\left\langle f\left( t\right) ,e_{i}\right\rangle
,\ \ \text{ }t\in \left[ a,b\right] ,\ i\in \left\{ 1,\dots ,n\right\} ,
\end{equation*}%
we deduce the desired result.
\end{proof}

In a different direction, we may state the following result as well \cite%
{RTCSSD2}.

\begin{corollary}
\label{RTxxc2.3.4}Let $f\in L\left( \left[ a,b\right] ;H\right) ,$ $\left\{
e_{i}\right\} _{i\in \left\{ 1,\dots ,n\right\} }$ an orthonormal family in $%
H$ and $r_{i}\in L^{2}\left( \left[ a,b\right] \right) ,$ $i\in \left\{
1,\dots ,n\right\} $ such that%
\begin{equation*}
\left\Vert f\left( t\right) -e_{i}\right\Vert \leq r_{i}\left( t\right) \ \
\ \text{for a.e. }t\in \left[ a,b\right] \text{ \ and \ }i\in \left\{
1,\dots ,n\right\} .
\end{equation*}%
Then we have the inequality%
\begin{equation}
\int_{a}^{b}\left\Vert f\left( t\right) \right\Vert dt\leq \frac{1}{\sqrt{n}}%
\left\Vert \int_{a}^{b}f\left( t\right) dt\right\Vert +\frac{1}{2n}%
\sum_{i=1}^{n}\left( \int_{a}^{b}r_{i}^{2}\left( t\right) dt\right) .
\label{RTxx2.3.15}
\end{equation}%
The equality holds in (\ref{RTxx2.3.15}) if and only if%
\begin{equation*}
\int_{a}^{b}\left\Vert f\left( t\right) \right\Vert dt\geq \frac{1}{2n}%
\sum_{i=1}^{n}\left( \int_{a}^{b}r_{i}^{2}\left( t\right) dt\right)
\end{equation*}%
and%
\begin{equation*}
\int_{a}^{b}f\left( t\right) dt=\left[ \int_{a}^{b}\left\Vert f\left(
t\right) \right\Vert dt-\frac{1}{n}\sum_{i=1}^{n}\left(
\int_{a}^{b}r_{i}^{2}\left( t\right) dt\right) \right] \sum_{i=1}^{n}e_{i}.
\end{equation*}
\end{corollary}

\begin{proof}
As in the proof of Corollary \ref{RTxxc2.2.4}, from (\ref{RTxx2.2.16}), we
deduce that%
\begin{equation}
\left\Vert f\left( t\right) \right\Vert -\func{Re}\left\langle f\left(
t\right) ,e_{i}\right\rangle \leq \frac{1}{2}r_{i}^{2}\left( t\right)
\label{RTxx2.3.16}
\end{equation}%
for a.e. $t\in \left[ a,b\right] $ and $i\in \left\{ 1,\dots ,n\right\} .$

Applying Theorem \ref{RTxxt2.3.1} for%
\begin{equation*}
M_{i}\left( t\right) :=\frac{1}{2}r_{i}^{2}\left( t\right) ,\ \ t\in \left[
a,b\right] ,\ i\in \left\{ 1,\dots ,n\right\} ,
\end{equation*}%
we get the desired result.
\end{proof}

Finally, the following result holds \cite{RTCSSD2}.

\begin{corollary}
\label{RTxxc2.3.5}Let $f\in L\left( \left[ a,b\right] ;H\right) ,$ $\left\{
e_{i}\right\} _{i\in \left\{ 1,\dots ,n\right\} }$ an orthonormal family in $%
H$, $M_{i},m_{i}:\left[ a,b\right] \rightarrow \lbrack 0,\infty )$ with $%
M_{i}\geq m_{i}$ a.e. on $\left[ a,b\right] $ and $\frac{\left(
M_{i}-m_{i}\right) ^{2}}{M_{i}+m_{i}}\in L\left[ a,b\right] ,$ and either%
\begin{multline}
\left\Vert f\left( t\right) -\frac{M_{i}\left( t\right) +m_{i}\left(
t\right) }{2}e_{i}\right\Vert \leq \frac{1}{2}\left[ M_{i}\left( t\right)
-m_{i}\left( t\right) \right] ^{2}\ \   \label{RTxx2.3.17} \\
\end{multline}%
or, equivalently,%
\begin{equation*}
\func{Re}\left\langle M_{i}\left( t\right) e_{i}-f\left( t\right) ,f\left(
t\right) -m_{i}\left( t\right) e_{i}\right\rangle \geq 0\ \ \
\end{equation*}%
for a.e. $t\in \left[ a,b\right] $ and any\ $i\in \left\{ 1,\dots ,n\right\}
$, then we have the inequality%
\begin{multline}
\int_{a}^{b}\left\Vert f\left( t\right) \right\Vert dt\leq \frac{1}{\sqrt{n}}%
\left\Vert \int_{a}^{b}f\left( t\right) dt\right\Vert  \label{RTxx2.3.18} \\
+\frac{1}{4n}\sum_{i=1}^{n}\left( \int_{a}^{b}\frac{\left[ M_{i}\left(
t\right) -m_{i}\left( t\right) \right] ^{2}}{M_{i}\left( t\right)
+m_{i}\left( t\right) }dt\right) .
\end{multline}%
The equality holds in (\ref{RTxx2.3.18}) if and only if%
\begin{equation*}
\int_{a}^{b}\left\Vert f\left( t\right) \right\Vert dt\geq \frac{1}{4n}%
\sum_{i=1}^{n}\left( \int_{a}^{b}\frac{\left[ M_{i}\left( t\right)
-m_{i}\left( t\right) \right] ^{2}}{M_{i}\left( t\right) +m_{i}\left(
t\right) }dt\right)
\end{equation*}%
and%
\begin{multline*}
\int_{a}^{b}f\left( t\right) dt \\
=\left( \int_{a}^{b}\left\Vert f\left( t\right) \right\Vert dt-\frac{1}{4n}%
\sum_{i=1}^{n}\left( \int_{a}^{b}\frac{\left[ M_{i}\left( t\right)
-m_{i}\left( t\right) \right] ^{2}}{M_{i}\left( t\right) +m_{i}\left(
t\right) }dt\right) \right) \sum_{i=1}^{n}e_{i}.
\end{multline*}
\end{corollary}

\begin{proof}
As in the proof of Corollary \ref{RTxxc2.2.5}, (\ref{RTxx2.3.17}), implies
that%
\begin{equation*}
\left\Vert f\left( t\right) \right\Vert -\func{Re}\left\langle f\left(
t\right) ,e_{i}\right\rangle \leq \frac{1}{4}\cdot \frac{\left[ M_{i}\left(
t\right) -m_{i}\left( t\right) \right] ^{2}}{M_{i}\left( t\right)
+m_{i}\left( t\right) }\
\end{equation*}%
for a.e. $t\in \left[ a,b\right] $ and $i\in \left\{ 1,\dots ,n\right\} .$

Applying Theorem \ref{RTxxt2.3.1} for%
\begin{equation*}
M_{i}\left( t\right) :=\frac{1}{4}\cdot \frac{\left[ M_{i}\left( t\right)
-m_{i}\left( t\right) \right] ^{2}}{M_{i}\left( t\right) +m_{i}\left(
t\right) },\text{ \ \ }t\in \left[ a,b\right] \text{, \ }i\in \left\{
1,\dots ,n\right\} ,
\end{equation*}%
we deduce the desired result.
\end{proof}

\section{Quadratic Reverses of the Triangle Inequality}

\subsection{Additive Reverses}

The following lemma holds \cite{RTCSSD3}.

\begin{lemma}[Dragomir, 2004]
\label{RTxxl3.2.1}Let $f\in L\left( \left[ a,b\right] ;H\right) $ be such
that there exists a function $k:\Delta \subset \mathbb{R}^{2}\rightarrow
\mathbb{R}$, $\Delta :=\left\{ \left( t,s\right) |a\leq t\leq s\leq
b\right\} $ with the property that $k\in L\left( \Delta \right) $ and%
\begin{equation}
\left( 0\leq \right) \left\Vert f\left( t\right) \right\Vert \left\Vert
f\left( s\right) \right\Vert -\func{Re}\left\langle f\left( t\right)
,f\left( s\right) \right\rangle \leq k\left( t,s\right) ,  \label{RTxx3.2.1}
\end{equation}%
for a.e. $\left( t,s\right) \in \Delta .$ Then we have the following
quadratic reverse of the continuous triangle inequality:%
\begin{equation}
\left( \int_{a}^{b}\left\Vert f\left( t\right) \right\Vert dt\right)
^{2}\leq \left\Vert \int_{a}^{b}f\left( t\right) dt\right\Vert
^{2}+2\iint_{\Delta }k\left( t,s\right) dtds.  \label{RTxx3.2.2}
\end{equation}%
The case of equality holds in (\ref{RTxx3.2.2}) if and only if it holds in (%
\ref{RTxx3.2.1}) for a.e. $\left( t,s\right) \in \Delta .$
\end{lemma}

\begin{proof}
We observe that the following identity holds%
\begin{align}
& \left( \int_{a}^{b}\left\Vert f\left( t\right) \right\Vert dt\right)
^{2}-\left\Vert \int_{a}^{b}f\left( t\right) dt\right\Vert ^{2}
\label{RTxx3.2.3} \\
& =\int_{a}^{b}\int_{a}^{b}\left\Vert f\left( t\right) \right\Vert
\left\Vert f\left( s\right) \right\Vert dtds-\left\langle
\int_{a}^{b}f\left( t\right) dt,\int_{a}^{b}f\left( s\right) ds\right\rangle
\notag \\
& =\int_{a}^{b}\int_{a}^{b}\left\Vert f\left( t\right) \right\Vert
\left\Vert f\left( s\right) \right\Vert dtds-\int_{a}^{b}\int_{a}^{b}\func{Re%
}\left\langle f\left( t\right) ,f\left( s\right) \right\rangle dtds  \notag
\\
& =\int_{a}^{b}\int_{a}^{b}\left[ \left\Vert f\left( t\right) \right\Vert
\left\Vert f\left( s\right) \right\Vert -\func{Re}\left\langle f\left(
t\right) ,f\left( s\right) \right\rangle \right] dtds:=I.  \notag
\end{align}%
Now, observe that for any $\left( t,s\right) \in \left[ a,b\right] \times %
\left[ a,b\right] ,$ we have%
\begin{multline*}
\qquad \left\Vert f\left( t\right) \right\Vert \left\Vert f\left(
s\right) \right\Vert -\func{Re}\left\langle f\left( t\right)
,f\left( s\right)
\right\rangle \\
=\left\Vert f\left( s\right) \right\Vert \left\Vert f\left(
t\right) \right\Vert -\func{Re}\left\langle f\left( s\right)
,f\left( t\right) \right\rangle\qquad
\end{multline*}%
and thus
\begin{equation}
I=2\iint_{\Delta }\left[ \left\Vert f\left( t\right) \right\Vert \left\Vert
f\left( s\right) \right\Vert -\func{Re}\left\langle f\left( t\right)
,f\left( s\right) \right\rangle \right] dtds.  \label{RTxx3.2.4}
\end{equation}%
Using the assumption (\ref{RTxx3.2.1}), we deduce%
\begin{equation*}
\iint_{\Delta }\left[ \left\Vert f\left( t\right) \right\Vert \left\Vert
f\left( s\right) \right\Vert -\func{Re}\left\langle f\left( t\right)
,f\left( s\right) \right\rangle \right] dtds\leq \iint_{\Delta }k\left(
t,s\right) dtds,
\end{equation*}%
and, by the identities (\ref{RTxx3.2.3}) and (\ref{RTxx3.2.4}), we deduce
the desired inequality (\ref{RTxx3.2.2}).

The case of equality is obvious and we omit the details.
\end{proof}

\begin{remark}
\label{RTxxr3.2.2}From (\ref{RTxx3.2.2}) one may deduce a coarser inequality
that can be useful in some applications. It is as follows:%
\begin{equation*}
\left( 0\leq \right) \int_{a}^{b}\left\Vert f\left( t\right) \right\Vert
dt-\left\Vert \int_{a}^{b}f\left( t\right) dt\right\Vert \leq \sqrt{2}\left(
\iint_{\Delta }k\left( t,s\right) dtds\right) ^{\frac{1}{2}}.
\end{equation*}
\end{remark}

\begin{remark}
\label{RTxxr3.2.3}If the condition (\ref{RTxx3.2.1}) is replaced with the
following refinement of the Schwarz inequality%
\begin{equation}
\left( 0\leq \right) k\left( t,s\right) \leq \left\Vert f\left( t\right)
\right\Vert \left\Vert f\left( s\right) \right\Vert -\func{Re}\left\langle
f\left( t\right) ,f\left( s\right) \right\rangle  \label{RTxx3.2.5}
\end{equation}%
for a.e. $\left( t,s\right) \in \Delta ,$ then the following refinement of
the quadratic triangle inequality is valid%
\begin{align}
\left( \int_{a}^{b}\left\Vert f\left( t\right) \right\Vert dt\right) ^{2}&
\geq \left\Vert \int_{a}^{b}f\left( t\right) dt\right\Vert
^{2}+2\iint_{\Delta }k\left( t,s\right) dtds  \label{RTxx3.2.6} \\
& \left( \geq \left\Vert \int_{a}^{b}f\left( t\right) dt\right\Vert
^{2}\right) .  \notag
\end{align}%
The equality holds in (\ref{RTxx3.2.6}) iff the case of equality holds in (%
\ref{RTxx3.2.5}) for a.e. $\left( t,s\right) \in \Delta .$
\end{remark}

The following result holds \cite{RTCSSD3}.

\begin{theorem}[Dragomir, 2004]
\label{RTxxt3.2.4}Let $f\in L\left( \left[ a,b\right] ;H\right) $ be such
that there exists $M\geq 1\geq m\geq 0$ such that either%
\begin{equation}
\func{Re}\left\langle Mf\left( s\right) -f\left( t\right) ,f\left( t\right)
-mf\left( s\right) \right\rangle \geq 0\text{ \ }  \label{RTxx3.2.7}
\end{equation}%
or, equivalently,%
\begin{equation}
\left\Vert f\left( t\right) -\frac{M+m}{2}f\left( s\right) \right\Vert \leq
\frac{1}{2}\left( M-m\right) \left\Vert f\left( s\right) \right\Vert \text{
\ }  \label{RTxx3.2.8}
\end{equation}%
for a.e. $\left( t,s\right) \in \Delta .$ Then we have the inequality:%
\begin{multline}
\left( \int_{a}^{b}\left\Vert f\left( t\right) \right\Vert dt\right)
^{2}\leq \left\Vert \int_{a}^{b}f\left( t\right) dt\right\Vert ^{2}
\label{RTxx3.2.9} \\
+\frac{1}{2}\cdot \frac{\left( M-m\right) ^{2}}{M+m}\int_{a}^{b}\left(
s-a\right) \left\Vert f\left( s\right) \right\Vert ^{2}ds.
\end{multline}%
The case of equality holds in (\ref{RTxx3.2.9}) if and only if%
\begin{equation}
\left\Vert f\left( t\right) \right\Vert \left\Vert f\left( s\right)
\right\Vert -\func{Re}\left\langle f\left( t\right) ,f\left( s\right)
\right\rangle =\frac{1}{4}\cdot \frac{\left( M-m\right) ^{2}}{M+m}\left\Vert
f\left( s\right) \right\Vert ^{2}  \label{RTxx3.2.10}
\end{equation}%
for a.e. $\left( t,s\right) \in \Delta .$
\end{theorem}

\begin{proof}
Taking the square in (\ref{RTxx3.2.8}), we get%
\begin{multline*}
\left\Vert f\left( t\right) \right\Vert ^{2}+\left( \frac{M+m}{2}\right)
^{2}\left\Vert f\left( s\right) \right\Vert ^{2} \\
\leq 2\func{Re}\left\langle f\left( t\right) ,\frac{M+m}{2}f\left( s\right)
\right\rangle +\frac{1}{4}\left( M-m\right) ^{2}\left\Vert f\left( s\right)
\right\Vert ^{2},
\end{multline*}%
for a.e. $\left( t,s\right) \in \Delta ,$ and obviously, since%
\begin{equation*}
2\left( \frac{M+m}{2}\right) \left\Vert f\left( t\right) \right\Vert
\left\Vert f\left( s\right) \right\Vert \leq \left\Vert f\left( t\right)
\right\Vert ^{2}+\left( \frac{M+m}{2}\right) ^{2}\left\Vert f\left( s\right)
\right\Vert ^{2},
\end{equation*}%
we deduce that%
\begin{multline*}
2\left( \frac{M+m}{2}\right) \left\Vert f\left( t\right) \right\Vert
\left\Vert f\left( s\right) \right\Vert \\
\leq 2\func{Re}\left\langle f\left( t\right) ,\frac{M+m}{2}f\left( s\right)
\right\rangle +\frac{1}{4}\left( M-m\right) ^{2}\left\Vert f\left( s\right)
\right\Vert ^{2},
\end{multline*}%
giving the much simpler inequality:%
\begin{equation}
\left\Vert f\left( t\right) \right\Vert \left\Vert f\left( s\right)
\right\Vert -\func{Re}\left\langle f\left( t\right) ,f\left( s\right)
\right\rangle \leq \frac{1}{4}\cdot \frac{\left( M-m\right) ^{2}}{M+m}%
\left\Vert f\left( s\right) \right\Vert ^{2}  \label{RTxx3.2.11}
\end{equation}%
for a.e. $\left( t,s\right) \in \Delta .$

Applying Lemma \ref{RTxxl3.2.1} for $k\left( t,s\right) :=\frac{1}{4}\cdot
\frac{\left( M-m\right) ^{2}}{M+m}\left\Vert f\left( s\right) \right\Vert
^{2},$ we deduce%
\begin{multline}
\left( \int_{a}^{b}\left\Vert f\left( t\right) \right\Vert dt\right)
^{2}\leq \left\Vert \int_{a}^{b}f\left( t\right) dt\right\Vert ^{2}
\label{RTxx3.2.12} \\
+\frac{1}{2}\cdot \frac{\left( M-m\right) ^{2}}{M+m}\iint_{\Delta
}\left\Vert f\left( s\right) \right\Vert ^{2}ds
\end{multline}%
with equality if and only if (\ref{RTxx3.2.11}) holds for a.e. $\left(
t,s\right) \in \Delta .$

Since%
\begin{equation*}
\iint_{\Delta }\left\Vert f\left( s\right) \right\Vert
^{2}ds=\int_{a}^{b}\left( \int_{a}^{s}\left\Vert f\left( s\right)
\right\Vert ^{2}dt\right) ds=\int_{a}^{b}\left( s-a\right) \left\Vert
f\left( s\right) \right\Vert ^{2}ds,
\end{equation*}%
then by (\ref{RTxx3.2.12}) we deduce the desired result (\ref{RTxx3.2.9}).
\end{proof}

Another result which is similar to the one above is incorporated in the
following theorem \cite{RTCSSD3}.

\begin{theorem}[Dragomir, 2004]
\label{RTxxt3.2.5}With the assumptions of Theorem \ref{RTxxt3.2.4}, we have%
\begin{multline}
\left( \int_{a}^{b}\left\Vert f\left( t\right) \right\Vert dt\right)
^{2}-\left\Vert \int_{a}^{b}f\left( t\right) dt\right\Vert ^{2}
\label{RTxx3.2.13} \\
\leq \frac{\left( \sqrt{M}-\sqrt{m}\right) ^{2}}{2\sqrt{Mm}}\left\Vert
\int_{a}^{b}f\left( t\right) dt\right\Vert ^{2}
\end{multline}%
or, equivalently,%
\begin{equation}
\int_{a}^{b}\left\Vert f\left( t\right) \right\Vert dt\leq \left( \frac{M+m}{%
2\sqrt{Mm}}\right) ^{\frac{1}{2}}\left\Vert \int_{a}^{b}f\left( t\right)
dt\right\Vert .  \label{RTxx3.2.14}
\end{equation}%
The case of equality holds in (\ref{RTxx3.2.13}) or (\ref{RTxx3.2.14}) if
and only if
\begin{equation}
\left\Vert f\left( t\right) \right\Vert \left\Vert f\left( s\right)
\right\Vert =\frac{M+m}{2\sqrt{Mm}}\func{Re}\left\langle f\left( t\right)
,f\left( s\right) \right\rangle ,  \label{RTxx3.2.14'}
\end{equation}%
for a.e. $\left( t,s\right) \in \Delta .$
\end{theorem}

\begin{proof}
From (\ref{RTxx3.2.7}), we deduce%
\begin{equation*}
\left\Vert f\left( t\right) \right\Vert ^{2}+Mm\left\Vert f\left( s\right)
\right\Vert ^{2}\leq \left( M+m\right) \func{Re}\left\langle f\left(
t\right) ,f\left( s\right) \right\rangle
\end{equation*}%
for a.e. $\left( t,s\right) \in \Delta .$ Dividing by $\sqrt{Mm}>0,$ we
deduce%
\begin{equation*}
\frac{\left\Vert f\left( t\right) \right\Vert ^{2}}{\sqrt{Mm}}+\sqrt{Mm}%
\left\Vert f\left( s\right) \right\Vert ^{2}\leq \frac{M+m}{\sqrt{Mm}}\func{%
Re}\left\langle f\left( t\right) ,f\left( s\right) \right\rangle
\end{equation*}%
and, obviously, since%
\begin{equation*}
2\left\Vert f\left( t\right) \right\Vert \left\Vert f\left( s\right)
\right\Vert \leq \frac{\left\Vert f\left( t\right) \right\Vert ^{2}}{\sqrt{Mm%
}}+\sqrt{Mm}\left\Vert f\left( s\right) \right\Vert ^{2},
\end{equation*}%
hence%
\begin{equation*}
\left\Vert f\left( t\right) \right\Vert \left\Vert f\left( s\right)
\right\Vert \leq \frac{M+m}{\sqrt{Mm}}\func{Re}\left\langle f\left( t\right)
,f\left( s\right) \right\rangle
\end{equation*}%
for a.e. $\left( t,s\right) \in \Delta ,$ giving%
\begin{equation*}
\left\Vert f\left( t\right) \right\Vert \left\Vert f\left( s\right)
\right\Vert -\func{Re}\left\langle f\left( t\right) ,f\left( s\right)
\right\rangle \leq \frac{\left( \sqrt{M}-\sqrt{m}\right) ^{2}}{2\sqrt{Mm}}%
\func{Re}\left\langle f\left( t\right) ,f\left( s\right) \right\rangle .
\end{equation*}%
Applying Lemma \ref{RTxxl3.2.1} for $k\left( t,s\right) :=\frac{\left( \sqrt{%
M}-\sqrt{m}\right) ^{2}}{\sqrt{Mm}}\func{Re}\left\langle f\left( t\right)
,f\left( s\right) \right\rangle ,$ we deduce%
\begin{multline}
\left( \int_{a}^{b}\left\Vert f\left( t\right) \right\Vert dt\right)
^{2}\leq \left\Vert \int_{a}^{b}f\left( t\right) dt\right\Vert ^{2}
\label{RTxx3.2.15} \\
+\frac{\left( \sqrt{M}-\sqrt{m}\right) ^{2}}{2\sqrt{Mm}}\func{Re}%
\left\langle f\left( t\right) ,f\left( s\right) \right\rangle .
\end{multline}%
On the other hand, since%
\begin{equation*}
\func{Re}\left\langle f\left( t\right) ,f\left( s\right) \right\rangle =%
\func{Re}\left\langle f\left( s\right) ,f\left( t\right) \right\rangle \text{
\ for any \ }\left( t,s\right) \in \left[ a,b\right] ^{2},
\end{equation*}%
hence%
\begin{align*}
\iint_{\Delta }\func{Re}\left\langle f\left( t\right) ,f\left( s\right)
\right\rangle dtds& =\frac{1}{2}\int_{a}^{b}\int_{a}^{b}\func{Re}%
\left\langle f\left( t\right) ,f\left( s\right) \right\rangle dtds \\
& =\frac{1}{2}\func{Re}\left\langle \int_{a}^{b}f\left( t\right)
dt,\int_{a}^{b}f\left( s\right) ds\right\rangle \\
& =\frac{1}{2}\left\Vert \int_{a}^{b}f\left( t\right) dt\right\Vert ^{2}
\end{align*}%
and thus, from (\ref{RTxx3.2.15}), we get (\ref{RTxx3.2.13}).

The equivalence between (\ref{RTxx3.2.13}) and (\ref{RTxx3.2.14}) is obvious
and we omit the details.
\end{proof}

\subsection{Related Results}

The following result also holds \cite{RTCSSD3}.

\begin{theorem}[Dragomir, 2004]
\label{RTxxt3.3.1}Let $f\in L\left( \left[ a,b\right] ;H\right) $ and $%
\gamma ,\Gamma \in \mathbb{R}$ be such that either%
\begin{equation}
\func{Re}\left\langle \Gamma f\left( s\right) -f\left( t\right) ,f\left(
t\right) -\gamma f\left( s\right) \right\rangle \geq 0\text{ \ }
\label{RTxx3.3.1}
\end{equation}%
or, equivalently,%
\begin{equation}
\left\Vert f\left( t\right) -\frac{\Gamma +\gamma }{2}f\left( s\right)
\right\Vert \leq \frac{1}{2}\left\vert \Gamma -\gamma \right\vert \left\Vert
f\left( s\right) \right\Vert \text{ }  \label{RTxx3.3.2}
\end{equation}%
\ for a.e. \ $\left( t,s\right) \in \Delta .$ Then we have the inequality:%
\begin{equation}
\int_{a}^{b}\left[ \left( b-s\right) +\gamma \Gamma \left( s-a\right) \right]
\left\Vert f\left( s\right) \right\Vert ^{2}ds\leq \frac{\Gamma +\gamma }{2}%
\left\Vert \int_{a}^{b}f\left( s\right) ds\right\Vert ^{2}.
\label{RTxx3.3.3}
\end{equation}%
The case of equality holds in (\ref{RTxx3.3.3}) if and only if the case of
equality holds in either (\ref{RTxx3.3.1}) or (\ref{RTxx3.3.2}) for a.e. $%
\left( t,s\right) \in \Delta $.
\end{theorem}

\begin{proof}
The inequality (\ref{RTxx3.3.1}) is obviously equivalent to%
\begin{equation}
\left\Vert f\left( t\right) \right\Vert ^{2}+\gamma \Gamma \left\Vert
f\left( s\right) \right\Vert ^{2}\leq \left( \Gamma +\gamma \right) \func{Re}%
\left\langle f\left( t\right) ,f\left( s\right) \right\rangle
\label{RTxx3.3.4}
\end{equation}%
for a.e. $\left( t,s\right) \in \Delta .$

Integrating (\ref{RTxx3.3.4}) on $\Delta ,$ we deduce%
\begin{multline}
\int_{a}^{b}\left( \int_{a}^{s}\left\Vert f\left( t\right) \right\Vert
^{2}dt\right) ds+\gamma \Gamma \int_{a}^{b}\left( \left\Vert f\left(
s\right) \right\Vert ^{2}\int_{a}^{s}dt\right) ds  \label{RTxx3.3.5} \\
=\left( \Gamma +\gamma \right) \int_{a}^{b}\left( \int_{a}^{s}\func{Re}%
\left\langle f\left( t\right) ,f\left( s\right) \right\rangle dt\right) ds.
\end{multline}%
It is easy to see, on integrating by parts, that%
\begin{align*}
\int_{a}^{b}\left( \int_{a}^{s}\left\Vert f\left( t\right) \right\Vert
^{2}dt\right) ds& =s\left. \int_{a}^{s}\left\Vert f\left( t\right)
\right\Vert ^{2}dt\right\vert _{a}^{b}-\int_{a}^{b}s\left\Vert f\left(
s\right) \right\Vert ^{2}ds \\
& =b\int_{a}^{s}\left\Vert f\left( s\right) \right\Vert
^{2}ds-\int_{a}^{b}s\left\Vert f\left( s\right) \right\Vert ^{2}ds \\
& =\int_{a}^{b}\left( b-s\right) \left\Vert f\left( s\right) \right\Vert
^{2}ds
\end{align*}%
and
\begin{equation*}
\int_{a}^{b}\left( \left\Vert f\left( s\right) \right\Vert
^{2}\int_{a}^{s}dt\right) ds=\int_{a}^{b}\left( s-a\right) \left\Vert
f\left( s\right) \right\Vert ^{2}ds.
\end{equation*}%
Since%
\begin{align*}
\frac{d}{ds}\left( \left\Vert \int_{a}^{b}f\left( t\right) dt\right\Vert
^{2}\right) & =\frac{d}{ds}\left\langle \int_{a}^{s}f\left( t\right)
dt,\int_{a}^{s}f\left( t\right) dt\right\rangle \\
& =\left\langle f\left( s\right) ,\int_{a}^{s}f\left( t\right)
dt\right\rangle +\left\langle \int_{a}^{s}f\left( t\right) dt,f\left(
s\right) \right\rangle \\
& =2\func{Re}\left\langle \int_{a}^{s}f\left( t\right) dt,f\left( s\right)
\right\rangle ,
\end{align*}%
hence%
\begin{align*}
\int_{a}^{b}\left( \int_{a}^{s}\func{Re}\left\langle f\left( t\right)
,f\left( s\right) \right\rangle dt\right) ds& =\int_{a}^{b}\func{Re}%
\left\langle \int_{a}^{s}f\left( t\right) dt,f\left( s\right) \right\rangle
ds \\
& =\frac{1}{2}\int_{a}^{b}\frac{d}{ds}\left( \left\Vert \int_{a}^{s}f\left(
t\right) dt\right\Vert ^{2}\right) ds \\
& =\frac{1}{2}\left\Vert \int_{a}^{b}f\left( t\right) dt\right\Vert ^{2}.
\end{align*}%
Utilising (\ref{RTxx3.3.5}), we deduce the desired inequality (\ref%
{RTxx3.3.3}).

The case of equality is obvious and we omit the details.
\end{proof}

\begin{remark}
Consider the function $\varphi \left( s\right) :=\left( b-s\right) +\gamma
\Gamma \left( s-a\right) ,$ $s\in \left[ a,b\right] .$ Obviously,%
\begin{equation*}
\varphi \left( s\right) =\left( \Gamma \gamma -1\right) s+b-\gamma \Gamma a.
\end{equation*}%
Observe that, if $\Gamma \gamma \geq 1,$ then%
\begin{equation*}
b-a=\varphi \left( a\right) \leq \varphi \left( s\right) \leq \varphi \left(
b\right) =\gamma \Gamma \left( b-a\right) ,\ \ \ \ s\in \left[ a,b\right]
\end{equation*}%
and, if $\Gamma \gamma <1,$ then%
\begin{equation*}
\gamma \Gamma \left( b-a\right) \leq \varphi \left( s\right) \leq b-a,\ \ \
\ s\in \left[ a,b\right] .
\end{equation*}
\end{remark}

Taking into account the above remark, we may state the following corollary
\cite{RTCSSD3}.

\begin{corollary}
\label{RTxxc3.3.2}Assume that $f,\gamma ,\Gamma $ are as in Theorem \ref%
{RTxxt3.3.1}.

\begin{enumerate}
\item[a)] If $\Gamma \gamma \geq 1,$ then we have the inequality%
\begin{equation*}
\left( b-a\right) \int_{a}^{b}\left\Vert f\left( s\right) \right\Vert
^{2}ds\leq \frac{\Gamma +\gamma }{2}\left\Vert \int_{a}^{b}f\left( s\right)
ds\right\Vert ^{2}.
\end{equation*}

\item[b)] If $0<\Gamma \gamma <1,$ then we have the inequality%
\begin{equation*}
\gamma \Gamma \left( b-a\right) \int_{a}^{b}\left\Vert f\left( s\right)
\right\Vert ^{2}ds\leq \frac{\Gamma +\gamma }{2}\left\Vert
\int_{a}^{b}f\left( s\right) ds\right\Vert ^{2}.
\end{equation*}
\end{enumerate}
\end{corollary}

\section{Refinements for Complex Spaces}

\subsection{The Case of a Unit Vector}

The following result holds \cite{RTCSSD4}.

\begin{theorem}[Dragomir, 2004]
\label{RTxxt4.2.1}Let $\left( H;\left\langle \cdot ,\cdot \right\rangle
\right) $ be a complex Hilbert space. If $f\in L\left( \left[ a,b\right]
;H\right) $ is such that there exists $k_{1},k_{2}\geq 0$ with%
\begin{equation}
k_{1}\left\Vert f\left( t\right) \right\Vert \leq \func{Re}\left\langle
f\left( t\right) ,e\right\rangle ,\ \ k_{2}\left\Vert f\left( t\right)
\right\Vert \leq \func{Im}\left\langle f\left( t\right) ,e\right\rangle
\label{RTxx4.2.1}
\end{equation}%
for a.e. $t\in \left[ a,b\right] ,$ where $e\in H,$ $\left\Vert e\right\Vert
=1,$ is given, then%
\begin{equation}
\sqrt{k_{1}^{2}+k_{2}^{2}}\int_{a}^{b}\left\Vert f\left( t\right)
\right\Vert dt\leq \left\Vert \int_{a}^{b}f\left( t\right) dt\right\Vert .
\label{RTxx4.2.2}
\end{equation}%
The case of equality holds in (\ref{RTxx4.2.2}) if and only if%
\begin{equation}
\int_{a}^{b}f\left( t\right) dt=\left( k_{1}+ik_{2}\right) \left(
\int_{a}^{b}\left\Vert f\left( t\right) \right\Vert dt\right) e.
\label{RTxx4.2.3}
\end{equation}
\end{theorem}

\begin{proof}
Using the Schwarz inequality $\left\Vert u\right\Vert \left\Vert
v\right\Vert \geq \left\vert \left\langle u,v\right\rangle \right\vert ,$ $%
u,v\in H;$ in the complex Hilbert space $\left( H;\left\langle \cdot ,\cdot
\right\rangle \right) ,$ we have%
\begin{align}
& \left\Vert \int_{a}^{b}f\left( t\right) dt\right\Vert ^{2}
\label{RTxx4.2.4}
 =\left\Vert \int_{a}^{b}f\left( t\right) dt\right\Vert ^{2}\left\Vert
e\right\Vert ^{2}  \\
& \geq \left\vert \left\langle \int_{a}^{b}f\left( t\right)
dt,e\right\rangle \right\vert ^{2}=\left\vert \int_{a}^{b}\left\langle
f\left( t\right) ,e\right\rangle dt\right\vert ^{2}  \notag \\
& =\left\vert \int_{a}^{b}\func{Re}\left\langle f\left( t\right)
,e\right\rangle dt+i\left( \int_{a}^{b}\func{Im}\left\langle f\left(
t\right) ,e\right\rangle dt\right) \right\vert ^{2}  \notag \\
& =\left( \int_{a}^{b}\func{Re}\left\langle f\left( t\right) ,e\right\rangle
dt\right) ^{2}+\left( \int_{a}^{b}\func{Im}\left\langle f\left( t\right)
,e\right\rangle dt\right) ^{2}.  \notag
\end{align}%
Now, on integrating (\ref{RTxx4.2.1}), we deduce%
\begin{align}
k_{1}\int_{a}^{b}\left\Vert f\left( t\right) \right\Vert dt& \leq
\int_{a}^{b}\func{Re}\left\langle f\left( t\right) ,e\right\rangle dt,
\label{RTxx4.2.5} \\
k_{2}\int_{a}^{b}\left\Vert f\left( t\right) \right\Vert dt& \leq
\int_{a}^{b}\func{Im}\left\langle f\left( t\right) ,e\right\rangle dt  \notag
\end{align}%
implying%
\begin{equation}
\left( \int_{a}^{b}\func{Re}\left\langle f\left( t\right) ,e\right\rangle
dt\right) ^{2}\geq k_{1}^{2}\left( \int_{a}^{b}\left\Vert f\left( t\right)
\right\Vert dt\right) ^{2}  \label{RTxx4.2.6}
\end{equation}%
and%
\begin{equation}
\left( \int_{a}^{b}\func{Im}\left\langle f\left( t\right) ,e\right\rangle
dt\right) ^{2}\geq k_{2}^{2}\left( \int_{a}^{b}\left\Vert f\left( t\right)
\right\Vert dt\right) ^{2}.  \label{RTxx4.2.7}
\end{equation}%
If we add (\ref{RTxx4.2.6}) and (\ref{RTxx4.2.7}) and use (\ref{RTxx4.2.4}),
we deduce the desired inequality (\ref{RTxx4.2.2}).

Further, if (\ref{RTxx4.2.3}) holds, then obviously%
\begin{align*}
\left\Vert \int_{a}^{b}f\left( t\right) dt\right\Vert & =\left\vert
k_{1}+ik_{2}\right\vert \left( \int_{a}^{b}\left\Vert f\left( t\right)
\right\Vert dt\right) \left\Vert e\right\Vert \\
& =\sqrt{k_{1}^{2}+k_{2}^{2}}\int_{a}^{b}\left\Vert f\left( t\right)
\right\Vert dt,
\end{align*}%
and the equality case holds in (\ref{RTxx4.2.2}).

Before we prove the reverse implication, let us observe that, for $x\in H$
and $e\in H,$ $\left\Vert e\right\Vert =1,$ the following identity is valid%
\begin{equation*}
\left\Vert x-\left\langle x,e\right\rangle e\right\Vert ^{2}=\left\Vert
x\right\Vert ^{2}-\left\vert \left\langle x,e\right\rangle \right\vert ^{2},
\end{equation*}%
therefore $\left\Vert x\right\Vert =\left\vert \left\langle x,e\right\rangle
\right\vert $ if and only if $x=\left\langle x,e\right\rangle e.$

If we assume that equality holds in (\ref{RTxx4.2.2}), then the case of
equality must hold in all the inequalities required in the argument used to
prove the inequality (\ref{RTxx4.2.2}). Therefore, we must have%
\begin{equation}
\left\Vert \int_{a}^{b}f\left( t\right) dt\right\Vert =\left\vert
\left\langle \int_{a}^{b}f\left( t\right) dt,e\right\rangle \right\vert
\label{RTxx4.2.7a}
\end{equation}%
and
\begin{equation}
k_{1}\left\Vert f\left( t\right) \right\Vert =\func{Re}\left\langle f\left(
t\right) ,e\right\rangle ,\ \ \ k_{2}\left\Vert f\left( t\right) \right\Vert
=\func{Im}\left\langle f\left( t\right) ,e\right\rangle  \label{RTxx4.2.8}
\end{equation}%
for a.e. $t\in \left[ a,b\right] .$

From (\ref{RTxx4.2.7a}) we deduce%
\begin{equation}
\int_{a}^{b}f\left( t\right) dt=\left\langle \int_{a}^{b}f\left( t\right)
dt,e\right\rangle e,  \label{RTxx4.2.9}
\end{equation}%
and from (\ref{RTxx4.2.8}), by multiplying the second equality with $i,$ the
imaginary unit, and integrating both equations on $\left[ a,b\right] ,$ we
deduce%
\begin{equation}
\left( k_{1}+ik_{2}\right) \int_{a}^{b}\left\Vert f\left( t\right)
\right\Vert dt=\left\langle \int_{a}^{b}f\left( t\right) dt,e\right\rangle .
\label{RTxx4.2.10}
\end{equation}%
Finally, by (\ref{RTxx4.2.9}) and (\ref{RTxx4.2.10}), we deduce the desired
equality (\ref{RTxx4.2.3}).
\end{proof}

The following corollary is of interest \cite{RTCSSD4}.

\begin{corollary}
\label{RTxxc4.2.2}Let $e$ be a unit vector in the complex Hilbert space $%
\left( H;\left\langle \cdot ,\cdot \right\rangle \right) $ and $\eta
_{1},\eta _{2}\in \left( 0,1\right) .$ If $f\in L\left( \left[ a,b\right]
;H\right) $ is such that%
\begin{equation}
\left\Vert f\left( t\right) -e\right\Vert \leq \eta _{1},\left\Vert f\left(
t\right) -ie\right\Vert \leq \eta _{2}\text{ }  \label{RTxx4.2.11a}
\end{equation}%
\ for a.e. \ $t\in \left[ a,b\right] ,$ then we have the inequality%
\begin{equation}
\sqrt{2-\eta _{1}^{2}-\eta _{2}^{2}}\int_{a}^{b}\left\Vert f\left( t\right)
\right\Vert dt\leq \left\Vert \int_{a}^{b}f\left( t\right) dt\right\Vert .
\label{RTxx4.2.12}
\end{equation}%
The case of equality holds in (\ref{RTxx4.2.12}) if and only if%
\begin{equation}
\int_{a}^{b}f\left( t\right) dt=\left( \sqrt{1-\eta _{1}^{2}}+i\sqrt{1-\eta
_{2}^{2}}\right) \left( \int_{a}^{b}\left\Vert f\left( t\right) \right\Vert
dt\right) e.  \label{RTxx4.2.13}
\end{equation}
\end{corollary}

\begin{proof}
From the first inequality in (\ref{RTxx4.2.11a}) we deduce, by taking the
square, that%
\begin{equation*}
\left\Vert f\left( t\right) \right\Vert ^{2}+1-\eta _{1}^{2}\leq 2\func{Re}%
\left\langle f\left( t\right) ,e\right\rangle ,
\end{equation*}%
implying%
\begin{equation}
\frac{\left\Vert f\left( t\right) \right\Vert ^{2}}{\sqrt{1-\eta _{1}^{2}}}+%
\sqrt{1-\eta _{1}^{2}}\leq \frac{2\func{Re}\left\langle f\left( t\right)
,e\right\rangle }{\sqrt{1-\eta _{1}^{2}}}  \label{RTxx4.2.14}
\end{equation}%
for a.e. $t\in \left[ a,b\right] .$

Since, obviously%
\begin{equation}
2\left\Vert f\left( t\right) \right\Vert \leq \frac{\left\Vert f\left(
t\right) \right\Vert ^{2}}{\sqrt{1-\eta _{1}^{2}}}+\sqrt{1-\eta _{1}^{2}},
\label{RTxx4.2.15}
\end{equation}%
hence, by (\ref{RTxx4.2.14}) and (\ref{RTxx4.2.15}) we get%
\begin{equation}
0\leq \sqrt{1-\eta _{1}^{2}}\left\Vert f\left( t\right) \right\Vert \leq
\func{Re}\left\langle f\left( t\right) ,e\right\rangle  \label{RTxx4.2.16}
\end{equation}%
for a.e. $t\in \left[ a,b\right] .$

From the second inequality in (\ref{RTxx4.2.11a}) we deduce%
\begin{equation*}
0\leq \sqrt{1-\eta _{2}^{2}}\left\Vert f\left( t\right) \right\Vert \leq
\func{Re}\left\langle f\left( t\right) ,ie\right\rangle
\end{equation*}%
for a.e. $t\in \left[ a,b\right] .$ Since%
\begin{equation*}
\func{Re}\left\langle f\left( t\right) ,ie\right\rangle =\func{Im}%
\left\langle f\left( t\right) ,e\right\rangle
\end{equation*}%
hence%
\begin{equation}
0\leq \sqrt{1-\eta _{2}^{2}}\left\Vert f\left( t\right) \right\Vert \leq
\func{Im}\left\langle f\left( t\right) ,e\right\rangle  \label{RTxx4.2.17}
\end{equation}%
for a.e. $t\in \left[ a,b\right] .$

Now, observe from (\ref{RTxx4.2.16}) and (\ref{RTxx4.2.17}), that the
condition (\ref{RTxx4.2.1}) of Theorem \ref{RTxxt4.2.1} is satisfied for $%
k_{1}=\sqrt{1-\eta _{1}^{2}},$ $k_{2}=\sqrt{1-\eta _{2}^{2}}\in \left(
0,1\right) ,$ and thus the corollary is proved.
\end{proof}

The following corollary may be stated as well \cite{RTCSSD4}.

\begin{corollary}
\label{RTxxc4.2.3}Let $e$ be a unit vector in the complex Hilbert space $%
\left( H;\left\langle \cdot ,\cdot \right\rangle \right) $ and $M_{1}\geq
m_{1}>0,$ $M_{2}\geq m_{2}>0.$ If $f\in L\left( \left[ a,b\right] ;H\right) $
is such that either%
\begin{align}
\func{Re}\left\langle M_{1}e-f\left( t\right) ,f\left( t\right)
-m_{1}e\right\rangle & \geq 0,\text{ }  \label{RTxx4.2.18} \\
\func{Re}\left\langle M_{2}ie-f\left( t\right) ,f\left( t\right)
-m_{2}ie\right\rangle & \geq 0  \notag
\end{align}%
or, equivalently,%
\begin{align}
\left\Vert f\left( t\right) -\frac{M_{1}+m_{1}}{2}e\right\Vert & \leq \frac{1%
}{2}\left( M_{1}-m_{1}\right) ,  \label{RTxx4.2.19} \\
\left\Vert f\left( t\right) -\frac{M_{2}+m_{2}}{2}ie\right\Vert & \leq \frac{%
1}{2}\left( M_{2}-m_{2}\right) ,  \notag
\end{align}%
for a.e. $t\in \left[ a,b\right] ,$ then we have the inequality%
\begin{multline}
2\left[ \frac{m_{1}M_{1}}{\left( M_{1}+m_{1}\right) ^{2}}+\frac{m_{2}M_{2}}{%
\left( M_{2}+m_{2}\right) ^{2}}\right] ^{\frac{1}{2}}\int_{a}^{b}\left\Vert
f\left( t\right) \right\Vert dt  \label{RTxx4.2.20} \\
\leq \left\Vert \int_{a}^{b}f\left( t\right) dt\right\Vert .\qquad
\end{multline}%
The equality holds in (\ref{RTxx4.2.20}) if and only if%
\begin{equation}
\int_{a}^{b}f\left( t\right) dt=2\left( \frac{\sqrt{m_{1}M_{1}}}{M_{1}+m_{1}}%
+i\frac{\sqrt{m_{2}M_{2}}}{M_{2}+m_{2}}\right) \left( \int_{a}^{b}\left\Vert
f\left( t\right) \right\Vert dt\right) e.  \label{RTxx4.2.21}
\end{equation}
\end{corollary}

\begin{proof}
From the first inequality in (\ref{RTxx4.2.18}), we get%
\begin{equation*}
\left\Vert f\left( t\right) \right\Vert ^{2}+m_{1}M_{1}\leq \left(
M_{1}+m_{1}\right) \func{Re}\left\langle f\left( t\right) ,e\right\rangle
\end{equation*}%
implying%
\begin{equation}
\frac{\left\Vert f\left( t\right) \right\Vert ^{2}}{\sqrt{m_{1}M_{1}}}+\sqrt{%
m_{1}M_{1}}\leq \frac{M_{1}+m_{1}}{\sqrt{m_{1}M_{1}}}\func{Re}\left\langle
f\left( t\right) ,e\right\rangle  \label{RTxx4.2.22}
\end{equation}%
for a.e. $t\in \left[ a,b\right] .$

Since, obviously,%
\begin{equation}
2\left\Vert f\left( t\right) \right\Vert \leq \frac{\left\Vert f\left(
t\right) \right\Vert ^{2}}{\sqrt{m_{1}M_{1}}}+\sqrt{m_{1}M_{1}},
\label{RTxx4.2.23}
\end{equation}%
hence, by (\ref{RTxx4.2.22}) and (\ref{RTxx4.2.23})%
\begin{equation}
0\leq \frac{2\sqrt{m_{1}M_{1}}}{M_{1}+m_{1}}\left\Vert f\left( t\right)
\right\Vert \leq \func{Re}\left\langle f\left( t\right) ,e\right\rangle
\label{RTxx4.2.24}
\end{equation}%
for a.e. $t\in \left[ a,b\right] .$

Using the same argument as in the proof of Corollary \ref{RTxxc4.2.2}, we
deduce the desired inequality. We omit the details.
\end{proof}

\subsection{The Case of Orthonormal Vectors}

The following result holds \cite{RTCSSD4}.

\begin{theorem}[Dragomir, 2004]
\label{RTxxt4.3.2}Let $\left\{ e_{1},\dots ,e_{n}\right\} $ be a family of
orthonormal vectors in the complex Hilbert space $\left( H;\left\langle
\cdot ,\cdot \right\rangle \right) $. If $k_{j},h_{j}\geq 0,$ $j\in \left\{
1,\dots ,n\right\} $ and $f\in L\left( \left[ a,b\right] ;H\right) $ are
such that%
\begin{equation}
k_{j}\left\Vert f\left( t\right) \right\Vert \leq \func{Re}\left\langle
f\left( t\right) ,e_{j}\right\rangle ,\quad h_{j}\left\Vert f\left( t\right)
\right\Vert \leq \func{Im}\left\langle f\left( t\right) ,e_{j}\right\rangle
\label{RTxx4.3.4}
\end{equation}%
for each $j\in \left\{ 1,\dots ,n\right\} $ and a.e. $t\in \left[ a,b\right]
,$ then%
\begin{equation}
\left[ \sum_{j=1}^{n}\left( k_{j}^{2}+h_{j}^{2}\right) \right] ^{\frac{1}{2}%
}\int_{a}^{b}\left\Vert f\left( t\right) \right\Vert dt\leq \left\Vert
\int_{a}^{b}f\left( t\right) dt\right\Vert .  \label{RTxx4.3.5}
\end{equation}%
The case of equality holds in (\ref{RTxx4.3.5}) if and only if%
\begin{equation}
\int_{a}^{b}f\left( t\right) dt=\left( \int_{a}^{b}\left\Vert f\left(
t\right) \right\Vert dt\right) \sum_{j=1}^{n}\left( k_{j}+ih_{j}\right)
e_{j}.  \label{RTxx4.3.6}
\end{equation}
\end{theorem}

\begin{proof}
Before we prove the theorem, let us recall that, if $x\in H$ and $%
e_{1},\dots ,e_{n}$ are orthonormal vectors, then the following identity
holds true:%
\begin{equation}
\left\Vert x-\sum_{j=1}^{n}\left\langle x,e_{j}\right\rangle
e_{j}\right\Vert ^{2}=\left\Vert x\right\Vert ^{2}-\sum_{j=1}^{n}\left\vert
\left\langle x,e_{j}\right\rangle \right\vert ^{2}.  \label{RTxx4.3.7}
\end{equation}%
As a consequence of this identity, we have the \textit{Bessel inequality}%
\begin{equation}
\sum_{j=1}^{n}\left\vert \left\langle x,e_{j}\right\rangle \right\vert
^{2}\leq \left\Vert x\right\Vert ^{2},x\in H,  \label{RTxx4.3.8}
\end{equation}%
in which, the case of equality holds if and only if
\begin{equation}
x=\sum_{j=1}^{n}\left\langle x,e_{j}\right\rangle e_{j}.  \label{RTxx4.3.9}
\end{equation}%
Now, applying Bessel's inequality for $x=\int_{a}^{b}f\left( t\right) dt,$
we have successively%
\begin{align}
& \left\Vert \int_{a}^{b}f\left( t\right) dt\right\Vert ^{2}
\label{RTxx4.3.10} \\
& \geq \sum_{j=1}^{n}\left\vert \left\langle \int_{a}^{b}f\left( t\right)
dt,e_{j}\right\rangle \right\vert ^{2}=\sum_{j=1}^{n}\left\vert
\int_{a}^{b}\left\langle f\left( t\right) ,e_{j}\right\rangle dt\right\vert
^{2}  \notag \displaybreak\\
& =\sum_{j=1}^{n}\left\vert \int_{a}^{b}\func{Re}\left\langle f\left(
t\right) ,e_{j}\right\rangle dt+i\left( \int_{a}^{b}\func{Im}\left\langle
f\left( t\right) ,e_{j}\right\rangle dt\right) \right\vert ^{2}  \notag \\
& =\sum_{j=1}^{n}\left[ \left( \int_{a}^{b}\func{Re}\left\langle f\left(
t\right) ,e_{j}\right\rangle dt\right) ^{2}+\left( \int_{a}^{b}\func{Im}%
\left\langle f\left( t\right) ,e_{j}\right\rangle dt\right) ^{2}\right] .
\notag
\end{align}%
Integrating (\ref{RTxx4.3.4}) on $\left[ a,b\right] ,$ we get%
\begin{equation}
\int_{a}^{b}\func{Re}\left\langle f\left( t\right) ,e_{j}\right\rangle
dt\geq k_{j}\int_{a}^{b}\left\Vert f\left( t\right) \right\Vert dt
\label{RTxx4.3.11}
\end{equation}%
and%
\begin{equation}
\int_{a}^{b}\func{Im}\left\langle f\left( t\right) ,e_{j}\right\rangle
dt\geq h_{j}\int_{a}^{b}\left\Vert f\left( t\right) \right\Vert dt
\label{RTxx4.3.12}
\end{equation}%
for each $j\in \left\{ 1,\dots ,n\right\} .$

Squaring and adding the above two inequalities (\ref{RTxx4.3.11}) and (\ref%
{RTxx4.3.12}), we deduce%
\begin{multline*}
\sum_{j=1}^{n}\left[ \left( \int_{a}^{b}\func{Re}\left\langle f\left(
t\right) ,e_{j}\right\rangle dt\right) ^{2}+\left( \int_{a}^{b}\func{Im}%
\left\langle f\left( t\right) ,e_{j}\right\rangle dt\right) ^{2}\right] \\
\geq \sum_{j=1}^{n}\left( k_{j}^{2}+h_{j}^{2}\right) \left(
\int_{a}^{b}\left\Vert f\left( t\right) \right\Vert dt\right) ^{2},
\end{multline*}%
which combined with (\ref{RTxx4.3.10}) will produce the desired inequality (%
\ref{RTxx4.3.5}).

Now, if (\ref{RTxx4.3.6}) holds true, then%
\begin{align*}
\left\Vert \int_{a}^{b}f\left( t\right) dt\right\Vert & =\left(
\int_{a}^{b}\left\Vert f\left( t\right) \right\Vert dt\right) \left\Vert
\sum_{j=1}^{n}\left( k_{j}+ih_{j}\right) e_{j}\right\Vert \\
& =\left( \int_{a}^{b}\left\Vert f\left( t\right) \right\Vert dt\right)
\left( \left\Vert \sum_{j=1}^{n}\left( k_{j}+ih_{j}\right) e_{j}\right\Vert
^{2}\right) ^{\frac{1}{2}} \\
& =\left( \int_{a}^{b}\left\Vert f\left( t\right) \right\Vert dt\right)
\left[ \sum_{j=1}^{n}\left( k_{j}^{2}+h_{j}^{2}\right) \right] ^{\frac{1}{2}%
},
\end{align*}%
and the case of equality holds in (\ref{RTxx4.3.5}).

Conversely, if the equality holds in (\ref{RTxx4.3.5}), then it must hold in
all the inequalities used to prove (\ref{RTxx4.3.5}) and therefore we must
have%
\begin{equation}
\left\Vert \int_{a}^{b}f\left( t\right) dt\right\Vert
^{2}=\sum_{j=1}^{n}\left\vert \left\langle \int_{a}^{b}f\left( t\right)
dt,e_{j}\right\rangle \right\vert ^{2}  \label{RTxx4.3.13}
\end{equation}%
and%
\begin{equation}
k_{j}\left\Vert f\left( t\right) \right\Vert =\func{Re}\left\langle f\left(
t\right) ,e_{j}\right\rangle \text{ \ and \ }h_{j}\left\Vert f\left(
t\right) \right\Vert =\func{Re}\left\langle f\left( t\right)
,e_{j}\right\rangle  \label{RTxx4.3.14}
\end{equation}%
for each $j\in \left\{ 1,\dots ,n\right\} $ and a.e. $t\in \left[ a,b\right]
.$

From (\ref{RTxx4.3.13}), on using the identity (\ref{RTxx4.3.9}), we deduce
that%
\begin{equation}
\int_{a}^{b}f\left( t\right) dt=\sum_{j=1}^{n}\left\langle
\int_{a}^{b}f\left( t\right) dt,e_{j}\right\rangle e_{j}.  \label{RTxx4.3.15}
\end{equation}%
Now, multiplying the second equality in (\ref{RTxx4.3.14}) with the
imaginary unit $i,$ integrating both inequalities on $\left[ a,b\right] $
and summing them up, we get%
\begin{equation}
\left( k_{j}+ih_{j}\right) \int_{a}^{b}\left\Vert f\left( t\right)
\right\Vert dt=\left\langle \int_{a}^{b}f\left( t\right)
dt,e_{j}\right\rangle  \label{RTxx4.3.16}
\end{equation}%
for each $j\in \left\{ 1,\dots ,n\right\} .$

Finally, utilising (\ref{RTxx4.3.15}) and (\ref{RTxx4.3.16}), we deduce (\ref%
{RTxx4.3.6}) and the theorem is proved.
\end{proof}

The following corollaries are of interest \cite{RTCSSD4}.

\begin{corollary}
\label{RTxxc4.3.2}Let $e_{1},\dots ,e_{m}$ be orthonormal vectors in the
complex Hilbert space $\left( H;\left\langle \cdot ,\cdot \right\rangle
\right) $ and $\rho _{k},\eta _{k}\in \left( 0,1\right) ,$ $k\in \left\{
1,\dots ,n\right\} .$ If $f\in L\left( \left[ a,b\right] ;H\right) $ is such
that%
\begin{equation*}
\left\Vert f\left( t\right) -e_{k}\right\Vert \leq \rho _{k},\qquad
\left\Vert f\left( t\right) -ie_{k}\right\Vert \leq \eta _{k}
\end{equation*}%
for each $k\in \left\{ 1,\dots ,n\right\} $ and for a.e. $t\in \left[ a,b%
\right] ,$ then we have the inequality%
\begin{equation}
\left[ \sum_{k=1}^{n}\left( 2-\rho _{k}^{2}-\eta _{k}^{2}\right) \right] ^{%
\frac{1}{2}}\int_{a}^{b}\left\Vert f\left( t\right) \right\Vert dt\leq
\left\Vert \int_{a}^{b}f\left( t\right) dt\right\Vert .  \label{RTxx4.3.17}
\end{equation}%
The case of equality holds in (\ref{RTxx4.3.17}) if and only if%
\begin{multline}
\quad \int_{a}^{b}f\left( t\right) dt\\=\left(
\int_{a}^{b}\left\Vert f\left(
t\right) \right\Vert dt\right) \sum_{k=1}^{n}\left( \sqrt{1-\rho _{k}^{2}}+i%
\sqrt{1-\eta _{k}^{2}}\right) e_{k}.  \label{RTxx4.3.18}\quad
\end{multline}
\end{corollary}

The proof follows by Theorem \ref{RTxxt4.3.2} and is similar to the one from
Corollary \ref{RTxxc4.2.2}. We omit the details.

Next, the following result may be stated \cite{RTCSSD4}:

\begin{corollary}
\label{RTxxc4.3.3}Let $e_{1},\dots ,e_{m}$ be as in Corollary \ref%
{RTxxc4.3.2} and $M_{k}\geq m_{k}>0,$ $N_{k}\geq n_{k}>0,$ $k\in \left\{
1,\dots ,n\right\} .$ If $f\in L\left( \left[ a,b\right] ;H\right) $ is such
that either%
\begin{align*}
\func{Re}\left\langle M_{k}e_{k}-f\left( t\right) ,f\left( t\right)
-m_{k}e_{k}\right\rangle & \geq 0,\  \\
\func{Re}\left\langle N_{k}ie_{k}-f\left( t\right) ,f\left( t\right)
-n_{k}ie_{k}\right\rangle & \geq 0
\end{align*}%
or, equivalently,%
\begin{align*}
\left\Vert f\left( t\right) -\frac{M_{k}+m_{k}}{2}e_{k}\right\Vert & \leq
\frac{1}{2}\left( M_{k}-m_{k}\right) ,\  \\
\left\Vert f\left( t\right) -\frac{N_{k}+n_{k}}{2}ie_{k}\right\Vert & \leq
\frac{1}{2}\left( N_{k}-n_{k}\right)
\end{align*}%
for each $k\in \left\{ 1,\dots ,n\right\} $ and a.e. $t\in \left[ a,b\right]
,$ then we have the inequality%
\begin{multline}
2\left\{ \sum_{k=1}^{m}\left[ \frac{m_{k}M_{k}}{\left( M_{k}+m_{k}\right)
^{2}}+\frac{n_{k}N_{k}}{\left( N_{k}+n_{k}\right) ^{2}}\right] \right\} ^{%
\frac{1}{2}}\int_{a}^{b}\left\Vert f\left( t\right) \right\Vert dt
\label{RTxx4.3.19} \\
\leq \left\Vert \int_{a}^{b}f\left( t\right) dt\right\Vert .
\end{multline}%
The case of equality holds in (\ref{RTxx4.3.19}) if and only if%
\begin{multline}
\quad \int_{a}^{b}f\left( t\right) dt=2\left(
\int_{a}^{b}\left\Vert f\left(
t\right) \right\Vert dt\right)   \label{RTxx4.3.20} \\
\times \sum_{k=1}^{n}\left( \frac{\sqrt{m_{k}M_{k}}}{M_{k}+m_{k}}+i\frac{%
\sqrt{n_{k}N_{k}}}{N_{k}+n_{k}}\right) e_{k}.\quad
\end{multline}
\end{corollary}

The proof employs Theorem \ref{RTxxt4.3.2} and is similar to the one in
Corollary \ref{RTxxc4.2.3}. We omit the details.

\section{Applications for Complex-Valued Functions}

The following proposition holds \cite{RTCSSD1}.

\begin{proposition}
\label{RTxxp1.4.1}If $f:\left[ a,b\right] \rightarrow \mathbb{C}$ is a
Lebesgue integrable function with the property that there exists a constant $%
K\geq 1$ such that%
\begin{equation}
\left\vert f\left( t\right) \right\vert \leq K\left[ \alpha \func{Re}f\left(
t\right) +\beta \func{Im}f\left( t\right) \right]  \label{RTxx1.4.1}
\end{equation}%
for a.e. $t\in \left[ a,b\right] ,$ where $\alpha ,\beta \in \mathbb{R}$, $%
\alpha ^{2}+\beta ^{2}=1$ are given, then we have the following reverse of
the continuous triangle inequality:%
\begin{equation}
\int_{a}^{b}\left\vert f\left( t\right) \right\vert dt\leq K\left\vert
\int_{a}^{b}f\left( t\right) dt\right\vert .  \label{RTxx1.4.2}
\end{equation}%
The case of equality holds in (\ref{RTxx1.4.2}) if and only if%
\begin{equation*}
\int_{a}^{b}f\left( t\right) dt=\frac{1}{K}\left( \alpha +i\beta \right)
\int_{a}^{b}\left\vert f\left( t\right) \right\vert dt.
\end{equation*}
\end{proposition}

The proof is obvious by Theorem \ref{RTxxt1.2.1}, and we omit the details.

\begin{remark}
\label{RTxxr1.4.2}If in the above Proposition \ref{RTxxp1.4.1} we choose $%
\alpha =1,$ $\beta =0,$ then the condition (\ref{RTxx1.4.1}) for $\func{Re}%
f\left( t\right) >0$ is equivalent to%
\begin{equation*}
\left[ \func{Re}f\left( t\right) \right] ^{2}+\left[ \func{Im}f\left(
t\right) \right] ^{2}\leq K^{2}\left[ \func{Re}f\left( t\right) \right] ^{2}
\end{equation*}%
or with the inequality:%
\begin{equation*}
\frac{\left\vert \func{Im}f\left( t\right) \right\vert }{\func{Re}f\left(
t\right) }\leq \sqrt{K^{2}-1}.
\end{equation*}%
Now, if we assume that%
\begin{equation}
\left\vert \arg f\left( t\right) \right\vert \leq \theta ,\ \ \ \theta \in
\left( 0,\frac{\pi }{2}\right) ,  \label{RTxx1.4.3}
\end{equation}%
then, for $\func{Re}f\left( t\right) >0,$%
\begin{equation*}
\left\vert \tan \left[ \arg f\left( t\right) \right] \right\vert =\frac{%
\left\vert \func{Im}f\left( t\right) \right\vert }{\func{Re}f\left( t\right)
}\leq \tan \theta ,
\end{equation*}%
and if we choose $K=\frac{1}{\cos \theta }>1,$ then%
\begin{equation*}
\sqrt{K^{2}-1}=\tan \theta ,
\end{equation*}%
and by Proposition \ref{RTxxp1.4.1}, we deduce%
\begin{equation}
\cos \theta \int_{a}^{b}\left\vert f\left( t\right) \right\vert dt\leq
\left\vert \int_{a}^{b}f\left( t\right) dt\right\vert ,  \label{RTxx1.4.4}
\end{equation}%
which is exactly the Karamata inequality (\ref{RTxx0.1.2}) from the
Introduction.
\end{remark}

Obviously, the result from Proposition \ref{RTxxp1.4.1} is more
comprehensive since for other values of $\left( \alpha ,\beta \right) \in
\mathbb{R}^{2}$ with $\alpha ^{2}+\beta ^{2}=1$ we can get different
sufficient conditions for the function $f$ such that the inequality (\ref%
{RTxx1.4.2}) holds true.

A different sufficient condition in terms of complex disks is incorporated
in the following proposition \cite{RTCSSD1}.

\begin{proposition}
\label{RTxxp1.4.3}Let $e=\alpha +i\beta $ with $\alpha ^{2}+\beta ^{2}=1,$ $%
r\in \left( 0,1\right) $ and $f:\left[ a,b\right] \rightarrow \mathbb{C}$ a
Lebesgue integrable function such that%
\begin{equation}
f\left( t\right) \in \bar{D}\left( e,r\right) :=\left\{ z\in \mathbb{C}|%
\text{ }\left\vert z-e\right\vert \leq r\right\} \text{ \ \ for a.e. }t\in %
\left[ a,b\right] .  \label{RTxx1.4.5}
\end{equation}%
Then we have the inequality%
\begin{equation}
\sqrt{1-r^{2}}\int_{a}^{b}\left\vert f\left( t\right) \right\vert dt\leq
\left\vert \int_{a}^{b}f\left( t\right) dt\right\vert .  \label{RTxx1.4.6}
\end{equation}%
The case of equality holds in (\ref{RTxx1.4.6}) if and only if%
\begin{equation*}
\int_{a}^{b}f\left( t\right) dt=\sqrt{1-r^{2}}\left( \alpha +i\beta \right)
\int_{a}^{b}\left\vert f\left( t\right) \right\vert dt.
\end{equation*}
\end{proposition}

The proof follows by Corollary \ref{RTxxc1.2.2} and we omit the details.

Further, we may state the following proposition as well \cite{RTCSSD1}.

\begin{proposition}
\label{RTxxp1.4.4}Let $e=\alpha +i\beta $ with $\alpha ^{2}+\beta ^{2}=1$
and $M\geq m>0.$ If $f:\left[ a,b\right] \rightarrow \mathbb{C}$ is such that%
\begin{equation}
\func{Re}\left[ \left( Me-f\left( t\right) \right) \left( \overline{f\left(
t\right) }-m\overline{e}\right) \right] \geq 0\text{ \ \ for a.e. }t\in %
\left[ a,b\right] ,  \label{RTxx1.4.7}
\end{equation}%
or, equivalently,%
\begin{equation}
\left\vert f\left( t\right) -\frac{M+m}{2}e\right\vert \leq \frac{1}{2}%
\left( M-m\right) \text{ \ \ for a.e. }t\in \left[ a,b\right] ,
\label{RTxx1.4.8}
\end{equation}%
then we have the inequality%
\begin{equation}
\frac{2\sqrt{mM}}{M+m}\int_{a}^{b}\left\vert f\left( t\right) \right\vert
dt\leq \left\vert \int_{a}^{b}f\left( t\right) dt\right\vert ,
\label{RTxx1.4.9}
\end{equation}%
or, equivalently,%
\begin{align}
(0& \leq )\int_{a}^{b}\left\vert f\left( t\right) \right\vert dt-\left\vert
\int_{a}^{b}f\left( t\right) dt\right\vert  \label{RTxx1.4.10} \\
& \leq \frac{\left( \sqrt{M}-\sqrt{m}\right) ^{2}}{M+m}\left\vert
\int_{a}^{b}f\left( t\right) dt\right\vert .  \notag
\end{align}%
The equality holds in (\ref{RTxx1.4.9}) (or in the second part of (\ref%
{RTxx1.4.10})) if and only if%
\begin{equation*}
\int_{a}^{b}f\left( t\right) dt=\frac{2\sqrt{mM}}{M+m}\left( \alpha +i\beta
\right) \int_{a}^{b}\left\vert f\left( t\right) \right\vert dt.
\end{equation*}
\end{proposition}

The proof follows by Corollary \ref{RTxxc1.2.3} and we omit the details.

\begin{remark}
\label{RTxxr1.4.5}Since%
\begin{align*}
Me-f\left( t\right) & =M\alpha -\func{Re}f\left( t\right) +i\left[ M\beta -%
\func{Im}f\left( t\right) \right] , \\
\overline{f\left( t\right) }-m\overline{e}& =\func{Re}f\left( t\right)
-m\alpha -i\left[ \func{Im}f\left( t\right) -m\beta \right]
\end{align*}%
hence%
\begin{multline}
\func{Re}\left[ \left( Me-f\left( t\right) \right) \left( \overline{f\left(
t\right) }-m\overline{e}\right) \right]  \label{RTxx1.4.11} \\
=\left[ M\alpha -\func{Re}f\left( t\right) \right] \left[ \func{Re}f\left(
t\right) -m\alpha \right] \\
+\left[ M\beta -\func{Im}f\left( t\right) \right] \left[ \func{Im}f\left(
t\right) -m\beta \right] .
\end{multline}%
It is obvious that, if%
\begin{equation}
m\alpha \leq \func{Re}f\left( t\right) \leq M\alpha \text{ \ \ for a.e. }%
t\in \left[ a,b\right] ,  \label{RTxx1.4.12}
\end{equation}%
and
\begin{equation}
m\beta \leq \func{Im}f\left( t\right) \leq M\beta \text{ \ \ for a.e. }t\in %
\left[ a,b\right] ,  \label{RTxx1.4.13}
\end{equation}%
then, by (\ref{RTxx1.4.11}),%
\begin{equation*}
\func{Re}\left[ \left( Me-f\left( t\right) \right) \left( \overline{f\left(
t\right) }-m\overline{e}\right) \right] \geq 0\text{ \ \ for a.e. }t\in %
\left[ a,b\right] ,
\end{equation*}%
and then either (\ref{RTxx1.4.9}) or (\ref{RTxx1.4.12}) hold true.
\end{remark}

We observe that the conditions (\ref{RTxx1.4.12}) and (\ref{RTxx1.4.13}) are
very easy to verify in practice and may be useful in various applications
where reverses of the continuous triangle inequality are required.

\begin{remark}
Similar results may be stated for functions $f:\left[ a,b\right] \rightarrow
\mathbb{R}^{n}$ or $f:\left[ a,b\right] \rightarrow H,$ with $H$ particular
instances of Hilbert spaces of significance in applications, but we leave
them to the interested reader.
\end{remark}

Let $e=\alpha +i\beta $ $\left( \alpha ,\beta \in \mathbb{R}\right) $ be a
complex number with the property that $\left\vert e\right\vert =1,$ i.e., $%
\alpha ^{2}+\beta ^{2}=1.$ The following proposition concerning a reverse of
the continuous triangle inequality for complex-valued functions may be
stated \cite{RTCSSD2}:

\begin{proposition}
\label{RTxxp.2.4.1} Let $f:\left[ a,b\right] \rightarrow \mathbb{C}$ be a
Lebesgue integrable function with the property that there exists a constant $%
\rho \in \left( 0,1\right) $ such that
\begin{equation}
\left\vert f\left( t\right) -e\right\vert \leq \rho \text{ for a.e. }t\in %
\left[ a,b\right] ,  \label{RTxxe.2.4.1}
\end{equation}%
where $e$ has been defined above. Then we have the following reverse of the
continuous triangle inequality
\begin{align}
(0& \leq )\int_{a}^{b}\left\vert f\left( t\right) \right\vert dt-\left\vert
\int_{a}^{b}f\left( t\right) dt\right\vert  \label{RTxxe.2.4.2} \\
& \leq \frac{\rho ^{2}}{\sqrt{1-\rho ^{2}}\left( 1+\sqrt{1-\rho ^{2}}\right)
}  \notag \\
& \qquad \times \left[ \alpha \int_{a}^{b}\func{Re}f\left( t\right) dt+\beta
\int_{a}^{b}\func{Im}f\left( t\right) dt\right] .  \notag
\end{align}
\end{proposition}

The proof follows by Corollary \ref{RTxxc2.2.2}, and the details are omitted.

On the other hand, the following result is perhaps more useful for
applications \cite{RTCSSD2}:

\begin{proposition}
\label{RTxxp.2.4.2} Assume that $f$ and $e$ are as in Proposition \ref%
{RTxxp.2.4.1}. If there exists the constants $M\geq m>0$ such that either
\begin{equation}
\func{Re}\left[ \left( Me-f\left( t\right) \right) \left( \overline{f\left(
t\right) }-m\overline{e}\right) \right] \geq 0  \label{RTxxe.2.4.3}
\end{equation}%
or, equivalently,
\begin{equation}
\left\vert f\left( t\right) -\frac{M+m}{2}e\right\vert \leq \frac{1}{2}%
\left( M-m\right)  \label{RTxxe.2.4.4}
\end{equation}%
for a.e. $t\in \left[ a,b\right] ,$ holds, then
\begin{align}
(0& \leq )\int_{a}^{b}\left\vert f\left( t\right) \right\vert dt-\left\vert
\int_{a}^{b}f\left( t\right) dt\right\vert  \label{RTxxe.2.4.5} \\
& \leq \frac{\left( \sqrt{M}-\sqrt{m}\right) ^{2}}{2\sqrt{Mm}}\left[ \alpha
\int_{a}^{b}\func{Re}f\left( t\right) dt+\beta \int_{a}^{b}\func{Im}f\left(
t\right) dt\right] .  \notag
\end{align}
\end{proposition}

The proof may be done on utilising Corollary \ref{RTxxc2.2.3}, but we omit
the details

Subsequently, on making use of Corollary \ref{RTxxc2.2.5}, one may state the
following result as well \cite{RTCSSD2}:

\begin{proposition}
\label{RTxxp.2.4.3} Let $f$ be as in Proposition \ref{RTxxp.2.4.1} and the
measurable functions $K,k:\left[ a,b\right] \rightarrow \lbrack 0,\infty )$
with the property that
\begin{equation*}
\frac{\left( K-k\right) ^{2}}{K+k}\in L\left[ a,b\right]
\end{equation*}%
and
\begin{equation*}
\alpha k\left( t\right) \leq \func{Re}f\left( t\right) \leq \alpha K\left(
t\right) \text{ and }\beta k\left( t\right) \leq \func{Im}f\left( t\right)
\leq \beta K\left( t\right)
\end{equation*}%
for a.e. $t\in \left[ a,b\right] ,$ where $\alpha ,\beta $ are assumed to be
positive and satisfying the condition $\alpha ^{2}+\beta ^{2}=1$. Then the
following reverse of the continuous triangle inequality is valid:
\begin{align*}
(0& \leq )\int_{a}^{b}\left\vert f\left( t\right) \right\vert dt-\left\vert
\int_{a}^{b}f\left( t\right) dt\right\vert \\
& \leq \frac{1}{4}\int_{a}^{b}\frac{\left[ K\left( t\right) -k\left(
t\right) \right] ^{2}}{K\left( t\right) +k\left( t\right) }dt.
\end{align*}%
The constant $\frac{1}{4}$ is best possible in the sense that it cannot be
replaced by a smaller quantity.
\end{proposition}

\begin{remark}
\label{RTxxr.2.4.2} One may realise that similar results can be stated if
the Corollaries \ref{RTxxc2.3.2}-\ref{RTxxc2.3.5}\ obtained above are used.
For the sake of brevity, we do not mention them here.
\end{remark}

Let $f:\left[ a,b\right] \rightarrow \mathbb{C}$ be a Lebesgue integrable
function and $M\geq 1\geq m\geq 0.$ The condition (\ref{RTxx3.2.7}) from
Theorem \ref{RTxxt3.2.4}, which plays a fundamental role in the results
obtained above, can be translated in this case as%
\begin{equation}
\func{Re}\left[ \left( Mf\left( s\right) -f\left( t\right) \right) \left(
\overline{f\left( t\right) }-m\overline{f\left( s\right) }\right) \right]
\geq 0  \label{RTxxe.3.4.1}
\end{equation}%
for a.e. $a\leq t\leq s\leq b.$

Since, obviously%
\begin{multline*}
\func{Re}\left[ \left( Mf\left( s\right) -f\left( t\right) \right) \left(
\overline{f\left( t\right) }-m\overline{f\left( s\right) }\right) \right] \\
=\left[ \left( M\func{Re}f\left( s\right) -\func{Re}f\left( t\right) \right)
\left( \func{Re}f\left( t\right) -m\func{Re}f\left( s\right) \right) \right]
\\
+\left[ \left( M\func{Im}f\left( s\right) -\func{Im}f\left( t\right) \right)
\left( \func{Im}f\left( t\right) -m\func{Im}f\left( s\right) \right) \right]
\end{multline*}%
hence a sufficient condition for the inequality in (\ref{RTxxe.3.4.1}) to
hold is%
\begin{equation}
m\func{Re}f\left( s\right) \leq \func{Re}f\left( t\right) \leq M\func{Re}%
f\left( s\right) \text{ }  \label{RTxxe.3.4.2}
\end{equation}%
and%
\begin{equation*}
m\func{Im}f\left( s\right) \leq \func{Im}f\left( t\right) \leq M\func{Im}%
f\left( s\right)
\end{equation*}%
for a.e. $a\leq t\leq s\leq b.$

Utilising Theorems \ref{RTxxt3.2.4}, \ref{RTxxt3.2.5} and \ref{RTxxt3.3.1}
we may state the following results incorporating quadratic reverses of the
continuous triangle inequality \cite{RTCSSD3}:

\begin{proposition}
\label{RTxxp.3.4.1} With the above assumptions for $f,M$ and $m,$ and if (%
\ref{RTxxe.3.4.1}) holds true, then we have the inequalities%
\begin{multline*}
\quad \left( \int_{a}^{b}\left\vert f\left( t\right) \right\vert
dt\right)
^{2}\leq \left\vert \int_{a}^{b}f\left( t\right) dt\right\vert ^{2} \\
+\frac{1}{2}\cdot \frac{\left( M-m\right)
^{2}}{M+m}\int_{a}^{b}\left( s-a\right) \left\vert f\left(
s\right) \right\vert ^{2}ds,\qquad
\end{multline*}%
\begin{equation*}
\int_{a}^{b}\left\vert f\left( t\right) \right\vert dt\leq \left( \frac{M+m}{%
2\sqrt{Mm}}\right) ^{\frac{1}{2}}\left\vert \int_{a}^{b}f\left( t\right)
dt\right\vert ,
\end{equation*}%
and%
\begin{equation*}
\int_{a}^{b}\left[ \left( b-s\right) +mM\left( s-a\right) \right] \left\vert
f\left( s\right) \right\vert ^{2}ds\leq \frac{M+m}{2}\left\vert
\int_{a}^{b}f\left( s\right) ds\right\vert ^{2}.
\end{equation*}
\end{proposition}

\begin{remark}
\label{RTxxr.3.4.1} One may wonder if there are functions satisfying the
condition (\ref{RTxxe.3.4.2}) above. It suffices to find examples of real
functions $\varphi :\left[ a,b\right] \rightarrow \mathbb{R}$ verifying the
following double inequality%
\begin{equation}
\gamma \varphi \left( s\right) \leq \varphi \left( t\right) \leq \Gamma
\varphi \left( s\right)  \label{RTxxe.3.4.3}
\end{equation}%
for some given $\gamma ,\Gamma $ with $0\leq \gamma \leq 1\leq \Gamma
<\infty $ for a.e. $a\leq t\leq s\leq b.$

For this purpose, consider $\psi :\left[ a,b\right] \rightarrow \mathbb{R}$
a differentiable function on $\left( a,b\right) $, continuous on $\left[ a,b%
\right] $ and with the property that there exists $\Theta \geq 0\geq \theta $
such that%
\begin{equation}
\theta \leq \psi ^{\prime }\left( u\right) \leq \Theta \text{ for any }u\in
\left( a,b\right) .  \label{RTxxe.3.4.4}
\end{equation}%
By Lagrange's mean value theorem, we have, for any $a\leq t\leq s\leq b$%
\begin{equation*}
\psi \left( s\right) -\psi \left( t\right) =\psi ^{\prime }\left( \xi
\right) \left( s-t\right)
\end{equation*}%
with $t\leq \xi \leq s.$ Therefore, for $a\leq t\leq s\leq b,$ by (\ref%
{RTxxe.3.4.4}), we have the inequality%
\begin{equation*}
\theta \left( b-a\right) \leq \theta \left( s-t\right) \leq \psi \left(
s\right) -\psi \left( t\right) \leq \Theta \left( s-t\right) \leq \Theta
\left( b-a\right) .
\end{equation*}%
If we choose the function $\varphi :\left[ a,b\right] \rightarrow \mathbb{R}$
given by%
\begin{equation*}
\varphi \left( t\right) :=\exp \left[ -\psi \left( t\right) \right] ,\text{ }%
t\in \left[ a,b\right] ,
\end{equation*}%
and $\gamma :=\exp \left[ \theta \left( b-a\right) \right] \leq 1,$ $\Gamma
:=\exp \left[ \Theta \left( b-a\right) \right] ,$ then (\ref{RTxxe.3.4.3})
holds true for any $a\leq t\leq s\leq b.$
\end{remark}

The following reverse of the continuous triangle inequality for
complex-valued functions that improves Karamata's result (\ref{RTxx0.1.1})
holds \cite{RTCSSD4}.

\begin{proposition}
\label{RTxxp4.4.1}Let $f\in L\left( \left[ a,b\right] ;\mathbb{C}\right) $
with the property that%
\begin{equation}
0\leq \varphi _{1}\leq \arg f\left( t\right) \leq \varphi _{2}<\frac{\pi }{2}
\label{RTxx4.4.1}
\end{equation}%
for a.e. $t\in \left[ a,b\right] .$ Then we have the inequality%
\begin{equation}
\sqrt{\sin ^{2}\varphi _{1}+\cos ^{2}\varphi _{2}}\int_{a}^{b}\left\vert
f\left( t\right) \right\vert dt\leq \left\vert \int_{a}^{b}f\left( t\right)
dt\right\vert .  \label{RTxx4.4.2}
\end{equation}%
The equality holds in (\ref{RTxx4.4.2}) if and only if%
\begin{equation}
\int_{a}^{b}f\left( t\right) dt=\left( \cos \varphi _{2}+i\sin \varphi
_{1}\right) \int_{a}^{b}\left\vert f\left( t\right) \right\vert dt.
\label{RTxx4.4.3}
\end{equation}
\end{proposition}

\begin{proof}
Let $f\left( t\right) =\func{Re}f\left( t\right) +i\func{Im}f\left( t\right)
.$ We may assume that $\func{Re}f\left( t\right) \geq 0,$ $\func{Im}f\left(
t\right) >0,$ for a.e. $t\in \left[ a,b\right] ,$ since, by (\ref{RTxx4.4.1}%
), $\frac{\func{Im}f\left( t\right) }{\func{Re}f\left( t\right) }=\tan \left[
\arg f\left( t\right) \right] \in \left[ 0,\frac{\pi }{2}\right) $, for a.e.
$t\in \left[ a,b\right] .$ By (\ref{RTxx4.4.1}), we obviously have%
\begin{equation*}
0\leq \tan ^{2}\varphi _{1}\leq \left[ \frac{\func{Im}f\left( t\right) }{%
\func{Re}f\left( t\right) }\right] ^{2}\leq \tan ^{2}\varphi _{2},
\end{equation*}%
for a.e. $t\in \left[ a,b\right] ,$ from where we get%
\begin{equation*}
\frac{\left[ \func{Im}f\left( t\right) \right] ^{2}+\left[ \func{Re}f\left(
t\right) \right] ^{2}}{\left[ \func{Re}f\left( t\right) \right] ^{2}}\leq
\frac{1}{\cos {}^{2}\varphi _{2}},
\end{equation*}%
for a.e. $t\in \left[ a,b\right] ,$ and%
\begin{equation*}
\frac{\left[ \func{Im}f\left( t\right) \right] ^{2}+\left[ \func{Re}f\left(
t\right) \right] ^{2}}{\left[ \func{Im}f\left( t\right) \right] ^{2}}\leq
\frac{1+\tan ^{2}\varphi _{1}}{\tan ^{2}\varphi _{1}}=\frac{1}{\sin \varphi
_{1}},
\end{equation*}%
for a.e. $t\in \left[ a,b\right] ,$ giving the simpler inequalities%
\begin{equation*}
\left\vert f\left( t\right) \right\vert \cos \varphi _{2}\leq \func{Re}%
\left( f\left( t\right) \right) ,\quad \left\vert f\left( t\right)
\right\vert \sin \varphi _{1}\leq \func{Im}\left( f\left( t\right) \right)
\end{equation*}%
for a.e. $t\in \left[ a,b\right] .$

Now, applying Theorem \ref{RTxxt4.2.1} for the complex Hilbert space $%
\mathbb{C}$ endowed with the inner product $\left\langle z,w\right\rangle
=z\cdot \bar{w} $ for $k_{1}=\cos \varphi _{2},$ $k_{2}=\sin \varphi _{1}$
and $e=1,$ we deduce the desired inequality (\ref{RTxx4.4.2}). The case of
equality is also obvious and we omit the details.
\end{proof}

Another result that has an obvious geometrical interpretation is the
following one \cite{RTCSSD4}.

\begin{proposition}
\label{RTxxp4.4.2}Let $e\in \mathbb{C}$ with $\left\vert e\right\vert =1$
and $\rho _{1},\rho _{2}\in \left( 0,1\right) .$ If $f\left( t\right) \in
L\left( \left[ a,b\right] ;\mathbb{C}\right) $ such that%
\begin{equation}
\left\vert f\left( t\right) -e\right\vert \leq \rho _{1},\quad \left\vert
f\left( t\right) -ie\right\vert \leq \rho _{2}\qquad \text{for a.e. }t\in %
\left[ a,b\right] ,  \label{RTxx4.4.5}
\end{equation}%
then we have the inequality%
\begin{equation}
\sqrt{2-\rho _{1}^{2}-\rho _{2}^{2}}\int_{a}^{b}\left\vert f\left( t\right)
\right\vert dt\leq \left\vert \int_{a}^{b}f\left( t\right) dt\right\vert ,
\label{RTxx4.4.6}
\end{equation}%
with equality if and only if%
\begin{equation}
\int_{a}^{b}f\left( t\right) dt=\left( \sqrt{1-\rho _{1}^{2}}+i\sqrt{1-\rho
_{2}^{2}}\right) \int_{a}^{b}\left\vert f\left( t\right) \right\vert dt\cdot
e.  \label{RTxx4.4.7}
\end{equation}
\end{proposition}

The proof is obvious by Corollary \ref{RTxxc4.2.2} applied for $H=\mathbb{C}$
and we omit the details.

\begin{remark}
If we choose $e=1,$ and for $\rho _{1},\rho _{2}\in \left( 0,1\right) $ we
define
\begin{equation*}
\bar{D}\left( 1,\rho _{1}\right) :=\left\{ z\in \mathbb{C}|\left\vert
z-1\right\vert \leq \rho _{1}\right\} ,\ \ \bar{D}\left( i,\rho _{2}\right)
:=\left\{ z\in \mathbb{C}|\left\vert z-i\right\vert \leq \rho _{2}\right\} ,
\end{equation*}%
then obviously the intersection domain%
\begin{equation*}
S_{\rho _{1},\rho _{2}}:=\bar{D}\left( 1,\rho _{1}\right) \cap \bar{D}\left(
i,\rho _{2}\right)
\end{equation*}%
is nonempty if and only if $\rho _{1}+\rho _{2}>\sqrt{2}.$

If $f\left( t\right) \in S_{\rho _{1},\rho _{2}}$ for a.e. $t\in \left[ a,b%
\right] ,$ then (\ref{RTxx4.4.6}) holds true. The equality holds in (\ref%
{RTxx4.4.6}) if and only if%
\begin{equation*}
\int_{a}^{b}f\left( t\right) dt=\left( \sqrt{1-\rho _{1}^{2}}+i\sqrt{1-\rho
_{2}^{2}}\right) \int_{a}^{b}\left\vert f\left( t\right) \right\vert dt.
\end{equation*}
\end{remark}

%


\chapter[CBS and Heisenberg Inequalities]{Reverses of the CBS and Heisenberg
Inequalities}\label{ch5}

\section{Introduction\label{Heisens1}}

Assume that $\left( K;\left\langle \cdot ,\cdot \right\rangle \right) $ is a
Hilbert space over the real or complex number field $\mathbb{K}$. If $\rho :%
\left[ a,b\right] \subset \mathbb{R}\rightarrow \lbrack 0,\infty )$ is a
Lebesgue integrable function with $\int_{a}^{b}\rho \left( t\right) dt=1,$
then we may consider the space $L_{\rho }^{2}\left( \left[ a,b\right]
;K\right) $ of all functions$f:\left[ a,b\right] \rightarrow K,$ that are
Bochner measurable and $\int_{a}^{b}\rho \left( t\right) \left\Vert f\left(
t\right) \right\Vert ^{2}dt<\infty .$ It is well known that $L_{\rho
}^{2}\left( \left[ a,b\right] ;K\right) $ endowed with the inner product $%
\left\langle \cdot ,\cdot \right\rangle _{\rho }$ defined by%
\begin{equation*}
\left\langle f,g\right\rangle _{\rho }:=\int_{a}^{b}\rho \left( t\right)
\left\langle f\left( t\right) ,g\left( t\right) \right\rangle dt
\end{equation*}%
and generating the norm%
\begin{equation*}
\left\Vert f\right\Vert _{\rho }:=\left( \int_{a}^{b}\rho \left( t\right)
\left\Vert f\left( t\right) \right\Vert ^{2}dt\right) ^{\frac{1}{2}},
\end{equation*}%
is a Hilbert space over $\mathbb{K}$.

The following integral inequality is known in the literature as the
Cauchy-Bunyakovsky-Schwarz (CBS) inequality%
\begin{multline}
\quad \int_{a}^{b}\rho \left( t\right) \left\Vert f\left( t\right)
\right\Vert ^{2}dt\int_{a}^{b}\rho \left( t\right) \left\Vert g\left(
t\right) \right\Vert ^{2}dt  \label{Heisen1.1} \\
\geq \left\vert \int_{a}^{b}\rho \left( t\right) \left\langle f\left(
t\right) ,g\left( t\right) \right\rangle dt\right\vert ^{2},\quad
\end{multline}%
provided $f,g\in L_{\rho }^{2}\left( \left[ a,b\right] ;K\right) .$

The case of equality holds in (\ref{Heisen1.1}) iff there exists a constant $%
\lambda \in \mathbb{K}$ such that $f\left( t\right) =\lambda g\left(
t\right) $ for a.e. $t\in \left[ a,b\right] .$

Another version of the (CBS) inequality for one vector-valued and one scalar
function is incorporated in:%
\begin{multline}
\quad\int_{a}^{b}\rho \left( t\right) \left\vert \alpha \left( t\right)
\right\vert ^{2}dt\int_{a}^{b}\rho \left( t\right) \left\Vert f\left(
t\right) \right\Vert ^{2}dt  \label{Heisen1.2} \\
\geq \left\Vert \int_{a}^{b}\rho \left( t\right) \alpha \left( t\right)
f\left( t\right) dt\right\Vert ^{2},\quad
\end{multline}%
provided $\alpha \in L_{\rho }^{2}\left( \left[ a,b\right] \right) $ and $%
f\in L_{\rho }^{2}\left( \left[ a,b\right] ;K\right) ,$ where $L_{\rho
}^{2}\left( \left[ a,b\right] \right) $ denotes the Hilbert space of scalar
functions $\alpha $ for which $\int_{a}^{b}\rho \left( t\right) \left\vert
\alpha \left( t\right) \right\vert ^{2}dt<\infty .$ The equality holds in (%
\ref{Heisen1.2}) iff there exists a vector $e\in K$ such that $f\left(
t\right) =\overline{\alpha \left( t\right) }e$ for a.e. $t\in \left[ a,b%
\right] .$

In this chapter some reverses of the inequalities (\ref{Heisen1.1}) and (\ref%
{Heisen1.2}) are given under various assumptions for the functions involved.
Natural applications for the Heisenberg inequality for vector-valued
functions in Hilbert spaces are also provided.

\section{Some Reverse Inequalities}

\subsection{The General Case}

The following result holds \cite{xxSSD1}.

\begin{theorem}[Dragomir, 2004]
\label{Heisent2.1}Let $f,g\in L_{\rho }^{2}\left( \left[ a,b\right]
;K\right) $ and $r>0$ be such that%
\begin{equation}
\left\Vert f\left( t\right) -g\left( t\right) \right\Vert \leq r\leq
\left\Vert g\left( t\right) \right\Vert  \label{Heisen2.1}
\end{equation}%
for a.e. $t\in \left[ a,b\right] .$ Then we have the inequalities:%
\begin{align}
0& \leq \int_{a}^{b}\rho \left( t\right) \left\Vert f\left( t\right)
\right\Vert ^{2}dt\int_{a}^{b}\rho \left( t\right) \left\Vert g\left(
t\right) \right\Vert ^{2}dt  \label{Heisen2.2} \\
& \qquad \qquad -\left\vert \int_{a}^{b}\rho \left( t\right) \left\langle
f\left( t\right) ,g\left( t\right) \right\rangle dt\right\vert ^{2}  \notag
\\
& \leq \int_{a}^{b}\rho \left( t\right) \left\Vert f\left( t\right)
\right\Vert ^{2}dt\int_{a}^{b}\rho \left( t\right) \left\Vert g\left(
t\right) \right\Vert ^{2}dt  \notag \\
& \qquad \qquad -\left[ \int_{a}^{b}\rho \left( t\right) \func{Re}%
\left\langle f\left( t\right) ,g\left( t\right) \right\rangle dt\right] ^{2}
\notag \\
& \leq r^{2}\int_{a}^{b}\rho \left( t\right) \left\Vert f\left( t\right)
\right\Vert ^{2}dt.  \notag
\end{align}%
The constant $C=1$ in front of $r^{2}$ is best possible in the sense that it
cannot be replaced by a smaller quantity.
\end{theorem}

\begin{proof}
We will use the following result obtained in \cite{xxSSD2}:

In the inner product space $\left( H;\left\langle \cdot ,\cdot \right\rangle
\right) ,$ if $x,y\in H$ and $r>0$ are such that $\left\Vert x-y\right\Vert
\leq r\leq \left\Vert y\right\Vert ,$ then%
\begin{align}
0& \leq \left\Vert x\right\Vert ^{2}\left\Vert y\right\Vert ^{2}-\left\vert
\left\langle x,y\right\rangle \right\vert ^{2}  \label{Heisen2.3} \\
& \leq \left\Vert x\right\Vert ^{2}\left\Vert y\right\Vert ^{2}-\left[ \func{%
Re}\left\langle x,y\right\rangle \right] ^{2}\leq r^{2}\left\Vert
x\right\Vert ^{2}.  \notag
\end{align}%
The constant $c=1$ in front of $r^{2}$ is best possible in the sense that it
cannot be replaced by a smaller quantity.

If (\ref{Heisen2.1}) holds, true, then%
\begin{equation*}
\left\Vert f-g\right\Vert _{\rho }^{2}=\int_{a}^{b}\rho \left( t\right)
\left\Vert f\left( t\right) -g\left( t\right) \right\Vert ^{2}dt\leq
r^{2}\int_{a}^{b}\rho \left( t\right) dt=r^{2}
\end{equation*}%
and%
\begin{equation*}
\left\Vert g\right\Vert _{\rho }^{2}=\int_{a}^{b}\rho \left( t\right)
\left\Vert g\left( t\right) \right\Vert ^{2}dt\geq r^{2}\int_{a}^{b}\rho
\left( t\right) dt=r^{2}
\end{equation*}%
and thus $\left\Vert f-g\right\Vert _{\rho }\leq r\leq \left\Vert
g\right\Vert _{\rho }.$ Applying the inequality (\ref{Heisen2.3}) for $%
\left( L_{\rho }^{2}\left( \left[ a,b\right] ;K\right) ,\left\langle \cdot
,\cdot \right\rangle _{p}\right) ,$ we deduce the desired inequality (\ref%
{Heisen2.2}).

If we choose $\rho \left( t\right) =\frac{1}{b-a},$ $f\left( t\right) =x,$ $%
g\left( t\right) =y,$ $x,y\in K,$ $t\in \left[ a,b\right] ,$ then from (\ref%
{Heisen2.2}) we recapture (\ref{Heisen2.3}) for which the constant $c=1$ in
front of $r^{2}$ is best possible.
\end{proof}

We next point out some general reverse inequalities for the second (CBS)
inequality (\ref{Heisen1.2}) \cite{xxSSD1}.

\begin{theorem}[Dragomir, 2004]
\label{Heisent2.2}Let $\alpha \in L_{\rho }^{2}\left( \left[ a,b\right]
\right) ,$\ $g\in L_{\rho }^{2}\left( \left[ a,b\right] ;K\right) $ and $%
a\in K,$ $r>0$ such that $\left\Vert a\right\Vert >r.$ If the following
condition holds%
\begin{equation}
\left\Vert g\left( t\right) -\bar{\alpha}\left( t\right) a\right\Vert \leq
r\left\vert \alpha \left( t\right) \right\vert  \label{Heisen2.4}
\end{equation}%
for a.e. $t\in \left[ a,b\right] ,$ (note that, if $\alpha \left( t\right)
\neq 0$ for a.e. $t\in \left[ a,b\right] ,$ then the condition (\ref%
{Heisen2.4}) is equivalent to%
\begin{equation}
\left\Vert \frac{g\left( t\right) }{\bar{\alpha}\left( t\right) }%
-a\right\Vert \leq r  \label{Heisen2.4'}
\end{equation}%
for a.e. $t\in \left[ a,b\right] ),$ then we have the following inequality%
\begin{align}
& \left( \int_{a}^{b}\rho \left( t\right) \left\vert \alpha \left( t\right)
\right\vert ^{2}dt\int_{a}^{b}\rho \left( t\right) \left\Vert g\left(
t\right) \right\Vert ^{2}dt\right) ^{\frac{1}{2}}  \label{Heisen2.5} \\
& \leq \frac{1}{\sqrt{\left\Vert a\right\Vert ^{2}-r^{2}}}\func{Re}%
\left\langle \int_{a}^{b}\rho \left( t\right) \alpha \left( t\right) g\left(
t\right) dt,a\right\rangle  \notag \\
& \leq \frac{\left\Vert a\right\Vert }{\sqrt{\left\Vert a\right\Vert
^{2}-r^{2}}}\left\Vert \int_{a}^{b}\rho \left( t\right) \alpha \left(
t\right) g\left( t\right) dt\right\Vert ;  \notag
\end{align}%
\begin{align}
0& \leq \left( \int_{a}^{b}\rho \left( t\right) \left\vert \alpha \left(
t\right) \right\vert ^{2}dt\int_{a}^{b}\rho \left( t\right) \left\Vert
g\left( t\right) \right\Vert ^{2}dt\right) ^{\frac{1}{2}}  \label{Heisen2.6}
\\
& \qquad \qquad -\left\Vert \int_{a}^{b}\rho \left( t\right) \alpha \left(
t\right) g\left( t\right) dt\right\Vert  \notag \\
& \leq \left( \int_{a}^{b}\rho \left( t\right) \left\vert \alpha \left(
t\right) \right\vert ^{2}dt\int_{a}^{b}\rho \left( t\right) \left\Vert
g\left( t\right) \right\Vert ^{2}dt\right) ^{\frac{1}{2}}  \notag \\
& \qquad \qquad -\func{Re}\left\langle \int_{a}^{b}\rho \left( t\right)
\alpha \left( t\right) g\left( t\right) dt,\frac{a}{\left\Vert a\right\Vert }%
\right\rangle  \notag \\
& \leq \frac{r^{2}}{\sqrt{\left\Vert a\right\Vert ^{2}-r^{2}}\left(
\left\Vert a\right\Vert +\sqrt{\left\Vert a\right\Vert ^{2}-r^{2}}\right) }
\notag \\
& \qquad \qquad \times \func{Re}\left\langle \int_{a}^{b}\rho \left(
t\right) \alpha \left( t\right) g\left( t\right) dt,\frac{a}{\left\Vert
a\right\Vert }\right\rangle  \notag \\
& \leq \frac{r^{2}}{\sqrt{\left\Vert a\right\Vert ^{2}-r^{2}}\left(
\left\Vert a\right\Vert +\sqrt{\left\Vert a\right\Vert ^{2}-r^{2}}\right) }%
\left\Vert \int_{a}^{b}\rho \left( t\right) \alpha \left( t\right) g\left(
t\right) dt\right\Vert ;  \notag
\end{align}%
\begin{align}
& \int_{a}^{b}\rho \left( t\right) \left\vert \alpha \left( t\right)
\right\vert ^{2}dt\int_{a}^{b}\rho \left( t\right) \left\Vert g\left(
t\right) \right\Vert ^{2}dt  \label{Heisen2.7} \\
& \leq \frac{1}{\left\Vert a\right\Vert ^{2}-r^{2}}\left[ \func{Re}%
\left\langle \int_{a}^{b}\rho \left( t\right) \alpha \left( t\right) g\left(
t\right) dt,a\right\rangle \right] ^{2}  \notag \\
& \leq \frac{\left\Vert a\right\Vert ^{2}}{\left\Vert a\right\Vert ^{2}-r^{2}%
}\left\Vert \int_{a}^{b}\rho \left( t\right) \alpha \left( t\right) g\left(
t\right) dt\right\Vert ^{2},  \notag
\end{align}%
and%
\begin{align}
0& \leq \int_{a}^{b}\rho \left( t\right) \left\vert \alpha \left( t\right)
\right\vert ^{2}dt\int_{a}^{b}\rho \left( t\right) \left\Vert g\left(
t\right) \right\Vert ^{2}dt  \label{Heisen2.8} \\
& \qquad \qquad -\left\Vert \int_{a}^{b}\rho \left( t\right) \alpha \left(
t\right) g\left( t\right) dt\right\Vert ^{2}  \notag \\
& \leq \int_{a}^{b}\rho \left( t\right) \left\vert \alpha \left( t\right)
\right\vert ^{2}dt\int_{a}^{b}\rho \left( t\right) \left\Vert g\left(
t\right) \right\Vert ^{2}dt  \notag \\
& \qquad \qquad -\left[ \func{Re}\left\langle \int_{a}^{b}\rho \left(
t\right) \alpha \left( t\right) g\left( t\right) dt,\frac{a}{\left\Vert
a\right\Vert }\right\rangle \right] ^{2}  \notag \\
& \leq \frac{r^{2}}{\left\Vert a\right\Vert ^{2}\left( \left\Vert
a\right\Vert ^{2}-r^{2}\right) }\left[ \func{Re}\left\langle
\int_{a}^{b}\rho \left( t\right) \alpha \left( t\right) g\left( t\right)
dt,a\right\rangle \right] ^{2}  \notag \\
& \leq \frac{r^{2}}{\left\Vert a\right\Vert ^{2}-r^{2}}\left\Vert
\int_{a}^{b}\rho \left( t\right) \alpha \left( t\right) g\left( t\right)
dt\right\Vert ^{2}.  \notag
\end{align}%
All the inequalities (\ref{Heisen2.5}) -- (\ref{Heisen2.8}) are sharp.
\end{theorem}

\begin{proof}
From (\ref{Heisen2.4}) we deduce%
\begin{equation*}
\left\Vert g\left( t\right) \right\Vert ^{2}-2\func{Re}\left\langle g\left(
t\right) ,\bar{\alpha}\left( t\right) a\right\rangle +\left\vert \alpha
\left( t\right) \right\vert ^{2}\left\Vert a\right\Vert ^{2}\leq \left\vert
\alpha \left( t\right) \right\vert ^{2}r^{2}
\end{equation*}%
for a.e. $t\in \left[ a,b\right] ,$ which is clearly equivalent to:%
\begin{equation}
\left\Vert g\left( t\right) \right\Vert ^{2}+\left( \left\Vert a\right\Vert
^{2}-r^{2}\right) \left\vert \alpha \left( t\right) \right\vert ^{2}\leq 2%
\func{Re}\left\langle \alpha \left( t\right) g\left( t\right) ,a\right\rangle
\label{Heisen2.9}
\end{equation}%
for a.e. $t\in \left[ a,b\right] .$

If we multiply (\ref{Heisen2.9}) by $\rho \left( t\right) \geq 0$ and
integrate over $t\in \left[ a,b\right] ,$ then we deduce%
\begin{multline}
\int_{a}^{b}\rho \left( t\right) \left\Vert g\left( t\right) \right\Vert
^{2}dt+\left( \left\Vert a\right\Vert ^{2}-r^{2}\right) \int_{a}^{b}\rho
\left( t\right) \left\vert \alpha \left( t\right) \right\vert ^{2}dt
\label{Heisen2.10} \\
\leq 2\func{Re}\left\langle \int_{a}^{b}\rho \left( t\right) \alpha \left(
t\right) g\left( t\right) dt,a\right\rangle .
\end{multline}%
Now, dividing (\ref{Heisen2.10}) by $\sqrt{\left\Vert a\right\Vert ^{2}-r^{2}%
}>0,$ we get%
\begin{multline}
\frac{1}{\sqrt{\left\Vert a\right\Vert ^{2}-r^{2}}}\int_{a}^{b}\rho \left(
t\right) \left\Vert g\left( t\right) \right\Vert ^{2}dt  \label{Heisen2.11}
\\
+\sqrt{\left\Vert a\right\Vert ^{2}-r^{2}}\int_{a}^{b}\rho \left( t\right)
\left\vert \alpha \left( t\right) \right\vert ^{2}dt \\
\leq \frac{2}{\sqrt{\left\Vert a\right\Vert ^{2}-r^{2}}}\func{Re}%
\left\langle \int_{a}^{b}\rho \left( t\right) \alpha \left( t\right) g\left(
t\right) dt,a\right\rangle .
\end{multline}%
On the other hand, by the elementary inequality%
\begin{equation*}
\frac{1}{\alpha }p+\alpha q\geq 2\sqrt{pq},\qquad \alpha >0,\ p,q\geq 0,
\end{equation*}%
we can state that%
\begin{multline}
2\sqrt{\int_{a}^{b}\rho \left( t\right) \left\vert \alpha \left( t\right)
\right\vert ^{2}dt}\cdot \sqrt{\int_{a}^{b}\rho \left( t\right) \left\Vert
g\left( t\right) \right\Vert ^{2}dt}  \label{Heisen2.12} \\
\leq \frac{1}{\sqrt{\left\Vert a\right\Vert ^{2}-r^{2}}}\int_{a}^{b}\rho
\left( t\right) \left\Vert g\left( t\right) \right\Vert ^{2}dt \\
+\sqrt{\left\Vert a\right\Vert ^{2}-r^{2}}\int_{a}^{b}\rho \left( t\right)
\left\vert \alpha \left( t\right) \right\vert ^{2}dt.
\end{multline}%
Making use of (\ref{Heisen2.11}) and (\ref{Heisen2.12}), we deduce the first
part of (\ref{Heisen2.5}).

The second part of (\ref{Heisen2.5}) is obvious by Schwarz's inequality%
\begin{equation*}
\func{Re}\left\langle \int_{a}^{b}\rho \left( t\right) \alpha \left(
t\right) g\left( t\right) dt,a\right\rangle \leq \left\Vert \int_{a}^{b}\rho
\left( t\right) \alpha \left( t\right) g\left( t\right) dt\right\Vert
\left\Vert a\right\Vert .
\end{equation*}%
If $\rho \left( t\right) =\frac{1}{b-a},\ \alpha \left( t\right) =1,$ $%
g\left( t\right) =x\in K,$ then, from (\ref{Heisen2.5}) we get%
\begin{equation*}
\left\Vert x\right\Vert \leq \frac{1}{\sqrt{\left\Vert a\right\Vert
^{2}-r^{2}}}\func{Re}\left\langle x,a\right\rangle \leq \frac{\left\Vert
x\right\Vert \left\Vert a\right\Vert }{\sqrt{\left\Vert a\right\Vert
^{2}-r^{2}}},
\end{equation*}%
provided $\left\Vert x-a\right\Vert \leq r<\left\Vert a\right\Vert ,$ $%
x,a\in K.$ The sharpness of this inequality has been shown in \cite{xxSSD2},
and we omit the details.

The other inequalities are obvious consequences of (\ref{Heisen2.5}) and we
omit the details.
\end{proof}

\subsection{Some Particular Cases}

It has been shown in \cite{xxSSD2} that, for $A,a\in \mathbb{K}$ $\left(
\mathbb{K}=\mathbb{C},\mathbb{R}\right) $ and $x,y\in H,$ where $\left(
H;\left\langle \cdot ,\cdot \right\rangle \right) $ is an inner product over
the real or complex number field $\mathbb{K}$, the following inequality holds%
\begin{align}
\left\Vert x\right\Vert \left\Vert y\right\Vert & \leq \frac{1}{2}\cdot
\frac{\func{Re}\left[ \left( \bar{A}+\bar{a}\right) \left\langle
x,y\right\rangle \right] }{\left[ \func{Re}\left( A\bar{a}\right) \right] ^{%
\frac{1}{2}}}  \label{Heisen3.1} \\
& \leq \frac{1}{2}\cdot \frac{\left\vert A+a\right\vert }{\left[ \func{Re}%
\left( A\bar{a}\right) \right] ^{\frac{1}{2}}}\left\vert \left\langle
x,y\right\rangle \right\vert  \notag
\end{align}%
provided $\func{Re}\left( A\bar{a}\right) >0$ and%
\begin{equation}
\func{Re}\left\langle Ay-x,x-ay\right\rangle \geq 0,  \label{Heisen3.2}
\end{equation}%
or, equivalently,%
\begin{equation}
\left\Vert x-\frac{a+A}{2}\cdot y\right\Vert \leq \frac{1}{2}\left\vert
A-a\right\vert \left\Vert y\right\Vert ,  \label{Heisen3.3}
\end{equation}%
holds. The constant $\frac{1}{2}$ is best possible in (\ref{Heisen3.1}).

From (\ref{Heisen3.1}), we can deduce the following results%
\begin{align}
0& \leq \left\Vert x\right\Vert \left\Vert y\right\Vert -\func{Re}%
\left\langle x,y\right\rangle  \label{Heisen3.4} \\
& \leq \frac{1}{2}\cdot \frac{\func{Re}\left[ \left( \bar{A}+\bar{a}-2\left[
\func{Re}\left( A\bar{a}\right) \right] ^{\frac{1}{2}}\right) \left\langle
x,y\right\rangle \right] }{\left[ \func{Re}\left( A\bar{a}\right) \right] ^{%
\frac{1}{2}}}  \notag \\
& \leq \frac{1}{2}\cdot \frac{\left\vert \bar{A}+\bar{a}-2\left[ \func{Re}%
\left( A\bar{a}\right) \right] ^{\frac{1}{2}}\right\vert }{\left[ \func{Re}%
\left( A\bar{a}\right) \right] ^{\frac{1}{2}}}\left\vert \left\langle
x,y\right\rangle \right\vert  \notag
\end{align}%
and%
\begin{align}
0& \leq \left\Vert x\right\Vert \left\Vert y\right\Vert -\left\vert
\left\langle x,y\right\rangle \right\vert  \label{Heisen3.5} \\
& \leq \frac{1}{2}\cdot \frac{\left\vert A+a\right\vert -2\left[ \func{Re}%
\left( A\bar{a}\right) \right] ^{\frac{1}{2}}}{\left[ \func{Re}\left( A\bar{a%
}\right) \right] ^{\frac{1}{2}}}\left\vert \left\langle x,y\right\rangle
\right\vert .  \notag
\end{align}%
If one assumes that $A=M,$ $a=m,$ $M\geq m>0,$ then, from (\ref{Heisen3.1}),
(\ref{Heisen3.4}) and (\ref{Heisen3.5}) we deduce the much simpler and more
useful results:%
\begin{equation}
\left\Vert x\right\Vert \left\Vert y\right\Vert \leq \frac{1}{2}\cdot \frac{%
M+m}{\sqrt{Mm}}\func{Re}\left\langle x,y\right\rangle ,  \label{Heisen3.6}
\end{equation}%
\begin{equation}
0\leq \left\Vert x\right\Vert \left\Vert y\right\Vert -\func{Re}\left\langle
x,y\right\rangle \leq \frac{1}{2}\cdot \frac{\left( \sqrt{M}-\sqrt{m}\right)
^{2}}{\sqrt{Mm}}\func{Re}\left\langle x,y\right\rangle  \label{Heisen3.7}
\end{equation}%
and%
\begin{equation}
0\leq \left\Vert x\right\Vert \left\Vert y\right\Vert -\left\vert
\left\langle x,y\right\rangle \right\vert \leq \frac{1}{2}\cdot \frac{\left(
\sqrt{M}-\sqrt{m}\right) ^{2}}{\sqrt{Mm}}\left\vert \left\langle
x,y\right\rangle \right\vert ,  \label{Heisen3.8}
\end{equation}%
provided%
\begin{equation*}
\func{Re}\left\langle My-x,x-my\right\rangle \geq 0
\end{equation*}%
or, equivalently%
\begin{equation}
\left\Vert x-\frac{M+m}{2}y\right\Vert \leq \frac{1}{2}\left( M-m\right)
\left\Vert y\right\Vert .  \label{Heisen3.9}
\end{equation}%
Squaring the second inequality in (\ref{Heisen3.1}), we can get the
following results as well:%
\begin{equation}
0\leq \left\Vert x\right\Vert ^{2}\left\Vert y\right\Vert ^{2}-\left\vert
\left\langle x,y\right\rangle \right\vert ^{2}\leq \frac{1}{4}\cdot \frac{%
\left\vert A-a\right\vert ^{2}}{\func{Re}\left( A\bar{a}\right) }\left\vert
\left\langle x,y\right\rangle \right\vert ^{2},  \label{Heisen3.10}
\end{equation}%
provided (\ref{Heisen3.2}) or (\ref{Heisen3.1}) holds. Here the constant $%
\frac{1}{4}$ is also best possible.

Using the above inequalities for vectors in inner product spaces, we are
able to state the following theorem concerning reverses of the (CBS)
integral inequality for vector-valued functions in Hilbert spaces \cite%
{xxSSD1}.

\begin{theorem}[Dragomir, 2004]
\label{Heisent3.1}Let $f,g\in L_{\rho }^{2}\left( \left[ a,b\right]
;K\right) $ and $\gamma ,\Gamma \in \mathbb{K}$ with $\func{Re}\left( \Gamma
\bar{\gamma}\right) >0.$ If%
\begin{equation}
\func{Re}\left\langle \Gamma g\left( t\right) -f\left( t\right) ,f\left(
t\right) -\gamma g\left( t\right) \right\rangle \geq 0  \label{Heisen3.11}
\end{equation}%
for a.e. $t\in \left[ a,b\right] ,$ or, equivalently,%
\begin{equation}
\left\Vert f\left( t\right) -\frac{\gamma +\Gamma }{2}\cdot g\left( t\right)
\right\Vert \leq \frac{1}{2}\left\vert \Gamma -\gamma \right\vert \left\Vert
g\left( t\right) \right\Vert  \label{Heisen3.12}
\end{equation}%
for a.e. $t\in \left[ a,b\right] ,$ then we have the inequalities%
\begin{gather}
\left( \int_{a}^{b}\rho \left( t\right) \left\Vert f\left( t\right)
\right\Vert ^{2}dt\int_{a}^{b}\rho \left( t\right) \left\Vert g\left(
t\right) \right\Vert ^{2}dt\right) ^{\frac{1}{2}}  \label{Heisen3.1a} \\
\leq \frac{1}{2}\cdot \frac{\func{Re}\left[ \left( \bar{\Gamma}+\bar{\gamma}%
\right) \int_{a}^{b}\rho \left( t\right) \left\langle f\left( t\right)
,g\left( t\right) \right\rangle dt\right] }{\left[ \func{Re}\left( \Gamma
\bar{\gamma}\right) \right] ^{\frac{1}{2}}}  \notag \\
\leq \frac{1}{2}\cdot \frac{\left\vert \Gamma +\gamma \right\vert }{\left[
\func{Re}\left( \Gamma \bar{\gamma}\right) \right] ^{\frac{1}{2}}}\left\vert
\int_{a}^{b}\rho \left( t\right) \left\langle f\left( t\right) ,g\left(
t\right) \right\rangle dt\right\vert ,  \notag
\end{gather}%
\begin{align}
0& \leq \left( \int_{a}^{b}\rho \left( t\right) \left\Vert f\left( t\right)
\right\Vert ^{2}dt\right) ^{\frac{1}{2}}\left( \int_{a}^{b}\rho \left(
t\right) \left\Vert g\left( t\right) \right\Vert ^{2}dt\right) ^{\frac{1}{2}}
\label{Heisen3.2a} \\
& \qquad \qquad -\int_{a}^{b}\rho \left( t\right) \func{Re}\left\langle
f\left( t\right) ,g\left( t\right) \right\rangle dt  \notag \\
& \leq \frac{1}{2}\cdot \frac{\func{Re}\left[ \left\{ \bar{\Gamma}+\bar{%
\gamma}-2\left[ \func{Re}\left( \Gamma \bar{\gamma}\right) \right] ^{\frac{1%
}{2}}\right\} \int_{a}^{b}\rho \left( t\right) \left\langle f\left( t\right)
,g\left( t\right) \right\rangle dt\right] }{\left[ \func{Re}\left( \Gamma
\bar{\gamma}\right) \right] ^{\frac{1}{2}}}  \notag \\
& \leq \frac{1}{2}\cdot \frac{\left\vert \bar{\Gamma}+\bar{\gamma}-2\left[
\func{Re}\left( \Gamma \bar{\gamma}\right) \right] ^{\frac{1}{2}}\right\vert
}{\left[ \func{Re}\left( \Gamma \bar{\gamma}\right) \right] ^{\frac{1}{2}}}%
\left\vert \int_{a}^{b}\rho \left( t\right) \left\langle f\left( t\right)
,g\left( t\right) \right\rangle dt\right\vert ,  \notag
\end{align}%
\begin{align}
0& \leq \left( \int_{a}^{b}\rho \left( t\right) \left\Vert f\left( t\right)
\right\Vert ^{2}dt\right) ^{\frac{1}{2}}\left( \int_{a}^{b}\rho \left(
t\right) \left\Vert g\left( t\right) \right\Vert ^{2}dt\right) ^{\frac{1}{2}}
\label{Heisen3.3a} \\
& \qquad \qquad -\left\vert \int_{a}^{b}\rho \left( t\right) \left\langle
f\left( t\right) ,g\left( t\right) \right\rangle dt\right\vert  \notag \\
& \leq \frac{1}{2}\cdot \frac{\left\vert \Gamma +\gamma \right\vert -2\left[
\func{Re}\left( \Gamma \bar{\gamma}\right) \right] ^{\frac{1}{2}}}{\left[
\func{Re}\left( \Gamma \bar{\gamma}\right) \right] ^{\frac{1}{2}}}\left\vert
\int_{a}^{b}\rho \left( t\right) \left\langle f\left( t\right) ,g\left(
t\right) \right\rangle dt\right\vert ,  \notag
\end{align}%
and%
\begin{align}
0& \leq \int_{a}^{b}\rho \left( t\right) \left\Vert f\left( t\right)
\right\Vert ^{2}dt\int_{a}^{b}\rho \left( t\right) \left\Vert g\left(
t\right) \right\Vert ^{2}dt  \label{Heisen3.4a} \\
& \qquad \qquad -\left\vert \int_{a}^{b}\rho \left( t\right) \left\langle
f\left( t\right) ,g\left( t\right) \right\rangle dt\right\vert ^{2}  \notag
\\
& \leq \frac{1}{4}\cdot \frac{\left\vert \Gamma -\gamma \right\vert ^{2}}{%
\func{Re}\left( \Gamma \bar{\gamma}\right) }\left\vert \int_{a}^{b}\rho
\left( t\right) \left\langle f\left( t\right) ,g\left( t\right)
\right\rangle dt\right\vert ^{2}.  \notag
\end{align}%
The constants $\frac{1}{2}$ and $\frac{1}{4}$ above are sharp.
\end{theorem}

In the case where $\Gamma ,\gamma $ are positive real numbers, the following
corollary incorporating more convenient reverses for the (CBS) integral
inequality, may be stated \cite{xxSSD1}.

\begin{corollary}
\label{Heisenc3.2}Let $f,g\in L_{\rho }^{2}\left( \left[ a,b\right]
;K\right) $ and $M\geq m>0.$ If%
\begin{equation}
\func{Re}\left\langle Mg\left( t\right) -f\left( t\right) ,f\left( t\right)
-mg\left( t\right) \right\rangle \geq 0  \label{Heisen3.5a}
\end{equation}%
for a.e. $t\in \left[ a,b\right] ,$ or, equivalently,%
\begin{equation}
\left\Vert f\left( t\right) -\frac{m+M}{2}\cdot g\left( t\right) \right\Vert
\leq \frac{1}{2}\left( M-m\right) \left\Vert g\left( t\right) \right\Vert
\label{Heisen3.6a}
\end{equation}%
for a.e. $t\in \left[ a,b\right] ,$ then we have the inequalities%
\begin{multline}
\left( \int_{a}^{b}\rho \left( t\right) \left\Vert f\left( t\right)
\right\Vert ^{2}dt\int_{a}^{b}\rho \left( t\right) \left\Vert g\left(
t\right) \right\Vert ^{2}dt\right) ^{\frac{1}{2}}  \label{Heisen3.7a} \\
\leq \frac{1}{2}\cdot \frac{M+m}{\sqrt{mM}}\int_{a}^{b}\rho \left( t\right)
\func{Re}\left\langle f\left( t\right) ,g\left( t\right) \right\rangle dt,
\end{multline}%
\begin{align}
0& \leq \left( \int_{a}^{b}\rho \left( t\right) \left\Vert f\left( t\right)
\right\Vert ^{2}dt\right) ^{\frac{1}{2}}\left( \int_{a}^{b}\rho \left(
t\right) \left\Vert g\left( t\right) \right\Vert ^{2}dt\right) ^{\frac{1}{2}}
\label{Heisen3.8a} \\
& \qquad \qquad -\int_{a}^{b}\rho \left( t\right) \func{Re}\left\langle
f\left( t\right) ,g\left( t\right) \right\rangle dt  \notag \\
& \leq \frac{1}{2}\cdot \frac{\left( \sqrt{M}-\sqrt{m}\right) ^{2}}{\sqrt{mM}%
}\int_{a}^{b}\rho \left( t\right) \func{Re}\left\langle f\left( t\right)
,g\left( t\right) \right\rangle dt,  \notag
\end{align}%
\begin{align}
0& \leq \left( \int_{a}^{b}\rho \left( t\right) \left\Vert f\left( t\right)
\right\Vert ^{2}dt\right) ^{\frac{1}{2}}\left( \int_{a}^{b}\rho \left(
t\right) \left\Vert g\left( t\right) \right\Vert ^{2}dt\right) ^{\frac{1}{2}}
\label{Heisen3.9a} \\
& \qquad \qquad -\left\vert \int_{a}^{b}\rho \left( t\right) \left\langle
f\left( t\right) ,g\left( t\right) \right\rangle dt\right\vert  \notag \\
& \leq \frac{1}{2}\cdot \frac{\left( \sqrt{M}-\sqrt{m}\right) ^{2}}{\sqrt{mM}%
}\left\vert \int_{a}^{b}\rho \left( t\right) \left\langle f\left( t\right)
,g\left( t\right) \right\rangle dt\right\vert ,  \notag
\end{align}%
and%
\begin{align}
0& \leq \int_{a}^{b}\rho \left( t\right) \left\Vert f\left( t\right)
\right\Vert ^{2}dt\int_{a}^{b}\rho \left( t\right) \left\Vert g\left(
t\right) \right\Vert ^{2}dt  \label{Heisen3.10a} \\
& \qquad \qquad -\left\vert \int_{a}^{b}\rho \left( t\right) \left\langle
f\left( t\right) ,g\left( t\right) \right\rangle dt\right\vert ^{2}  \notag
\\
& \leq \frac{1}{4}\cdot \frac{\left( M-m\right) ^{2}}{mM}\left\vert
\int_{a}^{b}\rho \left( t\right) \left\langle f\left( t\right) ,g\left(
t\right) \right\rangle dt\right\vert ^{2}.  \notag
\end{align}%
The constants $\frac{1}{2}$ and $\frac{1}{4}$ above are best possible.
\end{corollary}

On utilising the general result of Theorem \ref{Heisent2.2}, we are able to
state a number of interesting reverses for the (CBS) inequality in the case
when one function takes vector-values while the other is a scalar function
\cite{xxSSD1}.

\begin{theorem}[Dragomir, 2004]
\label{Heisent3.3}Let $\alpha \in L_{\rho }^{2}\left( \left[ a,b\right]
\right) ,$ $g\in L_{\rho }^{2}\left( \left[ a,b\right] ;K\right) ,$ $e\in K,$
$\left\Vert e\right\Vert =1,$ $\gamma ,\Gamma \in \mathbb{K}$ with $\func{Re}%
\left( \Gamma \bar{\gamma}\right) >0.$ If
\begin{equation}
\left\Vert g\left( t\right) -\bar{\alpha}\left( t\right) \cdot \frac{\Gamma
+\gamma }{2}e\right\Vert \leq \frac{1}{2}\left\vert \Gamma -\gamma
\right\vert \left\vert \alpha \left( t\right) \right\vert
\label{Heisen3.11a}
\end{equation}%
for a.e. $t\in \left[ a,b\right] ,$ or, equivalently%
\begin{equation}
\func{Re}\left\langle \Gamma \bar{\alpha}\left( t\right) e-g\left( t\right)
,g\left( t\right) -\gamma \bar{\alpha}\left( t\right) e\right\rangle \geq 0
\label{Heisen3.12a}
\end{equation}%
for a.e. $t\in \left[ a,b\right] ,$ (note that, if $\alpha \left( t\right)
\neq 0$ for a.e. $t\in \left[ a,b\right] ,$ then (\ref{Heisen3.11a}) is
equivalent to%
\begin{equation}
\left\Vert \frac{g\left( t\right) }{\overline{\alpha \left( t\right) }}-%
\frac{\Gamma +\gamma }{2}e\right\Vert \leq \frac{1}{2}\left\vert \Gamma
-\gamma \right\vert   \label{Heisen3.11'}
\end{equation}%
for a.e. $t\in \left[ a,b\right] ,$ and (\ref{Heisen3.12a}) is equivalent to
\begin{equation}
\func{Re}\left\langle \Gamma e-\frac{g\left( t\right) }{\overline{\alpha
\left( t\right) }},\frac{g\left( t\right) }{\overline{\alpha \left( t\right)
}}-\gamma e\right\rangle \geq 0  \label{Heisen3.12'}
\end{equation}%
for a.e. $t\in \left[ a,b\right] $), then the following reverse inequalities
are valid:%
\begin{align}
& \left( \int_{a}^{b}\rho \left( t\right) \left\vert \alpha \left( t\right)
\right\vert ^{2}dt\int_{a}^{b}\rho \left( t\right) \left\Vert g\left(
t\right) \right\Vert ^{2}dt\right) ^{\frac{1}{2}}  \label{Heisen3.13} \\
& \leq \frac{\func{Re}\left[ \left( \bar{\Gamma}+\bar{\gamma}\right)
\left\langle \int_{a}^{b}\rho \left( t\right) \alpha \left( t\right) g\left(
t\right) dt,e\right\rangle \right] }{2\left[ \func{Re}\left( \Gamma \bar{%
\gamma}\right) \right] ^{\frac{1}{2}}}  \notag \\
& \leq \frac{1}{2}\cdot \frac{\left\vert \Gamma +\gamma \right\vert }{\left[
\func{Re}\left( \Gamma \bar{\gamma}\right) \right] ^{\frac{1}{2}}}\left\Vert
\int_{a}^{b}\rho \left( t\right) \alpha \left( t\right) g\left( t\right)
dt\right\Vert ;  \notag
\end{align}%
\begin{align}
0& \leq \left( \int_{a}^{b}\rho \left( t\right) \left\vert \alpha \left(
t\right) \right\vert ^{2}dt\int_{a}^{b}\rho \left( t\right) \left\Vert
g\left( t\right) \right\Vert ^{2}dt\right) ^{\frac{1}{2}}  \label{Heisen3.14}
\\
& \qquad \qquad -\left\Vert \int_{a}^{b}\rho \left( t\right) \alpha \left(
t\right) g\left( t\right) dt\right\Vert   \notag \\
& \leq \left( \int_{a}^{b}\rho \left( t\right) \left\vert \alpha \left(
t\right) \right\vert ^{2}dt\int_{a}^{b}\rho \left( t\right) \left\Vert
g\left( t\right) \right\Vert ^{2}dt\right) ^{\frac{1}{2}}  \notag \\
& \qquad \qquad -\func{Re}\left[ \frac{\bar{\Gamma}+\bar{\gamma}}{\left\vert
\Gamma +\gamma \right\vert }\left\langle \int_{a}^{b}\rho \left( t\right)
\alpha \left( t\right) g\left( t\right) dt,e\right\rangle \right]   \notag
\end{align}%
\begin{align}
& \leq \frac{\left\vert \Gamma -\gamma \right\vert ^{2}}{2\sqrt{\func{Re}%
\left( \Gamma \bar{\gamma}\right) }\left( \left\vert \Gamma +\gamma
\right\vert +2\sqrt{\func{Re}\left( \Gamma \bar{\gamma}\right) }\right) }
\notag \\
& \qquad \qquad \times \func{Re}\left[ \frac{\bar{\Gamma}+\bar{\gamma}}{%
\left\vert \Gamma +\gamma \right\vert }\left\langle \int_{a}^{b}\rho \left(
t\right) \alpha \left( t\right) g\left( t\right) dt,e\right\rangle \right]
\notag \\
& \leq \frac{\left\vert \Gamma -\gamma \right\vert ^{2}}{2\sqrt{\func{Re}%
\left( \Gamma \bar{\gamma}\right) }\left( \left\vert \Gamma +\gamma
\right\vert +2\sqrt{\func{Re}\left( \Gamma \bar{\gamma}\right) }\right) }%
\left\Vert \int_{a}^{b}\rho \left( t\right) \alpha \left( t\right) g\left(
t\right) dt\right\Vert ;  \notag
\end{align}%
\begin{align}
& \int_{a}^{b}\rho \left( t\right) \left\vert \alpha \left( t\right)
\right\vert ^{2}dt\int_{a}^{b}\rho \left( t\right) \left\Vert g\left(
t\right) \right\Vert ^{2}dt  \label{Heisen3.15} \\
& \leq \frac{1}{4}\cdot \frac{1}{\func{Re}\left( \Gamma \bar{\gamma}\right) }%
\left[ \func{Re}\left( \left( \overline{\Gamma }+\overline{\gamma }\right)
\left\langle \int_{a}^{b}\rho \left( t\right) \alpha \left( t\right) g\left(
t\right) dt,e\right\rangle \right) \right] ^{2}  \notag \\
& \leq \frac{1}{4}\cdot \frac{\left\vert \Gamma +\gamma \right\vert ^{2}}{%
\func{Re}\left( \Gamma \bar{\gamma}\right) }\left\Vert \int_{a}^{b}\rho
\left( t\right) \alpha \left( t\right) g\left( t\right) dt\right\Vert ^{2}
\notag
\end{align}%
and%
\begin{align}
0& \leq \int_{a}^{b}\rho \left( t\right) \left\vert \alpha \left( t\right)
\right\vert ^{2}dt\int_{a}^{b}\rho \left( t\right) \left\Vert g\left(
t\right) \right\Vert ^{2}dt  \label{Heisen3.16} \\
& \qquad \qquad -\left\Vert \int_{a}^{b}\rho \left( t\right) \alpha \left(
t\right) g\left( t\right) dt\right\Vert ^{2}  \notag \\
& \leq \int_{a}^{b}\rho \left( t\right) \left\vert \alpha \left( t\right)
\right\vert ^{2}dt\int_{a}^{b}\rho \left( t\right) \left\Vert g\left(
t\right) \right\Vert ^{2}dt  \notag \\
& \qquad \qquad -\left[ \func{Re}\left( \frac{\overline{\Gamma }+\overline{%
\gamma }}{\left\vert \Gamma +\gamma \right\vert }\left\langle
\int_{a}^{b}\rho \left( t\right) \alpha \left( t\right) g\left( t\right)
dt,e\right\rangle \right) \right] ^{2}  \notag \\
& \leq \frac{1}{4}\cdot \frac{\left\vert \Gamma -\gamma \right\vert ^{2}}{%
\left\vert \Gamma +\gamma \right\vert ^{2}\func{Re}\left( \Gamma \bar{\gamma}%
\right) }  \notag \\
& \qquad \qquad \times \left[ \func{Re}\left( \left( \overline{\Gamma }+%
\overline{\gamma }\right) \left\langle \int_{a}^{b}\rho \left( t\right)
\alpha \left( t\right) g\left( t\right) dt,e\right\rangle \right) \right]
^{2}  \notag \\
& \leq \frac{1}{4}\cdot \frac{\left\vert \Gamma -\gamma \right\vert ^{2}}{%
\func{Re}\left( \Gamma \bar{\gamma}\right) }\left\Vert \int_{a}^{b}\rho
\left( t\right) \alpha \left( t\right) g\left( t\right) dt\right\Vert ^{2}.
\notag
\end{align}%
The constants $\frac{1}{2}$ and $\frac{1}{4}$ above are sharp.
\end{theorem}

In the particular case of positive constants, the following simpler version
of the above inequalities may be stated.

\begin{corollary}
\label{Heisenc3.4}Let $\alpha \in L_{\rho }^{2}\left( \left[ a,b\right]
\right) \backslash \left\{ 0\right\} ,$ $g\in L_{\rho }^{2}\left( \left[ a,b%
\right] ;K\right) ,$ $e\in K,$ $\left\Vert e\right\Vert =1$ and $M,m\in
\mathbb{R}$ with $M\geq m>0.$ If
\begin{equation}
\left\Vert \frac{g\left( t\right) }{\bar{\alpha}\left( t\right) }-\frac{M+m}{%
2}\cdot e\right\Vert \leq \frac{1}{2}\left( M-m\right)  \label{Heisen3.17}
\end{equation}%
for a.e. $t\in \left[ a,b\right] ,$ or, equivalently,%
\begin{equation}
\func{Re}\left\langle Me-\frac{g\left( t\right) }{\bar{\alpha}\left(
t\right) },\frac{g\left( t\right) }{\bar{\alpha}\left( t\right) }%
-me\right\rangle \geq 0  \label{Heisen3.18}
\end{equation}%
for a.e. $t\in \left[ a,b\right] ,$ then we have%
\begin{align}
& \left( \int_{a}^{b}\rho \left( t\right) \left\vert \alpha \left( t\right)
\right\vert ^{2}dt\int_{a}^{b}\rho \left( t\right) \left\Vert g\left(
t\right) \right\Vert ^{2}dt\right) ^{\frac{1}{2}}  \label{Heisen3.19} \\
& \leq \frac{1}{2}\cdot \frac{M+m}{\sqrt{Mm}}\func{Re}\left\langle
\int_{a}^{b}\rho \left( t\right) \alpha \left( t\right) g\left( t\right)
dt,e\right\rangle  \notag \\
& \leq \frac{1}{2}\cdot \frac{M+m}{\sqrt{Mm}}\left\Vert \int_{a}^{b}\rho
\left( t\right) \alpha \left( t\right) g\left( t\right) dt\right\Vert ;
\notag
\end{align}%
\begin{align}
0& \leq \left( \int_{a}^{b}\rho \left( t\right) \left\vert \alpha \left(
t\right) \right\vert ^{2}dt\int_{a}^{b}\rho \left( t\right) \left\Vert
g\left( t\right) \right\Vert ^{2}dt\right) ^{\frac{1}{2}}  \label{Heisen3.20}
\\
& \qquad \qquad -\left\Vert \int_{a}^{b}\rho \left( t\right) \alpha \left(
t\right) g\left( t\right) dt\right\Vert  \notag \\
& \leq \left( \int_{a}^{b}\rho \left( t\right) \left\vert \alpha \left(
t\right) \right\vert ^{2}dt\int_{a}^{b}\rho \left( t\right) \left\Vert
g\left( t\right) \right\Vert ^{2}dt\right) ^{\frac{1}{2}}  \notag \\
& \qquad \qquad -\func{Re}\left\langle \int_{a}^{b}\rho \left( t\right)
\alpha \left( t\right) g\left( t\right) dt,e\right\rangle  \notag \\
& \leq \frac{\left( \sqrt{M}-\sqrt{m}\right) ^{2}}{2\sqrt{Mm}}\func{Re}%
\left\langle \int_{a}^{b}\rho \left( t\right) \alpha \left( t\right) g\left(
t\right) dt,e\right\rangle  \notag \\
& \leq \frac{\left( \sqrt{M}-\sqrt{m}\right) ^{2}}{2\sqrt{Mm}}\left\Vert
\int_{a}^{b}\rho \left( t\right) \alpha \left( t\right) g\left( t\right)
dt\right\Vert  \notag
\end{align}%
\begin{align}
0& \leq \int_{a}^{b}\rho \left( t\right) \left\vert \alpha \left( t\right)
\right\vert ^{2}dt\int_{a}^{b}\rho \left( t\right) \left\Vert g\left(
t\right) \right\Vert ^{2}dt  \label{Heisen3.21} \\
& \leq \frac{1}{4}\cdot \frac{\left( M+m\right) ^{2}}{Mm}\left[ \func{Re}%
\left\langle \int_{a}^{b}\rho \left( t\right) \alpha \left( t\right) g\left(
t\right) dt,e\right\rangle \right] ^{2}  \notag \\
& \leq \frac{1}{4}\cdot \frac{\left( M+m\right) ^{2}}{Mm}\left\Vert
\int_{a}^{b}\rho \left( t\right) \alpha \left( t\right) g\left( t\right)
dt\right\Vert ^{2}  \notag
\end{align}%
and%
\begin{align}
0& \leq \int_{a}^{b}\rho \left( t\right) \left\vert \alpha \left( t\right)
\right\vert ^{2}dt\int_{a}^{b}\rho \left( t\right) \left\Vert g\left(
t\right) \right\Vert ^{2}dt  \label{Heisen3.22} \\
& \qquad \qquad -\left\Vert \int_{a}^{b}\rho \left( t\right) \alpha \left(
t\right) g\left( t\right) dt\right\Vert ^{2}  \notag \\
& \leq \int_{a}^{b}\rho \left( t\right) \left\vert \alpha \left( t\right)
\right\vert ^{2}dt\int_{a}^{b}\rho \left( t\right) \left\Vert g\left(
t\right) \right\Vert ^{2}dt  \notag \\
& \qquad \qquad -\left[ \func{Re}\left\langle \int_{a}^{b}\rho \left(
t\right) \alpha \left( t\right) g\left( t\right) dt,e\right\rangle \right]
^{2}  \notag \\
& \leq \frac{1}{4}\cdot \frac{\left( M-m\right) ^{2}}{Mm}\left[ \func{Re}%
\left\langle \int_{a}^{b}\rho \left( t\right) \alpha \left( t\right) g\left(
t\right) dt,e\right\rangle \right] ^{2}  \notag \\
& \leq \frac{1}{4}\cdot \frac{\left( M-m\right) ^{2}}{Mm}\left\Vert
\int_{a}^{b}\rho \left( t\right) \alpha \left( t\right) g\left( t\right)
dt\right\Vert ^{2}.  \notag
\end{align}%
The constants $\frac{1}{2}$ and $\frac{1}{4}$ above are sharp.
\end{corollary}

\subsection{Reverses of the Heisenberg Inequality}

It is well known that if $\left( H;\left\langle \cdot ,\cdot \right\rangle
\right) $ is a real or complex Hilbert space and $f:\left[ a,b\right]
\subset \mathbb{R\rightarrow }H$ is an \textit{absolutely continuous
vector-valued }function, then $f$ is differentiable almost everywhere on $%
\left[ a,b\right] ,$ the derivative $f^{\prime }:\left[ a,b\right]
\rightarrow H$ is Bochner integrable on $\left[ a,b\right] $ and%
\begin{equation}
f\left( t\right) =\int_{a}^{t}f^{\prime }\left( s\right) ds\qquad \text{for
any \ }t\in \left[ a,b\right] .  \label{Heisen4.1}
\end{equation}

The following theorem provides a version of the Heisenberg inequalities in
the general setting of Hilbert spaces \cite{xxSSD1}.

\begin{theorem}[Dragomir, 2004]
\label{Heisent4.1}Let $\varphi :\left[ a,b\right] \rightarrow H$ be an
absolutely continuous function with the property that $b\left\Vert \varphi
\left( b\right) \right\Vert ^{2}=a\left\Vert \varphi \left( a\right)
\right\Vert ^{2}.$ Then we have the inequality:%
\begin{equation}
\left( \int_{a}^{b}\left\Vert \varphi \left( t\right) \right\Vert
^{2}dt\right) ^{2}\leq 4\int_{a}^{b}t^{2}\left\Vert \varphi \left( t\right)
\right\Vert ^{2}dt\cdot \int_{a}^{b}\left\Vert \varphi ^{\prime }\left(
t\right) \right\Vert ^{2}dt.  \label{Heisen4.2}
\end{equation}%
The constant $4$ is best possible in the sense that it cannot be replaced by
a smaller quantity.
\end{theorem}

\begin{proof}
Integrating by parts, we have successively%
\begin{align}
& \int_{a}^{b}\left\Vert \varphi \left( t\right) \right\Vert ^{2}dt
\label{Heisen4.3} \\
& =t\left\Vert \varphi \left( t\right) \right\Vert ^{2}\bigg|%
_{a}^{b}-\int_{a}^{b}t\frac{d}{dt}\left( \left\Vert \varphi \left( t\right)
\right\Vert ^{2}\right) dt  \notag \\
& =b\left\Vert \varphi \left( b\right) \right\Vert ^{2}-a\left\Vert \varphi
\left( a\right) \right\Vert ^{2}-\int_{a}^{b}t\frac{d}{dt}\left\langle
\varphi \left( t\right) ,\varphi \left( t\right) \right\rangle dt  \notag \\
& =-\int_{a}^{b}t\left[ \left\langle \varphi ^{\prime }\left( t\right)
,\varphi \left( t\right) \right\rangle +\left\langle \varphi \left( t\right)
,\varphi ^{\prime }\left( t\right) \right\rangle \right] dt  \notag \\
& =-2\int_{a}^{b}t\func{Re}\left\langle \varphi ^{\prime }\left( t\right)
,\varphi \left( t\right) \right\rangle dt  \notag \\
& =2\int_{a}^{b}\func{Re}\left\langle \varphi ^{\prime }\left( t\right)
,\left( -t\right) \varphi \left( t\right) \right\rangle dt.  \notag
\end{align}%
If we apply the (CBS) integral inequality%
\begin{equation*}
\int_{a}^{b}\func{Re}\left\langle g\left( t\right) ,h\left( t\right)
\right\rangle dt\leq \left( \int_{a}^{b}\left\Vert g\left( t\right)
\right\Vert ^{2}dt\int_{a}^{b}\left\Vert h\left( t\right) \right\Vert
^{2}dt\right) ^{\frac{1}{2}}
\end{equation*}%
for $g\left( t\right) =\varphi ^{\prime }\left( t\right) ,$ $h\left(
t\right) =-t\varphi \left( t\right) ,$ $t\in \left[ a,b\right] ,$ then we
deduce the desired inequality (\ref{Heisen4.2}).

The fact that $4$ is the best possible constant in (\ref{Heisen4.2}) follows
from the fact that in the (CBS) inequality, the case of equality holds iff $%
g\left( t\right) =\lambda h\left( t\right) $ for a.e. $t\in \left[ a,b\right]
$ and $\lambda $ a given scalar in $\mathbb{K}$. We omit the details.
\end{proof}

For details on the classical Heisenberg inequality, see, for instance, \cite%
{xxHLP}.

The following reverse of the Heisenberg type inequality (\ref{Heisen4.2})
holds \cite{xxSSD1}.

\begin{theorem}[Dragomir, 2004]
\label{Heisent4.2}Assume that $\varphi :\left[ a,b\right] \rightarrow H$ is
as in the hypothesis of Theorem \ref{Heisent4.1}. In addition, if there
exists a $r>0$ such that%
\begin{equation}
\left\Vert \varphi ^{\prime }\left( t\right) -t\varphi \left( t\right)
\right\Vert \leq r\leq \left\Vert \varphi ^{\prime }\left( t\right)
\right\Vert  \label{Heisen4.4}
\end{equation}%
for a.e. $t\in \left[ a,b\right] ,$ then we have the inequalities%
\begin{align}
0& \leq \int_{a}^{b}t^{2}\left\Vert \varphi \left( t\right) \right\Vert
^{2}dt\int_{a}^{b}\left\Vert \varphi ^{\prime }\left( t\right) \right\Vert
^{2}dt-\frac{1}{4}\left( \int_{a}^{b}\left\Vert \varphi \left( t\right)
\right\Vert ^{2}dt\right) ^{2}  \label{Heisen4.5} \\
& \leq r^{2}\int_{a}^{b}t^{2}\left\Vert \varphi \left( t\right) \right\Vert
^{2}dt.  \notag
\end{align}
\end{theorem}

\begin{proof}
We observe, by the identity (\ref{Heisen4.3}), that%
\begin{equation}
\frac{1}{4}\left( \int_{a}^{b}\left\Vert \varphi \left( t\right) \right\Vert
^{2}dt\right) ^{2}=\left( \int_{a}^{b}\func{Re}\left\langle \varphi ^{\prime
}\left( t\right) ,t\varphi \left( t\right) \right\rangle dt\right) ^{2}.
\label{Heisen4.6}
\end{equation}%
Now, if we apply Theorem \ref{Heisent2.1} for the choices $f\left( t\right)
=t\varphi \left( t\right) ,$ $g\left( t\right) =\varphi ^{\prime }\left(
t\right) ,$ and $\rho \left( t\right) =\frac{1}{b-a},$ then by (\ref%
{Heisen2.2}) and (\ref{Heisen4.6}) we deduce the desired inequality (\ref%
{Heisen4.5}).
\end{proof}

\begin{remark}
\label{Heisenr4.3}Interchanging the place of $t\varphi \left( t\right) $
with $\varphi ^{\prime }\left( t\right) $ in Theorem \ref{Heisent4.2}, we
also have%
\begin{align}
0& \leq \int_{a}^{b}t^{2}\left\Vert \varphi \left( t\right) \right\Vert
^{2}dt\int_{a}^{b}\left\Vert \varphi ^{\prime }\left( t\right) \right\Vert
^{2}dt-\frac{1}{4}\left( \int_{a}^{b}\left\Vert \varphi \left( t\right)
\right\Vert ^{2}dt\right) ^{2}  \label{Heisen4.6a} \\
& \leq \rho ^{2}\int_{a}^{b}\left\Vert \varphi ^{\prime }\left( t\right)
\right\Vert ^{2}dt,  \notag
\end{align}%
provided%
\begin{equation*}
\left\Vert \varphi ^{\prime }\left( t\right) -t\varphi \left( t\right)
\right\Vert \leq \rho \leq \left\vert t\right\vert \left\Vert \varphi \left(
t\right) \right\Vert
\end{equation*}%
for a.e. $t\in \left[ a,b\right] ,$ where $\rho >0$ is a given positive
number.
\end{remark}

The following result also holds \cite{xxSSD1}.

\begin{theorem}[Dragomir, 2004]
\label{Heisent4.4}Assume that $\varphi :\left[ a,b\right] \rightarrow H$ is
as in the hypothesis of Theorem \ref{Heisent4.1}. In addition, if there
exists $M\geq m>0$ such that%
\begin{equation}
\func{Re}\left\langle Mt\varphi \left( t\right) -\varphi ^{\prime }\left(
t\right) ,\varphi ^{\prime }\left( t\right) -mt\varphi \left( t\right)
\right\rangle \geq 0  \label{Heisen4.7}
\end{equation}%
for a.e. $t\in \left[ a,b\right] ,$ or, equivalently,%
\begin{equation}
\left\Vert \varphi ^{\prime }\left( t\right) -\frac{M+m}{2}t\varphi \left(
t\right) \right\Vert \leq \frac{1}{2}\left( M-m\right) \left\vert
t\right\vert \left\Vert \varphi \left( t\right) \right\Vert
\label{Heisen4.8}
\end{equation}%
for a.e. $t\in \left[ a,b\right] ,$ then we have the inequalities%
\begin{multline}
\quad \int_{a}^{b}t^{2}\left\Vert \varphi \left( t\right)
\right\Vert ^{2}dt\int_{a}^{b}\left\Vert \varphi ^{\prime }\left(
t\right) \right\Vert
^{2}dt  \label{Heisen4.9} \\
\leq \frac{1}{16}\cdot \frac{\left( M+m\right) ^{2}}{Mm}\left(
\int_{a}^{b}\left\Vert \varphi \left( t\right) \right\Vert
^{2}dt\right) ^{2}\quad
\end{multline}%
and%
\begin{multline}
\int_{a}^{b}t^{2}\left\Vert \varphi \left( t\right) \right\Vert
^{2}dt\int_{a}^{b}\left\Vert \varphi ^{\prime }\left( t\right) \right\Vert
^{2}dt-\frac{1}{4}\left( \int_{a}^{b}\left\Vert \varphi \left( t\right)
\right\Vert ^{2}dt\right) ^{2}  \label{Heisen4.10} \\
\leq \frac{1}{16}\cdot \frac{\left( M-m\right) ^{2}}{Mm}\left(
\int_{a}^{b}\left\Vert \varphi \left( t\right) \right\Vert ^{2}dt\right) ^{2}
\end{multline}%
respectively.
\end{theorem}

\begin{proof}
We use Corollary \ref{Heisenc3.2} for the choices $f\left( t\right) =\varphi
^{\prime }\left( t\right) ,$ $g\left( t\right) =t\varphi \left( t\right) ,$ $%
\rho \left( t\right) =\frac{1}{b-a},$ to get%
\begin{multline*}
\quad \int_{a}^{b}\left\Vert \varphi ^{\prime }\left( t\right)
\right\Vert ^{2}dt\int_{a}^{b}t^{2}\left\Vert \varphi \left(
t\right) \right\Vert ^{2}dt
\\
\leq \frac{\left( M+m\right) ^{2}}{4Mm}\left( \int_{a}^{b}\func{Re}%
\left\langle \varphi ^{\prime }\left( t\right) ,t\varphi \left(
t\right) \right\rangle dt\right) ^{2}.\quad
\end{multline*}%
Since, by (\ref{Heisen4.6})%
\begin{equation*}
\left( \int_{a}^{b}\func{Re}\left\langle \varphi ^{\prime }\left( t\right)
,t\varphi \left( t\right) \right\rangle dt\right) ^{2}=\frac{1}{4}\left(
\int_{a}^{b}\left\Vert \varphi \left( t\right) \right\Vert ^{2}dt\right)
^{2},
\end{equation*}%
hence we deduce the desired result (\ref{Heisen4.9}).

The inequality (\ref{Heisen4.10}) follows from (\ref{Heisen4.9}), and we
omit the details.
\end{proof}

\begin{remark}
If one is interested in reverses for the Heisenberg inequality for scalar
valued functions, then all the other inequalities obtained above for one
scalar function may be applied as well. For the sake of brevity, we do not
list them here.
\end{remark}

\section{Other Reverses}

\subsection{The General Case}

The following result holds \cite{xxSEV2}.

\begin{theorem}[Dragomir, 2004]
\label{Hilbertt2.1}Let $f,g\in L_{\rho }^{2}\left( \left[ a,b\right]
;K\right) $ and $r>0$ be such that%
\begin{equation}
\left\Vert f\left( t\right) -g\left( t\right) \right\Vert \leq r
\label{Hilbert2.1}
\end{equation}%
for a.e. $t\in \left[ a,b\right] .$ Then we have the inequalities:%
\begin{align}
0& \leq \left( \int_{a}^{b}\rho \left( t\right) \left\Vert f\left( t\right)
\right\Vert ^{2}dt\int_{a}^{b}\rho \left( t\right) \left\Vert g\left(
t\right) \right\Vert ^{2}dt\right) ^{\frac{1}{2}}  \label{Hilbert2.2} \\
& \qquad \qquad \qquad -\left\vert \int_{a}^{b}\rho \left( t\right)
\left\langle f\left( t\right) ,g\left( t\right) \right\rangle dt\right\vert
\notag \\
& \leq \left( \int_{a}^{b}\rho \left( t\right) \left\Vert f\left( t\right)
\right\Vert ^{2}dt\int_{a}^{b}\rho \left( t\right) \left\Vert g\left(
t\right) \right\Vert ^{2}dt\right) ^{\frac{1}{2}}  \notag \\
& \qquad \qquad \qquad -\left\vert \int_{a}^{b}\rho \left( t\right) \func{Re}%
\left\langle f\left( t\right) ,g\left( t\right) \right\rangle dt\right\vert
\notag \\
& \leq \left( \int_{a}^{b}\rho \left( t\right) \left\Vert f\left( t\right)
\right\Vert ^{2}dt\int_{a}^{b}\rho \left( t\right) \left\Vert g\left(
t\right) \right\Vert ^{2}dt\right) ^{\frac{1}{2}}  \notag \\
& \qquad \qquad \qquad -\int_{a}^{b}\rho \left( t\right) \func{Re}%
\left\langle f\left( t\right) ,g\left( t\right) \right\rangle dt  \notag \\
& \leq \frac{1}{2}r^{2}.  \notag
\end{align}%
The constant $\frac{1}{2}$ in front of $r^{2}$ is best possible in the sense
that it cannot be replaced by a smaller quantity.
\end{theorem}

\begin{proof}
We will use the following result obtained in \cite{xxSSD2}:

In the inner product space $\left( H;\left\langle \cdot ,\cdot \right\rangle
\right) ,$ if $x,y\in H$ and $r>0$ are such that $\left\Vert x-y\right\Vert
\leq r,$ then%
\begin{align}
0& \leq \left\Vert x\right\Vert \left\Vert y\right\Vert -\left\vert
\left\langle x,y\right\rangle \right\vert \leq \left\Vert x\right\Vert
\left\Vert y\right\Vert -\left\vert \func{Re}\left\langle x,y\right\rangle
\right\vert  \label{Hilbert2.3a} \\
& \leq \left\Vert x\right\Vert \left\Vert y\right\Vert -\func{Re}%
\left\langle x,y\right\rangle \leq \frac{1}{2}r^{2}.  \notag
\end{align}%
The constant $\frac{1}{2}$ in front of $r^{2}$ is best possible in the sense
that it cannot be replaced by a smaller constant.

If (\ref{Hilbert2.1}) holds true, then%
\begin{equation*}
\left\Vert f-g\right\Vert _{\rho }^{2}=\int_{a}^{b}\rho \left( t\right)
\left\Vert f\left( t\right) -g\left( t\right) \right\Vert ^{2}dt\leq
r^{2}\int_{a}^{b}\rho \left( t\right) dt=r^{2}
\end{equation*}%
and thus $\left\Vert f-g\right\Vert _{\rho }\leq r.$

Applying the inequality (\ref{Hilbert2.3a}) for $\left( L_{\rho }^{2}\left( %
\left[ a,b\right] ;K\right) ,\left\langle \cdot ,\cdot \right\rangle
_{p}\right) ,$ we deduce the desired inequality (\ref{Hilbert2.2}).

If we choose $\rho \left( t\right) =\frac{1}{b-a},$ $f\left( t\right) =x,$ $%
g\left( t\right) =y,$ $x,y\in K,$ $t\in \left[ a,b\right] ,$ then from (\ref%
{Hilbert2.2}) we recapture (\ref{Hilbert2.3a}) for which the constant $\frac{%
1}{2}$ in front of $r^{2}$ is best possible.
\end{proof}

We next point out some general reverse inequalities for the second CBS
inequality (\ref{Heisen1.2})\cite{xxSEV2}.

\begin{theorem}[Dragomir, 2004]
\label{Hilbertt2.2}Let $\alpha \in L_{\rho }^{2}\left( \left[ a,b\right]
\right) ,$\ $g\in L_{\rho }^{2}\left( \left[ a,b\right] ;K\right) $ and $%
v\in K,\ r>0.$ If
\begin{equation}
\left\Vert \frac{g\left( t\right) }{\overline{\alpha \left( t\right) }}%
-v\right\Vert \leq r  \label{Hilbert2.4}
\end{equation}%
for a.e. $t\in \left[ a,b\right] ,$ then we have the inequality%
\begin{align}
0& \leq \left( \int_{a}^{b}\rho \left( t\right) \left\vert \alpha \left(
t\right) \right\vert ^{2}dt\int_{a}^{b}\rho \left( t\right) \left\Vert
g\left( t\right) \right\Vert ^{2}dt\right) ^{\frac{1}{2}}  \label{Hilbert2.5}
\\
& \qquad \qquad \qquad -\left\Vert \int_{a}^{b}\rho \left( t\right) \alpha
\left( t\right) g\left( t\right) dt\right\Vert  \notag \\
& \leq \left( \int_{a}^{b}\rho \left( t\right) \left\vert \alpha \left(
t\right) \right\vert ^{2}dt\int_{a}^{b}\rho \left( t\right) \left\Vert
g\left( t\right) \right\Vert ^{2}dt\right) ^{\frac{1}{2}}  \notag \\
& \qquad \qquad \qquad -\left\vert \left\langle \int_{a}^{b}\rho \left(
t\right) \alpha \left( t\right) g\left( t\right) dt,\frac{v}{\left\Vert
v\right\Vert }\right\rangle \right\vert  \notag \\
& \leq \left( \int_{a}^{b}\rho \left( t\right) \left\vert \alpha \left(
t\right) \right\vert ^{2}dt\int_{a}^{b}\rho \left( t\right) \left\Vert
g\left( t\right) \right\Vert ^{2}dt\right) ^{\frac{1}{2}}  \notag \\
& \qquad \qquad \qquad -\left\vert \func{Re}\left\langle \int_{a}^{b}\rho
\left( t\right) \alpha \left( t\right) g\left( t\right) dt,\frac{v}{%
\left\Vert v\right\Vert }\right\rangle \right\vert  \notag \\
& \leq \left( \int_{a}^{b}\rho \left( t\right) \left\vert \alpha \left(
t\right) \right\vert ^{2}dt\int_{a}^{b}\rho \left( t\right) \left\Vert
g\left( t\right) \right\Vert ^{2}dt\right) ^{\frac{1}{2}}  \notag \\
& \qquad \qquad \qquad -\func{Re}\left\langle \int_{a}^{b}\rho \left(
t\right) \alpha \left( t\right) g\left( t\right) dt,\frac{v}{\left\Vert
v\right\Vert }\right\rangle  \notag \\
& \leq \frac{1}{2}\cdot \frac{r^{2}}{\left\Vert v\right\Vert }%
\int_{a}^{b}\rho \left( t\right) \left\vert \alpha \left( t\right)
\right\vert ^{2}dt.  \notag
\end{align}%
The constant $\frac{1}{2}$ is best possible in the sense that it cannot be
replaced by a smaller quantity.
\end{theorem}

\begin{proof}
From (\ref{Hilbert2.4}) we deduce%
\begin{equation*}
\left\Vert g\left( t\right) \right\Vert ^{2}-2\func{Re}\left\langle \alpha
\left( t\right) g\left( t\right) ,v\right\rangle +\left\vert \alpha \left(
t\right) \right\vert ^{2}\left\Vert v\right\Vert ^{2}\leq r^{2}\left\vert
\alpha \left( t\right) \right\vert ^{2}
\end{equation*}%
which is clearly equivalent to%
\begin{equation}
\left\Vert g\left( t\right) \right\Vert ^{2}+\left\vert \alpha \left(
t\right) \right\vert ^{2}\left\Vert v\right\Vert ^{2}\leq 2\func{Re}%
\left\langle \alpha \left( t\right) g\left( t\right) ,v\right\rangle
+r^{2}\left\vert \alpha \left( t\right) \right\vert ^{2}.  \label{Hilbert2.6}
\end{equation}

If we multiply (\ref{Hilbert2.6}) by $\rho \left( t\right) \geq 0$ and
integrate over $t\in \left[ a,b\right] ,$ then we deduce%
\begin{multline}
\int_{a}^{b}\rho \left( t\right) \left\Vert g\left( t\right) \right\Vert
^{2}dt+\left\Vert v\right\Vert ^{2}\int_{a}^{b}\rho \left( t\right)
\left\vert \alpha \left( t\right) \right\vert ^{2}dt  \label{Hilbert2.7} \\
\leq 2\func{Re}\left\langle \int_{a}^{b}\rho \left( t\right) \alpha \left(
t\right) g\left( t\right) dt,v\right\rangle +r^{2}\int_{a}^{b}\rho \left(
t\right) \left\vert \alpha \left( t\right) \right\vert ^{2}dt.
\end{multline}%
Since, obviously%
\begin{multline}
2\left\Vert v\right\Vert \left( \int_{a}^{b}\rho \left( t\right) \left\vert
\alpha \left( t\right) \right\vert ^{2}dt\int_{a}^{b}\rho \left( t\right)
\left\Vert g\left( t\right) \right\Vert ^{2}dt\right) ^{\frac{1}{2}}
\label{Hilbert2.8} \\
\leq \int_{a}^{b}\rho \left( t\right) \left\Vert g\left( t\right)
\right\Vert ^{2}dt+\left\Vert v\right\Vert ^{2}\int_{a}^{b}\rho \left(
t\right) \left\vert \alpha \left( t\right) \right\vert ^{2}dt,
\end{multline}%
hence, by (\ref{Hilbert2.7}) and (\ref{Hilbert2.8}), we deduce%
\begin{multline*}
2\left\Vert v\right\Vert \left( \int_{a}^{b}\rho \left( t\right) \left\vert
\alpha \left( t\right) \right\vert ^{2}dt\int_{a}^{b}\rho \left( t\right)
\left\Vert g\left( t\right) \right\Vert ^{2}dt\right) ^{\frac{1}{2}} \\
\leq 2\func{Re}\left\langle \int_{a}^{b}\rho \left( t\right) \alpha \left(
t\right) g\left( t\right) dt,v\right\rangle +r^{2}\int_{a}^{b}\rho \left(
t\right) \left\vert \alpha \left( t\right) \right\vert ^{2}dt,
\end{multline*}%
which is clearly equivalent with the last inequality in (\ref{Hilbert2.5}).

The other inequalities are obvious and we omit the details.

Now, if $\rho \left( t\right) =\frac{1}{b-a},\ \alpha \left( t\right) =1,$ $%
g\left( t\right) =x,\ x\in K,$ then, by the last inequality in (\ref%
{Hilbert2.5}) we get%
\begin{equation*}
\left\Vert x\right\Vert \left\Vert v\right\Vert -\func{Re}\left\langle
x,v\right\rangle \leq \frac{1}{2}r^{2},
\end{equation*}%
provided $\left\Vert x-v\right\Vert \leq r,$ for which we know that (see
\cite{xxSSD2}), the constant $\frac{1}{2}$ is best possible.
\end{proof}

\subsection{Some Particular Cases of Interest}

It has been shown in \cite{xxSSD2} that, for $\gamma ,\Gamma \in \mathbb{K}$
$\left( \mathbb{K}=\mathbb{C}\text{ or }\mathbb{K}=\mathbb{R}\right) $ with $%
\Gamma \neq -\gamma $ and $x,y\in H,$ $\left( H;\left\langle \cdot ,\cdot
\right\rangle \right) $ is an inner product over the real or complex number
field $\mathbb{K}$, such that either%
\begin{equation}
\func{Re}\left\langle \Gamma y-x,x-\gamma y\right\rangle \geq 0,
\label{Hilbert3.1}
\end{equation}%
or, equivalently,%
\begin{equation}
\left\Vert x-\frac{\gamma +\Gamma }{2}\cdot y\right\Vert \leq \frac{1}{2}%
\left\vert \Gamma -\gamma \right\vert \left\Vert y\right\Vert ,
\label{Hilbert3.2}
\end{equation}%
holds, then one has the following reverse of Schwarz's inequality%
\begin{align}
0& \leq \left\Vert x\right\Vert \left\Vert y\right\Vert -\left\vert
\left\langle x,y\right\rangle \right\vert  \label{Hilbert3.3} \\
& \leq \left\Vert x\right\Vert \left\Vert y\right\Vert -\left\vert \func{Re}%
\left[ \frac{\bar{\Gamma}+\bar{\gamma}}{\left\vert \Gamma +\gamma
\right\vert }\left\langle x,y\right\rangle \right] \right\vert  \notag \\
& \leq \left\Vert x\right\Vert \left\Vert y\right\Vert -\func{Re}\left[
\frac{\bar{\Gamma}+\bar{\gamma}}{\left\vert \Gamma +\gamma \right\vert }%
\left\langle x,y\right\rangle \right]  \notag \\
& \leq \frac{1}{4}\cdot \frac{\left\vert \Gamma -\gamma \right\vert ^{2}}{%
\left\vert \Gamma +\gamma \right\vert }\left\Vert y\right\Vert ^{2}.  \notag
\end{align}%
The constant $\frac{1}{4}$ is best possible in (\ref{Hilbert3.3}) in the
sense that it cannot be replaced by a smaller constant.

If we assume \ that $\Gamma =M,$ $\gamma =m$ with $M\geq m>0,$ then from (%
\ref{Hilbert3.3}) we deduce the much simpler and more useful result%
\begin{align}
0& \leq \left\Vert x\right\Vert \left\Vert y\right\Vert -\left\vert
\left\langle x,y\right\rangle \right\vert \leq \left\Vert x\right\Vert
\left\Vert y\right\Vert -\left\vert \func{Re}\left\langle x,y\right\rangle
\right\vert  \label{Hilbert3.4} \\
& \leq \left\Vert x\right\Vert \left\Vert y\right\Vert -\func{Re}%
\left\langle x,y\right\rangle \leq \frac{1}{4}\cdot \frac{\left( M-m\right)
^{2}}{Mm}\left\Vert y\right\Vert ^{2},  \notag
\end{align}%
provided (\ref{Hilbert3.1}) or (\ref{Hilbert3.2}) holds true with $M$ and $m$
instead of $\Gamma $ and $\gamma .$

Using the above inequalities for vectors in inner product spaces, we are
able to state the following theorem concerning reverses of the CBS integral
inequality for vector-valued functions in Hilbert spaces \cite{xxSEV2}.

\begin{theorem}[Dragomir, 2004]
\label{Hilbertt3.1}Let $f,g\in L_{\rho }^{2}\left( \left[ a,b\right]
;K\right) $ and $\gamma ,\Gamma \in \mathbb{K}$ with $\Gamma \neq -\gamma .$
If%
\begin{equation}
\func{Re}\left\langle \Gamma g\left( t\right) -f\left( t\right) ,f\left(
t\right) -\gamma g\left( t\right) \right\rangle \geq 0  \label{Hilbert3.5}
\end{equation}%
for a.e. $t\in \left[ a,b\right] ,$ or, equivalently,%
\begin{equation}
\left\Vert f\left( t\right) -\frac{\gamma +\Gamma }{2}\cdot g\left( t\right)
\right\Vert \leq \frac{1}{2}\left\vert \Gamma -\gamma \right\vert \left\Vert
g\left( t\right) \right\Vert  \label{Hilbert3.6}
\end{equation}%
for a.e. $t\in \left[ a,b\right] ,$ then we have the inequalities%
\begin{align}
0& \leq \left( \int_{a}^{b}\rho \left( t\right) \left\Vert f\left( t\right)
\right\Vert ^{2}dt\int_{a}^{b}\rho \left( t\right) \left\Vert g\left(
t\right) \right\Vert ^{2}dt\right) ^{\frac{1}{2}}  \label{Hilbert3.7} \\
& \qquad \qquad \qquad -\left\vert \int_{a}^{b}\rho \left( t\right)
\left\langle f\left( t\right) ,g\left( t\right) \right\rangle dt\right\vert
\notag
\end{align}%
\begin{align}
& \leq \left( \int_{a}^{b}\rho \left( t\right) \left\Vert f\left( t\right)
\right\Vert ^{2}dt\int_{a}^{b}\rho \left( t\right) \left\Vert g\left(
t\right) \right\Vert ^{2}dt\right) ^{\frac{1}{2}}  \notag \\
& \qquad \qquad \qquad -\left\vert \func{Re}\left[ \frac{\bar{\Gamma}+\bar{%
\gamma}}{\left\vert \Gamma +\gamma \right\vert }\int_{a}^{b}\rho \left(
t\right) \left\langle f\left( t\right) ,g\left( t\right) \right\rangle dt%
\right] \right\vert  \notag \\
& \leq \left( \int_{a}^{b}\rho \left( t\right) \left\Vert f\left( t\right)
\right\Vert ^{2}dt\int_{a}^{b}\rho \left( t\right) \left\Vert g\left(
t\right) \right\Vert ^{2}dt\right) ^{\frac{1}{2}}  \notag \\
& \qquad \qquad \qquad -\func{Re}\left[ \frac{\bar{\Gamma}+\bar{\gamma}}{%
\left\vert \Gamma +\gamma \right\vert }\int_{a}^{b}\rho \left( t\right)
\left\langle f\left( t\right) ,g\left( t\right) \right\rangle dt\right]
\notag \\
& \leq \frac{1}{4}\cdot \frac{\left\vert \Gamma -\gamma \right\vert ^{2}}{%
\left\vert \Gamma +\gamma \right\vert }\int_{a}^{b}\rho \left( t\right)
\left\Vert g\left( t\right) \right\Vert ^{2}dt.  \notag
\end{align}%
The constant $\frac{1}{4}$ is best possible in (\ref{Hilbert3.7}).
\end{theorem}

\begin{proof}
Since, by (\ref{Hilbert3.5}),%
\begin{multline*}
\quad \func{Re}\left\langle \Gamma g-f,f-\gamma g\right\rangle _{\rho } \\
=\int_{a}^{b}\rho \left( t\right) \func{Re}\left\langle \Gamma
g\left( t\right) -f\left( t\right) ,f\left( t\right) -\gamma
g\left( t\right) \right\rangle dt\geq 0,\quad
\end{multline*}%
hence, by (\ref{Hilbert3.3}) applied for the Hilbert space $\left( L_{\rho
}^{2}\left( \left[ a,b\right] ;K\right) ;\left\langle \cdot ,\cdot
\right\rangle _{\rho }\right) ,$ we deduce the desired inequality (\ref%
{Hilbert3.7}).

The best constant follows by the fact that $\frac{1}{4}$ is a best constant
in (\ref{Hilbert3.7}) and we omit the details.
\end{proof}

\begin{corollary}
\label{Hilbertc3.2}Let $f,g\in L_{\rho }^{2}\left( \left[ a,b\right]
;K\right) $ and $M\geq m>0.$ If%
\begin{equation}
\func{Re}\left\langle Mg\left( t\right) -f\left( t\right) ,f\left( t\right)
-mg\left( t\right) \right\rangle \geq 0  \label{Hilbert3.8}
\end{equation}%
for a.e. $t\in \left[ a,b\right] ,$ or, equivalently,%
\begin{equation}
\left\Vert f\left( t\right) -\frac{m+M}{2}\cdot g\left( t\right) \right\Vert
\leq \frac{1}{2}\left( M-m\right) \left\Vert g\left( t\right) \right\Vert
\label{Hilbert3.9}
\end{equation}%
for a.e. $t\in \left[ a,b\right] ,$ then
\begin{align}
0& \leq \left( \int_{a}^{b}\rho \left( t\right) \left\Vert f\left( t\right)
\right\Vert ^{2}dt\int_{a}^{b}\rho \left( t\right) \left\Vert g\left(
t\right) \right\Vert ^{2}dt\right) ^{\frac{1}{2}}  \label{Hilbert3.10} \\
& \qquad \qquad \qquad -\left\vert \int_{a}^{b}\rho \left( t\right)
\left\langle f\left( t\right) ,g\left( t\right) \right\rangle dt\right\vert
\notag\displaybreak \\
& \leq \left( \int_{a}^{b}\rho \left( t\right) \left\Vert f\left( t\right)
\right\Vert ^{2}dt\int_{a}^{b}\rho \left( t\right) \left\Vert g\left(
t\right) \right\Vert ^{2}dt\right) ^{\frac{1}{2}}  \notag \\
& \qquad \qquad \qquad -\left\vert \int_{a}^{b}\rho \left( t\right) \func{Re}%
\left\langle f\left( t\right) ,g\left( t\right) \right\rangle dt\right\vert
\notag \\
& \leq \left( \int_{a}^{b}\rho \left( t\right) \left\Vert f\left( t\right)
\right\Vert ^{2}dt\int_{a}^{b}\rho \left( t\right) \left\Vert g\left(
t\right) \right\Vert ^{2}dt\right) ^{\frac{1}{2}}  \notag \\
& \qquad \qquad \qquad -\int_{a}^{b}\rho \left( t\right) \func{Re}%
\left\langle f\left( t\right) ,g\left( t\right) \right\rangle dt  \notag \\
& \leq \frac{1}{4}\cdot \frac{\left( M-m\right) ^{2}}{M+m}\int_{a}^{b}\rho
\left( t\right) \left\Vert g\left( t\right) \right\Vert ^{2}dt.  \notag
\end{align}%
The constant $\frac{1}{4}$ is best possible.
\end{corollary}

The case when a function is scalar is incorporated in the following theorem
\cite{xxSEV2}.

\begin{theorem}[Dragomir, 2004]
\label{Hilbertt3.3}Let $\alpha \in L_{\rho }^{2}\left( \left[ a,b\right]
\right) ,$ $g\in L_{\rho }^{2}\left( \left[ a,b\right] ;K\right) ,$ and $%
\gamma ,\Gamma \in \mathbb{K}$ with $\Gamma \neq -\gamma .$ If $e\in K,$ $%
\left\Vert e\right\Vert =1$ and%
\begin{equation}
\left\Vert \frac{g\left( t\right) }{\overline{\alpha \left( t\right) }}-%
\frac{\Gamma +\gamma }{2}e\right\Vert \leq \frac{1}{2}\left\vert \Gamma
-\gamma \right\vert  \label{Hilbert3.11}
\end{equation}%
for a.e. $t\in \left[ a,b\right] ,$ or, equivalently,%
\begin{equation}
\func{Re}\left\langle \Gamma e-\frac{g\left( t\right) }{\overline{\alpha
\left( t\right) }},\frac{g\left( t\right) }{\overline{\alpha \left( t\right)
}}-\gamma e\right\rangle \geq 0  \label{Hilbert3.12}
\end{equation}%
for a.e. $t\in \left[ a,b\right] ,$ then we have the inequalities
\begin{align}
0& \leq \left( \int_{a}^{b}\rho \left( t\right) \left\vert \alpha \left(
t\right) \right\vert ^{2}dt\int_{a}^{b}\rho \left( t\right) \left\Vert
g\left( t\right) \right\Vert ^{2}dt\right) ^{\frac{1}{2}}
\label{Hilbert3.13} \\
& \qquad \qquad \qquad -\left\Vert \int_{a}^{b}\rho \left( t\right) \alpha
\left( t\right) g\left( t\right) dt\right\Vert  \notag \\
& \leq \left( \int_{a}^{b}\rho \left( t\right) \left\vert \alpha \left(
t\right) \right\vert ^{2}dt\int_{a}^{b}\rho \left( t\right) \left\Vert
g\left( t\right) \right\Vert ^{2}dt\right) ^{\frac{1}{2}}  \notag \\
& \qquad \qquad \qquad -\left\vert \left\langle \int_{a}^{b}\rho \left(
t\right) \alpha \left( t\right) g\left( t\right) dt,e\right\rangle
\right\vert  \notag\displaybreak \\
& \leq \left( \int_{a}^{b}\rho \left( t\right) \left\vert \alpha
\left( t\right) \right\vert ^{2}dt\int_{a}^{b}\rho \left( t\right)
\left\Vert
g\left( t\right) \right\Vert ^{2}dt\right) ^{\frac{1}{2}}  \notag \\
& \qquad \qquad \qquad -\left\vert \func{Re}\left[ \frac{\bar{\Gamma}+\bar{%
\gamma}}{\left\vert \Gamma +\gamma \right\vert }\left\langle
\int_{a}^{b}\rho \left( t\right) \alpha \left( t\right) g\left( t\right)
dt,e\right\rangle \right] \right\vert  \notag \\
& \leq \left( \int_{a}^{b}\rho \left( t\right) \left\vert \alpha \left(
t\right) \right\vert ^{2}dt\int_{a}^{b}\rho \left( t\right) \left\Vert
g\left( t\right) \right\Vert ^{2}dt\right) ^{\frac{1}{2}}  \notag \\
& \qquad \qquad \qquad -\func{Re}\left[ \frac{\bar{\Gamma}+\bar{\gamma}}{%
\left\vert \Gamma +\gamma \right\vert }\left\langle \int_{a}^{b}\rho \left(
t\right) \alpha \left( t\right) g\left( t\right) dt,e\right\rangle \right]
\notag \\
& \leq \frac{1}{4}\cdot \frac{\left\vert \Gamma -\gamma \right\vert ^{2}}{%
\left\vert \Gamma +\gamma \right\vert }\int_{a}^{b}\rho \left( t\right)
\left\vert \alpha \left( t\right) \right\vert ^{2}dt.  \notag
\end{align}%
The constant $\frac{1}{4}$ is best possible in (\ref{Hilbert3.13}).
\end{theorem}

\begin{proof}
Follows by Theorem \ref{Hilbertt2.2} on choosing%
\begin{equation*}
v:=\frac{\Gamma +\gamma }{2}e\qquad \text{and}\qquad r:=\frac{1}{2}%
\left\vert \Gamma -\gamma \right\vert .
\end{equation*}%
We omit the details.
\end{proof}

\begin{corollary}
\label{Hilbertc3.4}Let $\alpha \in L_{\rho }^{2}\left( \left[ a,b\right]
\right) ,$ $g\in L_{\rho }^{2}\left( \left[ a,b\right] ;K\right) ,$ and $%
M\geq m>0.$ If $e\in K,$ $\left\Vert e\right\Vert =1$ and%
\begin{equation*}
\left\Vert \frac{g\left( t\right) }{\overline{\alpha \left( t\right) }}-%
\frac{M+m}{2}\cdot e\right\Vert \leq \frac{1}{2}\left( M-m\right)
\end{equation*}%
for a.e. $t\in \left[ a,b\right] ,$ or, equivalently,%
\begin{equation*}
\func{Re}\left\langle Me-\frac{g\left( t\right) }{\overline{\alpha \left(
t\right) }},\frac{g\left( t\right) }{\overline{\alpha \left( t\right) }}%
-me\right\rangle \geq 0
\end{equation*}%
for a.e. $t\in \left[ a,b\right] $, then we have the inequalities:%
\begin{align}
0& \leq \left( \int_{a}^{b}\rho \left( t\right) \left\vert \alpha \left(
t\right) \right\vert ^{2}dt\int_{a}^{b}\rho \left( t\right) \left\Vert
g\left( t\right) \right\Vert ^{2}dt\right) ^{\frac{1}{2}}
\label{Hilbert3.14} \\
& \qquad \qquad \qquad -\left\Vert \int_{a}^{b}\rho \left( t\right) \alpha
\left( t\right) g\left( t\right) dt\right\Vert  \notag \\
& \leq \left( \int_{a}^{b}\rho \left( t\right) \left\vert \alpha \left(
t\right) \right\vert ^{2}dt\int_{a}^{b}\rho \left( t\right) \left\Vert
g\left( t\right) \right\Vert ^{2}dt\right) ^{\frac{1}{2}}  \notag \\
& \qquad \qquad \qquad -\left\vert \left\langle \int_{a}^{b}\rho \left(
t\right) \alpha \left( t\right) g\left( t\right) dt,e\right\rangle
\right\vert  \notag\displaybreak \\
& \leq \left( \int_{a}^{b}\rho \left( t\right) \left\vert \alpha
\left( t\right) \right\vert ^{2}dt\int_{a}^{b}\rho \left( t\right)
\left\Vert
g\left( t\right) \right\Vert ^{2}dt\right) ^{\frac{1}{2}}  \notag \\
& \qquad \qquad \qquad -\left\vert \func{Re}\left\langle \int_{a}^{b}\rho
\left( t\right) \alpha \left( t\right) g\left( t\right) dt,e\right\rangle
\right\vert  \notag \\
& \leq \left( \int_{a}^{b}\rho \left( t\right) \left\vert \alpha \left(
t\right) \right\vert ^{2}dt\int_{a}^{b}\rho \left( t\right) \left\Vert
g\left( t\right) \right\Vert ^{2}dt\right) ^{\frac{1}{2}}  \notag \\
& \qquad \qquad \qquad -\func{Re}\left\langle \int_{a}^{b}\rho \left(
t\right) \alpha \left( t\right) g\left( t\right) dt,e\right\rangle  \notag \\
& \leq \frac{1}{4}\cdot \frac{\left( M-m\right) ^{2}}{M+m}\int_{a}^{b}\rho
\left( t\right) \left\vert \alpha \left( t\right) \right\vert ^{2}dt.  \notag
\end{align}%
The constant $\frac{1}{4}$ is best possible in (\ref{Hilbert3.14}).
\end{corollary}

\subsection{Applications for the Heisenberg Inequality}

The following reverse of the Heisenberg type inequality (\ref{Heisen4.2})
holds \cite{xxSEV2}.

\begin{theorem}[Dragomir, 2004]
\label{Hilbertt4.2}Assume that $\varphi :\left[ a,b\right] \rightarrow H$ is
as in the hypothesis of Theorem \ref{Heisent4.1}. In addition, if there
exists a $r>0$ such that%
\begin{equation}
\left\Vert \varphi ^{\prime }\left( t\right) +t\varphi \left( t\right)
\right\Vert \leq r  \label{Hilbert4.4}
\end{equation}%
for a.e. $t\in \left[ a,b\right] ,$ then we have the inequalities%
\begin{align}
0& \leq \left( \int_{a}^{b}t^{2}\left\Vert \varphi \left( t\right)
\right\Vert ^{2}dt\int_{a}^{b}\left\Vert \varphi ^{\prime }\left( t\right)
\right\Vert ^{2}dt\right) ^{\frac{1}{2}}-\frac{1}{2}\int_{a}^{b}\left\Vert
\varphi \left( t\right) \right\Vert ^{2}dt  \label{Hilbert4.5} \\
& \leq \frac{1}{2}r^{2}\left( b-a\right) .  \notag
\end{align}
\end{theorem}

\begin{proof}
We observe, by the identity (\ref{Heisen4.3}), that%
\begin{equation}
\int_{a}^{b}\func{Re}\left\langle \varphi ^{\prime }\left( t\right) ,\left(
-t\right) \varphi \left( t\right) \right\rangle dt=\frac{1}{2}%
\int_{a}^{b}\left\Vert \varphi \left( t\right) \right\Vert ^{2}dt.
\label{Hilbert4.6}
\end{equation}%
Now, if we apply Theorem \ref{Hilbertt2.1} for the choices $f\left( t\right)
=t\varphi \left( t\right) ,$ $g\left( t\right) =-t\varphi ^{\prime }\left(
t\right) ,$ $\rho \left( t\right) =\frac{1}{b-a},$ $t\in \left[ a,b\right] ,$
then we deduce the desired inequality (\ref{Hilbert4.5}).
\end{proof}

\begin{remark}
\label{Hilbertr4.3}It\ is interesting to remark that, from (\ref{Hilbert4.6}%
), we obviously have%
\begin{equation}
\frac{1}{2}\int_{a}^{b}\left\Vert \varphi \left( t\right) \right\Vert
^{2}dt=\left\vert \int_{a}^{b}\func{Re}\left\langle \varphi ^{\prime }\left(
t\right) ,t\varphi \left( t\right) \right\rangle dt\right\vert .
\label{Hilbert4.6a}
\end{equation}%
Now, if we apply the inequality (see (\ref{Hilbert2.2}))
\begin{equation*}
\int_{a}^{b}\left\Vert f\left( t\right) \right\Vert
^{2}dt\int_{a}^{b}\left\Vert g\left( t\right) \right\Vert ^{2}dt-\left\vert
\int_{a}^{b}\func{Re}\left\langle f\left( t\right) ,g\left( t\right)
\right\rangle dt\right\vert \leq \frac{1}{2}r^{2}\left( b-a\right) ,
\end{equation*}%
for the choices $f\left( t\right) =\varphi ^{\prime }\left( t\right) ,$ $%
g\left( t\right) =t\varphi \left( t\right) ,$ $t\in \left[ a,b\right] ,$
then we get the same inequality (\ref{Hilbert4.5}), but under the condition
\begin{equation}
\left\Vert \varphi ^{\prime }\left( t\right) -t\varphi \left( t\right)
\right\Vert \leq r  \label{Hilbert4.6.b}
\end{equation}%
for a.e. \ $t\in \left[ a,b\right] .$
\end{remark}

The following result holds as well \cite{xxSEV2}.

\begin{theorem}[Dragomir, 2004]
\label{Hilbertt4.4}Assume that $\varphi :\left[ a,b\right] \rightarrow H$ is
as in the hypothesis of Theorem \ref{Hilbertt4.2}. In addition, if there
exists $M\geq m>0$ such that%
\begin{equation}
\func{Re}\left\langle Mt\varphi \left( t\right) -\varphi ^{\prime }\left(
t\right) ,\varphi ^{\prime }\left( t\right) -mt\varphi \left( t\right)
\right\rangle \geq 0  \label{Hilbert4.7}
\end{equation}%
for a.e. $t\in \left[ a,b\right] ,$ or, equivalently,%
\begin{equation}
\left\Vert \varphi ^{\prime }\left( t\right) -\frac{M+m}{2}t\varphi \left(
t\right) \right\Vert \leq \frac{1}{2}\left( M-m\right) \left\vert
t\right\vert \left\Vert \varphi \left( t\right) \right\Vert \qquad
\label{Hilbert4.8}
\end{equation}%
for a.e. $t\in \left[ a,b\right] ,$ then we have the inequalities%
\begin{align}
0& \leq \left( \int_{a}^{b}t^{2}\left\Vert \varphi \left( t\right)
\right\Vert ^{2}dt\int_{a}^{b}\left\Vert \varphi ^{\prime }\left( t\right)
\right\Vert ^{2}dt\right) ^{\frac{1}{2}}-\frac{1}{2}\int_{a}^{b}\left\Vert
\varphi \left( t\right) \right\Vert ^{2}dt  \label{Hilbert4.9} \\
& \leq \frac{1}{4}\cdot \frac{\left( M-m\right) ^{2}}{M+m}%
\int_{a}^{b}t^{2}\left\Vert \varphi \left( t\right) \right\Vert ^{2}dt.
\notag
\end{align}
\end{theorem}

\begin{proof}
The proof follows by Corollary \ref{Hilbertc3.2} applied for the function $%
g\left( t\right) =t\varphi \left( t\right) $ and $f\left( t\right) =\varphi
^{\prime }\left( t\right) ,$ and on making use of the identity (\ref%
{Hilbert4.6a}). We omit the details.
\end{proof}

\begin{remark}
If one is interested in reverses for the Heisenberg inequality for real or
complex valued functions, then all the other inequalities obtained above for
one scalar and one vectorial function may be applied as well. For the sake
of brevity, we do not list them here.
\end{remark}

%

                      

\chapter[Other Inequalities]{Other Inequalities in Inner Product Spaces}

\label{ch6}

\section{Bounds for the Distance to Finite-Dimensional Subspaces}

\subsection{Introduction}

Let $\left( H;\left\langle \cdot ,\cdot \right\rangle \right) $ be an inner
product space over the real or complex number field $\mathbb{K}$, $\left\{
y_{1},\dots ,y_{n}\right\} $ a subset of $H$ and $G\left( y_{1},\dots
,y_{n}\right) $ the \textit{Gram matrix} of $\left\{ y_{1},\dots
,y_{n}\right\} $ where $\left( i,j\right) -$entry is $\left\langle
y_{i},y_{j}\right\rangle .$ The determinant of $G\left( y_{1},\dots
,y_{n}\right) $ is called the \textit{Gram determinant} of $\left\{
y_{1},\dots ,y_{n}\right\} $ and is denoted by $\Gamma \left( y_{1},\dots
,y_{n}\right) .$ Thus,%
\begin{equation*}
\Gamma \left( y_{1},\dots ,y_{n}\right) =\left\vert 
\begin{array}{c}
\left\langle y_{1},y_{1}\right\rangle \;\left\langle
y_{1},y_{2}\right\rangle \;\cdots \;\left\langle y_{1},y_{n}\right\rangle 
\\ 
\left\langle y_{2},y_{1}\right\rangle \;\left\langle
y_{2},y_{2}\right\rangle \;\cdots \;\left\langle y_{2},y_{n}\right\rangle 
\\ 
\cdots \cdots \cdots \cdots \cdots  \\ 
\left\langle y_{n},y_{1}\right\rangle \;\left\langle
y_{n},y_{2}\right\rangle \;\cdots \;\left\langle y_{n},y_{n}\right\rangle 
\end{array}%
\right\vert .
\end{equation*}

Following \cite[p. 129 -- 133]{DExx}, we state here some general results for
the Gram determinant that will be used in the sequel.

\begin{enumerate}
\item Let $\left\{ x_{1},\dots ,x_{n}\right\} \subset H.$ Then $\Gamma
\left( x_{1},\dots ,x_{n}\right) \neq 0$ if and only if $\left\{ x_{1},\dots
,x_{n}\right\} $ is linearly independent;

\item Let $M=span\left\{ x_{1},\dots ,x_{n}\right\} $ be $n-$dimensional in $%
H,$ i.e., $\{x_{1},\dots ,$ $x_{n}\}$ is linearly independent. Then for each 
$x\in H,$ the distance $d\left( x,M\right) $ from $x$ to the linear subspace 
$H$ has the representations%
\begin{equation}
d^{2}\left( x,M\right) =\frac{\Gamma \left( x_{1},\dots ,x_{n},x\right) }{%
\Gamma \left( x_{1},\dots ,x_{n}\right) }  \label{upper1.1}
\end{equation}%
and%
\begin{equation}
d^{2}\left( x,M\right) =\left\Vert x\right\Vert ^{2}-\beta ^{T}G^{-1}\beta ,
\label{upper1.2}
\end{equation}%
where $G=G\left( x_{1},\dots ,x_{n}\right) ,$ $G^{-1}$ is the inverse matrix
of $G$ and 
\begin{equation*}
\beta ^{T}=\left( \left\langle x,x_{1}\right\rangle ,\left\langle
x,x_{2}\right\rangle ,\dots ,\left\langle x,x_{n}\right\rangle \right) ,
\end{equation*}%
denotes the transpose of the column vector $\beta .$

Moreover, one has the simpler representation%
\begin{equation}
d^{2}\left( x,M\right) =\left\{ 
\begin{array}{ll}
\left\Vert x\right\Vert ^{2}-\frac{\left( \sum_{i=1}^{n}\left\vert
\left\langle x,x_{i}\right\rangle \right\vert ^{2}\right) ^{2}}{\left\Vert
\sum_{i=1}^{n}\left\langle x,x_{i}\right\rangle x_{i}\right\Vert ^{2}} & 
\text{if \ }x\notin M^{\perp }, \\ 
&  \\ 
\left\Vert x\right\Vert ^{2} & \text{if \ }x\in M^{\perp },%
\end{array}%
\right.  \label{upper1.3}
\end{equation}%
where $M^{\perp }$ denotes the orthogonal complement of $M.$

\item Let $\left\{ x_{1},\dots ,x_{n}\right\} $ be a set of nonzero vectors
in $H.$ Then%
\begin{equation}
0\leq \Gamma \left( x_{1},\dots ,x_{n}\right) \leq \left\Vert
x_{1}\right\Vert ^{2}\left\Vert x_{2}\right\Vert ^{2}\cdots \left\Vert
x_{n}\right\Vert ^{2}.  \label{upper1.4}
\end{equation}%
The equality holds on the left (respectively right) side of (\ref{upper1.4})
if and only if $\left\{ x_{1},\dots ,x_{n}\right\} $ is linearly dependent
(respectively orthogonal). The first inequality in (\ref{upper1.4}) is known
in the literature as \textit{Gram's inequality} while the second one is
known as \textit{Hadamard's inequality.}

\item If $\left\{ x_{1},\dots ,x_{n}\right\} $ is an orthonormal set in $H,$
i.e., $\left\langle x_{i},x_{j}\right\rangle =\delta _{ij},$ $i,j\in \left\{
1,\dots ,n\right\} ,$ where $\delta _{ij}$ is Kronecker's delta, then%
\begin{equation}
d^{2}\left( x,M\right) =\left\Vert x\right\Vert
^{2}-\sum_{i=1}^{n}\left\vert \left\langle x,x_{i}\right\rangle \right\vert
^{2}.  \label{upper1.5}
\end{equation}
\end{enumerate}

The following inequalities which involve Gram determinants may be stated as
well \cite[p. 597]{MPFxx}:%
\begin{equation}
\frac{\Gamma \left( x_{1},\dots ,x_{n}\right) }{\Gamma \left( x_{1},\dots
,x_{k}\right) }\leq \frac{\Gamma \left( x_{2},\dots ,x_{n}\right) }{\Gamma
\left( x_{1},\dots ,x_{k}\right) }\leq \cdots \leq \Gamma \left(
x_{k+1},\dots ,x_{n}\right) ,  \label{upper1.6}
\end{equation}%
\begin{equation}
\Gamma \left( x_{1},\dots ,x_{n}\right) \leq \Gamma \left( x_{1},\dots
,x_{k}\right) \Gamma \left( x_{k+1},\dots ,x_{n}\right)  \label{upper1.7}
\end{equation}%
and%
\begin{multline}
\quad \Gamma ^{\frac{1}{2}}\left( x_{1}+y_{1},x_{2},\dots ,x_{n}\right)
\label{upper1.8} \\
\leq \Gamma ^{\frac{1}{2}}\left( x_{1},x_{2},\dots ,x_{n}\right) +\Gamma ^{%
\frac{1}{2}}\left( y_{1},x_{2},\dots ,x_{n}\right) .\quad
\end{multline}

The main aim of this section is to point out some upper bounds for the
distance $d\left( x,M\right) $ in terms of the linearly independent vectors $%
\left\{ x_{1},\dots ,x_{n}\right\} $ that span $M$ and $x\notin M^{\perp },$
where $M^{\perp }$ is the orthogonal complement of $M$ in the inner product
space $\left( H;\left\langle \cdot ,\cdot \right\rangle \right) $.

As a by-product of this endeavour, some refinements of the generalisations
for Bessel's inequality due to several authors including: Boas, Bellman and
Bombieri are obtained. Refinements for the well known Hadamard's inequality
for Gram determinants are also derived.

\subsection{Upper Bounds for $d\left( x,M\right) $}

The following result may be stated \cite{SILV1xx}.

\begin{theorem}[Dragomir, 2005]
\label{uppert2.1}Let $\left\{ x_{1},\dots ,x_{n}\right\} $ be a linearly
independent system of vectors in $H$ and $M:=span\left\{ x_{1},\dots
,x_{n}\right\} .$ If $x\notin M^{\perp },$ then%
\begin{equation}
d^{2}\left( x,M\right) <\frac{\left\Vert x\right\Vert
^{2}\sum_{i=1}^{n}\left\Vert x_{i}\right\Vert ^{2}-\sum_{i=1}^{n}\left\vert
\left\langle x,x_{i}\right\rangle \right\vert ^{2}}{\sum_{i=1}^{n}\left\Vert
x_{i}\right\Vert ^{2}}  \label{upper2.1}
\end{equation}%
or, equivalently,%
\begin{multline}
\Gamma \left( x_{1},\dots ,x_{n},x\right)  \label{upper2.2} \\
<\frac{\left\Vert x\right\Vert ^{2}\sum_{i=1}^{n}\left\Vert x_{i}\right\Vert
^{2}-\sum_{i=1}^{n}\left\vert \left\langle x,x_{i}\right\rangle \right\vert
^{2}}{\sum_{i=1}^{n}\left\Vert x_{i}\right\Vert ^{2}}\cdot \Gamma \left(
x_{1},\dots ,x_{n}\right) .
\end{multline}
\end{theorem}

\begin{proof}
If we use the Cauchy-Bunyakovsky-Schwarz type inequality%
\begin{equation}
\left\Vert \sum_{i=1}^{n}\alpha _{i}y_{i}\right\Vert ^{2}\leq
\sum_{i=1}^{n}\left\vert \alpha _{i}\right\vert ^{2}\sum_{i=1}^{n}\left\Vert
y_{i}\right\Vert ^{2},  \label{upper2.3}
\end{equation}%
that can be easily deduced from the obvious identity%
\begin{equation}
\sum_{i=1}^{n}\left\vert \alpha _{i}\right\vert ^{2}\sum_{i=1}^{n}\left\Vert
y_{i}\right\Vert ^{2}-\left\Vert \sum_{i=1}^{n}\alpha _{i}y_{i}\right\Vert
^{2}=\frac{1}{2}\sum_{i,j=1}^{n}\left\Vert \overline{\alpha _{i}}x_{j}-%
\overline{\alpha _{j}}x_{i}\right\Vert ^{2},  \label{upper2.4}
\end{equation}%
we can state that%
\begin{equation}
\left\Vert \sum_{i=1}^{n}\left\langle x,x_{i}\right\rangle x_{i}\right\Vert
^{2}\leq \sum_{i=1}^{n}\left\vert \left\langle x,x_{i}\right\rangle
\right\vert ^{2}\sum_{i=1}^{n}\left\Vert x_{i}\right\Vert ^{2}.
\label{upper2.5}
\end{equation}%
Note that the equality case holds in (\ref{upper2.5}) if and only if, by (%
\ref{upper2.4}), 
\begin{equation}
\overline{\left\langle x,x_{i}\right\rangle }x_{j}=\overline{\left\langle
x,x_{i}\right\rangle }x_{i}  \label{upper2.6}
\end{equation}%
for each $i,j\in \left\{ 1,\dots ,n\right\} .$

Utilising the expression (\ref{upper1.3}) of the distance $d\left(
x,M\right) $, we have%
\begin{equation}
d^{2}\left( x,M\right) =\left\Vert x\right\Vert ^{2}-\frac{%
\sum_{i=1}^{n}\left\vert \left\langle x,x_{i}\right\rangle \right\vert
^{2}\sum_{i=1}^{n}\left\Vert x_{i}\right\Vert ^{2}}{\left\Vert
\sum_{i=1}^{n}\left\langle x,x_{i}\right\rangle x_{i}\right\Vert ^{2}}\cdot 
\frac{\sum_{i=1}^{n}\left\vert \left\langle x,x_{i}\right\rangle \right\vert
^{2}}{\sum_{i=1}^{n}\left\Vert x_{i}\right\Vert ^{2}}.  \label{upper2.7}
\end{equation}%
Since $\left\{ x_{1},\dots ,x_{n}\right\} $ are linearly independent, hence (%
\ref{upper2.6}) cannot be achieved and then we have strict inequality in (%
\ref{upper2.5}).

Finally, on using (\ref{upper2.5}) and (\ref{upper2.7}) we get the desired
result (\ref{upper2.1}).
\end{proof}

\begin{remark}
\label{upperr2.2}It is known that (see (\ref{upper1.4})) if not all $\left\{
x_{1},\dots ,x_{n}\right\} $ are orthogonal on each other, then the
following result, which is well known in the literature as Hadamard's
inequality holds:%
\begin{equation}
\Gamma \left( x_{1},\dots ,x_{n}\right) <\left\Vert x_{1}\right\Vert
^{2}\left\Vert x_{2}\right\Vert ^{2}\cdots \left\Vert x_{n}\right\Vert ^{2}.
\label{upper2.8}
\end{equation}%
Utilising the inequality (\ref{upper2.2}), we may write successively:%
\begin{align*}
\Gamma \left( x_{1},x_{2}\right) & \leq \frac{\left\Vert x_{1}\right\Vert
^{2}\left\Vert x_{2}\right\Vert ^{2}-\left\vert \left\langle
x_{2},x_{1}\right\rangle \right\vert ^{2}}{\left\Vert x_{1}\right\Vert ^{2}}%
\left\Vert x_{1}\right\Vert ^{2}\leq \left\Vert x_{1}\right\Vert
^{2}\left\Vert x_{2}\right\Vert ^{2}, \\
\Gamma \left( x_{1},x_{2},x_{3}\right) & <\frac{\left\Vert x_{3}\right\Vert
^{2}\sum_{i=1}^{2}\left\Vert x_{i}\right\Vert ^{2}-\sum_{i=1}^{2}\left\vert
\left\langle x_{3},x_{i}\right\rangle \right\vert ^{2}}{\sum_{i=1}^{2}\left%
\Vert x_{i}\right\Vert ^{2}}\Gamma \left( x_{1},x_{2}\right)  \\
& \leq \left\Vert x_{3}\right\Vert ^{2}\Gamma \left( x_{1},x_{2}\right)  \\
& \cdots \cdots \cdots \cdots \cdots \cdots \cdots \cdots \cdots \cdots
\cdots \cdots \cdots \cdots  \\
\Gamma \left( x_{1},\dots ,x_{n-1},x_{n}\right) & <\frac{\left\Vert
x_{n}\right\Vert ^{2}\sum_{i=1}^{n-1}\left\Vert x_{i}\right\Vert
^{2}-\sum_{i=1}^{n-1}\left\vert \left\langle x_{n},x_{i}\right\rangle
\right\vert ^{2}}{\sum_{i=1}^{n-1}\left\Vert x_{i}\right\Vert ^{2}} \\
& \qquad \qquad \times \Gamma \left( x_{1},\dots ,x_{n-1}\right)  \\
& \leq \left\Vert x_{n}\right\Vert ^{2}\Gamma \left( x_{1},\dots
,x_{n-1}\right) .
\end{align*}%
Multiplying the above inequalities, we deduce%
\begin{align}
& \Gamma \left( x_{1},\dots ,x_{n-1},x_{n}\right)   \label{upper2.9} \\
& <\left\Vert x_{1}\right\Vert ^{2}\prod_{k=2}^{n}\left( \left\Vert
x_{k}\right\Vert ^{2}-\frac{1}{\sum_{i=1}^{k-1}\left\Vert x_{i}\right\Vert
^{2}}\sum_{i=1}^{k-1}\left\vert \left\langle x_{k},x_{i}\right\rangle
\right\vert ^{2}\right)   \notag \\
& \leq \prod_{j=1}^{n}\left\Vert x_{j}\right\Vert ^{2},  \notag
\end{align}%
valid for a system of $n\geq 2$ linearly independent vectors which are not
orthogonal on each other.
\end{remark}

In \cite{DRA1xx}, the author has obtained the following inequality.

\begin{lemma}[Dragomir, 2004]
\label{upperl2.3}Let $z_{1},\dots ,z_{n}\in H$ and $\alpha _{1},\dots
,\alpha _{n}\in \mathbb{K}$. Then one has the inequalities:%
\begin{equation}
\left\Vert \sum_{i=1}^{n}\alpha _{i}z_{i}\right\Vert ^{2}  \label{upper2.10}
\end{equation}%
\begin{multline*}
\leq \left\{ 
\begin{array}{l}
\max\limits_{1\leq i\leq n}\left\vert \alpha _{i}\right\vert
^{2}\sum\limits_{i=1}^{n}\left\Vert z_{i}\right\Vert ^{2}; \\ 
\\ 
\left( \sum\limits_{i=1}^{n}\left\vert \alpha _{i}\right\vert ^{2\alpha
}\right) ^{\frac{1}{\alpha }}\left( \sum\limits_{i=1}^{n}\left\Vert
z_{i}\right\Vert ^{2\beta }\right) ^{\frac{1}{p}} \\ 
\hfill \text{where \ }\alpha >1,\ \frac{1}{\alpha }+\frac{1}{\beta }=1; \\ 
\\ 
\sum\limits_{i=1}^{n}\left\vert \alpha _{i}\right\vert
^{2}\max\limits_{1\leq i\leq n}\left\Vert z_{i}\right\Vert ^{2};%
\end{array}%
\right. \\
+\left\{ 
\begin{array}{l}
\max\limits_{1\leq i\neq j\leq n}\left\{ \left\vert \alpha _{i}\alpha
_{j}\right\vert \right\} \sum\limits_{1\leq i\neq j\leq n}\left\vert
\left\langle z_{i},z_{j}\right\rangle \right\vert ; \\ 
\\ 
\left[ \left( \sum\limits_{i=1}^{n}\left\vert \alpha _{i}\right\vert
^{\gamma }\right) ^{2}-\sum\limits_{i=1}^{n}\left\vert \alpha
_{i}\right\vert ^{2\gamma }\right] ^{\frac{1}{\gamma }}\left(
\sum\limits_{1\leq i\neq j\leq n}\left\vert \left\langle
z_{i},z_{j}\right\rangle \right\vert ^{\delta }\right) ^{\frac{1}{\delta }}
\\ 
\hfill \text{where \ }\gamma >1,\ \frac{1}{\gamma }+\frac{1}{\delta }=1; \\ 
\\ 
\left[ \left( \sum\limits_{i=1}^{n}\left\vert \alpha _{i}\right\vert \right)
^{2}-\sum\limits_{i=1}^{n}\left\vert \alpha _{i}\right\vert ^{2}\right]
\max\limits_{1\leq i\neq j\leq n}\left\vert \left\langle
z_{i},z_{j}\right\rangle \right\vert ;%
\end{array}%
\right.
\end{multline*}%
where any term in the first branch can be combined with each term from the
second branch giving 9 possible combinations.
\end{lemma}

Out of these, we select the following ones that are of relevance for further
consideration:%
\begin{align}
& \left\Vert \sum_{i=1}^{n}\alpha _{i}z_{i}\right\Vert ^{2}
\label{upper2.11} \\
& \leq \max\limits_{1\leq i\leq n}\left\Vert z_{i}\right\Vert
^{2}\sum\limits_{i=1}^{n}\left\vert \alpha _{i}\right\vert ^{2}  \notag \\
& \qquad \qquad +\max\limits_{1\leq i<j\leq n}\left\vert \left\langle
z_{i},z_{j}\right\rangle \right\vert \left[ \left(
\sum\limits_{i=1}^{n}\left\vert \alpha _{i}\right\vert \right)
^{2}-\sum\limits_{i=1}^{n}\left\vert \alpha _{i}\right\vert ^{2}\right] 
\notag \\
& \leq \sum\limits_{i=1}^{n}\left\vert \alpha _{i}\right\vert ^{2}\left(
\max\limits_{1\leq i\leq n}\left\Vert z_{i}\right\Vert ^{2}+\left(
n-1\right) \max\limits_{1\leq i<j\leq n}\left\vert \left\langle
z_{i},z_{j}\right\rangle \right\vert \right)  \notag
\end{align}%
and%
\begin{align}
& \left\Vert \sum_{i=1}^{n}\alpha _{i}z_{i}\right\Vert ^{2}
\label{upper2.12} \\
& \leq \max\limits_{1\leq i\leq n}\left\Vert z_{i}\right\Vert
^{2}\sum\limits_{i=1}^{n}\left\vert \alpha _{i}\right\vert ^{2}+\left[
\left( \sum\limits_{i=1}^{n}\left\vert \alpha _{i}\right\vert ^{2}\right)
^{2}-\sum\limits_{i=1}^{n}\left\vert \alpha _{i}\right\vert ^{4}\right]
^{1/2}  \notag \\
& \qquad \qquad \qquad \qquad \qquad \times \left( \sum\limits_{1\leq i\neq
j\leq n}\left\vert \left\langle z_{i},z_{j}\right\rangle \right\vert
^{2}\right) ^{\frac{1}{2}}  \notag \\
& \leq \sum\limits_{i=1}^{n}\left\vert \alpha _{i}\right\vert ^{2}\left[
\max\limits_{1\leq i\leq n}\left\Vert z_{i}\right\Vert ^{2}+\left(
\sum\limits_{1\leq i\neq j\leq n}\left\vert \left\langle
z_{i},z_{j}\right\rangle \right\vert ^{2}\right) ^{\frac{1}{2}}\right] . 
\notag
\end{align}%
Note that the last inequality in (\ref{upper2.11}) follows by the fact that%
\begin{equation*}
\left( \sum\limits_{i=1}^{n}\left\vert \alpha _{i}\right\vert \right)
^{2}\leq n\sum\limits_{i=1}^{n}\left\vert \alpha _{i}\right\vert ^{2},
\end{equation*}%
while the last inequality in (\ref{upper2.12}) is obvious.

Utilising the above inequalities (\ref{upper2.11}) and (\ref{upper2.12})
which provide alternatives to the Cauchy-Bunyakovsky-Schwarz inequality (\ref%
{upper2.3}), we can state the following results \cite{SILV1xx}.

\begin{theorem}[Dragomir, 2005]
\label{uppert2.4}Let $\left\{ x_{1},\dots ,x_{n}\right\} ,$ $M$ and $x$ be
as in Theorem \ref{uppert2.1}. Then%
\begin{multline}
d^{2}\left( x,M\right)   \label{upper2.13} \\
\leq \frac{\left\Vert x\right\Vert ^{2}\left[ \max\limits_{1\leq i\leq
n}\left\Vert x_{i}\right\Vert ^{2}+\left( \sum\limits_{1\leq i\neq j\leq
n}\left\vert \left\langle x_{i},x_{j}\right\rangle \right\vert ^{2}\right) ^{%
\frac{1}{2}}\right] -\sum\limits_{i=1}^{n}\left\vert \left\langle
x,x_{i}\right\rangle \right\vert ^{2}}{\max\limits_{1\leq i\leq n}\left\Vert
x_{i}\right\Vert ^{2}+\left( \sum\limits_{1\leq i\neq j\leq n}\left\vert
\left\langle x_{i},x_{j}\right\rangle \right\vert ^{2}\right) ^{\frac{1}{2}}}
\end{multline}%
or, equivalently,%
\begin{multline}
\Gamma \left( x_{1},\dots ,x_{n},x\right)   \label{upper2.14} \\
\leq \frac{\left\Vert x\right\Vert ^{2}\left[ \max\limits_{1\leq i\leq
n}\left\Vert x_{i}\right\Vert ^{2}+\left( \sum\limits_{1\leq i\neq j\leq
n}\left\vert \left\langle x_{i},x_{j}\right\rangle \right\vert ^{2}\right) ^{%
\frac{1}{2}}\right] -\sum\limits_{i=1}^{n}\left\vert \left\langle
x,x_{i}\right\rangle \right\vert ^{2}}{\max\limits_{1\leq i\leq n}\left\Vert
x_{i}\right\Vert ^{2}+\left( \sum\limits_{1\leq i\neq j\leq n}\left\vert
\left\langle x_{i},x_{j}\right\rangle \right\vert ^{2}\right) ^{\frac{1}{2}}}
\\
\times \Gamma \left( x_{1},\dots ,x_{n}\right) .
\end{multline}
\end{theorem}

\begin{proof}
Utilising the inequality (\ref{upper2.12}) for $\alpha _{i}=\left\langle
x,x_{i}\right\rangle $ and $z_{i}=x_{i},$ $i\in \left\{ 1,\dots ,n\right\} ,$
we can write:%
\begin{multline}
\left\Vert \sum_{i=1}^{n}\left\langle x,x_{i}\right\rangle x_{i}\right\Vert
^{2}  \label{upper2.15} \\
\leq \sum_{i=1}^{n}\left\vert \left\langle x,x_{i}\right\rangle \right\vert
^{2}\left[ \max\limits_{1\leq i\leq n}\left\Vert x_{i}\right\Vert
^{2}+\left( \sum_{1\leq i\neq j\leq n}\left\vert \left\langle
x_{i},x_{j}\right\rangle \right\vert ^{2}\right) ^{\frac{1}{2}}\right]
\end{multline}%
for any $x\in H.$

Now, since, by the representation formula (\ref{upper1.3})%
\begin{equation}
d^{2}\left( x,M\right) =\left\Vert x\right\Vert ^{2}-\frac{%
\sum_{i=1}^{n}\left\vert \left\langle x,x_{i}\right\rangle \right\vert ^{2}}{%
\left\Vert \sum_{i=1}^{n}\left\langle x,x_{i}\right\rangle x_{i}\right\Vert
^{2}}\cdot \sum_{i=1}^{n}\left\vert \left\langle x,x_{i}\right\rangle
\right\vert ^{2},  \label{upper2.16}
\end{equation}%
for $x\notin M^{\perp },$ hence, by (\ref{upper2.15}) and (\ref{upper2.16})
we deduce the desired result (\ref{upper2.13}).
\end{proof}

\begin{remark}
\label{upperr2.5}In 1941, R.P. Boas \cite{BOxx} and in 1944, R. Bellman \cite%
{BExx}, independent of each other, proved the following generalisation of
Bessel's inequality:%
\begin{equation}
\sum_{i=1}^{n}\left\vert \left\langle y,y_{i}\right\rangle \right\vert
^{2}\leq \left\Vert y\right\Vert ^{2}\left[ \max\limits_{1\leq i\leq
n}\left\Vert y_{i}\right\Vert ^{2}+\left( \sum_{1\leq i\neq j\leq
n}\left\vert \left\langle y_{i},y_{j}\right\rangle \right\vert ^{2}\right) ^{%
\frac{1}{2}}\right] ,  \label{upper2.17}
\end{equation}%
provided $y$ and $y_{i}$ $\left( i\in \left\{ 1,\dots ,n\right\} \right) $
are arbitrary vectors in the inner product space $\left( H;\left\langle
\cdot ,\cdot \right\rangle \right) .$ If $\left\{ y_{i}\right\} _{i\in
\left\{ 1,\dots ,n\right\} }$ are orthonormal, then (\ref{upper2.17})
reduces to Bessel's inequality.

In this respect, one may see (\ref{upper2.13}) as a refinement of the
Boas-Bellman result (\ref{upper2.17}).
\end{remark}

\begin{remark}
\label{upperr2.6}On making use of a similar argument to that utilised in
Remark \ref{upperr2.2}, one can obtain the following refinement of the
Hadamard inequality:%
\begin{align}
& \Gamma \left( x_{1},\dots ,x_{n}\right)  \label{upper2.18} \\
& \leq \left\Vert x_{1}\right\Vert ^{2}  \notag \\
& \times \prod_{k=2}^{n}\left( \left\Vert x_{k}\right\Vert ^{2}-\frac{%
\sum\limits_{i=1}^{k-1}\left\vert \left\langle x_{k},x_{i}\right\rangle
\right\vert ^{2}}{\max\limits_{1\leq i\leq k-1}\left\Vert x_{i}\right\Vert
^{2}+\left( \sum\limits_{1\leq i\neq j\leq k-1}\left\vert \left\langle
x_{i},x_{j}\right\rangle \right\vert ^{2}\right) ^{\frac{1}{2}}}\right) 
\notag \\
& \leq \prod_{j=1}^{n}\left\Vert x_{j}\right\Vert ^{2}.  \notag
\end{align}
\end{remark}

Further on, if we choose $\alpha _{i}=\left\langle x,x_{i}\right\rangle ,$ $%
z_{i}=x_{i},$ $i\in \left\{ 1,\dots ,n\right\} $ in (\ref{upper2.11}), then
we may state the inequality%
\begin{multline}
\left\Vert \sum_{i=1}^{n}\left\langle x,x_{i}\right\rangle x_{i}\right\Vert
^{2}  \label{upper2.19} \\
\leq \sum_{i=1}^{n}\left\vert \left\langle x,x_{i}\right\rangle \right\vert
^{2}\left( \max\limits_{1\leq i\leq n}\left\Vert x_{i}\right\Vert
^{2}+\left( n-1\right) \max_{1\leq i\neq j\leq n}\left\vert \left\langle
x_{i},x_{j}\right\rangle \right\vert \right) .
\end{multline}%
Utilising (\ref{upper2.19}) and (\ref{upper2.16}) we may state the following
result as well \cite{SILV1xx}:

\begin{theorem}[Dragomir, 2005]
\label{uppert2.7}Let $\left\{ x_{1},\dots ,x_{n}\right\} ,$ $M$ and $x$ be
as in Theorem \ref{uppert2.1}. Then%
\begin{multline}
d^{2}\left( x,M\right)   \label{upper2.20} \\
\leq \frac{\left\Vert x\right\Vert ^{2}\left[ \max\limits_{1\leq i\leq
n}\left\Vert x_{i}\right\Vert ^{2}+\left( n-1\right) \max\limits_{1\leq
i\neq j\leq n}\left\vert \left\langle x_{i},x_{j}\right\rangle \right\vert %
\right] -\sum_{i=1}^{n}\left\vert \left\langle x,x_{i}\right\rangle
\right\vert ^{2}}{\max\limits_{1\leq i\leq n}\left\Vert x_{i}\right\Vert
^{2}+\left( n-1\right) \max\limits_{1\leq i\neq j\leq n}\left\vert
\left\langle x_{i},x_{j}\right\rangle \right\vert }
\end{multline}%
or, equivalently,%
\begin{multline}
\Gamma \left( x_{1},\dots ,x_{n},x\right)   \label{upper2.21} \\
\leq \frac{\left\Vert x\right\Vert ^{2}\left[ \max\limits_{1\leq i\leq
n}\left\Vert x_{i}\right\Vert ^{2}+\left( n-1\right) \max\limits_{1\leq
i\neq j\leq n}\left\vert \left\langle x_{i},x_{j}\right\rangle \right\vert %
\right] -\sum_{i=1}^{n}\left\vert \left\langle x,x_{i}\right\rangle
\right\vert ^{2}}{\max\limits_{1\leq i\leq n}\left\Vert x_{i}\right\Vert
^{2}+\left( n-1\right) \max\limits_{1\leq i\neq j\leq n}\left\vert
\left\langle x_{i},x_{j}\right\rangle \right\vert } \\
\times \Gamma \left( x_{1},\dots ,x_{n}\right) .
\end{multline}
\end{theorem}

\begin{remark}
\label{upperr2.8}The above result (\ref{upper2.20}) provides a refinement
for the following generalisation of Bessel's inequality:%
\begin{equation}
\sum_{i=1}^{n}\left\vert \left\langle x,x_{i}\right\rangle \right\vert
^{2}\leq \left\Vert x\right\Vert ^{2}\left[ \max\limits_{1\leq i\leq
n}\left\Vert x_{i}\right\Vert ^{2}+\left( n-1\right) \max\limits_{1\leq
i\neq j\leq n}\left\vert \left\langle x_{i},x_{j}\right\rangle \right\vert %
\right] ,  \label{upper2.22}
\end{equation}%
obtained by the author in \cite{DRA1xx}.

One can also provide the corresponding refinement of Hadamard's inequality (%
\ref{upper1.4}) on using (\ref{upper2.21}), i.e., 
\begin{align}
& \Gamma \left( x_{1},\dots ,x_{n}\right)  \label{upper2.23} \\
& \leq \left\Vert x_{1}\right\Vert ^{2}  \notag \\
& \times \prod_{k=2}^{n}\left( \left\Vert x_{k}\right\Vert ^{2}-\frac{%
\sum\limits_{i=1}^{k-1}\left\vert \left\langle x_{k},x_{i}\right\rangle
\right\vert ^{2}}{\max\limits_{1\leq i\leq k-1}\left\Vert x_{i}\right\Vert
^{2}+\left( k-2\right) \max\limits_{1\leq i\neq j\leq k-1}\left\vert
\left\langle x_{i},x_{j}\right\rangle \right\vert }\right)  \notag \\
& \leq \prod_{j=1}^{n}\left\Vert x_{j}\right\Vert ^{2}.  \notag
\end{align}
\end{remark}

\subsection{Other Upper Bounds for $d\left( x,M\right) $}

In \cite[p. 140]{SSDxx} the author obtained the following inequality that is
similar to the Cauchy-Bunyakovsky-Schwarz result.

\begin{lemma}[Dragomir, 2004]
\label{upperl3.1}Let $z_{1},\dots ,z_{n}\in H$ and $\alpha _{1},\dots
,\alpha _{n}\in \mathbb{K}$. Then one has the inequalities:%
\begin{align}
\left\Vert \sum_{i=1}^{n}\alpha _{i}z_{i}\right\Vert ^{2}& \leq
\sum\limits_{i=1}^{n}\left\vert \alpha _{i}\right\vert
^{2}\sum\limits_{j=1}^{n}\left\vert \left\langle z_{i},z_{j}\right\rangle
\right\vert  \label{upper3.1} \\
& \leq \left\{ 
\begin{array}{l}
\sum\limits_{i=1}^{n}\left\vert \alpha _{i}\right\vert
^{2}\max\limits_{1\leq i\leq n}\left[ \sum\limits_{j=1}^{n}\left\vert
\left\langle z_{i},z_{j}\right\rangle \right\vert \right] ; \\ 
\\ 
\left( \sum\limits_{i=1}^{n}\left\vert \alpha _{i}\right\vert ^{2p}\right) ^{%
\frac{1}{p}}\left( \sum\limits_{i=1}^{n}\left(
\sum\limits_{j=1}^{n}\left\vert \left\langle z_{i},z_{j}\right\rangle
\right\vert \right) ^{q}\right) ^{\frac{1}{q}} \\ 
\hfill \text{where \ }p>1,\ \frac{1}{p}+\frac{1}{q}=1; \\ 
\\ 
\max\limits_{1\leq i\leq n}\left\vert \alpha _{i}\right\vert
^{2}\sum\limits_{i,j=1}^{n}\left\vert \left\langle z_{i},z_{j}\right\rangle
\right\vert .%
\end{array}%
\right.  \notag
\end{align}
\end{lemma}

We can state and prove now another upper bound for the distance $d\left(
x,M\right) $ as follows \cite{SILV1xx}.

\begin{theorem}[Dragomir, 2005]
\label{uppert3.2}Let $\left\{ x_{1},\dots ,x_{n}\right\} ,$ $M$ and $x$ be
as in Theorem \ref{uppert2.1}. Then%
\begin{equation}
d^{2}\left( x,M\right) \leq \frac{\left\Vert x\right\Vert
^{2}\max\limits_{1\leq i\leq n}\left[ \sum\limits_{j=1}^{n}\left\vert
\left\langle x_{i},x_{j}\right\rangle \right\vert \right] -\sum%
\limits_{i=1}^{n}\left\vert \left\langle x,x_{i}\right\rangle \right\vert
^{2}}{\max\limits_{1\leq i\leq n}\left[ \sum\limits_{j=1}^{n}\left\vert
\left\langle x_{i},x_{j}\right\rangle \right\vert \right] }  \label{upper3.2}
\end{equation}%
or, equivalently,%
\begin{multline}
\Gamma \left( x_{1},\dots ,x_{n},x\right)  \label{upper3.3} \\
\leq \frac{\left\Vert x\right\Vert ^{2}\max\limits_{1\leq i\leq n}\left[
\sum\limits_{j=1}^{n}\left\vert \left\langle x_{i},x_{j}\right\rangle
\right\vert \right] -\sum\limits_{i=1}^{n}\left\vert \left\langle
x,x_{i}\right\rangle \right\vert ^{2}}{\max\limits_{1\leq i\leq n}\left[
\sum\limits_{j=1}^{n}\left\vert \left\langle x_{i},x_{j}\right\rangle
\right\vert \right] }\cdot \Gamma \left( x_{1},\dots ,x_{n}\right) .
\end{multline}
\end{theorem}

\begin{proof}
Utilising the first branch in (\ref{upper3.1}) we may state that%
\begin{equation}
\left\Vert \sum_{i=1}^{n}\left\langle x,x_{i}\right\rangle x_{i}\right\Vert
^{2}\leq \sum_{i=1}^{n}\left\vert \left\langle x,x_{i}\right\rangle
\right\vert ^{2}\max\limits_{1\leq i\leq n}\left[ \sum\limits_{j=1}^{n}\left%
\vert \left\langle x_{i},x_{j}\right\rangle \right\vert \right]
\label{upper3.3a}
\end{equation}%
for any $x\in H.$

Now, since, by the representation formula (\ref{upper1.3}) we have%
\begin{equation}
d^{2}\left( x,M\right) =\left\Vert x\right\Vert ^{2}-\frac{%
\sum_{i=1}^{n}\left\vert \left\langle x,x_{i}\right\rangle \right\vert ^{2}}{%
\left\Vert \sum_{i=1}^{n}\left\langle x,x_{i}\right\rangle x_{i}\right\Vert
^{2}}\cdot \sum_{i=1}^{n}\left\vert \left\langle x,x_{i}\right\rangle
\right\vert ^{2},  \label{upper3.4}
\end{equation}%
for $x\notin M^{\perp },$ hence, by (\ref{upper3.3a}) and (\ref{upper3.4})
we deduce the desired result (\ref{upper3.2}).
\end{proof}

\begin{remark}
\label{upperr3.3}In 1971, E. Bombieri \cite{BOMxx} proved the following
generalisation of Bessel's inequality, however not stated in the general
form for inner products. The general version can be found for instance in 
\cite[p. 394]{MPFxx}. It reads as follows: if $y,y_{1},\dots ,y_{n}$ are
vectors in the inner product space $\left( H;\left\langle \cdot ,\cdot
\right\rangle \right) ,$ then%
\begin{equation}
\sum_{i=1}^{n}\left\vert \left\langle y,y_{i}\right\rangle \right\vert
^{2}\leq \left\Vert y\right\Vert ^{2}\max\limits_{1\leq i\leq n}\left\{
\sum\limits_{j=1}^{n}\left\vert \left\langle y_{i},y_{j}\right\rangle
\right\vert \right\} .  \label{upper3.5}
\end{equation}%
Obviously, when $\left\{ y_{1},\dots ,y_{n}\right\} $ are orthonormal, the
inequality (\ref{upper3.5}) produces Bessel's inequality.

In this respect, we may regard our result (\ref{upper3.2}) as a refinement
of the Bombieri inequality (\ref{upper3.5}).
\end{remark}

\begin{remark}
\label{upperr3.4}On making use of a similar argument to that in Remark \ref%
{upperr2.2}, we obtain the following refinement for the Hadamard inequality:%
\begin{align}
\Gamma \left( x_{1},\dots ,x_{n}\right) & \leq \left\Vert x_{1}\right\Vert
^{2}\prod_{k=2}^{n}\left[ \left\Vert x_{k}\right\Vert ^{2}-\frac{%
\sum\limits_{i=1}^{k-1}\left\vert \left\langle x_{k},x_{i}\right\rangle
\right\vert ^{2}}{\max\limits_{1\leq i\leq k-1}\left[ \sum%
\limits_{j=1}^{k-1}\left\vert \left\langle x_{i},x_{j}\right\rangle
\right\vert \right] }\right]  \label{upper3.6} \\
& \leq \prod_{j=1}^{n}\left\Vert x_{j}\right\Vert ^{2}.  \notag
\end{align}
\end{remark}

Another different Cauchy-Bunyakovsky-Schwarz type inequality is incorporated
in the following lemma \cite{DRA5xx}.

\begin{lemma}[Dragomir, 2004]
\label{upperl3.4}Let $z_{1},\dots ,z_{n}\in H$ and $\alpha _{1},\dots
,\alpha _{n}\in \mathbb{K}$. Then 
\begin{equation}
\left\Vert \sum_{i=1}^{n}\alpha _{i}z_{i}\right\Vert ^{2}\leq \left(
\sum\limits_{i=1}^{n}\left\vert \alpha _{i}\right\vert ^{p}\right) ^{\frac{2%
}{p}}\left( \sum\limits_{i,j=1}^{n}\left\vert \left\langle
z_{i},z_{j}\right\rangle \right\vert ^{q}\right) ^{\frac{1}{q}}
\label{upper3.7}
\end{equation}%
for $p>1,$ $\frac{1}{p}+\frac{1}{q}=1.$

If in (\ref{upper3.7}) we choose $p=q=2,$ then we get%
\begin{equation}
\left\Vert \sum_{i=1}^{n}\alpha _{i}z_{i}\right\Vert ^{2}\leq
\sum\limits_{i=1}^{n}\left\vert \alpha _{i}\right\vert ^{2}\left(
\sum\limits_{i,j=1}^{n}\left\vert \left\langle z_{i},z_{j}\right\rangle
\right\vert ^{2}\right) ^{\frac{1}{2}}.  \label{upper3.8}
\end{equation}
\end{lemma}

Based on (\ref{upper3.8}), we can state the following result that provides
yet another upper bound for the distance $d\left( x,M\right) $ \cite{SILV1xx}%
.

\begin{theorem}[Dragomir, 2005]
\label{uppert3.5}Let $\left\{ x_{1},\dots ,x_{n}\right\} ,$ $M$ and $x$ be
as in Theorem \ref{uppert2.1}. Then%
\begin{equation}
d^{2}\left( x,M\right) \leq \frac{\left\Vert x\right\Vert ^{2}\left(
\sum\limits_{i,j=1}^{n}\left\vert \left\langle x_{i},x_{j}\right\rangle
\right\vert ^{2}\right) ^{\frac{1}{2}}-\sum\limits_{i=1}^{n}\left\vert
\left\langle x,x_{i}\right\rangle \right\vert ^{2}}{\left(
\sum\limits_{i,j=1}^{n}\left\vert \left\langle x_{i},x_{j}\right\rangle
\right\vert ^{2}\right) ^{\frac{1}{2}}}  \label{upper3.9}
\end{equation}%
or, equivalently,%
\begin{multline}
\Gamma \left( x_{1},\dots ,x_{n},x\right)  \label{upper3.10} \\
\leq \frac{\left\Vert x\right\Vert ^{2}\left(
\sum\limits_{i,j=1}^{n}\left\vert \left\langle x_{i},x_{j}\right\rangle
\right\vert ^{2}\right) ^{\frac{1}{2}}-\sum\limits_{i=1}^{n}\left\vert
\left\langle x,x_{i}\right\rangle \right\vert ^{2}}{\left(
\sum\limits_{i,j=1}^{n}\left\vert \left\langle x_{i},x_{j}\right\rangle
\right\vert ^{2}\right) ^{\frac{1}{2}}}\cdot \Gamma \left( x_{1},\dots
,x_{n}\right) .
\end{multline}
\end{theorem}

Similar comments apply related to Hadamard's inequality. We omit the details.

\subsection{Some Conditional Bounds}

In the recent paper \cite{DRA4xx}, the author has established the following
reverse of the Bessel inequality.

Let $\left( H;\left\langle \cdot ,\cdot \right\rangle \right) $ be an inner
product space over the real or complex number field $\mathbb{K}$, $\left\{
e_{i}\right\} _{i\in I}$ a finite family of orthonormal vectors in $H,$ $%
\varphi _{i},\phi _{i}\in \mathbb{K}$, $i\in I$ and $x\in H.$ If%
\begin{equation}
\func{Re}\left\langle \sum_{i\in I}\phi _{i}e_{i}-x,x-\sum_{i\in I}\varphi
_{i}e_{i}\right\rangle \geq 0  \label{upper4.1}
\end{equation}%
or, equivalently,%
\begin{equation}
\left\Vert x-\sum_{i\in I}\frac{\varphi _{i}+\phi _{i}}{2}e_{i}\right\Vert
\leq \frac{1}{2}\left( \sum_{i\in I}\left\vert \phi _{i}-\varphi
_{i}\right\vert ^{2}\right) ^{\frac{1}{2}},  \label{upper4.2}
\end{equation}%
then%
\begin{equation}
\left( 0\leq \right) \left\Vert x\right\Vert ^{2}-\sum_{i\in I}\left\vert
\left\langle x,e_{i}\right\rangle \right\vert ^{2}\leq \frac{1}{4}\sum_{i\in
I}\left\vert \phi _{i}-\varphi _{i}\right\vert ^{2}.  \label{upper4.3}
\end{equation}%
The constant $\frac{1}{4}$ is best possible in the sense that it cannot be
replaced by a smaller constant \cite{SILV1xx}.

\begin{theorem}[Dragomir, 2005]
\label{uppert4.1}Let $\left\{ x_{1},\dots x_{n}\right\} $ be a linearly
independent system of vectors in $H$ and $M:=span\left\{ x_{1},\dots
x_{n}\right\} .$ If $\gamma _{i},$ $\Gamma _{i}\in \mathbb{K}$, $i\in
\left\{ 1,\dots ,n\right\} $ and $x\in H\backslash M^{\perp }$ is such that%
\begin{equation}
\func{Re}\left\langle \sum_{i=1}^{n}\Gamma
_{i}x_{i}-x,x-\sum_{i=1}^{n}\gamma _{i}x_{i}\right\rangle \geq 0,
\label{upper4.4}
\end{equation}%
then we have the bound%
\begin{equation}
d^{2}\left( x,M\right) \leq \frac{1}{4}\left\Vert \sum_{i=1}^{n}\left(
\Gamma _{i}-\gamma _{i}\right) x_{i}\right\Vert ^{2}  \label{upper4.5}
\end{equation}%
or, equivalently,%
\begin{equation}
\Gamma \left( x_{1},\dots ,x_{n},x\right) \leq \frac{1}{4}\left\Vert
\sum_{i=1}^{n}\left( \Gamma _{i}-\gamma _{i}\right) x_{i}\right\Vert
^{2}\Gamma \left( x_{1},\dots ,x_{n}\right) .  \label{upper4.6}
\end{equation}
\end{theorem}

\begin{proof}
It is easy to see that in an inner product space for any $x,z,Z\in H$ one has%
\begin{equation*}
\left\Vert x-\frac{z+Z}{2}\right\Vert ^{2}-\frac{1}{4}\left\Vert
Z-z\right\Vert ^{2}=\func{Re}\left\langle Z-x,x-z\right\rangle ,
\end{equation*}%
therefore, the condition (\ref{upper4.4}) is actually equivalent to%
\begin{equation}
\left\Vert x-\sum_{i=1}^{n}\frac{\Gamma _{i}+\gamma _{i}}{2}x_{i}\right\Vert
^{2}\leq \frac{1}{4}\left\Vert \sum_{i=1}^{n}\left( \Gamma _{i}-\gamma
_{i}\right) x_{i}\right\Vert ^{2}.  \label{upper4.7}
\end{equation}%
Now, obviously,%
\begin{equation}
d^{2}\left( x,M\right) =\inf_{y\in M}\left\Vert x-y\right\Vert ^{2}\leq
\left\Vert x-\sum_{i=1}^{n}\frac{\Gamma _{i}+\gamma _{i}}{2}x_{i}\right\Vert
^{2}  \label{upper4.8}
\end{equation}%
and thus, by (\ref{upper4.7}) and (\ref{upper4.8}) we deduce (\ref{upper4.5}%
).

The last inequality is obvious by the representation (\ref{upper1.2}).
\end{proof}

\begin{remark}
\label{upperr4.2}Utilising various Cauchy-Bunyakovsky-Schwarz type
inequalities we may obtain more convenient (although coarser) bounds for $%
d^{2}\left( x,M\right) .$ For instance, if we use the inequality (\ref%
{upper2.11}) we can state the inequality:%
\begin{multline*}
\left\Vert \sum_{i=1}^{n}\left( \Gamma _{i}-\gamma _{i}\right)
x_{i}\right\Vert ^{2} \\
\leq \sum_{i=1}^{n}\left\vert \Gamma _{i}-\gamma _{i}\right\vert ^{2}\left(
\max\limits_{1\leq i\leq n}\left\Vert x_{i}\right\Vert ^{2}+\left(
n-1\right) \max\limits_{1\leq i<j\leq n}\left\vert \left\langle
x_{i},x_{j}\right\rangle \right\vert \right) ,
\end{multline*}%
giving the bound:%
\begin{multline}
d^{2}\left( x,M\right) \leq \frac{1}{4}\sum_{i=1}^{n}\left\vert \Gamma
_{i}-\gamma _{i}\right\vert ^{2}  \label{upper4.9} \\
\times \left[ \max\limits_{1\leq i\leq n}\left\Vert x_{i}\right\Vert
^{2}+\left( n-1\right) \max\limits_{1\leq i<j\leq n}\left\vert \left\langle
x_{i},x_{j}\right\rangle \right\vert \right] ,
\end{multline}%
provided (\ref{upper4.4}) holds true.

Obviously, if $\left\{ x_{1},\dots ,x_{n}\right\} $ is an orthonormal family
in $H,$ then from (\ref{upper4.9}) we deduce the reverse of Bessel's
inequality incorporated in (\ref{upper4.3}).

If we use the inequality (\ref{upper2.12}), then we can state the inequality%
\begin{multline*}
\left\Vert \sum_{i=1}^{n}\left( \Gamma _{i}-\gamma _{i}\right)
x_{i}\right\Vert ^{2} \\
\leq \sum_{i=1}^{n}\left\vert \Gamma _{i}-\gamma _{i}\right\vert ^{2}\left[
\max\limits_{1\leq i\leq n}\left\Vert x_{i}\right\Vert ^{2}+\left(
\sum\limits_{1\leq i\neq j\leq n}\left\vert \left\langle
x_{i},x_{j}\right\rangle \right\vert ^{2}\right) ^{\frac{1}{2}}\right] ,
\end{multline*}%
giving the bound%
\begin{multline}
d^{2}\left( x,M\right) \leq \frac{1}{4}\sum_{i=1}^{n}\left\vert \Gamma
_{i}-\gamma _{i}\right\vert ^{2}  \label{upper4.10} \\
\times \left[ \max\limits_{1\leq i\leq n}\left\Vert x_{i}\right\Vert
^{2}+\left( \sum\limits_{1\leq i\neq j\leq n}\left\vert \left\langle
x_{i},x_{j}\right\rangle \right\vert ^{2}\right) ^{\frac{1}{2}}\right] ,
\end{multline}%
provided (\ref{upper4.4}) holds true.

In this case, when one assumes that $\left\{ x_{1},\dots ,x_{n}\right\} $ is
an orthonormal family of vectors, then (\ref{upper4.10}) reduces to (\ref%
{upper4.3}) as well.

Finally, on utilising the first branch of the inequality (\ref{upper3.1}),
we can state that%
\begin{equation}
d^{2}\left( x,M\right) \leq \frac{1}{4}\sum_{i=1}^{n}\left\vert \Gamma
_{i}-\gamma _{i}\right\vert ^{2}\max\limits_{1\leq i\leq n}\left[
\sum_{j=1}^{n}\left\vert \left\langle x_{i},x_{j}\right\rangle \right\vert %
\right] ,  \label{upper4.11}
\end{equation}%
provided (\ref{upper4.4}) holds true.

This inequality is also a generalisation of (\ref{upper4.3}).
\end{remark}

\section{Reversing the CBS Inequality for Sequences}

\subsection{Introduction}

Let $\left( H;\left\langle \cdot ,\cdot \right\rangle \right) $ be an inner
product space over the real or complex number field $\mathbb{K}$. One of the
most important inequalities in inner product spaces with numerous
applications, is the \textit{Schwarz inequality} 
\begin{equation}
\left\vert \left\langle x,y\right\rangle \right\vert ^{2}\leq \left\Vert
x\right\Vert ^{2}\left\Vert y\right\Vert ^{2},\ \ \ \ x,y\in H
\label{revcbs1.1}
\end{equation}%
or, equivalently, 
\begin{equation}
\left\vert \left\langle x,y\right\rangle \right\vert \leq \left\Vert
x\right\Vert \left\Vert y\right\Vert ,\ \ \ \ x,y\in H.  \label{revcbs1.2}
\end{equation}%
The case of equality holds iff there exists a scalar $\alpha \in \mathbb{K}$
such that $x=\alpha y.$

By a \textit{multiplicative reverse} of the Schwarz inequality we understand
an inequality of the form%
\begin{equation}
\left( 1\leq \right) \frac{\left\Vert x\right\Vert \left\Vert y\right\Vert }{%
\left\vert \left\langle x,y\right\rangle \right\vert }\leq k_{1}\text{ \ or
\ }\left( 1\leq \right) \frac{\left\Vert x\right\Vert ^{2}\left\Vert
y\right\Vert ^{2}}{\left\vert \left\langle x,y\right\rangle \right\vert ^{2}}%
\leq k_{2}  \label{revcbs1.3}
\end{equation}%
with appropriate $k_{1}$ and $k_{2}$ and under various assumptions for the
vectors $x$ and $y,$ while by an \textit{additive reverse} we understand an
inequality of the form%
\begin{align}
(0& \leq )\left\Vert x\right\Vert \left\Vert y\right\Vert -\left\vert
\left\langle x,y\right\rangle \right\vert \leq h_{1}\text{ \ or \ }
\label{revcbs1.4} \\
(0& \leq )\left\Vert x\right\Vert ^{2}\left\Vert y\right\Vert
^{2}-\left\vert \left\langle x,y\right\rangle \right\vert ^{2}\leq h_{2}. 
\notag
\end{align}

Similar definition apply when $\left\vert \left\langle x,y\right\rangle
\right\vert $ is replaced by $\func{Re}\left\langle x,y\right\rangle $ or $%
\left\vert \func{Re}\left\langle x,y\right\rangle \right\vert .$

The following recent reverses for the Schwarz inequality hold (see for
instance the monograph on line \cite[p. 20]{SSDxx}).

\begin{theorem}[Dragomir, 2004]
\label{revcbst1.1}Let $\left( H;\left\langle \cdot ,\cdot \right\rangle
\right) $ be an inner product space over the real or complex number field $%
\mathbb{K}$. If $x,y\in H$ and $r>0$ are such that 
\begin{equation}
\left\Vert x-y\right\Vert \leq r<\left\Vert y\right\Vert ,  \label{revcbs1.5}
\end{equation}%
then we have the following multiplicative reverse of the Schwarz inequality%
\begin{equation}
\left( 1\leq \right) \frac{\left\Vert x\right\Vert \left\Vert y\right\Vert }{%
\left\vert \left\langle x,y\right\rangle \right\vert }\leq \frac{\left\Vert
x\right\Vert \left\Vert y\right\Vert }{\func{Re}\left\langle
x,y\right\rangle }\leq \frac{\left\Vert y\right\Vert }{\sqrt{\left\Vert
y\right\Vert ^{2}-r^{2}}}  \label{revcbs1.6}
\end{equation}%
and the subsequent additive reverses%
\begin{align}
(0& \leq )\left\Vert x\right\Vert \left\Vert y\right\Vert -\left\vert
\left\langle x,y\right\rangle \right\vert \leq \left\Vert x\right\Vert
\left\Vert y\right\Vert -\func{Re}\left\langle x,y\right\rangle
\label{revcbs1.7} \\
& \leq \frac{r^{2}}{\sqrt{\left\Vert y\right\Vert ^{2}-r^{2}}\left(
\left\Vert y\right\Vert +\sqrt{\left\Vert y\right\Vert ^{2}-r^{2}}\right) }%
\func{Re}\left\langle x,y\right\rangle  \notag
\end{align}%
and%
\begin{align}
(0& \leq )\left\Vert x\right\Vert ^{2}\left\Vert y\right\Vert
^{2}-\left\vert \left\langle x,y\right\rangle \right\vert ^{2}
\label{revcbs1.8} \\
& \leq \left\Vert x\right\Vert ^{2}\left\Vert y\right\Vert ^{2}-\left[ \func{%
Re}\left\langle x,y\right\rangle \right] ^{2}  \notag \\
& \leq r^{2}\left\Vert x\right\Vert ^{2}.  \notag
\end{align}%
All the above inequalities are sharp.
\end{theorem}

Other additive reverses of the quadratic Schwarz's inequality are
incorporated in the following result \cite[p. 18-19]{SSDxx}.

\begin{theorem}[Dragomir, 2004]
\label{revcbst1.2}Let $x,y\in H$ and $a,A\in \mathbb{K}$. If%
\begin{equation}
\func{Re}\left\langle Ay-x,x-ay\right\rangle \geq 0  \label{revcbs1.9}
\end{equation}%
or, equivalently,%
\begin{equation}
\left\Vert x-\frac{a+A}{2}\cdot y\right\Vert \leq \frac{1}{2}\left\vert
A-a\right\vert \left\Vert y\right\Vert ,  \label{revcbs1.10}
\end{equation}%
then%
\begin{align}
(0& \leq )\left\Vert x\right\Vert ^{2}\left\Vert y\right\Vert
^{2}-\left\vert \left\langle x,y\right\rangle \right\vert ^{2}
\label{revcbs1.11} \\
& \leq \frac{1}{4}\left\vert A-a\right\vert ^{2}\left\Vert y\right\Vert
^{4}-\left\{ 
\begin{array}{l}
\left\vert \frac{A+a}{2}\left\Vert y\right\Vert ^{2}-\left\langle
x,y\right\rangle \right\vert ^{2} \\ 
\\ 
\left\Vert y\right\Vert ^{2}\func{Re}\left\langle Ay-x,x-ay\right\rangle%
\end{array}%
\right.  \notag \\
& \leq \frac{1}{4}\left\vert A-a\right\vert ^{2}\left\Vert y\right\Vert ^{4}.
\notag
\end{align}%
The constant $\frac{1}{4}$ is best possible in all inequalities.
\end{theorem}

If one were to assume more about the complex numbers $A$ and $a,$ then one
may state the following result as well \cite[p. 21-23]{SSDxx}.

\begin{theorem}[Dragomir, 2004]
\label{revcbst1.3}With the assumptions of Theorem \ref{revcbst1.2} and, if
in addition, $\func{Re}\left( A\bar{a}\right) >0,$ then%
\begin{equation}
\left\Vert x\right\Vert \left\Vert y\right\Vert \leq \frac{1}{2}\cdot \frac{%
\func{Re}\left[ \left( \bar{A}+\bar{a}\right) \left\langle x,y\right\rangle %
\right] }{\sqrt{\func{Re}\left( A\bar{a}\right) }}\leq \frac{1}{2}\cdot 
\frac{\left\vert A+a\right\vert }{\sqrt{\func{Re}\left( A\bar{a}\right) }}%
\left\vert \left\langle x,y\right\rangle \right\vert ,  \label{revcbs1.12}
\end{equation}%
\begin{align}
(0& \leq )\left\Vert x\right\Vert \left\Vert y\right\Vert -\func{Re}%
\left\langle x,y\right\rangle  \label{revcbs1.13} \\
& \leq \frac{1}{2}\cdot \frac{\func{Re}\left[ \left( \bar{A}+\bar{a}-2\sqrt{%
\func{Re}\left( A\bar{a}\right) }\right) \left\langle x,y\right\rangle %
\right] }{\sqrt{\func{Re}\left( A\bar{a}\right) }}  \notag
\end{align}%
and%
\begin{equation}
(0\leq )\left\Vert x\right\Vert ^{2}\left\Vert y\right\Vert ^{2}-\left\vert
\left\langle x,y\right\rangle \right\vert ^{2}\leq \frac{1}{4}\cdot \frac{%
\left\vert A-a\right\vert ^{2}}{\func{Re}\left( A\bar{a}\right) }\left\vert
\left\langle x,y\right\rangle \right\vert ^{2}.  \label{revcbs1.14}
\end{equation}%
The constants $\frac{1}{2}$ and $\frac{1}{4}$ are best possible.
\end{theorem}

\begin{remark}
\label{revcbsr1.4}If $A=M,$ $a=m$ and $M\geq m>0,$ then (\ref{revcbs1.12})
and (\ref{revcbs1.13}) may be written in a more convenient form as%
\begin{equation}
\left\Vert x\right\Vert \left\Vert y\right\Vert \leq \frac{M+m}{2\sqrt{mM}}%
\func{Re}\left\langle x,y\right\rangle  \label{revcbs1.15}
\end{equation}%
and%
\begin{equation}
(0\leq )\left\Vert x\right\Vert \left\Vert y\right\Vert -\func{Re}%
\left\langle x,y\right\rangle \leq \frac{\left( \sqrt{M}-\sqrt{m}\right) ^{2}%
}{2\sqrt{mM}}\func{Re}\left\langle x,y\right\rangle .  \label{revcbs1.16}
\end{equation}%
Here the constant $\frac{1}{2}$ is sharp in both inequalities.
\end{remark}

In this section several reverses for the Cauchy-Bunyakovsky-Schwarz (CBS)
inequality for sequences of vectors in Hilbert spaces are obtained.
Applications for bounding the distance to a finite-dimensional subspace and
in reversing the generalised triangle inequality are also given.

\subsection{Reverses of the $\left( CBS\right) -$Inequality for Two
Sequences in $\ell _{\mathbf{p}}^{2}\left( K\right) $}

Let $\left( K,\left\langle \cdot ,\cdot \right\rangle \right) $ be a Hilbert
space over $\mathbb{K}$, $p_{i}\geq 0,$ $i\in \mathbb{N}$ with $%
\sum_{i=1}^{\infty }p_{i}=1.$ Consider $\ell _{\mathbf{p}}^{2}\left(
K\right) $ as the space%
\begin{equation*}
\ell _{\mathbf{p}}^{2}\left( K\right) :=\left\{ x=\left( x_{i}\right) _{i\in 
\mathbb{N}}\left\vert x_{i}\in K,\ i\in \mathbb{N}\text{ \ and \ }%
\sum_{i=1}^{\infty }p_{i}\left\Vert x_{i}\right\Vert ^{2}<\infty \right.
\right\} .
\end{equation*}%
It is well known that $\ell _{\mathbf{p}}^{2}\left( K\right) $ endowed with
the inner product%
\begin{equation*}
\left\langle x,y\right\rangle _{\mathbf{p}}:=\sum_{i=1}^{\infty
}p_{i}\left\langle x_{i},y_{i}\right\rangle
\end{equation*}%
is a Hilbert space over $\mathbb{K}$. The norm $\left\Vert \cdot \right\Vert
_{\mathbf{p}}$ of $\ell _{\mathbf{p}}^{2}\left( K\right) $ is given by%
\begin{equation*}
\left\Vert x\right\Vert _{\mathbf{p}}:=\left( \sum_{i=1}^{\infty
}p_{i}\left\Vert x_{i}\right\Vert ^{2}\right) ^{\frac{1}{2}}.
\end{equation*}%
If $x,y\in \ell _{\mathbf{p}}^{2}\left( K\right) ,$ then the following
Cauchy-Bunyakovsky-Schwarz $\left( CBS\right) $ inequality holds true%
\begin{equation}
\sum_{i=1}^{\infty }p_{i}\left\Vert x_{i}\right\Vert ^{2}\sum_{i=1}^{\infty
}p_{i}\left\Vert y_{i}\right\Vert ^{2}\geq \left\vert \sum_{i=1}^{\infty
}p_{i}\left\langle x_{i},y_{i}\right\rangle \right\vert ^{2}
\label{revcbs2.1}
\end{equation}%
with equality iff there exists a $\lambda \in \mathbb{K}$ such that $%
x_{i}=\lambda y_{i}$ for each $i\in \mathbb{N}$.

This is an obvious consequence of the Schwarz inequality (\ref{revcbs1.1})
written for the inner product $\left\langle \cdot ,\cdot \right\rangle _{%
\mathbf{p}}$ defined on $\ell _{\mathbf{p}}^{2}\left( K\right) .$

The following proposition may be stated \cite{DRA3xx}.

\begin{proposition}
\label{revcbsp2.1}Let $x,y\in \ell _{\mathbf{p}}^{2}\left( K\right) $ and $%
r>0.$ Assume that%
\begin{equation}
\left\Vert x_{i}-y_{i}\right\Vert \leq r<\left\Vert y_{i}\right\Vert \text{
\ for each \ }i\in \mathbb{N}\text{.}  \label{revcbs2.2}
\end{equation}%
Then we have the inequality%
\begin{align}
(1& \leq )\frac{\left( \sum_{i=1}^{\infty }p_{i}\left\Vert x_{i}\right\Vert
^{2}\sum_{i=1}^{\infty }p_{i}\left\Vert y_{i}\right\Vert ^{2}\right) ^{\frac{%
1}{2}}}{\left\vert \sum_{i=1}^{\infty }p_{i}\left\langle
x_{i},y_{i}\right\rangle \right\vert }  \label{revcbs2.3} \\
& \leq \frac{\left( \sum_{i=1}^{\infty }p_{i}\left\Vert x_{i}\right\Vert
^{2}\sum_{i=1}^{\infty }p_{i}\left\Vert y_{i}\right\Vert ^{2}\right) ^{\frac{%
1}{2}}}{\sum_{i=1}^{\infty }p_{i}\func{Re}\left\langle
x_{i},y_{i}\right\rangle }  \notag \\
& \leq \frac{\left( \sum_{i=1}^{\infty }p_{i}\left\Vert y_{i}\right\Vert
^{2}\right) ^{\frac{1}{2}}}{\sqrt{\sum_{i=1}^{\infty }p_{i}\left\Vert
y_{i}\right\Vert ^{2}-r^{2}}},  \notag
\end{align}%
\begin{align}
(0& \leq )\left( \sum_{i=1}^{\infty }p_{i}\left\Vert x_{i}\right\Vert
^{2}\sum_{i=1}^{\infty }p_{i}\left\Vert y_{i}\right\Vert ^{2}\right) ^{\frac{%
1}{2}}-\left\vert \sum_{i=1}^{\infty }p_{i}\left\langle
x_{i},y_{i}\right\rangle \right\vert  \label{revcbs2.4} \\
& \leq \left( \sum_{i=1}^{\infty }p_{i}\left\Vert x_{i}\right\Vert
^{2}\sum_{i=1}^{\infty }p_{i}\left\Vert y_{i}\right\Vert ^{2}\right) ^{\frac{%
1}{2}}-\sum_{i=1}^{\infty }p_{i}\func{Re}\left\langle
x_{i},y_{i}\right\rangle  \notag \\
& \leq \frac{r^{2}\cdot \sum\limits_{i=1}^{\infty }p_{i}\func{Re}%
\left\langle x_{i},y_{i}\right\rangle }{\sqrt{\sum\limits_{i=1}^{\infty
}p_{i}\left\Vert y_{i}\right\Vert ^{2}-r^{2}}\left[ \left(
\sum\limits_{i=1}^{\infty }p_{i}\left\Vert y_{i}\right\Vert ^{2}\right) ^{%
\frac{1}{2}}+\sqrt{\sum\limits_{i=1}^{\infty }p_{i}\left\Vert
y_{i}\right\Vert ^{2}-r^{2}}\right] }  \notag
\end{align}%
and%
\begin{align}
(0& \leq )\sum_{i=1}^{\infty }p_{i}\left\Vert x_{i}\right\Vert
^{2}\sum_{i=1}^{\infty }p_{i}\left\Vert y_{i}\right\Vert ^{2}-\left\vert
\sum_{i=1}^{\infty }p_{i}\left\langle x_{i},y_{i}\right\rangle \right\vert
^{2}  \label{revcbs2.5} \\
& \leq \sum_{i=1}^{\infty }p_{i}\left\Vert x_{i}\right\Vert
^{2}\sum_{i=1}^{\infty }p_{i}\left\Vert y_{i}\right\Vert ^{2}-\left[
\sum_{i=1}^{\infty }p_{i}\func{Re}\left\langle x_{i},y_{i}\right\rangle %
\right] ^{2}  \notag \\
& \leq r^{2}\sum_{i=1}^{\infty }p_{i}\left\Vert x_{i}\right\Vert ^{2}. 
\notag
\end{align}
\end{proposition}

\begin{proof}
From (\ref{revcbs2.2}), we have%
\begin{equation*}
\left\Vert x-y\right\Vert _{\mathbf{p}}^{2}=\sum_{i=1}^{\infty
}p_{i}\left\Vert x_{i}-y_{i}\right\Vert ^{2}\leq r^{2}\sum_{i=1}^{\infty
}p_{i}\leq \sum_{i=1}^{\infty }p_{i}\left\Vert y_{i}\right\Vert
^{2}=\left\Vert y\right\Vert _{\mathbf{p}}^{2},
\end{equation*}%
giving $\left\Vert x-y\right\Vert _{\mathbf{p}}\leq r\leq \left\Vert
y\right\Vert _{\mathbf{p}}.$ Applying Theorem \ref{revcbst1.1} for $\ell _{%
\mathbf{p}}^{2}\left( K\right) $ and $\left\langle \cdot ,\cdot
\right\rangle _{\mathbf{p}},$ we deduce the desired inequality.
\end{proof}

The following proposition holds \cite{DRA3xx}.

\begin{proposition}
\label{revcbsp2.2}Let $x,y\in \ell _{\mathbf{p}}^{2}\left( K\right) $ and $%
a,A\in \mathbb{K}$. If%
\begin{equation}
\func{Re}\left\langle Ay_{i}-x_{i},x_{i}-ay_{i}\right\rangle \geq 0\quad 
\text{for each }i\in \mathbb{N}  \label{revcbs2.6}
\end{equation}%
or, equivalently,%
\begin{equation}
\left\Vert x_{i}-\frac{a+A}{2}y_{i}\right\Vert \leq \frac{1}{2}\left\vert
A-a\right\vert \left\Vert y_{i}\right\Vert \quad \text{for each }i\in 
\mathbb{N}  \label{revcbs2.7}
\end{equation}%
then%
\begin{align}
(0& \leq )\sum_{i=1}^{\infty }p_{i}\left\Vert x_{i}\right\Vert
^{2}\sum_{i=1}^{\infty }p_{i}\left\Vert y_{i}\right\Vert ^{2}-\left\vert
\sum_{i=1}^{\infty }p_{i}\left\langle x_{i},y_{i}\right\rangle \right\vert
^{2}  \label{revcbs2.8} \\
& \leq \frac{1}{4}\left\vert A-a\right\vert ^{2}\left( \sum_{i=1}^{\infty
}p_{i}\left\Vert y_{i}\right\Vert ^{2}\right) ^{2}  \notag \\
& \quad -\left\{ 
\begin{array}{c}
\left\vert \frac{A+a}{2}\sum_{i=1}^{\infty }p_{i}\left\Vert y_{i}\right\Vert
^{2}-\sum_{i=1}^{\infty }p_{i}\left\langle x_{i},y_{i}\right\rangle
\right\vert ^{2} \\ 
\\ 
\sum_{i=1}^{\infty }p_{i}\left\Vert y_{i}\right\Vert ^{2}\sum_{i=1}^{\infty
}p_{i}\func{Re}\left\langle Ay_{i}-x_{i},x_{i}-ay_{i}\right\rangle%
\end{array}%
\right.  \notag \\
& \leq \frac{1}{4}\left\vert A-a\right\vert ^{2}\left( \sum_{i=1}^{\infty
}p_{i}\left\Vert y_{i}\right\Vert ^{2}\right) ^{2}.  \notag
\end{align}
\end{proposition}

The proof follows by Theorem \ref{revcbst1.2}, we omit the details.

Finally, on using Theorem \ref{revcbst1.3}, we may state \cite{DRA3xx}:

\begin{proposition}
\label{revcbsp2.3}Assume that $x,y,a$ and $A$ are as in Proposition \ref%
{revcbsp2.2}. Moreover, if $\func{Re}\left( A\bar{a}\right) >0,$ then we
have the inequality:%
\begin{align}
& \left( \sum_{i=1}^{\infty }p_{i}\left\Vert x_{i}\right\Vert
^{2}\sum_{i=1}^{\infty }p_{i}\left\Vert y_{i}\right\Vert ^{2}\right) ^{\frac{%
1}{2}}  \label{revcbs2.9} \\
& \leq \frac{1}{2}\cdot \frac{\func{Re}\left[ \left( \bar{A}+\bar{a}\right)
\sum_{i=1}^{\infty }p_{i}\left\langle x_{i},y_{i}\right\rangle \right] }{%
\sqrt{\func{Re}\left( A\bar{a}\right) }}  \notag \\
& \leq \frac{1}{2}\cdot \frac{\left\vert A-a\right\vert }{\sqrt{\func{Re}%
\left( A\bar{a}\right) }}\left\vert \sum_{i=1}^{\infty }p_{i}\left\langle
x_{i},y_{i}\right\rangle \right\vert ,  \notag
\end{align}%
\begin{align}
(0& \leq )\left( \sum_{i=1}^{\infty }p_{i}\left\Vert x_{i}\right\Vert
^{2}\sum_{i=1}^{\infty }p_{i}\left\Vert y_{i}\right\Vert ^{2}\right) ^{\frac{%
1}{2}}-\sum_{i=1}^{\infty }p_{i}\func{Re}\left\langle
x_{i},y_{i}\right\rangle  \label{revcbs2.10} \\
& \leq \frac{1}{2}\cdot \frac{\func{Re}\left[ \left( \bar{A}+\bar{a}-2\sqrt{%
\func{Re}\left( A\bar{a}\right) }\right) \sum_{i=1}^{\infty
}p_{i}\left\langle x_{i},y_{i}\right\rangle \right] }{\sqrt{\func{Re}\left( A%
\bar{a}\right) }}  \notag
\end{align}%
and%
\begin{align}
(0& \leq )\sum_{i=1}^{\infty }p_{i}\left\Vert x_{i}\right\Vert
^{2}\sum_{i=1}^{\infty }p_{i}\left\Vert y_{i}\right\Vert ^{2}-\left\vert
\sum_{i=1}^{\infty }p_{i}\left\langle x_{i},y_{i}\right\rangle \right\vert
^{2}  \label{revcbs2.11} \\
& \leq \frac{1}{4}\cdot \frac{\left\vert A-a\right\vert ^{2}}{\func{Re}%
\left( A\bar{a}\right) }\left\vert \sum_{i=1}^{\infty }p_{i}\left\langle
x_{i},y_{i}\right\rangle \right\vert ^{2}.  \notag
\end{align}
\end{proposition}

\subsection{Reverses of the $\left( CBS\right) -$Inequality for Mixed
Sequences}

Let $\left( K,\left\langle \cdot ,\cdot \right\rangle \right) $ be a Hilbert
space over $\mathbb{K}$ and for $p_{i}\geq 0,$ $i\in \mathbb{N}$ with $%
\sum_{i=1}^{\infty }p_{i}=1,$ and $\ell _{\mathbf{p}}^{2}\left( K\right) $
the Hilbert space defined in the previous section.

If%
\begin{equation*}
\alpha \in \ell _{\mathbf{p}}^{2}\left( \mathbb{K}\right) :=\left\{ \alpha
=\left( \alpha _{i}\right) _{i\in \mathbb{N}}\left\vert \alpha _{i}\in 
\mathbb{K},\ i\in \mathbb{N}\text{ \ and \ }\sum_{i=1}^{\infty
}p_{i}\left\vert \alpha _{i}\right\vert ^{2}<\infty \right. \right\}
\end{equation*}%
and $x\in \ell _{\mathbf{p}}^{2}\left( K\right) ,$ then the following
Cauchy-Bunyakovsky-Schwarz $\left( CBS\right) $ inequality holds true:%
\begin{equation}
\sum_{i=1}^{\infty }p_{i}\left\vert \alpha _{i}\right\vert
^{2}\sum_{i=1}^{\infty }p_{i}\left\Vert x_{i}\right\Vert ^{2}\geq \left\Vert
\sum_{i=1}^{\infty }p_{i}\alpha _{i}x_{i}\right\Vert ^{2},  \label{revcbs3.1}
\end{equation}%
with equality if and only if there exists a vector $v\in K$ such that $x_{i}=%
\overline{\alpha _{i}}v$ for any $i\in \mathbb{N}$.

The inequality (\ref{revcbs3.1}) follows by the obvious identity%
\begin{multline*}
\qquad \sum_{i=1}^{n}p_{i}\left\vert \alpha _{i}\right\vert
^{2}\sum_{i=1}^{n}p_{i}\left\Vert x_{i}\right\Vert ^{2}-\left\Vert
\sum_{i=1}^{n}p_{i}\alpha _{i}x_{i}\right\Vert ^{2} \\
=\frac{1}{2}\sum_{i=1}^{n}\sum_{j=1}^{n}p_{i}p_{j}\left\Vert \overline{%
\alpha _{i}}x_{j}-\overline{\alpha _{j}}x_{i}\right\Vert ^{2},\qquad
\end{multline*}%
for any $n\in \mathbb{N}$, $n\geq 1.$

In the following we establish some reverses of the $\left( CBS\right) -$%
inequality in some of its various equivalent forms that will be specified
where they occur \cite{DRA3xx}.

\begin{theorem}[Dragomir, 2005]
\label{revcbst3.1}Let $\alpha \in \ell _{\mathbf{p}}^{2}\left( \mathbb{K}%
\right) ,$ $x\in \ell _{\mathbf{p}}^{2}\left( K\right) $ and $a\in K,$ $r>0$
such that $\left\Vert a\right\Vert >r.$ If the following condition holds%
\begin{equation}
\left\Vert x_{i}-\overline{\alpha _{i}}a\right\Vert \leq r\left\vert \alpha
_{i}\right\vert \quad \text{for each }i\in \mathbb{N},  \label{revcbs3.2}
\end{equation}%
(note that if $\alpha _{i}\neq 0$ for any $i\in \mathbb{N}$, then the
condition (\ref{revcbs3.2}) is equivalent to%
\begin{equation}
\left\Vert \frac{x_{i}}{\overline{\alpha _{i}}}-a\right\Vert \leq r\quad 
\text{for each }i\in \mathbb{N}),  \label{revcbs3.3}
\end{equation}%
then we have the following inequalities 
\begin{align}
\left( \sum_{i=1}^{\infty }p_{i}\left\vert \alpha _{i}\right\vert
^{2}\sum_{i=1}^{\infty }p_{i}\left\Vert x_{i}\right\Vert ^{2}\right) ^{\frac{%
1}{2}}& \leq \frac{1}{\sqrt{\left\Vert a\right\Vert ^{2}-r^{2}}}\func{Re}%
\left\langle \sum_{i=1}^{\infty }p_{i}\alpha _{i}x_{i},a\right\rangle
\label{revcbs3.4} \\
& \leq \frac{\left\Vert a\right\Vert }{\sqrt{\left\Vert a\right\Vert
^{2}-r^{2}}}\left\Vert \sum_{i=1}^{\infty }p_{i}\alpha _{i}x_{i}\right\Vert ;
\notag
\end{align}%
\begin{align}
0& \leq \left( \sum_{i=1}^{\infty }p_{i}\left\vert \alpha _{i}\right\vert
^{2}\sum_{i=1}^{\infty }p_{i}\left\Vert x_{i}\right\Vert ^{2}\right) ^{\frac{%
1}{2}}-\left\Vert \sum_{i=1}^{\infty }p_{i}\alpha _{i}x_{i}\right\Vert
\label{revcbs3.5} \\
& \leq \left( \sum_{i=1}^{\infty }p_{i}\left\vert \alpha _{i}\right\vert
^{2}\sum_{i=1}^{\infty }p_{i}\left\Vert x_{i}\right\Vert ^{2}\right) ^{\frac{%
1}{2}}-\func{Re}\left\langle \sum_{i=1}^{\infty }p_{i}\alpha _{i}x_{i},\frac{%
a}{\left\Vert a\right\Vert }\right\rangle  \notag
\end{align}%
\begin{align}
& \leq \frac{r^{2}}{\sqrt{\left\Vert a\right\Vert ^{2}-r^{2}}\left(
\left\Vert a\right\Vert +\sqrt{\left\Vert a\right\Vert ^{2}-r^{2}}\right) }%
\func{Re}\left\langle \sum_{i=1}^{\infty }p_{i}\alpha _{i}x_{i},\frac{a}{%
\left\Vert a\right\Vert }\right\rangle  \notag \\
& \leq \frac{r^{2}}{\sqrt{\left\Vert a\right\Vert ^{2}-r^{2}}\left(
\left\Vert a\right\Vert +\sqrt{\left\Vert a\right\Vert ^{2}-r^{2}}\right) }%
\left\Vert \sum_{i=1}^{\infty }p_{i}\alpha _{i}x_{i}\right\Vert ;  \notag
\end{align}%
\begin{align}
\sum_{i=1}^{\infty }p_{i}\left\vert \alpha _{i}\right\vert
^{2}\sum_{i=1}^{\infty }p_{i}\left\Vert x_{i}\right\Vert ^{2}& \leq \frac{1}{%
\left\Vert a\right\Vert ^{2}-r^{2}}\left[ \func{Re}\left\langle
\sum_{i=1}^{\infty }p_{i}\alpha _{i}x_{i},a\right\rangle \right] ^{2}
\label{revcbs3.6} \\
& \leq \frac{\left\Vert a\right\Vert ^{2}}{\left\Vert a\right\Vert ^{2}-r^{2}%
}\left\Vert \sum_{i=1}^{\infty }p_{i}\alpha _{i}x_{i}\right\Vert ^{2}  \notag
\end{align}%
and%
\begin{align}
0& \leq \sum_{i=1}^{\infty }p_{i}\left\vert \alpha _{i}\right\vert
^{2}\sum_{i=1}^{\infty }p_{i}\left\Vert x_{i}\right\Vert ^{2}-\left\Vert
\sum_{i=1}^{\infty }p_{i}\alpha _{i}x_{i}\right\Vert ^{2}  \label{revcbs3.7}
\\
& \leq \sum_{i=1}^{\infty }p_{i}\left\vert \alpha _{i}\right\vert
^{2}\sum_{i=1}^{\infty }p_{i}\left\Vert x_{i}\right\Vert ^{2}-\left[ \func{Re%
}\left\langle \sum_{i=1}^{\infty }p_{i}\alpha _{i}x_{i},\frac{a}{\left\Vert
a\right\Vert }\right\rangle \right] ^{2}  \notag \\
& \leq \frac{r^{2}}{\left\Vert a\right\Vert ^{2}\left( \left\Vert
a\right\Vert ^{2}-r^{2}\right) }\left[ \func{Re}\left\langle
\sum_{i=1}^{\infty }p_{i}\alpha _{i}x_{i},a\right\rangle \right] ^{2}  \notag
\\
& \leq \frac{r^{2}}{\left\Vert a\right\Vert ^{2}-r^{2}}\left\Vert
\sum_{i=1}^{\infty }p_{i}\alpha _{i}x_{i}\right\Vert ^{2}.  \notag
\end{align}%
All the inequalities in (\ref{revcbs3.4}) -- (\ref{revcbs3.7}) are sharp.
\end{theorem}

\begin{proof}
From (\ref{revcbs3.2}) we deduce%
\begin{equation*}
\left\Vert x_{i}\right\Vert ^{2}-2\func{Re}\left\langle x_{i},\overline{%
\alpha _{i}}a\right\rangle +\left\vert \alpha _{i}\right\vert ^{2}\left\Vert
a\right\Vert ^{2}\leq \left\vert \alpha _{i}\right\vert ^{2}r^{2}
\end{equation*}%
for any $i\in \mathbb{N}$, which is clearly equivalent to%
\begin{equation}
\left\Vert x_{i}\right\Vert ^{2}+\left( \left\Vert a\right\Vert
^{2}-r^{2}\right) \left\vert \alpha _{i}\right\vert ^{2}\leq 2\func{Re}%
\left\langle \alpha _{i}x_{i},a\right\rangle  \label{revcbs3.8}
\end{equation}%
for each $i\in \mathbb{N}$.

If we multiply (\ref{revcbs3.8}) by $p_{i}\geq 0$ and sum over $i\in \mathbb{%
N}$, then we deduce%
\begin{equation}
\sum_{i=1}^{\infty }p_{i}\left\Vert x_{i}\right\Vert ^{2}+\left( \left\Vert
a\right\Vert ^{2}-r^{2}\right) \sum_{i=1}^{\infty }p_{i}\left\vert \alpha
_{i}\right\vert ^{2}\leq 2\func{Re}\left\langle \sum_{i=1}^{\infty
}p_{i}\alpha _{i}x_{i},a\right\rangle .  \label{revcbs3.9}
\end{equation}%
Now, dividing (\ref{revcbs3.9}) by $\sqrt{\left\Vert a\right\Vert ^{2}-r^{2}}%
>0$ we get%
\begin{multline}
\frac{1}{\sqrt{\left\Vert a\right\Vert ^{2}-r^{2}}}\sum_{i=1}^{\infty
}p_{i}\left\Vert x_{i}\right\Vert ^{2}+\sqrt{\left\Vert a\right\Vert
^{2}-r^{2}}\sum_{i=1}^{\infty }p_{i}\left\vert \alpha _{i}\right\vert ^{2}
\label{revcbs3.10} \\
\leq \frac{2}{\sqrt{\left\Vert a\right\Vert ^{2}-r^{2}}}\func{Re}%
\left\langle \sum_{i=1}^{\infty }p_{i}\alpha _{i}x_{i},a\right\rangle .
\end{multline}%
On the other hand, by the elementary inequality%
\begin{equation*}
\frac{1}{\alpha }p+\alpha q\geq 2\sqrt{pq},\qquad \alpha >0,\ p,q\geq 0,
\end{equation*}%
we can state that:%
\begin{multline}
2\left[ \sum_{i=1}^{\infty }p_{i}\left\vert \alpha _{i}\right\vert
^{2}\sum_{i=1}^{\infty }p_{i}\left\Vert x_{i}\right\Vert ^{2}\right] ^{\frac{%
1}{2}}  \label{revcbs3.11} \\
\leq \frac{1}{\sqrt{\left\Vert a\right\Vert ^{2}-r^{2}}}\sum_{i=1}^{\infty
}p_{i}\left\Vert x_{i}\right\Vert ^{2}+\sqrt{\left\Vert a\right\Vert
^{2}-r^{2}}\sum_{i=1}^{\infty }p_{i}\left\vert \alpha _{i}\right\vert ^{2}.
\end{multline}%
Making use of (\ref{revcbs3.10}) and (\ref{revcbs3.11}), we deduce the first
part of (\ref{revcbs3.4}).

The second part is obvious by Schwarz's inequality%
\begin{equation*}
\func{Re}\left\langle \sum_{i=1}^{\infty }p_{i}\alpha
_{i}x_{i},a\right\rangle \leq \left\Vert \sum_{i=1}^{\infty }p_{i}\alpha
_{i}x_{i}\right\Vert \left\Vert a\right\Vert .
\end{equation*}

If $p_{1}=1,$ $x_{1}=x,$ $\alpha _{1}=1$ and $p_{i}=0,$ $\alpha _{i}=0,$ $%
x_{i}=0$ for $i\geq 2,$ then from (\ref{revcbs3.4}) we deduce the inequality%
\begin{equation*}
\left\Vert x\right\Vert \leq \frac{1}{\sqrt{\left\Vert a\right\Vert
^{2}-r^{2}}}\func{Re}\left\langle x,a\right\rangle \leq \frac{\left\Vert
x\right\Vert \left\Vert a\right\Vert }{\sqrt{\left\Vert a\right\Vert
^{2}-r^{2}}}
\end{equation*}%
provided $\left\Vert x-a\right\Vert \leq r<\left\Vert a\right\Vert ,$ $%
x,a\in K.$ The sharpness of this inequality has been shown in \cite[p. 20]%
{SSDxx}, and we omit the details.

The other inequalities are obvious consequences of (\ref{revcbs3.4}) and we
omit the details.
\end{proof}

The following corollary may be stated \cite{DRA3xx}.

\begin{corollary}
\label{revcbsc3.2}Let $\alpha \in \ell _{\mathbf{p}}^{2}\left( \mathbb{K}%
\right) ,$ $x\in \ell _{\mathbf{p}}^{2}\left( K\right) ,$ $e\in H,$ $%
\left\Vert e\right\Vert =1$ and $\varphi ,\phi \in \mathbb{K}$ with $\func{Re%
}\left( \phi \bar{\varphi}\right) >0.$ If%
\begin{equation}
\left\Vert x_{i}-\overline{\alpha _{i}}\cdot \frac{\varphi +\phi }{2}\cdot
e\right\Vert \leq \frac{1}{2}\left\vert \phi -\varphi \right\vert \left\vert
\alpha _{i}\right\vert  \label{revcbs3.12}
\end{equation}%
for each $i\in \mathbb{N}$, or, equivalently%
\begin{equation}
\func{Re}\left\langle \phi \overline{\alpha _{i}}e-x_{i},x_{i}-\varphi 
\overline{\alpha _{i}}e\right\rangle \geq 0  \label{revcbs3.13}
\end{equation}%
for each $i\in \mathbb{N}$, (note that, if $\alpha _{i}\neq 0$ for any $i\in 
\mathbb{N}$, then (\ref{revcbs3.12}) is equivalent to%
\begin{equation}
\left\Vert \frac{x_{i}}{\overline{\alpha _{i}}}-\frac{\varphi +\phi }{2}%
\cdot e\right\Vert \leq \frac{1}{2}\left\vert \phi -\varphi \right\vert
\label{revcbs3.14}
\end{equation}%
for each $i\in \mathbb{N}$ and (\ref{revcbs3.13}) is equivalent to%
\begin{equation*}
\func{Re}\left\langle \phi e-\frac{x_{i}}{\overline{\alpha _{i}}},\frac{x_{i}%
}{\overline{\alpha _{i}}}-\varphi e\right\rangle \geq 0
\end{equation*}%
for each $i\in \mathbb{N}$), then the following reverses of the $\left(
CBS\right) -$inequality are valid:%
\begin{align}
\left( \sum_{i=1}^{\infty }p_{i}\left\vert \alpha _{i}\right\vert
^{2}\sum_{i=1}^{\infty }p_{i}\left\Vert x_{i}\right\Vert ^{2}\right) ^{\frac{%
1}{2}}& \leq \frac{\func{Re}\left[ \left( \bar{\phi}+\bar{\varphi}\right)
\left\langle \sum_{i=1}^{\infty }p_{i}\alpha _{i}x_{i},e\right\rangle \right]
}{2\left[ \func{Re}\left( \phi \overline{\varphi }\right) \right] ^{\frac{1}{%
2}}}  \label{revcbs3.15} \\
& \leq \frac{1}{2}\cdot \frac{\left\vert \varphi +\phi \right\vert }{\left[ 
\func{Re}\left( \phi \overline{\varphi }\right) \right] ^{\frac{1}{2}}}%
\left\Vert \sum_{i=1}^{\infty }p_{i}\alpha _{i}x_{i}\right\Vert ;  \notag
\end{align}%
\begin{align}
0& \leq \left( \sum_{i=1}^{\infty }p_{i}\left\vert \alpha _{i}\right\vert
^{2}\sum_{i=1}^{\infty }p_{i}\left\Vert x_{i}\right\Vert ^{2}\right) ^{\frac{%
1}{2}}-\left\Vert \sum_{i=1}^{\infty }p_{i}\alpha _{i}x_{i}\right\Vert
\label{revcbs3.16} \\
& \leq \left( \sum_{i=1}^{\infty }p_{i}\left\vert \alpha _{i}\right\vert
^{2}\sum_{i=1}^{\infty }p_{i}\left\Vert x_{i}\right\Vert ^{2}\right) ^{\frac{%
1}{2}}  \notag \\
& \qquad -\func{Re}\left[ \frac{\bar{\phi}+\bar{\varphi}}{\left\vert \varphi
+\phi \right\vert }\left\langle \sum_{i=1}^{\infty }p_{i}\alpha
_{i}x_{i},e\right\rangle \right]  \notag \\
& \leq \frac{\left\vert \phi -\varphi \right\vert ^{2}}{2\sqrt{\func{Re}%
\left( \phi \varphi \right) }\left( \left\vert \varphi +\phi \right\vert +2%
\sqrt{\func{Re}\left( \phi \overline{\varphi }\right) }\right) }  \notag \\
& \qquad \times \func{Re}\left[ \frac{\bar{\phi}+\bar{\varphi}}{\left\vert
\varphi +\phi \right\vert }\left\langle \sum_{i=1}^{\infty }p_{i}\alpha
_{i}x_{i},e\right\rangle \right]  \notag \\
& \leq \frac{\left\vert \phi -\varphi \right\vert ^{2}}{2\sqrt{\func{Re}%
\left( \phi \varphi \right) }\left( \left\vert \varphi +\phi \right\vert +2%
\sqrt{\func{Re}\left( \phi \overline{\varphi }\right) }\right) }\left\Vert
\sum_{i=1}^{\infty }p_{i}\alpha _{i}x_{i}\right\Vert ;  \notag
\end{align}%
\begin{align}
& \sum_{i=1}^{\infty }p_{i}\left\vert \alpha _{i}\right\vert
^{2}\sum_{i=1}^{\infty }p_{i}\left\Vert x_{i}\right\Vert ^{2}
\label{revcbs3.17} \\
& \leq \frac{1}{4\func{Re}\left( \phi \bar{\varphi}\right) }\left[ \func{Re}%
\left\{ \left( \bar{\phi}+\bar{\varphi}\right) \left\langle
\sum_{i=1}^{\infty }p_{i}\alpha _{i}x_{i},e\right\rangle \right\} \right]
^{2}  \notag \\
& \leq \frac{1}{4}\cdot \frac{\left\vert \varphi +\phi \right\vert ^{2}}{%
\func{Re}\left( \phi \bar{\varphi}\right) }\left\Vert \sum_{i=1}^{\infty
}p_{i}\alpha _{i}x_{i}\right\Vert ^{2}  \notag
\end{align}%
and%
\begin{align}
0& \leq \sum_{i=1}^{\infty }p_{i}\left\vert \alpha _{i}\right\vert
^{2}\sum_{i=1}^{\infty }p_{i}\left\Vert x_{i}\right\Vert ^{2}-\left\Vert
\sum_{i=1}^{\infty }p_{i}\alpha _{i}x_{i}\right\Vert ^{2}  \label{revcbs3.18}
\\
& \leq \sum_{i=1}^{\infty }p_{i}\left\vert \alpha _{i}\right\vert
^{2}\sum_{i=1}^{\infty }p_{i}\left\Vert x_{i}\right\Vert ^{2}  \notag \\
& \qquad -\left[ \func{Re}\left\{ \frac{\bar{\phi}+\bar{\varphi}}{\left\vert
\varphi +\phi \right\vert }\left\langle \sum_{i=1}^{\infty }p_{i}\alpha
_{i}x_{i},e\right\rangle \right\} \right] ^{2}  \notag \\
& \leq \frac{\left\vert \phi -\varphi \right\vert ^{2}}{4\left\vert \phi
+\varphi \right\vert ^{2}\func{Re}\left( \phi \bar{\varphi}\right) }\left\{ 
\func{Re}\left[ \left( \bar{\phi}+\bar{\varphi}\right) \left\langle
\sum_{i=1}^{\infty }p_{i}\alpha _{i}x_{i},e\right\rangle \right] \right\}
^{2}  \notag \\
& \leq \frac{\left\vert \phi -\varphi \right\vert ^{2}}{4\func{Re}\left(
\phi \bar{\varphi}\right) }\left\Vert \sum_{i=1}^{\infty }p_{i}\alpha
_{i}x_{i}\right\Vert ^{2}.  \notag
\end{align}%
All the inequalities in (\ref{revcbs3.15}) -- (\ref{revcbs3.18}) are sharp.
\end{corollary}

\begin{remark}
\label{revcbsr2.3}We remark that if $M\geq m>0$ and for $\alpha \in \ell _{%
\mathbf{p}}^{2}\left( \mathbb{K}\right) ,$ $x\in \ell _{\mathbf{p}%
}^{2}\left( K\right) ,$ $e\in H$ with $\left\Vert e\right\Vert =1,$ one
would assume that either%
\begin{equation}
\left\Vert \frac{x_{i}}{\overline{\alpha _{i}}}-\frac{M+m}{2}\cdot
e\right\Vert \leq \frac{1}{2}\left( M-m\right)  \label{revcbs3.19}
\end{equation}%
for each $i\in \mathbb{N}$, or, equivalently%
\begin{equation}
\func{Re}\left\langle Me-\frac{x_{i}}{\overline{\alpha _{i}}},\frac{x_{i}}{%
\overline{\alpha _{i}}}-me\right\rangle \geq 0  \label{revcbs3.20}
\end{equation}%
for each $i\in \mathbb{N}$, then the following, much simpler reverses of the 
$\left( CBS\right) -$ inequality may be stated:%
\begin{align}
\left( \sum_{i=1}^{\infty }p_{i}\left\vert \alpha _{i}\right\vert
^{2}\sum_{i=1}^{\infty }p_{i}\left\Vert x_{i}\right\Vert ^{2}\right) ^{\frac{%
1}{2}}& \leq \frac{M+m}{2\sqrt{mM}}\func{Re}\left\langle \sum_{i=1}^{\infty
}p_{i}\alpha _{i}x_{i},e\right\rangle  \label{revcbs3.21} \\
& \leq \frac{M+m}{2\sqrt{mM}}\left\Vert \sum_{i=1}^{\infty }p_{i}\alpha
_{i}x_{i}\right\Vert ;  \notag
\end{align}%
\begin{align}
0& \leq \left( \sum_{i=1}^{\infty }p_{i}\left\vert \alpha _{i}\right\vert
^{2}\sum_{i=1}^{\infty }p_{i}\left\Vert x_{i}\right\Vert ^{2}\right) ^{\frac{%
1}{2}}-\left\Vert \sum_{i=1}^{\infty }p_{i}\alpha _{i}x_{i}\right\Vert
\label{revcbs3.22} \\
& \leq \left( \sum_{i=1}^{\infty }p_{i}\left\vert \alpha _{i}\right\vert
^{2}\sum_{i=1}^{\infty }p_{i}\left\Vert x_{i}\right\Vert ^{2}\right) ^{\frac{%
1}{2}}-\func{Re}\left\langle \sum_{i=1}^{\infty }p_{i}\alpha
_{i}x_{i},e\right\rangle  \notag \\
& \leq \frac{\left( M-m\right) ^{2}}{2\left( \sqrt{M}+\sqrt{m}\right) ^{2}%
\sqrt{mM}}\func{Re}\left\langle \sum_{i=1}^{\infty }p_{i}\alpha
_{i}x_{i},e\right\rangle  \notag \\
& \leq \frac{\left( M-m\right) ^{2}}{2\left( \sqrt{M}+\sqrt{m}\right) ^{2}%
\sqrt{mM}}\left\Vert \sum_{i=1}^{\infty }p_{i}\alpha _{i}x_{i}\right\Vert ; 
\notag
\end{align}%
\begin{align}
& \sum_{i=1}^{\infty }p_{i}\left\vert \alpha _{i}\right\vert
^{2}\sum_{i=1}^{\infty }p_{i}\left\Vert x_{i}\right\Vert ^{2}-\left\Vert
\sum_{i=1}^{\infty }p_{i}\alpha _{i}x_{i}\right\Vert ^{2}  \label{revcbs3.23}
\\
& \leq \frac{\left( M+m\right) ^{2}}{4mM}\left[ \func{Re}\left\langle
\sum_{i=1}^{\infty }p_{i}\alpha _{i}x_{i},e\right\rangle \right] ^{2}  \notag
\\
& \leq \frac{\left( M+m\right) ^{2}}{4mM}\left\Vert \sum_{i=1}^{\infty
}p_{i}\alpha _{i}x_{i}\right\Vert ^{2}  \notag
\end{align}%
and%
\begin{equation}
0\leq \sum_{i=1}^{\infty }p_{i}\left\vert \alpha _{i}\right\vert
^{2}\sum_{i=1}^{\infty }p_{i}\left\Vert x_{i}\right\Vert ^{2}-\left\Vert
\sum_{i=1}^{\infty }p_{i}\alpha _{i}x_{i}\right\Vert ^{2}  \label{revcbs3.24}
\end{equation}%
\begin{align*}
& \leq \sum_{i=1}^{\infty }p_{i}\left\vert \alpha _{i}\right\vert
^{2}\sum_{i=1}^{\infty }p_{i}\left\Vert x_{i}\right\Vert ^{2}-\left[ \func{Re%
}\left\langle \sum_{i=1}^{\infty }p_{i}\alpha _{i}x_{i},e\right\rangle %
\right] ^{2} \\
& \leq \frac{\left( M-m\right) ^{2}}{4mM}\left[ \func{Re}\left\langle
\sum_{i=1}^{\infty }p_{i}\alpha _{i}x_{i},e\right\rangle \right] ^{2} \\
& \leq \frac{\left( M-m\right) ^{2}}{4mM}\left\Vert \sum_{i=1}^{\infty
}p_{i}\alpha _{i}x_{i}\right\Vert ^{2}.
\end{align*}
\end{remark}

\subsection{Reverses for the Generalised Triangle Inequality}

In 1966, J.B. Diaz and F.T. Metcalf \cite{DMxx} proved the following reverse
of the generalised triangle inequality holding in an inner product space $%
\left( H;\left\langle \cdot ,\cdot \right\rangle \right) $ over the real or
complex number field $\mathbb{K}$:%
\begin{equation}
r\sum_{i=1}^{n}\left\Vert x_{i}\right\Vert \leq \left\Vert
\sum_{i=1}^{n}x_{i}\right\Vert  \label{revcbs4.1}
\end{equation}%
provided the vectors $x_{1},\dots ,x_{n}\in H\backslash \left\{ 0\right\} $
satisfy the assumption%
\begin{equation}
0\leq r\leq \frac{\func{Re}\left\langle x_{i},a\right\rangle }{\left\Vert
x_{i}\right\Vert },  \label{revcbs4.2}
\end{equation}%
where $a\in H$ and $\left\Vert a\right\Vert =1.$

In an attempt to diversify the assumptions for which such reverse results
hold, the author pointed out in \cite{DRA2xx} that%
\begin{equation}
\sqrt{1-\rho ^{2}}\sum_{i=1}^{n}\left\Vert x_{i}\right\Vert \leq \left\Vert
\sum_{i=1}^{n}x_{i}\right\Vert ,  \label{revcbs4.3}
\end{equation}%
where the vectors $x_{i,}i\in \left\{ 1,\dots ,n\right\} $ satisfy the
condition%
\begin{equation}
\left\Vert x_{i}-a\right\Vert \leq \rho ,\qquad i\in \left\{ 1,\dots
,n\right\}  \label{revcbs4.4}
\end{equation}%
where $a\in H,$ $\left\Vert a\right\Vert =1$ and $\rho \in \left( 0,1\right)
.$

If, for $M\geq m>0,$ the vectors $x_{i}\in H,$ $i\in \left\{ 1,\dots
,n\right\} $ verify either%
\begin{equation}
\func{Re}\left\langle Ma-x_{i},x_{i}-ma\right\rangle \geq 0,\qquad i\in
\left\{ 1,\dots ,n\right\} ,  \label{revcbs4.5}
\end{equation}%
or, equivalently,%
\begin{equation}
\left\Vert x_{i}-\frac{M+m}{2}\cdot a\right\Vert \leq \frac{1}{2}\left(
M-m\right) ,\qquad i\in \left\{ 1,\dots ,n\right\} ,  \label{revcbs4.6}
\end{equation}%
where $a\in H,$ $\left\Vert a\right\Vert =1$, then the following reverse of
the generalised triangle inequality may be stated as well \cite{DRA2xx}%
\begin{equation}
\frac{2\sqrt{mM}}{M+m}\sum_{i=1}^{n}\left\Vert x_{i}\right\Vert \leq
\left\Vert \sum_{i=1}^{n}x_{i}\right\Vert .  \label{revcbs4.7}
\end{equation}

Note that the inequalities (\ref{revcbs4.1}), (\ref{revcbs4.3}), and (\ref%
{revcbs4.7}) are sharp; necessary and sufficient equality conditions were
provided (see \cite{DMxx} and \cite{DRA2xx}).

It is obvious, from Theorem \ref{revcbst3.1}, that, if%
\begin{equation}
\left\Vert x_{i}-a\right\Vert \leq r,\quad \text{for }\quad i\in \left\{
1,\dots ,n\right\} ,  \label{revcbs4.8}
\end{equation}%
where $\left\Vert a\right\Vert >r,$ $a\in H$ and $x_{i}\in H$, $i\in \left\{
1,\dots ,n\right\} ,$ then one can state the inequalities%
\begin{align}
\sum_{i=1}^{n}\left\Vert x_{i}\right\Vert & \leq \sqrt{n}\left(
\sum_{i=1}^{n}\left\Vert x_{i}\right\Vert ^{2}\right) ^{\frac{1}{2}}
\label{revcbs4.9} \\
& \leq \frac{1}{\sqrt{\left\Vert a\right\Vert ^{2}-r^{2}}}\func{Re}%
\left\langle \sum_{i=1}^{n}x_{i},a\right\rangle  \notag \\
& \leq \frac{\left\Vert a\right\Vert }{\sqrt{\left\Vert a\right\Vert
^{2}-r^{2}}}\left\Vert \sum_{i=1}^{n}x_{i}\right\Vert  \notag
\end{align}%
and%
\begin{align}
0& \leq \sum_{i=1}^{n}\left\Vert x_{i}\right\Vert -\left\Vert
\sum_{i=1}^{n}x_{i}\right\Vert  \label{revcbs4.10} \\
& \leq \sqrt{n}\left( \sum_{i=1}^{n}\left\Vert x_{i}\right\Vert ^{2}\right)
^{\frac{1}{2}}-\left\Vert \sum_{i=1}^{n}x_{i}\right\Vert  \notag \\
& \leq \sqrt{n}\left( \sum_{i=1}^{n}\left\Vert x_{i}\right\Vert ^{2}\right)
^{\frac{1}{2}}-\func{Re}\left\langle \sum_{i=1}^{n}x_{i},\frac{a}{\left\Vert
a\right\Vert }\right\rangle  \notag \\
& \leq \frac{r^{2}}{\sqrt{\left\Vert a\right\Vert ^{2}-r^{2}}\left(
\left\Vert a\right\Vert +\sqrt{\left\Vert a\right\Vert ^{2}-r^{2}}\right) }%
\func{Re}\left\langle \sum_{i=1}^{n}x_{i},\frac{a}{\left\Vert a\right\Vert }%
\right\rangle  \notag \\
& \leq \frac{r^{2}}{\sqrt{\left\Vert a\right\Vert ^{2}-r^{2}}\left(
\left\Vert a\right\Vert +\sqrt{\left\Vert a\right\Vert ^{2}-r^{2}}\right) }%
\left\Vert \sum_{i=1}^{n}x_{i}\right\Vert .  \notag
\end{align}%
We note that for $\left\Vert a\right\Vert =1$ and $r\in \left( 0,1\right) ,$
the inequality (\ref{revcbs3.9}) becomes%
\begin{equation}
\sqrt{1-r^{2}}\sum_{i=1}^{n}\left\Vert x_{i}\right\Vert \leq \sqrt{\left(
1-r^{2}\right) n}\left( \sum_{i=1}^{n}\left\Vert x_{i}\right\Vert
^{2}\right) ^{\frac{1}{2}}  \label{revcbs4.11}
\end{equation}%
\begin{equation*}
\leq \func{Re}\left\langle \sum_{i=1}^{n}x_{i},a\right\rangle \leq
\left\Vert \sum_{i=1}^{n}x_{i}\right\Vert
\end{equation*}%
which is a refinement of (\ref{revcbs4.3}).

With the same assumptions for $a$ and $r,$ we have from (\ref{revcbs4.10})
the following additive reverse of the generalised triangle inequality:%
\begin{align}
0& \leq \sum_{i=1}^{n}\left\Vert x_{i}\right\Vert -\left\Vert
\sum_{i=1}^{n}x_{i}\right\Vert  \label{revcbs4.12} \\
& \leq \frac{r^{2}}{\sqrt{1-r^{2}}\left( 1+\sqrt{1-r^{2}}\right) }\func{Re}%
\left\langle \sum_{i=1}^{n}x_{i},a\right\rangle  \notag \\
& \leq \frac{r^{2}}{\sqrt{1-r^{2}}\left( 1+\sqrt{1-r^{2}}\right) }\left\Vert
\sum_{i=1}^{n}x_{i}\right\Vert .  \notag
\end{align}

We can obtain the following reverses of the generalised triangle inequality
from Corollary \ref{revcbsc3.2} when the assumptions are in terms of complex
numbers $\phi $ and $\varphi :$

If $\varphi ,\phi \in \mathbb{K}$ with $\func{Re}\left( \phi \bar{\varphi}%
\right) >0$ and $x_{i}\in H,$ $i\in \left\{ 1,\dots ,n\right\} ,$ $e\in H,$ $%
\left\Vert e\right\Vert =1$ are such that%
\begin{equation}
\left\Vert x_{i}-\frac{\varphi +\phi }{2}e\right\Vert \leq \frac{1}{2}%
\left\vert \phi -\varphi \right\vert \text{ \ for each \ }i\in \left\{
1,\dots ,n\right\} ,  \label{revcbs4.13}
\end{equation}%
or, equivalently,%
\begin{equation*}
\func{Re}\left\langle \phi e-x_{i},x_{i}-\varphi e\right\rangle \geq 0\text{
\ for each \ }i\in \left\{ 1,\dots ,n\right\} ,
\end{equation*}%
then we have the following reverses of the generalised triangle inequality:%
\begin{align}
\sum_{i=1}^{n}\left\Vert x_{i}\right\Vert & \leq \sqrt{n}\left(
\sum_{i=1}^{n}\left\Vert x_{i}\right\Vert ^{2}\right) ^{\frac{1}{2}}
\label{revcbs4.14} \\
& \leq \frac{\func{Re}\left[ \left( \bar{\phi}+\bar{\varphi}\right)
\left\langle \sum_{i=1}^{n}x_{i},e\right\rangle \right] }{2\sqrt{\func{Re}%
\left( \phi \bar{\varphi}\right) }}  \notag \\
& \leq \frac{1}{2}\cdot \frac{\left\vert \bar{\phi}+\bar{\varphi}\right\vert 
}{\sqrt{\func{Re}\left( \phi \bar{\varphi}\right) }}\left\Vert
\sum_{i=1}^{n}x_{i}\right\Vert  \notag
\end{align}%
and%
\begin{align}
0& \leq \sum_{i=1}^{n}\left\Vert x_{i}\right\Vert -\left\Vert
\sum_{i=1}^{n}x_{i}\right\Vert  \label{revcbs4.15} \\
& \leq \sqrt{n}\left( \sum_{i=1}^{n}\left\Vert x_{i}\right\Vert ^{2}\right)
^{\frac{1}{2}}-\left\Vert \sum_{i=1}^{n}x_{i}\right\Vert  \notag \\
& \leq \sqrt{n}\left( \sum_{i=1}^{n}\left\Vert x_{i}\right\Vert ^{2}\right)
^{\frac{1}{2}}-\func{Re}\left[ \frac{\left\vert \bar{\phi}+\bar{\varphi}%
\right\vert }{\sqrt{\func{Re}\left( \bar{\phi}\bar{\varphi}\right) }}%
\left\langle \sum_{i=1}^{n}x_{i},e\right\rangle \right]  \notag \\
& \leq \frac{\left\vert \phi -\varphi \right\vert ^{2}}{2\sqrt{\func{Re}%
\left( \phi \bar{\varphi}\right) }\left( \left\vert \phi +\varphi
\right\vert +2\sqrt{\func{Re}\left( \phi \bar{\varphi}\right) }\right) } 
\notag \\
& \qquad \times \func{Re}\left[ \frac{\bar{\phi}+\bar{\varphi}}{\left\vert 
\bar{\phi}+\bar{\varphi}\right\vert }\left\langle
\sum_{i=1}^{n}x_{i},e\right\rangle \right]  \notag \\
& \leq \frac{\left\vert \phi -\varphi \right\vert ^{2}}{2\sqrt{\func{Re}%
\left( \phi \bar{\varphi}\right) }\left( \left\vert \phi +\varphi
\right\vert +2\sqrt{\func{Re}\left( \phi \bar{\varphi}\right) }\right) }%
\left\Vert \sum_{i=1}^{n}x_{i}\right\Vert .  \notag
\end{align}

Obviously (\ref{revcbs4.14}) for $\phi =M,$ $\varphi =m,$ $M\geq m>0$
provides a refinement for (\ref{revcbs4.7}).

\subsection{Lower Bounds for the Distance to Finite-Dimensional Subspaces}

Let $\left( H;\left\langle \cdot ,\cdot \right\rangle \right) $ be an inner
product space over the real or complex number field $\mathbb{K}$, $\left\{
y_{1},\dots ,y_{n}\right\} $ a subset of $H$ and $G\left( y_{1},\dots
,y_{n}\right) $ the \textit{Gram matrix} of $\left\{ y_{1},\dots
,y_{n}\right\} $ where $\left( i,j\right) -$entry is $\left\langle
y_{i},y_{j}\right\rangle .$ The determinant of $G\left( y_{1},\dots
,y_{n}\right) $ is called the \textit{Gram determinant} of $\left\{
y_{1},\dots ,y_{n}\right\} $ and is denoted by $\Gamma \left( y_{1},\dots
,y_{n}\right) .$

Following \cite[p. 129 -- 133]{DExx}, we state here some general results for
the Gram determinant that will be used in the sequel:

\begin{enumerate}
\item Let $\left\{ x_{1},\dots ,x_{n}\right\} \subset H.$ Then $\Gamma
\left( x_{1},\dots ,x_{n}\right) \neq 0$ if and only if $\left\{ x_{1},\dots
,x_{n}\right\} $ is linearly independent;

\item Let $M=span\left\{ x_{1},\dots ,x_{n}\right\} $ be $n-$dimensional in $%
H,$ i.e., $\{x_{1},\dots ,$ $x_{n}\}$ is linearly independent. Then for each 
$x\in H,$ the distance $d\left( x,M\right) $ from $x$ to the linear subspace 
$H$ has the representations%
\begin{equation}
d^{2}\left( x,M\right) =\frac{\Gamma \left( x_{1},\dots ,x_{n},x\right) }{%
\Gamma \left( x_{1},\dots ,x_{n}\right) }  \label{revcbs5.2}
\end{equation}%
and%
\begin{equation}
d^{2}\left( x,M\right) =\left\{ 
\begin{array}{ll}
\left\Vert x\right\Vert ^{2}-\frac{\left( \sum_{i=1}^{n}\left\vert
\left\langle x,x_{i}\right\rangle \right\vert ^{2}\right) ^{2}}{\left\Vert
\sum_{i=1}^{n}\left\langle x,x_{i}\right\rangle x_{i}\right\Vert ^{2}} & 
\text{if \ }x\notin M^{\perp }, \\ 
&  \\ 
\left\Vert x\right\Vert ^{2} & \text{if \ }x\in M^{\perp },%
\end{array}%
\right.  \label{revcbs5.3}
\end{equation}%
where $M^{\perp }$ denotes the orthogonal complement of $M.$
\end{enumerate}

The following result may be stated \cite{DRA3xx}.

\begin{proposition}
\label{revcbsp5.1}Let $\left\{ x_{1},\dots ,x_{n}\right\} $ be a system of
linearly independent vectors, $M=span\left\{ x_{1},\dots ,x_{n}\right\} ,$ $%
x\in H\backslash M^{\perp },$ $a\in H,$ $r>0$ and $\left\Vert a\right\Vert
>r.$ If%
\begin{equation}
\left\Vert x_{i}-\overline{\left\langle x,x_{i}\right\rangle }a\right\Vert
\leq \left\vert \left\langle x,x_{i}\right\rangle \right\vert r\text{ \ for
each \ }i\in \left\{ 1,\dots ,n\right\} ,  \label{revcbs5.6}
\end{equation}%
(note that if $\left\langle x,x_{i}\right\rangle \neq 0$ for each $i\in
\left\{ 1,\dots ,n\right\} ,$ then (\ref{revcbs5.6}) can be written as%
\begin{equation}
\left\Vert \frac{x_{i}}{\overline{\left\langle x,x_{i}\right\rangle }}%
-a\right\Vert \leq r\text{ \ for each \ }i\in \left\{ 1,\dots ,n\right\} ),
\label{revcbs5.7}
\end{equation}%
then we have the inequality%
\begin{align}
d^{2}\left( x,M\right) & \geq \left\Vert x\right\Vert ^{2}-\frac{\left\Vert
a\right\Vert ^{2}}{\left\Vert a\right\Vert ^{2}-r^{2}}\cdot \frac{%
\sum_{i=1}^{n}\left\vert \left\langle x,x_{i}\right\rangle \right\vert ^{2}}{%
\sum_{i=1}^{n}\left\Vert x_{i}\right\Vert ^{2}}  \label{revcbs5.8} \\
& \geq 0.  \notag
\end{align}
\end{proposition}

\begin{proof}
Utilising (\ref{revcbs5.3}) we can state that%
\begin{equation}
d^{2}\left( x,M\right) =\left\Vert x\right\Vert ^{2}-\frac{%
\sum_{i=1}^{n}\left\vert \left\langle x,x_{i}\right\rangle \right\vert ^{2}}{%
\left\Vert \sum_{i=1}^{n}\left\langle x,x_{i}\right\rangle x_{i}\right\Vert
^{2}}\cdot \sum_{i=1}^{n}\left\vert \left\langle x,x_{i}\right\rangle
\right\vert ^{2}.  \label{revcbs5.9}
\end{equation}%
Also, by the inequality (\ref{revcbs3.6}) applied for $\alpha
_{i}=\left\langle x,x_{i}\right\rangle ,$ $p_{i}=\frac{1}{n},$ $i\in \left\{
1,\dots ,n\right\} ,$ we can state that%
\begin{equation}
\frac{\sum_{i=1}^{n}\left\vert \left\langle x,x_{i}\right\rangle \right\vert
^{2}}{\left\Vert \sum_{i=1}^{n}\left\langle x,x_{i}\right\rangle
x_{i}\right\Vert ^{2}}\leq \frac{\left\Vert a\right\Vert ^{2}}{\left\Vert
a\right\Vert ^{2}-r^{2}}\cdot \frac{1}{\sum_{i=1}^{n}\left\Vert
x_{i}\right\Vert ^{2}}  \label{revcbs5.10}
\end{equation}%
provided the condition (\ref{revcbs5.7}) holds true.

Combining (\ref{revcbs5.9}) with (\ref{revcbs5.10}) we deduce the first
inequality in (\ref{revcbs5.8}).

The last inequality is obvious since, by Schwarz's inequality%
\begin{equation*}
\left\Vert x\right\Vert ^{2}\sum_{i=1}^{n}\left\Vert x_{i}\right\Vert
^{2}\geq \sum_{i=1}^{n}\left\vert \left\langle x,x_{i}\right\rangle
\right\vert ^{2}\geq \frac{\left\Vert a\right\Vert ^{2}}{\left\Vert
a\right\Vert ^{2}-r^{2}}\sum_{i=1}^{n}\left\vert \left\langle
x,x_{i}\right\rangle \right\vert ^{2}.
\end{equation*}
\end{proof}

\begin{remark}
Utilising (\ref{revcbs5.2}), we can state the following result for Gram
determinants%
\begin{multline}
\Gamma \left( x_{1},\dots ,x_{n},x\right)  \label{revcbs5.11} \\
\geq \left[ \left\Vert x\right\Vert ^{2}-\frac{\left\Vert a\right\Vert ^{2}}{%
\left\Vert a\right\Vert ^{2}-r^{2}}\cdot \frac{\sum_{i=1}^{n}\left\vert
\left\langle x,x_{i}\right\rangle \right\vert ^{2}}{\sum_{i=1}^{n}\left\Vert
x_{i}\right\Vert ^{2}}\right] \Gamma \left( x_{1},\dots ,x_{n}\right) \geq 0
\end{multline}%
for $x\notin M^{\perp }$ and $x,x_{i},a$ and $r$ are as in Proposition \ref%
{revcbsp5.1}.
\end{remark}

The following corollary of Proposition \ref{revcbsp5.1} may be stated as
well \cite{DRA3xx}.

\begin{corollary}
\label{revcbsc5.3}Let $\left\{ x_{1},\dots ,x_{n}\right\} $ be a system of
linearly independent vectors, $M=span\left\{ x_{1},\dots ,x_{n}\right\} ,$ $%
x\in H\backslash M^{\perp }$ \ and $\phi ,\varphi \in K$ with $\func{Re}%
\left( \phi \bar{\varphi}\right) >0.$ If $e\in H,$ $\left\Vert e\right\Vert
=1$ and%
\begin{equation}
\left\Vert x_{i}-\overline{\left\langle x,x_{i}\right\rangle }\cdot \frac{%
\varphi +\phi }{2}e\right\Vert \leq \frac{1}{2}\left\vert \phi -\varphi
\right\vert \left\vert \left\langle x,x_{i}\right\rangle \right\vert
\label{revcbs5.12}
\end{equation}%
or, equivalently,%
\begin{equation*}
\func{Re}\left\langle \phi \overline{\cdot \left\langle x,x_{i}\right\rangle 
}e-x_{i},x_{i}-\varphi \cdot \overline{\left\langle x,x_{i}\right\rangle }%
e\right\rangle \geq 0,
\end{equation*}%
for each $i\in \left\{ 1,\dots ,n\right\} ,$ then%
\begin{equation}
d^{2}\left( x,M\right) \geq \left\Vert x\right\Vert ^{2}-\frac{1}{4}\cdot 
\frac{\left\vert \varphi +\phi \right\vert ^{2}}{\func{Re}\left( \phi \bar{%
\varphi}\right) }\cdot \frac{\sum_{i=1}^{n}\left\vert \left\langle
x,x_{i}\right\rangle \right\vert ^{2}}{\sum_{i=1}^{n}\left\Vert
x_{i}\right\Vert ^{2}}\geq 0,  \label{revcbs5.13}
\end{equation}%
or, equivalently,%
\begin{multline}
\Gamma \left( x_{1},\dots ,x_{n},x\right)  \label{revcbs5.14} \\
\geq \left[ \left\Vert x\right\Vert ^{2}-\frac{1}{4}\cdot \frac{\left\vert
\varphi +\phi \right\vert ^{2}}{\func{Re}\left( \phi \bar{\varphi}\right) }%
\cdot \frac{\sum_{i=1}^{n}\left\vert \left\langle x,x_{i}\right\rangle
\right\vert ^{2}}{\sum_{i=1}^{n}\left\Vert x_{i}\right\Vert ^{2}}\right]
\Gamma \left( x_{1},\dots ,x_{n}\right) \geq 0.
\end{multline}
\end{corollary}

\subsection{Applications for Fourier Coefficients}

Let $\left( H;\left\langle \cdot ,\cdot \right\rangle \right) $ be a Hilbert
space over the real or complex number field $\mathbb{K}$ and $\left\{
e_{i}\right\} _{i\in I}$ an \textit{orthornormal basis} for $H.$ Then (see
for instance \cite[p. 54 -- 61]{DExx})

\begin{enumerate}
\item[(i)] Every element $x\in H$ can be expanded in a \textit{Fourier
series, }i.e.,%
\begin{equation*}
x=\sum_{i\in I}\left\langle x,e_{i}\right\rangle e_{i},
\end{equation*}%
where $\left\langle x,e_{i}\right\rangle ,$ $i\in I$ are the \textit{Fourier
coefficients} of $x;$

\item[(ii)] (Parseval identity)%
\begin{equation*}
\left\Vert x\right\Vert ^{2}=\sum_{i\in I}\left\langle x,e_{i}\right\rangle
e_{i},\qquad x\in H;
\end{equation*}

\item[(iii)] (Extended Parseval identity)%
\begin{equation*}
\left\langle x,y\right\rangle =\sum_{i\in I}\left\langle
x,e_{i}\right\rangle \left\langle e_{i},y\right\rangle ,\qquad x,y\in H;
\end{equation*}

\item[(iv)] (Elements are uniquely determined by their Fourier coefficients)%
\begin{equation*}
\left\langle x,e_{i}\right\rangle =\left\langle y,e_{i}\right\rangle \text{
\ for every }i\in I\text{ \ implies that }x=y.
\end{equation*}
\end{enumerate}

Now, we must remark that all the results can be stated for $K=\mathbb{K}$
where $\mathbb{K}$ is the Hilbert space of complex (real) numbers endowed
with the usual norm and inner product.

Therefore, we can state the following proposition \cite{DRA3xx}.

\begin{proposition}
\label{revcbsp6.1}Let $\left( H;\left\langle \cdot ,\cdot \right\rangle
\right) $ be a Hilbert space over $\mathbb{K}$ and $\left\{ e_{i}\right\}
_{i\in I}$ an orthornormal base for $H.$ If $x,y\in H$ $\left( y\neq
0\right) ,$ $a\in \mathbb{K}$ $\left( \mathbb{C},\mathbb{R}\right) $ and $%
r>0 $ such that $\left\vert a\right\vert >r$ and%
\begin{equation}
\left\vert \frac{\left\langle x,e_{i}\right\rangle }{\left\langle
y,e_{i}\right\rangle }-a\right\vert \leq r\text{ \ for each \ }i\in I,
\label{revcbs6.1}
\end{equation}%
then we have the following reverse of the Schwarz inequality%
\begin{align}
\left\Vert x\right\Vert \left\Vert y\right\Vert & \leq \frac{1}{\sqrt{%
\left\vert a\right\vert ^{2}-r^{2}}}\func{Re}\left[ \bar{a}\cdot
\left\langle x,y\right\rangle \right]  \label{revcbs6.2} \\
& \leq \frac{\left\vert a\right\vert }{\sqrt{\left\vert a\right\vert
^{2}-r^{2}}}\left\vert \left\langle x,y\right\rangle \right\vert ;  \notag
\end{align}%
\begin{align}
(0& \leq )\left\Vert x\right\Vert \left\Vert y\right\Vert -\left\vert
\left\langle x,y\right\rangle \right\vert  \label{revcbs6.3} \\
& \leq \left\Vert x\right\Vert \left\Vert y\right\Vert -\func{Re}\left[ 
\frac{\bar{a}}{\left\vert a\right\vert }\cdot \left\langle x,y\right\rangle %
\right]  \notag \\
& \leq \frac{r^{2}}{\sqrt{\left\vert a\right\vert ^{2}-r^{2}}\left(
\left\vert a\right\vert +\sqrt{\left\vert a\right\vert ^{2}-r^{2}}\right) }%
\func{Re}\left[ \frac{\bar{a}}{\left\vert a\right\vert }\cdot \left\langle
x,y\right\rangle \right]  \notag \\
& \leq \frac{r^{2}}{\sqrt{\left\vert a\right\vert ^{2}-r^{2}}\left(
\left\vert a\right\vert +\sqrt{\left\vert a\right\vert ^{2}-r^{2}}\right) }%
\left\vert \left\langle x,y\right\rangle \right\vert ;  \notag
\end{align}%
\begin{align}
\left\Vert x\right\Vert ^{2}\left\Vert y\right\Vert ^{2}& \leq \frac{1}{%
\left\vert a\right\vert ^{2}-r^{2}}\left( \func{Re}\left[ \bar{a}\cdot
\left\langle x,y\right\rangle \right] \right) ^{2}  \label{revcbs6.4} \\
& \leq \frac{\left\vert a\right\vert ^{2}}{\left\vert a\right\vert ^{2}-r^{2}%
}\left\vert \left\langle x,y\right\rangle \right\vert ^{2}  \notag
\end{align}%
and%
\begin{align}
(0& \leq )\left\Vert x\right\Vert ^{2}\left\Vert y\right\Vert
^{2}-\left\vert \left\langle x,y\right\rangle \right\vert ^{2}
\label{revcbs6.5} \\
& \leq \left\Vert x\right\Vert ^{2}\left\Vert y\right\Vert ^{2}-\left( \func{%
Re}\left[ \frac{\bar{a}}{\left\vert a\right\vert }\cdot \left\langle
x,y\right\rangle \right] \right) ^{2}  \notag \\
& \leq \frac{r^{2}}{\left\vert a\right\vert ^{2}\left( \left\vert
a\right\vert ^{2}-r^{2}\right) }-\left( \func{Re}\left[ \frac{\bar{a}}{%
\left\vert a\right\vert }\cdot \left\langle x,y\right\rangle \right] \right)
^{2}  \notag \\
& \leq \frac{r^{2}}{\left\vert a\right\vert ^{2}-r^{2}}\left\vert
\left\langle x,y\right\rangle \right\vert .  \notag
\end{align}
\end{proposition}

The proof is similar to the one in Theorem \ref{revcbst3.1}, when instead of 
$x_{i}$ we take $\left\langle x,e_{i}\right\rangle ,$ instead of $\alpha
_{i} $ we take $\left\langle e_{i},y\right\rangle ,$ $\left\Vert \cdot
\right\Vert =\left\vert \cdot \right\vert ,$ $p_{i}=1,$ and we use the
Parseval identities mentioned above in (ii) and (iii). We omit the details.

The following result may be stated as well \cite{DRA3xx}.

\begin{proposition}
\label{revcbsp6.2}Let $\left( H;\left\langle \cdot ,\cdot \right\rangle
\right) $ be a Hilbert space over $\mathbb{K}$ and $\left\{ e_{i}\right\}
_{i\in I}$ an orthornormal base for $H.$ If $x,y\in H$ $\left( y\neq
0\right) ,$ $e,\varphi ,\phi \in \mathbb{K}$ with $\func{Re}\left( \phi \bar{%
\varphi}\right) >0,$ $\left\vert e\right\vert =1$ and, either%
\begin{equation}
\left\vert \frac{\left\langle x,e_{i}\right\rangle }{\left\langle
y,e_{i}\right\rangle }-\frac{\varphi +\phi }{2}\cdot e\right\vert \leq \frac{%
1}{2}\left\vert \phi -\varphi \right\vert  \label{revcbs6.6}
\end{equation}%
or, equivalently,%
\begin{equation}
\func{Re}\left[ \left( \phi e-\frac{\left\langle x,e_{i}\right\rangle }{%
\left\langle y,e_{i}\right\rangle }\right) \left( \frac{\left\langle
e_{i},x\right\rangle }{\left\langle e_{i},y\right\rangle }-\bar{\varphi}\bar{%
e}\right) \right] \geq 0  \label{revcbs6.7}
\end{equation}%
for each $i\in I,$ then the following reverses of the Schwarz inequality
hold:%
\begin{equation}
\left\Vert x\right\Vert \left\Vert y\right\Vert \leq \frac{\func{Re}\left[
\left( \bar{\phi}+\bar{\varphi}\right) \bar{e}\left\langle x,y\right\rangle %
\right] }{2\sqrt{\func{Re}\left( \phi \bar{\varphi}\right) }}\leq \frac{1}{2}%
\cdot \frac{\left\vert \varphi +\phi \right\vert }{\sqrt{\func{Re}\left(
\phi \bar{\varphi}\right) }}\left\vert \left\langle x,y\right\rangle
\right\vert ,  \label{revcbs6.8}
\end{equation}%
\begin{align}
(0& \leq )\left\Vert x\right\Vert \left\Vert y\right\Vert -\left\vert
\left\langle x,y\right\rangle \right\vert  \label{revcbs6.9} \\
& \leq \left\Vert x\right\Vert \left\Vert y\right\Vert -\func{Re}\left[ 
\frac{\left( \bar{\phi}+\bar{\varphi}\right) \bar{e}}{\left\vert \varphi
+\phi \right\vert }\left\langle x,y\right\rangle \right] \displaybreak 
\notag \\
& \leq \frac{\left\vert \phi -\varphi \right\vert ^{2}}{2\sqrt{\func{Re}%
\left( \phi \bar{\varphi}\right) }\left( \left\vert \varphi +\phi
\right\vert +2\sqrt{\func{Re}\left( \phi \bar{\varphi}\right) }\right) } 
\notag \\
& \qquad \times \func{Re}\left[ \frac{\left( \bar{\phi}+\bar{\varphi}\right) 
\bar{e}}{\left\vert \varphi +\phi \right\vert }\left\langle x,y\right\rangle %
\right]  \notag \\
& \leq \frac{\left\vert \phi -\varphi \right\vert ^{2}}{2\sqrt{\func{Re}%
\left( \phi \bar{\varphi}\right) }\left( \left\vert \varphi +\phi
\right\vert +2\sqrt{\func{Re}\left( \phi \bar{\varphi}\right) }\right) }%
\left\vert \left\langle x,y\right\rangle \right\vert  \notag
\end{align}%
and%
\begin{align}
(0& \leq )\left\Vert x\right\Vert ^{2}\left\Vert y\right\Vert
^{2}-\left\vert \left\langle x,y\right\rangle \right\vert ^{2}
\label{revcbs6.10} \\
& \leq \left\Vert x\right\Vert ^{2}\left\Vert y\right\Vert ^{2}-\left\{ 
\func{Re}\left[ \frac{\left( \bar{\phi}+\bar{\varphi}\right) \bar{e}}{%
\left\vert \varphi +\phi \right\vert }\left\langle x,y\right\rangle \right]
\right\} ^{2}  \notag \\
& \leq \frac{\left\vert \phi -\varphi \right\vert ^{2}}{4\left\vert \phi
+\varphi \right\vert ^{2}\func{Re}\left( \phi \bar{\varphi}\right) }\left\{ 
\func{Re}\left[ \left( \bar{\phi}+\bar{\varphi}\right) \bar{e}\left\langle
x,y\right\rangle \right] \right\} ^{2}  \notag \\
& \leq \frac{\left\vert \phi -\varphi \right\vert ^{2}}{4\func{Re}\left(
\phi \bar{\varphi}\right) }\left\vert \left\langle x,y\right\rangle
\right\vert ^{2}.  \notag
\end{align}
\end{proposition}

\begin{remark}
\label{revcbsr5.3}If $\phi =M\geq m=\varphi >0,$ then one may state simpler
inequalities from (\ref{revcbs6.8}) -- (\ref{revcbs6.10}). We omit the
details.
\end{remark}

\section{Other Reverses of the CBS Inequality}

\subsection{Introduction}

Let $\left( H;\left\langle \cdot ,\cdot \right\rangle \right) $ be an inner
product space over the real or complex number field $\mathbb{K}$.

The following reverses for the Schwarz inequality hold (see \cite{SSDR1}, or
the monograph on line \cite[p. 27]{SSDxx}).

\begin{theorem}[Dragomir, 2004]
\label{revt1}Let $\left( H;\left\langle \cdot ,\cdot \right\rangle \right) $
be an inner product space over the real or complex number field $\mathbb{K}$%
. If $x,a\in H$ and $r>0$ are such that%
\begin{equation}
x\in B\left( x,r\right) :=\left\{ z\in H|\left\Vert z-a\right\Vert \leq
r\right\} ,  \label{rev1.5}
\end{equation}%
then we have the inequalities%
\begin{align}
\left( 0\leq \right) \left\Vert x\right\Vert \left\Vert a\right\Vert
-\left\vert \left\langle x,a\right\rangle \right\vert & \leq \left\Vert
x\right\Vert \left\Vert a\right\Vert -\left\vert \func{Re}\left\langle
x,a\right\rangle \right\vert  \label{rev1.6} \\
& \leq \left\Vert x\right\Vert \left\Vert a\right\Vert -\func{Re}%
\left\langle x,a\right\rangle \leq \frac{1}{2}r^{2}.  \notag
\end{align}%
The constant $\frac{1}{2}$ is best possible in (\ref{rev1.5}) in the sense
that it cannot be replaced by a smaller quantity.
\end{theorem}

An additive version for the Schwarz inequality that may be more useful in
applications is incorporated in \cite{SSDR1} (see also \cite[p. 28]{SSDxx}).

\begin{theorem}[Dragomir, 2004]
\label{revt2}Let $\left( H;\left\langle \cdot ,\cdot \right\rangle \right) $
be an inner product space over $\mathbb{K}$ and $x,y\in H$ and $\gamma
,\Gamma \in \mathbb{K}$ with $\Gamma \neq -\gamma $ and either%
\begin{equation}
\func{Re}\left\langle \Gamma y-x,x-\gamma y\right\rangle \geq 0,
\label{rev1.7}
\end{equation}%
or, equivalently,%
\begin{equation}
\left\Vert x-\frac{\gamma +\Gamma }{2}y\right\Vert \leq \frac{1}{2}%
\left\vert \Gamma -\gamma \right\vert \left\Vert y\right\Vert  \label{rev1.8}
\end{equation}%
holds. Then we have the inequalities%
\begin{align}
0& \leq \left\Vert x\right\Vert \left\Vert y\right\Vert -\left\vert
\left\langle x,y\right\rangle \right\vert  \label{rev1.9} \\
& \leq \left\Vert x\right\Vert \left\Vert y\right\Vert -\left\vert \func{Re}%
\left[ \frac{\bar{\Gamma}+\bar{\gamma}}{\left\vert \Gamma +\gamma
\right\vert }\cdot \left\langle x,y\right\rangle \right] \right\vert  \notag
\\
& \leq \left\Vert x\right\Vert \left\Vert y\right\Vert -\func{Re}\left[ 
\frac{\bar{\Gamma}+\bar{\gamma}}{\left\vert \Gamma +\gamma \right\vert }%
\cdot \left\langle x,y\right\rangle \right]  \notag \\
& \leq \frac{1}{4}\cdot \frac{\left\vert \Gamma -\gamma \right\vert ^{2}}{%
\left\vert \Gamma +\gamma \right\vert }\left\Vert y\right\Vert ^{2}.  \notag
\end{align}%
The constant $\frac{1}{4}$ in the last inequality is best possible.
\end{theorem}

We remark that a simpler version of the above result may be stated if one
assumed that the scalars are real:

\begin{corollary}
\label{revc1}If $M\geq m>0,$ and either%
\begin{equation}
\func{Re}\left\langle My-x,x-my\right\rangle \geq 0,  \label{rev1.10}
\end{equation}%
or, equivalently,%
\begin{equation}
\left\Vert x-\frac{m+M}{2}y\right\Vert \leq \frac{1}{2}\left( M-m\right)
\left\Vert y\right\Vert  \label{rev1.11}
\end{equation}%
holds, then%
\begin{align}
0& \leq \left\Vert x\right\Vert \left\Vert y\right\Vert -\left\vert
\left\langle x,y\right\rangle \right\vert  \label{rev1.12} \\
& \leq \left\Vert x\right\Vert \left\Vert y\right\Vert -\left\vert \func{Re}%
\left\langle x,y\right\rangle \right\vert  \notag \\
& \leq \left\Vert x\right\Vert \left\Vert y\right\Vert -\func{Re}%
\left\langle x,y\right\rangle  \notag \\
& \leq \frac{1}{4}\cdot \frac{\left( M-m\right) ^{2}}{M+m}\left\Vert
y\right\Vert ^{2}.  \notag
\end{align}%
The constant $\frac{1}{4}$ is sharp.
\end{corollary}

Now, let $\left( K,\left\langle \cdot ,\cdot \right\rangle \right) $ be a
Hilbert space over $\mathbb{K}$, $p_{i}\geq 0,$ $i\in \mathbb{N}$ with $%
\sum_{i=1}^{\infty }p_{i}=1.$ Consider $\ell _{\mathbf{p}}^{2}\left(
K\right) $ as the space%
\begin{equation*}
\ell _{\mathbf{p}}^{2}\left( K\right) :=\left\{ x=\left( x_{i}\right)
|x_{i}\in K,\ i\in \mathbb{N}\text{ \ and \ }\sum_{i=1}^{\infty
}p_{i}\left\Vert x_{i}\right\Vert ^{2}<\infty \right\} .
\end{equation*}%
It is well known that $\ell _{\mathbf{p}}^{2}\left( K\right) $ endowed with
the inner product%
\begin{equation*}
\left\langle x,y\right\rangle _{\mathbf{p}}:=\sum_{i=1}^{\infty
}p_{i}\left\langle x_{i},y_{i}\right\rangle
\end{equation*}%
is a Hilbert space over $\mathbb{K}$. The norm $\left\Vert \cdot \right\Vert
_{\mathbf{p}}$ of $\ell _{\mathbf{p}}^{2}\left( K\right) $ is given by%
\begin{equation*}
\left\Vert x\right\Vert _{\mathbf{p}}:=\left( \sum_{i=1}^{\infty
}p_{i}\left\Vert x_{i}\right\Vert ^{2}\right) ^{\frac{1}{2}}.
\end{equation*}%
If $x,y\in \ell _{\mathbf{p}}^{2}\left( K\right) ,$ then the following
Cauchy-Bunyakovsky-Schwarz (CBS) inequality holds true:%
\begin{equation}
\sum_{i=1}^{\infty }p_{i}\left\Vert x_{i}\right\Vert ^{2}\sum_{i=1}^{\infty
}p_{i}\left\Vert y_{i}\right\Vert ^{2}\geq \left\vert \sum_{i=1}^{\infty
}p_{i}\left\langle x_{i},y_{i}\right\rangle \right\vert ^{2}  \label{rev1.13}
\end{equation}%
with equality iff there exists a $\lambda \in \mathbb{K}$ such that $%
x_{i}=\lambda y_{i}$ for each $i\in \mathbb{N}$.

If%
\begin{equation*}
\alpha \in \ell _{\mathbf{p}}^{2}\left( K\right) :=\left\{ \alpha =\left.
\left( \alpha _{i}\right) _{i\in \mathbb{N}}\right\vert \alpha _{i}\in 
\mathbb{K},\ i\in \mathbb{N}\text{ \ and \ }\sum_{i=1}^{\infty
}p_{i}\left\vert \alpha _{i}\right\vert ^{2}<\infty \right\}
\end{equation*}%
and $x\in \ell _{\mathbf{p}}^{2}\left( K\right) ,$ then the following
(CBS)-type inequality is also valid:%
\begin{equation}
\sum_{i=1}^{\infty }p_{i}\left\vert \alpha _{i}\right\vert
^{2}\sum_{i=1}^{\infty }p_{i}\left\Vert x_{i}\right\Vert ^{2}\geq \left\Vert
\sum_{i=1}^{\infty }p_{i}\alpha _{i}x_{i}\right\Vert ^{2}  \label{rev1.14}
\end{equation}%
with equality if and only if there exists a vector $v\in K$ such that $x_{i}=%
\overline{\alpha _{i}}v$ for each $i\in \mathbb{N}$.

In \cite{DRA3xx}, by the use of some preliminary results obtained in \cite%
{SSDR0}, various reverses for the (CBS)-type inequalities (\ref{rev1.13})
and (\ref{rev1.14}) for sequences of vectors in Hilbert spaces were
obtained. Applications for bounding the distance to a finite-dimensional
subspace and in reversing the generalised triangle inequality have also been
provided.

The aim of the present section is to provide different results by employing
some inequalities discovered in \cite{SSDR1}. Similar applications are
pointed out.

\subsection{Reverses of the (CBS)-Inequality for Two Sequences in $\ell _{%
\mathbf{p}}^{2}\left( K\right) $}

The following proposition may be stated \cite{5ab}.

\begin{proposition}
\label{revp2.1}Let $x,y\in \ell _{\mathbf{p}}^{2}\left( K\right) $ and $r>0.$
If%
\begin{equation}
\left\Vert x_{i}-y_{i}\right\Vert \leq r\text{ \ for each \ }i\in \mathbb{N}%
\text{,}  \label{rev2.1}
\end{equation}%
then%
\begin{align}
(0& \leq )\left( \sum_{i=1}^{\infty }p_{i}\left\Vert x_{i}\right\Vert
^{2}\sum_{i=1}^{\infty }p_{i}\left\Vert y_{i}\right\Vert ^{2}\right) ^{\frac{%
1}{2}}-\left\vert \sum_{i=1}^{\infty }p_{i}\left\langle
x_{i},y_{i}\right\rangle \right\vert  \label{rev2.2} \\
& \leq \left( \sum_{i=1}^{\infty }p_{i}\left\Vert x_{i}\right\Vert
^{2}\sum_{i=1}^{\infty }p_{i}\left\Vert y_{i}\right\Vert ^{2}\right) ^{\frac{%
1}{2}}-\left\vert \sum_{i=1}^{\infty }p_{i}\func{Re}\left\langle
x_{i},y_{i}\right\rangle \right\vert  \notag \\
& \leq \left( \sum_{i=1}^{\infty }p_{i}\left\Vert x_{i}\right\Vert
^{2}\sum_{i=1}^{\infty }p_{i}\left\Vert y_{i}\right\Vert ^{2}\right) ^{\frac{%
1}{2}}-\sum_{i=1}^{\infty }p_{i}\func{Re}\left\langle
x_{i},y_{i}\right\rangle  \notag \\
& \leq \frac{1}{2}r^{2}.  \notag
\end{align}%
The constant $\frac{1}{2}$ in front of $r^{2}$ is best possible in the sense
that it cannot be replaced by a smaller quantity.
\end{proposition}

\begin{proof}
If (\ref{rev2.1}) holds true, then%
\begin{equation*}
\left\Vert x-y\right\Vert _{\mathbf{p}}^{2}=\sum_{i=1}^{\infty
}p_{i}\left\Vert x_{i}-y_{i}\right\Vert ^{2}\leq r^{2}\sum_{i=1}^{\infty
}p_{i}=r^{2}
\end{equation*}%
and thus $\left\Vert x-y\right\Vert _{\mathbf{p}}\leq r$.

Applying the inequality (\ref{rev1.6}) for the inner product $\left( \ell _{%
\mathbf{p}}^{2}\left( K\right) ,\left\langle \cdot ,\cdot \right\rangle _{%
\mathbf{p}}\right) ,$ we deduce the desired result (\ref{rev2.2}).

The sharpness of the constant follows by Theorem \ref{revt1} and we omit the
details.
\end{proof}

The following result may be stated as well \cite{5ab}.

\begin{proposition}
\label{revp2.2}Let $x,y\in \ell _{\mathbf{p}}^{2}\left( K\right) $ and $%
\gamma ,\Gamma \in \mathbb{K}$ with $\Gamma \neq -\gamma .$ If either%
\begin{equation}
\func{Re}\left\langle \Gamma y_{i}-x_{i},x_{i}-\gamma y_{i}\right\rangle
\geq 0\ \text{\ for each}\ \ i\in \mathbb{N}  \label{rev2.3}
\end{equation}%
or, equivalently,%
\begin{equation}
\left\Vert x_{i}-\frac{\gamma +\Gamma }{2}y_{i}\right\Vert \leq \frac{1}{2}%
\left\vert \Gamma -\gamma \right\vert \left\Vert y_{i}\right\Vert \ \text{\
for each}\ \ i\in \mathbb{N}  \label{rev2.4}
\end{equation}%
holds, then:%
\begin{align}
(0& \leq )\left( \sum_{i=1}^{\infty }p_{i}\left\Vert x_{i}\right\Vert
^{2}\sum_{i=1}^{\infty }p_{i}\left\Vert y_{i}\right\Vert ^{2}\right) ^{\frac{%
1}{2}}-\left\vert \sum_{i=1}^{\infty }p_{i}\left\langle
x_{i},y_{i}\right\rangle \right\vert  \label{rev2.5} \\
& \leq \left( \sum_{i=1}^{\infty }p_{i}\left\Vert x_{i}\right\Vert
^{2}\sum_{i=1}^{\infty }p_{i}\left\Vert y_{i}\right\Vert ^{2}\right) ^{\frac{%
1}{2}}  \notag \\
& \qquad \qquad -\left\vert \func{Re}\left[ \frac{\bar{\Gamma}+\bar{\gamma}}{%
\left\vert \Gamma +\gamma \right\vert }\sum_{i=1}^{\infty }p_{i}\left\langle
x_{i},y_{i}\right\rangle \right] \right\vert  \notag \\
& \leq \left( \sum_{i=1}^{\infty }p_{i}\left\Vert x_{i}\right\Vert
^{2}\sum_{i=1}^{\infty }p_{i}\left\Vert y_{i}\right\Vert ^{2}\right) ^{\frac{%
1}{2}}  \notag \\
& \qquad \qquad -\func{Re}\left[ \frac{\bar{\Gamma}+\bar{\gamma}}{\left\vert
\Gamma +\gamma \right\vert }\sum_{i=1}^{\infty }p_{i}\left\langle
x_{i},y_{i}\right\rangle \right]  \notag \\
& \leq \frac{1}{4}\cdot \frac{\left\vert \Gamma -\gamma \right\vert ^{2}}{%
\left\vert \Gamma +\gamma \right\vert }\sum_{i=1}^{\infty }p_{i}\left\Vert
y_{i}\right\Vert ^{2}.  \notag
\end{align}%
The constant $\frac{1}{4}$ is best possible in (\ref{rev2.5}).
\end{proposition}

\begin{proof}
Since, by (\ref{rev2.3}),%
\begin{equation*}
\func{Re}\left\langle \Gamma y-x,x-\gamma y\right\rangle _{\mathbf{p}%
}=\sum_{i=1}^{\infty }p_{i}\func{Re}\left\langle \Gamma
y_{i}-x_{i},x_{i}-\gamma y_{i}\right\rangle \geq 0,
\end{equation*}%
hence, on applying the inequality (\ref{rev1.9}) for the Hilbert space 
\newline
$\left( \ell _{\mathbf{p}}^{2}\left( K\right) ,\left\langle \cdot ,\cdot
\right\rangle _{\mathbf{p}}\right) ,$ we deduce the desired inequality (\ref%
{rev2.5}).

The best constant follows by Theorem \ref{revt2} and we omit the details.
\end{proof}

\begin{corollary}
\label{revc2.3}If the conditions (\ref{rev2.3}) and (\ref{rev2.4}) hold for $%
\Gamma =M,$ $\gamma =m$ with $M\geq m>0,$ then%
\begin{align}
(0& \leq )\left( \sum_{i=1}^{\infty }p_{i}\left\Vert x_{i}\right\Vert
^{2}\sum_{i=1}^{\infty }p_{i}\left\Vert y_{i}\right\Vert ^{2}\right) ^{\frac{%
1}{2}}-\left\vert \sum_{i=1}^{\infty }p_{i}\left\langle
x_{i},y_{i}\right\rangle \right\vert  \label{rev2.6} \\
& \leq \left( \sum_{i=1}^{\infty }p_{i}\left\Vert x_{i}\right\Vert
^{2}\sum_{i=1}^{\infty }p_{i}\left\Vert y_{i}\right\Vert ^{2}\right) ^{\frac{%
1}{2}}-\left\vert \sum_{i=1}^{\infty }p_{i}\func{Re}\left\langle
x_{i},y_{i}\right\rangle \right\vert \displaybreak  \notag \\
& \leq \left( \sum_{i=1}^{\infty }p_{i}\left\Vert x_{i}\right\Vert
^{2}\sum_{i=1}^{\infty }p_{i}\left\Vert y_{i}\right\Vert ^{2}\right) ^{\frac{%
1}{2}}-\sum_{i=1}^{\infty }p_{i}\func{Re}\left\langle
x_{i},y_{i}\right\rangle  \notag \\
& \leq \frac{1}{4}\cdot \frac{\left( M-m\right) ^{2}}{M+m}\sum_{i=1}^{\infty
}p_{i}\left\Vert y_{i}\right\Vert ^{2}.  \notag
\end{align}%
The constant $\frac{1}{4}$ is best possible.
\end{corollary}

\subsection{Reverses of the (CBS)-Inequality for Mixed Sequences}

The following result holds \cite{5ab}:

\begin{theorem}[Dragomir, 2005]
\label{revt3.1}Let $\alpha \in \ell _{\mathbf{p}}^{2}\left( K\right) ,$ $%
x\in \ell _{\mathbf{p}}^{2}\left( K\right) $ and $v\in K\backslash \left\{
0\right\} ,$ $r>0.$ If%
\begin{equation}
\left\Vert x_{i}-\overline{\alpha _{i}}v\right\Vert \leq r\left\vert \alpha
_{i}\right\vert \ \text{\ for each}\ \ i\in \mathbb{N}  \label{rev3.1}
\end{equation}%
(note that if $\alpha _{i}\neq 0$ for any $i\in \mathbb{N}$, then the
condition (\ref{rev3.1}) is equivalent to the simpler one%
\begin{equation}
\left\Vert \frac{x_{i}}{\overline{\alpha _{i}}}-v\right\Vert \leq r\ \text{\
for each}\ \ i\in \mathbb{N}),  \label{rev3.2}
\end{equation}%
then%
\begin{align}
(0& \leq )\left( \sum_{i=1}^{\infty }p_{i}\left\vert \alpha _{i}\right\vert
^{2}\sum_{i=1}^{\infty }p_{i}\left\Vert x_{i}\right\Vert ^{2}\right) ^{\frac{%
1}{2}}-\left\Vert \sum_{i=1}^{\infty }p_{i}\alpha _{i}x_{i}\right\Vert 
\label{rev3.3} \\
& \leq \left( \sum_{i=1}^{\infty }p_{i}\left\vert \alpha _{i}\right\vert
^{2}\sum_{i=1}^{\infty }p_{i}\left\Vert x_{i}\right\Vert ^{2}\right) ^{\frac{%
1}{2}}-\left\vert \left\langle \sum_{i=1}^{\infty }p_{i}\alpha _{i}x_{i},%
\frac{v}{\left\Vert v\right\Vert }\right\rangle \right\vert   \notag \\
& \leq \left( \sum_{i=1}^{\infty }p_{i}\left\vert \alpha _{i}\right\vert
^{2}\sum_{i=1}^{\infty }p_{i}\left\Vert x_{i}\right\Vert ^{2}\right) ^{\frac{%
1}{2}}-\left\vert \func{Re}\left\langle \sum_{i=1}^{\infty }p_{i}\alpha
_{i}x_{i},\frac{v}{\left\Vert v\right\Vert }\right\rangle \right\vert  
\notag \\
& \leq \left( \sum_{i=1}^{\infty }p_{i}\left\vert \alpha _{i}\right\vert
^{2}\sum_{i=1}^{\infty }p_{i}\left\Vert x_{i}\right\Vert ^{2}\right) ^{\frac{%
1}{2}}-\func{Re}\left\langle \sum_{i=1}^{\infty }p_{i}\alpha _{i}x_{i},\frac{%
v}{\left\Vert v\right\Vert }\right\rangle   \notag \\
& \leq \frac{1}{2}\cdot \frac{r^{2}}{\left\Vert v\right\Vert }%
\sum_{i=1}^{\infty }p_{i}\left\vert \alpha _{i}\right\vert ^{2}.  \notag
\end{align}%
The constant $\frac{1}{2}$ is best possible in (\ref{rev3.3}).
\end{theorem}

\begin{proof}
From (\ref{rev3.1}) we deduce%
\begin{equation*}
\left\Vert x_{i}\right\Vert ^{2}-2\func{Re}\left\langle \alpha
_{i}x_{i},v\right\rangle +\left\vert \alpha _{i}\right\vert ^{2}\left\Vert
v\right\Vert ^{2}\leq r^{2}\left\vert \alpha _{i}\right\vert ^{2},
\end{equation*}%
which is clearly equivalent to%
\begin{equation}
\left\Vert x_{i}\right\Vert ^{2}+\left\vert \alpha _{i}\right\vert
^{2}\left\Vert v\right\Vert ^{2}\leq 2\func{Re}\left\langle \alpha
_{i}x_{i},v\right\rangle +r^{2}\left\vert \alpha _{i}\right\vert ^{2}
\label{rev3.4}
\end{equation}%
for each $i\in \mathbb{N}$.

If we multiply (\ref{rev3.4}) by $p_{i}\geq 0,$ $i\in \mathbb{N}$ and sum
over $i\in \mathbb{N}$, then we deduce%
\begin{multline}
\qquad \sum_{i=1}^{\infty }p_{i}\left\Vert x_{i}\right\Vert ^{2}+\left\Vert
v\right\Vert ^{2}\sum_{i=1}^{\infty }p_{i}\left\vert \alpha _{i}\right\vert
^{2}  \label{rev3.5} \\
\leq 2\func{Re}\left\langle \sum_{i=1}^{\infty }p_{i}\alpha
_{i}x_{i},v\right\rangle +r^{2}\sum_{i=1}^{\infty }p_{i}\left\vert \alpha
_{i}\right\vert ^{2}.\qquad
\end{multline}%
Since, obviously%
\begin{multline}
\qquad 2\left\Vert v\right\Vert \left( \sum_{i=1}^{\infty }p_{i}\left\vert
\alpha _{i}\right\vert ^{2}\sum_{i=1}^{\infty }p_{i}\left\Vert
x_{i}\right\Vert ^{2}\right) ^{\frac{1}{2}}  \label{rev3.6} \\
\leq \sum_{i=1}^{\infty }p_{i}\left\Vert x_{i}\right\Vert ^{2}+\left\Vert
v\right\Vert ^{2}\sum_{i=1}^{\infty }p_{i}\left\vert \alpha _{i}\right\vert
^{2},\qquad
\end{multline}%
hence, by (\ref{rev3.5}) and (\ref{rev3.6}), we deduce%
\begin{multline*}
\qquad 2\left\Vert v\right\Vert \left( \sum_{i=1}^{\infty }p_{i}\left\vert
\alpha _{i}\right\vert ^{2}\sum_{i=1}^{\infty }p_{i}\left\Vert
x_{i}\right\Vert ^{2}\right) ^{\frac{1}{2}} \\
\leq 2\func{Re}\left\langle \sum_{i=1}^{\infty }p_{i}\alpha
_{i}x_{i},v\right\rangle +r^{2}\sum_{i=1}^{\infty }p_{i}\left\vert \alpha
_{i}\right\vert ^{2},\qquad
\end{multline*}%
which is clearly equivalent to the last inequality in (\ref{rev3.3}).

The other inequalities are obvious.

The best constant follows by Theorem \ref{revt1}.
\end{proof}

The following corollary may be stated \cite{5ab}.

\begin{corollary}
\label{revc3.1}Let $\alpha \in \ell _{\mathbf{p}}^{2}\left( K\right) ,$ $%
x\in \ell _{\mathbf{p}}^{2}\left( K\right) ,$ $e\in H,$ $\left\Vert
e\right\Vert =1$ and $\gamma ,\Gamma \in \mathbb{K}$ with $\Gamma \neq
-\gamma .$ If%
\begin{equation}
\left\Vert x_{i}-\overline{\alpha _{i}}\frac{\gamma +\Gamma }{2}\cdot
e\right\Vert \leq \frac{1}{2}\left\vert \Gamma -\gamma \right\vert
\left\vert \alpha _{i}\right\vert  \label{rev3.7}
\end{equation}%
for each $i\in \mathbb{N}$, or, equivalently,%
\begin{equation}
\func{Re}\left\langle \Gamma \overline{\alpha _{i}}e-x_{i},x_{i}-\gamma 
\overline{\alpha _{i}}e\right\rangle  \label{rev3.8}
\end{equation}%
for each $i\in \mathbb{N}$ (note that, if $\alpha _{i}\neq 0$ for any $i\in 
\mathbb{N}$, then (\ref{rev3.7}) is equivalent to%
\begin{equation}
\left\Vert \frac{x_{i}}{\overline{\alpha _{i}}}-\frac{\gamma +\Gamma }{2}%
e\right\Vert \leq \frac{1}{2}\left\vert \Gamma -\gamma \right\vert
\label{rev3.9}
\end{equation}%
for each $i\in \mathbb{N}$ and (\ref{rev3.8}) is equivalent to%
\begin{equation}
\func{Re}\left\langle \Gamma e-\frac{x_{i}}{\overline{\alpha _{i}}},\frac{%
x_{i}}{\overline{\alpha _{i}}}-\gamma e\right\rangle \geq 0  \label{rev3.10}
\end{equation}%
for each $i\in \mathbb{N}$), then the following reverse of the
(CBS)-inequality is valid:%
\begin{align}
(0& \leq )\left( \sum_{i=1}^{\infty }p_{i}\left\vert \alpha _{i}\right\vert
^{2}\sum_{i=1}^{\infty }p_{i}\left\Vert x_{i}\right\Vert ^{2}\right) ^{\frac{%
1}{2}}-\left\Vert \sum_{i=1}^{\infty }p_{i}\alpha _{i}x_{i}\right\Vert
\label{rev3.11} \\
& \leq \left( \sum_{i=1}^{\infty }p_{i}\left\vert \alpha _{i}\right\vert
^{2}\sum_{i=1}^{\infty }p_{i}\left\Vert x_{i}\right\Vert ^{2}\right) ^{\frac{%
1}{2}}-\left\vert \left\langle \sum_{i=1}^{\infty }p_{i}\alpha
_{i}x_{i},e\right\rangle \right\vert  \notag \\
& \leq \left( \sum_{i=1}^{\infty }p_{i}\left\vert \alpha _{i}\right\vert
^{2}\sum_{i=1}^{\infty }p_{i}\left\Vert x_{i}\right\Vert ^{2}\right) ^{\frac{%
1}{2}}  \notag \\
& \qquad \qquad -\left\vert \func{Re}\left[ \frac{\bar{\Gamma}+\bar{\gamma}}{%
\left\vert \Gamma +\gamma \right\vert }\left\langle \sum_{i=1}^{\infty
}p_{i}\alpha _{i}x_{i},e\right\rangle \right] \right\vert  \notag \\
& \leq \left( \sum_{i=1}^{\infty }p_{i}\left\vert \alpha _{i}\right\vert
^{2}\sum_{i=1}^{\infty }p_{i}\left\Vert x_{i}\right\Vert ^{2}\right) ^{\frac{%
1}{2}}  \notag \\
& \qquad \qquad -\func{Re}\left[ \frac{\bar{\Gamma}+\bar{\gamma}}{\left\vert
\Gamma +\gamma \right\vert }\left\langle \sum_{i=1}^{\infty }p_{i}\alpha
_{i}x_{i},e\right\rangle \right]  \notag \\
& \leq \frac{1}{4}\cdot \frac{\left\vert \Gamma -\gamma \right\vert ^{2}}{%
\left\vert \Gamma +\gamma \right\vert }\sum_{i=1}^{\infty }p_{i}\left\vert
\alpha _{i}\right\vert ^{2}.  \notag
\end{align}%
The constant $\frac{1}{4}$ is best possible.
\end{corollary}

\begin{remark}
If $M\geq m>0,$ $\alpha _{i}\neq 0$ and for $e$ as above, either%
\begin{equation}
\left\Vert \frac{x_{i}}{\overline{\alpha _{i}}}-\frac{M+m}{2}e\right\Vert
\leq \frac{1}{2}\left( M-m\right) \ \text{\ for each}\ \ i\in \mathbb{N}
\label{rev3.12}
\end{equation}%
or, equivalently,%
\begin{equation*}
\func{Re}\left\langle Me-\frac{x_{i}}{\overline{\alpha _{i}}},\frac{x_{i}}{%
\overline{\alpha _{i}}}-me\right\rangle \geq 0\ \text{\ for each}\ \ i\in 
\mathbb{N}
\end{equation*}%
holds, then%
\begin{align*}
(0& \leq )\left( \sum_{i=1}^{\infty }p_{i}\left\vert \alpha _{i}\right\vert
^{2}\sum_{i=1}^{\infty }p_{i}\left\Vert x_{i}\right\Vert ^{2}\right) ^{\frac{%
1}{2}}-\left\Vert \sum_{i=1}^{\infty }p_{i}\alpha _{i}x_{i}\right\Vert %
\displaybreak \\
& \leq \left( \sum_{i=1}^{\infty }p_{i}\left\vert \alpha _{i}\right\vert
^{2}\sum_{i=1}^{\infty }p_{i}\left\Vert x_{i}\right\Vert ^{2}\right) ^{\frac{%
1}{2}}-\left\vert \left\langle \sum_{i=1}^{\infty }p_{i}\alpha
_{i}x_{i},e\right\rangle \right\vert \\
& \leq \left( \sum_{i=1}^{\infty }p_{i}\left\vert \alpha _{i}\right\vert
^{2}\sum_{i=1}^{\infty }p_{i}\left\Vert x_{i}\right\Vert ^{2}\right) ^{\frac{%
1}{2}}-\left\vert \func{Re}\left\langle \sum_{i=1}^{\infty }p_{i}\alpha
_{i}x_{i},e\right\rangle \right\vert \\
& \leq \left( \sum_{i=1}^{\infty }p_{i}\left\vert \alpha _{i}\right\vert
^{2}\sum_{i=1}^{\infty }p_{i}\left\Vert x_{i}\right\Vert ^{2}\right) ^{\frac{%
1}{2}}-\func{Re}\left\langle \sum_{i=1}^{\infty }p_{i}\alpha
_{i}x_{i},e\right\rangle \\
& \leq \frac{1}{4}\cdot \frac{\left( M-m\right) ^{2}}{M+m}\sum_{i=1}^{\infty
}p_{i}\left\vert \alpha _{i}\right\vert ^{2}.
\end{align*}%
The constant $\frac{1}{4}$ is best possible.
\end{remark}

\subsection{Reverses for the Generalised Triangle Inequality}

In 1966, Diaz and Metcalf \cite{DMxx} proved the following interesting
reverse of the generalised triangle inequality:%
\begin{equation}
r\sum_{i=1}^{\infty }\left\Vert x_{i}\right\Vert \leq \left\Vert
\sum_{i=1}^{\infty }x_{i}\right\Vert ,  \label{rev4.1}
\end{equation}%
provided the vectors $x_{1},\dots ,x_{n}\in H\backslash \left\{ 0\right\} $
satisfy the assumption%
\begin{equation}
0\leq r\leq \frac{\func{Re}\left\langle x_{i},a\right\rangle }{\left\Vert
x_{i}\right\Vert },\qquad i\in \left\{ 1,\dots ,n\right\} ,  \label{rev4.2}
\end{equation}%
where $a\in H,$ $\left\Vert a\right\Vert =1$ and $\left( H;\left\langle
\cdot ,\cdot \right\rangle \right) $ is a real or complex inner product
space.

In an attempt to provide other sufficient conditions for (\ref{rev4.1}) to
hold, the author pointed out in \cite{SEVER1} that%
\begin{equation}
\sqrt{1-\rho ^{2}}\sum_{i=1}^{\infty }\left\Vert x_{i}\right\Vert \leq
\left\Vert \sum_{i=1}^{\infty }x_{i}\right\Vert  \label{rev4.3}
\end{equation}%
where the vectors $x_{i},$ $i\in \left\{ 1,\dots ,n\right\} $ satisfy the
condition%
\begin{equation}
\left\Vert x_{i}-a\right\Vert \leq \rho ,\qquad i\in \left\{ 1,\dots
,n\right\} ,  \label{rev4.4}
\end{equation}%
where $r\in H,$ $\left\Vert a\right\Vert =1$ and $\rho \in \left( 0,1\right)
.$

Following \cite{SEVER1}, if $M\geq m>0$ and the vectors $x_{i}\in H,$ $i\in
\left\{ 1,\dots ,n\right\} $ verify either%
\begin{equation}
\func{Re}\left\langle Ma-x_{i},x_{i}-ma\right\rangle \geq 0,\qquad i\in
\left\{ 1,\dots ,n\right\} ,  \label{rev4.5}
\end{equation}%
or, equivalently,%
\begin{equation}
\left\Vert x_{i}-\frac{M+m}{2}\cdot a\right\Vert \leq \frac{1}{2}\left(
M-m\right) ,\qquad i\in \left\{ 1,\dots ,n\right\} ,  \label{rev4.6}
\end{equation}%
where $a\in H,$ $\left\Vert a\right\Vert =1,$ then%
\begin{equation}
\frac{2\sqrt{mM}}{M+m}\sum_{i=1}^{n}\left\Vert x_{i}\right\Vert \leq
\left\Vert \sum_{i=1}^{n}x_{i}\right\Vert .  \label{rev4.7}
\end{equation}%
It is obvious from Theorem \ref{revt3.1}, that, if%
\begin{equation}
\left\Vert x_{i}-v\right\Vert \leq r,\qquad \text{for }\quad i\in \left\{
1,\dots ,n\right\} ,  \label{rev4.7.1}
\end{equation}%
where $x_{i}\in H,$ $i\in \left\{ 1,\dots ,n\right\} ,$ $v\in H\backslash
\left\{ 0\right\} $ and $r>0,$ then we can state the inequality%
\begin{align}
(0& \leq )\left( \frac{1}{n}\sum_{i=1}^{n}\left\Vert x_{i}\right\Vert
^{2}\right) ^{\frac{1}{2}}-\left\Vert \frac{1}{n}\sum_{i=1}^{n}x_{i}\right%
\Vert  \label{rev4.8} \\
& \leq \left( \frac{1}{n}\sum_{i=1}^{n}\left\Vert x_{i}\right\Vert
^{2}\right) ^{\frac{1}{2}}-\left\vert \left\langle \frac{1}{n}%
\sum_{i=1}^{n}x_{i},\frac{v}{\left\Vert v\right\Vert }\right\rangle
\right\vert  \notag \\
& \leq \left( \frac{1}{n}\sum_{i=1}^{n}\left\Vert x_{i}\right\Vert
^{2}\right) ^{\frac{1}{2}}-\left\vert \func{Re}\left\langle \frac{1}{n}%
\sum_{i=1}^{n}x_{i},\frac{v}{\left\Vert v\right\Vert }\right\rangle
\right\vert  \notag \\
& \leq \left( \frac{1}{n}\sum_{i=1}^{n}\left\Vert x_{i}\right\Vert
^{2}\right) ^{\frac{1}{2}}-\func{Re}\left\langle \frac{1}{n}%
\sum_{i=1}^{n}x_{i},\frac{v}{\left\Vert v\right\Vert }\right\rangle  \notag
\\
& \leq \frac{1}{2}\cdot \frac{r^{2}}{\left\Vert v\right\Vert }.  \notag
\end{align}%
Since, by the (CBS)-inequality we have%
\begin{equation}
\frac{1}{n}\sum_{i=1}^{n}\left\Vert x_{i}\right\Vert \leq \left( \frac{1}{n}%
\sum_{i=1}^{n}\left\Vert x_{i}\right\Vert ^{2}\right) ^{\frac{1}{2}},
\label{rev4.9}
\end{equation}%
hence, by (\ref{rev4.8}) and (\ref{rev4.5}) we have \cite{5ab}:%
\begin{equation}
(0\leq )\sum_{i=1}^{n}\left\Vert x_{i}\right\Vert -\left\Vert
\sum_{i=1}^{n}x_{i}\right\Vert \leq \frac{1}{2}n\cdot \frac{r^{2}}{%
\left\Vert v\right\Vert }  \label{rev4.9.1}
\end{equation}%
provided that (\ref{rev4.7.1}) holds true.

Utilising Corollary \ref{revc3.1}, we may state that, if 
\begin{equation}
\left\Vert x_{i}-\frac{\gamma +\Gamma }{2}\cdot e\right\Vert \leq \frac{1}{2}%
\left\vert \Gamma -\gamma \right\vert ,\qquad i\in \left\{ 1,\dots
,n\right\} ,  \label{rev4.10}
\end{equation}%
or, equivalently,%
\begin{equation}
\func{Re}\left\langle \Gamma e-x_{i},x_{i}-\gamma e\right\rangle \geq
0,\qquad i\in \left\{ 1,\dots ,n\right\} ,  \label{rev4.11}
\end{equation}%
where $e\in H,$ $\left\Vert e\right\Vert =1,\gamma ,\Gamma \in \mathbb{K}$, $%
\Gamma \neq -\gamma $ and $x_{i}\in H,$ $i\in \left\{ 1,\dots ,n\right\} ,$
then%
\begin{align}
(0& \leq )\left( \frac{1}{n}\sum_{i=1}^{n}\left\Vert x_{i}\right\Vert
^{2}\right) ^{\frac{1}{2}}-\left\Vert \frac{1}{n}\sum_{i=1}^{n}x_{i}\right%
\Vert  \label{rev4.12} \\
& \leq \left( \frac{1}{n}\sum_{i=1}^{n}\left\Vert x_{i}\right\Vert
^{2}\right) ^{\frac{1}{2}}-\left\vert \left\langle \frac{1}{n}%
\sum_{i=1}^{n}x_{i},e\right\rangle \right\vert  \notag \\
& \leq \left( \frac{1}{n}\sum_{i=1}^{n}\left\Vert x_{i}\right\Vert
^{2}\right) ^{\frac{1}{2}}-\left\vert \func{Re}\left[ \frac{\bar{\Gamma}+%
\bar{\gamma}}{\left\vert \Gamma +\gamma \right\vert }\left\langle \frac{1}{n}%
\sum_{i=1}^{n}x_{i},e\right\rangle \right] \right\vert  \notag \\
& \leq \left( \frac{1}{n}\sum_{i=1}^{n}\left\Vert x_{i}\right\Vert
^{2}\right) ^{\frac{1}{2}}-\func{Re}\left[ \frac{\bar{\Gamma}+\bar{\gamma}}{%
\left\vert \Gamma +\gamma \right\vert }\left\langle \frac{1}{n}%
\sum_{i=1}^{n}x_{i},e\right\rangle \right]  \notag \\
& \leq \frac{1}{4}\cdot \frac{\left\vert \Gamma -\gamma \right\vert ^{2}}{%
\left\vert \Gamma +\gamma \right\vert }.  \notag
\end{align}%
Now, making use of (\ref{rev4.9}) and (\ref{rev4.12}) we can establish the
following additive reverse of the generalised triangle inequality \cite{5ab}%
\begin{equation}
(0\leq )\sum_{i=1}^{n}\left\Vert x_{i}\right\Vert -\left\Vert
\sum_{i=1}^{n}x_{i}\right\Vert \leq \frac{1}{4}n\cdot \frac{\left\vert
\Gamma -\gamma \right\vert ^{2}}{\left\vert \Gamma +\gamma \right\vert },
\label{rev4.13}
\end{equation}%
provided either (\ref{rev4.10}) or (\ref{rev4.11}) hold true.

\subsection{Applications for Fourier Coefficients}

Let $\left( H;\left\langle \cdot ,\cdot \right\rangle \right) $ be a Hilbert
space over the real or complex number field $\mathbb{K}$ and $\left\{
e_{i}\right\} _{i\in I}$ an \textit{orthonormal basis} for $H.$ Then (see
for instance \cite[p. 54 -- 61]{DExx}):

\begin{enumerate}
\item[(i)] Every element $x\in H$ can be expanded in a \textit{Fourier series%
}, i.e.,%
\begin{equation*}
x=\sum_{i\in I}\left\langle x,e_{i}\right\rangle e_{i},
\end{equation*}%
where $\left\langle x,e_{i}\right\rangle ,$ $i\in I$ are the \textit{Fourier
coefficients} of $x;$

\item[(ii)] (Parseval identity)%
\begin{equation*}
\left\Vert x\right\Vert ^{2}=\sum_{i\in I}\left\langle x,e_{i}\right\rangle
e_{i},\quad x\in H;
\end{equation*}

\item[(iii)] (Extended Parseval's identity)%
\begin{equation*}
\left\langle x,y\right\rangle =\sum_{i\in I}\left\langle
x,e_{i}\right\rangle \left\langle e_{i},y\right\rangle ,\qquad x,y\in H;
\end{equation*}

\item[(iv)] (Elements are uniquely determined by their Fourier coefficients)%
\begin{equation*}
\left\langle x,e_{i}\right\rangle =\left\langle y,e_{i}\right\rangle \text{
\ \ for every }i\in I\text{ \ implies that \ }x=y.
\end{equation*}
\end{enumerate}

We must remark that all the results from the second and third sections may
be stated for $K=\mathbb{K}$ where $\mathbb{K}$ is the Hilbert space of
complex (real) numbers endowed with the usual norm and inner product.

Therefore we can state the following reverses of the Schwarz inequality \cite%
{5ab}:

\begin{proposition}
\label{revp5.1}Let $\left( H;\left\langle \cdot ,\cdot \right\rangle \right) 
$ be a Hilbert space over $\mathbb{K}$ and $\left\{ e_{i}\right\} _{i\in I}$
an orthonormal base for $H.$ If $x,y\in H,$ $y\neq 0,$ $a\in \mathbb{K}$ $%
\left( \mathbb{C},\mathbb{R}\right) $ with $r>0$ such that%
\begin{equation}
\left\vert \frac{\left\langle x,e_{i}\right\rangle }{\left\langle
y,e_{i}\right\rangle }-a\right\vert \leq r\qquad \text{for each \ }i\in I,
\label{rev5.1}
\end{equation}%
then we have the following reverse of the Schwarz inequality:%
\begin{align}
(0& \leq )\left\Vert x\right\Vert \left\Vert y\right\Vert -\left\vert
\left\langle x,y\right\rangle \right\vert  \label{rev5.2} \\
& \leq \left\Vert x\right\Vert \left\Vert y\right\Vert -\left\vert \func{Re}%
\left[ \left\langle x,y\right\rangle \cdot \frac{\bar{a}}{\left\vert
a\right\vert }\right] \right\vert  \notag \\
& \leq \left\Vert x\right\Vert \left\Vert y\right\Vert -\func{Re}\left[
\left\langle x,y\right\rangle \cdot \frac{\bar{a}}{\left\vert a\right\vert }%
\right]  \notag \\
& \leq \frac{1}{2}\cdot \frac{r^{2}}{\left\vert a\right\vert }\left\Vert
y\right\Vert ^{2}.  \notag
\end{align}%
The constant $\frac{1}{2}$ is best possible in (\ref{rev5.2}).
\end{proposition}

The proof is similar to the one in Theorem \ref{revt3.1}, where instead of $%
x_{i}$ we take $\left\langle x,e_{i}\right\rangle $, instead of $\alpha _{i}$
we take $\left\langle e_{i},y\right\rangle ,$ $\left\Vert \cdot \right\Vert
=\left\vert \cdot \right\vert ,$ $p_{i}=1$ and use the Parseval identities
mentioned above in (ii) and (iii). We omit the details.

The following result may be stated as well \cite{5ab}.

\begin{proposition}
\label{revp5.2}Let $\left( H;\left\langle \cdot ,\cdot \right\rangle \right) 
$ be a Hilbert space over $\mathbb{K}$ and $\left\{ e_{i}\right\} _{i\in I}$
an orthonormal base for $H.$ If $x,y\in H,$ $y\neq 0,$ $e,\gamma ,\Gamma \in 
\mathbb{K}$ with $\left\vert e\right\vert =1,$ $\Gamma \neq -\gamma $ and%
\begin{equation}
\left\vert \frac{\left\langle x,e_{i}\right\rangle }{\left\langle
y,e_{i}\right\rangle }-\frac{\gamma +\Gamma }{2}\cdot e\right\vert \leq 
\frac{1}{2}\left\vert \Gamma -\gamma \right\vert  \label{rev5.3}
\end{equation}%
or equivalently,%
\begin{equation}
\func{Re}\left[ \left( \Gamma e-\frac{\left\langle x,e_{i}\right\rangle }{%
\left\langle y,e_{i}\right\rangle }\right) \left( \frac{\left\langle
e_{i},x\right\rangle }{\left\langle e_{i},y\right\rangle }-\bar{\gamma}\bar{e%
}\right) \right] \geq 0  \label{rev5.4}
\end{equation}%
for each $i\in I,$ then%
\begin{align}
(0& \leq )\left\Vert x\right\Vert \left\Vert y\right\Vert -\left\vert
\left\langle x,y\right\rangle \right\vert  \label{rev5.5} \\
& \leq \left\Vert x\right\Vert \left\Vert y\right\Vert -\left\vert \func{Re}%
\left[ \frac{\bar{\Gamma}+\bar{\gamma}}{\left\vert \Gamma +\gamma
\right\vert }\left\langle x,y\right\rangle \cdot \bar{e}\right] \right\vert 
\notag \\
& \leq \left\Vert x\right\Vert \left\Vert y\right\Vert -\func{Re}\left[ 
\frac{\bar{\Gamma}+\bar{\gamma}}{\left\vert \Gamma +\gamma \right\vert }%
\left\langle x,y\right\rangle \cdot \bar{e}\right]  \notag \\
& \leq \frac{1}{4}\cdot \frac{\left\vert \Gamma -\gamma \right\vert ^{2}}{%
\left\vert \Gamma +\gamma \right\vert }\left\Vert y\right\Vert ^{2}.  \notag
\end{align}%
The constant $\frac{1}{4}$ is best possible.
\end{proposition}

\begin{remark}
If $\Gamma =M\geq m=\gamma >0,$ then one may state simpler inequalities from
(\ref{rev5.5}). We omit the details.
\end{remark}

%

\end{document}